\documentclass[12pt,a4paper, reqno]{amsart}
\usepackage{amscd}
\usepackage{amssymb}
\usepackage[centering,text={14.5cm,22cm}]{geometry}
\usepackage{graphicx}
\usepackage{color}
\usepackage[all]{xy}
\usepackage{mathrsfs}
\usepackage{stmaryrd}
\usepackage{srcltx}
\usepackage{hyperref}
\usepackage{url}
\usepackage{enumerate}

\definecolor{shadecolor}{rgb}{1,0.9,0.7}

\usepackage[english]{babel}

\usepackage{amsmath, amsthm}
\usepackage{textcomp}
\usepackage{amscd}
\usepackage{amsfonts}
\usepackage{amssymb}
\usepackage{tikz-cd}
\usepackage{mathtools}
\usepackage{graphicx}
\usepackage{epigraph}
\usepackage{underscore}
\usepackage{color}
\usepackage{mathdots}
\usepackage{graphicx,wrapfig,lipsum}
\usepackage{pbox}
\usepackage{booktabs}
\usepackage{longtable}
\usepackage{bm}
\usepackage{csquotes}

\usepackage{hyperref}
\hypersetup{
  colorlinks   = true,          
  urlcolor     = violet,          
  linkcolor    = blue,          
  citecolor   = violet             
}
\hypersetup{linktocpage}

\usepackage{tikz}
\usetikzlibrary{external}
\usetikzlibrary{calc}

\usepackage{caption}
\usepackage{float}

\usepackage{pgfplots}
\pgfplotsset{compat=1.15}
\usepackage{mathrsfs}
\usetikzlibrary{arrows}

\usepackage{mathabx}

\newtheorem{theorem}{Theorem}[section]
\newtheorem{lemma}[theorem]{Lemma}
\newtheorem{fact}[theorem]{Proposition}
\newtheorem{cor}[theorem]{Corollary}
\newtheorem{conj}[theorem]{Conjecture}

\newtheorem*{theoremA}{Theorem A}
\newtheorem*{theoremB}{Theorem B}
\newtheorem*{theoremC}{Theorem C}

\theoremstyle{definition}
\newtheorem{defn}[theorem]{Definition}
\newtheorem{ass}[theorem]{Assumption}
\newtheorem{cons}[theorem]{Construction}

\newtheorem{ex}[theorem]{Example}

\newtheorem{basicsetup}[theorem]{Basic Setup}

\theoremstyle{remark}
\newtheorem{remark}[theorem]{Remark}
\newtheorem{remarks}[theorem]{Remarks}
\newtheorem{notation}[theorem]{Notation}
\newtheorem{notations}[theorem]{Notations}
\newtheorem{obs}[theorem]{Observations}
\newtheorem{observe}[theorem]{Observation}

\numberwithin{equation}{section}
\numberwithin{figure}{section}

\usepackage{pgfplots}
\pgfplotsset{compat=1.15}
\usepackage{mathrsfs}
\usetikzlibrary{arrows}

\newcommand {\kk} {\Bbbk}

\newcommand\btau{{\boldsymbol{\tau}}}

\newcommand\bsigma{{\boldsymbol{\sigma}}}

\DeclareMathOperator{\Spec}{Spec}

\DeclareMathOperator{\coker}{coker}

\DeclareMathOperator{\tors}{tors}

\DeclareMathOperator{\Conv}{Conv}

\DeclareMathOperator{\out}{out}

\DeclareMathOperator{\reduced}{red}

\DeclareMathOperator{\tot}{tot}

\DeclareMathOperator{\trop}{trop}

\DeclareMathOperator{\Proj}{Proj}

\DeclareMathOperator{\Pic}{Pic}

\DeclareMathOperator{\Sing}{Sing}

\DeclareMathOperator{\sat}{sat}

\DeclareMathOperator{\Spf}{Spf}

\DeclareMathOperator{\gp}{gp}

\DeclareMathOperator{\GS}{GS}

\DeclareMathOperator{\TD}{TD}

\DeclareMathOperator{\IMS}{IMS}

\DeclareMathOperator{\speci}{sp}

\DeclareMathOperator{\Id}{Id}

\DeclareMathOperator{\im}{im}

\DeclareMathOperator{\coar}{coar}

\DeclareMathOperator{\Int}{Int}

\DeclareMathOperator{\Aut}{Aut}

\DeclareMathOperator{\Nef}{Nef}

\DeclareMathOperator{\virt}{virt}

\DeclareMathOperator{\Hom}{Hom}

\DeclareMathOperator{\new}{new}

\DeclareMathOperator{\Div}{Div}

\DeclareMathOperator{\md}{mod}

\DeclareMathOperator{\rel}{rel}

\DeclareMathOperator{\dist}{dist}

\DeclareMathOperator{\ext}{ext}

\DeclareMathOperator{\MPA}{MPA}

\DeclareMathOperator{\Aff}{Aff}

\DeclareMathOperator{\numer}{num}

\DeclareMathOperator{\snc}{snc}

\DeclareMathOperator{\alg}{alg}

\DeclareMathOperator{\Rel}{Rel}

\DeclareMathOperator{\irrel}{irrel}

\DeclareMathOperator{\univ}{univ}

\DeclareMathOperator{\Griff}{Griff}

\DeclareMathOperator{\divisor}{div}

\DeclareMathOperator{\Cl}{Cl}

\DeclareMathOperator{\ind}{ind}

\DeclareMathOperator{\comp}{comp}

\DeclareMathOperator{\supp}{supp}

\DeclareMathOperator{\bMPA}{M\breve{P}A}

\DeclareMathOperator{\Num}{Num}

\usepackage{pstricks}

\usepackage{pstricks-add}

\mathtoolsset{showonlyrefs}

\definecolor{qqqqff}{rgb}{0,0,1}
\definecolor{ududff}{rgb}{0.30196078431372547,0.30196078431372547,1}
\definecolor{xdxdff}{rgb}{0.49019607843137253,0.49019607843137253,1}

\begin{document}
	
\title[Mirrors to toric degenerations via intrinsic mirror symmetry]
{Mirrors to toric degenerations via \\[0.05cm] intrinsic mirror symmetry}

\author[Evgeny Goncharov]{Evgeny Goncharov}

\dedicatory{\vspace{0.2cm} \normalfont\normalsize 7 August 2026}

\thanks{DPMMS, Centre for Mathematical Sciences,
Wilberforce Road, Cambridge, CB3 0WB, UK.
Email: \href{mailto:eg555@maths.cam.ac.uk}
{\texttt{eg555@maths.cam.ac.uk}}.}

\thanks{This article is a revised version of the author's PhD thesis at the University of Cambridge.}

\date{}

\maketitle

\vspace{-0.5cm}
	
\begin{abstract}
We explore the connection between two mirror constructions in Gross-Siebert mirror symmetry: toric degeneration mirror symmetry \cite{G1, GS2, models, GS1} and intrinsic mirror symmetry \cite{GS3, GS4}. After briefly exploring the case of degenerations of elliptic curves, we show that the Gross-Siebert mirror construction for minimal relative log Calabi-Yau degenerations generalizes that for divisorial toric degenerations $\bar{\mathfrak{X}} \to \mathcal{S}$ of K3-s that have a smooth generic fibre. We achieve this by constructing a resolution of $\bar{\mathfrak{X}} \to \mathcal{S}$ to a relative minimal log Calabi-Yau degeneration $\mathfrak{X} \to \mathcal{S}$  and comparing the algorithmic scattering diagram $\bar{\mathfrak{D}}$ giving rise to the toric degeneration mirror $\check{\bar{\mathfrak{X}}}$ and the canonical scattering diagram $\mathfrak{D}$ giving rise to the intrinsic mirror $\check{\mathfrak{X}}$. Moreover, we vastly expand the construction and obtain a correspondence between the restriction of the intrinsic mirror to the (numerical) minimal relative Gross-Siebert locus and the universal toric degeneration mirror. We also discuss generalizing the results to higher dimensions. In particular, we construct log smooth resolutions for a natural family of toric degenerations of Calabi-Yau threefolds.
\end{abstract}
	
\tableofcontents



\section{Introduction} \label{intro}

This project arose from the desire to understand the connection between various mirror symmetry constructions in the Gross-Siebert program. More precisely, we were interested in establishing a correspondence between the recently developed intrinsic mirror symmetry of \cite{GS3, GS4}, the toric degeneration mirror symmetry of \cite{G1, GS2, models, GS1}, and the classical Batyrev (hypersurfaces, \cite{B}) and Batyrev-Borisov (complete intersections, \cite{BB}) mirror symmetry.

\medskip

\noindent \textbf{Batyrev mirrors.} Batyrev and Batyrev-Borisov mirror symmetry are the most famous classical mirror symmetry constructions. They provide many examples of mirrors and incorporate the initial mirror symmetry observations (e.g. the construction of Greene and Plesser \cite{GP}).

Let $\Delta \subseteq \mathbb{R}^n$ be an $n$-dimensional lattice polytope with $0 \in \Int \Delta$. Then the \emph{dual} of $\Delta$ is defined as  $$\Delta^* = \left\{ x \in (\mathbb{R}^n)^* \ | \ \langle x, \Delta \rangle \geq -1 \right\}$$ and we also have $0 \in \Int \Delta^*$. We say that $\Delta$ is \emph{reflexive} if $\Delta^*$ is also a lattice polytope. One can check that $(\Delta^*)^* = \Delta$ so reflexive polytopes come in pairs $(\Delta, \Delta^*)$ (we say that a polytope $\Delta$ is \emph{self-dual} if $\Delta \cong \Delta^*$ as lattice polytopes). 

\begin{sloppypar}
Following \cite{B}, one can construct Batyrev mirrors as follows. A pair $(\Delta, \Delta^*)$ of reflexive polytopes gives rise to projective toric varieties $(\mathbb{P}_{\Delta}, \mathcal{O}_{\mathbb{P}_{\Delta}}(1))$, $(\mathbb{P}_{\Delta^*}, \mathcal{O}_{\mathbb{P}_{\Delta^*}}(1))$. Let $s \in \mathcal{O}_{\mathbb{P}_{\Delta}}(1), \ s^* \in \mathcal{O}_{\mathbb{P}_{\Delta^*}}(1)$ be two sections and define $$X_s := Z(s) \subseteq \mathbb{P}_{\Delta}, \quad X_{s^*} := Z(s^*) \subseteq \mathbb{P}_{\Delta^*}.$$ Then for a general ($\Delta$-\emph{regular} in the sense of \cite[Definition 3.1.1]{B}) choice of $s$ and $s^*$, we call $X_s$ and $X_{s^*}$ \emph{dual Batyrev mirrors}. In particular, they possess a duality of stringy Hodge numbers.
\end{sloppypar}

The setup for Batyrev-Borisov mirror symmetry of complete intersections is similar, but the combinatorics is more involved, see \cite{BB, G2} for details.

\medskip

\noindent \textbf{Toric degeneration mirrors.} Toric degenerations of Calabi-Yau varieties are a natural class to discuss mirror symmetry \cite{G1, GS2, models, GS1} and have been an object of study of the Gross-Siebert approach to mirror symmetry for the last few decades. Toric degeneration mirrors are double dual in the sense that for any toric degeneration $\bar{\mathfrak{X}} \to \mathcal{S}$ obtained via the reconstruction algorithm of \cite{GS1}, taking the mirror to the mirror degeneration (i.e. performing the reconstruction algorithm twice) gives back $\bar{\mathfrak{X}} \to \mathcal{S}$. Toric degeneration mirrors generalize most of the classical mirror symmetry constructions. In particular, they generalize the Batyrev mirror symmetry of hypersurfaces \cite{B} (and the Batyrev-Borisov mirror symmetry of complete intersections \cite{BB}) as we explain below (see \cite{G2} for details).

Let $\kk$ be an algebraically closed field of characteristic $0$ and let $\mathcal{S}:= \Spec R$ where $R$ is a complete\footnote{We only require completeness of $R$ for convenience of notation. For any non-complete $R$ and a toric degeneration $\bar{\mathfrak{X}} \to \Spec R$, the basechange by the canonical map $R \to \hat{R}$ to the completion defines a toric degeneration $\bar{\mathfrak{X}} \to \Spec \hat{R}$. Since the mirror toric degeneration only depends on the central fibre $\bar{\mathfrak{X}}_0$, replacing $\bar{\mathfrak{X}} \to \Spec R$ by $\bar{\mathfrak{X}} \to \Spec \hat{R}$ does not change the mirror. We shall be performing blowups of toric degenerations, but the blowup locus will always be supported on the central fibre as well. Therefore, the results of this thesis hold for non-complete $R$ by replacing relevant mentions of $\bar{\mathfrak{X}} \to \Spec R$ with mentions of the basechange $\bar{\mathfrak{X}} \to \Spec \hat{R}$.} discrete valuation $\kk$-algebra. 

\begin{remark} \label{dvrremark}
The reader may think of $R = \kk\llbracket t \rrbracket$, but we do not wish to restrict to this case. Note that for any complete discrete valuation $\kk$-algebra $R$, choosing a uniformizing parameter for $R$ defines a map $\kk\llbracket t \rrbracket \to R$. 
\end{remark}

Let $0 \in \mathcal{S}$ denote the closed point. Roughly, a \emph{toric degeneration} is a flat proper morphism $\bar{g}: \bar{\mathfrak{X}} \to \mathcal{S}$ of varieties such that the central (or \enquote{special}) fibre $\bar{\mathfrak{X}}_0$ is a union of toric varieties meeting along toric strata. Moreover, away from the \emph{singular} (or \emph{discriminant}) \emph{locus} $Z \subseteq \bar{\mathfrak{X}}$ of codimension $2$ and not containing any toric strata, $\bar{\mathfrak{X}} \to \mathcal{S}$ is described étale locally as a monomial map on an affine toric variety. We copy the precise technical definition from \cite[Definition 4.1]{GS2}. 

\begin{defn} \label{toricdegeneration}
A \emph{toric degeneration of Calabi-Yau
varieties} over $\mathcal{S}:=\Spec R$ is a normal algebraic space $\bar{\mathfrak{X}}$ along with a flat proper morphism
$\bar{\mathfrak{X}} \to \mathcal{S}$ satisfying the following properties:
\begin{enumerate}
    \item The generic fibre $\bar{\mathfrak{X}}_{\eta}$ is an irreducible normal variety over $\eta$.
    \item If $\nu:\tilde{\bar{\mathfrak{X}}}_0 \to \bar{\mathfrak{X}}_0$ is the normalization,
    then $\tilde{\bar{\mathfrak{X}}}_0$ is a disjoint union of toric varieties,
    the conductor scheme $C\subseteq\tilde{\bar{\mathfrak{X}}}_0$ is reduced
    and the map $C\to\nu(C)$ is unramified and generically
    two-to-one. The square
    \begin{equation} 
    \begin{tikzcd}
    {C} \arrow{r} \arrow{d} & \tilde{\bar{\mathfrak{X}}}_0  \arrow{d}{\nu} \\
    \nu(C) \arrow{r}{} & \bar{\mathfrak{X}}_0
    \end{tikzcd}
    \end{equation}
    is cartesian and cocartesian.
    \item $\bar{\mathfrak{X}}_0$ is Gorenstein, and the conductor locus $C$
    restricted to each irreducible component of $\tilde{\bar{\mathfrak{X}}}_0$ is the union
    of all toric Weil divisors. 
    \item There exists a closed subset $Z\subseteq \bar{\mathfrak{X}}$ of relative codimension $\ge 2$ such that $Z$ satisfies the following properties:
    \begin{enumerate}
        \item $Z$ does not contain the image under $\nu$ of any toric stratum of $\tilde{\bar{\mathfrak{X}}}_0$.
        \item For any geometric point $\bar x \to \bar{\mathfrak{X}} \setminus Z$,
        there is an étale neighbourhood $U_{\bar x}\to \bar{\mathfrak{X}} \setminus Z$ of
        $\bar x$, an affine toric variety $Y_{\bar x}$, a regular function
        $f_{\bar x}$ on $Y_{\bar x}$ given by a monomial, a choice of
        uniformizing parameter for $R$ giving a map $\kk[\mathbb{N}]\rightarrow R$, and
        a commutative diagram
        \begin{equation} 
        \begin{tikzcd}
        {U_{\bar x}} \arrow{r} \arrow{d}{f|_{U_{\bar x}}} & Y_{\bar x}  \arrow{d}{f_{\bar x}} \\
         \Spec R \arrow{r}{} & \Spec \kk[\mathbb{N}]
        \end{tikzcd}
        \end{equation}
        such that the induced map $U_{\bar x}\to \Spec R \times_{\Spec \kk[\mathbb{N}]}
        Y_{\bar x}$ is smooth. Furthermore, $f_{\bar x}$ vanishes on each
        toric divisor of $Y_{\bar x}$.
    \end{enumerate}
\end{enumerate}
\end{defn}

The central fibre $\bar{\mathfrak{X}}_0$ of a toric degeneration $\bar{\mathfrak{X}} \to \mathcal{S}$ is a union of toric varieties $\bar{D}_i, \ 1 \leq i \leq \bar{m}$. We let $\bar{D} := \bar{\mathfrak{X}}_0 = \bar{D}_1 + \dots + \bar{D}_{\bar{m}}$. We shall always endow $\bar{\mathfrak{X}} \to \mathcal{S}$ with \emph{divisorial log structures} $\mathcal{M}_{\bar{\mathfrak{X}}}$ on $\bar{\mathfrak{X}}$ and $\mathcal{M}_{\mathcal{S}}$ on $\mathcal{S}$ with divisors $\bar{D}$ and $0$ respectively. This makes $(\bar{\mathfrak{X}}, D) \to (\mathcal{S}, 0)$ a log morphism that we usually just denote $\bar{\mathfrak{X}} \to \mathcal{S}$. Moreover, condition (4)(b) of Definition \ref{toricdegeneration} implies that the log structure $\mathcal{M}_{\bar{\mathfrak{X}}}$ on $\bar{\mathfrak{X}}$ is \emph{fine and saturated} away from $Z$\footnote{That is, the stalks $\overline{\mathcal{M}}_{\bar{\mathfrak{X}}, \bar{x}}$ of the ghost sheaf $\overline{\mathcal{M}}_{\bar{\mathfrak{X}}}:=\mathcal{M}_{\bar{\mathfrak{X}}}/\mathcal{M}_{\bar{\mathfrak{X}}}^{\times}$ of the log structure $\mathcal{M}_{\bar{\mathfrak{X}}}$ are finitely generated, integral, saturated monoids for all geometric points $\bar x \to \bar{\mathfrak{X}} \setminus Z$. Note that by the definition of the ghost sheaf $\overline{\mathcal{M}}_{\bar{\mathfrak{X}}, \bar{x}}^{\times}$ is trivial so the stalks $\overline{\mathcal{M}}_{\bar{\mathfrak{X}}, \bar{x}}$ are actually toric and sharp monoids.} and that $\bar{\mathfrak{X}} \to \mathcal{S}$ is \emph{log smooth} away from $Z$.\footnote{These claims follow from the fact that the $Y_{\bar{x}}$ are toric varieties and from Kato's criterion (see \cite[Theorem (3.5)]{Kato}) respectively. See \cite[Proposition 4.6]{GS2} (using \cite[Definition 4.3]{GS2}) for details.} We shall always work with toric degenerations satisfying the following additional assumption (and just refer to them as toric degenerations).

\begin{ass} \label{toricdegenerationass}
Let $\bar{\mathfrak{X}} \to \mathcal{S}$ be a toric degeneration. We assume, in addition, that:
\begin{enumerate}
    \item $\bar{\mathfrak{X}}$ is a variety (not just an algebraic space).
    \item The morphism $\bar{\mathfrak{X}} \to \mathcal{S}$ is projective.
    \item $(\bar{\mathfrak{X}}, \bar{D})$ is a Zariski log scheme (it is enough to require that the irreducible components of $\bar{\mathfrak{X}}_0$ are normal) that is a \emph{simple log scheme} in the sense of \cite[Definition 2.1]{ACGS}.\footnote{In the construction of tropicalization in \cite[Section 2.1.4]{ACGS} (used to define simple log schemes) the authors assume that the log structure $\mathcal{M}_{\bar{\mathfrak{X}}}$ on $\bar{\mathfrak{X}}$ is fine and saturated. However, it is easy to extend the construction to the case of toric degenerations, see Section \ref{tropicalization}. The tropicalization of $\bar{\mathfrak{X}}$ naturally coincides with the tropicalization of $\bar{\mathfrak{X}} \setminus Z$ (see Remark \ref{tropremarks}(1)), so equivalently, we can require that $\bar{\mathfrak{X}} \setminus Z$ is a simple log scheme.}
    \item Any (possibly empty) intersection of the irreducible components of $\bar{D}$ is connected.
\end{enumerate}
\end{ass}

The conditions in Assumption \ref{toricdegenerationass} are necessary to fit with the setup for intrinsic mirrors (and log smooth degenerations) \cite{GS3, GS4} as well as the general framework for scattering diagrams of \cite{GHS}. As we explain in Section \ref{relax}, conditions (1), (3), and (4) are not critical, but removing them would require a significant revision of \cite{GS3, GS4, GHS}. Throughout the thesis, we shall keep in mind the following example.
	
\begin{ex} \label{runningexample}
Let $\bar{\mathfrak{X}}$ be defined as $$\bar{\mathfrak{X}}:=\left\{tf_4+x_0x_1x_2x_3\right\} \subseteq \mathbb{P}^3 \times \Spec \kk\llbracket t \rrbracket$$ where $f_4$ is a general homogeneous quartic polynomial, and let $\bar{\mathfrak{X}} \to \Spec \kk\llbracket t \rrbracket$ be the natural projection. Then $\bar{\mathfrak{X}}_0$ is the coordinate tetrahedron of four $\mathbb{P}^2$-s intersecting in $\mathbb{P}^1$-s. Let $$Z := \left\{t = f_4 = 0\right\} \cap \Sing\left(\bar{\mathfrak{X}}_0\right)$$ be the singular locus. Then $Z$ consists of $24$ \emph{(log) singular points}, $4$ on each $\mathbb{P}^1$, and it is easy to see that $\bar{\mathfrak{X}} \to \Spec \kk\llbracket t \rrbracket$ is a toric degeneration of K3-s satisfying Assumption \ref{toricdegenerationass}.
\end{ex}

The mirror to a toric degeneration $\bar{\mathfrak{X}} \to \mathcal{S}$ is obtained using the \emph{dual intersection complex} $\left(\bar B, \bar{\mathscr{P}}\right)$, which is a topological manifold $B$ with a natural polyhedral decomposition $\bar{\mathscr{P}}$. Conditions (3) and (4) of Assumption \ref{toricdegenerationass} ensure that $\left(\bar B, \bar{\mathscr{P}}\right)$ is a genuine polyhedral complex rather than a \emph{generalized polyhedral complex} of \cite[Section 2.1.3]{ACGS} (or, in the language of \cite{GS2}, a manifold with a \emph{toric polyhedral decomposition}, see \cite[Definitions 1.21 and 1.22]{GS2}). In particular, the cells of $\bar{\mathscr{P}}$ don't self-intersect, and an intersection of two cells of $\bar{\mathscr{P}}$ is a cell of $\bar{\mathscr{P}}$. As a result, $\left(\bar B, \bar{\mathscr{P}}\right)$ satisfies the conditions of \cite[Construction 1.1.1]{GHS}, and we will be able to use the conventions of \cite{GHS} (instead of the rather elaborate language of \cite{GS2, GS1}) to work with $\left(\bar B, \bar{\mathscr{P}}\right)$. The dual intersection complex $\left(\bar B, \bar{\mathscr{P}}\right)$ depends only  on the central fibre $\bar{\mathfrak{X}}_0$.

A toric degeneration $\bar{\mathfrak{X}} \to \mathcal{S}$ also gives rise to the structure of an \emph{integral affine manifold with singularities} (in codimension $\geq 2$) on $\left(\bar B, \bar{\mathscr{P}}\right)$. The \emph{toric degeneration mirror} to $\bar{\mathfrak{X}} \to \mathcal{S}$ is a toric degeneration $\check{\bar{\mathfrak{X}}} \to \Spec \kk\llbracket t \rrbracket$ that arises from $\left(\bar B, \bar{\mathscr{P}}\right)$ by an algorithmic procedure (the \enquote{reconstruction algorithm} of \cite{GS1}) encoded by an \emph{algorithmic scattering diagram} $\bar{\mathfrak{D}}$ on $\left(\bar B, \bar{\mathscr{P}}\right)$. We note that one can get a family over $\mathcal{S}$ by basechanging $\check{\bar{\mathfrak{X}}} \to \Spec \kk\llbracket t \rrbracket$ by the map of Remark \ref{dvrremark}. It will be more convenient to work with the universal $\check{\bar{\mathfrak{X}}} \to \Spec \kk\llbracket t \rrbracket$ but one can easily formulate similar results for $\check{\bar{\mathfrak{X}}} \to \mathcal{S}$ by applying the additional basechange. The construction of the toric degeneration mirror requires and depends on a choice of a relatively ample divisor $A$ on $\bar{\mathfrak{X}}$ that we call the \emph{polarization} of $\bar{\mathfrak{X}}$.

Following \cite{G2}, one can construct \emph{Batyrev degenerations} as follows. As before, we take sufficiently general sections $s \in \mathcal{O}_{\mathbb{P}_{\Delta}}(1), \ s^* \in \mathcal{O}_{\mathbb{P}_{\Delta^*}}(1)$ and let $s_0$ be given by $0 \in \Int \Delta$, $s^*_0$ be given by $0 \in \Int \Delta^*$. Then one can define\footnote{The construction works for any discrete valuation $\kk$-algebra $R$, but we choose to work with $R:=\kk\llbracket t \rrbracket$ for a direct comparison with the other mirror constructions. As before, one may get a family over any $R$ by basechanging via the map of Remark \ref{dvrremark}.} 
$$  \bar{\mathfrak{X}}_{\Delta} : = \left\{ ts+s_0 = 0 \right\} \subseteq \mathbb{P}_{\Delta} \times \Spec \kk\llbracket t \rrbracket, \quad \bar{\mathfrak{X}}_{\Delta^*} : = \left\{ ts^*+s^*_0 = 0 \right\} \subseteq \mathbb{P}_{\Delta^*} \times \Spec \kk\llbracket t \rrbracket.$$
Along with the natural projections to $\Spec \kk\llbracket t \rrbracket$, these give a pair of dual Batyrev degenerations. Note that the toric degeneration of Example \ref{runningexample} is a Batyrev degeneration by taking $\Delta$ to be the standard polytope for $\mathbb{P}^3$.
By \cite[Proposition 2.16]{G2} Batyrev degenerations are toric degenerations, and one can show that their dual intersection complexes are related by a \emph{discrete Legendre transform} (which is a purely combinatorial construction, see \cite[Section 1.4]{GS2}). It is possible to generalize the construction to other ample line bundles and to complete intersections, but one needs to take certain partial resolutions of $\mathbb{P}_{\Delta}$ and $\mathbb{P}_{\Delta^*}$ to achieve duality, see \cite[Theorem 3.25]{G2}. In this sense, toric degeneration mirror symmetry generalizes Batyrev (and Batyrev-Borisov) mirror symmetry. 

\medskip

\noindent \textbf{Intrinsic mirrors.} More recently in \cite{GHK} the authors have constructed a mirror to an arbitrary \emph{Looijenga pair} $(\mathfrak{X}, D)$, that is a smooth rational projective surface $\mathfrak{X}$ with $D \in |-K_{\mathfrak{X}}|$ a singular nodal curve. These assumptions ensure that $\mathfrak{X}  \setminus  D$ is a \emph{log Calabi-Yau (CY)} surface.\footnote{Conversely, all the minimal log CY surfaces $(\mathfrak{X}, D)$ in the sense of Gross-Siebert mirror symmetry are Looijenga pairs.} The mirror is a formal family $\check{\mathfrak{X}} \to \Spf{\widehat{\kk[P]}}$\footnote{Here $\widehat{\kk[P]}$ is the completion with respect to the maximal ideal of $\kk[P]$.} where $NE(\mathfrak{X}) \subseteq P \subseteq A_1(\mathfrak{X}, \mathbb{Z})$ is any finitely generated, saturated and sharp (having no non-trivial invertible elements) submonoid of the finitely generated free abelian group $A_1(\mathfrak{X}, \mathbb{Z})$ of $1$-cycles modulo algebraic or numerical equivalence (we sometimes write $A_1(\mathfrak{X}, \mathbb{Z})_{\alg}$ and $A_1(\mathfrak{X}, \mathbb{Z})_{\numer}$ respectively, note that the two groups coincide for a Looijenga pair) containing the submonoid $NE(\mathfrak{X})$ generated by the effective curve classes. 

The construction has been generalized to arbitrary log CY varieties $(\mathfrak{X}, D)$ (the \emph{absolute} setup) and to the \emph{relative} setup of a projective log smooth morphism $g: \mathfrak{X} \to \mathcal{S}$ where $\mathfrak{X}$ carries a Zariski log structure and $\mathcal{S}$ is a regular one-dimensional scheme over $\Spec \kk$ with a divisorial log structure coming from a single closed point $0 \in \mathcal{S}$ (in particular, $\mathcal{S}$ may be the spectrum of a discrete valuation $\kk$-algebra $R$) in \cite{GS3, GS4}. 

In general, $A_1(\mathfrak{X}, \mathbb{Z})$ is not a free group, and one needs to modify the conditions on $P$ slightly. One no longer requires that $P$ is sharp but that the group $P^{\times}$ of the invertible elements of $P$ coincides with the torsion part of $A_1(\mathfrak{X}, \mathbb{Z})$. Note that $A_1(\mathfrak{X}, \mathbb{Z})_{\numer}$ is always free. In the relative setup, one has similar assumptions on $P$ but should use a different group of curve classes that we denote $A_1(\mathfrak{X}_0, \mathbb{Z})$. This may be the group of $1$-cycles on $\mathfrak{X}_0$ or the group of relative $1$-cycles for $\mathfrak{X} \to \mathcal{S}$ modulo algebraic or numerical equivalence.

We shall always work in the relative setup. Similarly to the case of toric degenerations, let $D_i, \ 1 \leq i \leq m$ be the components of the reduced central fibre $\left(\mathfrak{X}_0 \right)_{\reduced}$\footnote{We will almost always deal with log smooth morphisms $\mathfrak{X} \to \mathcal{S}$ with $\mathfrak{X}_0$ reduced.} and let $D:= \left(\mathfrak{X}_0 \right)_{\reduced} = D_1 + \dots + D_m$. We assume that the log structure $\mathcal{M}_{\mathfrak{X}}$ on $\mathfrak{X}$ is divisorial with divisor $D$.\footnote{In \cite{GS4} one needs to require that the divisor $D$ giving the divisorial log structure on $\mathfrak{X}$ is simple normal crossings. However, this technical assumption can be removed (see Section \ref{technicalasssection}).} Note that the underlying morphism of schemes is a flat morphism\footnote{Indeed, a log smooth morphism is log flat, see \cite[Part IV, Proposition 4.1.2(1)]{Ogus}. Further, a log flat and \emph{integral} morphism is flat by \cite[Part IV, Proposition 4.3.5(1)]{Ogus}. But $\mathcal{S}$ is one-dimensional so $\mathfrak{X} \to \mathcal{S}$ is always integral by \cite[Part III, Proposition 2.5.3(3)]{Ogus}.} and the mirror extends to an algebraic family $\check{\mathfrak{X}} \to \Spec{\widehat{\kk[P]}}$.\footnote{This is true as long as $D = g^{-1}(0)$ set-theoretically, see \cite[Construction 1.19]{GS3}.} We also need to assume that any intersection of the irreducible components of $D$ is connected so that $\mathfrak{X} \to \mathcal{S}$ satisfies the same conditions as in Assumption \ref{toricdegenerationass}. Furthermore, we shall only work with projective log smooth morphisms $\mathfrak{X} \to \mathcal{S}$ that are \emph{minimal log CY}.\footnote{We refer to \cite{GS3} for the more general case of \emph{log CY} degenerations (see \cite[Definition 1.10]{GS3}).}

\begin{defn} \label{minlogcy}
$\mathfrak{X} \to \mathcal{S}$ is \emph{minimal log CY} if $K_{\mathfrak{X}} + D \equiv 0$ (numerically equivalent).
\end{defn}

The following lemma gives a logarithmic interpretation (used in \cite{GS3}) of Definition \ref{minlogcy}.

\begin{lemma}
Suppose that $D$ is simple normal crossings. Then we have $$c_1\left(\Omega^1_{\mathfrak{X}/\mathcal{S}}(\log D/\log 0)\right) \equiv K_{\mathfrak{X}} + D$$ where $c_1(\Omega^1_{\mathfrak{X}/\mathcal{S}}(\log D/\log 0))$ is the first Chern class of the relative sheaf of log differentials.
\end{lemma}

\begin{proof}
We calculate the first Chern class $c_1(\Omega^1_{\mathfrak{X}}(\log D))$ of the (absolute) sheaf of log differentials on $(\mathfrak{X}, D)$ first. From the definitions of the sheaf of log differentials and of the divisorial log structure, we have a short exact sequence 
$$0 \longrightarrow \Omega_{\mathfrak{X}}^1 \longrightarrow \Omega_{\mathfrak{X}}^1(\log D) \longrightarrow \oplus_i \mathcal{O}_{D_i} \longrightarrow 0$$ which implies that $$c_1(\Omega^1_{\mathfrak{X}}(\log D))= c_1(\Omega^1_{\mathfrak{X}}) + \sum_i c_1(\mathcal{O}_{D_i}) = K_{\mathfrak{X}} + \sum_i D_i = K_{\mathfrak{X}} + D.$$ Similarly, we have $$c_1(\Omega^1_{\mathcal{S}}(\log 0)) = K_{\mathcal{S}} + 0.$$
Now, since $g: \mathfrak{X} \to \mathcal{S}$ is log smooth, we have a logarithmic analogue of the classical short exact sequence (see \cite[Part III, Theorem 3.2.3]{Ogus}) 
$$0 \longrightarrow g^{*} \Omega_{\mathcal{S}}^1(\log 0) \longrightarrow \Omega_{\mathfrak{X}}^1(\log D) \longrightarrow \Omega^1_{\mathfrak{X}/\mathcal{S}}(\log D/\log 0) \longrightarrow 0$$ which implies that 
\begin{multline}
c_1\left(\Omega^1_{\mathfrak{X}/\mathcal{S}}(\log D/\log 0)\right) = c_1\left(\Omega^1_{\mathfrak{X}}(\log D)\right) - g^{*}\left( c_1\left(\Omega^1_{\mathcal{S}}(\log 0)\right)  \right) = \\ = K_{\mathfrak{X}} + D - g^{*} \left( K_{\mathcal{S}} + 0 \right).
\end{multline}
But $g^{*} \left( K_{\mathcal{S}} + 0 \right) \equiv 0$ since $\mathcal{S}$ is one-dimensional and the fibres of $g$ are numerically trivial. 
\end{proof}

As in the case of toric degenerations, $\mathfrak{X} \to \mathcal{S}$ gives rise to the structure of an integral affine manifold with singularities on the dual intersection complex $(B, \mathscr{P})$.\footnote{However, the location of the singularities is rather different from the toric degeneration case. The singularities of $(B, \mathscr{P})$ are contained in a union of cells of codimension at least $2$ whereas the singularities of $(\bar{B}, \bar{\mathscr{P}})$ are contained in the union of cells of the barycentric subdivision $\tilde{\bar{\mathscr{P}}}$ of $\bar{\mathscr{P}}$ not containing any vertex of $\bar{\mathscr{P}}$.} In \cite{GS3}, the mirror $\check{\mathfrak{X}} \to \Spec{\widehat{\kk[P]}}$ arises by defining an explicit product of the generators $\vartheta_p$ (that we call \emph{theta functions}) of the \emph{theta function ring} which correspond to rational points of the integral affine structure on $(B, \mathscr{P})$. In \cite{GS4} the mirror is obtained from the \emph{canonical scattering diagram} $\mathfrak{D}$ on $(B, \mathscr{P})$.\footnote{Strictly speaking, in \cite{GS3} the scattering diagram is defined on $(\mathbf{C}B, \mathbf{C}\mathscr{P})$, the cone over $(B, \mathscr{P})$. However, it is easy to modify the definition so that $\mathfrak{D}$ is defined on $(B, \mathscr{P})$ and gives rise to the same mirror, see Construction \ref{canonical scattering} and Remark \ref{diagramandcone}(1).} Both the product of the theta functions and the definition of $\mathfrak{D}$ use the theory of \emph{punctured log Gromov-Witten invariants} developed in \cite{ACGS, ACGS2} that allows, in particular, to make sense of negative contact orders. The two constructions agree by \cite[Theorem 6.1]{GS4}.

We shall call these mirrors \emph{intrinsic}. They are supposed to be a very general construction and generalize the previous approaches to the Gross-Siebert mirror symmetry. However, they are not expected to be double dual in general.\footnote{Although some forms of duality might work, see \cite[Section 0.4]{GHK} for a duality statement expected in the positive absolute surface case.} It is natural to investigate if the intrinsic mirror construction generalizes Batyrev-Borisov and (more generally) toric degeneration mirror symmetry.

\medskip

\noindent \textbf{Degenerations of elliptic curves.} The easiest case to consider is that of degenerations of elliptic curves since, in this case, all three mirrors can be computed explicitly. Consider a pair of dual $2$-dimensional reflexive polytopes $(\Delta, \Delta^*)$. As before, they give rise to dual Batyrev degenerations $\left( \mathfrak{X}_{\Delta}, \mathfrak{X}_{\Delta^{*}}\right)$ of elliptic curves that are both toric degenerations and log smooth.\footnote{Therefore, in the case of elliptic curves only, we drop the bar in the notation for a toric degeneration (and for related objects).} 

It is easy to explicitly compute both the toric degeneration mirror $\check{\mathfrak{X}}_{\TD} \to \Spec \kk\llbracket t \rrbracket = \Spec \widehat{\kk[\mathbb{N}]}$ (using e.g. the anticanonical polarization $A = -K_{\mathbb{P}_{\Delta}}$) and the intrinsic mirror $\check{\mathfrak{X}}_{\IMS} \to \Spec \widehat{\kk[NE(\mathfrak{X}_{\Delta, 0})]}$ (here $\mathfrak{X}_{\Delta, 0}$ is a wheel of irreducible rational curves so the submonoid $NE(\mathfrak{X}_{\Delta, 0}) \subseteq A_1(\mathfrak{X}_{\Delta,0}, \mathbb{Z})$ is finitely generated and we just let $P:=NE(\mathfrak{X}_{\Delta, 0})$) to $\mathfrak{X}_{\Delta}$, and to show that one has isomorphisms of the central fibres $$ \check{\mathfrak{X}}_{\IMS, 0}  \cong \check{\mathfrak{X}}_{\TD, 0} \cong \mathfrak{X}_{\Delta^*, 0} \cong  \partial \mathbb{P}_{\Delta^*}.$$
However, the isomorphism $ \check{\mathfrak{X}}_{\TD, 0} \cong \mathfrak{X}_{\Delta^*, 0}$ does not extend to the generic fibre and in most cases (unless $\mathbb{P}_{\Delta} \cong \mathbb{P}^2$ or $\mathbb{P}_{\Delta} \cong \mathbb{P}^1 \times \mathbb{P}^1$) there is no apparent connection between the two mirrors. 

This is, perhaps, not so surprising. By the technology of \cite{GHS}, both $\check{\mathfrak{X}}_{\TD}$ and $\check{\mathfrak{X}}_{\IMS}$ can be constructed by gluing together certain special thickenings of the irreducible components of $\check{\mathfrak{X}}_{\TD, 0}$. Viewing $\check{\mathfrak{X}}_{\TD}$ as a log deformation of the central fibre, one can prove that this deformation is, in fact, quite special and possesses well-behaved local models. In the language of \cite{models}, it is a \emph{divisorial log deformation} (see \cite[Definition 2.7]{models}).\footnote{Indeed, this follows from \cite[Corollary 2.18]{models} and the fact that simplicity (see \cite[Definition 1.60]{GS2}) of the dual intersection complex $(B, \mathscr{P})$ is a vacuous condition in dimension $1$.} Therefore, it has no reason to coincide with $\mathfrak{X}_{\Delta^*}$ (even though $\mathfrak{X}_{\Delta^*}$ has the same dual intersection complex), which is defined via toric geometry. 
	
On the other hand, in Chapter \ref{ell} we show that the constructions of the intrinsic and toric degeneration mirrors imply that for any toric degeneration $\mathfrak{X} \to \mathcal{S}$ of elliptic curves $\check{\mathfrak{X}}_{\TD} \to \Spec \kk\llbracket t \rrbracket$ is the basechange of $\check{\mathfrak{X}}_{\IMS} \to \Spec \widehat{\kk[NE(\mathfrak{X}_{0})]}$ by $$h: NE(\mathfrak{X}_{0}) \to \mathbb{N}, \quad \beta \mapsto A \cdot \beta$$ (inducing a map $\widehat{\kk[NE(\mathfrak{X}_{0})]} \to \kk\llbracket t \rrbracket$ by sending $t^{\beta}$ to $t^{h(\beta)} = t^{A \cdot \beta}$). Further, if one uses the universal version $\check{\mathfrak{X}}_{\TD} \to \Spec \widehat{\kk[NE(\mathfrak{X}_{0})]}$ of the toric degeneration mirror constructed in \cite[Appendix A.2]{GHS}, then $\check{\mathfrak{X}}_{\TD}  \cong \check{\mathfrak{X}}_{\IMS}$. One of the main objectives of this thesis is to generalize these observations about degenerations of elliptic curves to higher dimensions. We cover the case of degenerations of elliptic curves in more detail in Chapter \ref{ell}.

\medskip
	
\noindent \textbf{Resolution setup.} Suppose now that $(\Delta, \Delta^*)$ is a dual pair of $n$-dimensional reflexive polytopes for some $n \geq 3$. Again, they give rise to dual Batyrev degenerations $\left( \bar{\mathfrak{X}}_{\Delta}, \bar{\mathfrak{X}}_{\Delta^{*}}\right)$ of Calabi-Yau hypersurfaces of dimension $n-1$ that are toric degenerations. However, they are no longer log smooth, and to relate them to the log smooth setup, we need to resolve the log singularities. Assume that the generic fibre of $\bar{\mathfrak{X}}_{\Delta} \to \Spec \kk\llbracket t \rrbracket$ is smooth. This implies that $\bar{\mathfrak{X}}_{\Delta} \to \Spec \kk\llbracket t \rrbracket$ is log smooth away from a codimension $2$ subset $Z \subseteq \bar{\mathfrak{X}}_{\Delta, 0}$ of the central fibre. As we shall show in this thesis, the singularities of $Z$ can often be resolved (see Section \ref{admissible} for $n=3$ and Section \ref{resolvehd} for $n \geq 4$). Suppose that we have found a resolution $\pi: \mathfrak{X}_{\Delta} \to \bar{\mathfrak{X}}_{\Delta}$ of the log singularities such that the composition $\mathfrak{X}_{\Delta} \to \Spec \kk\llbracket t \rrbracket$ is a minimal log CY degeneration. Then it is natural to compare the dual Batyrev degeneration $\bar{\mathfrak{X}}_{\Delta^*}$, the toric degeneration mirror $\check{\bar{\mathfrak{X}}}_{\TD}$ to $\bar{\mathfrak{X}}_{\Delta}$, and the intrinsic mirror $\check{\mathfrak{X}}_{\IMS}$ to $\mathfrak{X}_{\Delta}$. 

As we shall discover, obtaining a \emph{small} (i.e. not contracting any divisors) resolution is not always possible, even in the case of $n=3$. If $\pi: \mathfrak{X}_{\Delta} \to \bar{\mathfrak{X}}_{\Delta}$ is not small, then $(B, \mathscr{P})$ is a non-trivial subdivision of $\left( \bar{B}, \bar{\mathscr{P}}\right)$ and as a result $\check{\mathfrak{X}}_{\IMS, 0}$ has more irreducible components than $\check{\bar{\mathfrak{X}}}_{\TD, 0}$. 
Assuming $\pi$ is small, one still has\footnote{The last isomorphism immediately follows from the Cox coordinate description of Batyrev degenerations (see Section \ref{batyrev}) via the anticanonical embedding, this is similar to the proof of a similar statement for degenerations of elliptic curves (see Proposition \ref{centralfibresprop}). Along with the results of \cite{G2}, this implies the second isomorphism (we use trivial gluing data to construct $\check{\bar{\mathfrak{X}}}_{\TD, 0}$). By the construction of the central fibre of the mirror family in \cite[Section 2.1]{GHS} (denoted $X_0$ in loc. cit.), it is clear that it only depends on the polyhedral manifold structure on the dual intersection complex. But the construction of the dual intersection complex as a height one slice of the tropicalization and the fact that $\pi: \mathfrak{X}_{\Delta} \to \bar{\mathfrak{X}}_{\Delta}$ is an isomorphism in codimension $1$ imply that $(B, \mathscr{P}) \cong \left( \bar{B}, \bar{\mathscr{P}}\right)$ as rational polyhedral complexes.}  
\begin{equation} \label{centralfibresequal}
\check{\mathfrak{X}}_{\IMS, 0}  \cong \check{\bar{\mathfrak{X}}}_{\TD, 0} \cong \bar{\mathfrak{X}}_{\Delta^*, 0} \cong  \partial \mathbb{P}_{\Delta^*}
\end{equation}
as in the $n=2$ case, but there is no reason to expect a more interesting connection between $\bar{\mathfrak{X}}_{\Delta^*}$ and $\check{\bar{\mathfrak{X}}}_{\TD}$, $\check{\mathfrak{X}}_{\IMS}$ as we don't have one even for $n = 2$. There is one notable exception. If $\mathbb{P}_{\Delta} \cong \mathbb{P}^3$, then one can use the symmetries of the dual intersection complex $\bar B \cong \partial \Delta^{*}$ of $\bar{\mathfrak{X}}_{\Delta}$ to restrict the form of the mirror equation and generalize the similar correspondence in the $n = 2$ case. In Appendix \ref{projcase}, we construct a rational map $\check{\bar{\mathfrak{X}}}_{\TD} \to \bar{\mathfrak{X}}_{\Delta^*}$ in the case that $\mathbb{P}_{\Delta} \cong \mathbb{P}^3$ using Cox coordinates \cite{cox} and conjecture a generalization to $\mathbb{P}^n, \ n \geq 4$. 

Since there is no reason to expect a (stronger than \eqref{centralfibresequal}) connection between $\bar{\mathfrak{X}}_{\Delta^*}$ and $\check{\bar{\mathfrak{X}}}_{\TD}$, $\check{\mathfrak{X}}_{\IMS}$ in general, we shall instead investigate if there is a connection between toric degeneration mirrors (that generalize Batyrev and Batyrev-Borisov mirror symmetry by \cite{G2}) and intrinsic mirrors.

\medskip 

\noindent \textbf{Main conjecture.}  We are going to be working with special toric degenerations. Roughly, we will define a toric degeneration $\bar{\mathfrak{X}} \to \mathcal{S}$ to be \emph{special} if it has a smooth generic fibre, one can construct the toric degeneration mirror $\check{\bar{\mathfrak{X}}} \to \Spec \kk\llbracket t \rrbracket$ (unless $\left(\bar{B}, \bar{\mathscr{P}}\right)$ is \emph{simple}\footnote{This is a technical maximal degeneracy assumption, see \cite[Definition 1.60]{GS2}.}, this is not automatic), and $\bar{\mathfrak{X}} \to \mathcal{S}$ admits good étale local models at points of the singular locus $Z \subseteq \bar{\mathfrak{X}}_0$. 
We defer the exact definition until Section \ref{specialtd} as it is rather technical. 

For every special toric degeneration $\bar{\mathfrak{X}} \to \mathcal{S}$ one should be able to define a class of \emph{strongly admissible} resolutions $\mathfrak{X} \to \bar{\mathfrak{X}}$ to minimal log CY degenerations $\mathfrak{X} \to \mathcal{S}$. Roughly, we require that a strongly admissible resolution $\pi: \mathfrak{X} \to \bar{\mathfrak{X}}$ is trivial away from strata of $\bar{\mathfrak{X}}_0$ of codimension at least one and étale locally at every point contained in a stratum of $\bar{\mathfrak{X}}_0$ of codimension at least one $\pi$ looks like one of the specified resolutions of the local model at that point. Further, the resolutions of the local models should be the same for two points in the same stratum and compatible with passing to deeper strata, and the dual intersection complex $\left(B, \mathscr{P}\right)$ of $\pi: \mathfrak{X} \to \mathcal{S}$ should be an integral subdivision of $\left(\bar{B}, \bar{\mathscr{P}}\right)$ (with a different affine structure). 

We say that a resolution $\pi: \mathfrak{X}' \to \bar{\mathfrak{X}}$ of a special toric degeneration $\bar{\mathfrak{X}} \to \mathcal{S}$ is \emph{admissible} if it factors as $\mathfrak{X}' \to \mathfrak{X} \to \bar{\mathfrak{X}}$ with $\mathfrak{X} \to \bar{\mathfrak{X}}$ strongly admissible and $\mathfrak{X}' \to \mathfrak{X}$ a logarithmic modification (i.e. a proper, representable, birational, log étale morphism). 

We are ready to state the main conjecture of this thesis. 
    
    \begin{conj} \label{conjecture}
        Let $\bar{\mathfrak{X}} \to \mathcal{S}$ be a special toric degeneration with polarization $A$. Then (possibly after a finite basechange) there exists an admissible resolution $\pi: \mathfrak{X} \to \bar{\mathfrak{X}}$ to a minimal log CY degeneration $\mathfrak{X} \to \mathcal{S}$.
        Moreover, for any admissible resolution $\pi$ the basechange of the intrinsic mirror $\check{\mathfrak{X}} \to \Spec \widehat{\kk[P]}$ by $h: P \to \mathbb{N}, \ \beta \mapsto \pi^{*}A \cdot \beta$ is isomorphic to the toric degeneration mirror $\check{\bar{\mathfrak{X}}} \to \Spec \kk\llbracket t \rrbracket$. 
    \end{conj}

There are a few things to clarify here:
\begin{enumerate}
    \item[1)] The definition of a \emph{special} toric degeneration should guarantee the existence of a (strongly) admissible $\pi: \mathfrak{X} \to \bar{\mathfrak{X}}$ and a $\pi$-ample effective simple normal crossings divisor $D'$ on $\mathfrak{X}$ supported on $D= \left(\mathfrak{X}_0 \right)_{\reduced}$. If $\bar{\mathfrak{X}} \to \mathcal{S}$ is a toric degeneration of K3-s, being special means that $\bar{\mathfrak{X}} \to \mathcal{S}$ is a divisorial log deformation in the sense of \cite[Definition 2.7]{models} (in particular, all the étale local models are toric) and has a smooth generic fibre.
    Basechange is also not necessary in this case. More generally, any $\bar{\mathfrak{X}} \to \mathcal{S}$ obtained via the reconstruction algorithm of \cite{GS1}, with a smooth generic fibre and with a simple $\left(\bar{B}, \bar{\mathscr{P}}\right)$ is special. This is a natural case to consider since such a $\bar{\mathfrak{X}} \to \mathcal{S}$ is a divisorial log deformation and the toric local models can be constructed from $\left(\bar{B}, \bar{\mathscr{P}}\right)$. In this thesis, unless $\bar{\mathfrak{X}} \to \mathcal{S}$ is a degeneration of K3-s, we shall mostly restrict our attention to this case.
    \item[2)] Note that a monoid homomorphism $h: P \to \mathbb{N}$ gives rise to a map $\widehat{\kk[P]} \to \widehat{\kk[\mathbb{N}]} = \kk\llbracket t \rrbracket$ only if the completion in $\widehat{\kk[P]}$ is taken with respect to an ideal $J \subseteq P$ with $h^{-1}(0) \subseteq P \setminus J$. The intrinsic mirror is defined over $\Spec \widehat{\kk[P]}$ where the completion in $\widehat{\kk[P]}$ is with respect to the maximal ideal. However, $h(\beta) = 0$ for any curve class $\beta$ of a curve contracted by $\pi: \mathfrak{X} \to \bar{\mathfrak{X}}$. As we explain in Section \ref{extendingintmirror}, the existence of a $\pi$-ample effective simple normal crossings divisor $D'$ supported on $D$ allows extending (for good choices of $P$ and under certain further restrictions on $D'$) the intrinsic mirror to a family over the completion of the toric stratum $\Spec \kk[K]$ corresponding to the face $K \subseteq P$ containing the classes of the contracted curves. This corresponds to taking the completion in $\widehat{\kk[P]}$ with respect to $J:=P \setminus K$. The map $h: P \to \mathbb{N}$ defines a basechange of this \emph{extended intrinsic mirror} which can be seen as the restriction to a specified one-parameter family through a point in $\Spec \kk[K]$. 
    \item[3)] Unless the dual intersection complex $\left(\bar B, \bar{\mathscr{P}}\right)$ of $\bar{\mathfrak{X}} \to \mathcal{S}$ is simple, there are some parameters (encoded in a choice of the \emph{initial slab functions} for $\left(\bar B, \bar{\mathscr{P}}\right)$) involved in the construction of the mirror $\check{\bar{\mathfrak{X}}}$ that the mirror depends on. In the above conjecture, the basechange should be isomorphic to the toric degeneration mirror for a certain specified choice of the initial slab functions.
    \item[4)] We will provide an exact definition of a strongly admissible resolution in two cases. For special toric degenerations of K3-s, we define a strongly admissible resolution in Definition \ref{stronglyadmissible}. In this case, one can often obtain a strongly admissible resolution by blowing up a sequence of irreducible components of the central fibre. In the case that $\bar{\mathfrak{X}} \to \mathcal{S}$ is a special toric degeneration of CY threefolds obtained via the reconstruction algorithm of \cite{GS1} and has a simple $\left(\bar{B}, \bar{\mathscr{P}}\right)$, we define a strongly admissible resolution in Definition \ref{stronglyadmissiblegeneral}. In both cases, one can define these resolutions tropically by subdividing $\left(\bar{B}, \bar{\mathscr{P}}\right)$ and using a certain PA-function on the subdivision to obtain the resolutions of the local models and glue them together. Under an additional assumption, the analysis generalizes to relative dimension $\geq 4$, see Section \ref{resolveindimatleast4}.
    \item[5)] Much of the literature on punctured log Gromov-Witten invariants and intrinsic mirrors has been written for the case that $D$ is a simple normal crossings divisor. However, having more flexibility in choosing a resolution is often desirable. Therefore, we will not restrict to this case and will argue that one may drop this requirement (see Section \ref{technicalasssection}). In all the relevant situations, we shall argue that there exists an admissible resolution $\mathfrak{X} \to \mathcal{S}$ such that $D$ is simple normal crossings. Therefore, Conjecture \ref{conjecture} still holds if we require in addition that $D$ is a simple normal crossings divisor. 
\end{enumerate}

The following shows that it is enough to prove Conjecture \ref{conjecture} for strongly admissible resolutions.

\begin{fact} \label{furtherresolution}
Suppose that Conjecture \ref{conjecture} holds for a toric degeneration $\bar{\mathfrak{X}} \to \mathcal{S}$ and a minimal log CY resolution $\pi: \mathfrak{X} \to \bar{\mathfrak{X}}$. Let $\pi': \mathfrak{X}' \to \mathfrak{X}$ be a further logarithmic modification. Then Conjecture \ref{conjecture} holds for $\bar{\mathfrak{X}} \to \mathcal{S}$ and the composed resolution $\mathfrak{X}' \to \mathcal{S}$.
\end{fact}

\begin{proof}
Let $P'$ be a monoid such that:
\begin{enumerate}
\item $NE(\mathfrak{X}'_0) \subseteq P \subseteq A_1(\mathfrak{X}'_0, \mathbb{Z})$.
\item $P'$ is finitely generated and saturated.
\item The group $(P')^{\times}$ of the invertible elements of $P'$ coincides with the torsion part of $A_1(\mathfrak{X}'_0, \mathbb{Z})$.
\end{enumerate}
Let $P$ be the saturation of $\pi'_{*}(P')$ and note that it satisfies the same conditions with respect to $\mathfrak{X}_0$. The result follows from the birational invariance of punctured log Gromov-Witten invariants recently proved by Samuel Johnston in \cite{J}. 

We claim that \cite[Corollary 1.6]{J} implies that the extended the intrinsic mirror $\check{\mathfrak{X}} \to \Spec \widehat{\kk[P]}$ is the basechange of the well-defined extended intrinsic mirror $\check{\mathfrak{X}}' \to \Spec \widehat{\kk[P']}$ by $\pi'_{*}: P' \to P$. Indeed, \cite[Corollary 1.6]{J} states that for any ideal $I \subseteq P$ with radical the maximal ideal the finite-order mirror $\check{\mathfrak{X}} \to \Spec \kk[P]/I$ is the basechange of the well-defined mirror over $\check{\mathfrak{X}}' \to \Spec \kk[P']/(\pi'_{*})^{-1}(I)$ via $\pi'_{*}: P' \to P$. The proof does not use the fact that the radical of $I$ is the maximal ideal apart from concluding that the mirror over $\kk[P]/I$ is well-defined. Therefore, using ideals $I$ with radical $J=P \setminus K$ corresponding to the curves contracted by $\pi$ (instead of ideals with radical the maximal ideal), we see that the extended intrinsic mirror $\check{\mathfrak{X}}' \to \Spec \widehat{\kk[P']}$ is well-defined if the extended intrinsic mirror $\check{\mathfrak{X}}' \to \Spec \widehat{\kk[P]}$ is. Further, the basechange claim holds since it holds modulo any ideal $I$ with radical $J$.

Let $\tilde{\pi}:= \pi' \circ \pi$. It suffices to show that the map $\tilde{h}: P' \to \mathbb{N}, \ \beta \mapsto \tilde{\pi}^{*}A \cdot \beta$ factors as the composition $\tilde{h} = h \circ \pi'_{*}$ with $h: P \to \mathbb{N}, \ \beta \mapsto \pi^{*}A \cdot \beta$ as before. Indeed, we have
\begin{multline}
h \circ \pi'_{*}(\beta) = h(\pi'_{*} \beta) = \pi^{*} A \cdot \pi'_{*} \beta = \\ = \pi_{*}'\left(\left(\pi'^{*} \circ \pi^{*} \right)\left(A\right) \cdot \beta \right) = \pi_{*}'\left(\tilde{\pi}^{*}A \cdot \beta \right) = \tilde{\pi}^{*}A \cdot \beta = \tilde{h}(\beta).
\end{multline}
Here we used the projection formula (see, e.g. \cite[Appendix A.1]{hartshorne}) for the third equality. 
\end{proof}

By \cite[Theorem 2.4.1 and Corollary 2.6.7]{AW}, there is a one-to-one correspondence between logarithmic modifications of $\mathfrak{X} \to \mathcal{S}$ and subdivisions of $(B, \mathscr{P})$.

\medskip

\noindent \textbf{Strategy and results.} A strategy to relate intrinsic mirrors and toric degeneration mirrors was first mentioned in \cite[Remark 2.15]{ann} of the paper that announced the development of the intrinsic mirrors viewpoint. The idea is to interpret intrinsic mirrors via scattering diagrams, compare the algorithmic scattering diagram $\bar{\mathfrak{D}}$ arising from a toric degeneration $\bar{\mathfrak{X}} \to \mathcal{S}$ and the canonical scattering diagram $\mathfrak{D}$ arising from the log smooth resolution $\mathfrak{X} \to \mathcal{S}$, and use the strong uniqueness properties of the inductive construction of $\bar{\mathfrak{D}}$ in \cite{GS1}. The first part was done in \cite{GS4}. In this thesis, we will prove the following.

\begin{theoremA}[Theorem \ref{proofconj}]
Conjecture \ref{conjecture} holds for special toric degenerations of K3-s. 
\end{theoremA}

We will also vastly expand the result. In Chapter \ref{gslocus}, we introduce the \emph{(numerical) minimal relative Gross-Siebert locus}  of the extended intrinsic mirror $\check{\mathfrak{X}} \to \Spec \widehat{\kk[P]}$. In the case that $\pi: \mathfrak{X} \to \bar{\mathfrak{X}}$ is small and  $A_1(\mathfrak{X}_0, \mathbb{Z})$ is the group of numerical curve classes on the central fibre, the minimal relative Gross-Siebert locus is just the completion of $\Spec \kk[K^{\gp}] \subseteq \Spec \kk[K]$. In general, one restricts to the completion of $\Spec \kk[E^{\gp}]$ where $E \subseteq K$ is the subface generated by $(-1)$-curves contracted by $\pi: \mathfrak{X} \to \bar{\mathfrak{X}}$. 

The definition of the relative Gross-Siebert locus is motivated by a similar construction in the absolute case of a log CY surface $(\mathfrak{X}, D)$ in \cite{GHK}. Such surfaces admit toric models $\pi: (\mathfrak{X}, D) \to (\bar{\mathfrak{X}}, \bar{D})$ (possibly after a toric blowup) where $(\bar{\mathfrak{X}}, \bar{D})$ is a toric variety with its boundary divisor and $\pi$ only contracts $(-1)$-curves. The authors show that the mirror $\check{\mathfrak{X}} \to \Spec \widehat{\kk[P]}$ extends to the completion of the stratum corresponding to the contracted curves and define the Gross-Siebert locus as the restriction to the completion of the large torus of that stratum. 

Our most general result can be summarized as follows (see Theorem \ref{compgen} for details). 

\begin{theoremB}[Theorem \ref{compgen}]
Let $\bar{\mathfrak{X}} \to \mathcal{S}$ be a special toric degeneration of K3-s and $\pi: \mathfrak{X} \to \bar{\mathfrak{X}}$ be a strongly admissible resolution. The basechange of Conjecture \ref{conjecture} extends to a correspondence between the restriction of the extended intrinsic mirror $\check{\mathfrak{X}} \to \Spec \widehat{\kk[P]}$ to the (numerical) minimal relative Gross-Siebert locus and a subfamily of the universal toric degeneration mirror of \cite[Theorem A.4.2]{GHS} varied in the free parameters of the initial slab functions and in gluing data.
\end{theoremB}

We also discuss in detail how to generalize the results to higher dimensions. In particular, we prove the following.

\begin{theoremC}[Theorem \ref{basechangstronglyadmissible}]
Let $\bar{\mathfrak{X}} \to \mathcal{S}$ be a special toric degeneration of CY threefolds obtained via the reconstruction algorithm of \cite{GS1} and with a simple $\left(\bar{B}, \bar{\mathscr{P}}\right)$. Then there exists an $L_{\bar{\mathfrak{X}}} \in \mathbb{Z}_{>0}$ such that for every $L \in \mathbb{Z}_{>0}$ with $L \geq L_{\bar{\mathfrak{X}}}$ the basechange $\bar{\mathfrak{X}}' \to \mathcal{S}$ of $\bar{\mathfrak{X}} \to \mathcal{S}$ by $R \to R, \ t \mapsto t^L$ (where $t$ is the uniformizer of $R$) admits a strongly admissible resolution $\pi: \mathfrak{X}' \to \bar{\mathfrak{X}}'$.
\end{theoremC}

In the course of the exposition, we shall often explain things in the basic case when $D$ is a simple normal crossings divisor and $\pi: \mathfrak{X} \to \bar{\mathfrak{X}}$ is small before explaining the general case. This is a natural case of interest in \cite{GHKS} (in the case of degenerations of K3-s).

\medskip

\noindent \textbf{Deformation philosophy.} We would like to mention a possible connection of this research to log deformation theory which may be an interesting topic to explore. Conjecture \ref{conjecture} suggests that the intrinsic mirror (at least in the case that $\pi: \mathfrak{X} \to \bar{\mathfrak{X}}$ is small) may be universal (or miniversal, etc.) for a suitable functor of log deformations of the central fibre of the mirror families (note from \eqref{centralfibresequal} that $\check{\mathfrak{X}}_{0}  \cong \check{\bar{\mathfrak{X}}}_{0}$ if $\pi$ is small).

In the case that $\left(\bar{B}, \bar{\mathscr{P}}\right)$ is simple, it is natural to consider the functor of divisorial log deformations since the toric degeneration mirrors (for any choice of polarization $A$) are deformations of this form. The functor of divisorial log deformations satisfies (see \cite[Appendix C]{RS}) the logarithmic analogue of Schlessinger's conditions \cite{sch} and therefore has a miniversal family by the results of \cite{K}. It is shown in \cite[Theorem 4.4]{RS} that this miniversal family is the completion of the universal toric degeneration of \cite[Theorem A.2.4]{GHS}. In general, there should be a corresponding theory of log deformations for toric degeneration mirrors since there are still well-behaved local models for the deformation (see \cite[Section 4.4.1]{GS1}) that can be \enquote{decomposed} into the divisorial local models. 

Intrinsic mirrors are not divisorial log deformations but one might be able to develop a suitable theory of log deformations that would include both toric degeneration and intrinsic mirrors. One should only consider deformations that are log smooth away from codimension $2$, and additional restrictions on the functor may be required. Indeed, all toric degenerations are log smooth away from the singular locus $Z$ of codimension $2$ and in Appendix \ref{log}, we also enhance $\check{\mathfrak{X}}\to \Spec \widehat{\kk[P]}$ to a log morphism that is log smooth away from codimension $2$. 

A result on the universality of the intrinsic mirror would connect the new theory back to the original hope of constructing mirrors in the Gross-Siebert program via a Bogomolov-Tian-Todorov type argument showing the smoothability of log CY spaces. We refer to \cite[Remark 2.19]{models} for a more detailed discussion. Smoothability of a general $\check{\bar{\mathfrak{X}}}_{0}$ (i.e. of a \emph{maximaly degenerate Calabi-Yau variety}) has recently been shown in \cite{CLM, FFR}.

\medskip

\noindent \textbf{Preliminaries.} We assume that the reader is familiar with logarithmic geometry but provide alternative interpretations where feasible. One can avoid the more technical discussion by restricting attention to the case that $D$ is a simple normal crossings divisor. We shall summarize the relevant facts about scattering diagrams, toric degenerations, and intrinsic mirror symmetry in Chapter \ref{setup}. We do not review punctured log Gromov-Witten theory, see \cite{ACGS, ACGS2}. We will not use punctured invariants directly and will refer to \cite{G3} and \cite{GHKS} for explicit computations.

\medskip

\noindent \textbf{Structure of the thesis.}  The rest of the thesis is organized as follows:

Chapter \ref{ell} is motivational, independent of the rest of the thesis, and covers the case of elliptic curves. In particular, in Proposition \ref{ellconjproof}, we prove Conjecture \ref{conjecture} for any toric degeneration of elliptic curves. 

In Chapter \ref{setup}, we review the setup in all dimensions, prove some auxiliary results (including a discussion of special toric degenerations in Section \ref{specialtd} and the construction of the extended intrinsic mirror following \cite{GHK} in Section \ref{extresults}), and set the stage for proving Conjecture \ref{conjecture} for toric degenerations of K3-s. 
We also explain the idea of the proof and give an overview of the results proved in Chapters \ref{proof} and \ref{gslocus}. 

Chapter \ref{proof} is the core of this thesis. It is devoted to the proof of Conjecture \ref{conjecture} for toric degenerations of K3-s, which involves constructing log smooth resolutions and introducing the notion of being an admissible resolution (Sections \ref{resolvesmall}, \ref{resolve}, and \ref{admissible}), interpreting the extended intrinsic mirror family via scattering diagrams (extending the results of \cite{GS4}) following \cite{GHKS} (Section \ref{scatint}), and relating the canonical scattering diagram $\mathfrak{D}$ giving rise to the intrinsic mirror and the algorithmic scattering diagram $\bar{\mathfrak{D}}$ giving rise to the toric degeneration mirror (Section \ref{main}). We prove Conjecture \ref{conjecture} in Theorem \ref{proofconj}.

In Chapter \ref{gslocus}, we gradually extend the correspondence to families over larger strata of the base of the extended intrinsic mirror. The most general result providing the correspondence between the (numerical) minimal relative Gross-Siebert locus of the extended intrinsic mirror (see Definition \ref{gslocusdef}) and the universal toric degeneration mirror of \cite[Theorem A.4.2]{GHS} is proved in Theorem \ref{compgen}.
We also briefly discuss our results, philosophy, and possible generalizations in Section \ref{discussion}. 

Chapter \ref{higherdimgen} discusses generalizing to higher dimensions. In particular, we construct admissible resolutions of special toric degenerations $\bar{\mathfrak{X}} \to \mathcal{S}$ of CY threefolds in the natural case that $\bar{\mathfrak{X}} \to \mathcal{S}$ is distinguished (i.e. can be obtained via the reconstruction algorithm of \cite{GS1}, see
Definition \ref{distinguished}) and has a simple dual intersection complex (see \cite[Definition 1.60]{GS2}). We also sketch a generalization to relative dimension $n \geq 4$ (modulo a combinatorial Conjecture \ref{combconjecture}) and reduce Conjecture \ref{conjecture} to Conjecture \ref{combconjecture} and a generalization of the results of Section \ref{scatint} (see Conjectures \ref{canonicalscathd} and \ref{scatdiaghigherpowershd}) in this case. Finally, we conjecture (under the same assumptions) a generalization of Theorem \ref{compgen} (the main result of Chapter \ref{gslocus}) to higher dimensions (see Conjecture \ref{conjecturehd}).

In Appendix \ref{projcase}, we prove a correspondence between the toric degeneration mirror to a Batyrev degeneration of K3-s in $\mathbb{P}^3$ (i.e. the toric degeneration of Example \ref{runningexample}) and the dual Batyrev degeneration by restricting the form of the equation for the toric degeneration mirror. We also conjecture a generalization to $\mathbb{P}^n, \ n \geq 4$.

In Appendix \ref{log}, we construct natural log structures on toric degeneration mirrors and intrinsic mirrors and note that all the basechanges considered in Chapters \ref{proof}, \ref{gslocus}, and \ref{higherdimgen} are also basechanges in the category of log schemes.

\medskip

\noindent \textbf{Related work.} 
The results of this thesis crucially depend on the theories of toric degeneration mirror symmetry \cite{G1, GS2, models, GS1}, intrinsic mirror symmetry \cite{GS3, GS4}, and punctured log Gromov-Witten invariants \cite{ACGS, ACGS2, G3}. 

Parts of this thesis are closely related to \cite{GHK} that carried out similar constructions for an arbitrary log CY surface $(\mathfrak{X}, D)$. The result was generalized to higher dimensions in \cite{AG} where the authors related the canonical scattering diagram of \cite{GS4} arising from a log CY variety $(\mathfrak{X}, D)$ that is a blowup of a toric variety along disjoint hypersurfaces in its toric boundary to the algorithmic scattering diagram arising from the toric variety. The results of this paper are relevant in higher dimensions (see the discussion after Conjecture \ref{toricresolvecount}). A related paper \cite{A} uses the results of \cite{AG} to calculate explicit equations of the mirrors to certain $(\mathfrak{X}, D)$. 

This thesis is also closely related to \cite{GHKS}, in particular in the basic case of a simple normal crossings $D$ and a small resolution $\pi: \mathfrak{X} \to \bar{\mathfrak{X}}$. The crucial part that we borrow (and slightly generalize) from this paper is the construction of a (finite) scattering diagram giving rise to the extended intrinsic mirror.

We shall use the results of \cite{G3} and \cite{GHKS} for explicit computations of punctured log Gromov-Witten invariants and the gluing formulas in the case of degenerations of K3-s. The gluing formula originally appeared in \cite{W} in the case of toric gluing strata (and arbitrary dimension).

\medskip

\noindent \textbf{Conventions.} Apart from Chapter \ref{ell}, where the notation would be redundant, we use the notation with a bar to denote the objects related to the toric degeneration $\bar{\mathfrak{X}} \to \mathcal{S}$: $\left(\bar B, \bar{\mathscr{P}}\right), \ \bar{\Delta}, \ \bar{\mathfrak{D}}$, etc. We use the notation without a bar to denote the similar objects related to the minimal log CY degeneration $\mathfrak{X} \to \mathcal{S}$ (that is usually a log smooth resolution of $\bar{\mathfrak{X}} \to \mathcal{S}$): $(B, \mathscr{P}), \ \Delta, \  \mathfrak{D}$, etc. We drop the bar for cells $\sigma \in \bar{\mathscr{P}}$ of the dual intersection complex $\left(\bar B, \bar{\mathscr{P}}\right)$ apart from places where that would lead to confusion. We shall always assume that the log structure on the total space $\mathfrak{X}$ of a log smooth $\mathfrak{X} \to \mathcal{S}$ is fine and saturated. 

We use the notation with a check to denote the mirror objects, e.g. $\check{\bar{\mathfrak{X}}}$ denotes the toric degeneration mirror to $\bar{\mathfrak{X}}$, $\check{\mathfrak{X}}$ denotes the intrinsic mirror to $\mathfrak{X}$, $(\check{B}, \check{\mathscr{P}})$ denotes the discrete Legendre transform of $(B, \mathscr{P})$, etc.

\begin{sloppypar}
$A_1(\mathfrak{X}_0, \mathbb{Z})$ stands for the relevant group of curve classes that may be the group of $1$-cycles on $\mathfrak{X}_0$ or the group of relative $1$-cycles for $\mathfrak{X} \to \mathcal{S}$  modulo algebraic or numerical equivalence (we sometimes write $A_1(\mathfrak{X}_0, \mathbb{Z})_{\alg}$, $A_1(\mathfrak{X}_0$, $\mathbb{Z})_{\numer}$,  $A_1(\mathfrak{X}/\mathcal{S}, \mathbb{Z})_{\alg}$, and $A_1(\mathfrak{X}/\mathcal{S}, \mathbb{Z})_{\numer}$ respectively) and $NE(\mathfrak{X}_0) \subseteq A_1(\mathfrak{X}_0, \mathbb{Z})$ stands for the submonoid generated by the effective curve classes. 
\end{sloppypar}

For $m \in M$ an element in a monoid, we use the notation $z^{m}$ for the corresponding element of the monoid ring $A [M]$ (here $A$ is any Noetherian ring). If $\beta \in P$ is a curve class in the monoid containing the effective curve classes, we shall use the notation $t^{\beta}$ instead. For any monoid ideal $I \subseteq P$, we abuse the notation by writing $I$ for the corresponding ideal in $A[P]$. We denote by $N$ the integer lattice of the relevant dimension (i.e. $N \cong \mathbb{Z}^n$ for some $n \geq 1$) and set $N_{\mathbb{R}}:= N \otimes_{\mathbb{Z}} \mathbb{R}$. We work over an algebraically closed field $\kk$ of characteristic $0$.

\medskip

\noindent \textbf{Acknowledgements.} I would like to thank my supervisor Mark Gross for suggesting this research direction, for our conversations on the topic, and for the continued advice and encouragement. I am also grateful to Bernd Siebert for communicating some ideas relevant to the project and to both Mark and Bernd for developing most of the theory that has made this project possible. I am thankful to Samuel Johnson for his research in \cite{J} that (among other things) allowed the extension of Corollary \ref{furtherresolution} and to Yixian Wu for our conversations on resolving toric degenerations and her research in \cite{W} that allowed to define an affine structure on $(B, \mathscr{P})$ in general. Last but not least, I am grateful to all the incredible people I have met through mathematics and to my family and friends for their unconditional love and support. 

The author was supported by the Cambridge University Department of Pure Mathematics and Mathematical Statistics and by Cambridge Trust.

\section{Degenerations of elliptic curves} \label{ell}

We shall discuss Batyrev degenerations of elliptic curves, compute the toric degeneration and intrinsic mirrors to a degeneration of elliptic curves, and compare the dual Batyrev degeneration, the toric degeneration mirror, and the intrinsic mirror. This chapter is motivational, and the results are not used in other parts of the thesis.

\subsection{Batyrev degenerations of elliptic curves} \label{batyrev}

There are $16$ (up to the action of $AGL(2, \mathbb{Z})$\footnote{The group of affine unimodular transformations consisting of translations by an integer vector and the linear transformations in $GL(2, \mathbb{Z})$.}) reflexive polytopes in dimension $2$ that form $6$ dual pairs and $4$ self-dual polytopes (see, e.g. \cite[Fig.1]{KOS} for the classification). We shall focus on the two pairs $(\Delta_1, \Delta_1^*)$ and $(\Delta_2, \Delta_2^*)$ as in Figure \ref{polytopes}. The computations for the other reflexive polytopes are similar.

\vspace{0.4cm}

\begin{figure}[H]
\begin{center}

\begin{minipage}{.2\textwidth}

\begin{tikzpicture}[scale=0.66]

\draw[fill] (5/4*0,5/4*0) circle (1.5pt) coordinate (m-0-0);
\draw[fill] (5/4*0,5/4*1) circle (1.5pt) coordinate (m-0-1);
\draw[fill] (5/4*1,5/4*0) circle (1.5pt) coordinate (m-1-0);
\draw[fill] (5/4*2,5/4*0) circle (1.5pt) coordinate (m-2-0);
\draw[fill] (5/4*1,5/4*1) circle (1.5pt) coordinate (m-1-1);
\draw[fill] (5/4*0,5/4*2) circle (1.5pt) coordinate (m-0-2);
\draw[fill] (5/4*3,5/4*0) circle (1.5pt) coordinate (m-3-0);
\draw[fill] (5/4*2,5/4*1) circle (1.5pt) coordinate (m-2-1);
\draw[fill] (5/4*1,5/4*2) circle (1.5pt) coordinate (m-1-2);
\draw[fill] (5/4*0,5/4*3) circle (1.5pt) coordinate (m-0-3);

\draw (m-0-0) -- (m-0-3) -- (m-3-0) -- cycle;
\draw(0.4,1.5) node{$D_0$};
\draw(1.5,0.5) node{$D_2$};
\draw(2.2,2.2) node{$D_1$};
\end{tikzpicture}
\captionsetup{labelformat=empty}
\caption{$\Delta_1$}
\end{minipage}
\hspace{0.2cm}
\begin{minipage}{.2\textwidth}
\begin{tikzpicture}

\draw[fill] (5/4*0,5/4*0) circle (1.5pt) coordinate (m-0-0);
\draw[fill] (5/4*1,5/4*1) circle (1.5pt) coordinate (m-1-1);
\draw[fill] (5/4*1,5/4*2) circle (1.5pt) coordinate (m-1-2);
\draw[fill] (5/4*2,5/4*1) circle (1.5pt) coordinate (m-2-1);

\draw (m-0-0) -- (m-1-2) -- (m-2-1) -- cycle;
\draw(0.3,1.5) node{$D_0$};
\draw(1.5,0.3) node{$D_2$};
\draw(2.1,2.1) node{$D_1$};
\draw(0.7,0.7) node{$A_2$};
\draw(1.9,1.3) node{$A_2$};
\draw(1.3,2.0) node{$A_2$};

\end{tikzpicture}
\captionsetup{labelformat=empty}
\caption{$\Delta_1^*$}
\end{minipage}
\hspace{1.5cm}
\begin{minipage}{.2\textwidth}
\begin{tikzpicture}

\foreach \x in {0,1,2}
   \foreach \y in {0,1,2} 
      \draw[fill] (5/4*\x,5/4*\y) circle (1.5pt) coordinate (m-\x-\y);

\draw (m-0-0) -- (m-0-2) -- (m-2-2) -- (m-2-0) -- cycle;
\draw(0.4,1.5) node{$D_0$};
\draw(1.5,0.3) node{$D_3$};
\draw(1.5,2.1) node{$D_1$};
\draw(2.1,1.5) node{$D_2$};

\end{tikzpicture}
\captionsetup{labelformat=empty}
\caption{$\Delta_2$}
\end{minipage}
\hspace{0.2cm}
\begin{minipage}{.2\textwidth}
   \begin{tikzpicture}

\draw[fill] (5/4*1,5/4*1) circle (1.5pt) coordinate (m-1-1);
\draw[fill] (5/4*0,5/4*1) circle (1.5pt) coordinate (m-0-1);
\draw[fill] (5/4*1,5/4*0) circle (1.5pt) coordinate (m-1-0);
\draw[fill] (5/4*1,5/4*2) circle (1.5pt) coordinate (m-1-2);
\draw[fill] (5/4*2,5/4*1) circle (1.5pt) coordinate (m-2-1);

\draw (m-1-0) -- (m-0-1) -- (m-1-2) -- (m-2-1) -- cycle;
\draw(0.4,0.4) node{$D_3$};
\draw(0.4,2.1) node{$D_0$};
\draw(2.1,0.4) node{$D_2$};
\draw(2.1,2.1) node{$D_1$};
\draw(2.0,1.3) node{$A_1$};
\draw(1.3,2.0) node{$A_1$};
\draw(1.3,0.5) node{$A_1$};
\draw(0.5,1.3) node{$A_1$};

\end{tikzpicture}
\captionsetup{labelformat=empty}
\caption{$\Delta_2^*$}
\end{minipage}
\end{center}
\setcounter{figure}{0}
\caption{Reflexive pairs $(\Delta_1, \Delta_1^*)$ and $(\Delta_2, \Delta_2^*)$.} \label{polytopes}
\end{figure}

Note that $\mathbb{P}_{\Delta_1} \cong \mathbb{P}^2$ and $\mathbb{P}_{\Delta_2} \cong \mathbb{P}^1 \times \mathbb{P}^1$ are smooth toric varieties, and $\mathbb{P}_{\Delta_1^*}$ and $\mathbb{P}_{\Delta_2^*}$ have $A_k, \ k \geq 1$ singularities. As in Chapter \ref{intro}, the polytopes of Figure \ref{polytopes} give rise to dual Batyrev degenerations $(\mathfrak{X}_{\Delta_1}, \mathfrak{X}_{\Delta_1^*})$ and $(\mathfrak{X}_{\Delta_2}, \mathfrak{X}_{\Delta_2^*})$. It is easy to read off information about the central fibre of $\mathfrak{X}_{\Delta}$ (which is just $\partial \mathbb{P}_{\Delta}$) from the polytope $\Delta$. Indeed, there is a one-to-one correspondence between components $D_i, \ 0 \leq i \leq m-1$ of the central fibre $\mathfrak{X}_{\Delta, 0}$ and edges $\delta_i$ of $\Delta$ (between intersections $D_{i-1} \cap D_i$ and vertices $\varepsilon_i$). The degree of $-K_{\mathbb{P}_{\Delta}}$ on $D_i$ is the integral length $d_i$ of $\delta_i$ and the singularity of the total space at $D_{i-1} \cap D_i$ is $A_{e_i}$ where $e_i = | \det (v_{\delta_{i-1}}, v_{\delta_{i}})|-1$ for $v_{\delta_{i}}$ the primitive generator of $\delta_i$ (we set $e_i:=0$ if the intersection $D_{i-1} \cap D_i$ is normal crossings). Here and later, we take the indices modulo $m$ so that $D_0 = D_m$, etc. For $\Delta^{*}$ we use notations $d_i^{*}$ and $e_i^{*}$. 
We have reflected this information in Figure \ref{polytopes}.

From the point of view of the Gross-Siebert mirror symmetry, knowing $d_i$ and $e_i$ corresponds to specifying the affine structure and kinks of a piecewise-linear function $\varphi$ on $\partial \Delta$ respectively. In the language of \cite{GS2}, $(\partial \Delta, e_i, d_i)$ is the \emph{intersection complex}\footnote{One can think of the intersection complex as the \enquote{real picture} of the central fibre with a polyhedral subdivision $\mathscr{P}$. Here $\mathscr{P}$ is the natural subdivision of $\partial \Delta$ by the edges and vertices of $\Delta$.} of $\mathfrak{X}_{\Delta}$, and $(\partial \Delta^*, e_i^*, d_i^*)$ is the dual intersection complex. Note that we have $e_i^* = d_i-1, \ d_i^* = e_i+1$. It is easy to see that this corresponds precisely to the \emph{discrete Legendre transform} (as defined in \cite[Section 1.4]{GS2}) between the intersection complex and its dual exchanging information about the affine structure and polarization.

A natural way to describe $\mathfrak{X}_{\Delta}$ explicitly is via the \emph{Cox coordinates} on $\mathbb{P}_{\Delta}$ (the construction works in all dimensions). We will briefly recall the setup, see \cite{cox} for more details. Consider a simplicial toric variety $\mathbb{P}_{\Delta}$ with $\dim \Delta = n$ and $m$ boundary divisors $D_{\rho}, \ \rho \subseteq \Delta$ corresponding to the maximal faces of $\Delta$. We can define its \emph{Cox ring} as $S(\mathbb{P}_{\Delta})=\kk[x_\rho \ | \ \rho \subseteq \Delta]$ with a natural grading as follows. For every divisor $D= \sum_{\rho} a_{\rho} D_{\rho}$ let $x^D : = \prod_{\rho} x_{\rho}^{a_{\rho}}$ be the corresponding element in $S(\mathbb{P}_{\Delta})$. We set $\deg (x^{D}): = \beta(D) \in A_{n-1}(\mathbb{P}_{\Delta})$ where $\beta$ is the map in the exact sequence 
\begin{equation} \label{exactseq}
\mathbb{Z}^m \longrightarrow \oplus_{\rho} \mathbb{Z} D_{\rho} \overset{\beta}{\longrightarrow} A_{n-1}(\mathbb{P}_{\Delta}, \mathbb{Z}) \longrightarrow 0
\end{equation}
taking a divisor to its class in the Chow group.

For each face $\sigma \subseteq \Delta$ let $$x^{\sigma} := \prod_{\sigma \not\subseteq \rho} x_{\rho}.$$ We call $S_{+}: = \langle x^{\sigma} \ | \ \sigma \subseteq \Delta \rangle$ the \emph{irrelevant ideal} and denote by $V(S_{+}) \subseteq \mathbb{A}^m$ the corresponding variety.

Let $G:=\Hom_{\mathbb{Z}} (A_{n-1}(\mathbb{P}_{\Delta}), \kk^{*})$. Dualizing \eqref{exactseq}, we obtain a natural action of $G$ on $\mathbb{A}^{m}$ leaving $V(S_{+})$ invariant (since it consists of coordinate subspaces). This action gives rise to a GIT description of toric varieties.

\begin{sloppypar}
\begin{theorem}[{\cite[Theorem 2.1(iii)]{cox}}]
$\mathbb{P}_{\Delta}$ is the geometric quotient $(\mathbb{A}^{m} \setminus V(S_{+}))/G$. 
\end{theorem}   
\end{sloppypar}

So one can describe $\mathbb{P}_{\Delta}$ by specifying the Cox coordinates $x_\rho, \  \rho \subseteq \Delta$ and the grading on the Cox ring $S(\mathbb{P}_{\Delta})$. One can also describe subvarieties of $\mathbb{P}_{\Delta}$ by giving expressions in terms of Cox coordinates on the cover. Given $D \in A_{n-1}(\mathbb{P}_{\Delta})$, the vanishing locus of the general section $s \in \mathcal{O}(D)$ is given by the sum of elements of the same degree in $S(\mathbb{P}_{\Delta})$: $$s = \sum_{D' \text{ s.t } \deg(x^{D'}) = \deg(x^D)}  x^{D'}$$ (where we abuse the notation by dropping the coefficients). The general section $s \in \mathcal{O}(-K_{\mathbb{P}_{\Delta}})$ admits a particularly nice description. We have $$s= \sum_{p \in \Delta(\mathbb{Z})} x^p$$ where the sum is over the lattice points of $\Delta$ and $$x^p :=\prod_{\rho} x_{\rho}^{\dist(p, D_\rho)}.$$ Here $\dist(p,D_{\rho}):=\langle p, n_{\rho} \rangle + a_{\rho}$ is the lattice distance between $p$ and the face $\langle m, n_{\rho} \rangle + a_{\rho} = 0$ of $\Delta$ corresponding to $D_{\rho}$. Note that in the case of a reflexive $\Delta$, there is just one lattice point $0$ in the interior of $\Delta$, and the section $s_0$ corresponding to $0$ is $\prod_{\rho} x_{\rho}$.

Using the above, we can compute explicit equations for the toric degenerations $\mathfrak{X}_{\Delta}$ (where $\Delta$ is a $2$-dimensional reflexive polytope) in the Cox coordinates. The polytopes of Figure \ref{polytopes} give rise to the toric degenerations in Figure \ref{gsequations} below.
The variables $x_i$ correspond to the $D_i$ in Figure \ref{polytopes}, and we list their degrees with respect to the standard basis of $A_{1} (\mathbb{P}_{\Delta}, \mathbb{Z}) \cong \mathbb{Z}^l \oplus \mathbb{Z}_k$ (for some $l \geq 1, \ k \geq 1$). To describe  $$\mathfrak{X}_{\Delta} = \left\{ ts+s_0 \right\} \subseteq \mathbb{P}_{\Delta} \times \Spec \kk\llbracket t \rrbracket$$ (with the anticanonical polarization) we give the Cox coordinates along with the expressions for $s_0$ and $s$ (we drop the $s_0$-monomial in $s$).

\begin{figure}[t]
\begin{center}

\begin{minipage}[t]{.4\textwidth}
\captionsetup{labelformat=empty}
\caption{\hspace{-3cm} $\mathfrak{X}_{\Delta_1}$}
$\begin{bmatrix}
& x_0 & x_1 & x_2  \\ \hline
\mathbb{Z} & 1 & 1 & 1 
\end{bmatrix}$ \\ \\ $s_0=x_0x_1x_2$ \\ $s=x_0^3+x_1^3+x_2^3+x_0x_1^2+x_0x_2^2+x_1x_0^2+x_1x_2^2+x_2x_0^2+x_2x_1^2$

\end{minipage}
\hspace{0.5cm}
\begin{minipage}[t]{.4\textwidth}
\captionsetup{labelformat=empty}
\caption{\hspace{-3cm} $\mathfrak{X}_{\Delta_1^*}$}
$\begin{bmatrix}
& x_0 & x_1 & x_2  \\ \hline
\mathbb{Z} & 1 & 1 & 1  \\ \mathbb{Z}_3 & 0 & 1 & 2 
\end{bmatrix}$ \\ \\ $s_0=x_0x_1x_2$ \\ $s=x_0^3+x_1^3+x_2^3$

\end{minipage}
\begin{minipage}[t]{.4\textwidth}
\captionsetup{labelformat=empty}
\caption{\hspace{-3cm} $\mathfrak{X}_{\Delta_2}$}
$\begin{bmatrix}
& x_0 & x_1 & x_2 & x_3  \\ \hline
\mathbb{Z} & 0 & 1 & 0 & 1 \\ \mathbb{Z} & 1 & 0 & 1 & 0 
\end{bmatrix}$ \\ \\ $s_0=x_0x_1x_2x_3$ \\ $s=x_0^2x_1^2+x_0^2x_3^2+x_1^2x_2^2+x_2^2x_3^2 + x_0x_2x_1^2 + x_0x_2x_3^2 + x_1 x_3 x_0^2 + x_1 x_3 x_2^2$

\end{minipage}
\hspace{0.5cm}
\begin{minipage}[t]{.4\textwidth}
\captionsetup{labelformat=empty}
\caption{\hspace{-3cm} $\mathfrak{X}_{\Delta_2^*}$}
$\begin{bmatrix}
& x_0 & x_1 & x_2 & x_3 \\ \hline
\mathbb{Z} & 0 & 1 & 0 & 1  \\ \mathbb{Z} & 1 & 0 & 1 & 0  \\ \mathbb{Z}_2 & 0 & 0 & 1 & 1 
\end{bmatrix}$ \\ \\ $s_0=x_0x_1x_2x_3$ \\ $s=x_0^2x_1^2+x_0^2x_3^2+x_1^2x_2^2+x_2^2x_3^2$

\end{minipage}
\end{center}
\setcounter{figure}{1}
\caption{Equations for the toric degenerations $\mathfrak{X}_{\Delta}$ in the Cox coordinates.} \label{gsequations}
\end{figure}

\subsection{The mirror to a degeneration of elliptic curves} \label{ellmirror}

We want to compute the toric degeneration mirrors $\check{\mathfrak{X}}_{\TD}$ and the intrinsic mirrors $\check{\mathfrak{X}}_{\IMS}$ to the Batyrev degenerations $\mathfrak{X}_{\Delta}$ of Figure \ref{gsequations}. More generally, let $\mathfrak{X} \to \mathcal{S}$ be a toric degeneration of elliptic curves in $\mathbb{P}_{\Delta} \times \mathcal{S} \to \mathcal{S}$.\footnote{Instead of an ambient $\mathbb{P}_{\Delta} \times \mathcal{S} \to \mathcal{S}$ one can specify a \emph{polarization} on $\mathfrak{X}$ (i.e. a relatively ample line bundle) defining the degrees $d_i$.} Then the central fibre is a cycle of rational curves $D_i, \ 0 \leq i \leq m-1$. As before, let $d_i$ specify the degrees of $-K_{\mathbb{P}_{\Delta}}$ on $D_i, \ 0 \leq i \leq m-1$ and  $e_i$ specify the singularities of the total space at the intersections $D_{i-1} \cap D_i$. By the duality discussed in Section \ref{batyrev}, the central fibre of the mirror $\check{\mathfrak{X}} \to \Spec \kk\llbracket t \rrbracket$ should again be a cycle of rational curves $D^{*}_i, \ 0 \leq i \leq m-1$, with $e_i^* = d_i-1, \ d_i^* = e_i+1$. We record this information in Figure \ref{centralfibres}.

\begin{figure}[H]
\begin{center}

\hspace{-1.5cm}
\begin{minipage}[t]{.4\textwidth}
\begin{tikzpicture}[line cap=round,line join=round,>=triangle 45,x=1cm,y=1cm, scale=0.45]
\clip(-1.127537653752384,0.8853926701456053) rectangle (14.613809545694705,10.158553320883723);
\draw  (5.175276826582198,10.068915317335827)-- (1.3552768265821977,4.908915317335829);
\draw  (3.575276826582198,9.46891531733583)-- (10.655276826582202,9.108915317335827);
\draw  (9.335276826582202,9.888915317335828)-- (12.935276826582204,4.808915317335828);
\draw  (12.895276826582204,6.308915317335828)-- (10.595276826582204,1.248915317335829);
\draw  (1.4952768265821983,6.408915317335829)-- (4.035276826582197,1.5289153173358274);
\begin{scriptsize}
\draw [fill=black] (4.686070477375847,1.3317724601929548) circle (1pt);
\draw [fill=black] (4.686070477375847,1.3317724601929548) circle (1pt);
\draw [fill=black] (4.686070477375847,1.3317724601929548) circle (1pt);
\draw [fill=black] (5.142567948159034,1.2546226259002484) circle (1pt);
\draw [fill=black] (5.630009078457307,1.1698502554135919) circle (1pt);
\draw [fill=black] (10.0593654363851,1.1910433480352545) circle (1pt);
\draw [fill=black] (9.677889769195147,1.1062709775485997) circle (1pt);
\draw [fill=black] (9.21164173151854,1.0214986070619423) circle (1pt);
\draw(2.7,8) node{\normalsize $D_0$};
\draw(5.1,8.6) node{\normalsize $A_{e_1}$};
\draw(2.8,5.5) node{\normalsize $A_{e_0}$};
\draw(2.2,3) node{\normalsize $D_{m-1}$};
\draw(7.3,9.7) node{\normalsize $D_1$};
\draw(9.7,8.5) node{\normalsize $A_{e_2}$};
\draw(11.8,7.5) node{\normalsize $D_2$};
\draw(11.6,5.5) node{\normalsize $A_{e_3}$};
\draw(12,3) node{\normalsize $D_3$};

\end{scriptsize}
\end{tikzpicture}
\captionsetup{labelformat=empty}
\caption{\hspace{2.5cm} $\mathfrak{X}_0$}
\end{minipage}
\hspace{1.5cm}
\begin{minipage}[t]{.4\textwidth}
\begin{tikzpicture}[line cap=round,line join=round,>=triangle 45,x=1cm,y=1cm, scale=0.45]
\clip(-1.127537653752384,0.8853926701456053) rectangle (14.613809545694705,10.158553320883723);
\draw  (5.175276826582198,10.068915317335827)-- (1.3552768265821977,4.908915317335829);
\draw  (3.575276826582198,9.46891531733583)-- (10.655276826582202,9.108915317335827);
\draw  (9.335276826582202,9.888915317335828)-- (12.935276826582204,4.808915317335828);
\draw  (12.895276826582204,6.308915317335828)-- (10.595276826582204,1.248915317335829);
\draw  (1.4952768265821983,6.408915317335829)-- (4.035276826582197,1.5289153173358274);
\begin{scriptsize}
\draw [fill=black] (4.686070477375847,1.3317724601929548) circle (1pt);
\draw [fill=black] (4.686070477375847,1.3317724601929548) circle (1pt);
\draw [fill=black] (4.686070477375847,1.3317724601929548) circle (1pt);
\draw [fill=black] (5.142567948159034,1.2546226259002484) circle (1pt);
\draw [fill=black] (5.630009078457307,1.1698502554135919) circle (1pt);
\draw [fill=black] (10.0593654363851,1.1910433480352545) circle (1pt);
\draw [fill=black] (9.677889769195147,1.1062709775485997) circle (1pt);
\draw [fill=black] (9.21164173151854,1.0214986070619423) circle (1pt);
\draw(2.7,8) node{\normalsize $D^*_0$};
\draw(5.1,8.6) node{\normalsize $A_{e^*_1}$};
\draw(2.8,5.5) node{\normalsize $A_{e^*_0}$};
\draw(2.2,3) node{\normalsize $D^*_{m-1}$};
\draw(7.6,9.7) node{\normalsize $D^*_1$};
\draw(9.5,8.5) node{\normalsize $A_{e^*_2}$};
\draw(11.8,7.5) node{\normalsize $D^*_2$};
\draw(11.6,5.5) node{\normalsize $A_{e^*_3}$};
\draw(12,3) node{\normalsize $D^*_3$};
\end{scriptsize}
\end{tikzpicture}
\captionsetup{labelformat=empty}
\caption{\hspace{2.5cm} $\check{\mathfrak{X}}_0$}
\end{minipage}
\end{center}
\setcounter{figure}{2}
\caption{The central fibres of a degeneration of elliptic curves and its mirror.} \label{centralfibres}
\end{figure}

We shall now explain how to construct the toric degeneration mirror $\check{\mathfrak{X}}_{\TD}$ and the intrinsic mirror $\check{\mathfrak{X}}_{\IMS}$ to $\mathfrak{X} \to \mathcal{S}$. We can make an explicit computation because the dual intersection complex $B$ of $\mathfrak{X} \to \mathcal{S}$ has an affine structure with no singularities. In this situation, the recipes of \cite{GS1} and \cite{GS4} for constructing mirrors reduce to a simple rule.

The construction that we give generalizes the classical construction (see, e.g. \cite[Chapter 8.4]{dirichlet} by Mark Gross) of the mirror to a degeneration of elliptic curves with no singularities in the total space. In \cite[Section 2]{theta}, the authors explain how one can explicitly compute mirrors to degenerations of abelian varieties (where the dual intersection complex $B$ is a torus) that includes our case, and \cite[Example 6.0.1]{GHS} implies that the mirrors that we construct are indeed the ones of \cite{GS1} (for the toric degeneration mirror) and \cite{GS4} (for the intrinsic mirror). The key to the construction is Mumford's description \cite{mumford} of a degeneration of abelian varieties of dimension $n$ via a toric cover $\mathfrak{X} \cong \hat{X}_{\Sigma} / \mathbb{Z}^n$ where $X_{\Sigma}$ is a toric variety that is not of finite type with a natural map $\pi: X_{\Sigma} \to \mathbb{A}^1$ and $\hat{X}_{\Sigma}$ is the completion of $X_{\Sigma}$ along $\pi^{-1}(0)$. We shall describe a tropical version of this for degenerations of elliptic curves. Unlike \cite[Chapter 8.4]{dirichlet}, we use the \emph{theta functions} of \cite[Section 3.3]{GHS} for the construction (and not the \emph{jagged path} theta functions that have been used historically, see \cite[Section 3.2]{GS4} or \cite[Section 4.5]{GHS}). They directly correspond to the classical theta functions on abelian varieties (see \cite[Example 6.0.1]{GHS}) but are much more general and admit a tropical construction via \emph{broken lines}.\footnote{We review the tropical construction of the theta functions in Section \ref{thetafunctions}.}
Our description will use the technology of \cite{GHS} as a black box. See Section \ref{scatteringsetup} and \cite{GHS} for the definitions of the terms we use when motivating the construction.

\subsubsection{The intrinsic mirror} \label{intellcurve}

We first construct the intrinsic mirror  $$\check{\mathfrak{X}}_{\IMS} \to \Spec \widehat{\kk[NE(\mathfrak{X}_0)]}.$$ The construction for the toric degeneration mirror $\check{\mathfrak{X}}_{\TD} \to \Spec \kk\llbracket t \rrbracket$ will be similar. The submonoid of $A_1(\mathfrak{X}_0, \mathbb{Z})$ generated by the effective curve classes is just $NE(\mathfrak{X}_0)= \oplus_i \mathbb{N}D_i$ (note that it is finitely generated) so we may indeed take the base of the family to be $$\Spec \widehat{\kk[NE(\mathfrak{X}_0)]} = \Spec \kk\left[\left[t^{D_0}, \dots, t^{D_{m-1}}\right]\right].$$ The intrinsic mirror can be viewed as a family
\begin{equation} \label{thetaring}
\check{\mathfrak{X}}_{\IMS} := \Proj \bigoplus_{p \in \mathbf{C}B(\mathbb{Z})} \kk\left[\left[t^{D_0}, \dots, t^{D_{m-1}}\right]\right] \vartheta_p \to \Spec \kk\left[\left[t^{D_0}, \dots, t^{D_{m-1}}\right]\right]
\end{equation}
where $\mathbf{C}B$ is the cone over the dual intersection complex $B$ of $\mathfrak{X}$ (equivalently, $\mathbf{C}B$ is the \emph{tropicalization} of $\mathfrak{X}$), $\mathbf{C}B(\mathbb{Z})$ is the set of integer points of $\mathbf{C}B$, and we need to define the products of the theta functions $\vartheta_p, \ p \in \mathbf{C}B(\mathbb{Z})$. For the purposes of this discussion one may think of $\vartheta_p, \ p \in \mathbf{C}B(\mathbb{Z})$ as the generators of the \emph{theta function ring} $$\bigoplus_{p \in \mathbf{C}B(\mathbb{Z})} \kk\left[\left[t^{D_0}, \dots, t^{D_{m-1}}\right]\right] \vartheta_p.$$

In our case, $\mathbf{C}B$ is just the cone over the $m$-gon with side lengths $d_i^*$. Let $d^*:=\sum_{i=0}^{m-1} d_i^*$ and consider a fan $\Sigma$ in $\mathbb{R}^2$ spanned by the countable collection of rays with primitive vectors $$\left(\sum_{i=0}^k d^*_i + jd^*, 1\right), \quad j \in \mathbb{Z}, \ 0 \leq k \leq m-1.$$ The support of $\Sigma$ is $\left\{\mathbb{R} \oplus \mathbb{R}_{>0} \right\} \cup (0, 0)$ and there is a natural action of $\mathbb{Z}$ on $\Sigma$ with $1 \in \mathbb{Z}$ acting by $(x, y) \to (x+d^*y, y)$. Then we have $\mathbf{C}B \cong \Sigma / \mathbb{Z}$ which corresponds to Mumford's description of $\mathfrak{X}$ as the quotient  $\mathfrak{X} \cong \hat{X}_{\Sigma} / \mathbb{Z}$. 

The integer points $\Sigma(\mathbb{Z}) = \left\{(p_1,p_2) \in \left\{\mathbb{Z} \oplus \mathbb{Z}_{>0}\right\} \cup (0,0) \right\}$ of $\Sigma$ define theta functions $\tilde{\vartheta}_{p_1,p_2}$ on the cover (viewed as generators of the theta function ring for $X_{\Sigma}$). We refer to the second term of the subscript (in either $\tilde{\vartheta}_{p_1,p_2}$ or $\vartheta_{p_1,p_2}$) as the \emph{degree} of the theta function (this corresponds to a grading on the theta function ring). Factoring  $\Sigma(\mathbb{Z})$ by the action of $\mathbb{Z}$ gives $\mathbf{C}B(\mathbb{Z})$ so there are exactly $d^*p_2$ theta functions $\vartheta_{0,p_2}, \dots, \vartheta_{d^*p_2-1, p_2}$ of degree $p_2$. 

We want to define a product $\vartheta_{p_1,p_2} \cdot \vartheta_{p'_1,p'_2}$ respecting the degrees. Let us first define products of theta functions $\tilde{\vartheta}_{p_1,p_2}$ on the cover and then descend it to a product of $\vartheta_{p_1,p_2}$. Denote the ray of $\Sigma$ with primitive generator $(i,1)$ by $\rho_i$ and denote the corresponding divisor of $X_{\Sigma}$ by $\tilde{D}_{\rho_i}$. In the language of \cite{GHS}, the fan $\Sigma$ defines a consistent scattering diagram with trivial wall functions. The triviality of the wall functions implies that the broken lines on $\Sigma$ are just rays supported on $|\Sigma|$ (and running off to infinity) and the broken line product formula for theta functions \cite[Theorem 3.5.1]{GHS} reduces to the following rule:
$$\tilde{\vartheta}_{p_1,p_2} \cdot \tilde{\vartheta}_{p'_1,p'_2}: =  \tilde{\vartheta}_{p_1 + p'_1, p_2+p'_2} t^{\deg(p_1, p_2, p'_1, p'_2)} \in \bigoplus_{p \in \Sigma(\mathbb{Z})} \kk\left[t^{\tilde{D}_{\rho_i}} \ | \ \rho_i \in \Sigma\right]\vartheta_p$$ 
where 
\begin{equation} \label{4pointdegree}
\deg(p_1, p_2, p'_1, p'_2) := \sum_{\rho_i} \left\langle (p_1,p_2), n_{\rho_i} \right\rangle \tilde{D}_{\rho_i} + \sum_{\rho_j} \left\langle (p'_1,p'_2), n_{\rho_j} \right\rangle \tilde{D}_{\rho_j}
\end{equation}
and the first sum is over all the rays intersecting the ray coming from infinity with direction vector $(p_1,p_2)$ and terminating at the point $(p_1+p'_1, p_2+p'_2)$, the second sum is over all the rays intersecting the ray coming from infinity with direction vector $(p'_1,p'_2)$ and terminating at the point $(p_1+p'_1, p_2+p'_2)$ (if the point $(p_1+p'_1, p_2+p'_2)$ is on a ray, we only include this ray in one of the sums), and $n_{\rho}$ is the primitive normal to $\rho$ pointing in the direction of the coming ray.

This product rule corresponds to a balanced tropical curve in $\Sigma$ where the term $t^{\langle(p_1,p_2), n_{\rho_i} \rangle \tilde{D}_{\rho_i}}$ comes from the wall-crossing automorphisms of \cite[(2.19)]{GHS} that use the canonically defined convex multi-valued piecewise-affine (MPA) function $\varphi$ on $\Sigma$ (with values in $NE(\mathfrak{X}_0)_{\mathbb{R}}^{\gp} \cong \mathbb{R}^n$) with kink $\tilde{D}_{\rho_i}$ at $\rho_i$. 

\begin{ex}
Suppose that $d_i^* = 2$ for $0 \leq i \leq m-1$ (we have $\mathfrak{X} = \mathfrak{X}_{\Delta_2^*}$ if $m = 4$). Then $$\tilde{\vartheta}_{-1, 1} \cdot \tilde{\vartheta}_{2,1} = \tilde{\vartheta}_{1,2}t^{\left\langle (-1, 1), (-1, 0) \right\rangle \tilde{D}_{\rho_{0}}} = \tilde{\vartheta}_{1,2}t^{\tilde{D}_{\rho_0}}$$ (see Figure \ref{tropics}).
\end{ex}

\begin{figure}[H]
\begin{center}
\begin{tikzpicture}[line cap=round,line join=round, scale = 1.5]
\clip(-0.7,-0.3) rectangle (8.7,4.3);
\draw [line width=1pt] (4,0) -- (4,6.451503796629907);
\draw [line width=1pt,domain=4:8.895642405250415] plot(\x,{(-4--1*\x)/2});
\draw [line width=1pt,domain=-0.9311121082279423:4] plot(\x,{(-4--1*\x)/-2});
\draw [line width=1pt,domain=4:8.895642405250415] plot(\x,{(-4--1*\x)/4});
\draw [line width=1pt,domain=-0.9311121082279423:4] plot(\x,{(-4--1*\x)/-4});
\draw [line width=2.8pt,color=qqqqff] (10.18268412099255,4.591096140853939)-- (5,2);
\draw [line width=2.8pt,color=qqqqff] (7.5204591310199485,3.260109969101327) -- (7.549997608949692,3.3782448214826823);
\draw [line width=2.8pt,color=qqqqff] (7.5204591310199485,3.260109969101327) -- (7.632686512042857,3.212851319371256);
\draw [line width=2.8pt,color=qqqqff] (2.6059973380551336,4.401678971280503)-- (5,2);
\draw [line width=2.8pt,color=qqqqff] (3.858945713373489,3.1447130485034838) -- (3.737517825701337,3.1355679339033404);
\draw [line width=2.8pt,color=qqqqff] (3.858945713373489,3.1447130485034838) -- (3.8684795123537943,3.266111037377162);
\draw [line width=2.8pt,color=qqqqff] (5,2)-- (4.409062660642468,0.8081383127076914);

\draw [line width=2.8pt,color=qqqqff] (4.382152405813356,0.7538703868289603) -- (4.344268278435678,0.8430078285088196);
\draw [line width=2.8pt,color=qqqqff] (4.382152405813356,0.7538703868289603) -- (4.476035214719449,0.777676033392335);
\begin{scriptsize}
\draw  [fill] (0,0) circle (1.5pt);
\draw  [fill] (1,0) circle (1.5pt);
\draw  [fill] (2,0) circle (1.5pt);
\draw  [fill] (3,0) circle (1.5pt);
\draw  [fill] (4,0) circle (1.5pt);
\draw  [fill] (5,0) circle (1.5pt);
\draw  [fill] (6,0) circle (1.5pt);
\draw  [fill] (7,0) circle (1.5pt);
\draw  [fill] (7,1) circle (1.5pt);
\draw  [fill] (6,1) circle (1.5pt);
\draw (6.069129682287555,1.2903157842613445) node {\normalsize $\tilde{\vartheta}_{2,1}$};
\draw  [fill] (5,1) circle (1.5pt);
\draw  [fill] (4,1) circle (1.5pt);
\draw  [fill] (3,1) circle (1.5pt);
\draw (3.0665102476136115,1.2903157842613445) node {\normalsize $\tilde{\vartheta}_{-1, 1}$};
\draw  [fill] (2,1) circle (1.5pt);
\draw  [fill] (1,1) circle (1.5pt);
\draw  [fill] (0,1) circle (1.5pt);
\draw  [fill] (0,2) circle (1.5pt);
\draw  [fill] (1,2) circle (1.5pt);
\draw [fill]  (2,2) circle (1.5pt);
\draw  [fill] (3,2) circle (1.5pt);
\draw  [fill] (4,2) circle (1.5pt);
\draw  [fill] (5,2) circle (1.5pt);
\draw (5.074126761119943,2.385318705428956) node {\normalsize $\tilde{\vartheta}_{1, 2}$};
\draw  [fill] (6,2) circle (1.5pt);
\draw  [fill] (7,2) circle (1.5pt);
\draw  [fill] (0,3) circle (1.5pt);
\draw  [fill] (1,3) circle (1.5pt);
\draw  [fill] (2,3) circle (1.5pt);
\draw [fill] (3,3) circle (1.5pt);
\draw [fill] (4,3) circle (1.5pt);
\draw [fill] (5,3) circle (1.5pt);
\draw [fill] (6,3) circle (1.5pt);
\draw [fill] (7,3) circle (1.5pt);
\draw [fill] (0,4) circle (1.5pt);
\draw [fill] (1,4) circle (1.5pt);
\draw [fill] (2,4) circle (1.5pt);
\draw [fill] (3,4) circle (1.5pt);
\draw [fill] (4,4) circle (1.5pt);
\draw [fill] (5,4) circle (1.5pt);
\draw [fill] (6,4) circle (1.5pt);
\draw [fill] (7,4) circle (1.5pt);
\draw [fill] (8,0) circle (1.5pt);
\draw [fill] (8,1) circle (1.5pt);
\draw [fill] (8,2) circle (1.5pt);
\draw [fill] (8,3) circle (1.5pt);
\draw [fill] (8,4) circle (1.5pt);
\draw [fill] (10.18268412099255,4.591096140853939) circle (2.5pt);
\draw (8.4, 3.3) node {\normalsize $(2,1)$};
\draw (2.8, 3.5) node {\normalsize $(-1,1)$};
\draw (2.8, 3.5) node {\normalsize $(-1,1)$};
\draw (4.35, 1.5) node {\normalsize $(1,2)$};
\draw (4.2, 4.2) node {\normalsize $\rho_0$};
\draw (-0.5, 2.05) node {\normalsize $\rho_{-2}$};
\draw (-0.5, 0.95) node {\normalsize $\rho_{-4}$};
\draw (8.5, 2.05) node {\normalsize $\rho_{2}$};
\draw (8.5, 0.95) node {\normalsize $\rho_{4}$};

\end{scriptsize}
\end{tikzpicture}

\end{center}
\caption{The tropical curve corresponding to the product $\tilde{\vartheta}_{-1, 1} \cdot \tilde{\vartheta}_{2,1} = \tilde{\vartheta}_{1,2}t^{\tilde{D}_{\rho_0}}$ when $d_i^* = 2$ for $0 \leq i \leq m-1$.} \label{tropics}
\end{figure}
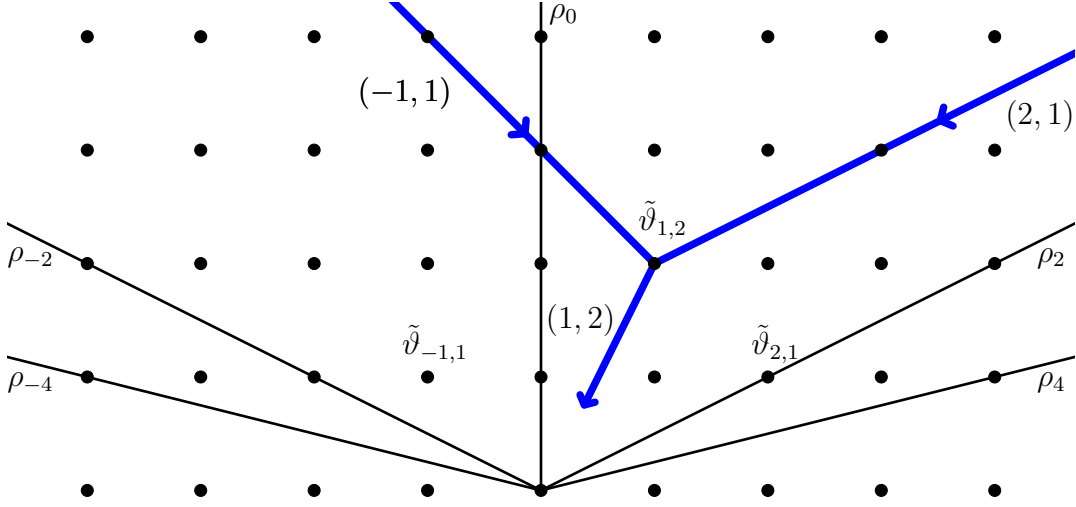

Note that all the $\tilde{D}_{\rho_i}$ with $i = \sum_{i=0}^k d^*_i + jd^*$ get identified with $D_k$ under the quotient map $\hat{X}_{\Sigma} \to \hat{X}_{\Sigma}/\mathbb{Z} \cong \mathfrak{X}$. Similarly, the theta functions $\tilde{\vartheta}_{p'_1 + p'_2 \alpha d^*,p'_2}$ for $ \alpha \in \mathbb{Z}$ get identified with $\vartheta_{p_1', p_2'}$. We define the product $\vartheta_{p_1,p_2} \cdot \vartheta_{p'_1,p'_2}$ by taking the sum over all the lifts of $\vartheta_{p_1', p_2'}$. This leads to the following formula:

\begin{multline} \label{productthetafs}
\vartheta_{p_1,p_2} \cdot \vartheta_{p'_1,p'_2} := \sum_{\alpha \in \mathbb{Z}} \vartheta_{p_1+p'_1+p'_2\alpha d^* \md d^*(p_2+p_2'), p_2+p_2'} t^{\deg(p_1, p'_1+p'_2\alpha d^*, p_2,p'_2)}  \\ \in \Proj \oplus_{p \in \mathbf{C}B(\mathbb{Z})} \kk\left[\left[t^{D_0}, \dots, t^{D_{m-1}}\right]\right] \vartheta_p   
\end{multline}
where $\deg(p_1, p_2, p'_1, p'_2)$ is as in \eqref{4pointdegree} with $\tilde{D}_{\rho_i}$ replaced by the corresponding $D_k$. Explicitly, we have: $$\deg(p_1, p_2, p'_1, p'_2) = \sum_{k=0}^{m-1} \deg_k (p_1, p_2, p'_1, p'_2) D_k$$ where
\begin{multline*}
\deg_k (p_1, p_2, p'_1, p'_2)= \left\langle (p_1,p_2), \sum_{j \text{ s.t }\frac{p_1}{p_2} < \sum_{i=0}^k d^*_i+jd^* \leq \frac{p_1+p'_1}{p_2+p'_2}} \left(-1, \sum_{i=0}^k d^*_i+jd^*\right) \right\rangle + \\ + \left\langle(p'_1,p'_2),  \sum_{j \text{ s.t }\frac{p_1+p'_1}{p_2+p'_2} < \sum_{i=0}^k d^*_i+jd^* < \frac{p'_1}{p'_2}} \left(1, -\sum_{i=0}^k d^*_i-jd^*\right) \right\rangle. 
\end{multline*}

Let $\vartheta_i := \vartheta_{i,1}$. One can compute the ideal of relations $I_{\rel}$ between $\vartheta_i, \ 0 \leq i \leq d^* - 1$ and represent

\begin{multline} \label{idealofrel}
\check{\mathfrak{X}}_{\IMS} = \Proj \bigoplus_{p \in \mathbf{C}B(\mathbb{Z})} \kk\left[\left[t^{D_0}, \dots, t^{D_{m-1}}\right]\right] \vartheta_p = \\ = \Proj \frac{\kk\left[\vartheta_0, \dots, \vartheta_{d^* - 1}\right]\left[\left[t^{D_0}, \dots, t^{D_{m-1}}\right]\right]}{I_{\rel}}
\end{multline}

as a subvariety of $\mathbb{P}^{d^{*}-1}_{\kk\left[\left[t^{D_0}, \dots, t^{D_{m-1}}\right]\right]}$ cut out by the ideal of relations $I_{\rel}$.

\subsubsection{The toric degeneration mirror}

To compute the toric degeneration mirror $\check{\mathfrak{X}}_{\TD} \to \Spec \kk\llbracket t \rrbracket$ we need to fix a relatively ample divisor (that we call \emph{polarization}) $A$ on $\mathfrak{X}$. The mirror can be defined via the theta functions as before:
$$\check{\mathfrak{X}}_{\TD} := \Proj \bigoplus_{p \in \mathbf{C}B(\mathbb{Z})} \kk\llbracket t \rrbracket \vartheta_p \to \Spec \kk\llbracket t \rrbracket.$$ The recipe to obtain the mirror is similar to that of Section \ref{intellcurve} except that now we use a different MPA function $\varphi_A$ on $\Sigma$ (with values in $\mathbb{N}^{\gp}_{\mathbb{R}} \cong \mathbb{R}$). Namely, we let the kink at $\rho_i$ be $A \cdot D_k > 0$ if $\tilde{D}_{\rho_i}$ is a lift of $D_k$. This corresponds to replacing $\tilde{D}_{\rho_i}$ and $D_k$ in the formulas by $A \cdot D_k$ and defines $\check{\mathfrak{X}}_{\TD}$ as a one-parameter family over $\Spec \kk\llbracket t \rrbracket$. 

As we are using the anticanoncal polarization $A = -K_{\mathbb{P}_{\Delta}}$, we have $-K_{\mathbb{P}_{\Delta}} \cdot D_k = d_k$. So we can recover the equations for the toric degeneration mirror $\check{\mathfrak{X}}_{\TD}$ from the equations for the intrinsic mirror $\check{\mathfrak{X}}_{\IMS}$ by replacing $D_k$ with $d_k$ in the equations for the mirror. Similarly to \eqref{idealofrel}, this represents
\begin{equation}
\check{\mathfrak{X}}_{\TD} = \Proj \bigoplus_{p \in \mathbf{C}B(\mathbb{Z})} \kk\llbracket t \rrbracket \vartheta_p = \Proj \frac{\kk[\vartheta_0, \dots, \vartheta_{d^* - 1}]\llbracket t \rrbracket}{I_{\rel}} 
\end{equation}
as a subvariety of $\mathbb{P}^{d^*-1}_{\kk\llbracket t \rrbracket}$. Setting $d_i = d^*_i = 1$ for $0 \leq i \leq m-1$, we recover the formula for the product of the theta functions of \cite[Chapter 8.4]{dirichlet}.\footnote{The setup in \cite[Chapter 8.4]{dirichlet} is slightly different and uses jagged path theta functions instead of the theta functions we use here. This is an equivalent viewpoint where one uses $\bigcup_{p_2 \in \mathbb{Z}_{\geq 0}}B\left(\frac{1}{p_2} \mathbb{Z}\right)$ instead of $\mathbf{C}B(\mathbb{Z})$ to parameterize the theta functions. The equivalence between the two viewpoints corresponds to replacing our $\vartheta_{p_1,p_2}$ with $\vartheta_{p_2, \frac{p_1}{p_2}}$ (in the notation of \cite[Chapter 8.4]{dirichlet}).}

We will only give explicit equations for toric degeneration mirrors for simplicity. Let $\check{\mathfrak{X}}_{\Delta_1, \TD}$ be the toric degeneration mirror to $\mathfrak{X}_{\Delta_1}$. It is easy to see that $\vartheta_0\vartheta_1\vartheta_2, \vartheta_0^3, \vartheta_1^3, \vartheta_2^3$ are all expressions in $\vartheta_{0,3}, \vartheta_{3,3}, \vartheta_{6,3}$ so we expect a unique relation between them of the form $$\lambda_1(t)\vartheta_0\vartheta_1\vartheta_2 = \lambda_2(t) (\vartheta_0^3 + \vartheta_1^3 + \vartheta_2^3)$$ for some $\lambda_1(t), \lambda_2(t) \in \kk\llbracket t \rrbracket$. Indeed, the mirror is described by a single equation $$\vartheta_0 \vartheta_1 \vartheta_2 = t^3\alpha(t)(\vartheta_0^3+\vartheta_1^3 + \vartheta_2^3)$$ in $\mathbb{P}_{\kk\llbracket t \rrbracket}^2$ where $$\alpha(t) =  1 - 5t^{9} + 28t^{18} - 150t^{27} + 794t^{36} - 4189t^{45} + 22092t^{54} + O(t^{55}).$$ In fact, the mirror is a Tate curve with a known $j$-invariant, and one could also compute $\alpha(t)$ from the $j$-invariant. More generally, similar considerations lead to the following. 

\begin{fact} \label{obstheta1}
Let $\check{\mathfrak{X}}_{\Delta, \TD}$ be the toric degeneration mirror to $\mathfrak{X}_{\Delta}$ for $\Delta$ a $2$-dimemsional reflexive polytope. Define a grading of $\mathbb{Z} \oplus \mathbb{Z}_{d^*}$ on $$\kk[\vartheta_0, \dots, \vartheta_{d^*-1}]\llbracket t \rrbracket$$ by setting $\deg(\vartheta_i)=(1,i)$. Unless $\Delta = \Delta_1$, the ideal $I_{\rel}$ defining $\check{\mathfrak{X}}_{\Delta, \TD}$ is generated by linear relations between any $3$ elements of degree $(2,i)$. If $\Delta = \Delta_1$, $I_{\rel}$ is generated by a single linear relation between the $4$ elements of degree $(3,0)$.
\end{fact}

Proposition \ref{obstheta1} follows from our definition of the product of the theta functions \eqref{productthetafs} by an explicit computation. We give equations for the toric degeneration mirrors (with polarization $-K_{\mathbb{P}_{\Delta}}$) to the degenerations of Figure \ref{gsequations} in Figure \ref{mirroreq} (here $\alpha(t), \ \beta(t), \ \gamma_i(t), \mu_j(t) \in \kk\llbracket t \rrbracket$ are certain power series with constant term $1$).

\begin{figure}[H]
\begin{center}

\begin{minipage}[t]{.3\textwidth}
\captionsetup{labelformat=empty}
\caption{$\check{\mathfrak{X}}_{\Delta_1, \TD}$}
1 equation in $\mathbb{P}^2_{\kk\llbracket t \rrbracket}$: \\ $\vartheta_0 \vartheta_1 \vartheta_2 = \\ t^3\alpha(t)(\vartheta_0^3+\vartheta_1^3 + \vartheta_2^3)$

\end{minipage}
\hspace{0.5cm}
\begin{minipage}[t]{.6\textwidth}
\captionsetup{labelformat=empty}
\caption{$\check{\mathfrak{X}}_{\Delta^*_1, \TD}$}
27 equations in $\mathbb{P}^8_{\kk\llbracket t \rrbracket}$: \\ $\vartheta_1 \vartheta_8  = t^4 \gamma_1(t)\vartheta_3 \vartheta_6 + t\gamma_2(t)\vartheta_0^2$ \\ $\vartheta_2 \vartheta_7 = t^2 \gamma_3(t)\vartheta_3 \vartheta_6 + t^2\gamma_4(t)\vartheta_0^2$ \\ $\vartheta_4 \vartheta_5 =  \gamma_5(t)\vartheta_3 \vartheta_6 + t^3\gamma_6(t)\vartheta_0^2$ \\ Changing $\vartheta_i \mapsto \vartheta_{i+3}, \vartheta_{i+6}$ in these gives 6 more equations. \\ Changing $\vartheta_i \mapsto \vartheta_{i+1}, \vartheta_{i+2}, \vartheta_{i+4}, \vartheta_{i+5}, \vartheta_{i+7}, \vartheta_{i+8}$ and multiplying the left side by $t$ gives the rest.

\end{minipage}
\end{center}
\begin{center}
\begin{minipage}[t]{.3\textwidth}
\captionsetup{labelformat=empty}
\caption{$\check{\mathfrak{X}}_{\Delta_2, \TD}$}

2 equations in $\mathbb{P}^3_{\kk\llbracket t \rrbracket}$: \\ $\vartheta_0 \vartheta_2 = t^2\beta(t)(\vartheta_1^2 + \vartheta_3^2)$ \\ $\vartheta_1 \vartheta_3 = t^2 \beta(t)(\vartheta_0^2+\vartheta_2^2)$

\end{minipage}
\hspace{0.4cm}
\begin{minipage}[t]{.6\textwidth}
\captionsetup{labelformat=empty}
\caption{$\check{\mathfrak{X}}_{\Delta^*_2, \TD}$}
20 equations in $\mathbb{P}^7_{\kk\llbracket t \rrbracket}$: \\ $\vartheta_3 \vartheta_5 = t^5\mu_1(t)\vartheta_0^2+t\mu_2(t)\vartheta_4^2$ \\ $\vartheta_2 \vartheta_6 = t^2 \mu_3(t)(\vartheta_0^2+\vartheta_4^2)$ \\ $\vartheta_1 \vartheta_7 = t \mu_2(t) \vartheta_0^2 + t^5 \mu_1(t) \vartheta_4^2$ \\  $\vartheta_4 \vartheta_6 = t^4\mu_1(t)\vartheta_1^2+\mu_2(t)\vartheta_5^2$ \\ $\vartheta_3 \vartheta_7 = t^2 \mu_3(t)(\vartheta_1^2+\vartheta_5^2)$ \\ $\vartheta_2 \vartheta_8 =  \mu_2(t) \vartheta_1^2 + t^4 \mu_1(t) \vartheta_5^2$ \\ Changing $\vartheta_i \mapsto \vartheta_{i+2}$ in these gives $6$ more equations. \\ We also have: \\ $\vartheta_3 \vartheta_6 = t^3 \mu_4(t) \vartheta_0 \vartheta_1 + t \mu_5(t) \vartheta_4 \vartheta_5$ \\ $\vartheta_2 \vartheta_7 = t \mu_5(t) \vartheta_0 \vartheta_1 + t^3 \mu_4(t) \vartheta_4 \vartheta_5$ \\ Changing $\vartheta_i \mapsto \vartheta_{i+1}, \vartheta_{i+2}, \vartheta_{i+3}$ gives the rest.

\end{minipage}
\end{center}
\setcounter{figure}{4}
\caption{Equations for the toric degeneration mirrors $\check{\mathfrak{X}}_{\Delta, \TD}$.} \label{mirroreq}
\end{figure}

\subsection{Comparing the mirror families}

We first note that Conjecture \ref{conjecture} tautologically holds for toric degenerations of elliptic curves.

\begin{fact} \label{ellconjproof}
Let $\mathfrak{X} \to \mathcal{S}$ be a toric degeneration of elliptic curves with polarization $A$. Then $\mathfrak{X} \to \mathcal{S}$ is a minimal log CY degeneration, and the basechange of the intrinsic mirror $\check{\mathfrak{X}}_{\IMS} \to \Spec \widehat{\kk[NE(\mathfrak{X}_0)]}$ by $h: NE(\mathfrak{X}) \to \mathbb{N}, \ \beta \mapsto A \cdot \beta$ is isomorphic to the toric degeneration mirror $\check{\mathfrak{X}}_{\TD} \to \Spec \kk\llbracket t \rrbracket$. Moreover, the universal toric degeneration mirror $\check{\mathfrak{X}}_{\TD} \to \Spec \widehat{\kk[NE(\mathfrak{X}_0)]}$ of \cite[Theorem A.2.4]{GHS} is isomorphic to $\check{\mathfrak{X}}_{\IMS} \to \Spec \widehat{\kk[NE(\mathfrak{X}_0)]}$.
\end{fact}

\begin{proof}
The first statement follows from our description of the mirrors in Section \ref{ellmirror}. $\mathfrak{X} \to \mathcal{S}$ satisfies $K_{\mathfrak{X}}+D = 0$ as a toric degeneration, so it is minimal log CY. The basechange corresponds exactly to replacing $D_k$ by $A \cdot D_k > 0$ (it is well-defined since $A$ is a relatively ample divisor).

For the second statement note that the dual intersection complex $B$ of $\mathfrak{X} \to \mathcal{S}$ has no singularities and that $\varphi_{-K_{\mathbb{P}_{\Delta}}}$ is a strictly convex MPA function since $-K_{\mathbb{P}_{\Delta}}$ is relatively ample. So we are in the setup of \cite[Theorem A.2.4]{GHS}. Moreover, $A_{\mathbb{P}} = \kk$ as there is only trivial gluing data in dimension $1$. The universal monoid $Q \subseteq \bMPA(B, \mathbb{N})^*$ coincides with $Q_0 =  \MPA(B, \mathbb{N})^* \cong \mathbb{N}^m \cong NE(\mathfrak{X}_0)$ by \cite[Proposition 1.2.9(a)]{GHS} and the fact that $\MPA(B, \mathbb{N}) = \bMPA(B, \mathbb{N})$ in dimension 1. The canonical version \cite[Theorem 5.2]{GS1} of the reconstruction algorithm corresponds to interpreting powers of $t$ as elements of $Q \cong NE(\mathfrak{X}_0)$ in exactly the same way as when constructing $\check{\mathfrak{X}}_{\IMS} \to \Spec \widehat{\kk[NE(\mathfrak{X}_0)]}$.
\end{proof}

The first part of Proposition \ref{ellconjproof} motivates Conjecture \ref{conjecture}. The second part inspired the generalizations of Chapter \ref{gslocus}.

The connection between $\check{\mathfrak{X}}_{\Delta, \TD}$ (and thus $\check{\mathfrak{X}}_{\Delta, \IMS}$ by the above) and $\mathfrak{X}_{\Delta^*}$ is more subtle. We use the Cox coordinate description \cite{coxrep} for rational maps between subvarieties of toric varieties. 

\begin{fact} \label{centralfibresprop}
    We have $\check{\mathfrak{X}}_{\Delta, \TD, 0} \cong \mathfrak{X}_{\Delta^*, 0} \cong \partial \mathbb{P}_{\Delta^*} $
\end{fact}

\begin{proof}
   The second identification is immediate from the Cox coordinate description of $\mathfrak{X}_{\Delta^*}$. 
   The first one is a consequence of the discrete Legendre transform. To describe the identification explicitly, note that we have natural bijections: $$\left\{\vartheta_i, \ 0 \leq i \leq m-1\right\} \Longleftrightarrow \left\{\mbox{sections } s^*_i \neq s^*_0 \mbox{ of } -K_{\mathbb{P}_{\Delta^*}} \right\} \Longleftrightarrow \partial \Delta^* (\mathbb{Z})$$
   The identification comes from considering the anticanonical embedding of $\mathfrak{X}_{\Delta^*}$ into $\mathbb{P}_{\kk\llbracket t \rrbracket}^m$ with coordinates $(\vartheta_0, \dots, \vartheta_{m-1}, \vartheta^*)$ and $\vartheta^*$ corresponding to $0 \in \Int \Delta$. In the Cox coordinates, it can be described as $\vartheta_i \mapsto s_i, \ \vartheta^* \mapsto s_0$. It is easy to see that the equations for the image of this embedding and the equations for $\check{\mathfrak{X}}_{\Delta, \TD}$ agree modulo the ideal $(t) \subseteq \kk\llbracket t \rrbracket$ (the image of the embedding has an additional equation $\vartheta^*=0$). This implies the first identification. Note that the isomorphism does not extend to the generic fibre.
\end{proof}

We slightly modify $\mathfrak{X}_{\Delta}$ by allowing more general deformations. Let $\eta(t) \in \kk\llbracket t \rrbracket$ and denote 
\begin{equation} \label{gentoricbat}
\mathfrak{X}^{\eta(t)}_{\Delta} := \left\{ \eta(t)s+s_0 =0 \right\} \subseteq \mathbb{P}_{\Delta} \times \Spec \kk\llbracket t \rrbracket.
\end{equation}

\begin{fact} \label{ldr}
If $\Delta = \Delta_1$ or $\Delta = \Delta_2$ (that is $\mathbb{P}_{\Delta} \cong \mathbb{P}^2$ or $\mathbb{P}_{\Delta} \cong \mathbb{P}^1 \times \mathbb{P}^1$), there is a rational map $\check{\mathfrak{X}}_{\Delta, \TD} \to \mathfrak{X}^{\eta(t)}_{\Delta^*}$ for a certain $\eta(t) \in \kk\llbracket t \rrbracket$. Explicitly:
\begin{enumerate}
    \item There is a rational map $\check{\mathfrak{X}}_{\Delta_1, \TD} \to \mathfrak{X}^{-t^3\alpha(t)}_{\Delta_1^*}$ given by $$x_0 \mapsto \vartheta_0, \ x_1 \mapsto \vartheta_1, \ x_2 \mapsto \vartheta_2$$ in the Cox coordinates.
    \item There is a rational map $\check{\mathfrak{X}}_{\Delta_2, \TD} \to \mathfrak{X}^{-t^4\beta^2(t)}_{\Delta_2^*}$ given by $$x_0 \mapsto \vartheta_0, \ x_1 \mapsto \vartheta_1, \  x_2 \mapsto \vartheta_2, \  x_3 \mapsto \vartheta_3$$ in the Cox coordinates.
\end{enumerate}
Here $\alpha(t), \beta(t) \in \kk\llbracket t \rrbracket$ are as in Figure \ref{mirroreq}.
\end{fact}

\begin{proof}
By the description of the image of a map given in the Cox coordinates \cite[Theorem 1.1]{coxrep}, it is enough to show that pulling back the equations for $\check{\mathfrak{X}}_{\Delta, \TD}$ to $S(\mathbb{P}_{\Delta^*})$ and taking the homogeneous part of the ideal generated by the pulled back equations with respect to the grading on $S(\mathbb{P}_{\Delta^*})$ gives the ideal $(\eta(t)s + s_0) \subseteq S(\mathbb{P}_{\Delta^*})$.

By comparing the equations in Figures \ref{gsequations} and \ref{mirroreq}, we see that in the first case, the pulled back equation is just $x_0x_1x_2 -t^3\alpha(t)(x_0^3+x_1^3+x_2^3)$ which proves (1). 

In the second case, the pulled back equations generate an ideal $$\left\langle x_0x_2 - t^2\beta(t)(x_1^2+x_3^2), x_1x_3-t^2\beta(t)(x_0^2+x_2^2) \right\rangle.$$ However, the pulled back equations are not homogeneous with respect to the grading on $S(\mathbb{P}_{\Delta_2^*})$. It is easy to check that the homogeneous part of this ideal is $$\left\langle x_0x_2x_2x_3 - t^4\beta^2(t)(x_0^2 x_1^2 + x_0^2 x_3^2 + x_1^2x_2^2 + x_2^2 x_3^2) \right\rangle $$ which proves (2).
\end{proof}

Note that the maps in the above proposition are not birational but rather quotients by an action of $\mathbb{Z}_3$ (in the first case) or $\mathbb{Z}_2$ (in the second case).

There are two reasons why we can construct these maps. First, in these examples, the number of degree one theta functions $\vartheta_i$ is the same as that of Cox coordinates $x_j$. Second, $\Delta_1$ and $\Delta_2$ have a lot of symmetries. One can explicitly check that for other $2$-dimensional reflexive polytopes, there are no rational maps given by linear expressions in terms of the theta functions in the Cox coordinates. So the results of Proposition \ref{ldr} are a coincidence and not a rule. 

In higher dimensions, one does not have a similar result for $\mathbb{P}_{\Delta} = (\mathbb{P}^1)^n$ since the number of $\vartheta_i$ and $x_j$ is no longer the same. However, in Appendix \ref{projcase}, we show that the result of Proposition \ref{ldr}(1) directly generalizes to the case of a Batyrev degeneration of K3-s in $\mathbb{P}^3$ (i.e. the case of Example \ref{runningexample}) and conjecture a generalization to the case of a Baryrev degeneration in $\mathbb{P}^n, \ n \geq 4$.

\section{Setup and preliminaries} \label{setup}

In this chapter, we review the setups, and adapt them to our needs, of scattering diagrams \cite{GHS, theta}, toric degenerations \cite{GS2, models, GS1}, and intrinsic mirrors \cite{ann, GS3, GS4}. We also discuss the resolution setup of a map $ \mathfrak{X} \overset{\pi}{\to}\bar{\mathfrak{X}} \to \mathcal{S}$ where $\bar{\mathfrak{X}} \to \mathcal{S}$ is a toric degeneration and the composed $\mathfrak{X} \to \mathcal{S}$ is a log smooth minimal log CY degeneration, and explain how $\pi: \mathfrak{X} \to \bar{\mathfrak{X}}$ gives rise to an extension of the intrinsic mirror $\check{\mathfrak{X}} \to \Spec \widehat{\kk[P]}$. We mostly keep the discussion general but sometimes specialize to degenerations of K3-s. Finally, we set up the scene for proving Conjecture \ref{conjecture} in the case of K3-s and give an overview of the results of Chapters \ref{proof} and \ref{gslocus}.

\subsection{Scattering diagrams} \label{scatteringsetup}

Both the toric degeneration and the intrinsic mirror constructions can be understood in terms of a combinatorial device of a \emph{scattering diagram} (also called \emph{wall structure} in \cite{GHS} and \cite{GS3}). The descriptions vary slightly throughout the literature. We stick to the conventions of \cite{GHS} except for imposing stronger conditions on the polyhedral (pseudo-)manifold $(B, \mathscr{P})$.

\subsubsection{Polyhedral manifold} \label{polymanifold}

We start with a \emph{polyhedral manifold} $(B, \mathscr{P})$. For us, such a manifold will usually arise as the dual intersection complex of a toric degeneration $\bar{\mathfrak{X}} \to \mathcal{S}$ or of its log smooth resolution $\mathfrak{X} \to \mathcal{S}$. 

\begin{defn} \label{polymani}
We say that a pair $(B, \mathscr{P})$ is a \emph{polyhedral manifold} of dimension $n$ if:
\begin{enumerate}
\item $B$ is a real topological manifold with $\partial B = 0$ and of pure dimension $n$.
\item $\mathscr{P}$ is a \emph{polyhedral decomposition} of $B$. That is, $\mathscr{P}$ is a set of bounded integral polyhedra along with a set of integral affine maps $\omega \to \tau$ identifying $\omega$ with a face of $\tau$ making $\mathscr{P}$ into a category such that:
\begin{enumerate}
\item For $\tau \in \mathscr{P}$ any face of $\tau$ occurs as a domain of an element of $\hom(\mathscr{P})$ with target $\tau$.
\item We have $$B =\lim_{\substack{\longrightarrow\\\tau\in\mathscr{P}}} \tau$$ with the colimit taken in the category of topological spaces.
\item For each $\tau \in \mathscr{P}$ the canonical map $\tau \to B$ is injective. That is, no cells self-intersect.
\item By abuse of notation we view the elements of $\mathscr{P}$ as subsets of $B$. We assume that a (possibly empty) intersection of any two cells of $\mathscr{P}$ is a cell of $\mathscr{P}$.
\end{enumerate}
\end{enumerate}
We refer to the cells of dimensions $0$, $1$, and $n$ as \emph{vertices}, \emph{edges}, and \emph{maximal cells} respectively. We use the notation $\mathscr{P}^{[k]}$ for the set of $k$-cells and $\mathscr{P}^{\max}:=\mathscr{P}^{[n]}$ for the set of maximal cells. We denote by $B(\mathbb{Z})$ the set of points of $B$ that are the images of integer points of some $\tau \in \mathscr{P}$ under $\tau \to B$. 

We say that two polyhedral manifolds $(B, \mathscr{P})$ and $\left(B', \mathscr{P}'\right)$ are \emph{isomorphic} if there exists an isomorphism $(B, \mathscr{P}) \cong \left(B', \mathscr{P}'\right)$ of polyhedral complexes that preserves the integral structure.

We also define a \emph{polyhedral manifold with boundary} as a pair $(B, \mathscr{P})$ satisfying all the conditions of (1) and (2) apart from requiring $\partial B = 0$. Unless explicitly mentioned otherwise, all the polyhedral manifolds in this thesis will be without boundary. 
\end{defn}

Our requirements on $(B, \mathscr{P})$ differ from both those of \cite{GHS} and those of \cite{GS2, GS1}, as we now explain. 

\begin{remarks} \label{polymaniremarks}
There are two directions in which Definition \ref{polymani} can be generalized:
\begin{enumerate}
\item Our notions are more restrictive than those of a \emph{polyhedral pseudomanifold} of \cite[Construction 1.1.1]{GHS} which allow $B$ not to be a manifold in codimension $\geq 2$ or have a non-empty boundary, and allow $\mathscr{P}$ to contain unbounded polyhedra. We include the additional assumptions since our manifolds will always satisfy them, leading to a simplified exposition. 
\item The dual intersection complex of a toric degeneration $\bar{\mathfrak{X}} \to \mathcal{S}$ satisfying Definition \ref{toricdegeneration} might not satisfy conditions (2)(c) and (2)(d) of Definition \ref{polymani}. Dropping these conditions leads to the notion of a \emph{toric polyhedral decomposition} of \cite[Definitions 1.21 and 1.22]{GS2}. In \cite{GS2, GS1}, the authors developed a sophisticated language to deal with such decompositions. It is easy to see that conditions (3) and (4) of Assumption \ref{toricdegenerationass} on $\bar{\mathfrak{X}} \to \mathcal{S}$ imply that its dual intersection complex satisfies conditions (2)(c) and (2)(d) of Definition \ref{polymani} respectively. 

One can likely generalize the results of this thesis to toric degenerations that don't satisfy conditions (3) and (4) of Assumption \ref{toricdegenerationass}. However, removing them would require a significant revision of \cite{GS3, GS4, GHS}. We don't feel that this generality is worth the added complexity of the language and prefer to work in the setup of \cite{GHS}.
\end{enumerate}
\end{remarks}

\subsubsection{Affine manifold} \label{affstr}

We assume that our polyhedral manifold $(B, \mathscr{P})$ is given a structure of an integral affine manifold with singularities. Recall that an \emph{affine manifold} is a topological manifold with transition maps in the group of affine transformations $$\Aff(N_{\mathbb{R}}) := N_{\mathbb{R}} \rtimes GL_n(\mathbb{R}) $$  of $N_{\mathbb{R}} = N \otimes_{\mathbb{Z}} \mathbb{R}$ for a fixed lattice $N \cong \mathbb{Z}^n$. We refer to such an atlas as an \emph{affine structure}. We say that the affine manifold is \emph{integral} if the transition maps lie in $$\Aff(N) :=N \rtimes GL_n(\mathbb{Z}).$$ 

Let $\tilde{\mathscr{P}}$ be the barycentric subdivision of $\mathscr{P}$ and let the \emph{discriminant} (or \emph{singular}) \emph{locus} $\Delta$ be the union of all $(n-2)$-dimensional cells of $\tilde{\mathscr{P}}$ that do not intersect the interiors of the maximal cells of $\mathscr{P}$. An \emph{integral affine manifold with singularities} (that we shall usually just call an \emph{affine manifold}) is a topological manifold $B$ with an integral affine structure on $B_0 := B \setminus \Delta$ (we write $i$ for the inclusion $i: B_0 \xhookrightarrow{} B$). We denote by $\underline{\rho} \subseteq \rho \in \mathscr{P}^{[n-1]}$ the connected components of $\rho \setminus \Delta$. We write $\underline{\rho} \in \tilde{\mathscr{P}}^{[n-1]}$ and call them \emph{slabs}.

\begin{remarks} \label{otherdelta}
In the future, we shall often need to modify the discriminant locus $\Delta$. 
\begin{enumerate}
\item Both for the dual intersection complex $\left( \bar{B}, \bar{\mathscr{P}} \right)$ of a toric degeneration $\bar{\mathfrak{X}} \to \mathcal{S}$ and for the dual intersection complex $(B, \mathscr{P})$ of a log CY degeneration $\mathfrak{X} \to \mathcal{S}$ we will actually choose a smaller discriminant locus (that will still be a union of submanifolds of codimension $\geq 2$), see Section \ref{tdaffstr} and Sections \ref{affstrsnc}, \ref{affstrgeneral} respectively. We will denote the discriminant locus of $\left( \bar{B}, \bar{\mathscr{P}} \right)$ by $\bar{\Delta}$ and denote the discriminant locus of $\left( B, \mathscr{P} \right)$ by $\Delta$ regardless of the exact situation.
\item In the case that we have a resolution $\pi: \mathfrak{X} \to \bar{\mathfrak{X}}$, the affine structure on $(B, \mathscr{P})$ will extend across each vertex corresponding to a toric component of the exceptional locus due to \cite[Proposition 2.3]{AG}. This corresponds to replacing $\Delta$ with a smaller discriminant locus and will enable us to define a PL-isomorphism $\Phi: B \to \bar{B}$ in general (see Constructions \ref{plmap} and \ref{plmapadmissible}).
\item In the definition of the affine structure on $\left( \bar{B}, \bar{\mathscr{P}} \right)$ we select $\bar{\Delta}$ in such a way that it does not contain any rational points of the integral structure (and we shall still denote the resulting subdivision of $\bar{\mathscr{P}}$ by $\tilde{\bar{\mathscr{P}}}$). 
\end{enumerate}
Neither choosing a smaller $\Delta$ nor modifying its precise location (as long as $\Delta$ respects the cell structure) changes the constructions in this section. We shall make explicit remarks when the exact location of $\Delta$ is important.
\end{remarks}

An integral affine manifold $B_0$ comes with a sheaf of integral (co-) tangent vectors $\Lambda = \Lambda_{B_0}$ (dually $\check{\Lambda} = \check{\Lambda}_{B_0}$). These are locally constant sheaves with stalks $\Lambda_x \cong \mathbb{Z}^n, \ \check{\Lambda}_x \cong \mathbb{Z}^n, \ x \in B_0$. Further, for $\tau \in \mathscr{P}$, we use the notation $\Lambda_{\tau}$ for the sheaf of integral tangent vectors on $\Int(\tau)$ (viewed as an integral affine manifold with boundary) or for the stalk of this sheaf at any $y \in \Int(\tau)$, depending on the context. If $\tau \subseteq \tau'$, we consider $\Lambda_{\tau}$ naturally as a subgroup of $\Lambda_{\tau'}$.

An important piece of data is the monodromy of the affine structure in the neighbourhood of $\rho \in \mathscr{P}^{[n-1]}$. Let $\underline{\rho}, \underline{\rho}' \subseteq \rho \in \mathscr{P}^{[n-1]}$ be two slabs and let $\sigma, \sigma' \in \mathscr{P}^{\max}$ be the maximal cells adjacent to $\rho$. Consider a simple loop $\gamma$ based at $x \in \Int \underline{\rho}$, passing successively into $\Int \sigma$, to $\Int \underline{\rho}'$, into $\Int \sigma'$, and back to $x$. Parallel transform along this path takes form
\begin{equation} \label{monodromy}
T_{\gamma}(m) = m + \check{d}_{\rho}(m) \cdot m_{\underline{\rho} \, \underline{\rho}'}, \quad m \in \Lambda_x
\end{equation}
where $m_{\underline{\rho} \, \underline{\rho}'} \in \Lambda_{\rho}$ and $\check{d}_{\rho} \in \check{\Lambda}_x$ is the generator of $\Lambda_{\rho}^{\perp} \subseteq \check{\Lambda}_{x}$ that takes non-negative values on $\sigma$. We call $m_{\underline{\rho} \, \underline{\rho}'}$ the monodromy vector. Note that $m_{\underline{\rho}'  \underline{\rho}} = -m_{\underline{\rho} \, \underline{\rho}'}$. 

We now describe a general procedure that we will use to define the structure of an integral affine manifold with singularities on a polyhedral manifold $(B, \mathscr{P})$.

\begin{cons} \label{affchartstool}
To define an integral affine structure on $B_0 =B \setminus \Delta$\footnote{Here, $\Delta$ should be understood as the union of all $(n-2)$-dimensional cells of $\tilde{\mathscr{P}}$ (with, possibly, the perturbed $\tilde{\mathscr{P}}$ of Remark \ref{otherdelta}(3)) that do not intersect the interiors of the maximal cells of $\mathscr{P}$. It will be easy to see when the affine structure defined via this construction extends over some codimension $2$ cells.} (viewed as a topological manifold only) it is enough to give:
\begin{enumerate}
\item To every 
\begin{equation} \label{polyh}
\left\{\sigma \ | \ \sigma \in \mathscr{P}^{\max} \right\}
\end{equation}
the structure of an integral polyhedron (so that for any $\rho \in \mathscr{P}^{[n-1]}$ such that $\rho = \sigma \cap \sigma'$ for $ \sigma, \sigma' \in \mathscr{P}^{\max}$ the structures of an integral polyhedron on $\rho$ induced from $\sigma$ and $\sigma'$ are the same). These define integral affine charts on 
\begin{equation} \label{polycharts}
\left\{\Int \sigma \ | \ \sigma \in \mathscr{P}^{\max} \right\}.
\end{equation}
Note that if $(B, \mathscr{P})$ is a polyhedral manifold, these structures are already provided by Definition \ref{polymani}. 
\item Structures of integral affine manifolds on every 
\begin{equation} \label{openstar}
\left\{W_v \setminus \Delta \ | \ v \in \mathscr{P}^{[0]} \right\}
\end{equation}
for 
\begin{equation} \label{openstardef}
W_v := \bigcup_{\left\{\tilde{\sigma} \in \tilde{\mathscr{P}} \ | \ v \subseteq \tilde{\sigma}\right\}} \Int \tilde{\sigma}
\end{equation}
that are compatible with the affine charts of \eqref{polycharts}.
\end{enumerate}

We often use the following \emph{refined description} of the affine structures on \eqref{openstar}.
\begin{enumerate}  
\item[$(2')$] To define the integral affine manifold structures on \eqref{openstar} it is enough to give, for any $\rho \in \mathscr{P}^{[n-1]}$ with $\rho = \sigma \cap \sigma'$ for $ \sigma, \sigma' \in \mathscr{P}^{\max}$ and any $v \subseteq \rho, \ v \in \mathscr{P}^{[0]}$, an integral piecewise-linear (PL) embedding 
\begin{equation} \label{affstrembeddings}
\psi_{\rho, v} : \sigma \cup \sigma' \to \mathbb{R}^n
\end{equation}
with $\psi_{\rho, v}(v) = 0$ compatible with the structures of affine polyhedra of \eqref{polyh}.
\end{enumerate}
\end{cons}

To see that any data of $(1)$ and $(2')$ in Construction \ref{affchartstool} defines a data of (2), note that for every $v \in \mathscr{P}^{[0]}$ the PL-embedding $\psi_{\rho, v} : \sigma \cup \sigma' \to \mathbb{R}^n$ of \eqref{affstrembeddings} restricts to a PL-embedding 
\begin{equation} \label{refinedembed}
\psi_{\underline{\rho}, v} : \tilde{\sigma} \cup \tilde{\sigma}' \to \mathbb{R}^n
\end{equation}
where $ \tilde{\sigma}, \tilde{\sigma}' \in \tilde{\mathscr{P}}^{\max}$ with $\tilde{\sigma} \cap \tilde{\sigma}' = \underline{\rho} \in \tilde{\mathscr{P}}^{[n-1]}$ are the unique cells of $\tilde{\mathscr{P}}$ such that $\underline{\rho} \subseteq \rho, \ \tilde{\sigma} \subseteq \sigma, \ \tilde{\sigma}' \subseteq \sigma'$, and $v \in \underline{\rho}$. Now the embedding $\psi_{\underline{\rho}, v}$ defines an affine chart on $\Int\left( \tilde{\sigma} \cup \tilde{\sigma}' \right)$ and the collection of such charts over all $\underline{\rho} \in \tilde{\mathscr{P}}^{[n-1]}$ with $v \in \underline{\rho}$, along with the charts on $$\left\{\Int \tilde{\sigma} \ | \ v \in \tilde{\sigma}, \   \tilde{\sigma} \in \tilde{\mathscr{P}}^{\max} \right\}$$ induced by restricting the charts of \eqref{polycharts}, defines the structure of an affine manifold on $W_v \setminus \Delta$. 

Conversely, assuming that for every $v \in \mathscr{P}^{[0]}$ there exist affine charts on $W_v \setminus \Delta$ covering $\Int\left( \tilde{\sigma} \cup \tilde{\sigma}' \right)$ for any $\tilde{\sigma}, \tilde{\sigma}' \in \tilde{\mathscr{P}}^{\max}$ with $\tilde{\sigma} \cap \tilde{\sigma}' = \underline{\rho} \in \tilde{\mathscr{P}}^{[n-1]}$ and $v \in \underline{\rho}$ (e.g. the affine structure on $W_v \setminus \Delta$ is given by a single chart), the data of (1) and (2) defines a data of $(2')$. Indeed, the affine charts induce PL-embeddings $\psi_{\underline{\rho}, v} : \tilde{\sigma} \cup \tilde{\sigma}' \to \mathbb{R}^n$ of \eqref{refinedembed}. Such a PL-embedding uniquely extends to a PL-embedding $\psi_{\rho, v} : \sigma \cup \sigma' \to \mathbb{R}^n$ of \eqref{affstrembeddings} for $\rho \in \mathscr{P}^{[n-1]}$ and $\sigma,\sigma' \in \mathscr{P}^{\max}$ the unique cells such that $\underline{\rho} \subseteq \rho, \ \tilde{\sigma} \subseteq \sigma, \ \tilde{\sigma}' \subseteq \sigma'$.

We will use Construction \ref{affchartstool} to define the affine structures on the dual intersection complex $\left( \bar{B}, \bar{\mathscr{P}} \right)$ of a toric degeneration $\bar{\mathfrak{X}} \to \mathcal{S}$ and the dual intersection complex $(B, \mathscr{P})$ of a log CY degeneration $\mathfrak{X} \to \mathcal{S}$ in Construction  \ref{affstrtoric} and Constructions \ref{affstrint}, \ref{affstrintgen} respectively. In the first case, we shall use the description of Construction \ref{affchartstool}(2). In the second case, we shall use the refined description of Construction \ref{affchartstool}$(2')$. We will use the above discussion on the equivalence of the two descriptions to define a PL-isomorphism $\Phi: B \to \bar{B}$ in increasing generality (see Constructions \ref{plisosmall}, \ref{plmap}, and \ref{plmapadmissible}).

\subsubsection{MPA function and the initial slab functions} \label{mpafunc}

The next ingredient of the setup is a convex multi-valued piecewise-affine (MPA) function $\varphi$ on $B$ with values in a \emph{toric} monoid. We fix a toric monoid $Q$, that is a finitely generated, integral, saturated monoid with $Q^{\gp}$ torsion-free.\footnote{Equivalently, we require $Q$ to be isomorphic to a finitely generated saturated submonoid of a finitely generated free abelian group. We do not require $Q$ to be \emph{sharp}, so it might have non-trivial invertible elements. So $Q$ can be realized as the set of integral points of some convex (but not necessarily strictly convex) rational polyhedral cone.} The choice of a $Q_{\mathbb{R}}^{\gp}$-valued MPA function (see \cite[Definition 1.2.5]{GHS}) $\varphi$ on $B$ consists of a choice of single-valued integral piecewise-linear (PL) functions $\varphi_U: U \to Q_{\mathbb{R}}^{\gp}$ in the charts of the affine structure well-defined up to linear functions.\footnote{This description is the reason that these functions are called multi-valued piecewise-linear (MPL) in \cite{GS4} and \cite{AG}.}

Let $\underline{\rho} \in \tilde{\mathscr{P}}^{[n-1]}$ be a slab contained in $\sigma, \sigma' \in \mathscr{P}^{\max}$ and let $\varphi_{U}: U \to Q_{\mathbb{R}}^{\gp}$ be the single-valued integral PL-function defined in a neighbourhood of $x \in \Int \underline{\rho}$. Let $n, n' \in \check{\Lambda}_x \otimes Q^{\gp}$ be the slopes of $\varphi_{U} |_{\sigma}$ and $\varphi_{U} |_{\sigma'}$. Then we may write 
\begin{equation} \label{kinks}
n'-n = \delta \cdot \kappa_{\underline{\rho}}
\end{equation}
where $\kappa_{\underline{\rho}} \in Q^{\gp}$ and $\delta: \Lambda_x \to \mathbb{Z}$ is the surjective map that vanishes on tangent vectors to $\underline{\rho}$ and is positive on tangent vectors pointing into $\sigma'$. Then $\kappa_{\underline{\rho}}$ (or $\kappa_{\underline{\rho}}(\varphi)$ if we want to emphasise the dependence on $\varphi$) is independent of $x \in \Int \underline{\rho}$ and we call it the \emph{kink} of $\varphi$ at $\underline{\rho}$. 

The collection of kinks $$\left\{ \kappa_{\underline{\rho}}(\varphi), \ \underline{\rho} \in \tilde{\mathscr{P}}^{[n-1]} \right\}$$ completely determines the MPA function $\varphi$ (see \cite[Proposition 1.2.6]{GHS})  and we will usually specify an MPA function by giving a collection of kinks. Unless otherwise specified, we require in addition that $\kappa_{\underline{\rho}} = \kappa_{\underline{\rho}'}=: \kappa_{\rho}$ for any $\underline{\rho}, \underline{\rho}' \subseteq \rho \in \mathscr{P}^{[n-1]}$. We fix an MPA function $\varphi$ and assume that it is \emph{convex}, that is $ \kappa_{\rho} \in Q \subseteq Q^{\gp}$ for all $\rho \in \mathscr{P}^{[n-1]}$. All MPA functions in this thesis will be convex. We say that $\varphi$ is \emph{strictly convex} if $ \kappa_{\rho} \in Q \setminus Q^{\times}$ for all $\rho \in \mathscr{P}^{[n-1]}$.

Fix a Noetherian ring $A$\footnote{We will always have $A = \kk$ until Chapter \ref{gslocus}.} and let $I_0 \subseteq Q$ be a monoid ideal. We assume that $\kappa_\rho \in I_0$ for any $\rho \in \mathscr{P}^{[n-1]}$.  We also fix an additional piece of discrete data $f_{\underline{\rho}} \in (A[Q]/I_0)[\Lambda_{\rho}]$ for each $\underline{\rho} \in \tilde{\mathscr{P}}^{[n-1]}$ such that for any $\underline{\rho}, \underline{\rho}' \subseteq \rho \in \mathscr{P}^{[n-1]}$ we have
\begin{equation} \label{comp}
f_{\underline{\rho}'} = z^{m_{\underline{\rho}'\underline{\rho} }}f_{\underline{\rho}}. 
\end{equation}
We call the elements of $$\left\{f_{\underline{\rho}} \ \Big| \ \underline{\rho} \in \tilde{\mathscr{P}}^{[n-1]} \right\}$$ the \emph{initial slab functions}.

\subsubsection{Monomials and the definition of a scattering diagram} \label{monomialsaffstr}
The choice of a $Q^{\gp}_{\mathbb{R}}$-valued 
MPA function $\varphi$ gives rise to a local system $\mathcal{P}$ (see \cite[Definition 1.2.12]{GHS}) on $B_0$ that fits into an exact sequence
\begin{equation} \label{monomialexactseq}
0\longrightarrow \underline{Q}^{\gp} \longrightarrow \mathcal{P} \longrightarrow \Lambda
\longrightarrow 0
\end{equation}
where $\underline{Q}^{\gp}$ denotes the constant sheaf with stalk $Q^{\gp}$. For every $m \in \mathcal{P}_x$, we write $\bar{m} \in \Lambda_x$ for the image of the projection. Moreover, each $x\in B_0$ gives a
submonoid $\mathcal{P}^+_x\subseteq \mathcal{P}_x$ (see \cite[Definition 1.2.13]{GHS}). We describe the monoids $\mathcal{P}^+_x$ and the effects of parallel transform on them explicitly. 

For $\sigma\in\mathscr{P}^{\max}$, $x\in\Int \sigma$, we have
\begin{equation} \label{posstalk}
\mathcal{P}^+_x= \Lambda_x \times Q.
\end{equation}

For $\rho\in\mathscr{P}^{[n-1]}$, $x\in \Int \rho \setminus \Delta$, we have
\begin{equation} \label{posstalk1}
\mathcal{P}^+_x=(\Lambda_{\rho}\oplus \mathbb{N} Z_+\oplus \mathbb{N} Z_-\oplus Q)/
\langle Z_++Z_-=\kappa_{\rho}\rangle.
\end{equation}
This description requires an ordering
$\sigma,\sigma'\in \mathscr{P}^{\max}$ of the maximal cells containing
$\rho$ and a choice of vector $\xi\in\Lambda_x$ pointing into
$\sigma$ and representing a generator of $\Lambda_{\sigma}/\Lambda_{\rho}$ (we call such a $\xi$ a \emph{normal generator}), so that $\bar{Z}_+ = \xi$ and $\bar{Z}_- = -\xi$.

For $x, x' \in \Int(\sigma), \ \sigma \in \mathscr{P}^{\max}$ the parallel transport $\mathcal{P}_x^+ \to \mathcal{P}_{x'}^+$ is trivial in the representation of (\ref{posstalk}). For $x\in \Int \rho \setminus \Delta$ and $y\in \Int \sigma$, $y'\in\Int \sigma'$, the parallel transports from $\mathcal{P}_x^+$ to $\mathcal{P}_y^+$ and $\mathcal{P}_{y'}^+$
in the local system $\mathcal{P}$ take form
\begin{equation} \label{pti}
(\lambda_{\rho},a Z_+,b Z_-,q) \mapsto \begin{cases} \big(\lambda_{\rho}+(a-b)\xi, q + b \kappa_{\rho}) \in \mathcal{P}_y^+ \\ \big(\lambda_{\rho}+(a-b)\xi, q + a \kappa_{\rho}) \in \mathcal{P}_{y'}^+  \end{cases}
\end{equation}
See the discussion of \cite[Section 2.2]{GHS} for more details.

More generally, let $x \in B_0$. Parallel transport inside the chart $U$ of the affine structure in the neighbourhood of $x \in \sigma, \ \sigma \in \mathscr{P}^{\max}$ identifies $\mathcal{P}_x$ with $\mathcal{P}_y$ for any $y \in \sigma \cap U$. Via this parallel transform we have 
\begin{equation} \label{stalkint}
\mathcal{P}_x^+ = \bigcap_{x \in \sigma \in \mathscr{P}^{\max}} \mathcal{P}_y^{+}.
\end{equation}

\begin{notation} \label{notation}
For $x\in\Int \sigma$, $\sigma\in\mathscr{P}_{\max}$, we will write a monomial
in $A[\mathcal{P}^+_x]=A[Q][\Lambda_x]$ either as $z^m$ for $m\in\mathcal{P}^+_x$
or as $t^q z^{\bar m}$
for $(\bar m, q)\in\Lambda_x \oplus Q$ via (\ref{posstalk}). Similarly, if $x\in \Int \rho \setminus \Delta, \ \rho\in \mathscr{P}^{[n-1]}$, 
we have a canonically defined submonoid
$\Lambda_{\rho} \oplus Q\subseteq \mathcal{P}^+_x$ via
\eqref{posstalk1} giving a subring
$A[Q][\Lambda_{\rho}]\subseteq A[\mathcal{P}^+_x]$. Again, we write the monomials in this subring as $t^qz^{\bar m}$ for $q\in Q, \ \bar m
\in\Lambda_{\rho}$. 
\end{notation}

For any monoid ideal $I \subseteq Q$ and any $x \in B_0$, we obtain an ideal $\mathcal{I}_x \subseteq A[\mathcal{P}_x^{+}]$. If $x \in \Int \sigma, \ \sigma \in \mathscr{P}^{\max}$, then $\mathcal{I}_x$ is defined via the inclusion of $A[Q]$ into $A[\mathcal{P}_x^{+}]$. More generally, using the description of \eqref{stalkint}, we can define 
\begin{equation} \label{idealint}
\mathcal{I}_x := \sum_{x \in \sigma \in \mathscr{P}^{\max}} \mathcal{I}_y \cap A[\mathcal{P}_y^{+}].
\end{equation}

Let $I$ be a monoid ideal with $\sqrt{I} = I_0$. 

\begin{defn} \label{wall}
A \emph{wall} on $(B,\mathscr{P})$ is a codimension one rational polyhedral subset $\mathfrak{p}$ of some $\sigma\in\mathscr{P}_{\max}$
along with a \emph{wall function}
\[
f_{\mathfrak{p}}=\sum_{m\in \mathcal{P}^+_x, \ \bar m\in\Lambda_{\mathfrak{p}}}c_m z^m\in A[\mathcal{P}_x^+]/\mathcal{I}_x
\]
for $x\in\Int \mathfrak{p}$. Let $y \in \mathfrak{p} \setminus \Delta$. If $c_m \neq 0$, we require that under the identification of $\mathcal{P}_y$ with $\mathcal{P}_x$ via
parallel transport inside $\sigma\setminus\Delta$ we have $m\in \mathcal{P}^+_y$. 
We further require that:
\begin{enumerate}
    \item If $\mathfrak{p} \cap \Int \sigma \neq \varnothing$, then $f_{\mathfrak{p}} = 1 \mod I_0$.
    \item If $\mathfrak{p} \subseteq \underline{\rho}$ for some $\underline{\rho} \in \tilde{\mathscr{P}}^{[n-1]}$, then $f_{\mathfrak{p}} = f_{\underline{\rho}} \mod I_0$.
\end{enumerate}
We denote a wall by either $(\mathfrak{p}, f_{\mathfrak{p}})$ or just $\mathfrak{p}$ and say that a wall is \emph{trivial} if $f_{\mathfrak{p}} = 1$.
\end{defn}

\begin{defn} \label{scatter}
A scattering diagram $\mathfrak{D}_I$ (or $(B, \mathfrak{D}_I)$ if we want to keep track of the affine manifold) on $(B,\mathscr{P})$ is a set of walls that is: 
\begin{enumerate}
    \item Finite.
    \item Forms the codimension one cells of a rational polyhedral decomposition $\mathscr{P}_{\mathfrak{D}_I}$ of $\mathscr{P}$ refining $\mathscr{P}$.
    \item Every $\underline{\mathcal{\rho}} \in \tilde{\mathscr{P}}^{[n-1]}$ is contained in a union of walls.\footnote{This is not required in \cite{GHS}, and is just a convention to make the notation easier, see Remark \ref{dec} and Definition \ref{combequiv}.} 
\end{enumerate}
We define the \emph{support} and \emph{singular locus} of $\mathfrak{D}_I$ as follows:
\begin{align*}
|\mathfrak{D}_I| &:= \bigcup_{\mathfrak{p}\in \mathfrak{D}_I} \mathfrak{p} \\
\Sing(\mathfrak{D}_{I}) &:= \Delta \cup \bigcup_{\mathfrak{p}, \mathfrak{p}' \in \mathfrak{D}_{I}} (\mathfrak{p} \cap \mathfrak{p}')
\end{align*}

We call elements $\mathfrak{u} \in \mathscr{P}_{\mathfrak{D}_I}^{\max}$ \emph{chambers}, elements $\mathfrak{b} \in \mathscr{P}_{\mathfrak{D}_I}^{[n-1]}$ such that $\mathfrak{b} \subseteq \underline{\rho}$ for $\rho \in \tilde{\mathscr{P}}^{[n-1]}$ \emph{slabs}\footnote{Note that we slightly abuse the notation here since we also call elements $\underline{\rho} \in \tilde{\mathscr{P}}^{[n-1]}$ slabs. The meaning will be clear from the context.} (and their wall functions \emph{slab functions}), and elements  $\mathfrak{j} \in \mathscr{P}_{\mathfrak{D}_I}^{[n-2]}$ \emph{joints}. The codimension $k$ of a joint $\mathfrak{j}$ or a wall $\mathfrak{p}$ is the codimension of the smallest cell of $\mathscr{P}$ containing $\mathfrak{j}$ or $\mathfrak{p}$ respectively. In particular, slabs are precisely the codimension $1$ walls.
\end{defn}

\begin{remark} \label{dec}
        The walls of \cite{GS1, GS4} often don't satisfy conditions (2) and (3) of Definition \ref{scatter}. Suppose that $\mathfrak{D}_{I}'$ is a scattering diagram in that sense and set \begin{align*}
|\mathfrak{D}'_I| &:= \bigcup_{\mathfrak{p}\in \mathfrak{D}'_I} \mathfrak{p} \cup \bigcup_{\underline{\rho}\in
\tilde{\mathscr{P}}^{[n-1]}} \underline{\rho}\\
\Sing(\mathfrak{D}'_{I}) &:= \Delta \cup \bigcup_{\mathfrak{p}, \mathfrak{p}' \in \mathfrak{D}'_{I}} (\mathfrak{p} \cap \mathfrak{p}')
\end{align*} where the union in $\Sing(\mathfrak{D}'_{I})$ is over all pairs of walls $\mathfrak{p},\mathfrak{p}'$ with
        $\mathfrak{p} \cap \mathfrak{p}'$ of codimension at least two. If $f_{\underline{\rho}} \neq 1$ for $\underline{\rho}\in
\tilde{\mathscr{P}}^{[n-1]}$, assume in addition that there are no two walls $\mathfrak{b}_1, \mathfrak{b}_2 \subseteq \underline{\rho}$ such that $f_{\mathfrak{b}_i} \neq 1 \mod I_0, \ i = 1, 2$  and $\dim \mathfrak{b}_1 \cap \mathfrak{b}_2 \cap \underline{\rho} = n-1$.\footnote{This assumption is necessary to ensure that step (2) below produces walls satisfying condition (2) of Definition \ref{wall}.} Now, do the following:
		\begin{enumerate}
			\item \begin{sloppypar}Decompose every wall ($\mathfrak{p}, f_{\mathfrak{p}})$ into walls $(\mathfrak{p}_1, f_{\mathfrak{p}_1}), \dots  ,(\mathfrak{p}_n, f_{\mathfrak{p}_n})$ where $\mathfrak{p}_1, \dots, \mathfrak{p}_n$ are the closures of the connected components of $\mathfrak{p} \setminus (\Sing(\mathfrak{D}_{I}') \cap \mathfrak{p})$ and $f_{\mathfrak{p}_i} := f_{\mathfrak{p}}$.\end{sloppypar}
			\item If after (1) there are multiple walls $(\mathfrak{p}_1, f_{\mathfrak{p}_1}), \dots  (\mathfrak{p}_m, f_{\mathfrak{p}_m})$ with the same support, replace them by one wall $(\mathfrak{p}', f_{\mathfrak{p}'})$ where $\mathfrak{p}'$ has the same support as $\mathfrak{p_i}$ and $f_{\mathfrak{p}'}:= \prod_{i = 1}^m f_{\mathfrak{p}_i}$.
			\item For every $\underline{\rho} \in \tilde{\mathscr{P}}^{[n-1]}$, introduce a slab $\mathfrak{b}$ for every connected component of $$\underline{\rho} \ \Big\backslash \left( \bigcup_{\mathfrak{p}\in \mathfrak{D}'_I} \mathfrak{p} \cup \Sing(\mathfrak{D}'_I)\right)$$ with slab function $f_{\mathfrak{b}} := f_{\underline{\rho}}$.
                \item Add trivial walls so that condition (2) of Definition \ref{scatter} is satisfied.
		\end{enumerate}
		The resulting scattering diagram $\mathfrak{D}_{I}$ satisfies Definition \ref{scatter}. We shall freely drop conditions (2) and (3) of Definition \ref{scatter} when convenient and assume that this procedure has taken place\footnote{The described procedure does not change the mirror family $\mathfrak{X}_{\mathfrak{D}_I}$ constructed from $\mathfrak{D}_I$. For \cite{GS1}, see the discussion of \cite[Appendix A.1]{GHS}. For \cite{GS4}, that is precisely how the construction goes, see \cite[Remark 3.4]{GS4}.} when using the setup of \cite{GHS}.
	\end{remark}

We will usually have to deal with a compatible collection of scattering diagrams over multiple ideals. 

\begin{defn} \label{compatibility}
The scattering diagrams $\mathfrak{D}_I$ and $\mathfrak{D}_{I'}$ for $I \subseteq I'$ are \emph{compatible} if for every wall $\mathfrak{p} \in \mathfrak{D}_{I'}$ we either have $\mathfrak{p} \in \mathfrak{D}_{I}$ and $f_{\mathfrak{p}, I} = f_{\mathfrak{p}, I'} \mod I$ or $\mathfrak{p} \in \mathfrak{D}_{I'} \setminus \mathfrak{D}_I$ and $f_{\mathfrak{p}} = 1 \mod I$.
\end{defn}

A scattering diagram $\mathfrak{D}_I$ gives rise to a canonical family $\check{\mathfrak{X}}_{\mathfrak{D}_I} \to \Spec A[Q]/I$ if $\mathfrak{D}_I$ satisfies the additional requirement of \emph{consistency}. We will now review the construction.

\subsubsection{Consistency in codimensions \texorpdfstring{$0$}{0} and \texorpdfstring{$1$}{1}, and the construction of \texorpdfstring{$\check{\mathfrak{X}}^{o}_{\mathfrak{D}_I}$}{Xo}} \label{gluedfamily}

First, we define a family $\check{X}_0 \to \Spec A[Q]/I_0$, which one should think of as the central fibre\footnote{It is indeed the central fibre of the family if $I_0 = \mathfrak{m}$ is the maximal ideal.} of $\check{\mathfrak{X}}_{\mathfrak{D}_I} \to \Spec A[Q]/I$. The construction is purely combinatorial. For any integral polyhedron $\sigma \in \mathbb{R}^k$ we define the cone over $\sigma$ as 
\begin{equation} \label{coneover}
  \mathbf{C}\sigma : = \mathbb{R}_{\geq 0} \cdot \left(\sigma \times \left\{1\right\}\right) \subseteq \mathbb{R} \times \mathbb{R}^k.
\end{equation}
Then for any $d \geq 0$ let 
\begin{equation} \label{ratpoints}
B\left(\frac{1}{d}\mathbb{Z}\right):= \bigcup_{\sigma \in \mathscr{P}^{\max}} \mathbf{C}\sigma \cap (\mathbb{Z}^n \times \left\{d \right\}).
\end{equation}
The notation is justified as $B\left(\frac{1}{d}\mathbb{Z}\right)$ is in bijection with the set of points of $B$ with denominator $d$ in some integral affine chart. Now let
\begin{equation} \label{ringcentralfibre}
\left(A[Q]/I_0\right)[B]:=\bigoplus_{d \in \mathbb{N}} \left(A[Q]/I_0\right)^{B\left(\frac{1}{d}\mathbb{Z}\right)}
\end{equation}
with basis elements $z^m, \ m \in B\left(\frac{1}{d}\mathbb{Z}\right)$ for some $d \in \mathbb{N}$. We define the multiplication as $z^m \cdot z^{m'} := z^{m'+m}$ if $m$ and $m'$ are in the same $\sigma \in \mathscr{P}^{\max}$ (the sum taken in $\mathbf{C}\sigma$) and $z^m \cdot z^{m'} := 0$ otherwise. This makes $\left(A[Q]/I_0\right)[B]$ into a ring and we have a natural $\mathbb{Z}$-grading by $d$. We define\footnote{In \cite{GHS}, this scheme is denoted $X_0$. We include a check here since we will use the construction to produce mirrors. The same remark applies to other objects constructed from $(B, \mathscr{P})$.} $$\check{X}_0 : = \Proj  \left(A[Q]/I_0\right)[B]$$ and assume that $\check{X}_0$ is projective.

We now want to deform $\check{X}_0$. We define rings $$R_{\sigma} := (A[Q]/I)[\Lambda_{\sigma}]$$ for every $\sigma \in \mathscr{P}^{\max}$ and $$R_{\underline{\rho}} := (A[Q]/I)[\Lambda_{\rho}][Z_+, Z_-]/(Z_+Z_- - f_{\underline{\rho}} \cdot z^{\kappa_{\rho}})$$ for every $ \underline{\rho} \in \tilde{\mathscr{P}}^{[n-1]}$, this arises from \eqref{posstalk} and \eqref{posstalk1}. 
We also define $$R_{\mathfrak{u}} := (A[Q]/I)[\Lambda_{\sigma}]$$ for every chamber $\mathfrak{u} \subseteq \sigma \in \mathscr{P}^{\max}$ and $$R_{\mathfrak{b}} := (A[Q]/I)[\Lambda_{\rho}][Z_+, Z_-]/(Z_+Z_- - f_{\mathfrak{b}} \cdot z^{\kappa_{\rho}})$$ for every slab $\mathfrak{b} \subseteq \underline{\rho} \in \tilde{\mathscr{P}}^{[n-1]}.$

For every $\mathfrak{b} \subseteq \mathfrak{u}$ we have \emph{localization homomorphisms} 
\begin{equation} \label{lochom}
\chi_{\mathfrak{b}, \mathfrak{u}}: R_{\mathfrak{b}} \to R_{\mathfrak{u}}
\end{equation} given by (this arises from \eqref{pti}): $$Z_+ \mapsto z^\xi, \quad Z_- \mapsto f_{\mathfrak{b}}z^{\kappa_{\rho}} z^{-\xi}$$ with $\xi \in \Lambda_{\sigma}$ the normal generator pointing into $\sigma$ as in \eqref{posstalk1} and other monomials identified using $\Lambda_{\rho} \subseteq \Lambda_{\sigma}$.

Let $\mathfrak{p}$ be a codimension $0$ wall separating two chambers $\mathfrak{u}, \mathfrak{u}' \subseteq \sigma \in \mathscr{P}^{\max}$. Let $n_{\mathfrak{p}}$ be the generator of $\check{\Lambda}_{\mathfrak{p}} \subseteq \check{\Lambda}_x$ for some $x \in \Int \mathfrak{p}$ that is non-negative on $\mathfrak{u}$. We have \emph{wall-crossing homomorphisms} 
\begin{equation} \label{wallcross}
\theta_{\mathfrak{p}}: R_{\mathfrak{u}} \to R_{\mathfrak{u}'}, \quad z^m \mapsto f_{\mathfrak{p}}^{\langle n_{\mathfrak{p}}, \bar{m} \rangle} z^m.
\end{equation}

\begin{defn} \label{codim0consistency}
Let $\mathfrak{j} \subseteq \sigma \in \mathscr{P}^{\max}$ be a codimension $0$ joint and suppose that it is contained in the set of walls $\mathfrak{p}_1, \dots, \mathfrak{p}_r$ ordered cyclically. $\mathfrak{D}_I$ is consistent around $\mathfrak{j}$ if $$\theta_{\gamma_{\mathfrak{j}}}: = \theta_{\mathfrak{p}_1} \circ \dots  \circ \theta_{\mathfrak{p}_r} = \Id$$ as an automorphism of $R_{\sigma}$.\footnote{The notation $\theta_{\gamma_{\mathfrak{j}}}$ comes from thinking of the composition as going around $\mathfrak{j}$ via a small loop $\gamma_{\mathfrak{j}}$ and multiplying by a wall-crossing homomorphism each time we cross a wall.} $\mathfrak{D}_I$ is \emph{consistent in codimension $0$} if it is consistent around every codimension $0$ joint.
\end{defn}

\begin{defn} \label{codim1consistency}
Let $\mathfrak{j} \subseteq \rho  \in \mathscr{P}^{[n-1]}$ be a codimension $1$ joint. Then we have unique $\sigma, \sigma' \in \mathscr{P}^{\max}$ such that $\rho = \sigma \cap \sigma'$ and unique slabs $\mathfrak{b}_1, \mathfrak{b}_2 \subseteq \rho$ such that $\mathfrak{j} = \mathfrak{b}_1 \cap \mathfrak{b}_2$. Suppose that $\mathfrak{j}$ is contained in a set of walls $\mathfrak{p}_1, \dots, \mathfrak{p}_r, \mathfrak{p}'_1, \dots, \mathfrak{p}'_s$ such that $\mathfrak{p}_1, \dots, \mathfrak{p}_r \subseteq \sigma$, $\mathfrak{p}'_1, \dots, \mathfrak{p}'_s \subseteq \sigma'$, and $\mathfrak{b}_1, \mathfrak{p}_1, \dots, \mathfrak{p}_r, \mathfrak{b}_2, \mathfrak{p}'_1, \dots, \mathfrak{p}'_s$ is a cyclic ordering around $\mathfrak{j}$. Then there are homomorphisms\footnote{Obtained by composing \eqref{lochom} with $\Id: R_{\mathfrak{u}} \to R_{\sigma}$ for any $\mathfrak{u} \subseteq \sigma \in \mathscr{P}^{\max}$.} $$\chi_{\mathfrak{b}_i, \sigma}: R_{\mathfrak{b}_i} \to R_{\sigma}, \quad \chi_{\mathfrak{b}_i, \sigma}: R_{\mathfrak{b}_i} \to R_{\sigma'}, \quad i=1,2$$ and compositions of wall-crossings 
$$
\begin{aligned}
\theta &:= \theta_{\mathfrak{p}_r} \circ \theta_{\mathfrak{p}_{r-1}} \circ \dots \circ \theta_{\mathfrak{p}_{1}}: \ &R_{\sigma} \to R_{\sigma} \\
\theta' &:= \theta_{\mathfrak{p}'_1} \circ \theta_{\mathfrak{p}'_{2}} \circ \dots \circ \theta_{\mathfrak{p}'_{s}}: \ &R_{\sigma'} \to R_{\sigma'}
\end{aligned}
$$ $\mathfrak{D}_I$ is consistent around $\mathfrak{j}$ if $$(\theta \times \theta')\left((\chi_{\mathfrak{b}_1, \sigma}, \chi_{\mathfrak{b}_1, \sigma'})(R_{\mathfrak{b}_1})\right) = (\chi_{\mathfrak{b}_2, \sigma}, \chi_{\mathfrak{b}_2, \sigma'})(R_{\mathfrak{b}_2}).$$ In this case we define 
\begin{equation}
    \theta_{\mathfrak{j}}: R_{\mathfrak{b}_1} \to R_{\mathfrak{b}_2}
\end{equation}
as the isomorphism induced by $\theta \times \theta'$.\footnote{It is easy to see that the homomorphisms $(\chi_{\mathfrak{b}_i, \sigma}, \chi_{\mathfrak{b}_i, \sigma'}), \ i=1,2$ are injective.} $\mathfrak{D}_I$ is \emph{consistent in codimension $1$} if it is consistent around every codimension $1$ joint. 
\end{defn}

Given a scattering diagram $\mathfrak{D}_I$ consistent in codimensions $0$ and $1$, \cite[Proposition 2.4.1]{GHS} constructs a family $\check{\mathfrak{X}}^o_{\mathfrak{D}_I} \to \Spec A[Q]/I$ by gluing together the various $\Spec R_{\mathfrak{u}}$ for $\mathfrak{u}$ a chamber of $\mathfrak{D}_I$ and $\Spec R_{\mathfrak{b}}$ for $\mathfrak{b}$ a slab of $\mathfrak{D}_I$ along the localization homomorphisms $\chi_{\mathfrak{b}, \mathfrak{u}}$ of \eqref{lochom}, the wall-crosssing homomorphisms $\theta_{\mathfrak{p}}$ of \eqref{wallcross} and the isomorphisms $\theta_{\mathfrak{j}}$ associated to crossing a codimension $1$ joint $\mathfrak{j}$. Consistency conditions ensure that the gluing is well-defined. The reduction of $\check{\mathfrak{X}}^o_{\mathfrak{D}_I}$ modulo $I_0$ is canonically isomorphic (see \cite[Proposition 2.4.4]{GHS}) to the complement of the union of codimension $2$ strata of $\check{X}_0$.

\subsubsection{Consistency in codimension \texorpdfstring{$2$}{2} and the extension to \texorpdfstring{$\check{\mathfrak{X}}_{\mathfrak{D}_I}$}{Xo}} \label{thetafunctions}

Consistency in codimension $2$ is more complicated than in codimensions $0$ and $1$ and roughly means that one has a well-defined basis $\vartheta_m$ of sections of $\Gamma(\check{\mathfrak{X}}_{\mathfrak{D}_I}^o, \mathcal{O}_{\check{\mathfrak{X}}_{\mathfrak{D}_I}^o})$ called \emph{theta functions}. We review the construction, see \cite[Sections 3.1 and 3.2]{GHS} for details.

First, we want to generalize the notion of a wall-crossing homomorphism \eqref{wallcross} to slabs. This is a bit more delicate. Let $\mathfrak{b}$ be a slab separating two chambers $\mathfrak{u} \subseteq \sigma \in \mathscr{P}^{\max}$ and $\mathfrak{u}' \subseteq \sigma' \in \mathscr{P}^{\max}$. Let $R_{\mathfrak{u}}^{\mathfrak{b}} \subseteq R_{\mathfrak{u}}$ be the $A[Q]/I$ subalgebra generated by $\Lambda_{\rho}$ and the $\chi_{\mathfrak{b}, \mathfrak{u}}(Z_{+})$ of \eqref{lochom}. Now define a \emph{slab-crossing homomorphism}
\begin{equation} \label{slabcrossing}
\theta_{\mathfrak{b}}:R_{\mathfrak{u}}^{\mathfrak{b}} \to R_{\mathfrak{u}'}
\end{equation}
as an identity on $\Lambda_{\rho}$ and by setting $$\theta_{\mathfrak{b}}(\chi_{\mathfrak{b}, \mathfrak{u}}(Z_{+})) := \chi_{\mathfrak{b}, \mathfrak{u}'}(Z_{-}) \cdot f_{\mathfrak{b}} \cdot z^{\kappa_{\rho}}.$$ 

Note that the expression $\theta_{\mathfrak{b}}(a z^{m})$ is defined for any $a z^{m} \in R_{\mathfrak{u}}$ that is defined at a point of $\Int(\mathfrak{u} \cap \mathfrak{u}')$ via parallel transport. 

Now we can define broken lines that can be roughly described as piecewise-linear paths in $B_0$ that bend when they cross a wall $\mathfrak{p}$ in a way specified by $f_{\mathfrak{p}}$. Formally, we have the following.

\begin{defn} \label{brokenlines}
A \emph{broken line} for a scattering diagram $\mathfrak{D}_I$ on $(B, \mathscr{P})$ is a proper continuous map $$\beta: (-\infty, 0] \to B \setminus \Sing(\mathfrak{D}_I)$$ along with:
\begin{itemize}
    \item A sequence $-\infty = t_0 < t_1 < \dots < t_r = 0$ for some $r \geq 1$ such that $\beta(t_i) \in |\mathfrak{D}_I|$ for $i \leq r-1$. We say that $p :=\beta(0)$ is the \emph{endpoint} of $\beta$.
    \item For each $1 \leq i \leq r$, an expression $a_iz^{m_i}$ with $a_i \in A[Q]/I, \ m_i \in \Lambda_{\beta(t)}$ for some $t \in (t_{i-1}, t_i)$, defined at all points of $\beta([t_{i-1}, t_i])$. We require $a_1 = 1$. We call $m_1$ the \emph{asymptotic monomial} of $\beta$ and introduce the notation: $$a_{\beta}:=a_r, \quad m_{\beta}:=m_r.$$
\end{itemize}
This information is subject to the following conditions:
\begin{enumerate}
    \item $\beta|_{(t_{i-1}, t_i)}$ is a non-constant affine map with image contained in a unique chamber $\mathfrak{u}$ and $\beta'(t) = -m_i$ for all $t \in (t_{i-1}, t_i)$. 
    \item If $\beta(t_i)$ for $1 \leq i \leq r-1$ lies in a codimension $0$ wall $\mathfrak{p} \subseteq \mathfrak{u} \cap \mathfrak{u}'$ for chambers $\mathfrak{u}, \mathfrak{u}' \subseteq \sigma \in \mathscr{P}^{\max}$ and $\beta|_{(t_{i-1}, t_i)}$ is contained in $\mathfrak{u}$, then we require that $a_{i+1}z^{m_{i+1}}$ is one of the $a_j z^{m_j}$ in $$\theta_{\mathfrak{p}}(a_i z^{m_i}) = \sum_j a_j z^{m_j}$$ where $\theta_{\mathfrak{p}}$ is the wall-crossing homomorphism of \eqref{wallcross}. 
    \item If $\beta(t_i)$ for $1 \leq i \leq r-1$ lies in a slab $\mathfrak{b} \subseteq \underline{\rho} \in \tilde{\mathscr{P}}^{\max}$, $\mathfrak{b} \subseteq \mathfrak{u} \cap \mathfrak{u}'$ is contained in the chambers $\mathfrak{u} \subseteq \sigma \in \mathscr{P}^{\max}, \ \mathfrak{u}' \subseteq \sigma' \in \mathscr{P}^{\max} $, and $\beta|_{(t_{i-1}, t_i)}$ is contained in $\mathfrak{u}$, we require that\footnote{This is well-defined since $\beta(t_i) \in \Int \mathfrak{b}$ as $\im(\beta)$ is disjoint from the joints.} $a_{i+1}z^{m_{i+1}}$ is one of the $a_j z^{m_j}$ in $$\theta_{\mathfrak{b}}(a_i z^{m_i}) = \sum_j a_j z^{m_j}$$ where $\theta_{\mathfrak{b}}$ is the slab-crossing homomorphism of \eqref{slabcrossing}. 
\end{enumerate}
In either (2) or (3), we say that $a_{i+1}z^{m_{i+1}}$ is the \emph{result of transport of $a_i z^{m_i}$ from $\mathfrak{u}$ to $\mathfrak{u}'$}.
\end{defn}

\begin{notation} \label{wallcrossnotation}
    Suppose that $\mathfrak{D}_I$ is a scattering diagram consistent in codimensions $0$ and $1$ and let $\mathfrak{u}, \mathfrak{u}' \in \mathscr{P}^{\max}$ be two adjacent chambers. We introduce a single notation $\theta_{\mathfrak{u}'\mathfrak{u}}$ for the \emph{change of chamber homomorphism} $\theta_{\mathfrak{p}}$ of \eqref{wallcross} or $\theta_{\mathfrak{b}}$ of \eqref{slabcrossing} when $\mathfrak{p} \subseteq \mathfrak{u} \cap \mathfrak{u}'$ (resp. $\mathfrak{b} \subseteq \mathfrak{u} \cap \mathfrak{u}'$).\footnote{Consistency in codimensions $0$ and $1$ implies that this does not depend on the choice of $\mathfrak{p} \subseteq \mathfrak{u} \cap \mathfrak{u}'$ (resp. $\mathfrak{b} \subseteq \mathfrak{u} \cap \mathfrak{u}'$). See the proof of \cite[Proposition 2.4.1]{GHS} and the reference to \cite[Lemma 2.30]{GS1} therein.}
\end{notation}

Let $\mathfrak{j}$ be a joint of codimension $2$ and let $\omega \in \mathscr{P}^{[n-2]}$ be the smallest cell containing $\mathfrak{j}$. Build a new affine manifold $(B_{\mathfrak{j}}, \mathscr{P}_{\mathfrak{j}})$ by replacing any $\tau \in \mathscr{P}$ with $\mathfrak{j} \subseteq \tau$ by the tangent wedge of $\omega$ in $\tau$. Note that the inclusion of faces $\tau \subseteq \tau'$ induces an inclusion of the corresponding tangent wedges. So $B_{\mathfrak{j}}$ is a local model of $(B, \mathscr{P})$ near $\mathfrak{j}$ all of whose cells are cones. Similarly, the scattering diagram $\mathfrak{D}_{I}$ induces a scattering diagram $\mathfrak{D}_{I, \mathfrak{j}}$ by considering only the walls containing $\mathfrak{j}$ and going over to the tangent wedges based at $\omega$ for the underlying polyhedral subsets of codimension one. Now, let $m$ be an asymptotic monomial on $(B_{\mathfrak{j}}, \mathscr{P}_{\mathfrak{j}})$. For a general point $p \in B_{\mathfrak{j}}$, say contained in a chamber $\mathfrak{u} \in \mathscr{P}_{\mathfrak{j},\mathfrak{D}_{I, \mathfrak{j}}}^{\max}$, define 
\begin{equation} \label{thetafunc}
\vartheta_m^{\mathfrak{j}}(p):=\sum_{\beta} a_{\beta} z^{m_{\beta}} \in R_{\mathfrak{u}}
\end{equation} where the sum is over all the broken lines on $(B_{\mathfrak{j}}, \mathscr{P}_{\mathfrak{j}})$ with asymptotic monomial $m$ and endpoint $p$.

\begin{defn} \label{consistency}
A scattering diagram $\mathfrak{D}_I$ is \emph{consistent along the codimension $2$ joint $\mathfrak{j}$} if the $\vartheta_m^{\mathfrak{j}}(p)$:
\begin{enumerate}
    \item Do not depend on the choice of a general point $p$ in the same chamber $\mathfrak{u}$.
    \item Are compatible with the change of chamber homomorphisms $\theta_{\mathfrak{u}' \mathfrak{u}}$ for $\mathfrak{D}_{I, \mathfrak{j}}$ as defined in Notation \ref{wallcrossnotation}.
\end{enumerate}
A scattering diagram is \emph{consistent in codimension $2$} if it is consistent along every codimension two joint. A scattering diagram is \emph{consistent} if it is consistent in codimensions $0$ (Definition \ref{codim0consistency}), $1$ (Definition \ref{codim1consistency}), and $2$.
\end{defn}

One can show that (see \cite[Theorem 3.3.1]{GHS}) if $\mathfrak{D}_I$ is consistent, then the local expressions $\vartheta_m^{\mathfrak{j}}(p)$ give rise to canonical globally defined functions $\vartheta_m \in \Gamma(\check{\mathfrak{X}}_{\mathfrak{D}_I}^o, \mathcal{O}_{\check{\mathfrak{X}}_{\mathfrak{D}_I}^o})$ that form an $A[Q]/I$ basis of $\Gamma(\check{\mathfrak{X}}_{\mathfrak{D}_I}^o, \mathcal{O}_{\check{\mathfrak{X}}_{\mathfrak{D}_I}^o})$. 

The construction of $\check{\mathfrak{X}}_{\mathfrak{D}_I} \to \Spec A[Q]/I$ proceeds as follows. One defines (see \cite[Definition 4.2.1]{GHS}) a new topological manifold with a polyhedral decomposition, the \emph{cone over $B$}, as a pair $(\mathbf{C}B, \mathbf{C}\mathscr{P})$ by taking the limit over $\left\{\mathbf{C}\sigma \ | \ \sigma \in \mathscr{P}\right\}$ with $\mathbf{C}\sigma$ defined as in \eqref{coneover}. Now, one can lift a consistent scattering diagram $(B, \mathfrak{D}_{I})$ to a consistent scattering diagram $(\mathbf{C}B, \mathbf{C}\mathfrak{D}_I)$ (see \cite[Definition 4.2.4 and Proposition 4.2.6]{GHS}). Following Section \ref{gluedfamily}, $\mathbf{C}\mathfrak{D}_I$ gives rise to a family $\check{\mathfrak{N}}_{\mathbf{C}\mathfrak{D}_I}^{o} \to \Spec A[Q]/I$. Then let $$\check{\mathfrak{X}}_{\mathfrak{D}_I} := \Proj \Gamma(\check{\mathfrak{N}}_{\mathbf{C}\mathfrak{D}_I}^o,\mathcal{O}_{\check{\mathfrak{N}}_{\mathbf{C}\mathfrak{D}_I}^o}) $$ along with the natural flat morphism to $\Spec A[Q]/I$. By \cite[Theorem 4.3.2(c)]{GHS}, there is a canonical embedding of $\check{\mathfrak{X}}_{\mathfrak{D}_I}^{o}$ into $\check{\mathfrak{X}}_{\mathfrak{D}_I}$ as an open dense subscheme and the reduction of $\check{\mathfrak{X}}_{\mathfrak{D}_I}$ modulo $I_0$ is canonically isomorphic to $\check{X}_0$. The product of the theta functions on $\mathbf{C}B$ forming the basis of $\Gamma(\check{\mathfrak{N}}_{\mathbf{C}\mathfrak{D}_I}^o, \mathcal{O}_{\check{\mathfrak{N}}_{\mathbf{C}\mathfrak{D}_I}^o})$ can be given explicitly in terms of the broken lines on $\mathbf{C}B$, see \cite[Theorem 3.5.1]{GHS}.

\subsubsection{Equivalence of scattering diagrams}

We need a way to detect whether two families $\check{\mathfrak{X}}_{(B, \mathfrak{D}_I)}$ and $\check{\mathfrak{X}}_{(B', \mathfrak{D}'_I)}$\footnote{Possibly on different affine manifolds but using the same general setup (i.e. the same $Q,  A,  I_0,  I$).} are isomorphic by either looking at the scattering diagrams or using the gluing setup of Section \ref{gluedfamily}.

\begin{defn} \label{equiv}
We say that two consistent scattering diagrams $(B, \mathfrak{D}_I)$ and $(B', \mathfrak{D}'_I)$ are \emph{equivalent} if the families $\check{\mathfrak{X}}^o_{(B, \mathfrak{D}_I)}$ and $\check{\mathfrak{X}}^o_{(B', \mathfrak{D}'_I)}$ constructed as in Section \ref{gluedfamily} are isomorphic. 
\end{defn}

Equivalence guarantees that the projective families constructed from consistent scattering diagrams coincide. 

\begin{fact} \label{help}
Suppose that $(B, \mathfrak{D}_I)$ and $(B', \mathfrak{D}'_I)$ are two equivalent consistent scattering diagrams and that there is an isomorphism 
$$\alpha:(B, \mathscr{P}) \to \left(B', \mathscr{P}'\right)$$ of polyhedral manifolds. Then $\check{\mathfrak{X}}_{(B, \mathfrak{D}_I)} \cong \check{\mathfrak{X}}_{(B', \mathfrak{D}'_I)}$.
\end{fact}

\begin{proof}
The construction of $\check{X}_0$ and the fact that for any $\sigma \in \mathscr{P}^{\max}$, $\sigma \cong \alpha(\sigma)$ as integral polyhedra imply that $\check{X}_0$ is the same using either scattering diagram. Equivalence says that the two families agree away from the codimension $2$ strata of $\check{X}_0$. By \cite[Proposition 2.1.6]{GHS}, both $\mathcal{O}_{\check{\mathfrak{X}}_{(B, \mathfrak{D}_I)}}$ and $\mathcal{O}_{\check{\mathfrak{X}}_{(B', \mathfrak{D}'_I)}}$ are sheaves on $\check{X}_0$ satisfying Serre's $S_2$ condition which are canonically isomorphic on $\check{X}_0^o$, the complement of the codimension $2$ strata of $\check{X}_0$. But now also $\check{\mathfrak{X}}_{(B, \mathfrak{D}_I)} \cong \check{\mathfrak{X}}_{(B', \mathfrak{D}'_I)}$ canonically since $$\mathcal{O}_{\check{\mathfrak{X}}_{(B, \mathfrak{D}_I)}} = i_{*} \mathcal{O}_{\check{\mathfrak{X}}^o_{(B, \mathfrak{D}_I)}} = i_{*} \mathcal{O}_{\check{\mathfrak{X}}^o_{(B', \mathfrak{D}'_I)}} = \mathcal{O}_{\check{\mathfrak{X}}_{(B', \mathfrak{D}'_I)}}$$ by the $S_2$ condition.
\end{proof}

\begin{defn} \label{combequiv}
Let $(B, \mathfrak{D}_I)$ be a scattering diagram (that does not necessarily satisfy conditions (2) and (3) of Definition \ref{scatter}, see the discussion of Remark \ref{dec}). 
If $x \in B\setminus\Sing(\mathfrak{D}_I)$, we define
$$
f_x:= \prod_{x\in\mathfrak{p}\in \mathfrak{D}_I} f_{\mathfrak{p}}.
$$
We say that two scattering diagrams $(B, \mathfrak{D}_I)$ and $(B, \mathfrak{D}'_I)$ (on the same affine manifold with singularities $B$) are \emph{combinatorially equivalent} if
$f_x=f'_x$ for all $x\in B\setminus (\Sing(\mathfrak{D}_I)\cup \Sing(\mathfrak{D}'_{I}))$.
\end{defn}

\begin{obs} \label{combimpliesequiv}
We make a few observations:
\begin{enumerate}
    \item Definition \ref{combequiv} agrees with the definition of equivalence of \cite[Definition 3.5]{GS4} and the definition of equivalence of \cite[Definition 3.3]{GS1}.\footnote{The notation of \cite{GS1} is rather different.}
    \item If $(B, \mathfrak{D}_I)$ and $(B, \mathfrak{D}'_I)$ are combinatorially equivalent and consistent, then they are equivalent in the sense of Definition \ref{equiv} (after applying the construction of Remark \ref{dec}). This is easy to see from the gluing construction of Section \ref{gluedfamily}.
    \item Adding a finite number of walls $\mathfrak{p}$ with $f_{\mathfrak{p}} = 1$ to a scattering diagram produces a combinatorially equivalent scattering diagram and does not change consistency (again, with the construction of Remark \ref{dec} assumed). So it produces an equivalent scattering diagram by (2).
\end{enumerate}
\end{obs}

\subsection{Toric degeneration setup} \label{tdsetup}

After slightly generalizing the notion of \emph{tropicalization} of a log scheme from \cite[Section 2.1.4]{ACGS}, we summarize the mirror reconstruction setup for toric degenerations. See \cite{G1} for an overview, \cite{GS2, GS1} for the full details, and Appendix \ref{projcase} for an explicit example. We also introduce and discuss the definition of being a \emph{special} toric degeneration (i.e. a degeneration for which we stated Conjecture \ref{conjecture}).

\subsubsection{Tropicalization} \label{tropicalization}
We recall the construction of the tropicalization of a Zariski log scheme. The construction of the tropicalization $\Sigma(X)$ of a Zariski log scheme $X$ of \cite[Section 2.1.4]{ACGS} requires that the log structure $\mathcal{M}_X$ on $X$ is fine and saturated. Since the total space $\bar{\mathfrak{X}}$ of a toric degeneration $\bar{\mathfrak{X}} \to \mathcal{S}$ only has a fine and saturated log structure away from the discriminant locus $Z$, we need to generalize the construction slightly.

Recall that \emph{logarithmic strata} of a fine and saturated Zariski log scheme $X$ are the connected components of the subsets where the ghost sheaf $\overline{\mathcal{M}}_{X}$ of the log structure $\mathcal{M}_{X}$ is constant. 

\begin{defn} \label{generalizedlogstrata}
Let $X$ be a Zariski log scheme such that the log structure $\mathcal{M}_{X}$ on $X$ is fine and saturated away from a subset $Z$ of codimension at least $2$. We say that $Y$ is a \emph{logarithmic stratum} of $X$ if:
\begin{enumerate}
    \item $Y \setminus Z$ is a logarithmic stratum $Y'$ of $X \setminus Z$ of the same dimension and the closure of $Y$ in $X$ coincides with the closure of $Y'$ in $X$.
    \item $\overline{\mathcal{M}}_{X}$ is constant on $Y \setminus Z$.
    \item $Y$ is a maximal subset satisfying (1) and (2).
\end{enumerate} 
\end{defn}

Note that Definition \ref{generalizedlogstrata} agrees with the usual definition if $Z = \varnothing$ and that the logarithmic strata of the total space $\bar{\mathfrak{X}}$ of a toric degeneration $\bar{\mathfrak{X}} \to \mathcal{S}$ agree with the toric strata. 

\begin{cons} \label{tropicalize}
We construct the \emph{tropicalization} $\Sigma(X)$ of a Zariski log scheme $X$ of finite type and with log structure fine and saturated away from a subset $Z$ of codimension at least $2$ following \cite[Section 2.1.4]{ACGS}. 
For the generic point $\eta$ of a logarithmic stratum of $X$ let $$\sigma_{\eta}:= \Hom(\overline{\mathcal{M}}_{X, \eta}, \mathbb{N})_{\mathbb{R}} \subseteq \Hom(\overline{\mathcal{M}}_{X, \eta}, \mathbb{Z})_{\mathbb{R}}$$ be the corresponding rational polyhedral cone. If $\eta$ is a specialization of $\eta'$, then there is a well-defined generization map $\overline{\mathcal{M}}_{X, \eta} \to \overline{\mathcal{M}}_{X, \eta'}$ since $X$ carries a Zariski log structure. Dualizing gives a face morphism $\sigma_{\eta'} \to \sigma_{\eta}$. 

Going over all the logarithmic strata of $X$ gives a diagram of rational polyhedral cones indexed by the strata with face morphisms. Thus, it defines a \emph{generalized cone complex} $\Sigma(X)$ with topological realization $$|\Sigma(X)| := \lim_{\substack{\longrightarrow\\\eta}} \sigma_{\eta}$$ (here the colimit is over all the generic points of logarithmic strata of $X$). It is clear from the construction that this agrees with \cite[Section 2.1.4]{ACGS} in the case that $Z = \varnothing$. 

The construction is functorial. Given a morphism of log schemes $f: X \to Y$, the map $f: f^{-1} \overline{\mathcal{M}}_Y \to \overline{\mathcal{M}}_X$ induces a morphism of generalized cone complexes 
\begin{equation} \label{tropmap}
f_{\trop}: \Sigma(X) \to \Sigma(Y).
\end{equation}

Note that $|\Sigma(X)|$ has a natural integral structure by setting $|\Sigma(X)|(\mathbb{Z})$ to be the set of points of $|\Sigma(X)|$ that are the images of integer points of some $\sigma_{\eta}$ under the canonical map $\sigma_{\eta} \to |\Sigma(X)|$. Moreover, for a morphism of log schemes $f: X \to Y$, the tropical morphism \eqref{tropmap} induces a map $|\Sigma(X)|(\mathbb{Z}) \to |\Sigma(Y)|(\mathbb{Z})$.

We refer to \cite[Section 2.1]{ACGS} for more details on the construction.
\end{cons}

\begin{remarks} \label{tropremarks}
\begin{enumerate}
    \item Note that if $\bar{\mathfrak{X}} \to \mathcal{S}$ is a toric degeneration, then the definitions immediately imply that $\Sigma(\bar{\mathfrak{X}}) \cong \Sigma(\bar{\mathfrak{X}} \setminus Z)$.
    \item Construction \ref{tropicalize} gives a unified view on tropicalizing toric degenerations $\bar{\mathfrak{X}} \to \mathcal{S}$ and log smooth degenerations $\mathfrak{X} \to \mathcal{S}$. However, it does not produce well-behaved tropicalizations in general. For example, suppose that $\mathfrak{X}' \to \mathcal{S}$ satisfies the assumptions of Definition \ref{toricdegeneration} of a toric degeneration but in condition 4(a) instead of requiring that $Z$ does not contain the image of a toric stratum under $\nu$, we let $Z$ be the union of (images of) toric strata of codimension $\geq 2$. Then $\Sigma(\mathfrak{X}')$ is two-dimensional and does not reflect the combinatorics of the central fibre $\mathfrak{X}'_0$ of $\mathfrak{X}'$. We shall only use Construction \ref{tropicalize} for toric degenerations and their (possibly partial) resolutions.
\end{enumerate}
\end{remarks}

\begin{defn} \label{simplelogscheme}
We say that $X$ is a \emph{simple log scheme} if for every $\sigma \in \Sigma(X)$ the canonical map $\sigma \to |\Sigma(X)|$ is injective.
\end{defn}

We shall always work with simple log schemes (recall that for a toric degeneration $\bar{\mathfrak{X}} \to \mathcal{S}$, Assumption \ref{toricdegenerationass}(3) is precisely this requirement). In our cases of interest, we will always have a one-to-one inclusion-reversing correspondence between cones of $\Sigma(X)$ and logarithmic strata of $X$. Further, the fact that $X$ is a simple log scheme implies that the topological realization $|\Sigma(X)|$ is a genuine cone complex, so we shall confuse it with $\Sigma(X)$ and write $\Sigma(X)$ for both from now on. Similarly, for a morphism $f:X \to Y$ of simple log schemes we shall confuse $f_{\trop}: \Sigma(X) \to \Sigma(Y)$ with the induced $f_{\trop}: |\Sigma(X)| \to |\Sigma(Y)|$.

\subsubsection{Dual intersection complex}
\label{tdaff}

Let $\bar{\mathfrak{X}} \overset{\bar{g}}{\to} \mathcal{S}$ be a toric degeneration of $n$-dimensional Calabi-Yaus in the sense of Definition \ref{toricdegeneration} and satisfying Assumption \ref{toricdegenerationass}. Let $\bar{D}_i, \ 1 \leq i \leq \bar{m}$ be the (toric) irreducible components of the central fibre $\bar{\mathfrak{X}}_0$. We endow $\bar{\mathfrak{X}}$ with the divisorial log structure given by $\bar{D}:= \bar{\mathfrak{X}}_0 = \bar{D}_1 + \dots + \bar{D}_{\bar{m}}$ and $\mathcal{S}$ with the divisorial log structure given by $0 \in \mathcal{S}$. Then $$\bar{g}: (\bar{\mathfrak{X}}, \bar{D}) \to (\mathcal{S}, 0)$$ is a log morphism, and we can apply the tropicalization functor $\Sigma$ of Construction \ref{tropicalize} to produce a map of rational cone complexes\footnote{Here we use Assumption \ref{toricdegenerationass}(3) to get a map of actual cone complexes and not the generalized cone complexes of Construction \ref{tropicalize}.}: $$\bar{g}_{\trop}: \Sigma \left(\bar{\mathfrak{X}}\right) \to \Sigma(\mathcal{S}) = \mathbb{R}_{\geq 0}.$$ We define the \emph{dual intersection complex} of $\bar{\mathfrak{X}} \to \mathcal{S}$ as 
\begin{equation} \label{dic}
\left(\bar B, \bar{\mathscr{P}}\right): = \bar{g}_{\trop}^{-1}(1)
\end{equation}
where the polyhedral structure $\bar{\mathscr{P}}$ comes from restricting the cones of $\Sigma \left(\bar{\mathfrak{X}}\right)$ to the fibre over $1 \in \mathbb{R}_{\geq 0}$. Clearly, we have $\mathbf{C}\bar{B} \cong \Sigma\left(\bar{\mathfrak{X}}\right)$ where $\mathbf{C}\bar{B}$ is the cone over $\bar{B}$ of Section \ref{thetafunctions}. 

\begin{fact} \label{dicpm}
The dual intersection complex $\left(\bar B, \bar{\mathscr{P}}\right)$ of $\bar{\mathfrak{X}} \to \mathcal{S}$ is a polyhedral manifold of dimension $n$ in the sense of Definition \ref{polymani}. 
\end{fact}

\begin{proof}
By \cite[Propositions 4.10]{GS2}, $\bar{B}$ is a real topological manifold with $\partial \bar{B} = 0$ and of pure dimension $n$. Therefore, $\left(\bar B, \bar{\mathscr{P}}\right)$ satisfies condition (1) of Definition \ref{polymani}. The fact that $\bar{\mathscr{P}}$ is a set of bounded integral polyhedra follows as in \cite[Lemma 4.9]{GS2} so we need to check conditions (2)(a)-(2)(d) of Definition \ref{polymani}. 

Condition (2)(a) of Definition \ref{polymani} is satisfied since the logarithmic strata of $\bar{\mathfrak{X}}$ agree with the toric strata. Condition (2)(b) is immediate from \eqref{dic} and Construction \ref{tropicalize} of the tropicalization functor. The analogue of condition (2)(c) for $\Sigma\left(\bar{\mathfrak{X}}\right)$ is exactly the requirement that $\bar{\mathfrak{X}}$ is a simple log scheme (satisfied by Assumption \ref{toricdegenerationass}(3)). So condition (2)(c) follows by \eqref{dic}. Finally, condition (2)(d) is a direct consequence of Assumption \ref{toricdegenerationass}(4) via tropicalization.
\end{proof}

\begin{notation} \label{otonot}
    There is an inclusion-reversing correspondence between the toric strata of $\bar{\mathfrak{X}}_0$ and the cells of $\bar{\mathscr{P}}$. We denote the stratum corresponding to a cell $\sigma \in \bar{\mathscr{P}}$ by $\bar{X}_{\sigma}$. We make an exception for vertices $v \in \bar{\mathscr{P}}^{[0]}$ where we denote the corresponding divisor by $\bar{D}_v$. 
\end{notation}

By tracing the definitions, it is easy to check that our description of $\left(\bar B, \bar{\mathscr{P}}\right)$ is equivalent to the construction of \cite[Section 4.1]{GS2}.

\begin{remark}
Under some simplifying assumptions\footnote{One needs to assume that the $f_{\bar x}$ of condition (4)(b) of Definition \ref{toricdegeneration} vanish precisely once along each toric divisor of $Y_{\bar x}$. \cite[Section 7]{G1} also assumes that the irreducible components of $\bar{\mathfrak{X}}_0$ are normal but that follows from the fact that $(\bar{\mathfrak{X}}, \bar{D})$ is a simple log scheme (see Assumption \ref{toricdegenerationass}(3)).}, it is possible to describe $\left(\bar B, \bar{\mathscr{P}}\right)$ without using logarithmic geometry. Condition (4)(b) of Definition \ref{toricdegeneration} of a toric degeneration implies that for every $0$-dimensional stratum $\bar x$ of $\bar{\mathfrak{X}}_0$ there exists an affine toric variety $Y_{\bar x}$ such that $\bar{\mathfrak{X}} \overset{\bar g}{\to} \mathcal{S}$ is étale locally isomorphic to $Y_{\bar x} \overset{f_{\bar x}}{\to} \Spec \kk\llbracket t \rrbracket$ near $\bar{x}$. Under the assumptions, the $Y_x$ are defined by cones over certain polyhedra $\sigma_x$. One can then glue together the $\sigma_x$ by integral affine transformations to obtain $\left(\bar B, \bar{\mathscr{P}}\right)$. We refer to \cite[Section 7]{G1} for details.
\end{remark}

\subsubsection{Affine structure} \label{tdaffstr}

We want to use the formalism of Section \ref{scatteringsetup} to construct the mirror family $\check{\bar{\mathfrak{X}}} \to \Spec \kk\llbracket t \rrbracket$ from $\left(\bar B, \bar{\mathscr{P}}\right)$. By Proposition \ref{dicpm}, $\left(\bar B, \bar{\mathscr{P}}\right)$ is a polyhedral manifold, and we need to give it the structure of an integral affine manifold with singularities.

We first choose the discriminant locus $\bar{\Delta}$. In \cite{GS2}, $\bar{\Delta}$ is the union of all the cells of the barycentric subdivision $\tilde{\bar{\mathscr{P}}}$ of $\bar{\mathscr{P}}$ which are not contained in any cells of $\bar{\mathscr{P}}^{\max}$ and do not contain any vertices $v \in \bar{\mathscr{P}}^{[0]}$. However, to construct the mirror family from $\left(\bar B, \bar{\mathscr{P}}\right)$ via a scattering diagram (following \cite{GS1}), one needs to choose a more general $\bar{\Delta}$. 

\begin{cons} \label{toricdelta}
Following \cite[Section 1.1]{GS1}, let $x_{\sigma} \in \Int(\sigma)$ for $\sigma \in \bar{\mathscr{P}}^{[i]}, \ 1 \leq i \leq n-1$ be any points and let $$\bar \Delta \left(\left\{ x_{\sigma} \right\}\right) := \bigcup_{\sigma_1 \subseteq \dots \subseteq \sigma_{n-1}} \Conv \left\{x_{\sigma_i} \ | \ 1 \leq i \leq n-1\right\}$$ where the union is over all chains of cells $\sigma_1 \subseteq \dots \subseteq \sigma_{n-1}$ with $\sigma_i \in \bar{\mathscr{P}}^{[i]}$. Then \cite[Lemma 1.3]{GS3} shows that for any sufficiently general choice of $x_{\sigma}$ for $\sigma \in \bar{\mathscr{P}}^{[i]}, \ 1 \leq i \leq n-1$, the discriminant locus $\Delta \left(\left\{ x_{\sigma} \right\}\right)$ contains no rational point. We let $\bar{\Delta} : = \bar \Delta \left(\left\{ x_{\sigma} \right\}\right)$ for such a choice. As in Remark \ref{otherdelta}(3), we let $\tilde{\bar{\mathscr{P}}}$ be the polyhedral subdecomposition of $\bar{\mathscr{P}}$ induced by $\bar \Delta$.
\end{cons}

The point of the deformation of Construction \ref{toricdelta} is that since the walls of any scattering diagram $\bar{\mathfrak{D}}_I$ on $\left(\bar B, \bar{\mathscr{P}}\right)$ are rationally defined, any joint $\mathfrak{j} \in \bar{\mathfrak{D}}_I$ of codimension $1$ either satisfies $\mathfrak{j} \not \subseteq \bar \Delta$ or is only contained in slabs. Note also that by construction, there are no codimension $2$ joints $\mathfrak{j}$ with $\mathfrak{j} \subseteq \bar \Delta$.

\begin{remark} \label{indofsinglocus}
By \cite[Remark 5.3]{GS1}, the mirror family constructed from $\left(\bar B, \bar{\mathscr{P}}\right)$ does not depend on the choice of $x_{\sigma}$ for $\sigma \in \bar{\mathscr{P}}^{[i]}, \ 1 \leq i \leq n-1$ in Construction \ref{toricdelta}. By the discussion of \cite[Appendix A.1]{GHS}, it is isomorphic to the family constructed using the barycentric choice for $\bar{\Delta}$ of \cite{GS2}.
\end{remark}

 We now give $\left(\bar B, \bar{\mathscr{P}}\right)$ an affine structure on $\bar{B}_0 :=\bar{B} \setminus \bar \Delta$. 

\begin{cons} \label{affstrtoric}
Following the general framework of Construction \ref{affchartstool}, note that all the $\sigma \in \bar{\mathscr{P}}^{\max}$ have natural structures of integral polyhedra via their inclusion into the cones of $\Sigma \left(\bar{\mathfrak{X}}\right)$. In other words, these are precisely the integral polyhedra given by the polyhedral manifold structure on $\left(\bar B, \bar{\mathscr{P}}\right)$. 

For any $v \in \bar{\mathscr{P}}^{[0]}$, let $\Sigma_v$ be the fan defining the toric variety $\bar{D}_v$. It is easy to show (see the proof of \cite[Proposition 4.10]{GS2}) that a neighbourhood of $v$ in the dual intersection complex $\left(\bar B, \bar{\mathscr{P}}\right)$ is homeomorphic to $\Sigma_v$. 
Let
$$W_v := \bigcup_{\left\{\tilde{\sigma} \in \tilde{\bar{\mathscr{P}}} \ | \ v \in \tilde{\sigma}\right\}} \Int \tilde{\sigma}$$ as in \eqref{openstardef}. Then there is a unique integral affine linear map $\psi_v:W_v \to \mathbb{R}^n$ compatible with the induced structures of integral polyhedra on all the $\tilde{\sigma} \in \tilde{\bar{\mathscr{P}}}^{\max}$ with $v \in \tilde{\sigma}$ that maps every $\tilde{\sigma} \in \tilde{\bar{\mathscr{P}}}^{\max}$ to the corresponding cone of $\Sigma_v$. 

As in Construction \ref{affchartstool}, this data gives an integral affine structure on $\bar{B}_0 =\bar{B} \setminus \bar \Delta$. The affine structure of Construction \ref{affchartstool} extends to the complement of $\bar \Delta$ since $\psi_v: W_v \to \mathbb{R}^n$ defines an affine manifold structure on the whole $W_v$ and not just the complement of codimension $2$ cells of $\tilde{\bar{\mathscr{P}}}$ containing $v$.

Note that, as in the discussion after Construction \ref{affchartstool}, this affine structure admits the refined description of Construction \ref{affchartstool}$(2')$. We shall freely use this description when we need to.
\end{cons}

\subsubsection{MPA function and the initial slab functions} \label{tdmpafunc}

We work over $A=\kk$, use the toric monoid $\mathbb{N}$ and let $I_0 = \mathfrak{m}:= \mathbb{N} \setminus \left\{0\right\}$ be the maximal ideal (which corresponds to $(t) \subseteq  \kk[t] = \kk[\mathbb{N}]$). We assume that $\bar{\mathfrak{X}}$ carries a $\bar g$-ample divisor $A$ and fix such a choice. We say that $A$ is the \emph{polarization} of $\bar{g}: \bar{\mathfrak{X}} \to \mathcal{S}$. The polarization $A$ gives rise to an $\mathbb{N}_{\mathbb{R}}^{\gp} = \mathbb{R}$-valued MPA function $\bar{\varphi}_A$ on $\left(\bar B, \bar{\mathscr{P}}\right)$ defined via its kinks by setting $$\bar{\kappa}_{\underline{\rho}} := \bar{X}_{\rho} \cdot A \in  \mathbb{Z} = \mathbb{N}^{\gp}$$ for every $\underline{\rho} \subseteq \tilde{\bar{\mathscr{P}}}^{[n-1]}$ such that $\underline{\rho} \subseteq \rho \in \bar{\mathscr{P}}^{[n-1]}$. Note that for any two $\underline{\rho}, \underline{\rho}' \subseteq \rho \in \bar{\mathscr{P}}^{[n-1]}$ we have $\bar{\kappa}_{\underline{\rho}} = \bar{\kappa}_{\underline{\rho}'} = : \bar{\kappa}_{\rho}$ as we required in Section \ref{mpafunc}. Since $A$ is $\bar{g}$-ample, we have $\bar{\kappa}_{\rho} \in \mathbb{N} \setminus \left\{0\right\}$ and $\bar{\varphi}_A$ is strictly convex.

Now, we specify any initial slab functions 
\begin{equation} \label{initialslabs}
\left\{f_{\underline{\rho}} \in \kk[\Lambda_{\rho}] \ \Big| \ \underline{\rho} \in \tilde{\bar{\mathscr{P}}}^{[n-1]} \right\}
\end{equation}
satisfying \eqref{comp}\footnote{This is the same equation as the one in \cite[Theorem 3.27]{GS2} for trivial gluing data by identifying $f_{\underline{\rho}}$ with $f_{\underline{\rho}, x}$ for any $x \in \Int \underline{\rho}$, following \cite[Appendix A.1]{GHS}.}, having no poles, and \emph{normalized} in the sense that their constant coefficient is $1$. Such a choice is always possible in dimensions $2$ and $3$ (this is claimed in \cite{GS2} after Example 4.28), but it is not known if it exists in general. This information determines a \emph{toric log CY} structure on $\check{\bar{X}}_0$ in the sense of \cite[Definition 4.3]{GS2}. If $\left(\bar B, \bar{\mathscr{P}}\right)$ is \emph{simple} in the sense of \cite[Definition 1.60]{GS2}\footnote{Not to confuse with the notion of simplicity of Definition \ref{simplelogscheme}.}, then there is a unique choice of slab functions, see Proposition \ref{prop65} or \cite[Lemma A.1.1]{GHS}.

\begin{remark} \label{tdcentralfibre}
An important element in \cite{GS2, GS3} is a choice of open gluing data $s$. Also, one can use the \emph{intersection complex}\footnote{Defined using the $g$-ample divisor $A$.} $(\check{\bar{B}}, \check{\bar{\mathscr{P}}})$ of $\bar{\mathfrak{X}}$ (which is the \emph{discrete Legendre transform} of $\left(\bar B, \bar{\mathscr{P}}\right)$ as defined in \cite[Section 1.4]{GS2}) to construct $\check{\bar{X}}_0$ instead of using $\left(\bar B, \bar{\mathscr{P}}\right)$. In the language of \cite{GS2}, we have $\check{\bar{X}}_0 = X_0(\check{\bar{B}}, \check{\bar{\mathscr{P}}}, 1) = \check{X}_0(\bar{B}, \bar{\mathscr{P}}, 1)$ where $1$ stands for trivial gluing data. We shall sometimes need to discuss gluing data, notably in Assumption \ref{techass}, but we do not use non-trivial gluing data for constructing mirrors until Section \ref{glue}. We will give an overview of gluing data (in the language of \cite{GHS}) in Section \ref{gluingdataintro}.
\end{remark}

\subsubsection{Toric degenerations of K3-s}

We now specialize to toric degenerations of K3-s. 

\begin{fact} \label{spheredic}
Suppose that $\bar{\mathfrak{X}} \to \mathcal{S}$ is a toric degeneration of K3-s. Then topologically, $\bar B$ is a sphere.
\end{fact}

\begin{sloppypar}
\begin{proof}
By Proposition \ref{dicpm}, $\bar B$ is a two-dimensional real manifold. So it is enough to show that it has genus zero, i.e. to check that $\dim H^1(\bar{B}, \kk) = 0$. By \cite[Proposition 4.6]{GS2} and \cite[Theorem 4.14]{GS2}, the central fibre $\bar{\mathfrak{X}}_0$ of $\bar{\mathfrak{X}} \to \mathcal{S}$ is of the form $\bar{\mathfrak{X}}_0 \cong X_0(\bar{B}, \bar{\mathscr{P}}, s)$ for some open gluing data $s$ (in the notations of Remark \ref{tdcentralfibre}). But now \cite[Proposition 2.37]{GS2} implies that $H^i(\bar{B}, \kk) \cong H^i(\bar{\mathfrak{X}}_0, \mathcal{O}_{\bar{\mathfrak{X}}_0})$ for all $i \geq 0$. So we have $H^0(\bar{\mathfrak{X}}_0, \mathcal{O}_{\bar{\mathfrak{X}}_0}) \cong H^0(\bar{B}, \kk) \cong \kk, \ H^2(\bar{\mathfrak{X}}_0, \mathcal{O}_{\bar{\mathfrak{X}}_0}) \cong H^2(\bar{B}, \kk) \cong \kk$, and $H^1(\bar{B}, \kk) \cong H^1(\bar{\mathfrak{X}}_0, \mathcal{O}_{\bar{\mathfrak{X}}_0})$. Since a toric degeneration $\bar{\mathfrak{X}} \to \mathcal{S}$ is proper and flat, the arithmetic genus is constant in the fibres of $\bar{\mathfrak{X}} \to \mathcal{S}$. Let $\bar{\mathfrak{X}}_{\eta}$ be the generic fibre of $\bar{\mathfrak{X}} \to \mathcal{S}$ and note that we have $\dim H^0(\bar{\mathfrak{X}}_0, \mathcal{O}_{\bar{\mathfrak{X}}_0}) = \dim H^0(\bar{\mathfrak{X}}_{\eta}, \mathcal{O}_{\bar{\mathfrak{X}}_{\eta}}) = 1$ and $ \dim H^2(\bar{\mathfrak{X}}_0, \mathcal{O}_{\bar{\mathfrak{X}}_0}) = \dim H^2(\bar{\mathfrak{X}}_{\eta}, \mathcal{O}_{\bar{\mathfrak{X}}_{\eta}}) = 1$. Therefore, we also have $ \dim H^1(\bar{\mathfrak{X}}_0, \mathcal{O}_{\bar{\mathfrak{X}}_0}) = \dim H^1(\bar{\mathfrak{X}}_{\eta}, \mathcal{O}_{\bar{\mathfrak{X}}_{\eta}})$. But $\dim H^1(\bar{\mathfrak{X}}_{\eta}, \mathcal{O}_{\bar{\mathfrak{X}}_{\eta}}) = 0$ since $\bar{\mathfrak{X}}_{\eta}$ is a K3-surface. So $\dim H^1(\bar{B}, \kk) = 0$ and $\bar{B}$ is a topological sphere.
\end{proof}
\end{sloppypar}

\begin{remark}
$\bar{B}$ is not always a topological sphere in higher dimensions. For example, products of spheres occur. However, $\bar{B}$ is a sphere if $\bar{\mathfrak{X}} \to \mathcal{S}$ is a Batyrev degeneration of hypersurfaces or a Batyrev-Borisov degeneration of complete intersections satisfying a certain natural assumption. We refer to \cite[Remark 2.17]{G2} and the references therein for details.
\end{remark}

The singular locus $\bar{\Delta}$ of $\bar{B}$ is a union of irrational points $$\bar \Delta = \left\{x_{\rho} \in \Int(\rho) \ | \ \rho \in \bar{\mathscr{P}}^{[1]}\right\}$$ and it is easy to see that the local affine monodromy around $x_{\rho}$ is conjugate to 
$\begin{pmatrix}
1 & r_{\rho} \\
0 & 1
\end{pmatrix}$ for some $r_{\rho} \in \mathbb{N}$. We will call such an $x_{\rho}$ an \emph{$r_{\rho}$-fold singularity} and call $r_{\rho}$ the \emph{index} of $x_{\rho}$. Simplicity of $\left(\bar{B}, \bar{\mathscr{P}} \right)$ (in the sense of \cite[Definition 1.60]{GS2}) corresponds to having $r_{\rho} \leq  1$ for all $\rho \in \bar{\mathscr{P}}^{[1]}$.
	
\begin{remark}
One can have $r_{\rho}=0$, in which case $\bar{\mathfrak{X}} \to \mathcal{S}$ is log smooth in a neighbourhood of $\bar{X}_{\rho}$. By \cite[Proposition 1.27]{GS2} the affine structure on $\bar{B}_0$ extends across $x_{\rho}$ in this case.
\end{remark}

The initial slab functions in dimension $2$ can be described as follows.
	
\begin{fact} \label{struct}
Let $\underline{\rho}, \underline{\rho}' \in \tilde{\bar{\mathscr{P}}}^{[1]}$ be two slabs with $ \underline{\rho}, \underline{\rho}' \subseteq \rho \in \bar{\mathscr{P}}^{[1]}$ and let $f_{\underline{\rho}}, f_{\underline{\rho}'}$ be the corresponding slab functions. Let $w_{\rho} := z^{m_{\rho}}$ where $m_{\rho}$ is the integral generator of $\Lambda_{\rho}$ that points towards the vertex endpoint of $\underline{\rho}'$. Then we have
\begin{equation} \label{slabfunctions}
\begin{aligned}f_{\underline{\rho}} = &1 + a_{\rho, 1}w_{\rho} + \dots + a_{\rho, r_{\rho}-1}w_{\rho}^{r_{\rho}-1} + w_{\rho}^{r_{\rho}}& \\ f_{\underline{\rho}'} = &1 + a_{\rho, r_{\rho}-1}w_{\rho}^{-1} + \dots + a_{\rho,1}w_{\rho}^{-r_{\rho}+1} + w_{\rho}^{-r_{\rho}}&
\end{aligned}
\end{equation}
for some fixed choice of constants $a_{\rho, i} \in \kk, \ 1 \leq i \leq r_{\rho} -1$.
\end{fact}

\begin{proof}
This follows immediately from the fact that the monodromy vector $m_{\underline{\rho}'\underline{\rho}}$ of \eqref{monodromy} is $m_{\underline{\rho}'\underline{\rho}} = -r_{\rho} m_{\rho}$, the compatibility of slab functions \eqref{comp}, and the normalization requirement.
\end{proof}

Note that there is a unique choice of slab functions in the case that $\left(\bar{B}, \bar{\mathscr{P}} \right)$ is simple.

\begin{ex} \label{exaffine}
Continuing with the setup of Example \ref{runningexample}, the dual intersection complex $\left( \bar{B}, \bar{\mathscr{P}}\right)$ of $\bar{\mathfrak{X}} \to \Spec \kk \llbracket t \rrbracket$ is the boundary of a tetrahedron with each face isomorphic to a standard triangle. The affine structure near each vertex makes the polyhedral decomposition look like the fan for $\mathbb{P}^2$. There is a singularity $x_{\rho}$ with index $r_{\rho} = 4$ at an irrational point of each edge $\rho \in \bar{\mathscr{P}}^{[1]}$ that subdivides $\rho$ into two slabs. One can choose any initial slab functions satisfying \eqref{slabfunctions}.
\end{ex}

\subsubsection{The scattering diagram \texorpdfstring{$\bar{\mathfrak{D}}$}{D}, the mirror family \texorpdfstring{$\check{\bar{\mathfrak{X}}}_{\bar{\mathfrak{D}}}$}{XD}, and uniqueness} 

Fixing a choice of the initial slab functions $$\left\{f_{\underline{\rho}} \in \kk[\Lambda_{\rho}] \ \Big| \ \underline{\rho} \in \tilde{\bar{\mathscr{P}}}^{[n-1]} \right\}$$ of \eqref{initialslabs}, we define the \emph{initial scattering diagram} $\bar{\mathfrak{D}}_0:=\bar{\mathfrak{D}}_{I_0} = \bar{\mathfrak{D}}_{(t)}$, whose only walls are the slabs $(\underline{\rho}, f_{\underline{\rho}})$ with support $\underline{\rho} \in \tilde{\bar{\mathscr{P}}}^{[n-1]}$ and the attached initial slab function $f_{\underline{\rho}}$. Consistency of $\bar{\mathfrak{D}}_0$ follows from equation \eqref{comp} relating $f_{\underline{\rho}}$ and $f_{\underline{\rho}'}$ for two slabs $\underline{\rho}, \underline{\rho}' \subseteq \rho \in \bar{\mathscr{P}}^{[n-1]}$. The main result of \cite{GS1}, proved in \cite[Proposition 3.9]{GS1}, is the following.

\begin{theorem} \label{reconstruction}
Suppose that $\check{\bar{X}}_0$ is locally rigid in the sense of \cite[Definition 1.26]{GS1}. Then we have the following.

\begin{enumerate}
    \begin{sloppypar}\item[]\textbf{Existence:} There exists a collection of scattering diagrams $\bar{\mathfrak{D}} = \left\{ \bar{\mathfrak{D}}_k, \ k \geq 0 \right\}$ where $$\bar{\mathfrak{D}}_k:=\bar{\mathfrak{D}}_{I^{k+1}_0} = \bar{\mathfrak{D}}_{(t^{k+1})},$$ such that $\bar{\mathfrak{D}}_k$ is compatible with $\bar{\mathfrak{D}}_{k-1}$ (in the sense of Definition \ref{compatibility}) for $k \geq 1$ and the $\bar{\mathfrak{D}}_k, \ k \geq 0$ are consistent in the sense of \cite[Definition 2.28]{GS1}. We sometimes refer to the whole collection $\bar{\mathfrak{D}}$ as an algorithmic scattering diagram.\end{sloppypar}
    \item[]\textbf{Uniqueness:} For any two sequences $\bar{\mathfrak{D}}, \ \bar{\mathfrak{D}}'$ of compatible scattering diagrams such that $\bar{\mathfrak{D}}_k, \ \bar{\mathfrak{D}}_k', \ k \geq 0$ are consistent in the sense of \cite[Definition 2.28]{GS1} and $\bar{\mathfrak{D}}_0$ is combinatorially equivalent to $\bar{\mathfrak{D}}'_0$, $\bar{\mathfrak{D}}_k$ is combinatorially equivalent to $\bar{\mathfrak{D}}'_k$ for any $k \geq 0$.
\end{enumerate}
\end{theorem}

The theorem is a far-reaching generalization of \cite[Theorem 6]{KS} and provides an explicit algorithm to produce $\bar{\mathfrak{D}}$. The most important property for us is uniqueness up to combinatorial equivalence. Indeed, combinatorially equivalent scattering diagrams produce isomorphic canonical families by Observation \ref{combimpliesequiv}(2) and Proposition \ref{help}.
\emph{Local rigidity} is a technical assumption that guarantees uniqueness up to equivalence. It is empty in dimension $2$. More generally, simplicity of $\left(\bar{B}, \bar{\mathscr{P}}\right)$ implies local rigidity of $\check{\bar{X}}_0$, see \cite[Remark 1.29]{GS1}.

As in Section \ref{thetafunctions}, we obtain families $$\check{\bar{\mathfrak{X}}}_{\bar{\mathfrak{D}}_k} \to \Spec \kk[t]/(t^{k+1})$$ for every $k \geq 0$ that form an inverse system since $\bar{\mathfrak{D}}_k$ is compatible with $\bar{\mathfrak{D}}_{k-1}$ for $k \geq 1$. Taking the limit over this system gives (see \cite[Proposition 2.42]{GS2} and \cite[Remark A.1.4]{GHS}) a toric degeneration
\begin{equation} \label{reconstructedfamily}
\check{\bar{\mathfrak{X}}}:=\check{\bar{\mathfrak{X}}}_{\bar{\mathfrak{D}}} \to \Spec \kk\llbracket t \rrbracket  
\end{equation}
that we call the \emph{toric degeneration mirror} to $\bar{\mathfrak{X}} \to \mathcal{S}$.

We shall not use consistency in the sense of \cite[Definition 2.28]{GS1} directly. Consistency in the sense of \cite[Definition 2.28]{GS1} implies consistency in the sense of Definition \ref{consistency} by \cite[Lemma A.1.2]{GHS} and the reference to \cite[Proposition 3.2]{CPS} therein. Conversely, we have the following.

\begin{fact} \label{consistencydiff}
Let $\bar{\mathfrak{D}}'_k$ be a scattering diagram on $\left(\bar B, \bar{\mathscr{P}}\right)$. It is consistent in the sense of \cite[Definition 2.28]{GS1} if the following two conditions are satisfied.
\begin{enumerate}
    \item $\bar{\mathfrak{D}}'_k$ is consistent in the sense of Definition \ref{consistency}.
    \item $\bar{\mathfrak{D}}'_k$ is consistent in the sense of \cite[Definition 2.28]{GS1} around every codimension one joint $\mathfrak{j} \subseteq \bar \Delta$.
\end{enumerate} 
\end{fact}

\begin{proof}
Suppose that $\bar{\mathfrak{D}}'_k$  satisfies (1) and (2). The notion of consistency in codimension $0$ is the same, so it is enough to prove consistency in the sense of \cite[Definition 2.28]{GS1} around all codimension $1$ and $2$ joints. 

Suppose that $\mathfrak{j} \in \bar{\mathscr{P}}_{\bar{\mathfrak{D}}'_k}$ is a codimension $1$ joint that is an intersection of two slabs $\mathfrak{j} = \mathfrak{b}_1 \cap \mathfrak{b}_2$. By (2), we may assume that $\mathfrak{j} \not \subseteq \bar \Delta$, which allows us to reinterpret consistency in codimension $1$ in the sense of Definition \ref{codim1consistency} in a way similar to Definition \ref{codim0consistency} for consistency in codimension $0$. We consider the ring $R'_{\mathfrak{b}_1}:= \kk[\mathcal{P}_x^+]/\mathcal{I}_x$ for $x \in \Int \mathfrak{b}_1$. Suppose that $\mathfrak{b}_1, \mathfrak{p}_1, \dots, \mathfrak{p}_r, \mathfrak{b}_2, \mathfrak{p}'_1, \dots, \mathfrak{p}'_s$ is a cyclic ordering of walls around $\mathfrak{j}$ as in Definition \ref{codim1consistency}. Since $\mathfrak{j} \not \subseteq \bar \Delta$, the definition of walls (Definition \ref{wall}) implies that we may canonically view $$f_{\mathfrak{b}_1}, f_{\mathfrak{p}_1}, \dots, f_{\mathfrak{p}_r}, f_{\mathfrak{b}_2}, f_{\mathfrak{p}'_1}, \dots, f_{\mathfrak{p}'_s} \in R'_{\mathfrak{b}_1}$$ as in \cite[Section 2.2.2]{AG}. Then define
$$\theta_{\mathfrak{p}}: R'_{\mathfrak{b}} \to R'_{\mathfrak{b}}, \quad z^m \mapsto f_{\mathfrak{p}}^{\langle n_{\mathfrak{p}}, \bar{m} \rangle} z^m$$ for any wall containing $\mathfrak{j}$ as in \eqref{wallcross} and define $\theta_{\gamma_{\mathfrak{j}}}: R'_{\mathfrak{b}} \to R'_{\mathfrak{b}}$ as the composition
$$\theta_{\gamma_{\mathfrak{j}}}: = \theta_{\mathfrak{b}_1} \circ \theta_{\mathfrak{p}_1} \circ \dots  \circ \theta_{\mathfrak{p}_r} \circ \theta_{\mathfrak{b}_2} \circ \theta_{\mathfrak{p}'_1} \circ \dots  \circ \theta_{\mathfrak{p}'_s}.$$
As in \cite[Remark 2.18]{AG}, consistency around $\mathfrak{j}$ is equivalent to having $\theta_{\gamma_{\mathfrak{j}}} = \Id$ as an automorphism of $R'_{\mathfrak{b}}$. But having $\theta_{\gamma_{\mathfrak{j}}} = \Id$ is also equivalent to consistency around $\mathfrak{j}$ in the sense of \cite[Definition 2.28]{GS1}.

Suppose $\mathfrak{j}$ is a codimension $2$ joint. By our choice of $\bar \Delta$ we have $\mathfrak{j} \not \subseteq \bar \Delta$. Then if $\mathfrak{p}$ is a wall containing $\mathfrak{j}$, $f_{\mathfrak{p}}$ is defined as an element of $\kk[\mathcal{P}_x^+]/\mathcal{I}_x$ for $x \in \mathfrak{j}$. By \eqref{stalkint} and \eqref{idealint}, we may actually view $f_{\mathfrak{p}} \in R'_{\mathfrak{b}}$ as above for a slab $\mathfrak{b}$ with $\mathfrak{j} \subseteq \mathfrak{b}$. We may then cyclically order the walls and define $\theta_{\gamma_{\mathfrak{j}}}$ as in the case of a codimension $1$ joint. As in the case of codimension $1$ joint, consistency around $\mathfrak{j}$ in the sense of \cite[Definition 2.28]{GS1} is equivalent to having $\theta_{\gamma_{\mathfrak{j}}} = \Id$. 

So it remains to show that $\theta_{\gamma_{\mathfrak{j}}} = \Id$ for a codimension $2$ joint $\mathfrak{j}$. We proceed as in the proof of \cite[Theorem 2.41]{AG}. Since $\mathfrak{j}$ is codimension $2$ we can pass to
$\left(\bar{B}_{\mathfrak{j}}, \bar{\mathscr{P}}_{\mathfrak{j}}\right)$ and $\bar{\mathfrak{D}}'_{k, \mathfrak{j}}$. Choose a general point $p \in \bar{\mathfrak{D}}'_{k, \mathfrak{j}} \setminus \mathfrak{j}$. Then by the construction of broken lines of Definition \ref{brokenlines} and the formula \eqref{thetafunc} for $\vartheta^{\mathfrak{j}}_{m}(p)$, we may view $\vartheta^{\mathfrak{j}}_{m}(p) \in R'_{\mathfrak{b}} = \kk[\mathcal{P}_x^+]/\mathcal{I}_x$ for $x \in \Int \mathfrak{b}$. Now, consistency along $\mathfrak{j}$ implies, by composing the wall-crossings homomorphisms around $\mathfrak{j}$, that $\theta_{\gamma_{\mathfrak{j}}}(\vartheta^{\mathfrak{j}}_{m}(p)) = \vartheta^{\mathfrak{j}}_{m}(p)$ for any asymptotic monomial $m$. But \cite[Proposition 3.2.9]{GHS} implies that $\vartheta^{\mathfrak{j}}_{m}(p)$ generate $R'_{\mathfrak{b}}$ as a $\kk[t]/(t^{k+1})$ module so $\theta_{\gamma_{\mathfrak{j}}}$ acts as the identity and we are done.
\end{proof}

\subsubsection{Special toric degenerations} \label{specialtd}

We shall now define the property of being a special toric degeneration in Conjecture \ref{conjecture}.  We need to use open gluing data (see Definition \ref{gddef} and \cite[Definition 2.25]{GS2}) to state the definition in full generality.
Consider a toric degeneration $\bar{\mathfrak{X}} \to \mathcal{S}$ with polarization $A$ and dual intersection complex $\left(\bar B, \bar{\mathscr{P}}\right)$. 
\begin{ass} \label{techass}
The central fibre $\bar{\mathfrak{X}}_0$ of $\bar{\mathfrak{X}} \to \mathcal{S}$ satisfies $\bar{\mathfrak{X}}_0 \cong X_0(\bar{B}, \bar{\mathscr{P}}, s)$ (using the same notation as in Remark \ref{tdcentralfibre}) for some open gluing data $s$ on $\left(\bar B, \bar{\mathscr{P}}\right)$ such that $o(s) = 1 \in H^2(\bar{B}, \kk^{\times})$ for the homomorphism $o$ defined in \cite[Theorem 2.34]{GS2}.
\end{ass}

This assumption is closely related to projectivity of $\bar{\mathfrak{X}}_0$ (which follows from Assumption \ref{toricdegenerationass}(2)) but unless $H^1(\bar B, \mathbb{Q}) = 0$ it is not equivalent. 
Apart from defining special toric degenerations, we shall use this assumption when generalizing Conjecture \ref{conjecture} in Chapter \ref{gslocus}. In the language of \cite{GS2} (see Remark \ref{tdcentralfibre}), Assumption \ref{techass} means that one can use the discrete Legendre transform $(\check{\bar{B}}, \check{\bar{\mathscr{P}}}, \check{\bar{\varphi}}_A)$ of $(\bar{B}, \bar{\mathscr{P}}, \bar{\varphi}_A)$ to construct $\bar{\mathfrak{X}}_0$ as $\bar{\mathfrak{X}}_0 \cong \check{X}_0(\check{\bar{B}}, \check{\bar{\mathscr{P}}}, \check{s})$ for some gluing data $\check{s}$ on $(\check{\bar{B}}, \check{\bar{\mathscr{P}}})$.

\begin{defn} \label{distinguished}
Suppose that the central fibre $\bar{\mathfrak{X}}_0$ of $\bar{\mathfrak{X}} \to \mathcal{S}$ satisfies Assumption \ref{techass} and is locally rigid in the sense of \cite[Definition 1.26]{GS2}. Then we say that $\bar{\mathfrak{X}} \to \mathcal{S}$ is \emph{distinguished} if it is isomorphic to the basechange of the family $\bar{\mathfrak{X}}_{\check{\bar{\mathfrak{D}}}} \to \Spec \kk\llbracket t \rrbracket$ of \eqref{reconstructedfamily} by the map $\kk\llbracket t \rrbracket \to R$ of Remark \ref{dvrremark}. Here $\check{\bar{\mathfrak{D}}}$ is a scattering diagram on $(\check{\bar{B}}, \check{\bar{\mathscr{P}}})$ constructed using Theorem \ref{reconstruction} and one uses the MPA function $\check{\bar{\varphi}}_A$ for the construction.
\end{defn}

Here Assumption \ref{techass} guarantees the existence of $\bar{\mathfrak{X}}_{\check{\bar{\mathfrak{D}}}} \to \Spec \kk\llbracket t \rrbracket$ as an algebraic (and not just formal) family.\footnote{We only stated Theorem \ref{reconstruction} for trivial gluing data. A similar statement is true for any choice of gluing data $s$ on $\left(\bar B, \bar{\mathscr{P}}\right)$, but one can only construct a formal family in general.} Indeed, the MPA-function $\check{\bar{\varphi}}_A$ gives rise to ample line bundles on the toric irreducible components of $\bar{\mathfrak{X}}_0$ and one can see that $o(s) \in H^2(\bar{B}, \kk^{\times})$ is the obstruction to gluing them to an ample line bundle on $\bar{\mathfrak{X}}_0$. Thus, Assumption \ref{techass} provides an ample line bundle $L$ on $\bar{\mathfrak{X}}_0$. One can show that $L$ extends to the formal degeneration and obtain an algebraic family via Grothendieck’s existence theorem. Alternatively, one can obtain the algebraic family directly by applying \cite[Theorem 5.2.19]{GHS}. See \cite[Remarks 5.2.15 and A.1.4]{GHS} for more details.

\begin{defn} \label{special}
A toric degeneration $\bar{\mathfrak{X}} \to \mathcal{S}$ is \emph{special} if it satisfies the following conditions:
\begin{enumerate}
    \item The generic fibre of $\bar{\mathfrak{X}} \to \mathcal{S}$ is smooth.
    \item There exists a toric log CY structure on $\check{\bar{X}}_0$\footnote{Recall from Section \ref{tdmpafunc} that such a structure always exists in dimensions $2$ and $3$ (this requires using trivial gluing data in dimension $3$) but it is not known if it exists in general.}  and it is locally rigid in the sense of \cite[Definition 1.26]{GS2}.
    \item One of the following conditions is satisfied:
    \begin{enumerate}
        \item $\bar{\mathfrak{X}} \to \mathcal{S}$ is a \emph{divisorial log deformation} of the central fibre $\bar{\mathfrak{X}}_0$ in the sense of \cite[Definition 2.7]{models}.\footnote{Strictly speaking, \cite[Definition 2.7]{models} does not apply in this situation since $R$ is not Artinian. However, the definition clearly extends to non-Artinian local $\kk$-algebras. Equivalently, we can say that $\bar{\mathfrak{X}} \to \Spec R$ is a divisorial log deformation if all the induced families $\bar{\mathfrak{X}} \to \Spec R/\mathfrak{m}^k$ for $k \geq 1$ are divisorial log deformations (here $\mathfrak{m}$ is the unique maximal ideal of $R$).}
        \item $\bar{\mathfrak{X}}_0$ satisfies Assumption \ref{techass} and is locally rigid. Moreover, $\bar{\mathfrak{X}} \to \mathcal{S}$ is distinguished in the sense of Definition \ref{distinguished}.
    \end{enumerate}
\end{enumerate}
\end{defn}

\begin{remarks} \label{specialass}
We explain the assumptions:
\begin{enumerate}
    \item We need this assumption to construct a log smooth resolution. Indeed, the log structure is trivial on the generic fibre, so log smoothness is equivalent to regular smoothness. Thus, we only need to resolve the singularities in the central fibre.
    \item This assumption is necessary to construct the toric degeneration mirror $\check{\bar{\mathfrak{X}}}_{\bar{\mathfrak{D}}} \to \Spec \kk\llbracket t \rrbracket$ using Theorem \ref{reconstruction}. Conjecture \ref{conjecture} would then provide a particular (possibly using a different choice of the initial slab functions) toric degeneration mirror.  
    \item These assumptions guarantee the existence of local models near points of the discriminant locus $Z$.
    \begin{enumerate}
        \item In this case, local models come from the local structure of the family (see \cite[Construction 2.1]{models}). Condition (1) further restricts the possible local models (see \cite[Proposition 2.2]{models}).
        \item The first two assumptions are necessary to construct $\bar{\mathfrak{X}}_{\check{\bar{\mathfrak{D}}}} \to \Spec \kk\llbracket t \rrbracket$ as an algebraic family. The local models come from the reconstruction algorithm of Theorem \ref{reconstruction} and can be \enquote{decomposed} into standard pieces that look like the local models in (a), see \cite[Sections 4.4.2 and 4.4.3]{GS2}. In fact, a distinguished toric degeneration is a divisorial log deformation if $(\check{\bar{B}}, \check{\bar{\mathscr{P}}})$ is simple, see \cite[Corollary 2.18]{models}.
    \end{enumerate}
\end{enumerate}
\end{remarks}

Although we have stated Conjecture \ref{conjecture} in full generality, in higher dimensions, we will mainly be interested in the case when $\bar{\mathfrak{X}} \to \mathcal{S}$ is a distinguished toric degeneration with $\left(\bar B, \bar{\mathscr{P}}\right)$ simple. We shall discuss this further in Chapter \ref{higherdimgen}.

\begin{fact} \label{distinguishedspecial}
A distinguished toric degeneration $\bar{\mathfrak{X}} \to \mathcal{S}$ with a simple dual intersection complex $\left(\bar B, \bar{\mathscr{P}}\right)$ satisfies conditions (2), (3)(a) and 3(b) of Definition \ref{special}. In particular, it is special if and only if it satisfies condition (1) of Definition \ref{special}.
\end{fact}

\begin{proof}
The simplicity of $\left(\bar B, \bar{\mathscr{P}}\right)$ implies that there is a unique way to put a toric log CY structure on $\check{\bar{X}}_0$ by \cite[Theorem 5.2]{GS2}. Then $\check{\bar{X}}_0$ is locally rigid by \cite[Remark 1.29]{GS1} since $\left(\bar B, \bar{\mathscr{P}}\right)$ is simple. So $\bar{\mathfrak{X}} \to \mathcal{S}$ satisfies condition (2) of Definition \ref{special}. $\bar{\mathfrak{X}} \to \mathcal{S}$ satisfies condition (3)(a) of Definition \ref{special} by \cite[Corollary 2.18]{models} and the fact that the basechange by the map $\kk\llbracket t \rrbracket \to R$ of Remark \ref{dvrremark} preserves the notion of being a divisorial log deformation. It satisfies condition (3)(b) since it is distinguished. Note that the discrete Legendre transform $(\check{\bar{B}}, \check{\bar{\mathscr{P}}})$ of $\left(\bar B, \bar{\mathscr{P}}\right)$ is also simple and so we don't need to check local rigidity of $\bar{\mathfrak{X}}_0$ in Definition \ref{distinguished}.
\end{proof}

We now specialize to the case of toric degenerations of K3-s. First, we would like to understand when a toric degeneration $\bar{\mathfrak{X}} \to \mathcal{S}$ of K3-s is a divisorial log deformation of the central fibre. \cite[Definition 2.7]{models} of being a divisorial log deformation provides étale local models at the points of the singular locus $Z$ of $\bar{\mathfrak{X}} \to \mathcal{S}$. Therefore, we need to investigate the local models for toric degenerations of K3-s in the neighbourhoods of one-dimensional strata.

Let $\rho \in \bar{\mathscr{P}}^{[1]}$ be an edge. The local model for the toric degeneration $\bar{\mathfrak{X}} \to \mathcal{S}$ in the neighbourhood of a codimension $1$ stratum $\bar{X}_{\rho}$ (one-dimensional for $\bar{\mathfrak{X}} \to \mathcal{S}$ a toric degeneration of K3-s) is given by (see \cite[(1.9)]{GS2} and \cite[Example 2.8]{models}): 
\begin{equation} \label{localpolynomials}
\left\{xy = t^lf_{\rho}(t,w_{\rho}) \right\} \subseteq \Spec \kk[x,y,w_{\rho}]\llbracket t \rrbracket
\end{equation}
(with the natural map to $\Spec \kk\llbracket t \rrbracket$) where $l \geq 1$ is the integral length of $\rho \in \bar{\mathscr{P}}^{[1]}$ and $f_{\rho}(t,w_{\rho}) \in \kk [w_{\rho}]\llbracket t \rrbracket$ is a polynomial in $w_{\rho}$ of degree $r_{\rho}$ (the index of the singularity $x_{\rho} \subseteq \rho$ of the affine structure) that is not divisible by $t$.  

\begin{lemma} \label{localmodelsk3}
Let $\bar{\mathfrak{X}} \to \mathcal{S}$ be a toric degeneration of K3-s and suppose that all the polynomials $f_{\rho}(t,w_{\rho}) \in  \kk[w_{\rho}]\llbracket t \rrbracket$ for $\rho \in \bar{\mathscr{P}}^{[1]}$ in \eqref{localpolynomials} are of the form 
\begin{equation} \label{localpolyexp}
f_{\rho}(t,w_{\rho}) = \prod_{i=1}^{r'_{\rho}} (w_{\rho}-\gamma^i_{\rho}(t))^{k_i} \alpha_{\rho}(t, w_{\rho})
\end{equation}
where $k_i \in \mathbb{Z}_{>0}$ are multiplicities of the roots,  $0 \leq \sum_{i=1}^{r'_{\rho}} k_i \leq r_{\rho}$, $\gamma^i_{\rho}(t) \in \kk \llbracket t \rrbracket$ are power series with $\gamma^i_{\rho}(0) \neq \gamma^j_{\rho} (0)$ for $1 \leq i < j \leq r_{\rho}'$, and $\alpha_{\rho}(t, w_{\rho}) \in \kk[w_{\rho}](\!( t )\!)$ is a polynomial in $w_{\rho}$ of degree $r_{\rho} - \sum_{i=1}^{r'_{\rho}} k_i$ with $\alpha_{\rho}(0, w_{\rho}) = \beta \in \kk \setminus \{0\}$. Then $\bar{\mathfrak{X}} \to \mathcal{S}$ is étale locally isomorphic to
\begin{equation} \label{localmodeldesired}
\left\{xy = t^lw_{\rho}^k \right\} \subseteq \Spec \kk[x,y,w_{\rho}]\llbracket t \rrbracket
\end{equation} (with the natural map to $\Spec \kk\llbracket t \rrbracket$) for some $ l \geq 1, \ k \geq 0$ in the neighbourhood of any point in a one-dimensional stratum of $\bar{\mathfrak{X}}_0$.
\end{lemma}

\begin{proof}
Fix a $\rho \in \bar{\mathscr{P}}^{[1]}$ and the corresponding stratum $\bar{X}_{\rho}$, and consider the local model 
\begin{equation} \label{localmodelproof}
\left\{xy = t^lf_{\rho}(t,w_{\rho}) \right\} \subseteq \Spec \kk[x,y,w_{\rho}]\llbracket t \rrbracket
\end{equation}
in the neighbourhood of $\bar{X}_{\rho}$. It is enough to check that $\bar{\mathfrak{X}} \to \mathcal{S}$ has the required local models at all the points of $\bar{X}_{\rho}$. 

The points of $\bar{X}_{\rho}$ correspond to setting $w_{\rho} = w^0_{\rho} \in \kk$ (and $x=y=0$). Fix an $1 \leq i \leq r'_{\rho}$ and let $w^0_{\rho} = \gamma_{\rho}^i(0)$. Consider a change of variables $$\kk[x,y,w_{\rho}]\llbracket t \rrbracket \to \kk[x,y,w_{\rho}]\llbracket t \rrbracket, \quad w_{\rho} \mapsto w_{\rho} - \gamma_{\rho}^i(t), \ (x, y, t) \mapsto (x,y,t).$$ This clearly determines an étale map under which the local model \eqref{localmodelproof} becomes
\begin{equation} \label{localmodelproof1}
\left\{xy = t^lw_{\rho}^{k_i} g_{\rho}(t,w_{\rho})\right\} \subseteq \Spec \kk[x,y,w_{\rho}]\llbracket t \rrbracket
\end{equation}
for a certain polynomial $g_{\rho}(t,w_{\rho}) \in \kk[w_{\rho}]\llbracket t \rrbracket$. Moreover, our assumptions ensure that $g_{\rho}(t, w_{\rho})$ is invertible in the neighbourhood of $0$. 

We need to check that \eqref{localmodeldesired} is étale locally isomorphic to \eqref{localmodelproof1} in the neighbourhood of $0$. Standard results about étale morphisms (see, e.g. \cite[Chapter I. 3.16]{étale}) say that a map of rings $C \to C[Z_1, \dots, Z_n]/(F_1, \dots, F_n)$ is étale if and only if $\det(\partial F_i/\partial Z_j)$ is a unit in $C[Z_1, \dots, Z_n]/(F_1, \dots, F_n)$. Consider the map $$\frac{\kk[x,y,w_{\rho}]\llbracket t \rrbracket[a]}{(xy - t^l w_{\rho}^{k_i}, g_{\rho}(t, w_{\rho})a -1)} \quad \longrightarrow \quad \frac{\kk[x,y,w_{\rho}]\llbracket t \rrbracket[a]}{(xy - t^l w_{\rho}^{k_i}, g_{\rho}(t, w_{\rho})a -1)}\left[ Z \right] \bigg/ \left(y-Za \right).$$ We have $\det(\partial F_i/\partial Z_j) = \partial (y-Za)/\partial Z = -a$ which is a unit. But the ring on the left-hand side is the localization of the coordinate ring for \eqref{localmodeldesired} (with $k:=k_i$) to $g(t, w_{\rho}) \neq 0$ and the ring on the right-hand side is the localization of the coordinate ring for \eqref{localmodelproof1} to $g_{\rho}(t, w_{\rho}) \neq 0$. Since $g(t, w_{\rho})$ is invertible in the neighbourhood of $0$, this implies that $\bar{\mathfrak{X}} \to \mathcal{S}$ has the required local model at the point corresponding to setting $w_{\rho} = \gamma_{\rho}^i(0)$. Note that this is compatible with the maps to $\Spec \kk \llbracket t \rrbracket$.

For every point of $\bar{X}_{\rho}$ corresponding to setting $w_{\rho} = w^0_{\rho}$ with $w^0_{\rho} \neq \gamma_{\rho}^i(0)$ for $1 \leq i \leq r'_{\rho}$, the same arguments as above show that $\bar{\mathfrak{X}} \to \mathcal{S}$ is étale locally isomorphic to 
$$\left\{xy = t^l \right\} \subseteq \Spec \kk[x,y,w_{\rho}]\llbracket t \rrbracket.$$
Therefore, $\bar{\mathfrak{X}} \to \mathcal{S}$ has the required local models at all the points of any one-dimensional stratum of $\bar{\mathfrak{X}}_0$.
\end{proof}

\begin{obs} \label{modelslocal}
Let $\bar{\mathfrak{X}} \to \mathcal{S}$ be a toric degeneration of K3-s satisfying the assumptions of Lemma \ref{localmodelsk3}.
\begin{enumerate}
    \item $\bar{\mathfrak{X}} \to \mathcal{S}$ is log smooth in the neighbourhood of any point with a local model of the form $\left\{xy = t^l \right\} \subseteq \Spec \kk[x,y,w_{\rho}]\llbracket t \rrbracket$ and the finite number of points with local models $\left\{xy = t^lw_{\rho}^k \right\} \subseteq \Spec \kk[x,y,w_{\rho}]\llbracket t \rrbracket$ for $k \geq 1$ form the singular locus $Z$ of $\bar{\mathfrak{X}} \to \mathcal{S}$.
    \item If the generic fibre of $\bar{\mathfrak{X}} \to \mathcal{S}$ is smooth, then we only have local models of the form $\left\{xy = t^l \right\} \subseteq \Spec \kk[x,y,w_{\rho}]\llbracket t \rrbracket$ and of the form $\left\{xy = t^lw_{\rho} \right\} \subseteq \Spec \kk[x,y,w_{\rho}]\llbracket t \rrbracket$ since for $k \geq 1$ the generic fibre of $\left\{xy = t^lw_{\rho}^k \right\} \subseteq \Spec \kk[x,y,w_{\rho}]\llbracket t \rrbracket$ has $A_{k-1}, \ k \geq 2$ singularities. This corresponds to $f_{\rho}(t, w_{\rho})$ not having any multiple roots $\gamma^i_{\rho}(t) \in \kk \llbracket t \rrbracket$ for all $\rho \in \bar{\mathscr{P}}^{[1]}$. 
\end{enumerate}
Note that (1) and (2) apply to any divisorial log deformation $\bar{\mathfrak{X}} \to \mathcal{S}$ of K3-s. This follows from \cite[Definition 2.7]{models} of being a divisorial log deformation, and the fact that all the étale local models in dimension $3$ given by \cite[Construction 2.1]{models} are of the form $\left\{xy = t^lw_{\rho}^k \right\} \subseteq \Spec \kk[x,y,w_{\rho}]\llbracket t \rrbracket$ for some $ l \geq 1, \ k \geq 0$ (see the proof of Proposition \ref{k3divisorial} below). 
\end{obs}

\begin{fact} \label{k3divisorial}
Let $\bar{\mathfrak{X}} \to \mathcal{S}$ be a toric degeneration of K3-s satisfying the assumptions of Lemma \ref{localmodelsk3}. Then $\bar{\mathfrak{X}} \to \mathcal{S}$ is a divisorial log deformation of the central fibre. 
\end{fact}

\begin{proof}
By \cite[Definition 2.7]{models} of being a divisorial log deformation, it is enough to check that there are étale local models given by \cite[Construction 2.1]{models} at all the points in the discriminant locus $Z$. Indeed, such a local model is an affine toric variety given by the cone over the convex hull 
\begin{equation} \label{convexhull}
\Conv\left( \bigcup_{i=0}^q \left(\Delta_i \times \left\{e_i\right\}\right) \right)
\end{equation}
where $\Delta_i$ are integral polytopes in a lattice $N$ (not necessarily of maximal dimension), $e_i$ for $ \ 0 \leq i \leq q$ are the standard generators of the second factor of $N' := N \oplus \mathbb{Z}^{q+1}$ and the convex hull is taken in $N'_{\mathbb{R}}$. The map to $\Spec \kk[t]$ is induced by $e_0^{*}$. Now, set $N := \mathbb{Z}, \ q := 1$, and let $\Delta_0$ be an interval of length $l$ and $\Delta_1$ be an interval of length $k$ (note that this is the general form of a local model in dimension $3$). The corresponding toric variety is a cone over the convex hull of $(0, 0, 1), (k, 0, 1), (0, 1, 0)$ and $(l, 1, 0)$, and it is easy to see that it is given by an equation of the form $$\{xy = t^l w^k\} \subseteq \Spec \kk[x,y,t,w].$$ Since the local models of \eqref{localmodeldesired} are given by the same equation (after completion in $t$), this implies the result.
\end{proof}

The assumptions of Lemma \ref{localmodelsk3} are rather general. Let $\bar{\kk}(\!( t )\!)$ be the field of Puiseux series. For a $\beta(t) \in \bar{\kk}(\!( t )\!)$ of the form $\beta(t) = \sum_{i=i_0}^{+\infty} a_i t^{\frac{i}{c}}$ for some $c \in \mathbb{Z}_{>0}$ and for $i_0 \in \mathbb{Z}$ the smallest integer such that the coefficient at $t^{\frac{i_0}{c}}$ is non-zero, we shall write $\beta(0) = a_0$ if $i_0 = 0$, $\beta(0) = 0$ if $i_0 > 0$, and $\beta(0) = \infty$ if $i_0 < 0$. We denote $\deg \beta(t):=\frac{i_0}{c}$.

Any $f_{\rho}(t, w_{\rho}) \in \kk[w_{\rho}]\llbracket t \rrbracket$ of \eqref{localpolynomials} can be written in the form
\begin{equation} \label{puiseux}
f_{\rho}(t,w_{\rho}) = b_{r_{\rho}}(t)\prod_{i=1}^{r''_{\rho}} (w_{\rho}-\beta^i_{\rho}(t))^{k_i}
\end{equation}
where $k_i \in \mathbb{Z}_{>0}$ are multiplicities of the roots, $ \sum_{i=1}^{r''_{\rho}} k_i = r_{\rho}$,  $b_{r_{\rho}}(t) \in \kk\llbracket t \rrbracket$ is the coefficient of $f_{\rho}(t, w_{\rho})$ at $w_{\rho}^{r_{\rho}}$,  and $\beta^i_{\rho}(t) \in \bar{\kk}(\!( t )\!)$ for $1 \leq i < j \leq r''_{\rho}$ are the distinct roots of $f_{\rho}(t, w_{\rho})$. 

\begin{cor} \label{puiseuxcor}
Let $\bar{\mathfrak{X}} \to \mathcal{S}$ be a toric degeneration of K3-s such that for any $\rho \in \bar{\mathscr{P}}^{[1]}$ in the expression \eqref{puiseux} for $f_{\rho}(t, w_{\rho}) \in \kk[w_{\rho}]\llbracket t \rrbracket$ we have $\beta^i_{\rho}(0) = \beta^j_{\rho} (0)$ for $1 \leq i < j \leq r''_{\rho}$ if and only if $\beta^i_{\rho}(0) = \beta^j_{\rho} (0) = \infty$. Then there exists a $c \in \mathbb{Z}_{>0}$ such that the basechange of $\bar{\mathfrak{X}} \to \mathcal{S}$ by $R \to R, \ t \mapsto t^c$ (here $t$ is the uniformizer of $R$) satisfies the assumptions of Lemma \ref{localmodelsk3}.
\end{cor}

\begin{proof}
Let $c \in \mathbb{Z}_{>0}$ be such that the $\beta_{\rho}^i(t)$ for all $\rho \in \bar{\mathscr{P}}^{[1]}$ and all $1 \leq i \leq r''_{\rho}$ are defined over $\kk(\!( t^{\frac{1}{c}} )\!)$. Then the basechange of $\bar{\mathfrak{X}} \to \mathcal{S}$ by $R \to R, \ t \mapsto t^c$ changes the local models in the neighbourhoods of a one-dimensional stratum $\bar{X}_{\rho}$ from  
$$\left\{xy = t^lf_{\rho}(t,w_{\rho}) \right\} \subseteq \Spec \kk[x,y,w_{\rho}]\llbracket t \rrbracket$$ for $f_{\rho}(t, w_{\rho}) \in \kk[w_{\rho}]\llbracket t \rrbracket$ of \eqref{puiseux} to 
$$\left\{xy = t^{c \cdot l}f_{\rho}(t^c,w_{\rho}) \right\} \subseteq \Spec \kk[x,y,w_{\rho}]\llbracket t \rrbracket$$
for $f_{\rho}(t^c, w_{\rho}) \in \kk[w_{\rho}]\llbracket t \rrbracket$ of the form
$$f_{\rho}(t^c,w_{\rho}) = b'_{r_{\rho}}(t)\prod_{i=1}^{r''_{\rho}} (w_{\rho}-\beta'^{i}_{\rho}(t))^{k_i}$$ where $b'_{r_{\rho}}(t) := b_{r_{\rho}}(t^c)$ is still defined over $\kk\llbracket t \rrbracket$ and the $\beta'^{i}_{\rho}(t) := \beta^{i}_{\rho}(t^c)$ for all $\rho \in \bar{\mathscr{P}}^{[1]}$ and all $1 \leq i \leq r''_{\rho}$ are now defined over $\kk(\!( t )\!)$. 

Suppose without loss of generality that $\beta_{\rho}^i(0) \neq \infty$ for $1 \leq i \leq r'_{\rho}$ and $\beta_{\rho}^i(0) = \infty$ for $r'_{\rho} < i \leq r''_{\rho}$ (for some $r'_{\rho}$ with $1 \leq r'_{\rho} \leq r''_{\rho}$ and suitably interpreted if $\beta_{\rho}^i(0) = \infty$ for $1 \leq i \leq r''_{\rho}$). Let $\gamma_{\rho}^i(t) := \beta'^{i}_{\rho}(t)$ for $1 \leq i \leq r'_{\rho}$ and let $$\alpha_{\rho}(t, w_{\rho}):= b'_{r_{\rho}}(t) \prod_{r'_{\rho} < i \leq r''_{\rho}} (w_{\rho}-\beta'^{i}_{\rho}(t))^{k_i}.$$ 
We claim that $\alpha_{\rho}(0, w_{\rho}) = \beta \in \kk \setminus \{0\}$. Indeed, $\kk[w_{\rho}]\llbracket t \rrbracket$ is a unique factorization domain and $$f_{\rho}(t^c,w_{\rho}) \in \kk[w_{\rho}]\llbracket t \rrbracket, \quad  \prod_{i =1}^{r'_{\rho}} (w_{\rho}-\beta'^{i}_{\rho}(t))^{k_i} \in \kk[w_{\rho}]\llbracket t \rrbracket.$$ So $\alpha_{\rho}(t, w_{\rho}) \in \kk[w_{\rho}]\llbracket t \rrbracket$ as well. Then we must have 
\begin{equation} \label{degreeineq}
\deg b'_{r_{\rho}}(t) \geq - \sum_{r'_{\rho} < i \leq r''_{\rho}}k_i \deg \beta'^{i}_{\rho}(t)
\end{equation}
(with the degree defined above) since the constant coefficient of $\alpha_{\rho}(t, w_{\rho})$ (viewed as a polynomial in $w_{\rho}$) has non-negative degree. It follows that all the non-constant coefficients of $\alpha_{\rho}(t, w_{\rho})$ are of some positive degrees. Therefore, 
$$\alpha_{\rho}(0, w_{\rho}) = b'_{r_{\rho}}(0) \prod_{r'_{\rho} < i \leq r''_{\rho}} (-\beta'^{i}_{\rho}(0))^{k_i} =: \beta \in \kk.$$ If $\beta = 0$, the inequality \eqref{degreeineq} is strict and the degree of the constant coefficient of $\alpha_{\rho}(t,w_{\rho})$ is positive. But then $\alpha_{\rho}(t,w_{\rho})$ is divisible by $t$. This implies that $f_{\rho}(t,w_{\rho})$ is also divisible by $t$, contradicting our assumptions. So $\alpha_{\rho}(0, w_{\rho}) = \beta \in \kk \setminus \{0\}$.

We have $$f_{\rho}(t^c,w_{\rho}) = \prod_{i=1}^{r'_{\rho}} (w_{\rho}-\gamma^i_{\rho}(t))^{k_i} \alpha_{\rho}(t, w_{\rho})$$ as in \eqref{localpolyexp} and $f_{\rho}(t^c,w_{\rho})$ satisfies the assumptions of Lemma \ref{localmodelsk3}. Since this holds for all $f_{\rho}(t^c,w_{\rho}), \ \rho \in \bar{\mathscr{P}}^{[1]}$, this implies the result.
\end{proof}

In the case of toric degenerations of K3-s, Definition \ref{special} of a special toric degeneration simplifies as follows.

\begin{fact} \label{whenspecialk3}
A toric degeneration $\bar{\mathfrak{X}} \to \mathcal{S}$ of K3-s is special if and only if it is a divisorial log deformation and the generic fibre of $\bar{\mathfrak{X}} \to \mathcal{S}$ is smooth. 
\end{fact}

\begin{proof}
Note that giving a toric log CY structure on $\check{\bar{X}}_0$ is equivalent to giving slab functions as in Proposition \ref{struct} and that the local rigidity condition is empty in dimension $2$. Therefore, condition (2) of Definition \ref{special} is empty in dimension $2$. The generic fibre of a special toric degeneration $\bar{\mathfrak{X}} \to \mathcal{S}$ is smooth by condition (1) of Definition \ref{special}. It is enough to check that a special toric degeneration $\bar{\mathfrak{X}} \to \mathcal{S}$ of K3-s is a divisorial log deformation. If $\bar{\mathfrak{X}} \to \mathcal{S}$ satisfies condition (3)(a) of Definition \ref{special}, then there is nothing to check. Condition (3)(b) of Definition \ref{special} says that $\bar{\mathfrak{X}} \to \mathcal{S}$ is distinguished (in the sense of Definition \ref{distinguished}). Therefore, it is enough to check that a distinguished toric degeneration $\bar{\mathfrak{X}} \to \mathcal{S}$ of K3-s is a divisorial log deformation.\footnote{Recall from Proposition \ref{distinguishedspecial} that this is true in any dimension provided that $\left(\bar B, \bar{\mathscr{P}}\right)$ is simple. The general claim only holds for toric degenerations of K3-s.} 

Let $\bar{\mathfrak{X}} \to \mathcal{S}$ be a distinguished toric degeneration of K3-s. It follows from \cite[Construction 2.7]{GS1} that for any $\rho \in \bar{\mathscr{P}}^{[1]}$ the polynomial $f_{\rho}(t,w_{\rho}) \in \kk[w_{\rho}]\llbracket t \rrbracket$ in the local model $$\left\{xy = t^lf_{\rho}(t,w_{\rho}) \right\} \subseteq \Spec \kk[x,y,w_{\rho}]\llbracket t \rrbracket$$ for $\bar{\mathfrak{X}} \to \mathcal{S}$ has constant coefficients, i.e. $f_{\rho}(t,w_{\rho})$ is of the form $f_{\rho}(t,w_{\rho}) = \sum_{j=0}^{r_{\rho}} b_j w_{\rho}^j$ for some $b_j \in \kk$ (this is no longer true in higher dimensions). But then we have 
$$f_{\rho}(t,w_{\rho}) = b_{r_{\rho}} \prod_{i=1}^{r''_{\rho}} (w_{\rho}-\beta^i_{\rho})^{k_i}$$ for $\beta^i_{\rho} \in \kk$ the distinct roots of $f_{\rho}(t,w_{\rho})$. Therefore, $f_{\rho}(t,w_{\rho})$ satisfies the assumptions of Lemma \ref{localmodelsk3} by taking $r'_{\rho} := r''_{\rho}$, $\gamma_{\rho}^i(t) := \beta_{\rho}^i$, and $\alpha_{\rho}(t, w_{\rho}):=b_{r_{\rho}}$. Since this holds for any $f_{\rho}(t,w_{\rho})$ with $ \rho \in \bar{\mathscr{P}}^{[1]}$, Proposition \ref{k3divisorial} implies that $\bar{\mathfrak{X}} \to \mathcal{S}$ is a divisorial log deformation.
\end{proof}

\begin{remarks} \label{basechangeconj}
We make a few remarks about extending Conjecture \ref{conjecture} to non-special toric degenerations of K3-s. 
\begin{enumerate}
\item Note that if Conjecture \ref{conjecture} holds for a toric degeneration $\bar{\mathfrak{X}} \to \mathcal{S}$ then it also holds for any toric degeneration $\bar{\mathfrak{X}}' \to \mathcal{S}$ such that $\bar{\mathfrak{X}} \to \mathcal{S}$ is the basechange of $\bar{\mathfrak{X}}' \to \mathcal{S}$ by $R \to R, \ t \mapsto t^c$ for some $c \in \mathbb{Z}_{>0}$ (since Conjecture \ref{conjecture} allows a finite basechange). Then our proof of Conjecture \ref{conjecture} for special toric degenerations $\bar{\mathfrak{X}} \to \mathcal{S}$ of K3-s in Chapter \ref{proof} along with Proposition \ref{whenspecialk3} above imply that Conjecture \ref{conjecture} holds for any toric degeneration $\bar{\mathfrak{X}} \to \mathcal{S}$ of K3-s satisfying the assumptions of Corollary \ref{puiseuxcor} and with a smooth generic fibre (i.e. with $k_i = 1$ in \eqref{puiseux} for all $1 \leq i \leq r''_{\rho}$). Such a degeneration need not be special, but unlike the case of special toric degenerations of K3-s, a basechange is necessary.
\item It would be interesting to see if one can construct resolutions of toric degenerations $\bar{\mathfrak{X}} \to \mathcal{S}$ of K3-s that are not divisorial log deformations. For instance, one can consider the case that $\bar{\mathfrak{X}} \to \mathcal{S}$ satisfies the assumptions of Lemma \ref{localmodelsk3} apart from requiring that the $\gamma^i_{\rho}(t) \in \kk \llbracket t \rrbracket$ of \eqref{localpolyexp} satisfy $\gamma^i_{\rho}(0) \neq \gamma^j_{\rho} (0)$ for $1 \leq i < j \leq r_{\rho}'$. The difficulty is that one can no longer use toric models to construct a resolution. If one finds such resolutions, Conjecture \ref{conjecture} for them should follow in the same way as in Chapter \ref{proof}. By the same argument as in (1), this would imply Conjecture \ref{conjecture} for all toric degenerations of K3-s with a smooth generic fibre. It is more likely that there is a way to relax (but not eliminate) the condition that a toric degeneration $\mathfrak{X} \to \mathcal{S}$ of K3-s has to be a divisorial log deformation in Conjecture \ref{conjecture}.
\end{enumerate}
Regardless of the (possible) extensions of Conjecture \ref{conjecture} to non-special degenerations of K3-s, in higher dimensions local models become much harder to control, so Definition \ref{special} of being special seems to be the most general way to guarantee well-behaved local models near points of the discriminant locus $Z$ (see Remark \ref{specialass}(3)).  
\end{remarks}

\subsection{Intrinsic setup} \label{intsetup}

We summarize two approaches to constructing intrinsic mirrors. See \cite{ann} for an overview, \cite{GS3} for the first construction, and \cite{GS4} for the second construction and equivalence of the two. 

Our exposition differs from this literature in two main aspects. First, we define the affine structure and the scattering diagram on the dual intersection complex $(B, \mathscr{P})$ of a log smooth degeneration $\mathfrak{X} \to \mathcal{S}$ directly instead of defining it on the cone $(\mathbf{C}B, \mathbf{C}\mathscr{P})$ over it. Second, \cite{GS4} is written under the simplifying assumption that the divisor $D$ is simple normal crossings which is also the main case of interest in \cite{GS3}. As we do not wish to restrict to this case, we extend the definition of the affine structure on $(B, \mathscr{P})$ to a more general situation. 

\subsubsection{Dual intersection complex} \label{intaff}

Let $\mathfrak{X} \overset{g}{\to} \mathcal{S}$ be a projective log smooth morphism of relative dimension $n$ where $\mathfrak{X}$ carries a fine and saturated Zariski log structure $\mathcal{M}_{\mathfrak{X}}$ and $\mathcal{S}$ is a regular one-dimensional scheme over $\Spec \kk$ with the divisorial log structure $\mathcal{M}_{\mathcal{S}}$ coming from a single closed point $0 \in \mathcal{S}$. In this thesis, we shall always have $\mathcal{S}:=\Spec R$ for $R$ a complete discrete valuation $\kk$-algebra as in Section \ref{tdsetup}.
Let $D_i, \ 1 \leq i \leq m$ be the components of the reduced central fibre $\left(\mathfrak{X}_0\right)_{\reduced}$. We will always assume that the log structure $\mathcal{M}_{\mathfrak{X}}$ on $\mathfrak{X}$ is divisorial with divisor $D:= \left(\mathfrak{X}_0 \right)_{\reduced} = D_1 + \dots + D_m$. So the log smooth morphism $g$ is of the form $$g: (\mathfrak{X}, D) \to (\mathcal{S}, 0).$$ We also require that $\mathfrak{X} \overset{g}{\to} \mathcal{S}$ is \emph{minimal log CY} (i.e. $K_{\mathfrak{X}} + D \equiv 0$, see Definition \ref{minlogcy}). 

Note that $(\mathfrak{X}, D)$ is Zariski and log smooth over $\Spec \kk$, so by \cite[Proposition 2.2]{ACGS}, it is a simple log scheme in the sense of Definition \ref{simplelogscheme}. As in the case of toric degenerations (see Section \ref{tdaff}), we apply the tropicalization functor $\Sigma$ of Construction \ref{tropicalize} to produce a map of rational cone complexes\footnote{It is a map of actual cone complexes (and not the generalized complexes of Construction \ref{tropicalize}) since $(\mathfrak{X}, D)$ is a simple log scheme.}: 
\begin{equation} \label{inttrop}
g_{\trop}: \Sigma (\mathfrak{X}) \to \Sigma(\mathcal{S}) =\mathbb{R}_{\geq 0}.
\end{equation}
We define the \emph{dual intersection complex} of $\mathfrak{X} \to \mathcal{S}$ as 
\begin{equation} \label{dicims}
(B, \mathscr{P}): = g_{\trop}^{-1}(1)
\end{equation}
where the polyhedral structure $\mathscr{P}$ comes from restricting the cones of $\Sigma \left(\mathfrak{X}\right)$ to the fibre over $1 \in \mathbb{R}_{\geq 0}$. Clearly, we have $\mathbf{C}B \cong \Sigma\left(\mathfrak{X}\right)$ where $\mathbf{C}B$ is the cone over $B$ of Section \ref{thetafunctions}.

\begin{ass} \label{manifoldass}
We assume that the dual intersection complex $\left( B, \mathscr{P}\right)$ of $g: \mathfrak{X} \to \mathcal{S}$ is a polyhedral manifold of dimension $n$ in the sense of Definition \ref{polymani}.
\end{ass}

\begin{remark} \label{nonreducedeness}
Note that if $g^{-1}(0)$ is not reduced, $\left( B, \mathscr{P}\right)$ is not a polyhedral manifold. Indeed, let $\rho_i \in \Sigma(\mathfrak{X})$ be the ray corresponding to a divisor $D_i$ of multiplicity $k > 1$ in $g^{-1}(0)$. Then Construction \ref{tropicalize} of tropicalization implies that $v_i:= \rho_i \cap g_{\trop}^{-1}(1)$ is not an integral point of $\rho_i$ (instead, we have $v_i \in \Sigma(\mathfrak{X})\left(\frac{1}{k}\mathbb{Z}\right)$). This is not a major issue as one may extend Definition \ref{polymani} to this case by allowing rational polyhedra in the construction instead of integral polyhedra. With this modification, all the constructions go through.

The morphism $\mathfrak{X} \to \mathcal{S}$ arising from an \emph{admissible} (see Definition \ref{admissibledef}) resolution $\mathfrak{X} \to \bar{\mathfrak{X}}$ in Conjecture \ref{conjecture}  may have a non-reduced central fibre. However, Proposition \ref{furtherresolution} implies that it is enough to prove Conjecture \ref{conjecture} for a more restricted class of \emph{strongly admissible} resolutions (see Definition \ref{stronglyadmissible}) which have a reduced central fibre.
\end{remark}

\begin{notation}
Similarly to the case of toric degenerations (see Notation \ref{otonot}), the fact that $(\mathfrak{X}, D)$ is a simple log scheme implies that there is an inclusion-reversing correspondence between the logarithmic strata of $\mathfrak{X}_0$ and the cells of $\mathscr{P}$. We denote the stratum corresponding to a cell $\sigma \in \mathscr{P}$ by $X_{\sigma}$. We make an exception for vertices $v \in \mathscr{P}^{[0]}$ where we denote the corresponding divisor by $D_v$.
\end{notation} 

\begin{remark} \label{coneremark1}
Both \cite{GS3} and \cite{GS4} work with $\mathbf{C}B \cong \Sigma (\mathfrak{X})$, the cone over $B$, directly (see Section \ref{thetafunctions}). To prove Conjecture \ref{conjecture}, we will need to relate, for a resolution $\pi: \mathfrak{X} \to \bar{\mathfrak{X}}$, the scattering diagram on $\bar{B}$ defining the toric degeneration mirror $\check{\bar{\mathfrak{X}}} \to \Spec \kk \llbracket t \rrbracket$ and a certain scattering diagram defining the intrinsic mirror $\check{\mathfrak{X}} \to \Spec \widehat{\kk[P]}$. So it is more convenient to work on $B$. Note that our notation differs from \cite{GS3,GS4} where $B$ denotes the manifold that we denote $\mathbf{C}B$.
\end{remark}

In the case that $D$ is a simple normal crossings divisor, there is a straightforward description of $\Sigma (\mathfrak{X})$ and $g_{\trop}$ in terms of $\Div_D(\mathfrak{X})^*_{\mathbb{R}}:= \Hom(\Div_D(\mathfrak{X}), \mathbb{R})$ where $\Div_D(\mathfrak{X}):=\bigoplus_{1 \leq i \leq m} \mathbb{Z} D_i$ is the group of divisors supported on $D$. 

\begin{cons} \label{tropicalizesnc}
Assume that $D$ is simple normal crossings. Then we have $\overline{\mathcal{M}}_{\mathfrak{X}} = \bigoplus_{1 \leq i \leq m} \mathbb{N}_{D_i}$ where $\mathbb{N}_{D_i}$ is the constant sheaf on $D_i$ with stalk $\mathbb{N}$. In particular, there is a natural isomoprhism $\Gamma\left(\mathfrak{X}, \overline{\mathcal{M}}_{\mathfrak{X}}^{\gp} \right) \cong \Div_{D}(\mathfrak{X})$. 

Assumption \ref{manifoldass} implies (by condition (2)(d) of Definition \ref{polymani}) that for
any index set $I\subseteq \{1,\ldots,m\}$ the (possibly empty)
\begin{equation} \label{DI}
D_I:=\bigcap_{i\in I} D_i
\end{equation}
is connected. Following \cite[Example 1.4]{GS3}, we can identify $\Sigma(\mathfrak{X})$
with the collection of cones
\[
\Sigma(\mathfrak{X}):=\left\{ \sum_{i\in I} \mathbb{R}_{\ge 0} D_i^*\,\bigg|\,
\hbox{$I\subseteq \{1,\ldots,s\}$ an index set with
$D_I \neq \varnothing$}\right\} \subseteq \Div_D(\mathfrak{X})^*_{\mathbb{R}}.
\]
with integral structure
\[
\Sigma(\mathfrak{X})(\mathbb{Z})=\left\{\sum a_i D_i^*\,\bigg|\,
a_i\in \mathbb{N}, \bigcap_{i: a_i>0} D_i \neq \varnothing\right\}.
\]
Under this identification, the map $g_{\trop}: \Sigma (\mathfrak{X}) \to \mathbb{R}_{\geq 0}$ is just the restriction of the linear map $\Div_D(\mathfrak{X})^*_{\mathbb{R}}\rightarrow \mathbb{R}$ given by evaluation
on the (not necessarily reduced) divisor $g^*(0)$.
\end{cons}

\subsubsection{Technical assumptions} \label{technicalasssection}
Consider a minimal log CY degeneration $\mathfrak{X} \to \mathcal{S}$ (with $(\mathfrak{X}, D)$ a Zariski log scheme) satisfying Assumption \ref{manifoldass}. 
The setups of \cite{GS3, GS4} have some additional technical assumptions that we now discuss. 

The intrinsic mirror constructions in \cite{GS3, GS4} rely heavily on the theory of \emph{punctured log Gromov-Witten invariants} developed in \cite{ACGS2} under the additional assumption that the ghost sheaf $\overline{\mathcal{M}}_{\mathfrak{X}}$ of the log structure $\mathcal{M}_{\mathfrak{X}}$ on $\mathfrak{X}$ is \emph{globally generated} in the sense that the natural map $\Gamma(\mathfrak{X}, \, \overline{\mathcal{M}}^{\gp}_{\mathfrak{X}} \otimes_{\mathbb{Z}} \mathbb{Q}) \to \overline{\mathcal{M}}^{\gp}_{\mathfrak{X}, x} \otimes_{\mathbb{Z}} \mathbb{Q}$ is surjective for every $x \in \mathfrak{X}$. For instance, this assumption holds in the case that $D$ is a simple normal crossings divisor.\footnote{This immediately follows from Construction \ref{tropicalizesnc} and \cite[Proposition 3.12]{ACGS2}.}

\begin{lemma} \label{topsphere}
Suppose that $\mathfrak{X} \to \mathcal{S}$ is a minimal log CY degeneration of K3-s satisfying Assumption \ref{manifoldass}. Then topologically, $B$ is a sphere.
\end{lemma}

\begin{proof}
We will only be interested in the case that $\mathfrak{X} \to \mathcal{S}$ is obtained by composing a toric degeneration $\bar{\mathfrak{X}} \to \mathcal{S}$ with an admissible (see Definition \ref{admissibledef}) resolution $\pi:\mathfrak{X} \to \bar{\mathfrak{X}}$. In this case, $(B, \mathscr{P})$ is a subdivision of the dual intersection complex $\left(\bar{B}, \bar{\mathscr{P}}\right)$ of $\bar{\mathfrak{X}} \to \mathcal{S}$ and the claim follows from Proposition \ref{spheredic}.

In general, rather than generalizing the proof of Proposition \ref{spheredic} and getting involved in the details of the proof of \cite[Proposition 2.37]{GS2}, we reduce to the classification of semi-stable degenerations of K3-s (see, e.g. \cite[Theorem 5.2]{F1}). 

First, note that $\mathfrak{X} \to \mathcal{S}$ is flat and proper. Indeed, $\mathfrak{X} \to \mathcal{S}$ is projective, thus proper. Since $\mathfrak{X} \to \mathcal{S}$ is log smooth, it is also log flat, see \cite[Part IV, Proposition 4.1.2(1)]{Ogus}. Further, a log flat and integral morphism is flat by \cite[Part IV, Proposition 4.3.5(1)]{Ogus}. But $\mathcal{S}$ is one-dimensional so $\mathfrak{X} \to \mathcal{S}$ is always integral by \cite[Part III, Proposition 2.5.3(3)]{Ogus}. So $\mathfrak{X} \to \mathcal{S}$ is flat. 

By Assumption \ref{manifoldass}, $(B, \mathscr{P})$ is a polyhedral manifold of dimension $2$. Recall from the introduction that by \cite[Theorem 2.4.1 and Corollary 2.6.7]{AW}, there is a one-to-one correspondence between logarithmic modifications of $\mathfrak{X} \to \mathcal{S}$ and subdivisions of $(B, \mathscr{P})$. Consider a subdivision $(\tilde{B}, \tilde{\mathscr{P}})$ of $(B, \mathscr{P})$ with all the cells $\sigma \in \tilde{\mathscr{P}}^{\max}$ standard triangles (it exists by e.g. Lemma \ref{pullingresolution}). Let $\tilde{\mathfrak{X}} \to \mathfrak{X}$ be the corresponding logarithmic modification. Then the composed $\tilde{\mathfrak{X}} \to \mathcal{S}$ is proper and flat (since a logarithmic modification is proper and log étale). Moreover, the central fibre $\tilde{\mathfrak{X}}_0$ of $\tilde{\mathfrak{X}} \to \mathcal{S}$ is simple normal crossings since $(\tilde{B}, \tilde{\mathscr{P}})$ is the dual intersection complex of $\tilde{\mathfrak{X}} \to \mathcal{S}$ and all the cells $\sigma \in \tilde{\mathscr{P}}^{\max}$ are standard triangles. Therefore, $\tilde{\mathfrak{X}} \to \mathcal{S}$ is a semi-stable degeneration of K3-s.\footnote{This construction is an instance of the semi-stable reduction theorem of \cite[Chapter II]{KKMS} which shows that every degeneration $X \to \mathcal{S}:= \Spec R$ admits a semi-stable model $X' \to \mathcal{S}'$ via a finite basechange $\mathcal{S}' \to \mathcal{S}$ and a blowup.} 

Since $(\tilde{B}, \tilde{\mathscr{P}})$ is a subdivision of $(B, \mathscr{P})$, it is enough to show that $\tilde{B}$ is a topological sphere. Indeed, since we have $\dim \tilde{B} = \dim B = 2$, the classification of semi-stable degenerations of K3-s (see, e.g. \cite[Theorem 5.2]{F1}) implies that $\tilde{\mathfrak{X}} \to \mathcal{S}$ is a Type III degeneration and $(\tilde{B}, \tilde{\mathscr{P}})$ is a triangulation of a sphere.
\end{proof}

\begin{fact} \label{propglobgen}
Let $\mathfrak{X} \to \mathcal{S}$ be a minimal log CY degeneration of K3-s satisfying Assumption \ref{manifoldass}. Then the ghost sheaf $\overline{\mathcal{M}}_{\mathfrak{X}}$ of the log structure $\mathcal{M}_{\mathfrak{X}}$ on $\mathfrak{X}$ is globally generated.
\end{fact}

\begin{proof}
By \cite[Proposition 3.12]{ACGS2}, global generation of $\overline{\mathcal{M}}_{\mathfrak{X}}$ is equivalent to having a piecewise-linear (PL) map $|\Sigma(\mathfrak{X})| \to \mathbb{R}^r$ for some $r \in \mathbb{N}$ that is injective when restricted to each cone $\sigma \in \Sigma(\mathfrak{X})$. It is enough to find a PL-embedding $B \to \mathbb{R}^r$. Indeed, given a PL-embedding $i:B \to \mathbb{R}^r$ we can define a PL-embedding $i: \Sigma(\mathfrak{X}) \cong \mathbf{C}B \to \mathbb{R}^{r+1}$ as follows. For any $\mathbf{C}\sigma \in \mathbf{C}\mathscr{P}^{[k]}$ and any $(a, x) \in \mathbf{C}\sigma = \mathbb{R}_{\geq 0} \cdot (\{1\}\times \sigma) \subseteq \mathbb{R} \times \mathbb{R}^k$ we set $i(\alpha)(a, x) := (a, \alpha(x)) \in \mathbb{R} \times \mathbb{R}^r$. 

$B$ is a topological sphere by Lemma \ref{topsphere}. Since $B$ has the structure of a polyhedral manifold of dimension $2$, by Steinitz's theorem (see, e.g. \cite{R}), $B$ is PL-isomorphic to the boundary complex of a polytope $\Delta \subseteq \mathbb{R}^3$. This gives rise to a PL-embedding $B \to \mathbb{R}^3$ with image $\Delta$. 
\end{proof}

In higher dimensions, the situation is more complicated. As the proof of Proposition \ref{propglobgen} shows, to prove that the ghost sheaf $\overline{\mathcal{M}}_{\mathfrak{X}}$ of the log structure $\mathcal{M}_{\mathfrak{X}}$ on $\mathfrak{X}$ is globally generated, it is enough to show that $B$ admits a PL-embedding $B \to \mathbb{R}^r$ for some $r \in \mathbb{N}$. It is well-known that every simplicial complex of dimension $n$ is \emph{geometrically realizable} in $\mathbb{R}^{2n+1}$. This implies the condition for any $\mathfrak{X} \to \mathcal{S}$ with all $\sigma \in \mathscr{P}^{\max}$ simplices (including the case that $D$ is simple normal crossings). However, already in dimension $3$, there exist polyhedral complexes (even if we only consider topological spheres) with no PL-embedding into $\mathbb{R}^r$ for any $r \in \mathbb{N}$ (see \cite{Ca, HZ}). There are also some positive embeddability results (see, e.g. \cite{PW, Z}).

During the writing of this thesis, the paper \cite{J} came out, showing the birational invariance of punctured log Gromov-Witten invariants. As in the case of ordinary stable logarithmic maps, this removes the need for the assumption that $\overline{\mathcal{M}}_{\mathfrak{X}}$ is globally generated in punctured log Gromov-Witten theory (see \cite{AC, GSlog} for stable logarithmic maps and \cite{ACMW} for the construction). Indeed, this assumption is only necessary to show \emph{boundedness} of $\mathscr{M}(\mathfrak{X}/\mathcal{S}, \beta)$ (the stack of all basic stable punctured maps to $\mathfrak{X}/\mathcal{S}$ of class $\beta$), i.e. to show that the natural map $\mathscr{M}(\mathfrak{X}/\mathcal{S}, \beta) \to \mathcal{S}$ is of finite type (see \cite[Section 3.3]{ACGS2} for details). The construction of \cite{ACMW} gives, for any $\mathfrak{X} \to \mathcal{S}$, a logarithmic modification $\mathcal{Y} \to \mathfrak{X}$ with $\overline{\mathcal{M}}_{\mathcal{Y}}$ globally generated. Then by the same arguments as in \cite[Section 5]{ACMW}, boundedness of $\mathscr{M}(\mathcal{Y}/\mathcal{S}, \beta)$ along with birational invariance of $\mathfrak{M}(\mathfrak{X}/\mathcal{S}, \beta)$ (the \emph{Artin stack} associated to $\mathscr{M}(\mathfrak{X}/\mathcal{S}, \beta)$) imply boundedness of $\mathscr{M}(\mathfrak{X}/\mathcal{S}, \beta)$.

The discussion above allows us to remove the assumption that $\overline{\mathcal{M}}_{\mathfrak{X}}$ is globally generated from \cite{GS3}. In \cite{GS3}, there are no other assumptions. We now discuss the assumptions of \cite{GS4}. First, \cite{GS4} assumes that the $D_I$ of \eqref{DI} is connected for
any index set $I\subseteq \{1,\ldots,m\}$. Moreover, \cite[Assumptions 1.2]{GS4} impose some conditions on $(\mathbf{C}B, \mathbf{C}\mathscr{P})$ and on the map $\mathfrak{X} \to \mathcal{S}$.

\begin{fact} \label{assofcanonical}
Let $\mathfrak{X} \to \mathcal{S}$ be a minimal log CY degeneration satisfying Assumption \ref{manifoldass}. Then $D_I$ is connected  for
any index set $I\subseteq \{1,\ldots,m\}$ and $\mathfrak{X} \to \mathcal{S}$ satisfies \cite[Assumptions 1.2]{GS4}.
\end{fact}

\begin{proof}
By Assumption \ref{manifoldass}, $\left( B, \mathscr{P}\right)$ is a polyhedral manifold in the sense of Definition \ref{polymani}. The first claim immediately follows from condition (2)(d) of Definition \ref{polymani}. Condition (2) of \cite[Assumptions 1.2]{GS4} is satisfied since for any $\rho \in \mathscr{P}$\footnote{Here $\rho \in \mathscr{P}$ is any cell, not necessarily of dimension $1$, following the notations of \cite[Assumptions 1.2]{GS4}.} the corresponding stratum $X_{\rho}$ is connected and the restriction $g|_{X_{\rho}}$ of $g: \mathfrak{X} \to \mathcal{S}$ is constant. So we only need to check that $(\mathfrak{X}, D)$ satisfies \cite[Assumptions 1.1]{GS4}. It satisfies condition (1) of \cite[Assumptions 1.1]{GS4} since $\left( B, \mathscr{P}\right)$ is of pure dimension $n$ by condition (1) of Definition \ref{polymani}. For conditions (2) and (3) of \cite[Assumptions 1.1]{GS4} note that since $\mathfrak{X} \to \mathcal{S}$ is minimal log CY, we have $K_{\mathfrak{X}}+D \equiv 0$ and all the strata $D_I$ for $I\subseteq \{1,\ldots,m\}$ are \emph{good} in the sense of \cite[Section 1.1]{GS4}. This immediately implies condition (2) of \cite[Assumptions 1.1]{GS4}. Moreover, using Construction \ref{tropicalize} of the tropicalization functor and the fact that $\left( B, \mathscr{P}\right)$ is a polyhedral manifold in the sense of Definition \ref{polymani}, it also implies that condition (3) of \cite[Assumptions 1.1]{GS4} is equivalent to the claim that the boundary of any polyhedron $\check{\sigma} \in \check{\mathscr{P}}^{[k]}, \ k > 1$ of the dual complex $\left(\check{B}, \check{\mathscr{P}}\right)$ of $\left( B, \mathscr{P}\right)$ is connected. This is indeed the case for any bounded polyhedron of dimension at least $2$.
\end{proof}

Finally, \cite{GS4} requires $D$ to be simple normal crossings. However, this assumption can easily be removed and only requires a generalization of the construction of the affine structure on $(B, \mathscr{P})$ that we describe in Section \ref{affstrgeneral} using \cite[Theorem 4.1]{W}. Indeed, one can check that all the arguments of \cite{GS4} work in this generalized setting (using boundedness of $\mathscr{M}(\mathfrak{X}/\mathcal{S}, \beta)$ in the case that $\overline{\mathcal{M}}_{\mathfrak{X}}$ is not globally generated, as discussed above).

\begin{remark}
There is a different way to see this generalization in the case that $\mathfrak{X} \to \mathcal{S}$ admits a logarithmic modification $\pi': \tilde{\mathfrak{X}} \to \mathfrak{X}$ with the divisor $(\pi')^{-1}D$ defining the pullback log structure on $\tilde{\mathfrak{X}}$ simple normal crossings and all the exceptional divisors toric varieties. In this case, the dual intersection complex $(\tilde{B}, \tilde{\mathscr{P}})$ of $\tilde{\mathfrak{X}} \to \mathcal{S}$ is a subdivision of $(B, \mathscr{P})$ (see \cite[Theorem 2.4.1 and Corollary 2.6.7]{AW})  and the affine structure on $(\tilde{B}, \tilde{\mathscr{P}})$ extends across each vertex corresponding to an exceptional divisor by \cite[Proposition 2.3]{AG}. This gives rise to an affine structure on $(B, \mathscr{P})$, and one can check that this is the same affine structure we describe in Section \ref{affstrgeneral}. Moreover, the \emph{canonical scattering diagram} (see Section \ref{canonicalmirror}) $\mathfrak{D}_I$ on $(B, \mathscr{P})$ is equivalent to the natural pushforward $\left((\pi')_{*} \tilde{\mathfrak{D}}\right)_{I}$ of the canonical scattering diagram $\tilde{\mathfrak{D}}_{(\pi')_{*}^{-1}(I)}$ on $(\tilde{B}, \tilde{\mathscr{P}})$ by an argument as in the proof of \cite[Corollary 1.7]{J}. 

In the case that $\mathfrak{X} \to \mathcal{S}$ is a strongly admissible resolution of a special toric degeneration of K3-s (see Definition \ref{stronglyadmissible}), $\mathfrak{X} \to \mathcal{S}$ always admits a logarithmic modification of the form above by Corollary \ref{furtherresolutionofadmissible}. Proposition \ref{furtherresolution} implies that it is enough to prove Conjecture \ref{conjecture} for strongly admissible resolutions. Therefore, this remark's description of the affine structure and the canonical scattering diagram on $(B, \mathscr{P})$ suffices for special toric degenerations of K3-s. 
\end{remark}

To summarize, we will only need to check that the degenerations $\mathfrak{X} \to \mathcal{S}$ we construct are minimal log CY and satisfy Assumption \ref{manifoldass}. One may wish to think of the case that $D$ is simple normal crossings to avoid some of the complexities of the general situation. We do not restrict to this case, allowing us to have more flexibility in Definition \ref{admissibledef} of an admissible resolution. We will show that every special toric degeneration $\bar{\mathfrak{X}} \to \mathcal{S}$ of K3-s admits an admissible resolution with $D$ simple normal crossings in Proposition \ref{strongadmissibleresolutionresult}. We generalize to special toric degenerations $\bar{\mathfrak{X}} \to \mathcal{S}$ of CY threefolds that are distinguished and have a simple dual intersection complex $\left(\bar{B}, \bar{\mathscr{P}}\right)$ in Corollary \ref{sncafterfurtherblowupgeneral} and Proposition \ref{basechangstronglyadmissiblesnc}. The general case is similar, see the discussion at the end of Section \ref{resolveindimatleast4}. 

\subsubsection{Affine structure in the simple normal crossings case} \label{affstrsnc}

We now need to put the structure of an affine manifold with singularities on $(B, \mathscr{P})$. We first describe the construction in the case that the central fibre $\mathfrak{X}_0$ of $\mathfrak{X} \to \mathcal{S}$ is reduced, $D$ is a simple normal crossings divisor and the $D_I$ of \eqref{DI} is connected for
any index set $I\subseteq \{1,\ldots,m\}$. It is easy to see that these assumptions imply that $\mathfrak{X} \to \mathcal{S}$ satisfies Assumption \ref{manifoldass}.

\begin{cons} \label{affstrint}
Suppose that the central fibre of $\mathfrak{X} \to \mathcal{S}$ is reduced and $D$ is simple normal crossings, and let $(B, \mathscr{P})$ be the dual intersection complex of $\mathfrak{X} \to \mathcal{S}$ which is a polyhedral manifold by Assumption \ref{manifoldass}. We let the discriminant locus $\Delta$ be the union of cells of codimension $2$:
\begin{equation} \label{discriminantlocusims}
\Delta := \bigcup_{\mathscr{P}^{[n-2]}} \sigma.
\end{equation}
We want to put an affine structure on $B_0 :=B \setminus  \Delta$.

Following the general framework of Construction \ref{affchartstool}, as in the case of toric degenerations, all the $\sigma \in\mathscr{P}^{\max}$ have natural structures of integral polyhedra via their inclusion into the cones of $\Sigma \left(\mathfrak{X}\right)$. In other words, these are precisely the integral polyhedra of the polyhedral manifold structure on $(B, \mathscr{P})$. Our assumptions and Construction \ref{tropicalizesnc} imply that all the $\sigma \in \mathscr{P}$ are standard simplices. 

Using the refined description of \ref{affchartstool}$(2')$, we now need to give, for any $\rho \in \mathscr{P}^{[n-1]}$ with $\rho = \sigma \cap \sigma'$ for $ \sigma, \sigma' \in \mathscr{P}^{\max}$ and any $v \subseteq \rho, \ v \in \mathscr{P}^{[0]}$, an integral PL-embedding 
\begin{equation} \label{intcharts}
\psi_{\rho, v}: \sigma \cup \sigma' \to \mathbb{R}^n  
\end{equation} 
with $\psi_{\rho, v}(v) = 0$ compatible with the structures of standard simplices on $\sigma \in \mathscr{P}^{\max}$.

We have $\rho = \langle v_{i} \ | \ 0 \leq i \leq n-1 \rangle$ for some $v_{i} \in \mathscr{P}^{[0]}$ with $v_0:=v$ (without loss of generality). Suppose that $\sigma = \langle \rho, \ v_{n} \rangle, \ \sigma' = \langle \rho, \ v'_{n} \rangle$ for some $v_n, \ v_n' \in \mathscr{P}^{[0]}$ and choose integral bases $(e_1, \dots, e_n)$ and $(e_1, \dots, e_{n-1}, e'_n)$ of $\mathbb{R}^n$ subject to the constraint that 
\begin{equation} \label{constraintaffine}
e_n + e'_n = - \sum_{i=1}^{n-1} \left( D_{v_{i}}  \cdot X_{\rho} \right)e_i.  
\end{equation}
Note that the intersection numbers are defined. Since $D$ is simple normal crossings, $D_{v_i}$ is a Cartier divisor and Assumption \ref{manifoldass} ensures that $X_{\rho}$ is proper. 

We define the embedding $\psi_{\rho, v_0}$ by requiring it to be linear on each cone and by setting 
\begin{equation} \label{mapaffstr}
\psi_{\rho, v_0}(v_{0}) = 0, \quad \psi_{\rho, v_0}(v_{i}) = e_i, \ 1 \leq i \leq n, \quad \psi_{\rho, v_0}(v'_{n}) = e'_n. 
\end{equation}

As in Construction \ref{affchartstool}, this data gives an integral affine structure on $B_0 =B \setminus \Delta$. The affine structure of Construction \ref{affchartstool} extends to the complement of $\Delta$ since the embeddings $\psi_{\rho, v_i}$ for $0 \leq i \leq n-1$ are compatible. This follows by an easy computation using Lemma \ref{intcentralfibrelemma} below, which gives a natural generalization of the triple point formula of \cite[p. 39]{P}. Alternatively, it follows from compatibility of the affine structure on $B$ with the construction of \cite[Section 1.3]{GS4} of the affine structure on $\mathbf{C}B$ under the natural inclusion $B \xhookrightarrow{} \mathbf{C}B$, that we show in Proposition \ref{indofchoices} below. 

Note that, as in the discussion after Construction \ref{affchartstool}, the affine structure on $B$ admits the description of Construction \ref{affchartstool}$(2)$ Namely, for every $v \in \mathscr{P}^{[0]}$, the restrictions of the embeddings $\psi_{v, \rho}$ for all the $\rho \in \mathscr{P}^{[n-1]}$ with $v \subseteq \rho$ along with the restrictions of the charts on the maximal cells $\sigma \in \mathscr{P}^{\max}$ define an affine structure on $W_v \setminus \Delta$ for $W_v$ as in \eqref{openstar}. We shall freely use this description when we need to.
\end{cons}

We now show that Construction \ref{affstrint} is compatible with the construction of \cite[Section 1.3]{GS4} of the affine structure on $\mathbf{C}B$. 

\begin{lemma} \label{intcentralfibrelemma}
In the setup of Construction \ref{affstrint}, we have 
\begin{equation} \label{intcentralfibre}
\sum_{i=0}^{n-1} \left( D_{v_{i}}  \cdot X_{\rho} \right)=-2.
\end{equation}
\end{lemma}

\begin{proof}
Note that $D$ is numerically equivalent to the trivial divisor, so we have 
\begin{equation} \label{intcentralfibre1}
\sum_{j=1}^m \left( D_j \cdot X_{\rho} \right) = \left(\sum_{j=1}^m D_j \right) \cdot X_{\rho} = D \cdot X_{\rho} = 0.
\end{equation}
Also, if $D_j \cap X_{\rho} \neq \varnothing$, we have $D_j \cdot X_{\rho} = 0$. Therefore \eqref{intcentralfibre1} implies that $$\sum_{i=0}^{n-1} \left( D_{v_{i}}  \cdot X_{\rho} \right) + D_{v_n} \cdot X_{\rho} + D_{v'_n} \cdot X_{\rho} = 0.$$ Now \eqref{intcentralfibre} follows by noting that $D$ being simple normal crossings implies that $D_{v_n} \cdot X_{\rho} = D_{v'_n} \cdot X_{\rho} = 1$.
\end{proof}

Compatibility with the construction of \cite[Section 1.3]{GS4} follows from an easy computation.

\begin{fact} \label{indofchoices}
The affine structure of Construction \ref{affstrint} is the affine structure on $B$ induced by the inclusion $B \hookrightarrow \mathbf{C}B$ (using the affine structure of \cite[Section 1.3]{GS4} on $\mathbf{C}B$).
\end{fact}

\begin{proof}
In the conventions of Construction \ref{affstrint}, the affine structure of \cite[Section 1.3]{GS4} on $\mathbf{C}B$ is defined via the natural charts on $$\left\{\Int \left(\mathbf{C}\sigma \right) \ | \ \mathbf{C}\sigma \in \mathbf{C}\mathscr{P}^{\max} \right\}$$ induced by Construction \ref{tropicalize} of tropicalization and the charts on $$\left\{ \Int \left(\mathbf{C}\sigma \cup \mathbf{C}\sigma' \right) \ | \ \mathbf{C}\sigma, \mathbf{C}\sigma' \in \mathbf{C}\mathscr{P}^{\max}, \ \mathbf{C}\sigma \cap \mathbf{C}\sigma' = \mathbf{C}\rho \in \mathbf{C}\mathscr{P}^{[(n+1)-1]} \right\}$$ given by selecting integral bases $(\bar{e}_0, \bar{e}_1, \dots, \bar{e}_n)$ and $(\bar{e}_0, \bar{e}_1, \dots, \bar{e}_{n-1}, \bar{e}'_n)$ of $\mathbb{R}^{n+1}$ subject to the constraint that 
\begin{equation} \label{constraintaffinecone}
\bar{e}_n + \bar{e}'_n = - \sum_{i=0}^{n-1} \left( D_{v_{i}}  \cdot X_{\rho} \right)\bar{e}_i,  
\end{equation}
and defining the embedding 
\begin{equation} \label{intchartscone}
\psi_{\mathbf{C\rho}}: \mathbf{C}\sigma \cup \mathbf{C}\sigma' \to \mathbb{R}^{n+1}  
\end{equation}
by requiring it to be linear on each cone and by setting 
\begin{equation} \label{mapaffstrcone}
\psi_{\mathbf{C\rho}}(v_{i}) = \bar{e}_i, \ 0 \leq i \leq n, \quad \psi_{\mathbf{C\rho}}(v'_{n}) = \bar{e}'_n.  
\end{equation}

Construction \ref{tropicalize} of tropicalization implies that the charts on the maximal cells of $\mathbf{C}B$ induce the corresponding charts on $B$ defined via Construction \ref{affstrint}. Passing from a chart of $\mathbf{C}B$ on $\Int \left(\mathbf{C}\sigma \cup \mathbf{C}\sigma' \right)$ to a chart of $B$ on $\Int (\sigma \cup \sigma')$ sending $v_0$ to $0 \in \mathbb{R}^n$ corresponds to setting $$e_i := \bar{e}_i - \bar{e}_0, \ 1 \leq i \leq n, \quad e'_n := \bar{e}'_n - \bar{e}_0.$$ This recovers \eqref{mapaffstr} from \eqref{mapaffstrcone} . Also, we have 
\begin{multline}
e_n + e'_n = \bar{e}_n - \bar{e}_0 + \bar{e}'_n - \bar{e_0} = - \sum_{i=0}^{n-1} \left( D_{v_{i}}  \cdot X_{\rho} \right)\bar{e}_i - 2\bar{e}_0 = \\ = -\sum_{i=0}^{n-1} \left( D_{v_{i}}  \cdot X_{\rho} \right)\left(\bar{e}_i - \bar{e}_0\right) - \left(\sum_{i=0}^{n-1} \left( D_{v_{i}}  \cdot X_{\rho}\right) + 2 \right) \bar{e}_0  =  - \sum_{i=1}^{n-1} \left( D_{v_{i}}  \cdot X_{\rho} \right)e_i
\end{multline}
(here we used \eqref{constraintaffinecone} for the second equality and Lemma \ref{intcentralfibrelemma} for the last one) which recovers \eqref{constraintaffine}. 

So the charts of \cite[Section 1.3]{GS4} of the affine structure on $\mathbf{C}B$ induce by inclusion the charts of Construction \ref{affstrint} of the affine structure on $B$, which implies the claim. 
\end{proof}

\begin{remark} \label{logremark}
The affine structure of \cite[Section 1.3]{GS4} on $\mathbf{C}B$ is natural from the logarithmic point of view. For every $\mathbf{C}\rho \in \mathbf{C}\mathscr{P}^{[(n+1)-1]}$, the log scheme $X_{\mathbf{C}\rho}$ with the log structure induced from $(\mathfrak{X}, D)$ by inclusion is isomorphic as a log scheme to the stratum $X_{\Sigma_{\mathbf{C}\rho}, \psi_{\mathbf{C}\rho}(\mathbf{C}\rho)}$ of the toric variety $X_{\Sigma_{\mathbf{C}\rho}}$ (with the standard toric log structure) constructed from the fan $\Sigma_{\mathbf{C}\rho}$ that consists of $\psi_{\mathbf{C}\rho}(\mathbf{C}\sigma), \ \psi_{\mathbf{C}\rho}(\mathbf{C}\sigma')$, and their faces. This follows from a direct computation of the sheaves defining the log structures induced on $X_{\mathbf{C}\rho}$ and $X_{\Sigma_{\mathbf{C}\rho}, \psi_{\mathbf{C}\rho}(\mathbf{C}\rho)}$ by inclusion. A similar statement holds in the absolute case of a log CY variety $(\mathfrak{X}, D)$. We refer to \cite[Lemma 1.10]{GS4} for details.
\end{remark}

We now specialize to the case $n=2$ corresponding to minimal log CY degenerations of K3-s and show that in this case, the affine structure of Construction \ref{affstrint} recovers the affine structure of \cite[Example 1.1.3(2)]{GHS}. 

\begin{observe} \label{affstrdim2}
If $n=2$, then the singular locus $\Delta$ of the affine structure on $(B, \mathscr{P})$ is just the set of $v \in \mathscr{P}^{[0]}$, all the maximal cells are standard triangles, and for $\rho = \langle v_0, v_1 \rangle, \ \sigma = \langle v_0, v_1, v_2 \rangle, \ \sigma' = \langle v_0, v_1, v'_2 \rangle$, the definition \eqref{mapaffstr} of the embedding $\psi_{\rho, v_0}: \sigma \cup \sigma' \to \mathbb{R}^2$ reduces to a simple rule:
\begin{multline}
$$\psi_{\rho, v_0}(v_0) = (0, 0), \quad \psi_{\rho, v_0}(v_1) = (1, 0), \quad \psi_{\rho, v_0}(v_2) = (0, 1), \\ \psi_{\rho, v_0}(v'_2) = ( - D_{v_1} \cdot X_{\langle v_0, v_1 \rangle}, -1),
\end{multline}
and $\psi_{\rho, v_0}$ linear on each cone. 
\end{observe}

\begin{notation} \label{intnotation}
Consider subvarieties $Y_1, \ Y_2$ and $W$ of $\mathfrak{X}$ such that $Y_1 \subseteq W \subseteq \mathfrak{X}$ and $Y_2 \subseteq W \subseteq \mathfrak{X}$. In what follows, we shall denote by $(Y_1 \cdot Y_2)_W$ the intersection product of $Y_1$ and $Y_2$ computed in $W$.
\end{notation}

\begin{fact}  \label{affinestrghs}
In the case $n=2$, Construction \ref{affstrint} recovers the affine structure of \cite[Example 1.1.3(2)]{GHS}.
\end{fact}

\begin{proof}
The definitions of charts are the same for the interiors of the maximal cells. The definition of $\psi_{\rho, v}: \sigma \cup \sigma' \to \mathbb{R}^2$ in \cite[Example 1.1.3(2)]{GHS} is the same up to permuting the basis vectors and using $\left( X^2_{\langle v_0, v_1 \rangle}\right)_{D_{v_0}}$ instead of $D_{v_1} \cdot X_{\langle v_0, v_1 \rangle}$. So it is enough to show that $$\left( X^2_{\langle v_0, v_1 \rangle}\right)_{D_{v_0}} =  D_{v_1} \cdot X_{\langle v_0, v_1 \rangle}.$$
Indeed, consider the embedding $i: D_{v_0} \xhookrightarrow{} \mathfrak{X}$. We have 
\begin{multline}
\left( X^2_{\langle v_0, v_1 \rangle}\right)_{D_{v_0}} = i_{*} \left( X_{\langle v_0, v_1 \rangle} \cdot X_{\langle v_0, v_1 \rangle}\right) = \\ = i_{*}\left( i^*\left(D_{v_1}\right) \cdot  X_{\langle v_0, v_1 \rangle} \right) = D_{v_1} \cdot i_{*}( X_{\langle v_0, v_1 \rangle})  = D_{v_1} \cdot X_{\langle v_0, v_1 \rangle}.
\end{multline}
Here we used the projection formula (see, e.g. \cite[Appendix A.1]{hartshorne}) for the third equality. 
\end{proof}

\begin{remark} \label{triplepointformula}
The argument in the proof of Proposition \ref{affinestrghs} also implies that \eqref{intcentralfibre} is the natural generalization of the triple point formula $$\left( X_{\langle v_0, v_1 \rangle}\right)^2_{D_{v_0}} + \left( X_{\langle v_0, v_1 \rangle}\right)^2_{D_{v_{1}}} = -2$$ of \cite[p. 39]{P} in the case $n=2$.
\end{remark}

We can generalize the argument of Proposition \ref{affinestrghs} to show that the intersection numbers $D_{v_i} \cdot X_{\rho}$ in \eqref{constraintaffine} and \eqref{intcentralfibre} can be computed in strata of the central fibre of the desired dimension. 

\begin{fact} \label{intnumbersarbdimension}
Let $D_{v_k}$ be a divisor corresponding to a fixed $v_k$ in $$\rho = \langle v_{i} \ | \ 0 \leq i \leq n-1 \rangle$$ and let $$X_{J} := \bigcap_{j \in J} D_{v_j}, \text{ where } J \subseteq \left\{0, \dots, n-1\right\}, \ J \cap \left\{k\right\} = \varnothing$$ be a stratum of the central fibre. Then we have $$\left( \left( X_J \cap D_{v_k} \right) \cdot X_{\rho}\right)_{X_{J}} = D_{v_k} \cdot X_{\rho}.$$
\end{fact}

\begin{proof}
The result follows from the projection formula (see, e.g. \cite[Appendix A.1]{hartshorne}) in the same way as in the proof of Proposition \ref{affinestrghs} using the embedding $i: X_{J} \xhookrightarrow{} \mathfrak{X}.$
\end{proof}

\begin{sloppypar}
This allows us to rewrite the sums $\sum_{i=1}^{n-1} \left( D_{v_{i}}  \cdot X_{\rho} \right)e_i$ in \eqref{constraintaffine} and $\sum_{i=0}^{n-1} \left( D_{v_{i}}  \cdot X_{\rho} \right)$ in \eqref{intcentralfibre} in terms of intersection numbers in certain strata of the central fibre of chosen codimension. In particular, this, perhaps, gives a more natural way to view \eqref{intcentralfibre} as the generalization of the triple point formula (see Remark \ref{triplepointformula}).
\end{sloppypar}

\begin{cor} \label{compinsmallerstrata}
For any $1 \leq m \leq n-2$ and any $0 \leq i \leq n-1$ let $$J_i^m := \left\{i -1 \mod (n-1), \ \dots, \ i-m \mod (n-1)\right\}.$$ 
\begin{enumerate}
    \item For any $m$, we have $$\sum_{i=1}^{n-1} \left( D_{v_{i}}  \cdot X_{\rho} \right)e_i = \sum_{i=1}^{n-1} \left( \left( X_{J_i^m} \cap D_{v_i} \right) \cdot X_{\rho}\right)_{X_{J_i^m}} e_i$$ which is a sum of intersection numbers in strata of the central fibre of codimension $m$ in $\mathfrak{X}$.
    \item In particular, we have $$\sum_{i=1}^{n-1} \left( D_{v_{i}}  \cdot X_{\rho} \right)e_i = \sum_{i=1}^{n-1} \left( \left( D_{v_{i-1}} \cap D_{v_i} \right) \cdot X_{\rho}\right)_{D_{v_{i-1}}} e_i = \sum_{i=1}^{n-1} \left(  X_{\rho}^2\right)_{X_{\langle \rho \setminus v_{i} \rangle}} e_i.$$
    \item Similarly, we have
    $$\sum_{i=1}^{n-1} \left( D_{v_{i}}  \cdot X_{\rho} \right)e_i = \sum_{i=1}^{n-1} \left( \left( D_{v_{0}} \cap D_{v_i} \right) \cdot X_{\rho}\right)_{D_{v_{0}}} e_i$$
    \item Similar statements hold for $\sum_{i=0}^{n-1} \left( D_{v_{i}}  \cdot X_{\rho} \right)$.
\end{enumerate}
\end{cor}

\begin{proof}
For (1), apply Proposition \ref{intnumbersarbdimension} to every intersection product in the sum with $J = J_i^m$. Now (2) follows by setting $m = 1$ for the first equality and $m = n-2$ for the second. For (3), apply Proposition \ref{intnumbersarbdimension} to every intersection product in the sum with $J = \{0\}$. The statements for (4) follow by the same arguments.
\end{proof}

The following observation gives a different way to view Construction \ref{affstrint} that is better suited to generalizing the construction of the affine structure on $\left(\bar{B}, \bar{\mathscr{P}}\right)$ to arbitrary minimal log CY degenerations $\mathfrak{X} \to \mathcal{S}$ satisfying Assumption \ref{manifoldass}.

\begin{observe} \label{altaffinecons}
Corollary \ref{compinsmallerstrata}(3) gives an alternative way to view Construction \ref{affstrint} giving $B$ the structure of an affine manifold with singularities. For any $v \in \mathscr{P}^{[0]}$, the component $D_v$ of the central fibre with the log structure induced by the inclusion $D_v \xhookrightarrow{} \mathfrak{X}$ is a log CY variety $(D_v, \partial D_v)$ where the divisorial log structure is given by $$\partial D_v := \bigcup_{v \subseteq \rho \in \mathscr{P}^{[n-1]}} X_{\rho}$$
(indeed, $K_{\mathfrak{X}}+D \equiv 0$ implies that $K_{D_v} + \partial D_v \equiv 0$ by adjunction). In the case that $\mathfrak{X} \to \mathcal{S}$ is a degeneration of K3-s, $(D_v, \partial D_v)$ is a Looijenga pair (e.g. by the classification of surfaces).

Now, \cite[Section 1.3]{GS4} gives an affine structure on $\Sigma(D_v) \setminus \Delta_v$ (where $\Delta_v$ is the union of cells of $\Sigma(D_v)$ of codimension $2$) in the same way as in the description of the affine structure on $\mathbf{C}B$ in the proof of Proposition \ref{indofchoices}. This gives PL-embeddings 
\begin{equation} \label{locembed}
\psi_{\rho_v}: \sigma_v \cup \sigma'_v \to \mathbb{R}^n
\end{equation}
for any maximal cones $\sigma_v,  \sigma'_v \in \Sigma(D_v)$ with intersection a codimension $1$ cone $\rho_v \in \Sigma(D_v)$. 

It is clear from Construction \ref{tropicalize} of tropicalization and definition \eqref{dicims} of the dual intersection complex $(B, \mathscr{P})$ that there is a one-to-one correspondence between cones of $\Sigma(D_v)$ and cells of $\mathscr{P}$ containing $v$. Moreover, it respects the integral structure. Let $b(\sigma_v), \ b(\sigma'_v) \in \mathscr{P}^{\max}$ be the maximal cells corresponding to $\sigma_v, \ \sigma'_v \in \Sigma(D_v)$ and let $b(\rho_v) \in \mathscr{P}^{[n-1]}$ be the cell corresponding to $\rho_v \in \Sigma(D_v)$. Then there are natural inclusions $b(\sigma_v) \subseteq \sigma_v, \ b(\sigma'_v) \subseteq \sigma'_v$, and $b(\rho_v) \subseteq \rho_v$ compatible with inclusions of faces. We can now define an embedding 
\begin{equation} \label{indembed}
\psi_{b(\rho_v), v}: b(\sigma_v) \cup b(\sigma'_v) \to \mathbb{R}^n
\end{equation}
as the restriction of $\psi_{\rho_v}$ and doing this over all the $v \in \mathscr{P}^{[0]}$ and all the codimension $1$ cells $\rho_v$ of $\Sigma(D_v)$ defines a data of the refined description of charts of Construction \ref{affchartstool}$(2')$. 

It is straightforward to check using the equality of Corollary \ref{compinsmallerstrata}(3) that \eqref{indembed} defines the same embedding as \eqref{mapaffstr}. Therefore, an alternative way to define the structure of an affine manifold with singularities on $B$ using Construction \ref{affchartstool} is to define the structures of integral polyhedra on $\sigma \in \mathscr{P}^{\max}$ via the tropicalization $\Sigma(\mathfrak{X})$ as before and to define the affine structure on every $W_v \setminus \Delta, \ v \in \mathscr{P}^{[0]}$ to be the one induced by the affine structure on $\Sigma(D_v) \setminus \Delta_v$. 
\end{observe}

Finally, we note what happens if a component $D_{v}, \ v \in \mathscr{P}^{[0]}$ of $\mathfrak{X}_0$ is toric and observe that Construction \ref{affstrint} can be seen as a generalization of Construction \ref{affstrtoric} of an affine structure on the dual intersection complex $\left( \bar{B}, \bar{\mathscr{P}}\right)$ of a toric degeneration $\bar{\mathfrak{X}} \to \mathcal{S}$ with all $\sigma \in \bar{\mathscr{P}}^{\max}$ standard simplices.

\begin{remarks} \label{toriccomponent}
Suppose that $D_{v}, \ v \in \mathscr{P}^{[0]}$ is a toric variety, and its induced log stratification coincides with its toric stratification.
\begin{enumerate}
    \item Then the embedding $\psi_{\rho, v}$ of Construction \ref{affstrint} maps $\sigma \cup \sigma'$ to the union of the corresponding cones of the fan $\Sigma_{D_{v}}$ defining $D_{v}$. This follows from Observation \ref{altaffinecons}, Corollary \ref{compinsmallerstrata}(3), and standard toric geometry. This is actually the motivation behind Construction \ref{affstrint}.
    \item By \cite[Proposition 2.3]{AG} the affine structure on $W_v \setminus \Delta$ extends to an affine structure on the whole $W_v$. It follows from (1) that it is the same affine structure on $W_v$ as that of Construction \ref{affstrtoric}. A similar statement holds for toric strata $D_{\tau}, \tau \in \mathscr{P}$ of lower dimension.
    \item In particular, the recipe of Construction \ref{affstrint} (viewed via Observation \ref{altaffinecons}, i.e. replacing the right hand side of \eqref{constraintaffine} with the expression of Corollary \ref{compinsmallerstrata}(3)) recovers Construction \ref{affstrtoric} if applied to the dual intersection complex $\left( \bar{B}, \bar{\mathscr{P}}\right)$ of a toric degeneration $\bar{\mathfrak{X}} \to \mathcal{S}$ with all $\sigma \in \bar{\mathscr{P}}^{\max}$ standard simplices. However, note that there is no generalization of the triple point formula of Lemma \ref{intcentralfibrelemma} in this case, so the affine structure is defined on the complement of the discriminant locus $\bar{\Delta}$ of Construction \ref{affstrint}.
\end{enumerate}
\end{remarks}

\subsubsection{Affine structure in general} \label{affstrgeneral}

We now describe the structure of an affine manifold with singularities on $B$ in the general case of a minimal log CY degeneration $\mathfrak{X} \to \mathcal{S}$ satisfying Assumption \ref{manifoldass}. Since the components $D_i, \ 1 \leq i \leq m$ of the central fibre $\mathfrak{X}_0$ are not necessarily Cartier, we can't use intersection theory for the general description. Instead, we will argue as in Observation \ref{altaffinecons} using logarithmic geometry to define the embeddings \eqref{locembed}.

\begin{cons} \label{affstrintgen}
Let $\mathfrak{X} \to \mathcal{S}$ be a minimal log CY degeneration satisfying Assumption \ref{manifoldass} and let $(B, \mathscr{P})$ be the dual intersection complex of $\mathfrak{X} \to \mathcal{S}$. We let the discriminant locus $\Delta$ be the union of codimension $2$ cells \eqref{discriminantlocusims} of $\mathscr{P}$ as before and want to put the affine structure on $B_0 := B \setminus \Delta$ using the framework of Construction \ref{affchartstool}. As before, we let the structures of integral polyhedra on $\sigma \in \mathscr{P}^{\max}$ be given by their inclusions into cones of $\Sigma(\mathfrak{X})$. By Construction \ref{affchartstool}(2), it suffices to give structures of integral affine manifolds on every $\left\{W_v \setminus \Delta \ | \ v \in \mathscr{P}^{[0]} \right\}$
for $W_v$ as in \eqref{openstar}, compatible with the integral structures on the maximal cells. 

As in Observation \ref{altaffinecons}, for any $v \in \mathscr{P}^{[0]}$, we endow $D_v$ with the log structure induced by the inclusion $D_v \xhookrightarrow{} \mathfrak{X}$, making it a log CY variety. Arguing as in Observation \ref{altaffinecons}, it is enough to give, for every $v \in \mathscr{P}^{[0]}$, an affine structure on $\Sigma(D_v) \setminus \Delta_v$ respecting the integral structure. We let the affine structure on the interiors of the maximal cones be given via the images of those cones in $\Sigma(D_v)$ as before.\footnote{By Construction \ref{tropicalize} of tropicalization and definition \eqref{dicims} of the dual intersection complex $(B, \mathscr{P})$, this is the only definition compatible with the structures of integral polyhedra on $\sigma \in \mathscr{P}^{\max}$.} Now, it is enough to give PL-embeddings $\psi_{\rho_v}: \sigma_v \cup \sigma'_v \to \mathbb{R}^n$ for any maximal cones $\sigma_v, \ \sigma'_v \in \Sigma(D_v)$ with intersection a codimension $1$ cone $\rho_v \in \Sigma(D_v)$. Following the logic of Section \ref{affstrsnc}, these embeddings should be natural from the logarithmic point of view in the same sense as in Remark \ref{logremark}.

Since $\mathfrak{X} \to \mathcal{S}$ and $(B, \mathscr{P})$ satisfy \cite[Assumptions 1.2]{GS4} by Proposition \ref{assofcanonical}, it follows that $D_v$ and $\Sigma(D_v)$ satisfy \cite[Assumptions 1.1]{GS4}. 
Let $\rho_v \in \Sigma(D_v)$ be a codimension $1$ cone (of the form $\rho_v = \sigma_v \cap \sigma'_v$ for maximal cones $\sigma_v, \sigma'_v \in \Sigma(D_v)$). Assumption \ref{manifoldass} implies that for any stratum $D_I$ of \eqref{DI} with $X_{\rho_v} \cap D_I$ non-empty, $X_{\rho_v} \cap D_I$ is either $X_{\rho_v}$ itself or one of the zero-dimensional strata $X_{\sigma_v}, X_{\sigma'_v} \subseteq X_{\rho_v}$. As in the proof of \cite[Proposition 1.3]{GS4}, the fact that $K_{D_v} + \partial D_v \equiv 0$ implies, by repeatedly applying adjunction, that $K_{X_{\rho_v}} + \partial X_{\rho_v} \equiv K_{X_{\rho_v}} +  X_{\sigma_v} + X_{\sigma'_v} \equiv 0$ (using the notations of Observation \ref{altaffinecons}). In particular, $(X_{\rho_v}, \partial X_{\rho_v})$ is itself a log CY variety. Since $X_{\rho_v}$ is one-dimensional, it is a compact non-singular curve. But $K_{X_{\rho_v}} + \partial X_{\rho_v} \equiv 0$ so the only possibility is that $X_{\rho_v} \cong \mathbb{P}^1$ and $X_{\sigma_v}, X_{\sigma'_v}$ are the torus-fixed points of $X_{\rho_v}$.
In particular, $X_{\rho_v}$ is a toric variety, and the induced log stratification of $X_{\rho_v}$ is the same as the toric stratification. Under these assumptions, \cite[Theorem 4.1]{W} uses charts of the log structure $\mathcal{M}_{X_{\rho_v}}$ to provide a toric variety $X_{\Sigma_{\rho_v}}$ and a cone $\tilde{\sigma}_{\rho_v} \in \Sigma_{\rho_v}$ such that $X_{\rho_v}$ (with the pullback log structure induced from $D_v$ by inclusion) is isomorphic as a log scheme to the toric stratum $X_{\Sigma_{\rho_v}, \tilde{\sigma}_{\rho_v}}$ (with the log structure induced from $X_{\Sigma_{\rho_v}}$).\footnote{\cite[Theorem 4.1]{W} also requires that $\overline{\mathcal{M}}_{X_{\rho_v}}$ is globally generated. Since $X_{\rho_v} \cong \mathbb{P}^1$, $\Sigma(X_{\rho_v})$ admits a PL-embedding into $\mathbb{R}$ so $\overline{\mathcal{M}}_{X_{\rho_v}}$ is globally generated by \cite[Proposition 3.12]{ACGS2}.} 

More explicitly, in the case of $X_{\rho_v}$ the construction in the proof of \cite[Theorem 4.1]{W} produces a fan $\Sigma_{X_{\rho_v}}$ with $|\Sigma_{X_{\rho_v}}| \subseteq \mathbb{R}^n$ consisting of two cones $\tilde{\sigma}_v, \tilde{\sigma}'_v$ of dimension $n$ and $\tilde{\sigma}_{\rho_v} = \tilde{\sigma}_v \cap \tilde{\sigma}'_v$. The $n$-dimensional cones correspond to the torus-fixed points of $X_{\rho} \cong \mathbb{P}^1$, we let $\tilde{\sigma}_v$ (resp. $\tilde{\sigma}'_v$) be the cone corresponding to the point $X_{\sigma_v}$ (resp. $X_{\sigma'_v}$). The construction of $\Sigma_{X_{\rho_v}}$ is performed via charts of the log structure $\mathcal{M}_{X_{\rho_v}}$ which ensures that the integral structure on $\tilde{\sigma}_v$ (resp. $\tilde{\sigma}'_v$) is the same as on $\sigma_v$ (resp. $\sigma'_v$). We define
\begin{equation} \label{embeddinglocalgen}
\psi_{\rho_v}: \sigma_v \cup \sigma'_v \to \mathbb{R}^n
\end{equation}
to be the unique PL-embedding identifying $\sigma_v$ with $\tilde{\sigma_v}$, $\sigma'_v$ with $\tilde{\sigma'_v}$ and respecting the integral structures. This concludes the construction of the affine structure on $\Sigma(D_v) \setminus \Delta_v$ and defining such structures for all $v \in \mathscr{P}^{[0]}$ defines an affine structure on $B_0 = B \setminus \Delta$. 
\end{cons}

\begin{remarks} \label{genaffstrremarks}
One needs to check a few things regarding Construction \ref{affstrintgen}.
\begin{enumerate}
    \item One needs to show that the structure of an affine manifold with singularities on $B$ defined via Construction \ref{affchartstool} extends to the complement of the smaller discriminant locus $\Delta$ of \eqref{discriminantlocusims}. As in the case of Construction \ref{affstrint}, one can show that the affine structure on $B$ is the one induced from $\mathbf{C}B$ via the natural inclusion $B \xhookrightarrow{} \mathbf{C}B$. Here the affine structure on $\mathbf{C}B \cong \Sigma(\mathfrak{X})$ is defined in the same way as the affine structures on $\Sigma(D_v), \ v \in \mathscr{P}^{[0]}$ in Construction \ref{affstrintgen}. Compatibility under the inclusion is a consequence of the construction in the proof of \cite[Theorem 4.1]{W} and the fact that the $D_v, \ v \in \mathscr{P}^{[0]}$ carry the pullback log structures. 
    \item One also needs to check that in the case that $\mathfrak{X} \to \mathcal{S}$ has a reduced central fibre, and $D$ is simple normal crossings, Construction \ref{affstrintgen} reduces to Construction \ref{affstrint}. It suffices to show that the PL-embedding \eqref{embeddinglocalgen} recovers \eqref{locembed} under these assumptions. This is a consequence of the construction in the proof of \cite[Theorem 4.1]{W} and the explicit presentation of $\mathcal{M}_{X_{\rho_{v}}}$ in the proof of \cite[Lemma 1.10]{GS4}.
    \item Finally, in the case that $D_{v}, \ v \in \mathscr{P}^{[0]}$ is a toric variety, and its induced log stratification coincides with its toric stratification, the analogues of Remarks \ref{toriccomponent} hold. This follows from standard toric geometry. In particular, the affine structure extends across $v$ to the whole $W_v$ in this case, and applying Construction \ref{affstrintgen} to the dual intersection complex $\left( \bar{B}, \bar{\mathscr{P}} \right)$ of a toric degeneration $\bar{\mathfrak{X}} \to \mathcal{S}$ recovers Construction \ref{affstrtoric} of the affine structure with discriminant locus $\bar{\Delta}$.
\end{enumerate}
\end{remarks}

\subsubsection{The group \texorpdfstring{$A_1(\mathfrak{X}_0, \mathbb{Z})$}{A1}  and the monoid \texorpdfstring{$P$}{P}}
Let $A_1(\mathfrak{X}_0, \mathbb{Z})$ be the relevant group of curve classes. That may be the group of $1$-cycles on $\mathfrak{X}_0$ or the group of relative $1$-cycles for $\mathfrak{X} \to \mathcal{S}$  modulo algebraic or numerical equivalence (we sometimes write $A_1(\mathfrak{X}_0, \mathbb{Z})_{\alg}$, $ A_1(\mathfrak{X}_0, \mathbb{Z})_{\numer}$, $ A_1(\mathfrak{X}/\mathcal{S}, \mathbb{Z})_{\alg}$, and $A_1(\mathfrak{X}/\mathcal{S}, \mathbb{Z})_{\numer}$ respectively). Out of the possible choices for $A_1(\mathfrak{X}_0, \mathbb{Z})$, the group $ A_1(\mathfrak{X}/\mathcal{S}, \mathbb{Z})_{\numer}$ fits best with the extension of the intrinsic mirror results of Section \ref{extresults} and $A_1(\mathfrak{X}_0, \mathbb{Z})_{\numer}$ fits best with the generalizations of Chapter \ref{gslocus}. We have canonical splittings as follows. 

\begin{cons} \label{splittings}
We have a commutative diagram with exact rows and columns:
\begin{equation} 
    \begin{tikzcd}
    & & 0 \arrow{d} & 0 \arrow{d} & \\
    & & \Rel(\mathfrak{X}/\mathcal{S})_{\alg} \arrow{d} & \Rel(\mathfrak{X}/\mathcal{S})_{\numer} \arrow{d} & \\
    0 \arrow{r} & \Griff(\mathfrak{X}/\mathcal{S}) \arrow{r} & A_1(\mathfrak{X}/\mathcal{S}, \mathbb{Z})_{\alg} \arrow{r}{a} \arrow{d}{c} & A_1(\mathfrak{X}/\mathcal{S}, \mathbb{Z})_{\numer} \arrow{r} \arrow{d}{d} & 0 \\
    0 \arrow{r} & \Griff(\mathfrak{X}_0) \arrow{r} & A_1(\mathfrak{X}_0, \mathbb{Z})_{\alg} \arrow{r}{b} \arrow{d} & A_1(\mathfrak{X}_0, \mathbb{Z})_{\numer} \arrow{r} \arrow{d} & 0 \\
    & & 0 & 0 & 
    \end{tikzcd}
\end{equation}
Here $\Griff(\mathfrak{X}_0) := \frac{Z_1(\mathfrak{X}_0)_{\numer}}{Z_1(\mathfrak{X}_0)_{\alg}}$ and $\Griff(\mathfrak{X}/\mathcal{S}) := \frac{Z_1(\mathfrak{X}/\mathcal{S})_{\numer}}{Z_1(\mathfrak{X}/\mathcal{S})_{\alg}}$ are the Griffiths groups of $1$-cycles that are numerically but not algebraically equivalent to $0$.\footnote{A Griffiths group of a scheme $X$ is more commonly defined as $\Griff(X) :=\frac{Z_1(X)_{\hom}}{Z_1(X)_{\alg}}$ where $Z_1(X)_{\hom}$ is the group of $1$-cycles homologically equivalent to $0$ for a particular choice of a Weil cohomology theory. Unless we assume standard conjecture D, our groups may be larger.} With this definition, the exactness of the rows is clear. Since $\mathcal{S}$ is the spectrum of a discrete valuation $\kk$-algebra, every relative $1$-cycle for $\mathfrak{X} \to \mathcal{S}$ is defined over $\mathfrak{X}_0$, which gives rise to the surjections $c$ and $d$. Moreover, we have $Z_1(\mathfrak{X}_0)_{\alg} \subseteq Z_1(\mathfrak{X}/\mathcal{S})_{\alg}$ and $Z_1(\mathfrak{X}_0)_{\numer} \subseteq Z_1(\mathfrak{X}/\mathcal{S})_{\numer}$ so exactness of the columns follows once we set $\Rel(\mathfrak{X}/\mathcal{S})_{\alg} := \frac{Z_1(\mathfrak{X}_0)_{\alg}}{Z_1(\mathfrak{X}/\mathcal{S})_{\alg}}$ and $\Rel(\mathfrak{X}/\mathcal{S})_{\numer} := \frac{Z_1(\mathfrak{X}_0)_{\numer}}{Z_1(\mathfrak{X}/\mathcal{S})_{\numer}}$.

Now, note that $A_1(\mathfrak{X}_0, \mathbb{Z})_{\numer}$ and $A_1(\mathfrak{X}/\mathcal{S}, \mathbb{Z})_{\numer}$ are finitely generated free abelian groups. Therefore, the maps $a$, $b$, and $d$ admit splittings. Fixing such splittings, we have 
\begin{equation} \label{standardsplittings}
\begin{aligned}
A_1(\mathfrak{X}_0, \mathbb{Z})_{\alg} &\cong A_1(\mathfrak{X}_0, \mathbb{Z})_{\numer} \oplus \Griff(\mathfrak{X}_0) \\
A_1(\mathfrak{X}/\mathcal{S}, \mathbb{Z})_{\numer} &\cong A_1(\mathfrak{X}_0, \mathbb{Z})_{\numer} \oplus \Rel(\mathfrak{X}/\mathcal{S})_{\numer} \\
A_1(\mathfrak{X}/\mathcal{S}, \mathbb{Z})_{\alg} &\cong A_1(\mathfrak{X}_0, \mathbb{Z})_{\numer} \oplus \Rel(\mathfrak{X}/\mathcal{S})_{\numer} \oplus \Griff(\mathfrak{X}/\mathcal{S}) 
\end{aligned}
\end{equation}
Generally, we will write $A_1(\mathfrak{X}_0, \mathbb{Z}) = A_1(\mathfrak{X}_0, \mathbb{Z})_{\numer} \oplus G$ for the splitting and we will write $A_1(\mathfrak{X}_0, \mathbb{Z})_{\numer} \subseteq A_1(\mathfrak{X}_0, \mathbb{Z})$ for the inclusion as the first factor.
\end{cons}

We further assume that $A_1(\mathfrak{X}_0, \mathbb{Z})$ is a finitely generated abelian group. This is always the case for $ A_1(\mathfrak{X}_0, \mathbb{Z})_{\numer}$ and $ A_1(\mathfrak{X}/\mathcal{S}, \mathbb{Z})_{\numer}$. For $A_1(\mathfrak{X}_0, \mathbb{Z})_{\alg}$ and $ A_1(\mathfrak{X}/\mathcal{S}, \mathbb{Z})_{\alg}$ this is equivalent to requiring that $\Griff(\mathfrak{X}_0)$ and $\Griff(\mathfrak{X}/\mathcal{S})$ respectively are finitely generated. In many examples, this is not the case. For instance, the Griffiths group of a general Calabi-Yau threefold is not finitely generated even after factoring by the torsion subgroup (see \cite{cl} for the case of a quintic threefold, and \cite{V} for the general case). Note that $\Griff(\mathfrak{X}_0) = 0$ if $\mathfrak{X} \to \mathcal{S}$ is a degeneration of K3-s.\footnote{Since algebraic and numerical equivalence coincide for divisors, see \cite{L} (in French).} 

We allow various choices for $A_1(\mathfrak{X}_0, \mathbb{Z})$ to offer more flexibility in the construction. The reader may wish to just think of $A_1(\mathfrak{X}_0, \mathbb{Z}) = A_1(\mathfrak{X}_0, \mathbb{Z})_{\numer}$. Note that for any choice of $A_1(\mathfrak{X}_0, \mathbb{Z})$ we have a natural pairing $$A_1(\mathfrak{X}_0, \mathbb{Z}) \times \Pic(\mathfrak{X}) \to \mathbb{Z}, \ (\beta, \mathcal{L}) \mapsto \deg_{\mathcal{L}}(\beta)$$ which allows us to compute intersection numbers of curve classes with Cartier divisors. 

\begin{remark}
One may consider more general groups $H_2(\mathfrak{X}/\mathcal{S})$ of degree data than the choices of $A_1(\mathfrak{X}_0, \mathbb{Z})$ we have described. See \cite[Basic Setup 1.6]{GS3} for a discussion on this. In particular, one can use other equivalence relations (as long as they are at least as coarse as algebraic equivalence). One can also use $H_2(\mathfrak{X}, \mathbb{Z})$ if working over $\kk = \mathbb{C}$. For our purposes, one can take $H_2(\mathfrak{X}/\mathcal{S})$ to be any finitely generated abelian group of curve classes that admits a splitting $H_2(\mathfrak{X}/\mathcal{S}) \cong A_1(\mathfrak{X}_0, \mathbb{Z})_{\numer} \oplus G$.
\end{remark}

Let $NE(\mathfrak{X}_0) \subseteq A_1(\mathfrak{X}_0, \mathbb{Z})$ be the submonoid generated by the effective curve classes (again, we specify $NE(\mathfrak{X}_0)_{\alg}$, $NE(\mathfrak{X}_0)_{\numer}$, $NE(\mathfrak{X}/\mathcal{S})_{\alg}$, $NE(\mathfrak{X}/\mathcal{S})_{\numer}$ when the distinction is important). 

\begin{remark}
In the situation that we have a map $\pi: \mathfrak{X} \to \bar{\mathfrak{X}}$ resolving a toric degeneration $\bar{\mathfrak{X}} \to \mathcal{S}$ to a minimal log CY degeneration $\mathfrak{X} \to \mathcal{S}$ and such that the cone $NE(\pi) \subseteq NE(\mathfrak{X}_0)$ of curves contracted by $\pi$ is finitely generated (our resolutions will always satisfy this requirement), the group $A_1(\mathfrak{X}_0, \mathbb{Z})$ is finitely generated if and only if the corresponding $A_1(\bar{\mathfrak{X}}_0, \mathbb{Z})$ is finitely generated.

Therefore, it is enough to assume that $A_1(\bar{\mathfrak{X}}_0, \mathbb{Z})$ is finitely generated for the toric degeneration $\bar{\mathfrak{X}} \to \mathcal{S}$ we work with. The groups $ A_1(\bar{\mathfrak{X}}_0, \mathbb{Z})_{\numer}$ and $ A_1(\bar{\mathfrak{X}}/\mathcal{S}, \mathbb{Z})_{\numer}$ are always finitely generated. The group $A_1(\bar{\mathfrak{X}}_0, \mathbb{Z})_{\alg}$ is also finitely generated in this case since there is a natural surjective map $$\bigoplus_i A_1(\bar{D}_i, \mathbb{Z})_{\alg} \to A_1(\bar{\mathfrak{X}}_0, \mathbb{Z})_{\alg}$$ and the groups $A_1(\bar{D}_i, \mathbb{Z})_{\alg}$ are finitely generated because the irreducible components $\bar{D}_i, \ 1 \leq i \leq \bar{m}$ of $\bar{\mathfrak{X}}_0$ are toric. We refrain from making comments on finite generation of $A_1(\mathfrak{X}_0, \mathbb{Z})$ from now on.
\end{remark}

The monoid $NE(\mathfrak{X}_0) \subseteq A_1(\mathfrak{X}_0, \mathbb{Z})$ is usually not finitely generated (see \cite[Example 1.23]{KM} for some examples) so we choose a larger finitely generated monoid. 

\begin{defn} \label{monoidass}
We choose a \emph{base monoid} $P$ such that:
\begin{enumerate}
\item $NE(\mathfrak{X}_0) \subseteq P \subseteq A_1(\mathfrak{X}_0, \mathbb{Z})$.
\item $P$ is finitely generated and saturated.
\item The group $P^{\times}$ of the invertible elements of $P$ coincides with the torsion part of $A_1(\mathfrak{X}_0, \mathbb{Z})$.
\end{enumerate}
\end{defn}

Note that if $f:C \to \mathfrak{X}/\mathcal{S}$ is a stable map with $f_{*}[C] \in P^{\times}$, then $f_{*}[C] = 0$.

\begin{remarks} Some remarks on this setup are in order.
\begin{enumerate}
    \item A monoid $P$ satisfying Definition \ref{monoidass} is not necessarily toric since $P^{\gp} = A_1(\mathfrak{X}_0, \mathbb{Z})$ is not always torsion-free. Therefore, Definition \ref{monoidass} allows slightly more general monoids than the setup of Section \ref{mpafunc}. Since $A_1(\mathfrak{X}_0, \mathbb{Z})_{\numer}$ and $A_1(\mathfrak{X}/\mathcal{S}, \mathbb{Z})_{\numer}$ are always free, these choices guarantee that $P$ is a toric and \emph{sharp} (i.e. $P^{\times} = 0$) monoid.
    \item The construction of an extended intrinsic mirror in Section \ref{extresults} allows to relax Condition (3) of Definition \ref{monoidass} under some conditions (see Remark \ref{extmonoidchange}).
\end{enumerate}
\end{remarks}

We prove a quick lemma that shows how to produce monoids satisfying Definition \ref{monoidass}.

\begin{lemma} \label{mirrorexttrivial}
Suppose that $P'$ is a monoid satisfying conditions (1) and (2) of Definition \ref{monoidass} and such that $A_1(\mathfrak{X}_0, \mathbb{Z})_{\tors} \subseteq (P')^{\times}$. Then there exists a monoid $P$ with $NE(\mathfrak{X}_0) \subseteq P \subseteq P' \subseteq A_1(\mathfrak{X}_0, \mathbb{Z})$ and satisfying Definition \ref{monoidass}. In particular, there exist monoids $P$ satisfying Definition \ref{monoidass}.
\end{lemma}

\begin{proof}
Since $\mathfrak{X}_0$ is projective, $NE(\mathfrak{X}_0)_{\mathbb{R}}  \subseteq A_1(\mathfrak{X}_0, \mathbb{R})$ is a strictly convex cone of maximal dimension $k$. Let $NE(\mathfrak{X}_0)\widecheck{_{\mathbb{R}}} \subseteq A_1(\mathfrak{X}_0, \mathbb{R})\,\widecheck{}$ be the dual cone and let $e_0, \dots, e_k \in \Hom(NE(\mathfrak{X}_0), \mathbb{Z})$ be linearly independent as elements of $A_1(\mathfrak{X}_0, \mathbb{R})\,\widecheck{}$. Then we have $$\Hom(NE(\mathfrak{X}_0), \mathbb{Z}) \supseteq \langle e_0, \dots, e_k \rangle_{\mathbb{Z}}^{\sat}$$ (we use the notation $M^{\sat}$ for the saturation of a monoid $M$), so 
\begin{multline}
NE(\mathfrak{X}_0) \subseteq \Hom(\Hom(NE(\mathfrak{X}_0), \mathbb{Z}), \mathbb{Z}) \oplus A_1(\mathfrak{X}_0, \mathbb{Z})_{\tors} \subseteq \\ \subseteq \left(\langle e_0, \dots, e_k \rangle_{\mathbb{Z}}^{\sat} \right)\widecheck{} \  \oplus A_1(\mathfrak{X}_0, \mathbb{Z})_{\tors} 
\end{multline}
where the first inclusion is since $NE(\mathfrak{X}_0)$ is integral. It follows that $$NE(\mathfrak{X}_0) \subseteq P:=\left( \left(\langle e_0, \dots, e_k \rangle_{\mathbb{Z}}^{\sat} \right)\widecheck{} \cap P'/P'_{\tors} \right)\  \oplus A_1(\mathfrak{X}_0, \mathbb{Z})_{\tors} $$ (where the intersection is computed in $A_1(\mathfrak{X}_0, \mathbb{Z})/A_1(\mathfrak{X}_0, \mathbb{Z})_{\tors}$). It is easy to check that $P$ satisfies all the necessary conditions. The last claim follows by taking $P' = A_1(\mathfrak{X}_0, \mathbb{Z})$. 
\end{proof}

\subsubsection{MPA function and the initial slab functions} \label{mpaintrinsic}
We work over $A=\kk$. We let $I_0 = \mathfrak{m} := P \setminus P^{\times}$ be the maximal ideal and define a $P_{\mathbb{R}}^{\gp} = A_1(\mathfrak{X}_0, \mathbb{R})$-valued convex MPA function $\varphi$ on $B$ via its kinks\footnote{If $P$ has torsion, then there are a few subtleties. The kinks as defined in \eqref{kinks} only lie in $P^{\gp}/P_{\tors}^{\gp}$. However, if $P$ has torsion, it is more natural to view an MPA function as a collection of kinks in $P^{\gp} = A_1(\mathfrak{X}_0, \mathbb{Z})$. This only affects the definition of the sheaf $\mathcal{P}$ of \eqref{monomialexactseq}. See \cite[Footnotes 5 and 7]{GS4} for more details.} by setting $$\kappa_{\underline{\rho}} := [X_{\rho}] \in A_1(\mathfrak{X}_0, \mathbb{Z}) = P^{\gp}$$ for every $\underline{\rho} \subseteq \tilde{\mathscr{P}}^{[n-1]}$ with $\underline{\rho} \subseteq \rho \subseteq \mathscr{P}^{[n-1]}.$\footnote{By \cite[Remark 1.17]{GS4}, this agrees with the restriction of the MPA function defined in \cite{GS4} on $\mathbf{C}B$.} Note that for any two slabs $\underline{\rho}, \underline{\rho}' \subseteq \rho \in \mathscr{P}^{[n-1]}$ we have $\kappa_{\underline{\rho}} = \kappa_{\underline{\rho}'} = : \kappa_{\rho}$ as we required in Section \ref{mpafunc}. The curve class makes sense since Assumption \ref{manifoldass} implies that $X_{\rho}$ is proper. Further, the curve class is not invertible, so we actually have  $\kappa_{\rho} \in P \setminus P^{\times}$ and $\varphi$ is strictly convex.

We set the initial slab functions to be trivial: $$\left\{f_{\underline{\rho}} = 1 \ \Big| \ f_{\underline{\rho}} \in \kk[\Lambda_{\rho}], \  \underline{\rho} \in \tilde{\mathscr{P}}^{[n-1]} \right\}.$$

\subsubsection{The mirror via three-pointed punctured maps} \label{3pointmirror}
In the paper \cite{GS3}, the mirror to $\mathfrak{X} \overset{g}{\to} \mathcal{S}$ is constructed without using scattering diagrams by defining the products in the theta function ring for $\mathbf{C}B$ directly, using the theory of \emph{punctured log Gromov-Witten invariants} developed in \cite{ACGS2}. 

\begin{notation} \label{logpairing}
For any element $s \in \Gamma\left(\mathfrak{X}, \overline{\mathcal{M}}_{\mathfrak{X}}^{\gp}\right)$, $\eta \in \mathfrak{X}$ a generic point of a logarithmic stratum, and $u \in \Hom\left(\overline{\mathcal{M}}_{\mathfrak{X}, \eta}, \mathbb{Z}\right)$ an integral tangent vector to the cone $\sigma_{\eta} \in \Sigma(\mathfrak{X})$, we have a natural evaluation of $s$ on $u$ that we denote $\langle s, u \rangle$ (this defines a bilinear pairing). In particular, by Construction \ref{tropicalize} of tropicalization we have a well-defined $\langle s, p \rangle$ for any $s \in \Gamma\left(\mathfrak{X}, \overline{\mathcal{M}}_{\mathfrak{X}}^{\gp}\right)$ and any $p \in \Sigma(\mathfrak{X})(\mathbb{Z}) = \mathbf{C}B(\mathbb{Z})$. If $s \in \Gamma\left(\mathfrak{X}, \overline{\mathcal{M}}_{\mathfrak{X}}\right) \subseteq \Gamma\left(\mathfrak{X}, \overline{\mathcal{M}}_{\mathfrak{X}}^{\gp}\right)$, then $\langle s, p \rangle \geq 0$ for any $p \in \mathbf{C}B(\mathbb{Z})$.

In the case that $D$ is a simple normal crossings divisor, we have a natural isomorphism $\Gamma\left( \mathfrak{X}, \overline{\mathcal{M}}_{\mathfrak{X}}\right) \cong \Div_{D}(\mathfrak{X})$ (see Construction \ref{tropicalizesnc}). Under this isomorphism, the above pairing corresponds to the canonical pairing between $\Div_{D}(\mathfrak{X})$ and $\Div_{D}(\mathfrak{X})^*$.
\end{notation}

Let $P$ be as in Definition \ref{monoidass} and let $\widehat{\kk [P]}$ be the completion of $\kk[P]$ with respect to the maximal ideal $\mathfrak{m} := P \setminus P^{\times}$. We define a graded ring (that we call the \emph{theta function ring}): 
\begin{equation} \label{ringofdef}
\widehat{R} = \bigoplus_{p \in \mathbf{C}B(\mathbb{Z})} \widehat{\kk [P]} \vartheta_p.
\end{equation}
Here we regard the theta functions $\vartheta_p, \ p \in \mathbf{C}B$ as just the generators of $\widehat{R}$. 

Let $\rho \in \Gamma \left( \mathfrak{X}, \overline{\mathcal{M}}_{\mathfrak{X}}\right)$ be the pullback of $1 \in \Gamma\left(\mathcal{S},\overline{\mathcal{M}}_{\mathcal{S}}\right)$.
The grading on \eqref{ringofdef} comes from setting $\deg \vartheta_p:= \langle \rho, p \rangle$. Note that $\deg \vartheta_p$ for $p \in \mathbf{C} \sigma$ is equal to the $d$ such that $p \in B\left(\frac{1}{d} \mathbb{Z} \right)$ (as defined in \eqref{ratpoints}). It remains to define a product rule on $\hat{R}$, i.e. to specify the structure constants $\alpha_{p_1 p_2 r} \in \widehat{ \kk [P]}$ in the product of the theta functions
\begin{equation} \label{prodtheta}
\vartheta_{p_1} \cdot \vartheta_{p_2} = \sum_{r \in \mathbf{C}B(\mathbb{Z})} \alpha_{p_1 p_2 r} \vartheta_r,
\end{equation}
such that only a finite number of $\alpha_{p_1 p_2 r} \in \widehat{ \kk [P]}$ are non-zero and $\deg \vartheta_r = \deg \vartheta_{p_1} + \deg \vartheta_{p_2}$ for any $\vartheta_{r}$ with $\alpha_{p_1 p_2 r} \neq 0$. In fact, it is enough to require that the product rule respects the grading since $B$ is compact, so $B \left( \frac{1}{d} \mathbb{Z} \right)$ is finite. 

We can expand
\begin{equation} \label{sumalpha}
\alpha_{p_1 p_2 r}= \sum_{\beta \in P} N_{p_1 p_2 r}^{\beta} t^{\beta}
\end{equation}
where the sum is over all classes $\beta \in P$ of stable maps to $\mathfrak{X}/\mathcal{S}$ and $N_{p_1 p_2 r} \in \kk$ are the structure constants we need to define. Note that the sum is finite modulo any $\mathfrak{m}^k, \ k \geq 1$ since $P \setminus \mathfrak{m}^k$ is finite.  So we can view $\alpha_{p_1 p_2 r} \in \widehat{\kk [P]}$ by taking the inverse limit over $\widehat{\kk [P]}/\mathfrak{m}^k, \ k \geq 1$.

\begin{remark} \label{fewerideals}
In \cite{GS3}, the inverse limit is taken over all monoid ideals $I \subseteq P$ with $P \setminus I$ is finite. However, our description is equivalent. Indeed, it is easy to check that $P \setminus I$ is finite if and only if $\sqrt{I} = \mathfrak{m}$. So it is enough to show that any $I$ with $\sqrt{I} = \mathfrak{m}$ is contained in $\mathfrak{m}^k$ for some $k \geq 1$. But since $\kk [P]$ is Noetherian, any ideal is contained in some power of its radical. 
\end{remark}

The structure constants $N_{p_1 p_2 r}^{\beta} \in \mathbb{Q}$ are rationally defined as certain \emph{punctured log Gromov-Witten invariants}, as introduced in \cite{ACGS2}. A point $p \in \mathbf{C}B$ can be seen as imposing tangency conditions on stable or punctured log curves, we refer to \cite[Section 2.3]{ACGS}. In the case that $D$ is a simple normal crossings divisor, the structure constants $N_{p_1 p_2 r}^{\beta} \in \mathbb{Q}$ are virtual counts of genus zero three-pointed punctured stable maps of class $\beta$ with marked points $x_1, x_2, x_{\out}$, having contact order $\langle D_i, p_1 \rangle$ with $D_i$ at $x_1$, contact order $\langle D_i, p_2 \rangle$ with $D_i$ at $x_2$, and contact order $- \langle D_i, r \rangle$ with $D_i$ at $x_{\out}$ (here we make use of Notation \ref{logpairing}), with a certain logarithmically imposed point constraint.  

Note that $- \langle D_i, r \rangle$ is negative unless $r=0$, but one can make sense of negative contact orders using punctured theory. Negative contact orders correspond to \emph{punctures} of a punctured curve, and if the contact order with $D_i$ at the puncture $x_{\out}$ is negative, the component of the curve containing $x_{\out}$ is contained in $D_i$.  

The logarithmic point constraint at the schematic level means that we only count punctured maps that map $x_{\out}$ to an arbitrary fixed point $z$ in the interior of the stratum $X_{\sigma}$ such that the corresponding $\mathbf{C}\sigma$ is the smallest cell containing $r$. 

Associated to a choice of data $\beta, p_1, p_2, r$ and $z$, one can define a moduli space $\mathscr{M}(\mathfrak{X}, \beta, z)$ that is a Deligne-Mumford stack of virtual dimension $0$\footnote{This is only true in the case that $\mathfrak{X} \to \mathcal{S}$ is minimal log CY. More generally, the virtual dimension of this stack is $(K_{\mathfrak{X}}+D) \cdot \beta$, see \cite[Construction 1.8]{GS3}.}, and set $$N_{p_1 p_2 r}^{\beta} := \deg[\mathscr{M}(\mathfrak{X}, \beta, z)]^{\virt}.$$ The construction of $\mathscr{M}(\mathfrak{X}, \beta, z)$ is quite technical, and we do not review it here. We refer to \cite[Section 3]{GS3} for more details on defining the invariants $N_{p_1 p_2 r}^{\beta}$. 

We now explain that the product rule \eqref{prodtheta} for theta functions respects the grading on the theta function ring $\hat{R}$. Consider the natural exact sequence 
$$0 \longrightarrow \mathcal{O}_{\mathfrak{X}}^{\times} \longrightarrow \mathcal{M}_{\mathfrak{X}}^{\gp}  \longrightarrow  \overline{\mathcal{M}}_{\mathfrak{X}}^{\gp} \longrightarrow 0.$$
The associated long exact sequence in cohomology gives a connecting homomorphism 
\begin{equation} \label{connectinghom}
\tau: \Gamma\left(\mathfrak{X}, \overline{\mathcal{M}}_{\mathfrak{X}}^{\gp}\right) \to H^1(\mathfrak{X}, \mathcal{O}_{\mathfrak{X}}^{\times}) = \Pic(\mathfrak{X}).
\end{equation}
For $s \in \Gamma\left(\mathfrak{X}, \overline{\mathcal{M}}_{\mathfrak{X}}^{\gp}\right)$ we denote $\mathcal{L}_s:= \tau(s)$. Explicitly, $\mathcal{L}_s$ is the line bundle associated to the $\mathcal{O}_{\mathfrak{X}}^{\times}$-torsor that is the preimage of $s$ under $\mathcal{M}_{\mathfrak{X}}^{\gp} \to \overline{\mathcal{M}}_{\mathfrak{X}}^{\gp}$.

\begin{lemma} \label{puncturedordersgeneral}
For any stable punctured map  $f: C \to \mathfrak{X}/\mathcal{S}, \ f \in \mathscr{M}(\mathfrak{X}, \beta, z)$ and any $s \in \Gamma\left(\mathfrak{X}, \overline{\mathcal{M}}_{\mathfrak{X}}^{\gp}\right)$ we have (using Notation \ref{logpairing}) $$\deg_{\mathcal{L}_s}(\beta) = \langle s, p_1 \rangle + \langle s, p_2 \rangle - \langle s, r \rangle .$$
\end{lemma}

\begin{proof}
This is \cite[Proposition 1.13]{GS3} specialized to $f \in \mathscr{M}(\mathfrak{X}, \beta, z)$.
\end{proof}

The fact that the product rule \eqref{prodtheta} respects the grading on $\hat{R}$ follows from Lemma \ref{puncturedordersgeneral} and the fact that $D$ is numerically equivalent to the trivial divisor. We note the following consequence of Lemma \ref{puncturedordersgeneral} for future use.

\begin{cor} \label{puncturedorders}
Suppose that $D$ is simple normal crossings. For any stable punctured map  $f: C \to \mathfrak{X}/\mathcal{S}, \ f \in \mathscr{M}(\mathfrak{X}, \beta, z)$ and any divisor $D'$ supported on $D$ we have $$\beta \cdot D' = \langle D', p_1 \rangle + \langle D', p_2 \rangle - \langle D', r \rangle .$$
\end{cor}

\begin{proof}
This is \cite[Corollary 1.14]{GS3} specialized to $f \in \mathscr{M}(\mathfrak{X}, \beta, z)$. Following the proof of \cite[Corollary 1.14]{GS3}, the result follows from Lemma \ref{puncturedordersgeneral} and the fact that $\overline{\mathcal{M}}_{\mathfrak{X}} = \bigoplus_{1 \leq i \leq m} \mathbb{N}_{D_i}$ where $\mathbb{N}_{D_i}$ is the constant sheaf on $D_i$ with stalk $\mathbb{N}$. Indeed, for $s$ given by $1 \in \mathbb{N}_{D_i}$ we have $\mathcal{L}_s = \mathcal{O}_{\mathfrak{X}}(-D_i)$, which implies the claim.
\end{proof}

One of the main results of \cite{GS3} (that we formulate here in the limit case) is as follows.
\begin{theorem}[{\cite[Theorem 1.9]{GS3}}]
The structure constants $N_{p_1 p_2 r}^{\beta}$ define via \eqref{prodtheta} and \eqref{sumalpha} a commutative, associative $\widehat{\kk [P]}$ algebra structure on $\widehat{R}$ with unit $\vartheta_0$.
\end{theorem}
We define the \emph{intrinsic mirror} as the family  
\begin{equation} \label{intmirror}
\check{\mathfrak{X}} : = \Proj \widehat{R} \to \Spec \widehat{\kk [P]}.
\end{equation}

Note that it is clear from the construction that if $P^{\times} = A_1(\mathfrak{X}_0, \mathbb{Z})_{\tors}$ is non-trivial, then the intrinsic mirror is a disjoint union of $|A_1(\mathfrak{X}_0, \mathbb{Z})_{\tors}|$ isomorphic copies.

\begin{obs} \label{heuristics}
We make two observations concerning the choice of the base monoid $P$:
\begin{enumerate}
    \item For any $NE(\mathfrak{X}_0) \subseteq P \subseteq P' \subseteq A_1(\mathfrak{X}_0, \mathbb{Z})$ where both $P$ and $P'$ satisfy Definition \ref{monoidass}, the mirror constructed using $P'$ is the basechange of the mirror constructed using $P$ by $P \hookrightarrow P'$. So the intrinsic mirror is independent of the choice of the base monoid in this sense. 
    \item The product \eqref{sumalpha} is defined over $NE(\mathfrak{X}_0)$. However, $NE(\mathfrak{X}_0)$ is not finitely generated in general (implying that $\widehat{\kk [NE(\mathfrak{X}_0)]}$ is not Noetherian), so we prefer to work over a (usually) larger finitely generated monoid $P$. We will sometimes think of the family over $NE(\mathfrak{X}_0)$ as a useful heuristic.  
\end{enumerate}
\end{obs}

\subsubsection{The mirror via the canonical scattering diagram \texorpdfstring{$\mathfrak{D}_I$}{D}} \label{canonicalmirror}

The paper \cite{GS4} gives a scattering diagram interpretation of the intrinsic mirror family. We have already described the affine manifold with singularities $B$, the monoid $P$, the ideal $I_0 = \mathfrak{m}:= P \setminus P^{\times}$, the MPA function $\varphi$ and the (trivial) initial slab functions. Let $I \subseteq P$ be any monoid ideal with $P \setminus I$ finite (we are mostly interested in $I = \mathfrak{m}^k, \ k \geq 1$). We will now describe a finite set of walls on $B$ (and thus define a scattering diagram, see Remark \ref{dec}) following \cite[Sections 2 and 3.2]{GS4}. See loc. cit. for details.

The walls are defined via punctured theory governed by tropical geometry. We review the notations of \cite{ACGS, ACGS2} for tropical maps to $(\Sigma(\mathfrak{X}), \Sigma \mathscr{P})$.\footnote{Here $\Sigma\mathscr{P}$ is the natural cone structure on $\Sigma(\mathfrak{X})$ provided by Construction \ref{tropicalize} of tropicalization. Naturally, we have $(\Sigma(\mathfrak{X}), \Sigma \mathscr{P}) \cong (\mathbf{C}B, \mathbf{C}\mathscr{P})$ but the notation $(\Sigma(\mathfrak{X}), \Sigma \mathscr{P})$ fits better with the conventions of \cite{GS4}. The cone structure $\Sigma \mathscr{P}$ is denoted by $\mathscr{P}$ in \cite{GS4}.} We consider genus zero graphs $G$ with sets of vertices $V(G)$, edges $E(G)$, and legs $L(G)$. Legs have only one vertex, correspond to marked and punctured points of a punctured curve, and are rays in the marked case and
compact line segments in the punctured case. An abstract tropical curve of genus zero over a rational polyhedral cone $\omega$ with integral structure $\Lambda_{\omega}$ is specified by the data $G: = (G, \ell)$ where $\ell:E(G) \to \Hom (\omega \cap \Lambda_{\omega}, \mathbb{N})$ determines the edge lengths. 

One can associate to $(G, \ell)$ and a rational polyhedral cone $\omega$ a generalized cone complex (a diagram in the category of rational polyhedral cones) along with a morphism of generalized cone complexes $\Gamma(G, \ell) \to \omega$. The fibre over $s\in\Int(\omega)$ is a metric graph with underlying graph
$G$ and the affine edge length of $E\in E(G)$ equal to $\ell(E)(s)\in \mathbb{R}_{\ge 0}$. 
Associated to each vertex $v\in V(G)$
of $G$ is a copy $\omega_v$ of $\omega$
in $\Gamma(G,\ell)$. Associated to each edge or leg $E\in E(G) \cup L(G)$
is a cone $\omega_E \in \Gamma(G,\ell)$ with
$\omega_E\subseteq \omega\times\mathbb{R}_{\ge 0}$ and the map to $\omega$ is
given by projection onto the first coordinate.

\begin{defn} \label{types}

A \emph{family of tropical maps
to $\Sigma(\mathfrak{X})$} over $\omega$ is a morphism of cone complexes
\[
h:\Gamma(G,\ell)\rightarrow \Sigma(\mathfrak{X}).
\]
If $s\in\Int(\omega)$, we may view the fibre of
$\Gamma(G,\ell)\rightarrow \omega$ over $s$ as a metric graph with underlying graph $G$, and
write
\[
h_s:G\rightarrow \Sigma(\mathfrak{X})
\] for the corresponding tropical map. 

The \emph{type} of $h$ consists of the data
$\tau:=(G,\bsigma,\mathbf{u})$ where
\[
\bsigma:V(G)\cup E(G)\cup L(G)\rightarrow \Sigma\mathscr{P}
\]
associates to $x\in V(G)\cup E(G)\cup L(G)$ the minimal
cone of $\Sigma \mathscr{P}$ containing $h(\omega_x)$. Further,
$\mathbf{u}$ associates to each (oriented) edge or leg $E\in E(G)\cup L(G)$
the corresponding \emph{contact order} $\mathbf{u}(E)\in \Lambda_{\bsigma(E)}$,
the image of the tangent vector $(0,1)\in \Lambda_{\omega_E}=
\Lambda_{\omega}\oplus\mathbb{Z}$ under the map $h$. 

We say a type $\tau$ is \emph{realizable} if there exists a family of tropical maps to $\Sigma(\mathfrak{X})$ of type $\tau$. If $\tau$ is realizable, then there exists a \emph{universal family} $h = h_{\tau}:\Gamma(G,\ell)\rightarrow \Sigma(\mathfrak{X})$ of tropical maps of type $\tau$ over a certain cone that we also denote by $\tau$. Generally, we write $h$ rather than $h_{\tau}$ when unambiguous.

A \emph{decorated type} is data $\btau=(\tau,\mathbf{A})$
where $\mathbf{A}:V(G)\rightarrow A_1(\mathfrak{X}_0, \mathbb{Z})$ associates a curve class to
each vertex of $G$. The \emph{total curve class} of ${\bf A}$
is $A=\sum_{v\in V(G)} {\bf A}(v)$.

For a decorated type $\btau$, we define $\Aut(\btau)$ as the group of automorphisms of the underlying graph $G$ preserving $\bsigma, \mathbf{u}$, and $\mathbf{A}$.

\end{defn}

Using the fact that tropical maps in the above sense arise as tropicalizations of punctured maps, one can define the notion of marking a punctured map by $\btau$ and the Deligne-Mumford stack $\mathscr{M}(\mathfrak{X}, \btau)$ of $\btau$-marked punctured maps to $\mathfrak{X}$, see \cite[Definition 3.8]{ACGS2}. The punctured theory here is rather technical. We do not review it and refer to \cite{ACGS2}.  

The walls of the canonical scattering diagram arise from a particular family of types $\btau$ such that the corresponding stack $\mathscr{M}(\mathfrak{X},\btau)$ is of virtual dimension $0$. First, the tropicalizations of stable punctured maps satisfy a certain balancing condition we impose on the types. 

\begin{defn} \label{balanced}
A type $\tau$ of a tropical map to $\Sigma(\mathfrak{X})$ is \emph{balanced} if:
\begin{enumerate}
\item
For each $v\in V(G)$ with $\bsigma(v)\in \Sigma\mathscr{P}$ a codimension zero or one
cone, a balancing condition holds at $v$. Namely, if $E_1,\ldots,E_m$ are the legs
and edges adjacent to $v$ and oriented away from $v$,
we interpret the contact orders $\mathbf{u}(E_i)$
as elements of $\Lambda_{h_s(v)}$ and require that
$$
\sum_{i=1}^m \mathbf{u}(E_i)=0.
$$
\item $\tau$ induces a type of a balanced tropical
map to $\Sigma(\mathcal{S}) = \mathbb{R}_{\geq 0}$ by composing with the $g_{\trop}$ of \eqref{inttrop}.
\end{enumerate}
\end{defn}

We impose some further restrictions that ensure that we get well-defined walls.

\begin{defn} \label{walltype}
A \emph{wall type} is a type $\tau=(G,\bsigma,{\mathbf u})$ 
of a tropical map to $\Sigma(\mathfrak{X})$ such that:
\begin{enumerate}
\item
$G$ is a genus zero graph with $L(G)=\left\{L_{\out}\right\}$ with $u_{\tau}:={\bf u}(L_{\out})\neq 0$.
\item $\tau$ is realizable and balanced.
\item Let $h:\Gamma(G,\ell)\rightarrow \Sigma(\mathfrak{X})$ be the corresponding
universal family of tropical maps, and let $\tau_{\out}\in \Gamma(G,\ell)$ be
the cone corresponding to $L_{\out}$. Then\footnote{There is a notational disparity with \cite[Definition 3.6]{GS4} here that might seem confusing. The reason is simply that in \cite{GS4}, $n$ is the dimension of $\mathfrak{X}$ and not the relative dimension of $\mathfrak{X} \to \mathcal{S}$.} $\dim\tau=n-1$ and
$\dim h(\tau_{\out})=n$. 
\end{enumerate}
A \emph{decorated wall type} is a decorated type $\btau=(\tau,{\bf A})$
with $\tau$ a wall type.

Now $h|_{\tau_{\out}}:\tau_{\out}\rightarrow \sigma$
induces a morphism
\[
h_*: \Lambda_{\tau_{\out}}\rightarrow \Lambda_{\sigma},
\]
and we define
\begin{equation}
\label{eq:ktau def}
k_{\tau}:=|\coker(h_*)_{\tors}|=
|\Lambda_{h(\tau_{\out})}/h_*(\Lambda_{\tau_{\out}})|.
\end{equation}
\end{defn}

\begin{remark} \label{dicinstead}
Note that the universal family $h:\Gamma(G,\ell)\rightarrow \Sigma(\mathfrak{X})$ of tropical maps in Definition \ref{walltype} carries the same information as its restriction $h_B:\Gamma(G,\ell)_B \rightarrow B$ to the dual intersection complex $B = g_{\trop}^{-1}(1) \subseteq \Sigma(\mathfrak{X})$. Indeed, this immediately follows from the fact that $\tau$ is balanced and from condition (2) of Definition \ref{balanced}. This phenomenon was first observed in \cite[Section 2.5.3]{ACGS}. 

It is often easier to perform explicit computations by studying tropical maps to $B$ instead of tropical maps to $\Sigma(\mathfrak{X})$ since this lowers $\dim \tau$ and $\dim h(\tau_{\out})$ by $1$. In particular, in the case that $\mathfrak{X} \to \mathcal{S}$ is a minimal log CY degeneration of K3-s, a wall type $\tau$ defines a single tropical curve $h_B:\Gamma(G,\ell)_B \rightarrow B$ that is \emph{rigid} in the sense of \cite[Definition 3.6]{ACGS} (since we require $\dim \tau = n-1$ in the definition of a wall type).
\end{remark}

We need the following definition to define the canonical scattering diagram on $B$ and not on $\mathbf{C}B$ (as in \cite[Construction 3.13]{GS4}). 

\begin{defn} \label{index}
Let $\mathbf{C}\mathfrak{p}$ be a cone over a codimension one rational polyhedral subset $\mathfrak{p}$ of some $\sigma \in \mathscr{P}^{\max}$. Then there is a map $\Lambda_{\mathbf{C}\mathfrak{p}} \to \mathbb{Z}$ (projection to the first component, see \eqref{coneover}) that is not necessarily surjective. The image is of the form $\ind(\mathbf{C}\mathfrak{p}) \cdot \mathbb{Z}$ for some $\ind(\mathbf{C}\mathfrak{p}) \in \mathbb{N}$. Following \cite[Definition 4.2.2]{GHS}, we call $\ind(\mathbf{C}\mathfrak{p})$ the \emph{index} of $\mathbf{C}\mathfrak{p}$. 

For $\tau$ a wall type we denote $\ind(\tau):= \ind(h(\tau_{\out}))$ and we set $\ind(\btau):=\ind(\tau)$ for a decorated wall type $\btau$ with underlying wall type $\tau$.
\end{defn}

\begin{observe} \label{nnmn}
Note that $\ind(\mathbf{C}\mathfrak{p}) = 1$ if $\mathfrak{p} \subseteq \rho \in \mathscr{P}^{[n-1]}$ since in this case $\rho$ contains integral points.
\end{observe}

We are now ready to define the canonical scattering diagram.

\begin{cons} \label{canonical scattering}
Fix a decorated wall type $\btau=(\tau,{\bf A})$, and let $A=\sum_{v\in V(G)}
{\bf A}(v)$ be the total curve class.
Then $\mathscr{M}(\mathfrak{X},\btau)$ is of virtual dimension $0$ (see the proof of \cite[Lemma 3.9]{GS4}) and one can define
\[
W_{\btau}:=\frac{\deg [\mathscr{M}(\mathfrak{X},\btau)]^{\virt}}{ |\Aut(\btau)|}.
\]
For every $\btau$ such that $W_{\btau} \neq 0$ we define a wall
\begin{equation} \label{refinedwalleq}
\mathfrak{p}_{\btau}:=
\left(h(\tau_{\out}) \cap g_{\trop}^{-1}(1), \ \exp(k_{\tau} W_{\btau}t^A z^{-u_{\tau}})^{\ind(\btau)}\right)
\end{equation} 
Here we view $t^A z^{-u_{\tau}}$ as a monomial in $\kk[\mathcal{P}^+_x]$
for $x\in \Int(h(\tau_{\out}))$
using Notation \ref{notation}. Note that
$(\kk[P]/I)[\Lambda_{h(\tau_{\out})}] \subseteq \kk[\mathcal{P}^+_x]/\mathcal{I}_x$, and $\exp(k_{\tau} W_{\btau}t^A z^{-u_{\tau}})^{\ind(\btau)}$ makes sense as an element in $(\kk[P]/I)[\Lambda_{h(\tau_{\out})}]$ by removing all the monomials that are zero in $(\kk[P]/I)[\Lambda_{h(\tau_{\out})}]$.
We define the \emph{canonical scattering diagram} as follows.
\begin{equation} \label{canonicalscattering}
\mathfrak{D}_I:=\left\{\mathfrak{p}_{\btau}\,\Big|\,
\substack{\hbox{$\btau$ an isomorphism class of 
a decorated wall}\\\hbox{type with total curve class lying in 
$P\setminus I$}}\right\}.
\end{equation}
\end{cons}

The following theorem is the first major result of \cite{GS4}.

\begin{theorem}[{\cite[Theorem 5.2]{GS4}}] \label{canonicalconsistency}
$\mathfrak{D}_I$ is a consistent scattering diagram.
\end{theorem}

\begin{remarks} \label{diagramandcone}
Some remarks are in order:
\begin{enumerate}
    \item The scattering diagram of \cite{GS4} is defined on $\mathbf{C}B$, rather than $B$. Our definition ensures (see \cite[Definition 4.2.4]{GHS}) that the lifting of the scattering diagram of Construction \ref{canonical scattering} to $\mathbf{C}B$ is the scattering diagram of \cite[Construction 3.13]{GS4}.
    \item \cite[Theorem 5.2]{GS4} proves consistency of the lifted scattering diagram on $\mathbf{C}B$. However, by reversing the proof of \cite[Proposition 4.2.6]{GHS}, one can see that consistency of a scattering diagram on $B$ is equivalent to consistency of the lifted scattering diagram on $\mathbf{C}B$. 
    \item \cite[Theorem 5.2]{GS4} proves consistency of the undecorated version \cite[Definition 3.10]{GS4} of the canonical scattering diagram (using $\tau$ instead of $\btau$). However, the two scattering diagrams are combinatorially equivalent, see \cite[Construction 3.13]{GS4}.
\end{enumerate}
\end{remarks}

We now define 
\begin{equation} \label{intmirrormodI}
\check{\mathfrak{X}}_{\mathfrak{D}_I} \to \Spec \kk[P]/I
\end{equation}
following the recipe of Section \ref{thetafunctions}. The scattering diagram $\mathfrak{D}_{\mathfrak{m}^k}$ is compatible with $\mathfrak{D}_{\mathfrak{m}^{k-1}}$ (in the sense of Definition \ref{compatibility}) for $k \geq 1$  so the families $\check{\mathfrak{X}}_{\mathfrak{D}_{\mathfrak{m}^k}} \to \Spec \kk[P]/\mathfrak{m}^k$ for $k \geq 1$ form an inverse system and taking the limit over this system gives $$\check{\mathfrak{X}}_{\mathfrak{D}} \to \Spec \widehat{\kk[P]}$$ where $\widehat{\kk[P]}$ is the completion with respect to $\mathfrak{m}$ as in Section \ref{3pointmirror}. The second major result of \cite{GS4} (formulated in the limit case) is proved by comparing the two product formulas for theta functions: the product formula \eqref{prodtheta} defined via three-pointed punctured maps and the broken line product formula on $\mathbf{C}B$ (see \cite[Theorem 3.5.1]{GHS}).

\begin{theorem}[{\cite[Theorem 6.1]{GS4}}]
The family $\check{\mathfrak{X}}_{\mathfrak{D}} \to \Spec \widehat{\kk[P]}$ is isomorphic to the intrinsic mirror family $\check{\mathfrak{X}} \to \Spec \widehat{\kk [P]}$ of \eqref{intmirror}.
\end{theorem}

So we may use the constructions of the intrinsic mirror of this section and of Section \ref{3pointmirror} interchangeably. 

\begin{remark}
For any ideal $I$ with $\sqrt{I} = \mathfrak{m}$ we can talk about the \emph{intrinsic mirror modulo $I$}, which is just the family of \eqref{intmirrormodI}. Using the setup of Section \ref{3pointmirror}, this corresponds to viewing the coefficients $\alpha_{p_1 p_2 r}$ defined by \eqref{sumalpha} modulo $I$ and replacing $\widehat{R}$ in \eqref{ringofdef} with 
\begin{equation}\label{thetaringmodI}
R_I := \bigoplus_{p \in \mathbf{C}B(\mathbb{Z})} \left(\kk [P]/I \right) \vartheta_p.
\end{equation}
\end{remark}

\subsection{Resolution setup} \label{extresults}

Suppose that we have a projective resolution $ \pi: \mathfrak{X} \to \bar{\mathfrak{X}}$ where $\bar{g}: \bar{\mathfrak{X}} \to \mathcal{S}$ is a toric degeneration with polarization $A$ and the composed $g: \mathfrak{X} \to \mathcal{S}$ is a log smooth minimal log CY degeneration satisfying Assumption \ref{manifoldass}. Note that since $(\bar{\mathfrak{X}}, \bar{D})$ is a Zariski log scheme (by Assumption \ref{toricdegenerationass}(3)), so is $(\mathfrak{X}, D)$. Projectivity of $\mathfrak{X} \to \mathcal{S}$ follows from projectivity of $\bar{\mathfrak{X}} \to \mathcal{S}$ (i.e. Assumption \ref{toricdegenerationass}(2)) and projectivity of the resolution $ \pi: \mathfrak{X} \to \bar{\mathfrak{X}}$. 

In this section, we will explain that one may choose the monoid $P$ of Definition \ref{monoidass} so that the curves contracted by $\pi$ are contained in a face $K$ of $P$. We will also study extensions of the intrinsic mirror $\check{\mathfrak{X}} \to \Spec \widehat{\kk[P]}$ and show that under certain assumptions on the resolution $\pi: \mathfrak{X} \to \bar{\mathfrak{X}}$ one can extend $\check{\mathfrak{X}} \to \Spec \widehat{\kk[P]}$ to a family over the completion of $\Spec \kk[K]$.

\subsubsection{Restricting to well-chosen monoids \texorpdfstring{$P$}{P}}

The map $\pi: \mathfrak{X} \to \bar{\mathfrak{X}}$ contracts some curve classes. 

\begin{defn} \label{wellchosendef}
We say that a monoid $P$ satisfying Definition \ref{monoidass} is \emph{well-chosen} if there exists a face $K$ of $P$ such that for any curve $C \subseteq \mathfrak{X}$ we have $$C \in K \cap NE(\mathfrak{X}_0) \text{ if and only if } \dim(\pi(C)) = 0.$$ We always assume that $K$ is the minimal face of $P$ with this property.
\end{defn}

\begin{fact} \label{wellchosen}
For any monoid $P'$ satisfying the assumptions of Definition \ref{monoidass} there exists a well-chosen monoid $P$ with $NE(\mathfrak{X}_0) \subseteq P \subseteq P' \subseteq A_1(\mathfrak{X}_0, \mathbb{Z})$.
\end{fact}

\begin{proof}
Note that $\pi^{*}A$ is $g$-nef and that a curve  $C \subseteq \mathfrak{X}_0$ is contracted by $\pi: \mathfrak{X} \to \bar{\mathfrak{X}}$ if and only if $\pi^{*}A \cdot C = 0$. We let $P = P' \cap \{\beta \in A_1(\mathfrak{X}_0, \mathbb{Z}) \ | \ \pi^{*} A \cdot \beta \geq 0\}$ and let $K = P' \cap \{\beta \in A_1(\mathfrak{X}_0, \mathbb{Z}) \ | \ \pi^{*} A \cdot \beta = 0\}$. Then clearly, $K$ is a face of $P$, and $P$ is well-chosen. 
\end{proof}

\begin{fact} \label{enoughforwellchosen}
Suppose that Conjecture \ref{conjecture} holds for a toric degeneration $\bar{\mathfrak{X}} \to \mathcal{S}$, a choice of resolution $\pi: \mathfrak{X} \to \bar{\mathfrak{X}}$ and all well-chosen monoids $P$. Then it holds for $\bar{\mathfrak{X}} \to \mathcal{S}$, $\pi$, and all monoids $P'$ satisfying Definition \ref{monoidass}. 
\end{fact}

\begin{proof}
Let $P'$ be any monoid satisfying Definition \ref{monoidass} and let $P$ be the well-chosen monoid with $NE(\mathfrak{X}_0) \subseteq P \subseteq P' \subseteq A_1(\mathfrak{X}_0, \mathbb{Z})$ given by Proposition \ref{wellchosen}. By Observation \ref{heuristics}(1), $\check{\mathfrak{X}} \to \Spec \widehat{\kk[P']}$ is the basechange of $\check{\mathfrak{X}} \to \Spec \widehat{\kk[P]}$ by the inclusion $P \xhookrightarrow{} P'$. If the completion in $\check{\mathfrak{X}} \to \Spec \widehat{\kk[P]}$ is taken with respect to an ideal $J$, then we can take the completion in $\check{\mathfrak{X}} \to \Spec \widehat{\kk[P']}$ with respect to the ideal generated by the image of $J$ under $P \xhookrightarrow{} P'$. 
It is enough to check that the basechange of the intrinsic mirror $\check{\mathfrak{X}} \to \Spec \widehat{\kk[P']}$ by $h': P' \to \mathbb{N}, \ \beta \mapsto \pi^{*}A \cdot \beta$ coincides with the basechange of $\check{\mathfrak{X}} \to \Spec \widehat{\kk[P]}$ by $h: P \to \mathbb{N}, \ \beta \mapsto \pi^{*}A \cdot \beta$.  Indeed, the two basechanges agree since $h: P \to \mathbb{N}$ factors as $P \xhookrightarrow{} P' \overset{h'} \to \mathbb{N}$.
\end{proof}

\subsubsection{The minimal ideal for extension}

In Sections \ref{3pointmirror} and \ref{canonicalmirror}, we defined the intrinsic mirror $\check{\mathfrak{X}} \to \Spec \widehat{\kk [P]}$ to a minimal log CY degeneration $\mathfrak{X} \to \mathcal{S}$ satisfying Assumption \ref{manifoldass} to lie over the completion of the smallest toric stratum $\Spec \kk [P]/\mathfrak{m}$ (here $\mathfrak{m} = P \setminus P^{\times}$ is the maximal ideal of $P$). However, the mirror can often be defined over the completion of a larger union of toric strata $\Spec \kk[P]/J$ for some monoid ideal $J \subseteq \mathfrak{m} \subseteq P$. This phenomenon was first observed in \cite[Section 6]{GHK} in the case of a log CY surface $(\mathfrak{X}, D)$. Our analysis in this section is similar but simpler since, as explained in Section \ref{3pointmirror}, the sum $\vartheta_{p_1} \cdot \vartheta_{p_2} = \sum_{r \in \mathbf{C}B(\mathbb{Z})} \alpha_{p_1 p_2 r} \vartheta_r$ in \eqref{prodtheta} is always a finite sum.  

\begin{sloppypar}
Note that for any monoid ideal $J \subseteq P$ such that the sums $\alpha_{p_1 p_2 r}= \sum_{\beta \in P \setminus J} N_{p_1 p_2 r}^{\beta} t^{\beta}$ of \eqref{sumalpha} are finite, the intrinsic mirror is defined over $\kk [P]/J$ via \eqref{prodtheta} and \eqref{sumalpha}. Let 
\begin{equation} \label{idealcollection}
\mathcal{A}:=\left\{
\substack{\hbox{monoid ideals  $J \subseteq P$ such that  $\alpha_{p_1 p_2 r}= \sum_{\beta \in P \setminus J} N_{p_1 p_2 r}^{\beta} t^{\beta}$}\\\hbox{is finite for all   $p_1, p_2, r \in \mathbf{C}B$. }}\right\}.
\end{equation}
\end{sloppypar}

\begin{lemma} \label{p1}
		If $J \in \mathcal{A}$ and $J \subset J'$, then $J' \in \mathcal{A}$. In addition, $\mathcal{A}$ is closed under finite intersections. 
	\end{lemma}
	
	\begin{proof}
		The first statement follows since any $\beta \in P \setminus J'$ is contained in $P \setminus J$. For the second statement let $J_1, J_2 \in \mathcal{A}$. Then we know that $\sum_{\beta \in P \setminus J_1} N_{p_1 p_2 r}^{\beta} t^{\beta}$ and $\sum_{\beta \in P \setminus J_2} N_{p_1 p_2 r}^{\beta} t^{\beta}$ are finite sums. But then $$\sum_{\beta \in P \setminus (J_1 \cap J_2)} N_{p_1 p_2 r}^{\beta} t^{\beta} = \sum_{\beta \in (P \setminus J_1) \cup (P \setminus J_2)} N_{p_1 p_2 r}^{\beta} t^{\beta}$$ is also a finite sum.
	\end{proof}

\begin{fact} \label{imin}
There is a unique radical monomial ideal $I_{\min}$ such that if $I_{\min} \subseteq \sqrt{J}$, then $J \in \mathcal{A}$.
\end{fact}

\begin{proof}
We have $J \in \mathcal{A}$ for any $J$ with $\sqrt{J}=\mathfrak{m}$ (since $\sqrt{J} = \mathfrak{m}$ if and only if $P \setminus J$ is finite). Note that a radical monomial ideal is the complement of a union of faces of $P$, so there are only a finite number of such ideals. Suppose $I_1, I_2$ are two radical ideals such that if $I_i \subseteq \sqrt{J}$, then $J \in \mathcal{A}$ for $i=1,2$. Note that any ideal $J$ with $I_1 \cap I_2 \subseteq \sqrt{J}$ can be written as $J_1 \cap J_2$ with $I_i \subseteq \sqrt{J_i}$ (Indeed, we can use the primary decomposition of $J$. If $J = \cap_k \mathfrak{p}_k$ is an
intersection of primary ideals, necessarily the prime ideal $\sqrt{\mathfrak{p}_k}$ contains either $I_1$ or $I_2$ for each $k$. Then let $J_1$ be the intersection of those $\mathfrak{p}_k$ whose radical contains $I_1$ and $J_2$
be the intersection of those $\mathfrak{p}_k$ whose radical contains $I_2$). Lemma \ref{p1} implies that $J \in \mathcal{A}$. This shows the existence of $I_{\min}$.
\end{proof}

\subsubsection{Extending the intrinsic mirror family: simple normal crossings case.} \label{extendsnccase}

We want to understand $I_{\min}$ better. Assume for now that $D$ is simple normal crossings. We shall treat the general case in the next section. It is easy to prove the following using the same reasoning as in \cite[Remark 1.22]{GS3}.

\begin{fact} \label{P1}
    For any contraction $\mathfrak{X} \overset{\pi'}{\to} \mathfrak{X}'$ (not necessarily to a toric degeneration), let $P$ be a well-chosen monoid and $K \subseteq P$ be the face containing the contracted curves. Suppose that there exists a $\pi'$-ample $D'$ supported on $D$. Then $ P \setminus K \in \mathcal{A}$.
\end{fact}
	
\begin{proof}
	For any curve class $\beta \in P$, the intersection number $\beta \cdot D'$ is determined by $p, q, r$ by Corollary \ref{puncturedorders}. Because $D'$ is $\pi'$-ample, the relative Hilbert scheme is of finite type, and there are only a finite number of effective curve classes $\beta \in K$ with the given value of $\beta \cdot D'$. So every $\alpha_{p_1 p_2 r}$ is a finite sum and $P \setminus K \in \mathcal{A}$. 
\end{proof}

We would like to further show that $I_{\min} \subseteq P \setminus K$. We note as in \cite[Construction 1.24]{GS3} that for any monoid ideal $J \subseteq P$, there is a natural grading (corresponding to a torus action on the intrinsic mirror modulo $J$) of $R_J$ (as defined in \eqref{thetaringmodI}) by $\Div_D(\mathfrak{X})$, the group of divisors supported on $D$. Indeed, define the $D_i$-component of the degree by $\deg_{D_i} \vartheta_p := \langle D_i, p \rangle$ for $p \in \mathbf{C}B(\mathbb{Z})$ (using Notation \ref{logpairing}) and $\deg_{D_i} t^\beta := D_i \cdot \beta$ for $t^\beta \in \kk [P]/J$. Additivity of the degrees follows from Corollary \ref{puncturedorders}. 

Now, let $(e_{D_1}, \dots, e_{D_m})$ be the basis for $\mathbb{Z}^{\Div_D(\mathfrak{X})} \cong \mathbb{Z}^m$. We define a character map 
\begin{equation} \label{charactermap}
w: \mathbf{C}B(\mathbb{Z}) \times P \setminus J \to \mathbb{Z}^m, \quad (p, \beta) \mapsto \sum_{i=1}^m \langle D_i, p \rangle e_{D_i} + \sum_{i=1}^m (D_i \cdot \beta)e_{D_i}.
\end{equation}

Now, we need a technical lemma.

\begin{lemma} \label{monoidlemma}
Let $Q$ be a toric monoid with $Q^{\times}$ finite, $I \subseteq Q$ a monoid ideal with $\sqrt{I}$ a prime ideal, and $$w: Q \setminus I \to \mathbb{Z}^m$$ a linear map to a finitely generated free abelian group. Then $w$ has finite fibres if and only if $$\ker w \cap \left(Q \setminus \sqrt{I} \right) = Q^{\times}.$$
\end{lemma}

\begin{proof}
Since $Q^{\times}$ is finite, it is torsion so $w(q) = 0$ for any $q \in Q^{\times}$ by linearity of $w$. Let $\alpha: Q \to Q/Q^{\times}$ be the canonical quotient map. Then $\alpha(I)$ is an ideal of $Q/Q^{\times}$, $w:Q \setminus I \to \mathbb{Z}^m$ factors through $\left(Q/Q^{\times}\right) \setminus \alpha(I)$, and the fibres of $Q\setminus I \to \left(Q/Q^{\times}\right) \setminus \alpha(I)$ have cardinality $|Q^{\times}|$. So it suffices to show the claim for a monoid $Q$ that is toric and sharp.

A toric and sharp $Q$ can be viewed as the set of integral points of a strictly convex rational polyhedral cone. So there exists a unique minimal system of generators $v_1, \dots, v_k \in Q$ such that every $q \in Q$ is of the form $q = \sum_{i=1}^k a_i v_i$ with $a_i \in \mathbb{N}$ for all $1 \leq i \leq k$. Every face of $Q$ is generated by a subset of $ \left\{v_1, \dots, v_k \right\}$  so we have two subsets of generators $\left\{\alpha_1, \dots, \alpha_{k_1} \right\}$ and $\left\{\beta_1, \dots, \beta_{k_2} \right\}$ such that $$Q \setminus \sqrt{I} = \langle \alpha_1, \dots, \alpha_{k_1} \rangle, \quad Q = \langle \alpha_1, \dots, \alpha_{k_1}, \beta_1, \dots, \beta_{k_2} \rangle = \langle v_1, \dots, v_k \rangle.$$

The fibre of $w: Q \setminus I \to \mathbb{Z}^m$ over $d \in \mathbb{Z}^m$ is the set of $ \left\{ q   \in Q \setminus I \ | \ w(q) = d \right\}$. That is, it is the set of solutions of the system of linear equations 
\begin{equation} \label{lineardioff}
w(q) = a_1w(\alpha_1) + \dots + a_{k_1}w(\alpha_{k_1}) + b_1 w(\beta_1) + \dots + b_{k_2} w(\beta_{k_2}) = d
\end{equation}
with $a_i, b_j \in \mathbb{N}$ and $q = a_1 \alpha_1 + \dots + a_{k_1} \alpha_{k_1} + b_1 \beta_1 + \dots + b_{k_2} \beta_{k_2} \in Q \setminus I$ a general element.

Since $Q \setminus I$ is the complement of an ideal, we have $b_j \beta_j \in Q \setminus I$ for all $1 \leq j \leq k_2$. Also, we have $b_j \beta_j \notin Q \setminus \sqrt{I}$ by construction so $b_j \beta_j \in \sqrt{I} \setminus I$. But then $b_j \beta_j \in I$ for large enough $b_j \in \mathbb{N}$ so there are a finite number of choices for $b_j, \ 1 \leq j \leq k_2$. 

Therefore, it is enough to show that there are a finite number of solutions of \eqref{lineardioff} for any fixed choices of $b_j = b_j^0 \in \mathbb{N}, \ 1 \leq j \leq k_2$ if and only if $\ker w \cap \left(Q \setminus \sqrt{I} \right) = 0$. This amounts to checking that the system of equations $$ a_1w(\alpha_1) + \dots + a_{k_1}w(\alpha_{k_1}) = \tilde{d}$$ with $\tilde{d} := d -  b^0_1 w(\beta_1) - \dots - b^0_{k_2} w(\beta_{k_2})$ has finitely many solutions with $a_i \in \mathbb{N}$. Every equation in this system is a linear Diophantine equation with integer coefficients. It is classical (see, e.g. \cite[Lemma 1]{CF}) that such an equation has finitely many solutions in non-negative integers if and only if the corresponding homogeneous equation has finitely many solutions in non-negative integers. So it is enough to check that $$ a_1w(\alpha_1) + \dots + a_{k_1}w(\alpha_{k_1}) = 0$$ has finitely many solutions with $a_i \in \mathbb{N}$. But this is equivalent to $\ker w \cap \left(Q \setminus \sqrt{I} \right) = 0$ by construction.
\end{proof}

\begin{fact} \label{P2}
In the setup of Proposition \ref{P1}, suppose in addition that:
\begin{enumerate}
    \item $D' = \sum_{i=1}^{m} a_i D_i$ is an effective divisor supported on $D$ and such that $D_i$ is $\pi'$-nef for any $i$ with $a_i = 0$ (e.g. $a_i > 0$ for $1 \leq i \leq m$).
    \item $K \cap NE(\mathfrak{X}_0) = K \cap NE(\mathfrak{X}_0)_{\numer}$ (under the splittings of \eqref{standardsplittings}), i.e. $\pi'$ only contracts numerical classes of curves $C \subseteq \mathfrak{X}_0$.
\end{enumerate} Then we have $I_{\min} \subseteq P \setminus K$.
\end{fact}

\begin{proof}
We first want to show that we can choose a monoid $P$ and a face $K$ so that $D' \cdot \beta > 0$ for any $\beta \in K \setminus P^{\times}$. 

\textbf{Step 1.} Suppose without loss of generality that $D' = \sum_{i=1}^{m'} a_i D_i$ with $a_i > 0$. Let $D_{\rel} : = \sum_{i=1}^{m'} D_i$ and $D_{\irrel}:=\sum_{i=m'+1}^{m} D_i$ so that $D = D_{\rel} + D_{\irrel}$. Let $\Delta : = D_{\rel} - \varepsilon D'$ for some rational $ 0 < \varepsilon \ll 1$, then $\lfloor \Delta \rfloor = 0$ (here we define $\lfloor \sum_{i=1}^m b_i D_i \rfloor := \sum_{i=1}^m \lfloor b_i \rfloor D_i$), and $\Delta$ is effective and simple normal crossings (as it is supported on a simple normal crossings $D$). Also, $K_{\mathfrak{X}} + \Delta = -D_{\irrel} - \varepsilon D'$ is $\mathbb{Q}$-Cartier (some rational multiple of it is Cartier) since all the $D_i, \ 1 \leq i \leq m$ are Cartier and $\varepsilon$ is rational. Under these conditions, $(\mathfrak{X}, \Delta)$ is \emph{dlt} (\emph{divisorial log terminal}, see \cite[Definition 2.37]{KM}) by \cite[Proposition 2.40]{KM}. Further, since $\lfloor \Delta \rfloor = 0$, it is actually \emph{klt} (\emph{Kawamata log
terminal}, see \cite[Definition 2.34]{KM}) by \cite[Proposition 2.41]{KM}. 

\textbf{Step 2.} Suppose that $A_1(\mathfrak{X}_0, \mathbb{Z}):= A_1(\mathfrak{X}/\mathcal{S}, \mathbb{Z})_{\numer}$. It is clear from assumption (1) that $K_{\mathfrak{X}} + \Delta = -D_{\irrel} - \varepsilon D'$ is negative on $K \cap NE(\mathfrak{X}/\mathcal{S})_{\numer}$ so $K \cap NE(\mathfrak{X}/\mathcal{S})_{\numer}$ is a $K_{\mathfrak{X}} + \Delta$-negative extremal face of $NE(\mathfrak{X}/\mathcal{S})_{\numer}$. By the relative cone theorem \cite[Theorem 3.25]{KM}, $NE(\mathfrak{X}/\mathcal{S})_{\numer}$ is rational polyhedral near $K \cap NE(\mathfrak{X}/\mathcal{S})_{\numer}$. So there exists a monoid $P'$ with $NE(\mathfrak{X}/\mathcal{S})_{\numer} \subseteq P' \subseteq P \subseteq A_1(\mathfrak{X}/\mathcal{S}, \mathbb{Z})_{\numer}$, satisfying Definition \ref{monoidass}, and such that $K \cap NE(\mathfrak{X}/\mathcal{S})_{\numer}$ is a face of $P'$. By Observation \ref{heuristics}(1), $\check{\mathfrak{X}} \to \Spec \widehat{\kk[P']}$ is the basechange of $\check{\mathfrak{X}} \to \Spec \widehat{\kk[P']}$ by the inclusion $P' \xhookrightarrow{} P$, so (replacing $P$ by $P'$) we can assume that $K = K \cap NE(\mathfrak{X}/\mathcal{S})_{\numer}$. Then $D'$ is positive on $K \setminus P^{\times} = K \setminus \{0\}$.

Using the splittings of Construction \ref{splittings}, the same claim is true for other choices of $A_1(\mathfrak{X}_0, \mathbb{Z})$. Indeed, we have $NE(\mathfrak{X}_0)_{\numer} = NE(\mathfrak{X}_0) \cap A_1(\mathfrak{X}_0, \mathbb{Z})_{\numer}$ under the splittings. Fix a choice of $P'$ as above (for any corresponding choice of $P$ with $NE(\mathfrak{X}/\mathcal{S})_{\numer} \subseteq  P \subseteq A_1(\mathfrak{X}/\mathcal{S}, \mathbb{Z})_{\numer}$). Then if using $A_1(\mathfrak{X}_0, \mathbb{Z})_{\numer}$ (and a choice of $P$ and $K$), we can take the corresponding restricted monoid and face to be $P' \cap P \cap A_1(\mathfrak{X}_0, \mathbb{Z})_{\numer}$ and $P' \cap K \cap NE(\mathfrak{X}_0)_{\numer}$ respectively. Now, if $A_1(\mathfrak{X}_0, \mathbb{Z}) \cong A_1(\mathfrak{X}_0, \mathbb{Z})_{\numer} \oplus G$, then assumption (2) implies that we can take the restricted monoid and face to be $\left(\left(P' \cap A_1(\mathfrak{X}_0, \mathbb{Z})_{\numer} \right) \oplus G \right) \cap P$ and $\left(\left(P' \cap A_1(\mathfrak{X}_0, \mathbb{Z})_{\numer} \right) \oplus G \right) \cap \left(K \cap NE(\mathfrak{X}_0)\right)$ respectively. 

\textbf{Step 3.} Assume that we have a choice of $P$ and $K$ such that $D'$ is positive on $K \setminus P^{\times}$ and let $J \subseteq P$ be a monoid ideal with $\sqrt{J} = P \setminus K$. Combining together \eqref{prodtheta} and \eqref{sumalpha}, we have
$$
\vartheta_{p_1} \cdot \vartheta_{p_2} = \sum_{r \in \mathbf{C}B(\mathbb{Z}), \ \beta \in P \setminus J} N_{p_1 p_2 r}^{\beta} t^{\beta} \vartheta_r
$$ and all the $t^{\beta}\vartheta_r$ with $N_{p_1 p_2 r}^{\beta} \neq 0$ are of the same $\Div_D(\mathfrak{X})$-degree by additivity of the degrees. 

It is enough to show that the character map $w$ of \eqref{charactermap} has finite fibres. Clearly, it suffices to consider the fibres of $\mathbf{C}\sigma (\mathbb{Z}) \times P \setminus J \to \mathbb{Z}^m$ for every $\mathbf{C}\sigma \in \mathbf{C}\mathscr{P}^{\max}$. Fix a $\mathbf{C}\sigma \in \mathbf{C}\mathscr{P}^{\max}$ and note that $\mathbf{C}\sigma(\mathbb{Z}) \times P$ is a toric monoid with $(\mathbf{C}\sigma(\mathbb{Z}) \times P)^{\times} = P^{\times}$ finite, and that $w$ is linear on $\mathbf{C}\sigma(\mathbb{Z}) \times P$. Also $\mathbf{C}\sigma (\mathbb{Z}) \times J$ is an ideal of $\mathbf{C}\sigma(\mathbb{Z}) \times P$ with prime radical $\sqrt{\mathbf{C}\sigma (\mathbb{Z}) \times J} = \mathbf{C}\sigma \times \sqrt{J}$ that is the complement of the face $\mathbf{C}\sigma(\mathbb{Z}) \times K$ of $\mathbf{C}\sigma(\mathbb{Z}) \times P$.

\textbf{Step 4.} We are in the setup of Lemma \ref{monoidlemma}, so it is enough to check that $\ker w \cap (\mathbf{C}\sigma(\mathbb{Z}) \times K) = P^{\times}$. Suppose that $(r, \beta) \in \ker w \cap (\mathbf{C}\sigma(\mathbb{Z}) \times K)$. We have $\mathbf{C}\sigma = \langle v_1, \dots, v_{n+1} \rangle$ for $v_i \in \mathscr{P}^{[0]}, \ 1 \leq i \leq n+1$ the generators of the rays of $\mathbf{C}\sigma$. Then $r = \sum_{i=1}^{n+1} a_i v_i$ for some $a_i \in \mathbb{N}, \ 1 \leq i \leq n+1$. Now $$w(r, \beta) = \sum_{i=1}^{n+1} a_i e_{D_{v_i}} + \sum_{j=1}^m (D_{j} \cdot \beta) e_{D_{j}} = 0,$$ so we have $D_j \cdot \beta \leq 0$ for all $1 \leq j \leq m$. In particular, $D' \cdot \beta \leq 0$ since $D'$ is effective. But $D'$ is positive on $K \setminus P^{\times}$ so $\beta \in P^{\times}$. Therefore, $D_j \cdot \beta = 0$ for all $1 \leq j \leq m$ and we have $\sum_{i=1}^{n+1} a_i e_{D_{v_i}} = 0$. But then $a_i = 0$ for all $1 \leq i \leq n+1$. So $r=0$ and $\ker w \cap (\mathbf{C}\sigma(\mathbb{Z}) \times K) = P^{\times}$. 
\end{proof}

\begin{remark} \label{irreldivisor}
Note from the proof that it is enough to require that $D_{\irrel}$ is $\pi'$-nef instead of the second part of condition (1) of Proposition \ref{P2}. 
\end{remark}

\subsubsection{Extending the intrinsic mirror family: general case} \label{generalextensionstuff}

We now want to generalize Propositions \ref{P1} and \ref{P2} to the case when the divisor $D$ is not simple normal crossings. Note that the components $D_i, \ 1 \leq i \leq m$ are not generally Cartier, and we can't use Corollary \ref{puncturedorders}. One can still argue similarly to Section \ref{extendsnccase} in the case that $D_i, \ 1 \leq i \leq m$ are $\mathbb{Q}$-Cartier (e.g. if $\mathfrak{X}$ is $\mathbb{Q}$-factorial) and all the cells $\sigma \in \mathscr{P}^{\max}$ are simplices. We shall consider the general case. 

In general, it is natural to use Lemma \ref{puncturedordersgeneral} instead of Corollary \ref{puncturedorders}. Lemma \ref{puncturedordersgeneral} applies to line bundles that lie in the image of the map $\tau: \Gamma\left(\mathfrak{X}, \overline{\mathcal{M}}_{\mathfrak{X}}^{\gp}\right) \to \Pic(\mathfrak{X}), \ s \mapsto \mathcal{L}_s$ of \eqref{connectinghom}. Accordingly, we are going to use the group $\Gamma\left(\mathfrak{X}, \overline{\mathcal{M}}_{\mathfrak{X}}^{\gp}\right)$ instead of the group $\Div_D(\mathfrak{X})$ we used in the previous section. 

First, we want to understand $\Gamma\left(\mathfrak{X}, \overline{\mathcal{M}}_{\mathfrak{X}}^{\gp}\right)$ and $\tau$ better. Let $PA(B):=PA(B, \mathbb{Z})$ and $PL(\mathbf{C}B):=PL(\mathbf{C}B, \mathbb{Z})$ be the groups of integral piecewise-affine (PA) functions on $B$ and integral piecewise-linear (PL) functions on $\mathbf{C}B$ respectively.\footnote{Not to confuse with the multi-valued functions of Section \ref{mpafunc}.} 

\begin{fact} \label{logpicard}
Let $\mathfrak{X} \to \mathcal{S}$ be a log smooth minimal log CY degeneration (with the log structure $\mathcal{M}_{\mathfrak{X}}$ on $\mathfrak{X}$ fine, saturated, and Zariski as usual). Then: 
\begin{enumerate}
    \item We have canonical isomorphisms $$ PA(B) \cong PL(\mathbf{C}B) \cong \Gamma\left(\mathfrak{X}, \overline{\mathcal{M}}_{\mathfrak{X}}^{\gp}\right).$$
    \item The natural map $$\divisor: \ PA(B) \cong  \Gamma\left(\mathfrak{X}, \overline{\mathcal{M}}_{\mathfrak{X}}^{\gp}\right) \overset{\tau}{\to} \Pic(\mathfrak{X}) \xhookrightarrow{} \Cl(\mathfrak{X})$$ is of the form 
    \begin{equation} \label{divmap}
    \divisor: \alpha \mapsto \sum_{v \in \mathscr{P}^{[0]}} \alpha(v) D_v.
    \end{equation}
\end{enumerate}
\end{fact}

\begin{proof} 
We first explain that $PA(B) \cong PL(\mathbf{C}B)$. We have a natural inclusion $i: PA(B) \xhookrightarrow{} PL(\mathbf{C}B)$ by extending an $\alpha \in PA(B)$ to a PL-function on $\mathbf{C}B$ as follows. For any $\mathbf{C}\sigma \in \mathbf{C}\mathscr{P}^{[k]}$ and any $(a, x) \in \mathbf{C}\sigma = \mathbb{R}_{\geq 0} \cdot (\sigma \times \{1\}) \subseteq \mathbb{R} \times \mathbb{R}^k$, we define $i(\alpha)(a, x) := a \cdot \alpha(x)$. Then clearly, $i(\alpha)$ is a PL-function on $\mathbf{C}B$. Conversely, a PL-function $\beta$ on $\mathbf{C}B$ induces a PA-function $\beta|_{B}$ on $B$ by restriction. 

The isomorphism $PL(\mathbf{C}B) \cong \Gamma\left(\mathfrak{X}, \overline{\mathcal{M}}_{\mathfrak{X}}^{\gp}\right)$ can be viewed as a global version of \cite[Lemma 5.15]{GS2} which applies in the context of toric degenerations. 
Let $\mathbf{C} \sigma \in \mathbf{C}\mathscr{P}$ be a cell of $\mathbf{C}B$ (here we include the minimal cell $0 \in \mathbf{C}\mathscr{P}^{[0]}$) and let  $\eta_{\sigma} \in \mathfrak{X}$ be the generic point of the corresponding logarithmic stratum $X_{\mathbf{C}\sigma} \subseteq \mathfrak{X}$. Then specifying an integral PL-function on $\mathbf{C}\sigma = \Hom\left(\overline{\mathcal{M}}_{\mathfrak{X}, \eta_{\sigma}}, \mathbb{R}_{\geq 0}\right)$ is equivalent to specifying an element of $\Hom \left( \Hom\left(\overline{\mathcal{M}}_{\mathfrak{X}, \eta_{\sigma}}, \mathbb{N} \right), \mathbb{Z} \right).$ But under the assumptions on the log structure, we have canonical isomorphisms: 
\begin{multline} \label{dualisopl}
\Hom \left( \Hom\left(\overline{\mathcal{M}}_{\mathfrak{X}, \eta_{\sigma}}, \mathbb{N} \right), \mathbb{Z} \right) \cong \Hom \left( \Hom\left(\overline{\mathcal{M}}_{\mathfrak{X}, \eta_{\sigma}}, \mathbb{N} \right)^{\gp}, \mathbb{Z} \right) \cong \\ \cong \Hom \left( \Hom\left(\overline{\mathcal{M}}_{\mathfrak{X}, \eta_{\sigma}}, \mathbb{Z} \right), \mathbb{Z} \right) \cong \Hom \left( \Hom\left(\overline{\mathcal{M}}_{\mathfrak{X}, \eta_{\sigma}}^{\gp}, \mathbb{Z} \right), \mathbb{Z} \right) \cong \overline{\mathcal{M}}_{\mathfrak{X}, \eta_{\sigma}}^{\gp}
\end{multline}
where the second isomorphism is since $\overline{\mathcal{M}}_{\mathfrak{X}, \eta_{\sigma}}$ is sharp and the last isomorphism is since $\overline{\mathcal{M}}_{\mathfrak{X}, \eta_{\sigma}}$ is toric so $\overline{\mathcal{M}}_{\mathfrak{X}, \eta_{\sigma}}^{\gp}$ is torsion-free. 

Now, an element $\beta \in PL(\mathbf{C}B)$ specifies linear functions $\beta_{\mathbf{C}\sigma}, \ \mathbf{C} \sigma \in \mathbf{C}\mathscr{P}$ on the cones of $\mathbf{C}B$ such that for any $\mathbf{C} \tau = \mathbf{C} \sigma \cap \mathbf{C} \sigma'$ with $\mathbf{C} \sigma, \mathbf{C} \sigma', \mathbf{C} \tau \in \mathbf{C}\mathscr{P}$  we have $\beta_{\mathbf{C} \tau} = \beta_{\mathbf{C} \sigma}|_{\mathbf{C} \tau} = \beta_{\mathbf{C} \sigma'}|_{\mathbf{C} \tau}$. This implies that the corresponding stalks of $\overline{\mathcal{M}}_{\mathfrak{X}}^{\gp}$ agree under generization maps (well-defined since $\mathcal{M}_{\mathfrak{X}}$ is Zariski) $$\overline{\mathcal{M}}_{\mathfrak{X}, \eta_{\sigma}}^{\gp} \xrightarrow{} \overline{\mathcal{M}}_{\mathfrak{X}, \eta_{\tau}}^{\gp} \xleftarrow{} \overline{\mathcal{M}}_{\mathfrak{X}, \eta_{\sigma'}}^{\gp}$$ and hence give rise to a global section of $\overline{\mathcal{M}}_{\mathfrak{X}}^{\gp}$.

Conversely, evaluating a section $s$ of $\overline{\mathcal{M}}_{\mathfrak{X}}^{\gp}$ on all $u \in \Hom\left(\overline{\mathcal{M}}_{\mathfrak{X}, \eta_{\sigma}}, \mathbb{N} \right) \subseteq \Hom\left(\overline{\mathcal{M}}_{\mathfrak{X}, \eta_{\sigma}}, \mathbb{Z} \right)$ specifies a linear function $\beta_{\mathbf{C}\sigma}$ on $\mathbf{C} \sigma \in \mathbf{C}\mathscr{P}$. For any $\mathbf{C} \tau = \mathbf{C} \sigma \cap \mathbf{C} \sigma'$ with $\mathbf{C} \sigma, \mathbf{C} \sigma', \mathbf{C} \tau \in \mathbf{C}\mathscr{P}$  we have $\beta_{\mathbf{C} \tau} = \beta_{\mathbf{C} \sigma}|_{\mathbf{C} \tau} = \beta_{\mathbf{C} \sigma'}|_{\mathbf{C} \tau}$ since the stalks $s_{\sigma} \in \overline{\mathcal{M}}_{\mathfrak{X}, \eta_{\sigma}}^{\gp}$ of $s$ that we use to evaluate the section agree under generization. So the linear functions $\beta_{\mathbf{C}\sigma}, \ \mathbf{C} \sigma \in \mathbf{C}\mathscr{P}$ glue to a $\beta \in PL(\mathbf{C}B)$. This shows (1).

From the isomorphisms $PA(B) \cong PL(\mathbf{C}B) \cong \Gamma\left(\mathfrak{X}, \overline{\mathcal{M}}_{\mathfrak{X}}^{\gp}\right)$ we see that the section $s \in \Gamma\left(\mathfrak{X}, \overline{\mathcal{M}}_{\mathfrak{X}}^{\gp}\right)$ corresponding to $i(\alpha) \in PL(\mathbf{C}B)$ has stalks $s_{v} = \alpha(v) \in \overline{\mathcal{M}}_{\mathfrak{X}, \eta_v}^{\gp} \cong \mathbb{Z}$ at the generic points $\eta_v$ of $D_v, \ v \in \mathscr{P}^{[0]}$ (and a stalk $s_0 = 0 \in \overline{\mathcal{M}}_{\mathfrak{X}, \eta_0}^{\gp} = 0$ at the generic point $\eta_0$ of the locus of triviality of the log structure $\mathcal{M}_{\mathfrak{X}}$). Now the constructions of $\tau: \Gamma\left(\mathfrak{X}, \overline{\mathcal{M}}_{\mathfrak{X}}^{\gp}\right) \to \Pic(\mathfrak{X})$ and $\Pic(\mathfrak{X}) \xhookrightarrow{} \Cl(\mathfrak{X})$ imply (2). 
\end{proof}

\begin{defn} \label{pagenerated}
We say that a (Weil) divisor $D'$ on $\mathfrak{X}$ is PA-generated if it is a Cartier divisor of the form $D = \divisor(\alpha)$ for some $\alpha \in PA(B)$. Note that \eqref{divmap} implies that all PA-generated divisors are supported on $D$. 
\end{defn}

\begin{remark} 
If the divisor $D$ is simple normal crossings, then we have $\Gamma\left(\mathfrak{X}, \overline{\mathcal{M}}_{\mathfrak{X}}^{\gp}\right) \cong \Div_D(\mathfrak{X})$ (see Construction \ref{tropicalizesnc}). So every $D'$ supported on $D$ is PA-generated.
\end{remark}

\begin{cor} \label{DisPLgen}
In the assumptions of Proposition \ref{logpicard}, the divisor $D$ is PA-generated. 
\end{cor}

\begin{proof}
This is an immediate consequence of Proposition \ref{logpicard} by taking the constant function $1 \in PA(B)$.
\end{proof}

We are ready to generalize Propositions \ref{P1} and \ref{P2}.

\begin{fact} \label{P1general}
    For any contraction $\mathfrak{X} \overset{\pi'}{\to} \mathfrak{X}'$ (not necessarily to a toric degeneration), let $P$ be a well-chosen monoid and $K \subseteq P$ be the face containing the contracted curves. Suppose that there exists a $\pi'$-ample PA-generated divisor $D'$ on $\mathfrak{X}$. Then $ P \setminus K \in \mathcal{A}$.
\end{fact}

\begin{proof}
The proof is the same as that of Proposition \ref{P1} using Lemma \ref{puncturedordersgeneral} instead of Corollary \ref{puncturedorders}.
\end{proof}

For any monoid ideal $J \subseteq P$, we have a natural grading of $R_J$ (as defined in \eqref{thetaringmodI}) by $\Hom\left(\Gamma\left(\mathfrak{X}, \overline{\mathcal{M}}_{\mathfrak{X}}^{\gp}\right), \mathbb{Z}\right)$ via the maps
$$
\mathbf{C}B(\mathbb{Z}) \to \Hom\left(\Gamma\left(\mathfrak{X}, \overline{\mathcal{M}}_{\mathfrak{X}}^{\gp}\right), \mathbb{Z}\right), \quad p \mapsto (s \mapsto \langle s, p \rangle)
$$
(using Notation \ref{logpairing}) and 
\begin{equation} \label{degonP}
P \to \Hom\left(\Gamma\left(\mathfrak{X}, \overline{\mathcal{M}}_{\mathfrak{X}}^{\gp}\right), \mathbb{Z}\right), \quad \beta \mapsto \left(s \mapsto  \deg_{\mathcal{L}_s}(\beta)\right).
\end{equation}
Additivity of the degrees follows from Lemma \ref{puncturedordersgeneral}.

\begin{remark}
Note that we have defined a grading by $\Hom\left(\Gamma\left(\mathfrak{X}, \overline{\mathcal{M}}_{\mathfrak{X}}^{\gp}\right), \mathbb{Z}\right)$ rather than by $\Gamma\left(\mathfrak{X}, \overline{\mathcal{M}}_{\mathfrak{X}}^{\gp}\right)$. The reason is that the $\Hom\left(\Gamma\left(\mathfrak{X}, \overline{\mathcal{M}}_{\mathfrak{X}}^{\gp}\right), \mathbb{Z}\right)$-grading we defined falls into the general framework for gradings on $R_J$ (corresponding to a torus action on the intrinsic mirror modulo $J$) of \cite[Section 4.4]{GHS}. This is the natural way to define gradings (in particular, it guarantees the functoriality of the torus action).

\begin{sloppypar}
According to \cite[Section 4.4]{GHS}, to define a grading by an arbitrary finitely generated free abelian group $\Gamma$, it is enough to define maps $\delta_B: \Hom(PA(B), \mathbb{Z}) \to \Gamma$ and $\delta_P : P \to \Gamma$ that fit into a commutative diagram \cite[(4.8)]{GHS}\footnote{Unlike \cite[(4.8)]{GHS}, we are using $PA(B) \cong PL(\mathbf{C}B)$ instead of $PL(B)$ since we are in the projective case, see \cite[Remark 4.4.4]{GHS}.}:
    \begin{equation} \label{commtorus}
        \begin{tikzcd}
        Q_0 \arrow{r}{g} \arrow{d}{h} & \Hom(PA(B), \mathbb{Z})  \arrow{d}{\delta_B} \\
         P \arrow{r}{\delta_P} & \Gamma
        \end{tikzcd}
    \end{equation}
Here $Q_0:=\Hom(MPA(B, \mathbb{N}), \mathbb{N})$, $g$ is the natural map $$g: Q_0 \xhookrightarrow{} \Hom(MPA(B, \mathbb{Z}), \mathbb{Z}) \to \Hom(PA(B), \mathbb{Z})$$ and $h: Q_0 \to P$ is the map defined by the MPA function $\varphi$ (see \cite[Proposition 1.2.9(b)]{GHS}). One also needs to check that the wall functions $f_{\mathfrak{p}}$ of $\mathfrak{D}_J$ are homogeneous of degree $0$ with respect to the resulting $\Gamma$-degree. 
\end{sloppypar}
    
We can define the $\Hom\left(\Gamma\left(\mathfrak{X}, \overline{\mathcal{M}}_{\mathfrak{X}}^{\gp}\right), \mathbb{Z}\right)$-grading by taking $\delta_B$ to be the dual of the canonical isomorphism $PA(B) \cong \Gamma\left(\mathfrak{X}, \overline{\mathcal{M}}_{\mathfrak{X}}^{\gp}\right)$ and taking $\delta_P$ to be the map of \eqref{degonP}. It is straightforward to check that \eqref{commtorus} is commutative and the induced $\Hom\left(\Gamma\left(\mathfrak{X}, \overline{\mathcal{M}}_{\mathfrak{X}}^{\gp}\right), \mathbb{Z}\right)$-grading is the same grading we defined above. The fact that the wall functions $f_{\mathfrak{p}}$ of $\mathfrak{D}_J$ are homogeneous of degree $0$ follows from the analogue of Lemma \ref{puncturedordersgeneral} for $f \in \mathscr{M}(\mathfrak{X}, \btau)$ (see \cite[Proposition 1.13]{GS3}).

In the case that $D$ is simple normal crossings (so that $\Gamma\left(\mathfrak{X}, \overline{\mathcal{M}}_{\mathfrak{X}}^{\gp}\right) \cong \Div_D(\mathfrak{X})$, see Construction \ref{tropicalizesnc}), our grading recovers the grading by $\Hom(PL(\mathbf{C}B), \mathbb{Z})$ of \cite[Example 4.4.6]{GHS} and \cite[Construction 3.16]{GS4} (using the canonical isomorphism $PL(\mathbf{C}B) \cong \Gamma\left(\mathfrak{X}, \overline{\mathcal{M}}_{\mathfrak{X}}^{\gp}\right)$). Moreover, the irreducible components $D_i, \ 1 \leq i \leq m$ of $\mathfrak{X}_0$ define a natural basis of $\Div_D(\mathfrak{X})$ that gives rise to the canonical isomorphism $$\Div_D(\mathfrak{X}) \xrightarrow{\cong} \Hom(\Div_D(\mathfrak{X}), \mathbb{Z}), \quad D_i \to D^*_i.$$ Under this isomorphism, our grading is just the grading by $\Div_D(\mathfrak{X})$ from the previous section.
\end{remark}

Suppose that the rank of $\Gamma\left(\mathfrak{X}, \overline{\mathcal{M}}_{\mathfrak{X}}^{\gp}\right)$ is $l$ and let $s_i \in \Gamma\left(\mathfrak{X}, \overline{\mathcal{M}}_{\mathfrak{X}}\right) \subseteq \Gamma\left(\mathfrak{X}, \overline{\mathcal{M}}_{\mathfrak{X}}^{\gp}\right)$, $ 1 \leq i \leq l$ be any linearly independent sections (not necessarily forming a $\mathbb{Z}$-basis), such that $\langle s_1, \cdots, s_l \rangle_{\mathbb{Z}} \otimes_{\mathbb{Z}} \mathbb{Q} = \Gamma\left(\mathfrak{X}, \overline{\mathcal{M}}_{\mathfrak{X}}^{\gp}\right) \otimes_{\mathbb{Z}} \mathbb{Q}$. Let $(e_1, \dots, e_{l})$ be a basis for $\mathbb{Z}^{l}$. We define a character map 
\begin{equation} \label{charactermapgeneral}
w_s: \mathbf{C}B(\mathbb{Z}) \times P \setminus J \to \mathbb{Z}^{l}, \quad (p, \beta) \mapsto \sum_{i=1}^{l} \langle s_i, p \rangle e_{i} + \sum_{i=1}^l \deg_{\mathcal{L}_{s_i}}(\beta) e_{i}.
\end{equation}

\begin{fact} \label{P2general}
In the setup of Proposition \ref{P1general}, suppose in addition that:
\begin{enumerate}
    \item The divisor $D'$ is simple normal crossings and effective.
    \item The divisor $D_{\irrel}$ (defined as in the proof of Proposition \ref{P2}) is $\mathbb{Q}$-Cartier and $\pi'$-nef.
    \item $K \cap NE(\mathfrak{X}_0) = K \cap NE(\mathfrak{X}_0)_{\numer}$ (under the splittings of \eqref{standardsplittings}), i.e. $\pi'$ only contracts numerical classes of curves $C \subseteq \mathfrak{X}_0$.
\end{enumerate} Then we have $I_{\min} \subseteq P \setminus K$.
\end{fact}

\begin{proof}
The proof that $(\mathfrak{X}, \Delta)$ is klt is as in Proposition \ref{P2} using assumptions (1) and (2), Remark \ref{irreldivisor}, and the fact that $D'$ is PA-generated, thus Cartier. We continue as in the proof of Proposition \ref{P2} (using assumption (3)) until we need to use the character map (Step 3). 

Let $s_1 \in \Gamma\left(\mathfrak{X}, \overline{\mathcal{M}}^{\gp}_{\mathfrak{X}}\right)$ be the global section corresponding to $D'$. The fact that $D'$ is effective (by condition (1)) implies that $s_1 \in \Gamma\left(\mathfrak{X}, \overline{\mathcal{M}}_{\mathfrak{X}}\right) \subseteq \Gamma\left(\mathfrak{X}, \overline{\mathcal{M}}^{\gp}_{\mathfrak{X}}\right)$. Indeed, by Proposition \ref{logpicard}(2), the PL-function $\beta \in PL(\mathbf{C}B)$ corresponding to $D'$ is non-negative on $\mathbf{C}B$. Following the proof of Proposition \ref{logpicard}(1), $\beta$ specifies a compatible collection $\beta_{\mathbf{C}\sigma}, \ \mathbf{C} \sigma \in \mathbf{C}\mathscr{P}$ of non-negative linear functions on the cones of $\mathbf{C}B$. Non-negativity of $\beta_{\mathbf{C}\sigma}$ implies that it corresponds to a stalk of $\overline{\mathcal{M}}_{\mathfrak{X}, \eta_{\sigma}} \subseteq \overline{\mathcal{M}}_{\mathfrak{X}, \eta_{\sigma}}^{\gp}$ under the isomorphism of \eqref{dualisopl}. The collection of such stalks over all the generic points $\eta_{\sigma} \in \mathfrak{X}$ of logarithmic strata $X_{\mathbf{C}\sigma} \subseteq \mathfrak{X}$ is compatible under generization maps and defines the global section $s_1 \in \Gamma\left(\mathfrak{X}, \overline{\mathcal{M}}_{\mathfrak{X}}\right) \subseteq \Gamma\left(\mathfrak{X}, \overline{\mathcal{M}}^{\gp}_{\mathfrak{X}}\right)$. 

Choose any $s_{2}, \dots, s_{l} \in \Gamma\left(\mathfrak{X}, \overline{\mathcal{M}}_{\mathfrak{X}}\right) \subseteq \Gamma\left(\mathfrak{X}, \overline{\mathcal{M}}_{\mathfrak{X}}^{\gp}\right)$ with $\langle s_1, \cdots, s_l \rangle_{\mathbb{Z}} \otimes_{\mathbb{Z}} \mathbb{Q} = \Gamma\left(\mathfrak{X}, \overline{\mathcal{M}}_{\mathfrak{X}}^{\gp}\right) \otimes_{\mathbb{Z}} \mathbb{Q}$ (note that it is not always possible to choose $s_{2}, \dots, s_{l} \in \Gamma\left(\mathfrak{X}, \overline{\mathcal{M}}_{\mathfrak{X}}\right)$, such that $\langle s_1, \cdots, s_l \rangle_{\mathbb{Z}} \cong \Gamma\left(\mathfrak{X}, \overline{\mathcal{M}}_{\mathfrak{X}}^{\gp} \right)$) and let $(e_1, \dots, e_{l})$ be a basis for $\mathbb{Z}^{l}$.
We now have the character map $w_s$ of \eqref{charactermapgeneral}.

Similarly to the proof of Proposition \ref{P2} (using duality and the fact that $\langle s_1, \cdots, s_l \rangle_{\mathbb{Z}} \otimes_{\mathbb{Z}} \mathbb{Q} = \Gamma\left(\mathfrak{X}, \overline{\mathcal{M}}_{\mathfrak{X}}^{\gp}\right) \otimes_{\mathbb{Z}} \mathbb{Q}$), it is enough to check that $w_s$ has finite fibres. As in the proof of Proposition \ref{P2}, it suffices to show that $\ker w_s \cap (\mathbf{C}\sigma(\mathbb{Z}) \times K) = P^{\times}$ for every $\sigma \in \mathscr{P}^{\max}$. Suppose that $(r, \beta) \in \ker w_s \cap (\mathbf{C}\sigma(\mathbb{Z}) \times K)$ and let $D'_i$ be the divisor corresponding to $\mathcal{L}_{s_i}$ for all $2 \leq i \leq l$). We have  $$w_s(r, \beta) = \sum_{i=1}^l \langle s_i, r \rangle e_i + \left(D' \cdot \beta \right) e_1 + \sum_{i=2}^l (D_i' \cdot \beta) e_{i} = 0$$ Note that $\langle s_i, r \rangle \geq 0$ for all $1 \leq i \leq l$ since $s_i \in \Gamma\left(\mathfrak{X}, \overline{\mathcal{M}}_{\mathfrak{X}}\right)$ by the above. So $D' \cdot \beta \leq 0$ and $D_i' \cdot \beta \leq 0$ for all $2 \leq i \leq l$. But $D'$ is positive on $K \setminus P^{\times}$, so $\beta \in P^{\times}$. But then 
$D' \cdot \beta = 0$ and $D_i' \cdot \beta = 0$ for all $2 \leq i \leq l$, so we have $\sum_{i=1}^l \langle s_i, r \rangle e_i = 0$. So $\langle s_i, r \rangle = 0$ for all $1 \leq i \leq l$. Since the log structure $\mathcal{M}_{\mathfrak{X}}$ is fine and saturated (in particular, $\Gamma\left(\mathfrak{X}, \overline{\mathcal{M}}_{\mathfrak{X}}^{\gp}\right)$ is torsion-free) and $\langle s_1, \cdots, s_l \rangle_{\mathbb{Z}} \otimes_{\mathbb{Z}} \mathbb{Q} = \Gamma\left(\mathfrak{X}, \overline{\mathcal{M}}_{\mathfrak{X}}^{\gp}\right) \otimes_{\mathbb{Z}} \mathbb{Q}$, for every section $s \in \Gamma\left(\mathfrak{X}, \overline{\mathcal{M}}_{\mathfrak{X}}^{\gp}\right) \setminus \{0\}$ we have $\alpha s \in \langle s_1, \cdots, s_l \rangle_{\mathbb{Z}}$ for some $\alpha \in \mathbb{Z} \setminus \{0\}$. So $\alpha \langle s, r \rangle = \langle \alpha s , r \rangle = 0$ and $\langle s, r \rangle = 0$. Therefore, we have $\langle s, r \rangle = 0$ for all $s \in \Gamma\left(\mathfrak{X}, \overline{\mathcal{M}}_{\mathfrak{X}}^{\gp}\right)$. This implies that $r=0$, so $\ker w_s \cap (\mathbf{C}\sigma(\mathbb{Z}) \times K) = P^{\times}$.
\end{proof}

\subsubsection{Extension corresponding to the resolution \texorpdfstring{$\pi: \mathfrak{X} \to \bar{\mathfrak{X}}$}{pi}} \label{extendingintmirror}

Going back to the setup at the start of this section, suppose that the resolution $\pi: \mathfrak{X} \to \bar{\mathfrak{X}}$ of a toric degeneration $\bar{\mathfrak{X}} \to \mathcal{S}$ to a log smooth minimal log CY degeneration $\mathfrak{X} \to \mathcal{S}$ satisfies the assumptions of Proposition \ref{P2general} (or the assumptions of Proposition \ref{P2} if $D$ is simple normal crossings). We want to use $\pi: \mathfrak{X} \to \bar{\mathfrak{X}}$ to extend the intrinsic mirror $\check{\mathfrak{X}} \to \Spec \widehat{\kk [P]}$.

Let $I_{\min}$ be the ideal of Proposition \ref{imin} for $\mathfrak{X} \to \mathcal{S}$ and let $J \subseteq P$ be any radical ideal with $I_{\min} \subseteq J$. For any $J'$ with $\sqrt{J'} = J$ we have $J' \in \mathcal{A}$ by Proposition \ref{imin}. So the intrinsic mirror is defined over $J'$. Taking the inverse limit over all the $J'$  with $\sqrt{J'} = J$ we obtain the \emph{extended intrinsic mirror} 
\begin{equation} \label{extendedmirror}
\check{\mathfrak{X}} \to \Spec \widehat{\kk [P]}_J
\end{equation}
where we denote by $\widehat{\kk [P]}_J$ the completion of $\kk [P]$ with respect to $J$.\footnote{ For example, if $P=\mathbb{N}^k$ and $J=\mathbb{N}^{k-l}$ for some $1 \leq l < k$, then $\widehat{\kk [P]}_{J} = \kk [t_1, \dots, t_l]\llbracket t_{l+1}, \dots, t_k \rrbracket$ where the formal elements correspond to the generators of $J$.} As in Remark \ref{fewerideals}, it is actually enough to take the limit over the families $$\check{\mathfrak{X}}_{J^{k+1}} \to \Spec \kk [P]/J^{k+1}$$ for $ k \geq 0$.

Proposition \ref{P2general} (or Proposition \ref{P2} if $D$ is simple normal crossings) implies that for a well-chosen monoid $P$ and a face $K \subseteq P$ containing the contracted curve classes, we may take the ideal $J$ to be $J:=P \setminus K$ (which is a prime, thus radical, ideal of $P$) in which case \eqref{extendedmirror} is a family over the completion of $\Spec \kk[K]$. We shall only use the extension \eqref{extendedmirror} in this context.

\begin{remark} \label{extmonoidchange}
We explain a slightly different way to view the extension \eqref{extendedmirror} that will be useful in Chapter \ref{gslocus}. For any well-chosen monoid $P$ and an extension using $J:= P \setminus K$, the intrinsic mirror is well-defined over the monoid $P + K^{\gp}$. Note that $P + K^{\gp}$ does not satisfy condition (3) of Definition \ref{monoidass} on the base monoid. As a result, the maximal ideal $\mathfrak{m}_{\ext}$ of $P + K^{\gp}$ has complement $K^{\gp}$ which is not finite. However, we have $\mathfrak{m}_{\ext} = J$ as subsets of $A_1(\mathfrak{X}_0, \mathbb{Z})$ which implies that the intrinsic mirror $\check{\mathfrak{X}} \to \Spec \widehat{\kk[P + K^{\gp}]}$ is well-defined via \eqref{prodtheta} and \eqref{sumalpha}. This mirror may also be viewed as the basechange of the extended intrinsic mirror of \eqref{extendedmirror} via the natural inclusion $P \xhookrightarrow{} P + K^{\gp}$. We still call $\check{\mathfrak{X}} \to \Spec \widehat{\kk[P + K^{\gp}]}$ the \emph{extended intrinsic mirror}.  
\end{remark}

\begin{remark} \label{scatforext}
Even though one can use $\pi: \mathfrak{X} \to \bar{\mathfrak{X}}$ to extend the intrinsic mirror of Section \ref{3pointmirror}, it is not obvious that one can give a scattering diagram interpretation for the extended intrinsic mirror. Indeed, the recipe of Construction \ref{canonical scattering} produces infinitely many walls (even for $I = J$). By grouping certain walls together, we will provide a scattering diagram interpretation of $\check{\mathfrak{X}} \to \Spec \widehat{\kk [P]}_J$ in Section \ref{scatint}. Our analysis will generalize \cite[Section 5.3]{GHKS}, which covers the case when the resolution $\pi: \mathfrak{X} \to \bar{\mathfrak{X}}$ factors as a composition of a small contraction and a log étale blowup (e.g. $\pi: \mathfrak{X} \to \bar{\mathfrak{X}}$ is small). On one hand, we will require our resolution to be of a specific form (i.e. toric, integral, and homogeneous, see Definition \ref{toricresolutions}). On the other hand, we won't need to require the existence of a $\pi$-ample PA-generated divisor $D'$ (or assumptions (1) and (2) of Proposition \ref{P2general}), see Remark \ref{stronglyadmissibleweakening}(1).
\end{remark}

\subsection{Overview of the results} \label{problem}

We will now give a brief overview of the results that we are going to prove in Chapters \ref{proof} and \ref{gslocus}, using Example \ref{runningexample} as a working model. First, we set up the stage for proving Conjecture \ref{conjecture} by summarizing the setups of Sections \ref{tdsetup}, \ref{intsetup}, and \ref{extresults}. 

\begin{basicsetup} \label{setupforproof}
Let $\bar{\mathfrak{X}} \overset{\bar{g}}{\to} \mathcal{S}$ be a special toric degeneration with polarization $A$ and let $\pi: \mathfrak{X} \to \bar{\mathfrak{X}}$ be a projective resolution to a minimal log CY degeneration $\mathfrak{X} \overset{g}{\to} \mathcal{S}$. We assume that $\bar{\mathfrak{X}} \overset{\bar{g}}{\to} \mathcal{S}$ satisfies Assumption \ref{toricdegenerationass} and $\mathfrak{X} \overset{g}{\to} \mathcal{S}$ satisfies Assumption \ref{manifoldass}.

Let $P$ be a well-chosen monoid (see Definition \ref{wellchosendef}) with $K \subseteq P$ the face containing the classes of the contracted curves, and let  
$J: = P \setminus K$. We assume that $\pi: \mathfrak{X} \to \bar{\mathfrak{X}}$ satisfies the assumptions of Proposition \ref{P2general} (or Proposition \ref{P2} if $D$ is simple normal crossings). Fix also a choice of the initial slab functions $$\left\{f_{\underline{\rho}} \in \kk[\Lambda_{\rho}] \ \Big| \ \underline{\rho} \in \tilde{\bar{\mathscr{P}}}^{[n-1]} \right\}$$  of \eqref{initialslabs} for the toric degeneration that defines a toric log CY structure on $\check{\bar{X}}_0$. 

We have defined the algorithmic scattering diagram $\bar{\mathfrak{D}} = \left\{ \bar{\mathfrak{D}}_k, \ k \geq 0 \right\}$ in Theorem \ref{reconstruction} and defined the toric degeneration mirror $$\check{\bar{\mathfrak{X}}} \to \Spec \kk \llbracket t \rrbracket$$ to $\bar{\mathfrak{X}} \to \mathcal{S}$ as the limit over $\check{\bar{\mathfrak{X}}}_{\bar{\mathfrak{D}}_k} \to \Spec \kk[t]/(t^{k+1})$ for $k \geq 0$ in \eqref{reconstructedfamily}. We also defined the extended intrinsic mirror $$\check{\mathfrak{X}} \to \Spec \widehat{\kk [P]}_J$$ to $\mathfrak{X} \to \mathcal{S}$ as the limit over $\check{\mathfrak{X}}_{J^{k+1}} \to \Spec \kk [P]/J^{k+1}$ for $ k \geq 0$ in \eqref{extendedmirror}. Note that the basechange (of Conjecture \ref{conjecture}) of the extended intrinsic mirror by $P \to \mathbb{N}, \ \beta \mapsto \pi^{*} A \cdot \beta$ makes sense.
\end{basicsetup}

Chapter \ref{proof} is devoted to the proof of Conjecture \ref{conjecture} for toric degenerations of K3-s. Consider the toric degeneration $\bar{\mathfrak{X}} \to \Spec \kk \llbracket t \rrbracket$ of Example \ref{runningexample} (which clearly satisfies Assumption \ref{toricdegenerationass}). Recall that we have $$\bar{\mathfrak{X}}:= \left\{tf_4+x_0x_1x_2x_3\right\} \subseteq \mathbb{P}^3 \times \Spec \kk\llbracket t \rrbracket$$ and the generic fibre of $\bar{\mathfrak{X}} \to \Spec \kk \llbracket t \rrbracket$ is smooth for a general choice of $f_4$. Moreover, genericity of $f_4$ implies, by Lemma \ref{localmodelsk3} and Proposition \ref{k3divisorial}, that $\bar{\mathfrak{X}} \to \Spec \kk \llbracket t \rrbracket$ is a divisorial log deformation. So it is special by Proposition \ref{whenspecialk3}. 

The dual intersection complex $\left(\bar B, \bar{\mathscr{P}}\right)$ of $\bar{\mathfrak{X}} \to \Spec \kk \llbracket t \rrbracket$ is the boundary of a tetrahedron with all faces standard triangles, see Example \ref{exaffine}. The singularities $x_{\rho} \in \Int(\rho)$ of the affine structure are at the irrational points of the edges $\rho \in \bar{\mathscr{P}}^{[1]}$ and have monodromy index $r_{\rho} = 4$. Every $x_{\rho}, \ \rho \in \bar{\mathscr{P}}^{[1]}$ corresponds to $4$ ordinary double point singularities on $\bar{X}_{\rho}$ with local models $\left\{xy = tw_{\rho} \right\} \subseteq \Spec \kk[x,y,w_{\rho}]\llbracket t \rrbracket$. They can be resolved by blowing up the components $\bar{D}_1, \bar{D}_2, \bar{D}_3, \bar{D}_4$ of the central fibre $\bar{\mathfrak{X}}_0$ in any order. Suppose that we blow up $\bar{D}_i$ before $\bar{D}_{i+1}$. We obtain a small resolution $\pi: \mathfrak{X} \to \bar{\mathfrak{X}}$ to a minimal log CY degeneration $\mathfrak{X} \to \Spec \kk \llbracket t \rrbracket$. 

The dual intersection complex $(B, \mathscr{P})$ of $\mathfrak{X} \to \Spec \kk \llbracket t \rrbracket$ is isomorphic to $\left(\bar B, \bar{\mathscr{P}}\right)$ as a polyhedral manifold but has a different affine structure with singularities contained at the vertices. The irreducible components $D_i, \ 1 \leq i \leq 4$ of the central fibre $\mathfrak{X}_0$ of $\mathfrak{X} \to \Spec \kk \llbracket t \rrbracket$ are the strict transforms of $\bar{D}_i, \ 1 \leq i \leq 4$ and the divisor $D = D_1 + D_2 + D_3 + D_4$ is clearly simple normal crossings. Since $(B, \mathscr{P}) \cong \left(\bar B, \bar{\mathscr{P}}\right)$ as polyhedral manifolds and $\left(\bar B, \bar{\mathscr{P}}\right)$ satisfies Assumption \ref{manifoldass}, so does $(B, \mathscr{P})$. 

A divisor of the form $D' = a_1 D_1 + a_2 D_2 + a_3 D_3 + a_4 D_4$ for some $a_i \in \mathbb{Z}$ is $\pi$-ample as long as $a_1 < a_2 < a_3 < a_4$. Further assuming $a_1 > 0$, the degeneration $\mathfrak{X} \to \mathcal{S}$, the resolution $\pi: \mathfrak{X} \to \bar{\mathfrak{X}}$, and the effective $\pi$-ample divisor $D'$ satisfy the assumptions of Proposition \ref{P2}. The exceptional locus of $\pi: \mathfrak{X} \to \bar{\mathfrak{X}}$ is the union of $24$ disjoint curves 
\begin{equation} \label{24curves}
\left\{E_{\rho, k} \ | \ 1 \leq k \leq 4, \ \rho \in \bar{\mathscr{P}}^{[1]} \right\},
\end{equation}
one for each singularity of $\bar{\mathfrak{X}}_0$. Let $P$ be a well-chosen monoid (obtained via Proposition \ref{wellchosen}) with a face $K$ containing the contracted curves. 

In the case that all the maximal cells $\sigma \in \bar{\mathscr{P}}^{\max}$ of the dual intersection complex $\left(\bar{B}, \bar{\mathscr{P}} \right)$ of a special toric degeneration $\mathfrak{X} \to \mathcal{S}$ of K3-s are standard triangles, we construct a resolution $\pi: \mathfrak{X} \to \bar{\mathfrak{X}}$ (and prove that it has the required properties) to a minimal log CY degeneration $\mathfrak{X} \to \mathcal{S}$ in Section \ref{resolvesmall} by blowing up a sequence of irreducible components of the central fibre. Generalizing Example \ref{runningexample}, we actually construct a small resolution $\pi: \mathfrak{X} \to \bar{\mathfrak{X}}$ in this case. We generalize the construction to a larger class of $\sigma \in \bar{\mathscr{P}}^{\max}$ in Section \ref{resolve}. To perform the construction in general and have more control over the resolution, we need a tropical interpretation. We give this interpretation in Section \ref{admissible}, introducing admissible resolutions in Definition \ref{admissibledef} and showing that any special toric degeneration $\bar{\mathfrak{X}} \to \mathcal{S}$ of K3-s admits a (strongly) admissible resolution $\pi: \mathfrak{X} \to \bar{\mathfrak{X}}$ in Proposition \ref{strongadmissibleresolutionresult}.
 We define a piecewise-linear (PL) isomorphism $\Phi: B \to \bar B$ in these three cases in Constructions \ref{plisosmall}, \ref{plmap}, and \ref{plmapadmissible} respectively.

In \eqref{choiceslab}, we define the initial slab functions (of the form required by \eqref{slabfunctions}) for the toric degeneration $\bar{\mathfrak{X}} \to \mathcal{S}$. In the case of Example \ref{runningexample}, for $\underline{\rho}, \underline{\rho}' \in \tilde{\bar{\mathscr{P}}}^{[1]}$ two slabs with $\underline{\rho}, \underline{\rho}' \subseteq \rho \in \bar{\mathscr{P}}^{[1]}$ and $w_{\rho}$ a choice of coordinate on $\bar{X}_{\rho}$, we have 
\begin{equation} \label{exinitslabfunc}
f_{\underline{\rho}} = (1+ w_{\rho})^4, \quad f_{\underline{\rho}'} = (1+ w_{\rho}^{-1})^4.
\end{equation}
We are now in Basic Setup \ref{setupforproof} for proving Conjecture \ref{conjecture}. 

In Section \ref{scatint}, we give a scattering diagram interpretation of the extended intrinsic mirror
\begin{equation} \label{over1}
\check{\mathfrak{X}} \to \Spec \widehat{\kk [P]}_J.    
\end{equation}

\begin{sloppypar}
Namely, in Constructions \ref{newwalls} and \ref{extendcanonicalscat}, we produce a collection $\mathfrak{D}^J: = \left\{ \mathfrak{D}_{J^{k+1}}, \ k \geq 0  \right\}$ of compatible scattering diagrams such that every $\mathfrak{D}_{J^{k+1}}$ is consistent modulo $J^{k+1}$ and $\check{\mathfrak{X}}_{\mathfrak{D}_{J^{k+1}}} \to \Spec \kk [P]/J^{k+1}$ is isomorphic to the family $\check{\mathfrak{X}}_{J^{k+1}} \to \Spec \kk [P]/J^{k+1}$ obtained by reducing \eqref{over1} modulo $J^{k+1}$. This involves computing the walls modulo $J$ using the recipe of Construction \ref{canonical scattering}, see Proposition \ref{diag}. There are infinitely many such walls, but one can replace all the walls with the same support with one wall, with the wall function the infinite product of the wall functions (that turns out to be polynomial). This defines the scattering diagram $\mathfrak{D}_J$. In the case of Example \ref{runningexample}, $\mathfrak{D}_J$ consists of six slabs $\mathfrak{b}_{\rho}$ supported on the edges $\rho \in \mathscr{P}^{[1]}$, with attached functions $$f_{\mathfrak{b}_{\rho}} = \prod_{k=1}^4 (1+ t^{E_{\rho, k}}w_{\rho}).$$ Note that these look similar to the initial slab function $f_{\underline{\rho}}$ of \eqref{exinitslabfunc} for the toric degeneration.   
\end{sloppypar}

We prove Conjecture \ref{conjecture} in Section \ref{main} by relating $\mathfrak{D}^J$ and $\bar{\mathfrak{D}}$. The basechange of Conjecture \ref{conjecture} corresponds to interpreting scattering diagrams $\mathfrak{D}_{J^{k+1}}, \ k \geq 0$ as scattering diagrams $\mathfrak{D}_{k}, \ k \geq 0$ with base monoid $\mathbb{N}$, see Construction \ref{ccc}. In the case of Example \ref{runningexample}, the basechange just sets all the $E_{\rho, k}$ of \eqref{24curves} to $0$. Now, using the PL-isomorphism $\Phi: B \to \bar B$, in Construction \ref{chs} we define the image of
the scattering diagram $\mathfrak{D}_k$ on $(B, \mathscr{P})$ as a scattering diagram $\Phi(\mathfrak{D}_k)$ on $\left(\bar B, \bar{\mathscr{P}}\right)$. The image of every codimension $0$ wall $\mathfrak{p} \in \mathfrak{D}_k$ that is not a slab is a codimension $0$ wall $\Phi(\mathfrak{p}) \in \Phi(\mathfrak{D}_k)$ with the same attached wall function. Defining the image of a slab $\mathfrak{b} \subseteq \rho \in \mathscr{P}^{[1]}$ is more complicated. The image is either one or two slabs depending on whether $\Phi(\mathfrak{b}) \subseteq \Phi(\rho)$ contains the singularity $x_{\Phi(\rho)}$ of the affine structure on $\left(\bar{B}, \bar{\mathscr{P}}\right)$. We also need to adjust the slab functions to account for the monodromy around the singularities.

By construction, $\Phi(\mathfrak{D}_0)$ is (combinatorially) equivalent to  $\bar{\mathfrak{D}}_0$. We show that $\Phi(\mathfrak{D}_k)$ is consistent for all $k \geq 0$ in Proposition \ref{consis} and that $\Phi(\mathfrak{D}_k)$ is equivalent to $\mathfrak{D}_k$ in Proposition \ref{equivalence}. We finish the argument in Theorem \ref{proofconj}. The uniqueness of the reconstruction algorithm of Theorem \ref{reconstruction} implies that $\Phi(\mathfrak{D}_k)$ is equivalent to $\bar{\mathfrak{D}}_k$ for all $k \geq 0$. So Conjecture \ref{conjecture} follows by Proposition \ref{help}.

In Chapter \ref{gslocus}, we vastly generalize the conjecture. Thinking again of Example \ref{runningexample} (or any situation where the resolution is small) and $A_1(\mathfrak{X}_0, \mathbb{Z}):= A_1(\mathfrak{X}_0, \mathbb{Z})_{\numer}$, let us explain the heuristics, treating the monoid $P$ as $NE(\mathfrak{X}_0)$ (see Observation \ref{heuristics}(2)). Since the resolution $\pi: \mathfrak{X} \to \bar{\mathfrak{X}}$ is small, it induces a splitting of $f_{*}:NE(\mathfrak{X}_0) \to  NE\left(\bar{\mathfrak{X}}_0\right)$ by sending every curve $C \in \bar{\mathfrak{X}}_0$ to the scheme-theoretic preimage $f^{-1}(C)$. So we have 
\begin{equation} \label{canonicalsplitting}
NE(\mathfrak{X}_0) = K \oplus \pi^{-1} NE\left(\bar{\mathfrak{X}}_0\right),
\end{equation}
inducing a canonical isomorphism 
\begin{equation} \label{canonicalsplittingring}
\widehat{\kk [NE(\mathfrak{X}_0)]} \cong \kk [K] \llbracket NE\left(\bar{\mathfrak{X}}_0\right) \rrbracket = \kk[t^{E_{\rho, k}}]\llbracket NE\left(\bar{\mathfrak{X}}_0\right) \rrbracket
\end{equation}
(where the notation means that the completion is only with respect to the second factor). Now the restriction of $\Spec \kk[t^{E_{\rho, k}}]$ to $\left\{t^{E_{\rho, k}} \neq 0, \ 1 \leq k \leq 4, \ \rho \in \bar{\mathscr{P}}\right\}$ can be seen as the parameter space of toric log CY structures on $\check{\bar{X}}_0$ of a certain form (we will make this more precise in Proposition \ref{fibres} and Remark \ref{fibresremarks}(2)). Here the log structure on $\check{\bar{X}}_0$ is induced from the natural log structure on $\check{\mathfrak{X}} \to \Spec \kk[t^{E_{\rho, k}}]\llbracket NE\left(\bar{\mathfrak{X}}_0\right) \rrbracket$ of Appendix \ref{log} by inclusion. By analogy with \cite{GHK}, we call the family $\check{\mathfrak{X}} \to \Spec \kk[t^{\pm E_{\rho, k}}]\llbracket NE\left(\bar{\mathfrak{X}}_0\right) \rrbracket$ the \emph{(numerical) minimal relative Gross-Siebert locus} (see Definition \ref{gslocusdef}). We formalize these heuristics using Noetherian families in Section \ref{gensetup}. 

Since the fibres over the points with $\left\{t^{E_{\rho, k}} \neq 0, \ 1 \leq k \leq 4, \ \rho \in \bar{\mathscr{P}}\right\}$ are toric log CY, they can be deformed to toric degenerations. We recover these deformations from the minimal relative Gross-Siebert locus, performing the construction universally in polarization $A$, choices of the initial slab functions, and in a family of gluing data for toric degenerations. Note that the fibres over the points of $\Spec \kk[t^{E_{\rho, k}}]\llbracket NE\left(\bar{\mathfrak{X}}_0\right) \rrbracket$ that are not in the base of the minimal relative Gross-Siebert locus are not toric log CY (so such fibres can't define central fibres of toric degenerations) since the singularities of the log structure fall into deeper strata, see the discussion of Section \ref{gslocusstuff}.   

Figure \ref{pictureextentions} illustrates the various extension results. Theorem \ref{proofconj} constructs a one-dimensional family, specified by the polarization $A$, through the point where $t^{E_{\rho, k}} =1$ for all $E_{\rho,k}, \ \rho \in \bar{\mathscr{P}}^{[1]}, \ 1 \leq k \leq r_{\rho}$ (since $\pi^{*} A \cdot E_{\rho, k} = 0$). It corresponds to the black line in Figure \ref{pictureextentions}. In Section \ref{universal}, we use the universal version of the toric degeneration mirror (see \cite[Theorem A.2.4]{GHS}) to remove the dependence on $A$. We recover the toric degeneration mirror over a certain universal finitely generated monoid containing $NE\left(\bar{\mathfrak{X}}_0\right)$. This mirror corresponds to the red sphere in Figure \ref{pictureextentions}, and we prove the result in Proposition \ref{proofuniversal}.

\begin{figure}[t]
\begin{center}
\begin{tikzpicture}
\node [
    above right,
    inner sep=0] (image) at (0,0) {\includegraphics[width=\textwidth]{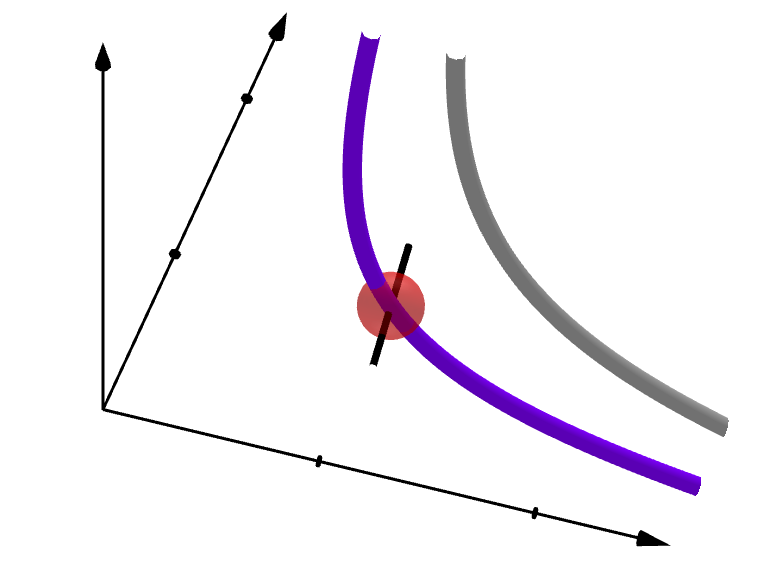}};
 
\begin{scope}[
x={($0.1*(image.south east)$)},
y={($0.1*(image.north west)$)}]

    \node at (1.5,9.5){\small $t^{E}, \ E \in NE\left(\bar{\mathfrak{X}}_0\right)$};
    \node at (8.85,0.7){\small $t^{E_{\rho, k}}$};
    \node at (4.05,9.5){\small $t^{E_{\rho, k}}$};
    \node at (3.3,4){\small \boxed{\text{Theorem \ref{compgen}}}};
 
    \draw[latex-, very thick,dashed] (4.85, 3.65) -- (3,1.5) 
        node[left,black,fill=white]{\small \boxed{\text{Theorem \ref{proofconj}}}};

    \draw[latex-, very thick,dashed] 
        (5.65,4.7) -- (8,4.7)
         node[right,black,fill=white]{\small \boxed{\text{Proposition \ref{proofuniversal}}}};
         
    \draw[latex-, very thick,dashed] 
        (4.85,6.5) -- (8,6.5)
         node[right,black,fill=white]{\small \boxed{\text{Proposition \ref{prooffree}}}};
    
    \draw[latex-, very thick,dashed] 
        (6.25,8) -- (8,8)
         node[right,black,fill=white]{\small \boxed{\text{Proposition \ref{proofglue}}}};
\end{scope}
 
\end{tikzpicture}
\end{center}
\caption{Extensions of Conjecture \ref{conjecture}.} \label{pictureextentions}
\end{figure}
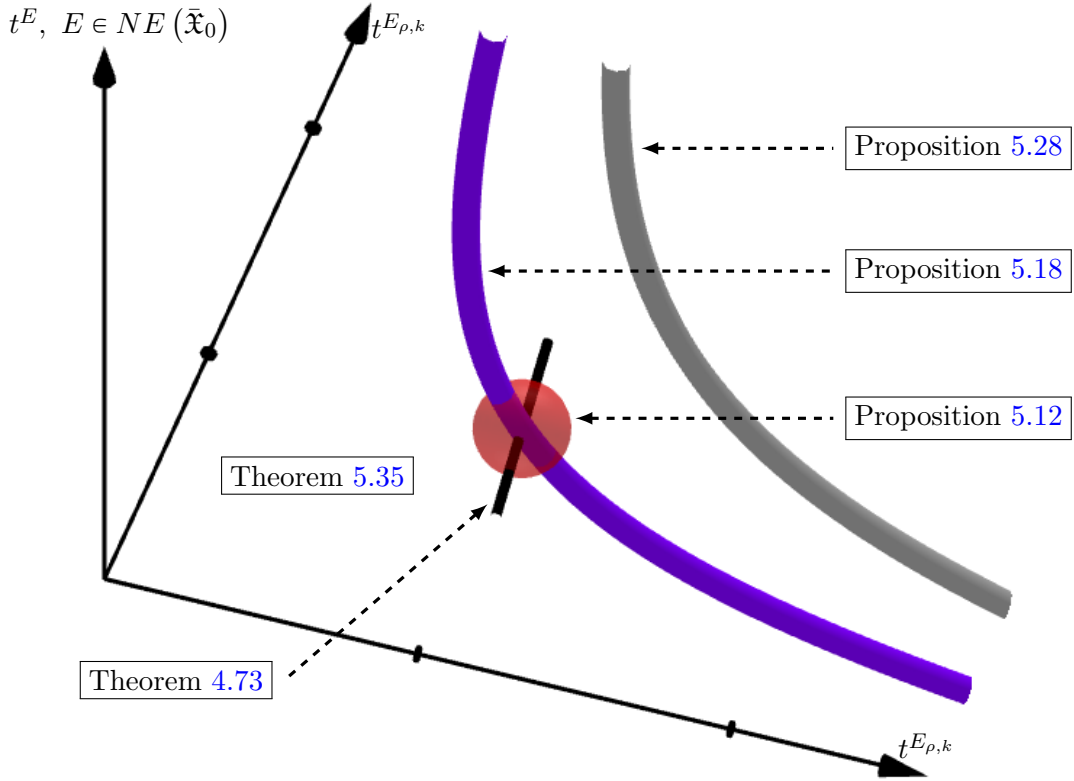

In Section \ref{freeparam}, we vary the point of $\kk[t^{E_{\rho, k}}]$ to recover the toric degeneration mirrors constructed using other choices of the initial slab functions satisfying Proposition \ref{struct}. Since the local rigidity condition is empty in dimension $2$, one can interpret the free coefficients $a_{\rho, i}$ in the slab functions as new variables and construct the toric degeneration mirror using the ring $A := \kk [a_{\rho, i}] = \kk [\mathbb{N}^K]$ for $K := \sum_{\rho \in \bar{\mathscr{P}}^{[1]}, r_{\rho} > 0} (r_{\rho} -1)$ instead of $A = \kk$  (see \cite[Theorem A.4.2]{GHS}). We obtain a correspondence between this toric degeneration mirror and the restriction of the extended intrinsic mirror to the family over the subvariety $\left\{\prod_{k=1}^{r_{\rho}} t^{E_{\rho, k}} = 1, \ \rho \in \bar{\mathscr{P}}^{[1]} \right\}$ of $\Spec \widehat{\kk [NE(\mathfrak{X}_0)]} = \Spec \kk[t^{E_{\rho, k}}]\llbracket NE\left(\bar{\mathfrak{X}}_0\right) \rrbracket$. This mirror corresponds to the blue tubular (corresponding to the formal terms) neighbourhood in Figure \ref{pictureextentions}, and we prove the result in Proposition \ref{prooffree}. Note that this result is more general than Proposition \ref{proofuniversal} only if the dual intersection complex $\left( \bar{B}, \bar{\mathscr{P}}\right)$ of $\bar{\mathfrak{X}} \to \mathcal{S}$ is not simple.

Finally, in Section \ref{glue}, we study the fibres of the extended intrinsic mirror $\check{\mathfrak{X}} \to \Spec \kk[t^{E_{\rho, k}}]\llbracket NE\left(\bar{\mathfrak{X}}_0\right) \rrbracket$ over points with arbitrary choices of $\prod_{k=1}^4 t^{E_{\rho, k}} \in  \kk^{\times}$ for $\rho \in \bar{\mathscr{P}}^{[1]}$ (such that this system of equations is consistent\footnote{In general, one can have relations between the $E_{\rho, k}$, see the discussion after \eqref{egp}.}). These produce slab functions of the form not covered by Proposition \ref{struct} (i.e. non-normalized), so understanding them requires introducing gluing data for toric degenerations into the picture. We generalize the correspondence of Proposition \ref{prooffree} to arbitrary fixed $\left\{\prod_{k=1}^{r_{\rho}} t^{E_{\rho, k}} \in \kk^{\times}, \ \rho \in \bar{\mathscr{P}}^{[1]} \right\}$ in Proposition \ref{proofglue}, corresponding to the grey tubular neighbourhood in Figure \ref{pictureextentions}. In Theorem \ref{compgen}, we obtain the result in complete generality, showing a correspondence between the minimal relative Gross-Siebert locus and a certain subfamily in the gluing data of the universal (varied in the free parameters of the initial slab functions and in gluing data) toric degeneration mirror of \cite[Theorem A.4.2]{GHS}.

\section[Proof of Conjecture 1.7 for toric degenerations of K3-s]{Proof of Conjecture \ref{conjecture} for toric degenerations of K3-s} \label{proof}

In this chapter, we prove Conjecture \ref{conjecture} following the plan of Section \ref{problem}.   

\subsection[Small resolution in the generically snc case]{Small resolution in the generically simple normal crossings case} \label{resolvesmall}

Let $\bar{\mathfrak{X}} \to \mathcal{S}$ be a special toric degeneration of K3-s (satisfying Assumption \ref{toricdegenerationass}) with polarization $A$ and dual intersection complex $\left(\bar B, \bar{\mathscr{P}}\right)$. As in Section \ref{tdaff}, let $\bar{D}_i, \ 1 \leq i \leq \bar{m}$ be the irreducible components of the central fibre $\bar{\mathfrak{X}}_0$. We endow $\bar{\mathfrak{X}}$ with the divisorial log structure given by $\bar{D}:= \bar{\mathfrak{X}}_0 = \bar{D}_1 + \dots + \bar{D}_{\bar{m}}$. We also endow $\mathcal{S}$ with the divisorial log structure given by $0 \in \mathcal{S}$ making $\bar{\mathfrak{X}} \to \mathcal{S}$ a log morphism. By Proposition \ref{whenspecialk3}, being special is equivalent to $\bar{\mathfrak{X}} \to \mathcal{S}$ being a divisorial log deformation with a smooth generic fibre. 

In this section, we assume in addition that all the maximal cells $\sigma \in \bar{\mathscr{P}}^{\max}$ are standard triangles. This assumption corresponds to requiring that $\bar D$ is simple normal crossings away from a finite number of singularities at the double curves $\bar{X}_{\rho}, \ \rho \in \bar{\mathscr{P}}^{[1]}$. In particular, Example \ref{runningexample} satisfies this requirement. We are going to construct a small resolution $\pi: \mathfrak{X} \to \bar{\mathfrak{X}}$ to a minimal log CY degeneration $\mathfrak{X} \to \mathcal{S}$ by blowing up a sequence of irreducible components of the central fibre $\bar{\mathfrak{X}}_0$. We shall generalize the construction to a larger class of $\sigma \in \bar{\mathscr{P}}^{\max}$ in Section \ref{resolve} and treat the case of general $\sigma \in \bar{\mathscr{P}}^{\max}$ in Section \ref{admissible}.

\subsubsection{Resolving local models} \label{locmod}

By Observation \ref{modelslocal}(2), a toric degeneration $\bar{\mathfrak{X}} \to \mathcal{S}$ that is a divisorial log deformation with a smooth generic fibre has local models of the form $\left\{xy = t^l \right\} \subseteq \Spec \kk[x,y,w_{\rho}]\llbracket t \rrbracket$ and $\left\{xy = t^lw_{\rho} \right\} \subseteq \Spec \kk[x,y,w_{\rho}]\llbracket t \rrbracket$ (with the natural map to $\Spec \kk\llbracket t \rrbracket$) at the non-singular and singular points of the double curve $\bar{X}_{\rho}, \rho \in \bar{\mathscr{P}}^{[1]}$ respectively, where $l$ is the integral length of $\rho$. Because of the assumption on the maximal cells, we have $l = 1$ in all the local models. 

The local model $\left\{xy = t \right\} \subseteq \Spec \kk[x,y,w_{\rho}]\llbracket t \rrbracket$ is simple normal crossings, and the map to $\Spec \kk\llbracket t \rrbracket$ is log smooth (here the induced divisorial log structure in the model is given by $\left\{t=0\right\}$). So we only need to resolve the log singularities with local models of the form $\left\{xy = tw_{\rho} \right\} \subseteq \Spec \kk[x,y,w_{\rho}]\llbracket t \rrbracket$, which are \emph{ordinary double point} singularities (ODP-s).

The local model $\left\{xy = tw_{\rho} \right\} \subseteq \Spec \kk[x,y,w_{\rho}]\llbracket t \rrbracket$ can be resolved to a simple normal crossings and log smooth degeneration by blowing up one of the divisors in the central fibre $\left\{t=0\right\}$, that is by blowing up $\{x=t=0\}$ or $\{y=t=0\}$. Torically, $\left\{xy = tw_{\rho} \right\} \subseteq \Spec \kk[x,y,w_{\rho}]\llbracket t \rrbracket$ is (the completion in $t$ of) the affine variety defined by the cone over the square, that is the convex hull of $(0, 0, 1), (1, 0, 1), (0, 1, 1), (1, 1, 1)$. Vertices of the square correspond to toric divisors, and the two blowups correspond to the two ways to subdivide the square, see Figure \ref{coneoversquare}.

\begin{figure}[H]
\begin{center}

\begin{minipage}{.3\textwidth}
\begin{tikzpicture}[line cap=round,line join=round,>=triangle 45,x=1cm,y=1cm, scale = 3]
\clip(-0.2,-0.2) rectangle (1.25,1.25);
\draw  (0,1)-- (1,1);
\draw  (1,1)-- (1,0);
\draw  (0,0)-- (1,0);
\draw  (0,0)-- (0,1);
\draw  (0,0)-- (1,1);
\draw [fill] (0,0) circle (0.5pt);
\draw [fill] (0,1) circle (0.5pt);
\draw [fill] (1,0) circle (0.5pt);
\draw [fill] (1,1) circle (0.5pt);
\draw(0.1, -0.1) node{\scriptsize $x=t=0$};
\draw(0.88, -0.1) node{\scriptsize $y=t=0$};
\draw(0.1, 1.1) node{\scriptsize $x=w_{\rho}=0$};
\draw(0.88, 1.1) node{\scriptsize $y=w_{\rho}=0$};

\end{tikzpicture}
\captionsetup{labelformat=empty}
\caption{Blowup of $x=t=0$}
\end{minipage}
\begin{minipage}{.3\textwidth}
\begin{tikzpicture}[line cap=round,line join=round,>=triangle 45,x=1cm,y=1cm, scale = 3]
\clip(-0.2,-0.2) rectangle (1.25,1.25);
\draw  (0,1)-- (1,1);
\draw  (1,1)-- (1,0);
\draw  (0,0)-- (1,0);
\draw  (0,0)-- (0,1);
\draw [fill] (0,0) circle (0.5pt);
\draw [fill] (0,1) circle (0.5pt);
\draw [fill] (1,0) circle (0.5pt);
\draw [fill] (1,1) circle (0.5pt);
\draw(0.1, -0.1) node{\scriptsize $x=t=0$};
\draw(0.88, -0.1) node{\scriptsize $y=t=0$};
\draw(0.1, 1.1) node{\scriptsize $x=w_{\rho}=0$};
\draw(0.88, 1.1) node{\scriptsize $y=w_{\rho}=0$};

\end{tikzpicture}
\captionsetup{labelformat=empty}
\caption{\ \ \ ODP \\ singularity}
\end{minipage}
\begin{minipage}{.3\textwidth}
\begin{tikzpicture}[line cap=round,line join=round,>=triangle 45,x=1cm,y=1cm, scale = 3]
\clip(-0.2,-0.2) rectangle (1.25,1.25);
\draw  (0,1)-- (1,1);
\draw  (1,1)-- (1,0);
\draw  (0,0)-- (1,0);
\draw  (0,0)-- (0,1);
\draw  (0,1)-- (1,0);
\draw [fill] (0,0) circle (0.5pt);
\draw [fill] (0,1) circle (0.5pt);
\draw [fill] (1,0) circle (0.5pt);
\draw [fill] (1,1) circle (0.5pt);
\draw(0.1, -0.1) node{\scriptsize $x=t=0$};
\draw(0.88, -0.1) node{\scriptsize $y=t=0$};
\draw(0.1, 1.1) node{\scriptsize $x=w_{\rho}=0$};
\draw(0.88, 1.1) node{\scriptsize $y=w_{\rho}=0$};
\end{tikzpicture}
\captionsetup{labelformat=empty}
\caption{Blowup of $y=t=0$}
\end{minipage}
\end{center}
\setcounter{figure}{0}
\caption{Two resolutions of an ODP.} \label{coneoversquare}
\end{figure}

\begin{figure}[t]
\begin{center}

\begin{tikzpicture}[line cap=round,line join=round,>=triangle 45,x=1cm,y=1cm]

\begin{scope}[scale = 2]
 \clip(-1.5068464964661663,-0.6875061121455801) rectangle (1.5626116514137933,1.7468917292764545);
\draw (0,1)-- (0,0);
\draw (0,1)-- (0.5,1.5);
\draw (0,1)-- (-0.5,1.5);
\draw (0,0)-- (-0.5,-0.5);
\draw (0,0)-- (0.5,-0.5);
\begin{scriptsize}

\draw [color=red, line width=1pt] (0,0.5)-- ++(-6pt,6pt);
\draw(-0.25, 0.6) node{$E$};

\end{scriptsize}
\end{scope}

\begin{scope}[xshift=4.5cm, scale = 2]
 \clip(-1.5068464964661663,-0.6875061121455801) rectangle (1.5626116514137933,1.7468917292764545);
\draw (0,1)-- (0,0);
\draw (0,1)-- (0.5,1.5);
\draw (0,1)-- (-0.5,1.5);
\draw (0,0)-- (-0.5,-0.5);
\draw (0,0)-- (0.5,-0.5);
\begin{scriptsize}
\draw [color=blue, line width = 1pt] (0, 0.5)-- ++(-1.5pt,-1.5pt) -- ++(3pt,3pt) ++(-3pt,0) -- ++(3pt,-3pt);
\draw(0.5, 0.3) node{$y=t=0$};
\draw(-0.5, 0.3) node{$x=t=0$};
\end{scriptsize}
\end{scope}

\begin{scope}[xshift=9cm, scale = 2]
 \clip(-1.5068464964661663,-0.6875061121455801) rectangle (1.5626116514137933,1.7468917292764545);
\draw (0,1)-- (0,0);
\draw (0,1)-- (0.5,1.5);
\draw (0,1)-- (-0.5,1.5);
\draw (0,0)-- (-0.5,-0.5);
\draw (0,0)-- (0.5,-0.5);
\begin{scriptsize}

\draw [color=red, line width=1pt] (0,0.5)-- ++(6pt,6pt);
\draw(0.25, 0.6) node{$E$};

\end{scriptsize}
\end{scope}

\draw[->] (7.5,1.3) -- (6,1.3);
\draw[->] (1.5,1.3) -- (3,1.3);

\end{tikzpicture}

\end{center}
\caption{Two blowups of an ODP.} \label{odpresolve}
\end{figure}
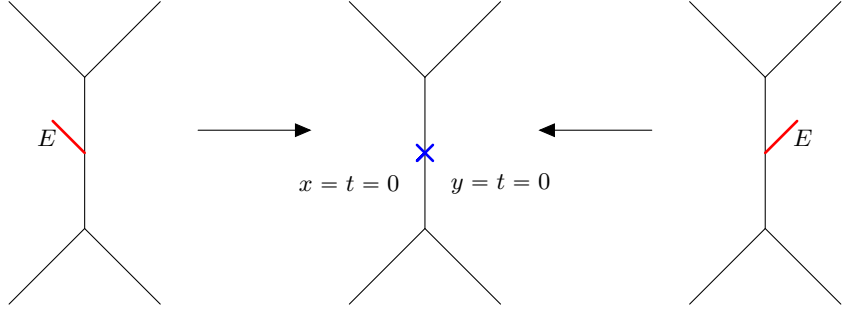

The exceptional locus of the blowup is a single curve $E$ meeting $\left\{x=y=t=0\right\}$ at one point and contained in the divisor that was blown up. See Figure \ref{odpresolve} below for a sketch of the central fibres of the local model and the resolutions corresponding to Figure \ref{coneoversquare}.

\subsubsection{Blowing up a divisor} \label{divisorblowup}

We want to obtain a global resolution of $\bar{\mathfrak{X}} \to \mathcal{S}$ by blowing up a sequence of irreducible components of the central fibre $\bar{\mathfrak{X}}_0$. Let $\bar{D}_i$ for some $ 1 \leq i \leq \bar m$ be a component of $\bar{\mathfrak{X}}_0$. Then blowing up $\bar{D}_i$ corresponds to blowing up either $\{x=t=0\}$ or $\{y=t=0\}$ in the local models of all the singular points contained in every double curve $\bar{X}_{\rho} \subseteq \bar{D}_i$ (with a compatible choice of local models for the singular points contained in the same $\bar{X}_{\rho}, \rho \in \bar{\mathscr{P}}^{[1]}$). Note that the number of singular points in $\bar{X}_{\rho}, \rho \in \bar{\mathscr{P}}^{[1]}$ is equal to $r_{\rho}$ (the index of the singularity $x_{\rho} \subseteq \Int(\rho)$ of the affine structure). The exceptional locus is a disjoint union of exceptional curves contained in the strict transform of $\bar{D}_i$, intersecting the adjacent divisors at one point. Every exceptional curve arises as the curve $E$ in Section \ref{locmod} using the local model at the corresponding singular point.

\begin{notations} \label{blowupnotations}
We introduce the following notations:
\begin{enumerate}
    \item We denote the strict transform of $\bar{D}_i$ (after blowing it up) by $D_i$. We denote the strict transforms of the other divisors $\bar{D}_k$ (for $1 \leq k \leq \bar{m}$ and $k \neq i$) by $\bar{D}_k$.
    \item We denote the exceptional curves contained in $D_i$ and intersecting $\bar{D}_j$ by $E_{k}^{ij}$ for $1 \leq k \leq r_{\rho}$. Alternatively, if the divisor that the curve is contained in is not important, we will use the notation $E_{\rho, k}$ for $\bar{X}_{\rho} := \bar{D}_i \cap \bar{D}_j$. Note that $E_{k}^{ij}$ does not intersect any irreducible components of the central fibre of the resolution apart from $D_i$ and $\bar{D}_j$.
\end{enumerate}
We will use similar notations for further blowups of irreducible components of (the proper transform of) the central fibre $\bar{\mathfrak{X}}_0$ of $\bar{\mathfrak{X}} \to \mathcal{S}$. That is, we shall denote the divisors that have been blown up without a bar and use the notation $E_{k}^{ij}$ (or $E_{\rho, k}$) for the exceptional curves of further blowups.
\end{notations}

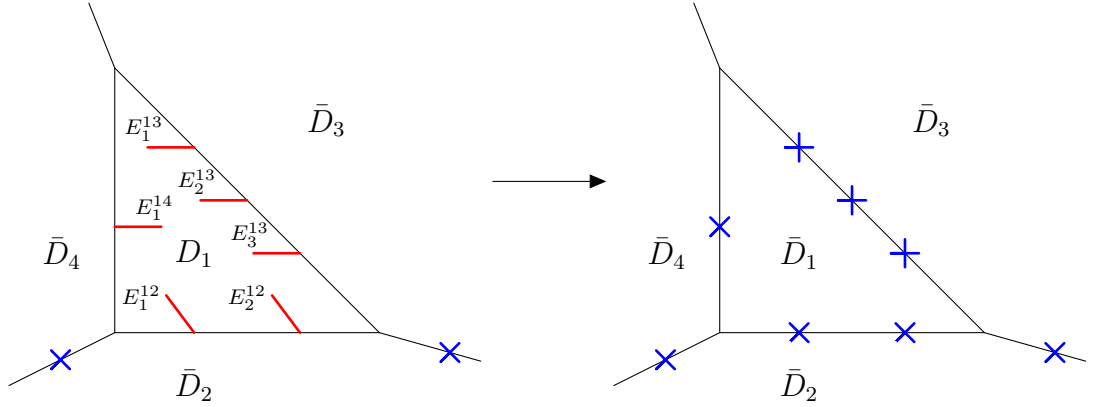
\begin{figure}[t]
\begin{center}
\hspace{-1cm}
\begin{tikzpicture}[line cap=round,line join=round,>=triangle 45,x=1cm,y=1cm]

\begin{scope}[scale = 3.5]
  \clip(-0.6422904971084227,-0.38192239585650195) rectangle (1.5701289445169364,1.3727550923291256);
\draw (0,1)-- (1,0);
\draw (0,1)-- (0,0);
\draw (0,0)-- (1,0);
\draw (-0.4,-0.2)-- (0,0);
\draw (-0.09786989626372086,1.2449016467603005)-- (0,1);
\draw  (1,0)-- (1.3830218834654358,-0.10779353114155361);
\begin{scriptsize}

\draw [color=blue, line width = 1pt] (-0.20489333052568517,-0.10244666526284259)-- ++(-1pt,-1pt) -- ++(2pt,2pt) ++(-2pt,0) -- ++(2pt,-2pt);
\draw [color=blue, line width = 1pt] (1.2664045376373911,-0.074974008180028)-- ++(-1pt,-1pt) -- ++(2pt,2pt) ++(-2pt,0) -- ++(2pt,-2pt);
\draw [color=red, line width=1pt] (0.3,0)-- ++(-3pt,4pt);
\draw [color=red, line width=1pt] (0.7,0)-- ++(-3pt,4pt);
\draw [color=red, line width=1pt] (0,0.4)-- ++(5pt,0pt);
\draw [color=red, line width=1pt] (0.7,0.3)-- ++(-5pt,0pt);
\draw [color=red, line width=1pt] (0.3,0.7)-- ++(-5pt,0pt);
\draw [color=red, line width=1pt] (0.5,0.5)-- ++(-5pt,0pt);
\draw(0.3, 0.3) node{\normalsize $D_1$};
\draw(0.3, -0.2) node{\normalsize $\bar{D}_2$};
\draw(0.8, 0.8) node{\normalsize $\bar{D}_3$};
\draw(-0.2, 0.3) node{\normalsize $\bar{D}_4$};
\draw(0.1, 0.13) node{\tiny $E^{12}_{1}$};
\draw(0.5, 0.13) node{\tiny $E^{12}_{ 2}$};
\draw(0.15, 0.47) node{\tiny $E^{14}_{ 1}$};
\draw(0.51, 0.37) node{\tiny $E^{13}_{ 3}$};
\draw(0.31, 0.57) node{\tiny $E^{13}_{ 2}$};
\draw(0.11, 0.77) node{\tiny $E^{13}_{ 1}$};
\end{scriptsize}
\end{scope}

\begin{scope}[xshift=8cm, scale = 3.5]
  \clip(-0.6422904971084227,-0.38192239585650195) rectangle (1.5701289445169364,1.3727550923291256);
\draw (0,1)-- (1,0);
\draw (0,1)-- (0,0);
\draw (0,0)-- (1,0);
\draw (-0.4,-0.2)-- (0,0);
\draw (-0.09786989626372086,1.2449016467603005)-- (0,1);
\draw  (1,0)-- (1.3830218834654358,-0.10779353114155361);
\begin{scriptsize}
\draw [color=blue, line width = 1pt] (0.3,0)-- ++(-1pt,-1pt) -- ++(2pt,2pt) ++(-2pt,0) -- ++(2pt,-2pt);
\draw [color=blue, line width = 1pt] (0.7,0)-- ++(-1pt,-1pt) -- ++(2pt,2pt) ++(-2pt,0) -- ++(2pt,-2pt);
\draw [color=blue, line width = 1pt] (0.7,0.3)-- ++(-1.5pt,0 pt) -- ++(3pt,0 pt) ++(-1.5pt,-1.5pt) -- ++(0 pt,3pt);
\draw [color=blue, line width = 1pt] (0.5,0.5)-- ++(-1.5pt,0 pt) -- ++(3pt,0 pt) ++(-1.5pt,-1.5pt) -- ++(0 pt,3pt);
\draw [color=blue, line width = 1pt] (0.3,0.7)-- ++(-1.5pt,0 pt) -- ++(3pt,0 pt) ++(-1.5pt,-1.5pt) -- ++(0 pt,3pt);
\draw [color=blue, line width = 1pt] (0,0.4)-- ++(-1pt,-1pt) -- ++(2pt,2pt) ++(-2pt,0) -- ++(2pt,-2pt);
\draw [color=blue, line width = 1pt] (-0.20489333052568517,-0.10244666526284259)-- ++(-1pt,-1pt) -- ++(2pt,2pt) ++(-2pt,0) -- ++(2pt,-2pt);
\draw [color=blue, line width = 1pt] (1.2664045376373911,-0.074974008180028)-- ++(-1pt,-1pt) -- ++(2pt,2pt) ++(-2pt,0) -- ++(2pt,-2pt);
\draw(0.3, 0.3) node{\normalsize $\bar{D}_1$};
\draw(0.3, -0.2) node{\normalsize $\bar{D}_2$};
\draw(0.8, 0.8) node{\normalsize $\bar{D}_3$};
\draw(-0.2, 0.3) node{\normalsize $\bar{D}_4$};
\end{scriptsize}
\end{scope}

\draw[->] (5,2) -- (6.5,2);

\end{tikzpicture}

\end{center}
\caption{Blowing up the divisor $\bar{D}_1$. The number of singularities on every $\bar{X}_{\rho} = \bar{D}_1 \cap \bar{D}_j, \ j=1,2,3$ corresponds to $r_{\rho}$.} \label{blowingadivisorsmall}
\end{figure}

Note that after blowing up $\bar{D}_i$, the deformation becomes log smooth and simple normal crossings in a neighbourhood of $D_i$. We give an example of a divisor blowup in Figure \ref{blowingadivisorsmall}.

We explain the effect of blowing up $\bar{D}_i$ on the dual intersection complex $\left( \bar{B}, \bar{\mathscr{P}} \right)$. Let $g': \mathfrak{X}' \to \mathcal{S}$ be the degeneration obtained after blowing up $\bar{D}_i$. 
Similarly to the toric degeneration and log smooth cases (see \eqref{dic} and \eqref{dicims} respectively), we define the dual intersection complex $\left(B', \mathscr{P}'\right)$ of $\mathfrak{X}' \overset{g'}{\to} \mathcal{S}$ as $(g')^{-1}_{\trop}(1)$ for $g'_{\trop}: \Sigma(\mathfrak{X}') \to \mathbb{R}_{\geq 0}$ the tropicalization of $g'$ of Construction \ref{tropicalize}. As usual, the polyhedral structure $\mathscr{P}'$ comes from restricting the cones of $\Sigma \left(\mathfrak{X}'\right)$ to the fibre over $1 \in \mathbb{R}_{\geq 0}$. The fact that all the maximal cells $\sigma \in \bar{\mathscr{P}}^{\max}$ are standard triangles implies that $\left(B', \mathscr{P}'\right) \cong \left(\bar{B}, \bar{\mathscr{P}}\right)$ as polyhedral manifolds. Indeed, the local toric model for a codimension $2$ stratum $\bar{X}_{\sigma}$ (with $\bar{v}_{i} \subseteq \sigma$) is given by the cone over $\sigma \in \bar{\mathscr{P}}^{\max}$. So it is smooth, and the blowup of the toric divisor corresponding to $\bar{D}_{i}$ is trivial in the local model. Therefore, blowing up $\bar{D}_i$ does not subdivide any cells of $\bar{\mathscr{P}}$. In particular, since $\left(\bar{B}, \bar{\mathscr{P}}\right)$ is a polyhedral manifold of dimension $2$ in the sense of Definition \ref{polymani} (see Proposition \ref{dicpm}), so is $\left(B', \mathscr{P}'\right)$. 

\begin{notation} \label{blowupdicnotations}
We denote by $v_i \in \mathscr{P}'^{[0]}$ the vertex corresponding to $D_i$ and by $\bar{v}_k \in \mathscr{P}'^{[0]}$ the vertices corresponding to $\bar{D}_k$ (for $1 \leq k \leq \bar{m}$ and $k \neq i$). We use similar notations for $\bar{\mathfrak{X}}_0$ and for further blowups of irreducible components of $\mathfrak{X}'_{0}$.
\end{notation}

Now, we can give $\left(B', \mathscr{P}'\right)$ the structure of an affine manifold with singularities using Construction \ref{affstrint} (viewed via Observation \ref{altaffinecons}, i.e. replacing the right-hand side of \eqref{constraintaffine} with the expression of Corollary \ref{compinsmallerstrata}(3)). Note from Remark \ref{toriccomponent}(2) that this defines the usual toric charts of Construction \ref{affstrtoric} on $W_{\bar{v}_k}$ for $1 \leq k \leq \bar{m}, \ k \neq i$ and the intrinsic affine structure on $W_{v_i} \setminus v_i$. Therefore, at the level of affine manifolds, the blowup of $\bar{D}_i$ can be visualized as pulling the singularities $x_{\rho} \in \Int(\rho)$ for $\rho \in \bar{\mathscr{P}}^{[1]}$ with $\bar{v}_i \subseteq \rho$ (and with $r_\rho \neq 0$ so that there is an $x_{\rho}$) into $\bar{v}_i$. We give the transformation corresponding to the blowup of Figure \ref{blowingadivisorsmall} in Figure \ref{blowingadivisorsmallaff}. Here and later, if we draw a singularity at an interior point $x_{\rho}$ of an edge $\rho$, we always assume that the point is not rational. 

\begin{figure}[t]
\begin{center}

\begin{tikzpicture}[line cap=round,line join=round,>=triangle 45,x=1cm,y=1cm]

\hspace{-1cm}
\begin{scope}[scale = 2]
\clip(-2.1,-1.7) rectangle (2.1,1.7);
\draw  (0,0)-- (-1,0);
\draw  (0,0)-- (0,-1);
\draw  (0,0)-- (1,1);
\draw  (-1,0)-- (1,1);
\draw  (1,1)-- (0,-1);
\draw  (-1,0)-- (0,-1);
\draw  (0,-1)-- (0,-1.5);
\draw  (1,1)-- (1.5,1.5);
\draw  (-1,0)-- (-1.5,0);
\draw [fill] (0,0) circle (1pt);
\draw [fill] (-1,0) circle (1pt);
\draw [fill] (0,-1) circle (1pt);
\draw [fill] (1,1) circle (1pt);
\begin{scriptsize}
\draw [color=blue, line width = 1pt] (0,0)-- ++(-1.5pt,-1.5pt) -- ++(3pt,3pt) ++(-3pt,0) -- ++(3pt,-3pt);
\draw [color=blue, line width = 1pt] (-0.5,-0.5)-- ++(-2.25pt,0 pt) -- ++(4.5pt,0 pt) ++(-2.25pt,-2.25pt) -- ++(0 pt,4.5pt);
\draw [color=blue, line width = 1pt] (0.5,0)-- ++(-2.25pt,0 pt) -- ++(4.5pt,0 pt) ++(-2.25pt,-2.25pt) -- ++(0 pt,4.5pt);
\draw [color=blue, line width = 1pt] (1.3,1.3)-- ++(-2.25pt,0 pt) -- ++(4.5pt,0 pt) ++(-2.25pt,-2.25pt) -- ++(0 pt,4.5pt);
\draw [color=blue, line width = 1pt] (0,-1.3)-- ++(-1.5pt,-1.5pt) -- ++(3pt,3pt) ++(-3pt,0) -- ++(3pt,-3pt);
\draw(0.2, 0) node{\normalsize $v_1$};
\draw(1.2, 0.9) node{\normalsize $\bar{v}_2$};
\draw(0.25, -1) node{\normalsize $\bar{v}_3$};
\draw(-1, 0.25) node{\normalsize $\bar{v}_4$};
\draw(0.2, 1.4) node{\normalsize $\left(B', \mathscr{P}'\right)$};
\end{scriptsize}
\end{scope}

\begin{scope}[xshift=8cm, scale = 2]
\clip(-2.1,-1.7) rectangle (2.1,1.7);
\draw  (0,0)-- (-1,0);
\draw  (0,0)-- (0,-1);
\draw  (0,0)-- (1,1);
\draw  (-1,0)-- (1,1);
\draw  (1,1)-- (0,-1);
\draw  (-1,0)-- (0,-1);
\draw  (0,-1)-- (0,-1.5);
\draw  (1,1)-- (1.5,1.5);
\draw  (-1,0)-- (-1.5,0);
\draw [fill] (0,0) circle (1pt);
\draw [fill] (-1,0) circle (1pt);
\draw [fill] (0,-1) circle (1pt);
\draw [fill] (1,1) circle (1pt);
\begin{scriptsize}
\draw [color=blue, line width = 1pt] (-0.5,0)-- ++(-1.5pt,-1.5pt) -- ++(3pt,3pt) ++(-3pt,0) -- ++(3pt,-3pt);
\draw [color=blue, line width = 1pt] (0.5,0.5)-- ++(-2.25pt,0 pt) -- ++(4.5pt,0 pt) ++(-2.25pt,-2.25pt) -- ++(0 pt,4.5pt);
\draw [color=blue, line width = 1pt] (0,-0.5)-- ++(-1.5pt,-1.5pt) -- ++(3pt,3pt) ++(-3pt,0) -- ++(3pt,-3pt);
\draw [color=blue, line width = 1pt] (-0.5,-0.5)-- ++(-2.25pt,0 pt) -- ++(4.5pt,0 pt) ++(-2.25pt,-2.25pt) -- ++(0 pt,4.5pt);
\draw [color=blue, line width = 1pt] (0.5,0)-- ++(-2.25pt,0 pt) -- ++(4.5pt,0 pt) ++(-2.25pt,-2.25pt) -- ++(0 pt,4.5pt);
\draw [color=blue, line width = 1pt] (1.3,1.3)-- ++(-2.25pt,0 pt) -- ++(4.5pt,0 pt) ++(-2.25pt,-2.25pt) -- ++(0 pt,4.5pt);
\draw [color=blue, line width = 1pt] (0,-1.3)-- ++(-1.5pt,-1.5pt) -- ++(3pt,3pt) ++(-3pt,0) -- ++(3pt,-3pt);
\draw(0.2, 0) node{\normalsize $\bar{v}_1$};
\draw(1.2, 0.9) node{\normalsize $\bar{v}_2$};
\draw(0.25, -1) node{\normalsize $\bar{v}_3$};
\draw(-1, 0.25) node{\normalsize $\bar{v}_4$};
\draw(0.2, 1.4) node{\normalsize $\left(\bar B, \bar{\mathscr{P}}\right)$};
\end{scriptsize}
\end{scope}

\draw[->] (2.5,0) -- (4,0);

\end{tikzpicture}
\end{center}
\caption{Transformation of $\left(\bar B, \bar{\mathscr{P}}\right)$ corresponding to blowing up $\bar{D}_1$.} \label{blowingadivisorsmallaff}
\end{figure}

Blowing up the component $\bar{D}_i$ of the central fibre $\bar{\mathfrak{X}}_0$ of $\bar{\mathfrak{X}} \to \mathcal{S}$ resolves all the singularities contained in every double curve $\bar{X}_{\rho} \subseteq \bar{D}_i$ and gives rise to a partial resolution $\mathfrak{X}' \to \mathcal{S}$. We may now blow up any other component $\bar{D}_k$ (for $1 \leq k \leq \bar{m}, \ k \neq i$) of the central fibre $\mathfrak{X}'_{0}$ of $\mathfrak{X}' \to \mathcal{S}$, resolving all the (remaining) singularities in the one-dimensional strata of $\bar{D}_k$. Moreover, the effect of this blowup on the dual intersection complex $\left( B', \mathscr{P}' \right)$ of $\mathfrak{X}' \to \mathcal{S}$ is similar to the description we give above. 

\begin{defn} \label{genlogsmooth}
We call the degenerations obtained from $\bar{\mathfrak{X}} \to \mathcal{S}$ by blowing up a sequence of irreducible components of the central fibre $\bar{\mathfrak{X}}_0$ (and its proper transforms after some blowups) \emph{generically log smooth} partial resolutions. Indeed, they are log smooth away from a subset of codimension $2$ of the singularities of $\bar{\mathfrak{X}} \to \mathcal{S}$ that have not been resolved.
\end{defn}

Note that any generically log smooth partial resolution $\mathfrak{X}'' \to \mathcal{S}$ of $\bar{\mathfrak{X}} \to \mathcal{S}$ has a well-defined dual intersection complex $\left(B'', \mathscr{P}'' \right)$ by a description similar to the above. Under the assumption that all the $\sigma \in \bar{\mathscr{P}}^{\max}$ are standard triangles, we have  $\left(B'', \mathscr{P}'' \right) \cong \left( \bar{B}, \bar{\mathscr{P}} \right)$ as polyhedral manifolds and $\left(B'', \mathscr{P}'' \right)$ has the structure of an affine manifold with singularities via Construction \ref{affstrint} (viewed via Observation \ref{altaffinecons}).

\subsubsection{A global resolution} \label{globalresolvesmall}

Blow up all the irreducible components of the central fibre $\bar{\mathfrak{X}}_0$ of $\bar{\mathfrak{X}} \to \mathcal{S}$ in any order. By the analysis of the previous section, this gives a resolution $\pi: \mathfrak{X} \to \bar{\mathfrak{X}}$ to a log smooth degeneration $\mathfrak{X} \to \mathcal{S}$. Here the log structure on $\mathfrak{X}$ is the divisorial log structure given by $D = D_1 + \dots + D_m$ where $D_i, \ 1 \leq i \leq m$ are the irreducible components of the central fibre $\mathfrak{X}_0$ (note that we have $m=\bar{m}$). Moreover, the exceptional locus of $\pi: \mathfrak{X} \to \bar{\mathfrak{X}}$ is the disjoint union of curves (using Notation \ref{blowupnotations}(2)) 
\begin{equation} \label{contcurves}
\left\{ E_{\rho,k}  \ | \ \rho \in \bar{\mathscr{P}}^{[1]}, \ 1 \leq k \leq r_{\rho}\right\},
\end{equation}
so $\pi$ is small.
The singularities of the affine structure on the dual intersection complex $(B, \mathscr{P})$ of $\mathfrak{X} \to \mathcal{S}$ are at the vertices (any singularity at an edge is pulled into a vertex at some point) and one may view the resolution tropically as a composition of transformations as in Figure \ref{blowingadivisorsmallaff}. Note that we have $(B, \mathscr{P}) \cong \left(\bar B, \bar{\mathscr{P}}\right)$ as polyhedral manifolds.

\begin{fact} \label{sncmlcy}
	$\mathfrak{X} \to \mathcal{S}$ is minimal log CY and $D$ is simple normal crossings. The dual intersection complex $(B, \mathscr{P})$ of $\mathfrak{X} \to \mathcal{S}$ satisfies Assumption \ref{manifoldass}.
\end{fact}
	
\begin{proof}
	    
To show that $D$ is simple normal crossings, it is enough to check the double and triple intersections of the components of $\mathfrak{X}_0$. $D$ is simple normal crossings at all the triple intersections by our assumption on the maximal cells $\sigma \in \bar{\mathscr{P}}^{\max}$ since $(B, \mathscr{P}) \cong \left(\bar B, \bar{\mathscr{P}}\right)$ implies that all the maximal cells $\sigma \in \mathscr{P}^{\max}$ are standard triangles as well. But $\sigma \in \mathscr{P}^{\max}$ being a standard triangle is equivalent to $D$ being simple normal crossings at the point $X_{\sigma}$. $D$ is simple normal crossings at all the points of the double curves $X_{\rho}, \ \rho \in \mathscr{P}^{[1]}$ since the local models at such points are of the form $\{xy=t \} \subseteq \Spec \kk[x,y,w_{\rho}]\llbracket t \rrbracket$. This comes from the strict transform of the divisorial local model at points that are not the intersection of $X_{\rho}$ with the curves $E_{\rho, k}, \ 1 \leq k \leq r_{\rho}$. At the intersections, the local model is given by an affine chart of the resolution of an ODP singularity. 
	
To show that $\mathfrak{X} \to \mathcal{S}$ is minimal log CY note that $K_{\bar{\mathfrak{X}}} \equiv -(\bar{D}_1+\dots + \bar{D}_{\bar{m}})$ since $\bar{\mathfrak{X}} \to \mathcal{S}$ is a toric degeneration. So $K_{\bar{\mathfrak{X}}} + \bar{D} \equiv 0$. But since $\pi: \mathfrak{X} \to \bar{\mathfrak{X}}$ is small, we have $$K_{\mathfrak{X}} \equiv \pi^{*}(K_{\bar{\mathfrak{X}}}) \equiv  -\left(\pi^{*}(\bar{D}_1) + \dots + \pi^{*}(\bar{D}_{\bar{m}})\right)  \equiv -(D_1+\dots +D_{\bar{m}}) \equiv -D,$$ so $\mathfrak{X} \to \mathcal{S}$ is minimal log CY (see Definition \ref{minlogcy}).

We have $(B, \mathscr{P}) \cong \left(\bar B, \bar{\mathscr{P}}\right)$ as polyhedral manifolds and $ \left(\bar B, \bar{\mathscr{P}}\right)$ is a polyhedral manifold of dimension $2$ in the sense of Definition \ref{polymani}. Therefore, $(B, \mathscr{P})$ is also a polyhedral manifold of dimension $2$ in the sense of Definition \ref{polymani}. Hence, $(B, \mathscr{P})$ satisfies Assumption \ref{manifoldass}.
\end{proof}

Let $P$ be a well-chosen monoid in the sense of Definition \ref{wellchosendef}. Explicitly, $P$ satisfies the following properties:
\begin{enumerate}
    \item $NE(\mathfrak{X}_0) \subseteq P \subseteq A_1(\mathfrak{X}_0, \mathbb{Z})$.
    \item $P$ is finitely generated and saturated.
    \item The group $P^{\times}$ of the invertible elements of $P$ coincides with the torsion part of $A_1(\mathfrak{X}_0, \mathbb{Z})$.
    \item There exists a (minimal) face $K$ of $P$ containing the classes of the contracted curves \eqref{contcurves}. 
\end{enumerate}
Note that Proposition \ref{wellchosen} implies (using Lemma \ref{mirrorexttrivial}) that such a monoid exists.

\begin{fact} \label{ampledivisorsmall}
The resolution $\pi: \mathfrak{X} \to \bar{\mathfrak{X}}$ satisfies the assumptions of Proposition \ref{P2}. Explicitly:
\begin{enumerate}
    \item There exists an effective $\pi$-ample divisor $D' = \sum_{i=1}^{m} a_i D_i$ such that $D_i$ is $\pi$-nef for any $i$ with $a_i = 0$.
    \item $K \cap NE(\mathfrak{X}_0) = K \cap NE(\mathfrak{X}_0)_{\numer}$ (under the splittings of \eqref{standardsplittings}), i.e. $\pi$ only contracts numerical classes of curves $C \subseteq \mathfrak{X}_0$.
\end{enumerate}
\end{fact}

\begin{proof}
It is clear that the contracted curves $E_{\rho, k}$ of \eqref{contcurves} define non-zero classes in $A_1(\mathfrak{X}_0, \mathbb{Z})_{\numer}$. Therefore, assumption (2) is satisfied. To show that $\pi: \mathfrak{X} \to \bar{\mathfrak{X}}$ satisfies assumption (1) it is enough to find a $\pi$-ample divisor $D' = \sum_{i=1}^{m} a_i D_i$ with $a_i > 0$ for $1 \leq i \leq m$.

Without loss of generality, assume that $\bar{D}_i$ is blown up before $\bar{D}_j$ for $1 \leq i < j \leq \bar{m}$. The cone $NE(\pi) \subseteq NE(\mathfrak{X}_0)$ of curves contracted by $\pi$ is finitely generated by the $E_{\rho, k}$ of \eqref{contcurves} so $NE(\pi) = \overline{NE(\pi)}$. By the (relative) Kleiman's criterion for ampleness, the divisor $D'$ is $\pi$-ample if and only if $D' \cdot E^{ij}_{ k} > 0$ for all $1 \leq i < j \leq \bar{m}$ with $\bar{D}_i \cap \bar{D}_j = \bar{X}_{\rho}, \ \rho \in \bar{\mathscr{P}}^{[1]}$ and all $1 \leq k \leq r_{\rho}$. The intersection of $D_j$ and $E^{ij}_{k}$ is transversal so $D_j \cdot E^{ij}_{k} = 1$. Further, we have $$(D_1 + \dots + D_m) \cdot E_k^{ij} = D \cdot E_k^{ij} = 0$$ since $D$ is numerically equivalent to the trivial divisor. We also have $D_l \cdot E_k^{ij} = 0$ for any $1 \leq l \leq m$ such that $l \neq i,j$ (since $D_l \cap E_k^{ij} = \varnothing$ in this case). This implies that $D_i \cdot E^{ij}_{k} = -1$.   So any divisor $D' = \sum_{i=1}^{m} a_i D_i$ with $a_i < a_j$ for $1 \leq i < j \leq m$ is $\pi$-ample. Further requiring $a_i > 0$ for $1 \leq i \leq m$, the resolution $\pi: \mathfrak{X} \to \bar{\mathfrak{X}}$ satisfies assumption (1).
\end{proof}

We define the initial slab functions $$\left\{f_{\underline{\rho}} \in \kk[\Lambda_{\rho}] \ \Big| \ \underline{\rho} \in \tilde{\bar{\mathscr{P}}}^{[n-1]} \right\}$$ for the toric degeneration $\bar{\mathfrak{X}} \to \mathcal{S}$ by setting
\begin{equation} \label{choiceslab}
f_{\underline{\rho}} := (1+w_{\rho})^{r_{\rho}}, \quad f_{\underline{\rho}'} := z^{m_{\underline{\rho}'\underline{\rho} }}f_{\underline{\rho}} = (1+w_{\rho}^{-1})^{r_{\rho}}
\end{equation}
for $\underline{\rho}, \underline{\rho}' \in \tilde{\bar{\mathscr{P}}}^{[1]}$ two slabs with $\underline{\rho}, \underline{\rho}' \subseteq \rho \in \bar{\mathscr{P}}^{[1]}$ and $w_{\rho} := z^{m_{\rho}}$ where $m_{\rho}$ is the integral generator of $\Lambda_{\rho}$ that points towards the vertex endpoint of $\underline{\rho}'$. If $\rho = \langle \bar{v}_i, \bar{v}_j \rangle$ for $\bar{v}_i, \bar{v_j} \in \bar{\mathscr{P}}^{[0]}$ and $\bar{D}_{i}$ is blown up before $\bar{D}_{j}$, then we assume that $\underline{\rho}, \underline{\rho}' \in \tilde{\bar{\mathscr{P}}}^{[1]}$ are chosen so that the vertex endpoint of $\underline{\rho}'$ is $\bar{v}_j$. 
Note that these are of the form required by \eqref{slabfunctions} in Proposition \ref{struct} that describes the possible initial slab functions. We have now defined all the necessary data in Basic Setup \ref{setupforproof}.

\subsubsection{The PL-isomorphism \texorpdfstring{$\Phi: B \to \bar{B}$}{Phi}}

The resolution $\pi: \mathfrak{X} \to \bar{\mathfrak{X}}$ gives rise to a piecewise-linear (PL) isomorphism $\Phi: (B, \mathscr{P}) \to \left(\bar{B}, \bar{\mathscr{P}}\right)$, linear on the maximal cells of $\mathscr{P}$. We shall use this isomorphism to relate the canonical scattering diagram $\mathfrak{D}$ on $(B, \mathscr{P})$ to the algorithmic scattering diagram $\bar{\mathfrak{D}}$ on $\left(\bar{B}, \bar{\mathscr{P}}\right)$.

\begin{cons} \label{plisosmall}
We define the PL-isomorphism $\Phi: (B, \mathscr{P}) \to \left(\bar{B}, \bar{\mathscr{P}}\right)$, linear on the maximal cells of $\mathscr{P}$, by using the descriptions of the affine structures with singularities on $(B, \mathscr{P})$ and $\left(\bar{B}, \bar{\mathscr{P}}\right)$ of Constructions \ref{affstrint} and \ref{affstrtoric} respectively. Recall that we have $(B, \mathscr{P}) \cong \left(\bar B, \bar{\mathscr{P}}\right)$ as polyhedral manifolds. For every $\tau \in \mathscr{P}$, let $\bar{\tau} \in \bar{\mathscr{P}}$ be the corresponding cell of $\bar{\mathscr{P}}$.

For all the maximal cells $\sigma \in \mathscr{P}^{\max}$, there are canonical identity maps $\Id: \sigma \to \bar{\sigma}$ compatible with each other. Therefore, it suffices to give PL-maps $\Phi_v: W_v \to W_{\bar{v}}$, linear on the maximal cells of $W_v$ (using the notation of \eqref{openstardef} for the open stars of $v$ and $\bar{v}$) and compatible with the identifications on the maximal cells. The fact that $(B, \mathscr{P}) \cong \left(\bar B, \bar{\mathscr{P}}\right)$ as polyhedral manifolds implies that there exist unique such maps $\Phi_v$. Indeed, using the refined description of Construction \ref{affchartstool}$(2')$, fix a $\rho \in \mathscr{P}^{[n-1]}$ such that $\rho = \sigma \cap \sigma'$ for $\sigma, \sigma' \in \mathscr{P}^{\max}$ and with $v \subseteq \rho$. Consider the embedding  $\psi_{\rho, v} : \sigma \cup \sigma' \to \mathbb{R}^2$ of Construction \ref{affstrint} and the corresponding embedding $\bar{\psi}_{\bar{\rho}, \bar{v}} : \bar{\sigma} \cup \bar{\sigma}' \to \mathbb{R}^2$ of Construction \ref{affstrtoric}. Then there is a unique piecewise-linear, linear on $\im(\Int(\sigma))$ and $\im(\Int(\sigma'))$, and sending the generator of $\im(\rho)$ to the generator of $\im(\bar{\rho})$, identification $\Phi_{\rho, v}: \im (\psi_{\rho, v}) \to \im (\bar{\psi}_{\bar \rho, \bar{v}})$ and the collection $\{\Phi_{\rho, v} \ | \ \rho \in \mathscr{P}^{[1]}, v \subseteq \rho \}$ of such identifications defines $\Phi_v$. We have defined a PL-isomorphism $\Phi: (B, \mathscr{P}) \to \left(\bar{B}, \bar{\mathscr{P}}\right)$, linear on the maximal cells of $\mathscr{P}$.

Since we assume that $D$ is simple normal crossings in this section, we can write an easy formula for $\Phi_{\rho, v}$. Indeed, we have an explicit description for $\psi_{\rho, v}$ given by Observation \ref{affstrdim2} and an explicit description for $\bar{\psi}_{\bar{\rho}, \bar{v}}$ given by Observation \ref{affstrdim2} and Remark \ref{toriccomponent}(3). Using the notations of Observation \ref{affstrdim2} and Notation \ref{intnotation} (and Corollary \ref{compinsmallerstrata}(3)), we just send 
\begin{multline*}
(0, 0) \to (0, 0), \ (1, 0) \to (1, 0), \ (0, 1) \to (0, 1), \\ \left(-  \left(X^2_{\langle v_0, v_1 \rangle}\right)_{D_{v_0}}, -1\right) \to \left(-  \left(\bar{X}^2_{\langle \bar{v}_0, \bar{v}_1 \rangle}\right)_{\bar{D}_{\bar{v}_0}}, -1\right),    
\end{multline*}
 and extend linearly. 

For any $\tau \in \mathscr{P}$, we write $\Phi(\tau) \in \bar{\mathscr{P}}$ for the image of $\tau$ under $\Phi$. Similarly, the image of any rational polyhedral subset $\mathfrak{p} \subseteq \sigma \in \mathscr{P}^{\max}$ of codimension $1$ is a well-defined rational polyhedral subset $\Phi(\mathfrak{p})\subseteq \Phi(\sigma) \in \bar{\mathscr{P}}^{\max}$ of codimension $1$. Note that $\Phi$ induces a canonical isomorphism $\mathcal{P}_x^{+} \cong \mathcal{\bar{P}}_{\Phi(x)}^{+}$ (see Section \ref{monomialsaffstr} for the definition of $\mathcal{P}_x^{+}$) for every $x \in \Int(\sigma), \ \sigma \in \mathscr{P}^{\max}$, so we also have an identification of monomials. We denote the image of a monomial $m \in \mathcal{P}_x^{+}$ by $\Phi(m) \in \mathcal{\bar{P}}_{\Phi(x)}^{+}$ and set $\Phi(z^m):=z^{\Phi(m)}$. We define the image under $\Phi$ of any codimension $0$ wall $(\mathfrak{p}, f_{\mathfrak{p}})$ on $(B, \mathscr{P})$ as a codimension $0$ wall $\Phi(\mathfrak{p}):= (\Phi(\mathfrak{p}), \Phi(f_{\mathfrak{p}}))$ on $\left(\bar B, \bar{\mathscr{P}}\right)$.

The image of any rational polyhedral subset $\mathfrak{b} \subseteq \rho \in \mathscr{P}^{[1]}$ (of full dimension $1$) is a well-defined rational polyhedral subset $\Phi(\mathfrak{b})\subseteq \Phi(\rho) \in \bar{\mathscr{P}}^{[1]}$ (of full dimension $1$). For $x \in \Int(\rho), \ \rho \in \mathscr{P}^{[1]}$ there is no canonical isomorphism between $\mathcal{P}_x^{+}$ and $\mathcal{\bar{P}}_{\Phi(x)}^{+}$ since $\Phi$ is only piecewise-linear in the neighbourhood of $\rho$. For a slab $\mathfrak{b} = (\mathfrak{b}, f_{\mathfrak{b}})$ on $(B, \mathscr{P})$ we will always have $f_{\mathfrak{b}} \in \kk[P][\Lambda_{\rho}] = \kk[P][z^{m_{\rho}}] \subseteq  \kk[\mathcal{P}^+_x]$ (via Notation \ref{notation}) where $m_{\rho}$ is an integral generator of $\Lambda_{\rho}$. We assume that $m_{\rho}$ is chosen to point in the same direction as $m_{\Phi(\rho)}$ under the identification of $\Lambda_{\rho}$ and $\Lambda_{\Phi(\rho)}$ (induced by the identification of $\rho$ and $\Phi(\rho)$ via the PL-isomorphism). We define $\Phi(m_{\rho}) := m_{\Phi(\rho)}$. As above, we set $\Phi(z^{m_{\rho}}):= z^{\Phi(m_{\rho})} = z^{m_{\Phi(\rho)}}$ and define $\Phi(\mathfrak{b}):=(\Phi(\mathfrak{b}), \Phi(f_{\mathfrak{b}}))$. Note that $\Phi(\mathfrak{b})$ is not a slab if there is a singularity $x_{\Phi(\rho)}$ contained in $\Phi(\mathfrak{b})$. For now, we just regard $\Phi(\mathfrak{b})$ as a codimension $1$ polyhedral subset of $\bar{B}$ with an attached function. We will define the slabs (one if $x_{\Phi(\rho)} \notin \Int(\Phi(\mathfrak{b}))$ and two if $x_{\Phi(\rho)} \in \Int(\Phi(\mathfrak{b}))$) on $\left(\bar B, \bar{\mathscr{P}}\right)$ corresponding to $\mathfrak{b}$ in Construction \ref{chs}.
\end{cons}

\subsection{An extension to non-small resolutions} \label{resolve}

We continue with the general setup of Section \ref{resolvesmall}. As in Section \ref{resolvesmall}, we are going to construct a resolution $\pi: \mathfrak{X} \to \bar{\mathfrak{X}}$ of a special (i.e. a divisorial log deformation with a smooth generic fibre) toric degeneration $\bar{\mathfrak{X}} \to \mathcal{S}$ of K3-s to a minimal log CY degeneration $\mathfrak{X} \to \mathcal{S}$ by blowing up a sequence of irreducible components of the central fibre $\bar{\mathfrak{X}}_0$ (and irreducible components of their proper transforms after some of the blowups). We relax the requirement on the maximal cells of $\bar{\mathscr{P}}$ and unless all the $\sigma \in \bar{\mathscr{P}}^{\max}$ are standard triangles or standard squares (i.e. as in Figure \ref{admisspolytopes}(2) with $l=k=1$), the resolution $\pi: \mathfrak{X} \to \bar{\mathfrak{X}}$ will not be small. 

We no longer assume that the maximal cells $\sigma \in \bar{\mathscr{P}}^{\max}$ of the dual intersection complex $\left( \bar{B}, \bar{\mathscr{P}} \right)$ are standard triangles. Instead, we require that every $\sigma \in \bar{\mathscr{P}}^{\max}$ is, up to the action of $AGL(2, \mathbb{Z})$\footnote{The group of affine unimodular transformations consisting of translations by an integer vector and the linear transformations in $GL(2, \mathbb{Z})$.}, one of the following (see Figure \ref{admisspolytopes}):
\begin{enumerate}
    \item A \enquote{lattice triangle of height $1$}, i.e. the convex hull of $(0,0), (0,1)$, and $(l,0)$ for some $l \in \mathbb{Z}_{\geq 1}$ (if $l=1$, we get a standard triangle).
    \item A \enquote{lattice trapezoid of height $1$}, i.e. the convex hull of $(0,0), (0,1), (l,0)$, and $(k,1)$ for some $l,k \in \mathbb{Z}_{\geq1}$.
    \item A \enquote{lattice equilateral triangle}, i.e. the convex hull of $(0,0), (0,k)$, and $(k,0)$ for some $k \in \mathbb{Z}_{\geq 1}$ (if $k=1$, we get a standard triangle).
    \item A \enquote{lattice rectangle}, i.e. the convex hull of $(0,0), (0,k), (l,0)$, and $(l,k)$ for some $l,k \in \mathbb{Z}_{\geq1}$ (if $l=1$ or $k=1$, this is also an example of (2)).
    \item The \enquote{minimal lattice hexagon}, i.e. the convex hull of $(0,0), (1,0), (2, 1),$ $(2,2)$, $(1,2)$, and $(0,1)$.
\end{enumerate}

\begin{figure}[H]
\begin{center}
\begin{minipage}{.3\textwidth}
\begin{tikzpicture}[scale=0.7]

\draw[fill] (5/4*0,5/4*0) circle (1.5pt) coordinate (m-0-0);
\draw[fill] (5/4*0,5/4*1) circle (1.5pt) coordinate (m-0-1);
\draw[fill] (5/4*1,5/4*0) circle (1.5pt) coordinate (m-1-0);
\draw[fill] (5/4*3,5/4*0) circle (1.5pt) coordinate (m-3-0);
\draw[fill] (5/4*2,5/4*0) circle (1.5pt) coordinate (m-2-0);

\draw (m-0-0) -- (m-3-0) -- (m-0-1) -- cycle;

\end{tikzpicture}
\captionsetup{labelformat=empty}
\caption{(1)}
\end{minipage}
\hspace{1cm}
\begin{minipage}{.3\textwidth}
\begin{tikzpicture}[scale=0.7]

\draw[fill] (5/4*0,5/4*0) circle (1.5pt) coordinate (m-0-0);
\draw[fill] (5/4*0,5/4*1) circle (1.5pt) coordinate (m-0-1);
\draw[fill] (5/4*1,5/4*0) circle (1.5pt) coordinate (m-1-0);
\draw[fill] (5/4*3,5/4*0) circle (1.5pt) coordinate (m-3-0);
\draw[fill] (5/4*2,5/4*0) circle (1.5pt) coordinate (m-2-0);
\draw[fill] (5/4*4,5/4*0) circle (1.5pt) coordinate (m-4-0);
\draw[fill] (5/4*1,5/4*1) circle (1.5pt) coordinate (m-1-1);
\draw[fill] (5/4*2,5/4*1) circle (1.5pt) coordinate (m-2-1);

\draw (m-0-0) -- (m-4-0) -- (m-2-1) -- (m-0-1) -- cycle;

\end{tikzpicture}
\captionsetup{labelformat=empty}
\caption{(2)}
\end{minipage}

\begin{minipage}{.3\textwidth}

\begin{tikzpicture}[scale=0.5]

\draw[fill] (5/4*0,5/4*0) circle (1.5pt) coordinate (m-0-0);
\draw[fill] (5/4*0,5/4*1) circle (1.5pt) coordinate (m-0-1);
\draw[fill] (5/4*1,5/4*0) circle (1.5pt) coordinate (m-1-0);
\draw[fill] (5/4*2,5/4*0) circle (1.5pt) coordinate (m-2-0);
\draw[fill] (5/4*1,5/4*1) circle (1.5pt) coordinate (m-1-1);
\draw[fill] (5/4*0,5/4*2) circle (1.5pt) coordinate (m-0-2);
\draw[fill] (5/4*3,5/4*0) circle (1.5pt) coordinate (m-3-0);
\draw[fill] (5/4*2,5/4*1) circle (1.5pt) coordinate (m-2-1);
\draw[fill] (5/4*1,5/4*2) circle (1.5pt) coordinate (m-1-2);
\draw[fill] (5/4*0,5/4*3) circle (1.5pt) coordinate (m-0-3);
\draw[fill] (5/4*4,5/4*0) circle (1.5pt) coordinate (m-4-0);
\draw[fill] (5/4*3,5/4*1) circle (1.5pt) coordinate (m-3-1);
\draw[fill] (5/4*2,5/4*2) circle (1.5pt) coordinate (m-2-2);
\draw[fill] (5/4*1,5/4*3) circle (1.5pt) coordinate (m-1-3);
\draw[fill] (5/4*0,5/4*4) circle (1.5pt) coordinate (m-0-4);

\draw (m-0-0) -- (m-0-4) -- (m-4-0) -- cycle;
\end{tikzpicture}
\captionsetup{labelformat=empty}
\caption{(3)}
\end{minipage}
\hspace{0.5cm}
\begin{minipage}{.3\textwidth}
\begin{tikzpicture}[scale=0.5]

\draw[fill] (5/4*0,5/4*0) circle (1.5pt) coordinate (m-0-0);
\draw[fill] (5/4*0,5/4*1) circle (1.5pt) coordinate (m-0-1);
\draw[fill] (5/4*1,5/4*0) circle (1.5pt) coordinate (m-1-0);
\draw[fill] (5/4*2,5/4*0) circle (1.5pt) coordinate (m-2-0);
\draw[fill] (5/4*1,5/4*1) circle (1.5pt) coordinate (m-1-1);
\draw[fill] (5/4*0,5/4*2) circle (1.5pt) coordinate (m-0-2);
\draw[fill] (5/4*3,5/4*0) circle (1.5pt) coordinate (m-3-0);
\draw[fill] (5/4*2,5/4*1) circle (1.5pt) coordinate (m-2-1);
\draw[fill] (5/4*1,5/4*2) circle (1.5pt) coordinate (m-1-2);
\draw[fill] (5/4*0,5/4*3) circle (1.5pt) coordinate (m-0-3);
\draw[fill] (5/4*4,5/4*0) circle (1.5pt) coordinate (m-4-0);
\draw[fill] (5/4*3,5/4*1) circle (1.5pt) coordinate (m-3-1);
\draw[fill] (5/4*2,5/4*2) circle (1.5pt) coordinate (m-2-2);
\draw[fill] (5/4*1,5/4*3) circle (1.5pt) coordinate (m-1-3);
\draw[fill] (5/4*2,5/4*3) circle (1.5pt) coordinate (m-2-3);
\draw[fill] (5/4*3,5/4*3) circle (1.5pt) coordinate (m-3-3);
\draw[fill] (5/4*4,5/4*3) circle (1.5pt) coordinate (m-4-3);
\draw[fill] (5/4*3,5/4*2) circle (1.5pt) coordinate (m-3-2);
\draw[fill] (5/4*4,5/4*2) circle (1.5pt) coordinate (m-4-2);
\draw[fill] (5/4*4,5/4*1) circle (1.5pt) coordinate (m-4-1);

\draw (m-0-0) -- (m-4-0) -- (m-4-3) -- (m-0-3) -- cycle;

\end{tikzpicture}
\captionsetup{labelformat=empty}
\caption{(4)}
\end{minipage}
\begin{minipage}{.2\textwidth}
\begin{tikzpicture}[scale=0.7]

\draw[fill] (5/4*0,5/4*0) circle (1.5pt) coordinate (m-0-0);
\draw[fill] (5/4*0,5/4*1) circle (1.5pt) coordinate (m-0-1);
\draw[fill] (5/4*1,5/4*0) circle (1.5pt) coordinate (m-1-0);
\draw[fill] (5/4*1,5/4*1) circle (1.5pt) coordinate (m-1-1);
\draw[fill] (5/4*2,5/4*1) circle (1.5pt) coordinate (m-2-1);
\draw[fill] (5/4*1,5/4*2) circle (1.5pt) coordinate (m-1-2);
\draw[fill] (5/4*2,5/4*2) circle (1.5pt) coordinate (m-2-2);

\draw (m-0-0) -- (m-1-0) -- (m-2-1) -- (m-2-2) -- (m-1-2) -- (m-0-1) -- cycle;

\end{tikzpicture}
\captionsetup{labelformat=empty}
\caption{(5)}
\end{minipage}
\end{center}
\setcounter{figure}{4}
\caption{Types of $\sigma \in \bar{\mathscr{P}}^{\max}$ in Section \ref{resolve}.} 
\label{admisspolytopes}
\end{figure}

The motivation for restricting the types of $\sigma \in \bar{\mathscr{P}}^{\max}$ to Figure \ref{admisspolytopes} is as follows. Blowing up a component of $\bar{\mathfrak{X}} \to \mathcal{S}$ corresponding to a vertex $\bar{v} \in \bar{\mathscr{P}}^{[0]}$ induces subdivisions of the maximal cells $\sigma \in \bar{\mathscr{P}}^{\max}$ with $\bar{v} \subseteq \sigma$. To ensure that the dual intersection complex of the partial resolution is a polyhedral manifold, we need to require that the subdivision of every such $\sigma \in \bar{\mathscr{P}}^{\max}$ is integral and only subdivides edges of $\sigma$ adjacent to $\bar{v}$. The types of $\sigma \in \bar{\mathscr{P}}^{\max}$ in Figure \ref{admisspolytopes} form a class satisfying these properties (see Figure \ref{admisspolytopesresolve} below). We expect that it is the largest such class. 

The subdivision of every $\sigma \in \bar{\mathscr{P}}^{\max}$ with $\bar{v} \subseteq \sigma$ can be understood via the local toric model for the codimension $2$ stratum $\bar{X}_{\sigma}$. We shall now discuss these models in detail. We will treat the case of general $\sigma \in \bar{\mathscr{P}}^{\max}$ in Section \ref{admissible}.

\subsubsection{Local models in codimension \texorpdfstring{$2$}{2}} \label{toricmodelcodim2}

The local toric model for a codimension $2$ stratum $\bar{X}_{\sigma}, \ \sigma \in \bar{\mathscr{P}}^{\max}$ is given by the cone over $\sigma$. To understand the effect of blowing up a component $\bar{D}_i$ of $\bar{\mathfrak{X}}_0$ on the dual intersection complex $\left( \bar{B}, \bar{\mathscr{P}} \right)$, we need to first study the effect of blowing up a toric divisor of the toric variety  $X_{\mathbf{C}\sigma}$ defined by the cone $\mathbf{C} \sigma$.

\begin{lemma} \label{blowuplemma}
Let $X_{\mathbf{C}\sigma}$ be the toric variety defined by the cone $\mathbf{C}\sigma$. Then blowing up the toric divisor corresponding to a ray $\mathbf{C}v$ of $\mathbf{C}\sigma$ (here $v \subseteq \sigma$ is a vertex) corresponds to the subdivision of $\mathbf{C}\sigma$ by the bend locus of the PL-function $$\varphi_{D_{\mathbf{C}v}} :=\min \left\{m \ | \ z^m \in I_{D_{\mathbf{C}v}} 
 \right\}.$$
Here $I_{D_{\mathbf{C}v}} \subseteq \kk[\widecheck{\mathbf{C}\sigma} \cap N]$ is the ideal of the divisor $D_{\mathbf{C}v}$. Note that it is enough to take the minimum over the (finitely many) generators of $I_{D_{\mathbf{C}v}}$.
\end{lemma}

\begin{proof}
The claim follows from the discussion of \cite[Section 3]{T}. Indeed, consider the subdivision $\Sigma$ of $\mathbf{C}\sigma$ given by the bend locus of $\varphi_{D_{\mathbf{C}v}}$. Then, by construction, $\varphi_{D_{\mathbf{C}v}}$ satisfies the requirements (1) and (2) of \cite[Section 3]{T} on the support function for $\Sigma$ and the generators $z^{m_1}, \dots, z^{m_k}$ of $I_{D_{\mathbf{C}v}}$ define the Cartier data $m_1, \dots, m_k$ for $\Sigma$ and $\varphi_{D_{\mathbf{C}v}}$. Then, by the recipe of \cite[Section 3]{T}, the ideal of the blowup is defined by $\langle z^{m_1}, \dots, z^{m_k} \rangle \subseteq \kk[\widecheck{\mathbf{C}\sigma} \cap N]$, which is precisely $I_{D_{\mathbf{C}v}}$.
\end{proof}

\begin{cor} \label{subwellbehaved}
For $\sigma \in \bar{\mathscr{P}}^{\max}$ of one of the types in Figure \ref{admisspolytopes}, consider the subdivision of $\sigma$ (induced from a subdivision of $\mathbf{C}\sigma$) that corresponds to blowing up the toric divisor of $X_{\mathbf{C\sigma}}$ corresponding to a vertex $v$. This subdivision is integral, only subdivides edges of $\sigma$ adjacent to $v$, and all the maximal cells of the subdivision are of one of the types in Figure \ref{admisspolytopes}.
\end{cor}

\begin{proof}
We have discussed the cases of $\sigma$ a standard triangle or a square (i.e. the convex hull of $(0, 0), (1, 0), (0, 1)$, and $(1,1)$) in Section \ref{resolvesmall}. The other cases follow by a direct computation using Lemma \ref{blowuplemma} and the fact that the ideal $I_{D_{\mathbf{C}v}}$ is generated by, for $m_1, \dots, m_k$ a choice of generators of $\widecheck{\mathbf{C}\sigma} \cap N$,  $$\left\{z^{m_i} \ | \ m_i \notin \mathbf{C}v^{\perp}\cap N, \ 1 \leq i \leq k  \right\}.$$ 
\end{proof}

In Figure \ref{admisspolytopesresolve}, we give the subdivisions of $\sigma \in \bar{\mathscr{P}}^{\max}$ for the types of $\sigma \in \bar{\mathscr{P}}^{\max}$ in Figure \ref{admisspolytopes} (the subdivision is given in blue and the vertex corresponding to the blown up divisor is highlighted in red).\footnote{Figure \ref{admisspolytopesresolve} does not cover the cases when $\sigma \in \bar{\mathscr{P}}^{\max}$ is either a standard triangle or square, or of the type of Figure \ref{admisspolytopes}(2) with $l=1$. If $\sigma \in \bar{\mathscr{P}}^{\max}$ is a standard triangle, the subdivision is trivial. If it is a standard square, the subdivision is as in Figure \ref{coneoversquare}. Finally, if $\sigma \in \bar{\mathscr{P}}^{\max}$ is of the type of Figure \ref{admisspolytopes}(2) with $l=1$, it is subdivided by the diagonal connecting $(0,0)$ and $(k,1)$.}

\begin{figure}[t]
\begin{center}
\begin{minipage}{.3\textwidth}
\begin{tikzpicture}[scale=0.7]

\draw[fill, color = red] (5/4*0,5/4*0) circle (3pt) coordinate (m-0-0);
\draw[fill] (5/4*0,5/4*1) circle (1.5pt) coordinate (m-0-1);
\draw[fill] (5/4*1,5/4*0) circle (1.5pt) coordinate (m-1-0);
\draw[fill] (5/4*3,5/4*0) circle (1.5pt) coordinate (m-3-0);
\draw[fill] (5/4*2,5/4*0) circle (1.5pt) coordinate (m-2-0);

\draw (m-0-0) -- (m-3-0) -- (m-0-1) -- cycle;

\draw [color=blue, line width=1pt] (m-0-1)--(m-2-0);

\end{tikzpicture}
\captionsetup{labelformat=empty}
\caption{(1)}
\end{minipage}
\hspace{1cm}
\begin{minipage}{.3\textwidth}
\begin{tikzpicture}[scale=0.7]

\draw[fill, color = red] (5/4*0,5/4*0) circle (3pt) coordinate (m-0-0);
\draw[fill] (5/4*0,5/4*1) circle (1.5pt) coordinate (m-0-1);
\draw[fill] (5/4*1,5/4*0) circle (1.5pt) coordinate (m-1-0);
\draw[fill] (5/4*3,5/4*0) circle (1.5pt) coordinate (m-3-0);
\draw[fill] (5/4*2,5/4*0) circle (1.5pt) coordinate (m-2-0);
\draw[fill] (5/4*4,5/4*0) circle (1.5pt) coordinate (m-4-0);
\draw[fill] (5/4*1,5/4*1) circle (1.5pt) coordinate (m-1-1);
\draw[fill] (5/4*2,5/4*1) circle (1.5pt) coordinate (m-2-1);

\draw (m-0-0) -- (m-4-0) -- (m-2-1) -- (m-0-1) -- cycle;

\draw [color=blue, line width=1pt] (m-2-1)--(m-3-0);

\end{tikzpicture}
\captionsetup{labelformat=empty}
\caption{(2)}
\end{minipage}

\begin{minipage}{.3\textwidth}

\begin{tikzpicture}[scale=0.5]

\draw[fill, color = red] (5/4*0,5/4*0) circle (4pt) coordinate (m-0-0);
\draw[fill] (5/4*0,5/4*1) circle (1.5pt) coordinate (m-0-1);
\draw[fill] (5/4*1,5/4*0) circle (1.5pt) coordinate (m-1-0);
\draw[fill] (5/4*2,5/4*0) circle (1.5pt) coordinate (m-2-0);
\draw[fill] (5/4*1,5/4*1) circle (1.5pt) coordinate (m-1-1);
\draw[fill] (5/4*0,5/4*2) circle (1.5pt) coordinate (m-0-2);
\draw[fill] (5/4*3,5/4*0) circle (1.5pt) coordinate (m-3-0);
\draw[fill] (5/4*2,5/4*1) circle (1.5pt) coordinate (m-2-1);
\draw[fill] (5/4*1,5/4*2) circle (1.5pt) coordinate (m-1-2);
\draw[fill] (5/4*0,5/4*3) circle (1.5pt) coordinate (m-0-3);
\draw[fill] (5/4*4,5/4*0) circle (1.5pt) coordinate (m-4-0);
\draw[fill] (5/4*3,5/4*1) circle (1.5pt) coordinate (m-3-1);
\draw[fill] (5/4*2,5/4*2) circle (1.5pt) coordinate (m-2-2);
\draw[fill] (5/4*1,5/4*3) circle (1.5pt) coordinate (m-1-3);
\draw[fill] (5/4*0,5/4*4) circle (1.5pt) coordinate (m-0-4);

\draw [color=blue, line width=1pt] (m-0-3)--(m-3-0);

\draw (m-0-0) -- (m-0-4) -- (m-4-0) -- cycle;
\end{tikzpicture}
\captionsetup{labelformat=empty}
\caption{(3)}
\end{minipage}
\hspace{0.5cm}
\begin{minipage}{.3\textwidth}
\begin{tikzpicture}[scale=0.5]

\draw[fill, color = red] (5/4*0,5/4*0) circle (4pt) coordinate (m-0-0);
\draw[fill] (5/4*0,5/4*1) circle (1.5pt) coordinate (m-0-1);
\draw[fill] (5/4*1,5/4*0) circle (1.5pt) coordinate (m-1-0);
\draw[fill] (5/4*2,5/4*0) circle (1.5pt) coordinate (m-2-0);
\draw[fill] (5/4*1,5/4*1) circle (1.5pt) coordinate (m-1-1);
\draw[fill] (5/4*0,5/4*2) circle (1.5pt) coordinate (m-0-2);
\draw[fill] (5/4*3,5/4*0) circle (1.5pt) coordinate (m-3-0);
\draw[fill] (5/4*2,5/4*1) circle (1.5pt) coordinate (m-2-1);
\draw[fill] (5/4*1,5/4*2) circle (1.5pt) coordinate (m-1-2);
\draw[fill] (5/4*0,5/4*3) circle (1.5pt) coordinate (m-0-3);
\draw[fill] (5/4*4,5/4*0) circle (1.5pt) coordinate (m-4-0);
\draw[fill] (5/4*3,5/4*1) circle (1.5pt) coordinate (m-3-1);
\draw[fill] (5/4*2,5/4*2) circle (1.5pt) coordinate (m-2-2);
\draw[fill] (5/4*1,5/4*3) circle (1.5pt) coordinate (m-1-3);
\draw[fill] (5/4*2,5/4*3) circle (1.5pt) coordinate (m-2-3);
\draw[fill] (5/4*3,5/4*3) circle (1.5pt) coordinate (m-3-3);
\draw[fill] (5/4*4,5/4*3) circle (1.5pt) coordinate (m-4-3);
\draw[fill] (5/4*3,5/4*2) circle (1.5pt) coordinate (m-3-2);
\draw[fill] (5/4*4,5/4*2) circle (1.5pt) coordinate (m-4-2);
\draw[fill] (5/4*4,5/4*1) circle (1.5pt) coordinate (m-4-1);

\draw [color=blue, line width=1pt] (m-0-2)--(m-3-2);
\draw [color=blue, line width=1pt] (m-3-0)--(m-3-2);
\draw [color=blue, line width=1pt] (m-4-3)--(m-3-2);

\draw (m-0-0) -- (m-4-0) -- (m-4-3) -- (m-0-3) -- cycle;

\end{tikzpicture}
\captionsetup{labelformat=empty}
\caption{(4)}
\end{minipage}
\begin{minipage}{.2\textwidth}
\begin{tikzpicture}[scale=0.7]

\draw[fill, color = red] (5/4*0,5/4*0) circle (3pt) coordinate (m-0-0);
\draw[fill] (5/4*0,5/4*1) circle (1.5pt) coordinate (m-0-1);
\draw[fill] (5/4*1,5/4*0) circle (1.5pt) coordinate (m-1-0);
\draw[fill] (5/4*1,5/4*1) circle (1.5pt) coordinate (m-1-1);
\draw[fill] (5/4*2,5/4*1) circle (1.5pt) coordinate (m-2-1);
\draw[fill] (5/4*1,5/4*2) circle (1.5pt) coordinate (m-1-2);
\draw[fill] (5/4*2,5/4*2) circle (1.5pt) coordinate (m-2-2);

\draw (m-0-0) -- (m-1-0) -- (m-2-1) -- (m-2-2) -- (m-1-2) -- (m-0-1) -- cycle;

\draw [color=blue, line width=1pt] (m-0-0)--(m-1-1);
\draw [color=blue, line width=1pt] (m-2-1)--(m-1-1);
\draw [color=blue, line width=1pt] (m-1-2)--(m-1-1);
\draw [color=blue, line width=1pt] (m-2-2)--(m-1-1);

\end{tikzpicture}
\captionsetup{labelformat=empty}
\caption{(5)}
\end{minipage}
\end{center}
\setcounter{figure}{5}
\caption{Subdivisions corresponding to divisor blowups.} 
\label{admisspolytopesresolve}
\end{figure}

\begin{remark} \label{moregeneralsubdivisions}
The reason that we restrict to $\sigma \in \bar{\mathscr{P}}^{\max}$ of Figure \ref{admisspolytopes} in this section is that the subdivision of a more general $\sigma \in \bar{\mathscr{P}}^{\max}$ induced by blowing up the divisor corresponding to a vertex $v \subseteq \sigma$ could be non-integral or subdivide edges that are not adjacent to $v$. Globally, a blowup of an irreducible component $\bar{D}_i$ of $\bar{\mathfrak{X}}_0$ with $\bar{v}_i \subseteq \sigma \in \bar{\mathscr{P}}^{\max}$ would produce a subdivision $\left(B', \mathscr{P}'\right)$ of $\left(\bar{B}, \bar{\mathscr{P}}\right)$ that no longer satisfies Assumption \ref{manifoldass}. See Figure \ref{pathexamples} for examples of these phenomena (we use the same notations as in Figure \ref{admisspolytopesresolve}).
\end{remark}

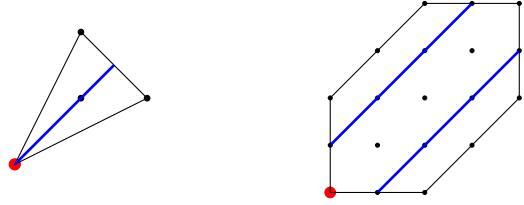
\begin{figure}[H]
\begin{center}
\begin{minipage}{.2\textwidth}
\begin{tikzpicture}[scale=0.7]

\draw[fill, color = red] (5/4*0,5/4*0) circle (3pt) coordinate (m-0-0);
\draw[fill] (5/4*2,5/4*1) circle (1.5pt) coordinate (m-2-1);
\draw[fill] (5/4*1,5/4*2) circle (1.5pt) coordinate (m-1-2);
\draw[fill] (5/4*1,5/4*1) circle (1.5pt) coordinate (m-1-1);

\draw (m-0-0) -- (m-1-2) -- (m-2-1) -- cycle;

\draw [color=blue, line width=1pt] (m-0-0)--(5/4*1.5,5/4*1.5);

\end{tikzpicture}
\end{minipage}
\hspace{1cm}
\begin{minipage}{.3\textwidth}
\begin{tikzpicture}[scale=0.5]

\draw[fill, color = red] (5/4*0,5/4*0) circle (4pt) coordinate (m-0-0);
\draw[fill] (5/4*0,5/4*1) circle (1.5pt) coordinate (m-0-1);
\draw[fill] (5/4*1,5/4*0) circle (1.5pt) coordinate (m-1-0);
\draw[fill] (5/4*2,5/4*0) circle (1.5pt) coordinate (m-2-0);
\draw[fill] (5/4*1,5/4*1) circle (1.5pt) coordinate (m-1-1);
\draw[fill] (5/4*0,5/4*2) circle (1.5pt) coordinate (m-0-2);
\draw[fill] (5/4*2,5/4*1) circle (1.5pt) coordinate (m-2-1);
\draw[fill] (5/4*1,5/4*2) circle (1.5pt) coordinate (m-1-2);
\draw[fill] (5/4*3,5/4*1) circle (1.5pt) coordinate (m-3-1);
\draw[fill] (5/4*2,5/4*2) circle (1.5pt) coordinate (m-2-2);
\draw[fill] (5/4*1,5/4*3) circle (1.5pt) coordinate (m-1-3);
\draw[fill] (5/4*2,5/4*3) circle (1.5pt) coordinate (m-2-3);
\draw[fill] (5/4*3,5/4*3) circle (1.5pt) coordinate (m-3-3);
\draw[fill] (5/4*4,5/4*3) circle (1.5pt) coordinate (m-4-3);
\draw[fill] (5/4*3,5/4*2) circle (1.5pt) coordinate (m-3-2);
\draw[fill] (5/4*4,5/4*2) circle (1.5pt) coordinate (m-4-2);
\draw[fill] (5/4*3,5/4*4) circle (1.5pt) coordinate (m-3-4);
\draw[fill] (5/4*2,5/4*4) circle (1.5pt) coordinate (m-2-4);
\draw[fill] (5/4*4,5/4*4) circle (1.5pt) coordinate (m-4-4);

\draw [color=blue, line width=1pt] (m-0-1)--(m-3-4);
\draw [color=blue, line width=1pt] (m-1-0)--(m-4-3);


\draw (m-0-0) -- (m-2-0) -- (m-4-2) -- (m-4-4) -- (m-2-4) -- (m-0-2) -- cycle;

\end{tikzpicture}
\end{minipage}
\end{center}
\setcounter{figure}{6}
\caption{Subdivisions of more general $\sigma \in \bar{\mathscr{P}}^{\max}$.} 
\label{pathexamples}
\end{figure}

We expect that the types of $\sigma \in \bar{\mathscr{P}}^{\max}$ in Figure \ref{admisspolytopes} define the most general class for which Corollary \ref{subwellbehaved} holds.

\subsubsection{Resolving local models in codimension \texorpdfstring{$1$}{1}} \label{localmodelgen}

As in Section \ref{locmod}, by Observation \ref{modelslocal}(2), a toric degeneration $\bar{\mathfrak{X}} \to \mathcal{S}$ that is a divisorial log deformation with a smooth generic fibre has local models of the form $\left\{xy = t^l \right\} \subseteq \Spec \kk[x,y,w_{\rho}]\llbracket t \rrbracket$ and $\left\{xy = t^lw_{\rho} \right\} \subseteq \Spec \kk[x,y,w_{\rho}]\llbracket t \rrbracket$ (with the natural map to $\Spec \kk\llbracket t \rrbracket$) at the non-singular and singular points of the double curve $\bar{X}_{\rho}, \rho \in \bar{\mathscr{P}}^{[1]}$ respectively, where $l$ is the integral length of $\rho$.

The local model $\left\{xy = t^l \right\} \subseteq \Spec \kk[x,y,w_{\rho}]\llbracket t \rrbracket$ is not simple normal crossings unless $l =1$ (it is a line of $A_{l-1}$ singularities), but the map to $\Spec \kk\llbracket t \rrbracket$ is log smooth (here the induced divisorial log structure in the model is given by $\left\{t=0\right\}$). Suppose that $l>1$. Blowing up $\left\{x = t =0 \right\}$ (without loss of generality) produces a chart that is simple normal crossings and a chart where the local model takes form $\left\{xy = t^{l-1} \right\} \subseteq \Spec \kk[x,y,w_{\rho}]\llbracket t \rrbracket$ (here we still denote the strict transform of $\left\{x = t =0 \right\}$ by $\left\{x = t =0 \right\}$ and denote the \enquote{new} divisor by $\left\{x'=t=0 \right\}$). We can then blow up $\left\{x=t=0 \right\}$ again or we can blow up $\left\{x' = t =0 \right\}$. Repeating this procedure $l-1$ times (with choices of divisors to blow up at every step) produces a simple normal crossings degeneration. 

\begin{sloppypar}
\begin{notation} \label{divisorlocalnotation}
We will follow the convention that for every successive blowup, the strict transform of $\left\{x=t=0 \right\}$ (resp. $\left\{x'=t=0 \right\}$) will be denoted by $\left\{x=t=0 \right\}$ (resp. $\left\{x'=t=0 \right\}$) and the \enquote{new} divisor contained in the proper transform will be denoted by $\left\{x'=t=0 \right\}$ (resp. $\left\{x=t=0 \right\}$).    
\end{notation}   
\end{sloppypar}

Torically, $\left\{xy = t^l \right\} \subseteq \Spec \kk[x,y,w_{\rho}]\llbracket t \rrbracket$ is (the completion in $t$ of) the affine variety defined by the cone over the triangle that is the convex hull of $(0, 0, 1)$, $(0, 1, 1)$, and $ (l, 0, 1)$ (i.e. it is of the same type as Figure \ref{admisspolytopesresolve}(1)). We assume that $(0, 0, 1)$ corresponds to $\{x = t = 0\}$, $(l, 0, 1)$ corresponds to $\{y = t = 0\}$, and that we first blow up $\{x = t = 0\}$. Suppose that we have blown up $\{x = t = 0\}$ $i_x$ times and $\{x' = t = 0\}$ $i_{x'}$ times (with $i_x + i_{x'} \leq l-1$), using Notation \ref{divisorlocalnotation}. Then this partial resolution corresponds to the subdivision of the triangle by the line segments connecting $(0, 1, 1)$ to $(1, 0, 1), (2, 0, 1), \dots, (i_{x'}, 0, 1)$ and to $(l-1, 0, 1), \dots, (l-i_x, 0, 1)$. Note that after $l-1$ blowups, we arrive at the 
same resolution regardless of the choices of divisors that we blow up at each step. See Figure \ref{model1} for these subdivisions in the case that $l=4$.

\begin{figure}[t]
\begin{center}

\begin{minipage}{.3\textwidth}
\begin{tikzpicture}[line cap=round,line join=round,>=triangle 45,x=1cm,y=1cm, scale = 3]
\clip(-0.2,-0.2) rectangle (1.25,1.25);
\draw  (0,0)-- (1,0);
\draw  (0,0)-- (0,1);
\draw  (0,1)-- (1,0);
\draw [fill] (0,0) circle (0.5pt);
\draw [fill] (0.25,0) circle (0.5pt);
\draw [fill] (0.5,0) circle (0.5pt);
\draw [fill] (0.75,0) circle (0.5pt);
\draw [fill] (0,1) circle (0.5pt);
\draw [fill] (1,0) circle (0.5pt);
\draw(0, -0.1) node{\tiny $x=t=0$};
\draw(1, -0.1) node{\tiny $y=t=0$};
\draw(0, 1.07) node{\tiny $w_{\rho}$};
\end{tikzpicture}
\captionsetup{labelformat=empty}
\caption{$xy=t^4$ \textcolor{white}{llllll}}
\end{minipage}
\hspace{0.3cm}
\begin{minipage}{.3\textwidth}
\begin{tikzpicture}[line cap=round,line join=round,>=triangle 45,x=1cm,y=1cm, scale = 3]
\clip(-0.2,-0.2) rectangle (1.25,1.25);
\draw  (0,0)-- (1,0);
\draw  (0,0)-- (0,1);
\draw  (0,1)-- (1,0);
\draw  (0,1)-- (0.75,0);
\draw [fill] (0,0) circle (0.5pt);
\draw [fill] (0.25,0) circle (0.5pt);
\draw [fill] (0.5,0) circle (0.5pt);
\draw [fill] (0.75,0) circle (0.5pt);
\draw [fill] (0,1) circle (0.5pt);
\draw [fill] (1,0) circle (0.5pt);
\draw(0, -0.1) node{\tiny $x=t=0$};
\draw(1, -0.1) node{\tiny $y=t=0$};
\draw(0, 1.07) node{\tiny $w_{\rho}$};

\end{tikzpicture}
\captionsetup{labelformat=empty}
\caption{Blowup of $x=t=0$}
\end{minipage}
\hspace{0.3cm}
\begin{minipage}{.3\textwidth}
\begin{tikzpicture}[line cap=round,line join=round,>=triangle 45,x=1cm,y=1cm, scale = 3]
\clip(-0.2,-0.2) rectangle (1.25,1.25);
\draw  (0,0)-- (1,0);
\draw  (0,0)-- (0,1);
\draw  (0,1)-- (1,0);
\draw  (0,1)-- (0.25,0);
\draw  (0,1)-- (0.5,0);
\draw  (0,1)-- (0.75,0);
\draw [fill] (0,0) circle (0.5pt);
\draw [fill] (0.25,0) circle (0.5pt);
\draw [fill] (0.5,0) circle (0.5pt);
\draw [fill] (0.75,0) circle (0.5pt);
\draw [fill] (0,1) circle (0.5pt);
\draw [fill] (1,0) circle (0.5pt);
\draw(0, -0.1) node{\tiny $x=t=0$};
\draw(1, -0.1) node{\tiny $y=t=0$};
\draw(0, 1.07) node{\tiny $w_{\rho}$};

\end{tikzpicture}
\captionsetup{labelformat=empty}
\caption{\ \ \ After \\ 3 blowups}
\end{minipage}
\end{center}
\setcounter{figure}{7}
\caption{Resolution of $\left\{xy = t^4 \right\} \subseteq \Spec \kk[x,y,w_{\rho}]\llbracket t \rrbracket$.} \label{model1}
\end{figure}

The exceptional locus after all the $l-1$ blowups is a chain of $l-1$ new divisors $D_1, \dots, D_{l-1}$ that separate the original divisors $\{x=t=0\}$ and $\{y=t=0\}$, see Figure \ref{asingresolve} for a sketch of the central fibres of the local model and the resolution corresponding to Figure \ref{model1}. The new divisors are $\mathbb{P}^1$-bundles with exceptional curves $F_i, \ 1 \leq i \leq l-1$.

\begin{figure}[t]
\begin{center}

\begin{tikzpicture}[line cap=round,line join=round,>=triangle 45,x=1cm,y=1cm]

\begin{scope}[scale = 2]
\draw (0,1)-- (0,0);
\draw (0,1)-- (-0.5,1.5);
\draw (0,0)-- (-0.5,-0.5);
\draw (0,1)-- (0.4,1);
\draw (0.4,1)-- (0.4,0);
\draw (0.4,0)-- (0,0);
\draw (0.4,0)-- (0.8,0);
\draw (0.8,0)-- (0.8,1);
\draw (0.8,1)-- (0.4,1);
\draw (1.2,1)-- (1.2,0);
\draw (1.2,0)-- (0.8,0);
\draw (1.2,1)-- (0.8,1);

\draw (1.2,1)-- (1.7,1.5);
\draw (1.2,0)-- (1.7,-0.5);
\begin{scriptsize}

\end{scriptsize}
\begin{scriptsize}

\draw(1.7, 0.3) node{$y=t=0$};
\draw(-0.5, 0.3) node{$x=t=0$};
\draw(0.2, 0.3) node{$D_1$};
\draw(0.6, 0.3) node{$D_2$};
\draw(1, 0.3) node{$D_3$};
\draw [color=red, line width=1pt] (0,0.6)--(0.4,0.5);
\draw [color=red, line width=1pt] (0.4,0.5)--(0.8,0.6);
\draw [color=red, line width=1pt] (0.8,0.6)--(1.2,0.5);
\draw(1, 0.65) node{$F_3$};
\draw(0.6, 0.65) node{$F_2$};
\draw(0.2, 0.65) node{$F_1$};

\end{scriptsize}
\end{scope}

\begin{scope}[xshift=7cm, scale = 2]
 \clip(-1.5068464964661663,-0.6875061121455801) rectangle (1.5626116514137933,1.7468917292764545);
\draw (0,1)-- (0,0);
\draw (0,1)-- (0.5,1.5);
\draw (0,1)-- (-0.5,1.5);
\draw (0,0)-- (-0.5,-0.5);
\draw (0,0)-- (0.5,-0.5);
\begin{scriptsize}
\draw(0.5, 0.3) node{$y=t=0$};
\draw(-0.5, 0.3) node{$x=t=0$};
\end{scriptsize}
\end{scope}

\draw[->] (4,1.3) -- (5.5,1.3);

\end{tikzpicture}

\end{center}
\caption{Blowing up $\left\{xy = t^4 \right\} \subseteq \Spec \kk[x,y,w_{\rho}]\llbracket t \rrbracket$.} \label{asingresolve}
\end{figure}

Similarly, the local model $\left\{xy = t^lw_{\rho} \right\} \subseteq \Spec \kk[x,y,w_{\rho}]\llbracket t \rrbracket$ can be resolved to a simple normal crossings and log smooth degeneration by blowing up the divisors $\{x=t = 0\}$ or $\{x'=t = 0\}$ at each step $l$ times (using Notation \ref{divisorlocalnotation}). Indeed, suppose that $l>1$ (otherwise, we are in the case of an ODP singularity of Section \ref{locmod}). Blowing up $\{x=t=0\}$ (without loss of generality) produces a chart that is simple normal crossings and a chart where the local model takes form $\left\{xy = t^{l-1}w_{\rho} \right\} \subseteq \Spec \kk[x,y,w_{\rho}]\llbracket t \rrbracket$. We can then blow up $\left\{x=t=0 \right\}$ again or we can blow up $\left\{x' = t =0 \right\}$. Repeating this procedure $l-1$ times (with choices of divisors to blow up at every step) produces a degeneration that is simple normal crossings away from an ODP singularity as in Section \ref{locmod}. We then blow up one of the divisors in the local model of the ODP singularity to obtain a log smooth and simple normal crossings degeneration. 

Torically, $\left\{xy = t^l w_{\rho} \right\} \subseteq \Spec \kk[x,y,w_{\rho}]\llbracket t \rrbracket$ is (the completion in $t$ of) the affine variety defined by the cone over the trapezoid that is the convex hull of $(0, 0, 1), (0, 1, 1), (l, 0, 1)$, and $ (1, 1, 1)$ (i.e. it is of the same type as Figure \ref{admisspolytopesresolve}(2)). Again, we assume that $(0, 0, 1)$ corresponds to $\{x = t = 0\}$, $(l, 0, 1)$ corresponds to $\{y = t = 0\}$, and that we first blow up $\{x = t = 0\}$. Suppose that we have blown up $\{x = t = 0\}$ $i_x$ times and $\{x' = t = 0\}$ $i_{x'}$ times (with $i_x + i_{x'} \leq l-1$), using Notation \ref{divisorlocalnotation}. Then this partial resolution corresponds to the subdivision of the trapezoid by the segments connecting $(0, 1, 1)$ to $(1, 0, 1), (2, 0, 1), \dots, (i_{x'}, 0, 1)$ and connecting $(1,1,1)$ to $(l-1, 0, 1), \dots, (l-i_x, 0, 1)$. If $i_x+i_{x'}= l-1$, the final blowup resolving the ODP singularity subdivides the square that is the convex hull of $(0, 1, 1), (1,1,1), (i_{x'},0,1), (l-i_{x}, 0,1)$ by the diagonal connecting $(i_{x'},0, 1)$ to $(1,1,1)$ if we blow up $\{x = t = 0\}$ and by the diagonal connecting $(l-i_x,0, 1)$ to $(0,1,1)$ if we blow up $\{x' = t = 0\}$.

Note that, unlike the resolution of $\left\{xy = t^l \right\} \subseteq \Spec \kk[x,y,w_{\rho}]\llbracket t \rrbracket$ described above, the final resolution depends on the choices of divisors we blow up. However, it does not depend on the order of the blowups and only depends on the total number $i_x^{\tot}$ of times that we blow up $\{x = t = 0\}$ and the total number $i_{x'}^{\tot}$ of times that we blow up $\{x' = t = 0\}$ (here $i_x^{\tot}+ i_{x'}^{\tot} = l$).
See Figure \ref{model2} for the subdivisions in the case that $l=4, \ i_x^{\tot} = 2, \ i_{x'}^{\tot} = 2$.

\begin{figure}[H]
\begin{center}

\begin{minipage}{.3\textwidth}
\begin{tikzpicture}[line cap=round,line join=round,>=triangle 45,x=1cm,y=1cm, scale = 3]
\clip(-0.22,-0.22) rectangle (1.25,1.25);
\draw  (0,0)-- (1,0);
\draw  (0,0)-- (0,1);
\draw  (0,1)-- (0.25,1);
\draw  (0.25,1)-- (1,0);
\draw [fill] (0,0) circle (0.5pt);
\draw [fill] (0.25,0) circle (0.5pt);
\draw [fill] (0.5,0) circle (0.5pt);
\draw [fill] (0.75,0) circle (0.5pt);
\draw [fill] (0,1) circle (0.5pt);
\draw [fill] (1,0) circle (0.5pt);
\draw [fill] (0.25,1) circle (0.5pt);
\draw(0, -0.1) node{\tiny $x=t=0$};
\draw(1, -0.1) node{\tiny $y=t=0$};
\draw(0.025, 1.07) node{\tiny $x=w_{\rho} =0$};
\draw(0.575, 1.07) node{\tiny $y=w_{\rho} =0$};
\end{tikzpicture}
\captionsetup{labelformat=empty}
\caption{$xy=t^4 w_{\rho}$ \textcolor{white}{llllll}}
\end{minipage}
\hspace{0.3cm}
\begin{minipage}{.3\textwidth}
\begin{tikzpicture}[line cap=round,line join=round,>=triangle 45,x=1cm,y=1cm, scale = 3]
\clip(-0.22,-0.22) rectangle (1.25,1.25);
\draw  (0,0)-- (1,0);
\draw  (0,0)-- (0,1);
\draw  (0,1)-- (0.25,1);
\draw  (0.25,1)-- (1,0);
\draw  (0,1)-- (0.25,0);
\draw  (0,1)-- (0.5,0);
\draw  (0.25,1)-- (0.75,0);
\draw [fill] (0,0) circle (0.5pt);
\draw [fill] (0.25,0) circle (0.5pt);
\draw [fill] (0.5,0) circle (0.5pt);
\draw [fill] (0.75,0) circle (0.5pt);
\draw [fill] (0,1) circle (0.5pt);
\draw [fill] (1,0) circle (0.5pt);
\draw [fill] (0.25,1) circle (0.5pt);
\draw(0, -0.1) node{\tiny $x=t=0$};
\draw(1, -0.1) node{\tiny $y=t=0$};
\draw(0.025, 1.07) node{\tiny $x=w_{\rho} =0$};
\draw(0.575, 1.07) node{\tiny $y=w_{\rho} =0$};

\end{tikzpicture}
\captionsetup{labelformat=empty}
\caption{\ \ \ After \\ 3 blowups}
\end{minipage}
\hspace{0.3cm}
\begin{minipage}{.3\textwidth}
\begin{tikzpicture}[line cap=round,line join=round,>=triangle 45,x=1cm,y=1cm, scale = 3]
\clip(-0.22,-0.22) rectangle (1.25,1.25);
\draw  (0,0)-- (1,0);
\draw  (0,0)-- (0,1);
\draw  (0,1)-- (0.25,1);
\draw  (0.25,1)-- (1,0);
\draw  (0,1)-- (0.25,0);
\draw  (0,1)-- (0.5,0);
\draw  (0.25,1)-- (0.75,0);
\draw  (0.25,1)-- (0.5,0);
\draw [fill] (0,0) circle (0.5pt);
\draw [fill] (0.25,0) circle (0.5pt);
\draw [fill] (0.5,0) circle (0.5pt);
\draw [fill] (0.75,0) circle (0.5pt);
\draw [fill] (0,1) circle (0.5pt);
\draw [fill] (1,0) circle (0.5pt);
\draw [fill] (0.25,1) circle (0.5pt);
\draw(0, -0.1) node{\tiny $x=t=0$};
\draw(1, -0.1) node{\tiny $y=t=0$};
\draw(0.025, 1.07) node{\tiny $x=w_{\rho} =0$};
\draw(0.575, 1.07) node{\tiny $y=w_{\rho} =0$};

\end{tikzpicture}
\captionsetup{labelformat=empty}
\caption{\ \ \ After \\ 4 blowups}
\end{minipage}
\end{center}
\setcounter{figure}{9}
\caption{Resolution of $\left\{xy = t^4 w_{\rho} \right\} \subseteq \Spec \kk[x,y,w_{\rho}]\llbracket t \rrbracket$.} \label{model2}
\end{figure}

The exceptional locus after all the $l$ blowups is a chain of $l-1$ new divisors $D_1, \dots, D_{l-1}$ that separate the strict transforms of the original divisors $\{x=t=0\}$ and $\{y=t=0\}$, along with an exceptional curve $E$ (possibly contained in one of these divisors). 

The divisors $D_{i}, \ 1 \leq i \leq l-1$ with $i \neq i_{x'}^{\tot}$ are $\mathbb{P}^1$-bundles with exceptional curves $F_i$. The last divisor blown up is $D_{i_{x'}^{\tot}}$ (interpreted as $\{x=t=0\}$ if $i_{x'}^{\tot} = 0$ and as $\{y=t=0\}$ if $i_{x'}^{\tot} = l$). If $i_{x'}^{\tot} \neq  0, l$ then $D_{i_{x'}^{\tot}}$ is the strict transform of a $\mathbb{P}^1$-bundle (with an exceptional curve $F_{i_{x'}^{\tot}}$) under the last blowup and has an exceptional curve $E$ intersecting a curve of class $F_{i_{x'}^{\tot}} - E$ (the strict transform of the exceptional curve of the $\mathbb{P}^1$-bundle meeting the singularity) at one point.\footnote{Note that even though the resolution does not depend on the order of the blowups, our notations for the exceptional curves reflect which divisor was blown up last.} If $i_{x'}^{\tot} = 0$, then we have an exceptional curve $E$ contained in the strict transform of $\{x = t = 0\}$ that meets $D_1$ at one point. If $i_{x'}^{\tot} = l$, then we have an exceptional curve $E$ contained in the strict transform of $\{y = t = 0\}$ that meets $D_{l-1}$ at one point. See Figure \ref{asingresolve1} for a sketch of the central fibres of the local model and the resolution corresponding to Figure \ref{model2}.

\begin{figure}[H]
\begin{center}

\begin{tikzpicture}[line cap=round,line join=round,>=triangle 45,x=1cm,y=1cm]

\begin{scope}[scale = 2]
\draw (0,1)-- (0,0);
\draw (0,1)-- (-0.5,1.5);
\draw (0,0)-- (-0.5,-0.5);
\draw (0,1)-- (0.4,1);
\draw (0.4,1)-- (0.4,0);
\draw (0.4,0)-- (0,0);
\draw (0.4,0)-- (0.8,0);
\draw (1.1,0)-- (1.1,1);
\draw (0.8,1)-- (0.4,1);
\draw (1.5,1)-- (1.5,0);
\draw (1.5,0)-- (0.8,0);
\draw (1.5,1)-- (0.8,1);

\draw (1.5,1)-- (2,1.5);
\draw (1.5,0)-- (2,-0.5);
\begin{scriptsize}

\end{scriptsize}
\begin{scriptsize}

\draw(2, 0.3) node{$y=t=0$};
\draw(-0.5, 0.3) node{$x=t=0$};
\draw(0.2, 0.3) node{$D_1$};
\draw(0.75, 0.3) node{$D_2$};
\draw(1.3, 0.3) node{$D_3$};
\draw [color=red, line width=1pt] (1.1,0.5)-- ++(-9pt,7.5pt);
\draw(1.01, 0.7) node{$E$};
\draw [color=red, line width=1pt] (0,0.6)--(0.4,0.5);
\draw [color=red, line width=1pt] (0.4,0.5)--(0.87,0.75);
\draw [color=red, line width=1pt] (1.1,0.5)--(1.5,0.6);
\draw(1.3, 0.65) node{$F_3$};
\draw(0.2, 0.65) node{$F_1$};

\draw[latex-,dashed] 
        (0.55,0.65) -- (0.55, 1.3)
         node[fill=white]{$F_2 -E$};

\end{scriptsize}
\end{scope}

\begin{scope}[xshift=7cm, scale = 2]
 \clip(-1.5068464964661663,-0.6875061121455801) rectangle (1.5626116514137933,1.7468917292764545);
\draw (0,1)-- (0,0);
\draw (0,1)-- (0.5,1.5);
\draw (0,1)-- (-0.5,1.5);
\draw (0,0)-- (-0.5,-0.5);
\draw (0,0)-- (0.5,-0.5);
\begin{scriptsize}
\draw [color=blue, line width = 1pt] (0, 0.5)-- ++(-1.5pt,-1.5pt) -- ++(3pt,3pt) ++(-3pt,0) -- ++(3pt,-3pt);
\draw(0.5, 0.3) node{$y=t=0$};
\draw(-0.5, 0.3) node{$x=t=0$};
\end{scriptsize}
\end{scope}

\draw[->] (4,1.3) -- (5.5,1.3);

\end{tikzpicture}

\end{center}
\caption{Blowing up $\left\{xy = t^4 w_{\rho} \right\} \subseteq \Spec \kk[x,y,w_{\rho}]\llbracket t \rrbracket$.} \label{asingresolve1}
\end{figure}

\begin{sloppypar}
\begin{remark}
Note that performing a similar construction in the case of the more general divisorial singularity with local model $\left\{xy = t^lw_{\rho}^k \right\} \subseteq \Spec \kk[x,y,w_{\rho}]\llbracket t \rrbracket$ for $k \geq 2$ (see Observation \ref{modelslocal}(1)) is not possible since there is no way to resolve the singularity in the generic fibre of this local model by blowing up irreducible components of the central fibre. This is why we require that a special toric degeneration $\bar{\mathfrak{X}} \to \mathcal{S}$ has a smooth generic fibre (see condition (1) of Definition \ref{special}).
\end{remark}
\end{sloppypar}

\subsubsection{Blowing up a divisor} \label{divisorblowupgeneral}

We want to obtain a global resolution of $\bar{\mathfrak{X}} \to \mathcal{S}$ by blowing up a sequence of irreducible components of the central fibre $\bar{\mathfrak{X}}_0$ (and irreducible components of the proper transform of $\bar{\mathfrak{X}}_0$ under the blowups). Let $\bar{D}_i$ for some $  1 \leq i \leq \bar{m}$ be a component of the central fibre. As in Section \ref{divisorblowup}, blowing up $\bar{D}_i$ corresponds to blowing up either $\{x=t=0\}$ or $\{y=t=0\}$ in the local models of all the (singular and non-singular) points contained in every double curve $\bar{X}_{\rho} \subseteq \bar{D}_i$ (with a compatible choice for the points contained in the same $\bar{X}_{\rho}, \rho \in \bar{\mathscr{P}}^{[1]}$). Note that the number of singular points in $\bar{X}_{\rho}, \rho \in \bar{\mathscr{P}}^{[1]}$ is equal to $r_{\rho}$ (the index of the singularity $x_{\rho} \subseteq \Int(\rho)$ of the affine structure). 

\begin{notation}
From now on, we shall denote the integral length of an edge $\rho \in \bar{\mathscr{P}}^{[1]}$ by $l_{\rho}$.    
\end{notation}

Fix a $\rho \in \bar{\mathscr{P}}^{[1]}$ with $\bar{X}_{\rho} \subseteq \bar{D}_i$. If $l_{\rho} = 1$, then the proper transform of each singular point of $\bar{X}_{\rho} \subseteq \bar{D}_i$ is an exceptional curve contained in the strict transform of $\bar{D}_i$ as in Section \ref{divisorblowup}. If $l_{\rho}>1$, then the resolution is not small on $\bar{X}_{\rho}$ and the proper transform of $\bar{X}_{\rho}$ gives a $\mathbb{P}^1$-bundle. This is a globalization of the first blowup in the local models of Section \ref{localmodelgen}.

\begin{notations} \label{blowupnotationsgeneral}
We expand Notations \ref{blowupnotations} to this case.
\begin{enumerate}
    \item We use Notations \ref{blowupnotations} for the strict transforms of the irreducible components of $\bar{\mathfrak{X}}_0$. However, we denote the strict transform of $\bar{D}_i$ by $\bar{D}_i$ (and not $D_i$) unless the degeneration becomes log smooth in the neighbourhood of the strict transform after blowing up $\bar{D}_i$ (we will also use the bar if speaking of the general situation). 
    \item We use Notations \ref{blowupnotations} for the exceptional curves $E_k^{ij}$ (or $E_{\rho, k}$) contained in the strict transform of $\bar{D}_i$.
    \item For every $\rho \in \bar{\mathscr{P}}^{[1]}$ with $l_{\rho}>1$ and $\bar{X}_{\rho} = \bar{D}_i \cap \bar{D}_j$, we denote by $\bar{D}^{ij}_1$ (or by $D^{ij}_1$ if the degeneration is log smooth in the neighbourhood of that component) the $\mathbb{P}^1$-bundle contained in the proper transform of $\bar{X}_{\rho}$ and we denote its exceptional curve by $F_1^{ij}$. Alternatively, we use the notations $\bar{D}_{\rho, 1}$ and $F_{\rho, 1}$.
\end{enumerate}
\end{notations}

Suppose that $\bar{X}_{\rho} \subseteq \bar{D}_i$ and $l_{\rho} > 1$. Then the points of $\bar{D}_{\rho, 1} \cap \bar{D}_i$ have local models of the form $\left\{xy = t^{l_{\rho} - 1} \right\} \subseteq \Spec \kk[x,y,w_{\rho}]\llbracket t \rrbracket$ at the strict transforms of non-singular points of $\bar{X}_{\rho}$ and local models of the form $\left\{xy = t^{l_{\rho} - 1}w_{\rho} \right\} \subseteq \Spec \kk[x,y,w_{\rho}]\llbracket t \rrbracket$ at strict transforms of the singularities of $\bar{X}_{\rho}$. Therefore, we have improved the singularities.

Let $g': \mathfrak{X}' \to \mathcal{S}$ be the generically log smooth (in the sense of Definition \ref{genlogsmooth}) degeneration obtained after blowing up $\bar{D}_i$. To better understand the global structure of the central fibre $\mathfrak{X}'_{0}$ of $g'$, we can use the local models at the $0$-dimensional strata of $\bar{D}_i$. 

Consider a stratum $\bar{X}_{\sigma} \subseteq \bar{D}_i$ with $\sigma \in \bar{\mathscr{P}}^{\max}$. By condition (4)(a) in Definition \ref{toricdegeneration} of a toric degeneration, the point $\bar{X}_{\sigma}$ does not lie in the singular locus $Z$ of $\bar{\mathfrak{X}} \to \mathcal{S}$. This implies, via Construction \ref{tropicalize} of tropicalization, that there is a local toric model at $\bar{X}_{\sigma}$ given by the cone over the polyhedron $\sigma$. The blowup of $\bar{D}_i$ corresponds to blowing up the divisor corresponding to $\bar{v}_i$ (following Notation \ref{blowupdicnotations}) in the local model. We have described such local models and their blowups in Section \ref{toricmodelcodim2} (for the allowed class of $\sigma \in \bar{\mathscr{P}}^{\max}$). Note from Figure \ref{admisspolytopesresolve} that blowing up $\bar{D}_i$ can give rise to new toric irreducible components and new curve classes in the central fibre.

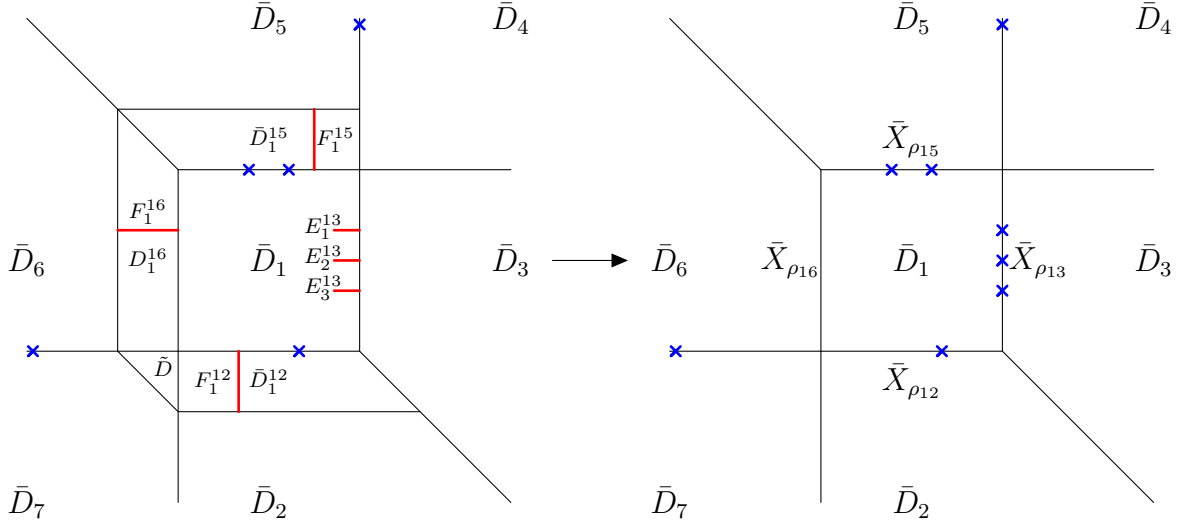
\begin{figure}[t]
\begin{center}

\hspace{-0.5cm}
\begin{tikzpicture}[line cap=round,line join=round,>=triangle 45,x=1cm,y=1cm]

\begin{scope}[scale = 0.8]
\clip(-4.3,-4.3) rectangle (4.3,4.3);
\draw  (-1.5,1.5)-- (1.5,1.5);
\draw  (1.5,1.5)-- (1.5,-1.5);
\draw  (-1.5,1.5)-- (-1.5,-1.5);
\draw  (-1.5,-1.5)-- (1.5,-1.5);
\draw  (-1.5,1.5)-- (-4,4);
\draw  (1.5,1.5)-- (1.5,4);
\draw  (1.5,1.5)-- (4,1.5);
\draw  (1.5,-1.5)-- (4,-4);
\draw  (-1.5,-1.5)-- (-1.5,-4);
\draw  (-1.5,-1.5)-- (-4,-1.5);
\draw  (-2.5, 2.5)-- (1.5, 2.5);
\draw  (-2.5, 2.5)-- (-2.5, -1.5);
\draw  (-1.5, -2.5)-- (2.5, -2.5);
\draw  (-1.5, -2.5)-- (-2.5, -1.5);

\begin{scriptsize}
\draw(0,0) node{\normalsize $\bar{D}_1$};
\draw(0,-4) node{\normalsize $\bar{D}_2$};
\draw(4,0) node{\normalsize $\bar{D}_3$};
\draw(4,4) node{\normalsize $\bar{D}_4$};
\draw(0,4) node{\normalsize $\bar{D}_5$};
\draw(-4,0) node{\normalsize $\bar{D}_6$};
\draw(-4,-4) node{\normalsize $\bar{D}_7$};
\draw [color=red, line width=1pt] (1.5,0.5)-- ++(-12pt,0pt);
\draw(0.9, 0.55) node{\tiny $E_1^{13}$};
\draw [color=red, line width=1pt] (1.5,0)-- ++(-12pt,0pt);
\draw(0.9, 0.05) node{\tiny $E_2^{13}$};
\draw [color=red, line width=1pt] (1.5,-0.5)-- ++(-12pt,0pt);
\draw(0.9, -0.45) node{\tiny $E_3^{13}$};

\draw [color=blue, line width = 1pt] (0.33, 1.5)-- ++(-2.25pt,-2.25pt) -- ++(4.5pt,4.5pt) ++(-4.5pt,0) -- ++(4.5pt,-4.5pt);
\draw [color=blue, line width = 1pt] (-0.33, 1.5)-- ++(-2.25pt,-2.25pt) -- ++(4.5pt,4.5pt) ++(-4.5pt,0) -- ++(4.5pt,-4.5pt);

\draw [color=red, line width=1pt] (0.75,1.5)-- (0.75,2.5);
\draw(1.1, 2) node{\tiny $F_1^{15}$};
\draw(0, 2) node{\tiny $\bar{D}_1^{15}$};

\draw [color=red, line width=1pt] (-2.5,0.5)-- (-1.5,0.5);
\draw(-2, 0.8) node{\tiny $F_1^{16}$};
\draw(-2, 0) node{\tiny $D_1^{16}$};
\draw(-1.75, -1.75) node{\tiny $\tilde{D}$};

\draw [color=blue, line width = 1pt] (0.5, -1.5)-- ++(-2.25pt,-2.25pt) -- ++(4.5pt,4.5pt) ++(-4.5pt,0) -- ++(4.5pt,-4.5pt);
\draw [color=red, line width=1pt] (-0.5,-1.5)-- (-0.5,-2.5);
\draw(-0.92, -2) node{\tiny $F_1^{12}$};
\draw(0, -2) node{\tiny $\bar{D}_1^{12}$};

\draw [color=blue, line width = 1pt] (1.5, 3.9)-- ++(-2.25pt,-2.25pt) -- ++(4.5pt,4.5pt) ++(-4.5pt,0) -- ++(4.5pt,-4.5pt);

\draw [color=blue, line width = 1pt] (-3.9, -1.5)-- ++(-2.25pt,-2.25pt) -- ++(4.5pt,4.5pt) ++(-4.5pt,0) -- ++(4.5pt,-4.5pt);

\end{scriptsize}
\end{scope}

\begin{scope}[xshift=8.5cm, scale = 0.8]
\clip(-4.3,-4.3) rectangle (4.3,4.3);
\draw  (-1.5,1.5)-- (1.5,1.5);
\draw  (1.5,1.5)-- (1.5,-1.5);
\draw  (-1.5,1.5)-- (-1.5,-1.5);
\draw  (-1.5,-1.5)-- (1.5,-1.5);
\draw  (-1.5,1.5)-- (-4,4);
\draw  (1.5,1.5)-- (1.5,4);
\draw  (1.5,1.5)-- (4,1.5);
\draw  (1.5,-1.5)-- (4,-4);
\draw  (-1.5,-1.5)-- (-4,-1.5);
\draw  (-1.5,-1.5)-- (-1.5,-4);
\begin{scriptsize}
\draw [color=blue, line width = 1pt] (1.5, 0.5)-- ++(-2.25pt,-2.25pt) -- ++(4.5pt,4.5pt) ++(-4.5pt,0) -- ++(4.5pt,-4.5pt);
\draw [color=blue, line width = 1pt] (1.5, 0)-- ++(-2.25pt,-2.25pt) -- ++(4.5pt,4.5pt) ++(-4.5pt,0) -- ++(4.5pt,-4.5pt);
\draw [color=blue, line width = 1pt] (1.5, -0.5)-- ++(-2.25pt,-2.25pt) -- ++(4.5pt,4.5pt) ++(-4.5pt,0) -- ++(4.5pt,-4.5pt);

\draw [color=blue, line width = 1pt] (0.33, 1.5)-- ++(-2.25pt,-2.25pt) -- ++(4.5pt,4.5pt) ++(-4.5pt,0) -- ++(4.5pt,-4.5pt);
\draw [color=blue, line width = 1pt] (-0.33, 1.5)-- ++(-2.25pt,-2.25pt) -- ++(4.5pt,4.5pt) ++(-4.5pt,0) -- ++(4.5pt,-4.5pt);

\draw [color=blue, line width = 1pt] (0.5, -1.5)-- ++(-2.25pt,-2.25pt) -- ++(4.5pt,4.5pt) ++(-4.5pt,0) -- ++(4.5pt,-4.5pt);

\draw(0,0) node{\normalsize $\bar{D}_1$};
\draw(0,-4) node{\normalsize $\bar{D}_2$};
\draw(4,0) node{\normalsize $\bar{D}_3$};
\draw(4,4) node{\normalsize $\bar{D}_4$};
\draw(0,4) node{\normalsize $\bar{D}_5$};
\draw(-4,0) node{\normalsize $\bar{D}_6$};
\draw(-4,-4) node{\normalsize $\bar{D}_7$};

\draw(2.1,0) node{\normalsize $\bar{X}_{\rho_{13}}$};
\draw(0,2) node{\normalsize $\bar{X}_{\rho_{15}}$};
\draw(0,-2) node{\normalsize $\bar{X}_{\rho_{12}}$};
\draw(-2,0) node{\normalsize $\bar{X}_{\rho_{16}}$};

\draw [color=blue, line width = 1pt] (1.5, 3.9)-- ++(-2.25pt,-2.25pt) -- ++(4.5pt,4.5pt) ++(-4.5pt,0) -- ++(4.5pt,-4.5pt);

\draw [color=blue, line width = 1pt] (-3.9, -1.5)-- ++(-2.25pt,-2.25pt) -- ++(4.5pt,4.5pt) ++(-4.5pt,0) -- ++(4.5pt,-4.5pt);

\end{scriptsize}
\end{scope}

\draw[->] (3.75,0) -- (4.75,0);

\end{tikzpicture}

\end{center}
\caption{Blowing up the divisor $\bar{D}_1$. Here $l_{\rho_{13}}=1$, $l_{\rho_{12}}=3$, and $l_{\rho_{15}}=l_{\rho_{16}}=2$ (one can only see that $l_{\rho_{12}}, l_{\rho_{15}}, l_{\rho_{16}} > 1$ from this sketch).} \label{blowingadivisor}
\end{figure}

We give an example of a divisor blowup in Figure \ref{blowingadivisor} (obtained using the analysis of the local models of Sections \ref{toricmodelcodim2} and \ref{localmodelgen}). The blowup resolves the ODP singularities of $\bar{X}_{\rho_{13}}$ producing exceptional curves $E_{1}^{13}, E_{2}^{13}, E_3^{13}$ and improves the other singularities, producing three $\mathbb{P}^1$-bundles $\bar{D}_1^{12}, \bar{D}_1^{15}, \bar{D}_{1}^{16}$ with exceptional curves $F_1^{12}, F_1^{15}, F_{1}^{16}$ respectively (it also introduces an additional toric irreducible component, isomorphic to $\mathbb{P}^2$, that we denote by $\tilde{D}$).

We explain the effect of blowing up $\bar{D}_i$ on the dual intersection complex $\left( \bar{B}, \bar{\mathscr{P}} \right)$. As usual, we define the dual intersection complex $\left(B', \mathscr{P}'\right)$ of $\mathfrak{X}' \overset{g'}{\to} \mathcal{S}$ as $(g')^{-1}_{\trop}(1)$ for $g'_{\trop}: \Sigma(\mathfrak{X}') \to \mathbb{R}_{\geq 0}$ the tropicalization of $g'$, with the polyhedral structure $\mathscr{P}'$ obtained by restricting the cones of $\Sigma \left(\mathfrak{X}'\right)$ to the fibre over $1 \in \mathbb{R}_{\geq 0}$. 

\begin{notation} \label{vertexnot}
Expanding Notation \ref{blowupdicnotations}, we denote by $\bar{v}_i \in \mathscr{P}'^{[0]}$ (resp. $v_i \in \mathscr{P}'^{[0]}$) the vertex corresponding to $\bar{D}_i$ (resp. $D_i$) and by $\bar{v}_k \in \mathscr{P}'^{[0]}$ the vertices corresponding to $\bar{D}_k$ (for $1 \leq k \leq \bar{m}, \ k \neq i$). For every $\rho \in \bar{\mathscr{P}}^{[1]}$ with $l_{\rho}>1$ and $\bar{X}_{\rho} = \bar{D}_i \cap \bar{D}_j$, we denote by $\bar{v}_{1}^{ij} \in \mathscr{P}'^{[0]}$ the vertex corresponding to $\bar{D}^{ij}_1$. Alternatively, we denote it by $\bar{v}_{\rho, 1}$.
\end{notation}

Unlike Section \ref{divisorblowup}, it is not the case that $\left(B', \mathscr{P}'\right) \cong \left(\bar{B}, \bar{\mathscr{P}}\right)$ as polyhedral manifolds unless all the maximal cells $\sigma \in \bar{\mathscr{P}}^{\max}$ with $\bar{v}_i \subseteq \sigma$ are standard triangles. Instead, $\left(B', \mathscr{P}'\right)$ is a natural subdivision of $\left(\bar{B}, \bar{\mathscr{P}}\right)$ and the analysis of the local models in codimension $2$ of Section \ref{toricmodelcodim2} implies that the cells $\sigma \in \bar{\mathscr{P}}^{\max}$ are subdivided according to Figure \ref{admisspolytopesresolve}. Since $\left(\bar{B}, \bar{\mathscr{P}}\right)$ is a polyhedral manifold of dimension $2$ in the sense of Definition \ref{polymani} (see Proposition \ref{dicpm}), this implies that $\left(B', \mathscr{P}'\right)$ is also a polyhedral manifold of dimension $2$ in the sense of Definition \ref{polymani}. 

Now, we can give $\left(B', \mathscr{P}'\right)$ the structure of an affine manifold with singularities using Construction \ref{affstrintgen}. Note from Remark \ref{genaffstrremarks}(3) that this defines the usual toric charts of Construction \ref{affstrtoric} on $W_{\bar{v}_k}$ for $1 \leq k \leq \bar{m}, \ k \neq i$ and the intrinsic affine structure on $W_{\bar{v}_i} \setminus \bar{v}_i$. We also have an intrinsic affine structure on $W_{\bar{v}_{\rho,1}} \setminus \bar{v}_{\rho,1}$ for all $\rho \in \bar{\mathscr{P}}^{[1]}$ with $\bar{X}_{\rho} \subseteq \bar{D}_i$ and $l_{\rho}>1$. By Remark \ref{genaffstrremarks}(3), it extends to the whole $W_{\bar{v}_{\rho,1}}$ since $\bar{D}_{\rho,1}$ is a toric variety and its log stratification (in the sense of Definition \ref{generalizedlogstrata}) coincides with its toric stratification. In this case, it is easy to see directly that in Construction \ref{affstrintgen} the embeddings $\psi_{\rho_{\bar{v}_{\rho,1}}}$ of \eqref{embeddinglocalgen} over the codimension $1$ cones $\rho_{\bar{v}_{\rho,1}} \in \Sigma(D_{\bar{v}_{\rho,1}})$ glue to a global embedding $\Sigma(D_{\bar{v}_{\rho,1}}) \to \mathbb{R}^2$ (inducing a chart on $W_{\bar{v}_{\rho,1}}$). Similarly, the intrinsic affine structure on $W_{\bar{v}_i} \setminus \bar{v}_i$ extends to the whole $W_{\bar{v}_i}$ in the case that $l_{\rho} > 1$ for all $\rho \in \bar{\mathscr{P}}^{[1]}$ with $\bar{X}_{\rho} \subseteq \bar{D}_i$.

This analysis, along with the description of the blowups in the local models of Section \ref{localmodelgen}, implies that at the level of the affine manifolds, the blowup of $\bar{D}_i$ can be visualized as follows. First, the dual intersection complex $\left( \bar{B}, \bar{\mathscr{P}}\right)$ is subdivided to $\left( B', \mathscr{P}' \right)$ as described above. Let $\rho \in \bar{\mathscr{P}}^{[1]}$ be an edge with $\bar{v}_i \subseteq \rho$ and $r_\rho \neq 0$ (so that there is a singularity $x_{\rho} \in \Int(\rho)$). If $l_{\rho} = 1$, then $x_{\rho}$ is pulled into $\bar{v}_i$ as in Section \ref{divisorblowup}. If $l_{\rho} > 1$ then $\rho$ is subdivided into two edges $\rho_{\bar{v}_1}, \rho' \in \mathscr{P}'^{[1]}$ with $\bar{v}_i \subseteq \rho_{\bar{v}_i}$ and $l_{\rho_{\bar{v}_i}} = l-1, \ l_{\rho'} = 1$. The singularity $x_{\rho}$ is, again, pulled into $\bar{v}_i$, i.e. it now lies in $\Int(\rho_{\bar{v}_i})$ (the monodromy of the new singularity is still $r_{\rho}$). 
We give the transformation corresponding to the blowup of Figure \ref{blowingadivisor} in Figure \ref{blowingadivisoraff}.

\begin{figure}[t]
\begin{center}

\begin{tikzpicture}[line cap=round,line join=round,>=triangle 45,x=1cm,y=1cm]

\begin{scope}[scale = 1.2]
\clip(-2.5,-3.5) rectangle (2,2.5);
\draw (0,2)-- (1,2);
\draw (1,2)-- (1,0);
\draw (1,0)-- (0,-3);
\draw (0,-3)-- (-2,-3);
\draw (-2,0)-- (-2,-3);
\draw (-2,0)-- (0,2);
\draw (0,2)-- (0,0);
\draw (0,0)-- (1,0);
\draw (0,0)-- (0,-3);
\draw (0,0)-- (-2,0);
\draw (1,2)-- (-1,0);
\draw (-1,0)-- (-1,-2);
\draw (-1,-2)-- (-2,-3);
\draw (0,-2)-- (-1,-2);
\draw (0,-2)-- (1,0);
\begin{scriptsize}
\draw [fill] (0,0) circle (1pt);
\draw [fill] (0,1) circle (1pt);
\draw [fill] (1,0) circle (1pt);
\draw [fill] (1,1) circle (1pt);
\draw [fill] (0,2) circle (1pt);
\draw [fill] (1,2) circle (1pt);
\draw [fill] (-1,0) circle (1pt);
\draw [fill] (-2,0) circle (1pt);
\draw [fill] (-1,1) circle (1pt);
\draw [fill] (0,-1) circle (1pt);
\draw [fill] (0,-2) circle (1pt);
\draw [fill] (0,-3) circle (1pt);
\draw [fill] (-1,-1) circle (1pt);

\draw [fill] (-2,-1) circle (1pt);
\draw [fill] (-2,-2) circle (1pt);
\draw [fill] (-2,-3) circle (1pt);
\draw [fill] (-1,-2) circle (1pt);
\draw [fill] (-1,-3) circle (1pt);

\draw [color=blue, line width = 1pt] (0,0)-- ++(-2.25pt,-2.25pt) -- ++(4.5pt,4.5pt) ++(-4.5pt,0) -- ++(4.5pt,-4.5pt);
\draw [color=blue, line width = 1pt] (0,0.5)-- ++(-2.25pt,-2.25pt) -- ++(4.5pt,4.5pt) ++(-4.5pt,0) -- ++(4.5pt,-4.5pt);
\draw [color=blue, line width = 1pt] (0,-1.5)-- ++(-2.25pt,-2.25pt) -- ++(4.5pt,4.5pt) ++(-4.5pt,0) -- ++(4.5pt,-4.5pt);

\draw [color=blue, line width = 1pt] (0.5,2)-- ++(-2.25pt,-2.25pt) -- ++(4.5pt,4.5pt) ++(-4.5pt,0) -- ++(4.5pt,-4.5pt);
\draw [color=blue, line width = 1pt] (-2, -1.5)-- ++(-2.25pt,-2.25pt) -- ++(4.5pt,4.5pt) ++(-4.5pt,0) -- ++(4.5pt,-4.5pt);

\draw(-0.25, 0.2) node{\normalsize $\bar{v}_1$};
\draw(0, -3.2) node{\normalsize $\bar{v}_2$};
\draw(1.2, 0.2) node{\normalsize $\bar{v}_3$};
\draw(1.2, 2.15) node{\normalsize $\bar{v}_4$};
\draw(0.2, 2.15) node{\normalsize $\bar{v}_5$};
\draw(-2.2, 0.2) node{\normalsize $\bar{v}_6$};
\draw(-2.2, -3.15) node{\normalsize $\bar{v}_7$};

\draw(-1.2, 0.2) node{\tiny $v_1^{16}$};
\draw(0.2, 0.85) node{\tiny $\bar{v}_1^{15}$};
\draw(-0.2, -1.85) node{\tiny $\bar{v}_1^{12}$};
\draw(-1.5, 2) node{\normalsize $\left(B', \mathscr{P}'\right)$};

\end{scriptsize}
\end{scope}

\begin{scope}[xshift=8cm, scale = 1.2]
\clip(-3.5,-3.5) rectangle (2,2.5);
\draw (0,2)-- (1,2);
\draw (1,2)-- (1,0);
\draw (1,0)-- (0,-3);
\draw (0,-3)-- (-2,-3);
\draw (-2,0)-- (-2,-3);
\draw (-2,0)-- (0,2);
\draw (0,2)-- (0,0);
\draw (0,0)-- (1,0);
\draw (0,0)-- (0,-3);
\draw (0,0)-- (-2,0);
\begin{scriptsize}
\draw [fill] (0,0) circle (1pt);
\draw [fill] (0,1) circle (1pt);
\draw [fill] (1,0) circle (1pt);
\draw [fill] (1,1) circle (1pt);
\draw [fill] (0,2) circle (1pt);
\draw [fill] (1,2) circle (1pt);
\draw [fill] (-1,0) circle (1pt);
\draw [fill] (-2,0) circle (1pt);
\draw [fill] (-1,1) circle (1pt);
\draw [fill] (0,-1) circle (1pt);
\draw [fill] (0,-2) circle (1pt);
\draw [fill] (0,-3) circle (1pt);
\draw [fill] (-1,-1) circle (1pt);

\draw [fill] (-2,-1) circle (1pt);
\draw [fill] (-2,-2) circle (1pt);
\draw [fill] (-2,-3) circle (1pt);
\draw [fill] (-1,-2) circle (1pt);
\draw [fill] (-1,-3) circle (1pt);

\draw [color=blue, line width = 1pt] (0.5,0)-- ++(-2.25pt,-2.25pt) -- ++(4.5pt,4.5pt) ++(-4.5pt,0) -- ++(4.5pt,-4.5pt);
\draw [color=blue, line width = 1pt] (0,1.5)-- ++(-2.25pt,-2.25pt) -- ++(4.5pt,4.5pt) ++(-4.5pt,0) -- ++(4.5pt,-4.5pt);
\draw [color=blue, line width = 1pt] (0,-2.5)-- ++(-2.25pt,-2.25pt) -- ++(4.5pt,4.5pt) ++(-4.5pt,0) -- ++(4.5pt,-4.5pt);

\draw [color=blue, line width = 1pt] (0.5,2)-- ++(-2.25pt,-2.25pt) -- ++(4.5pt,4.5pt) ++(-4.5pt,0) -- ++(4.5pt,-4.5pt);
\draw [color=blue, line width = 1pt] (-2, -1.5)-- ++(-2.25pt,-2.25pt) -- ++(4.5pt,4.5pt) ++(-4.5pt,0) -- ++(4.5pt,-4.5pt);

\draw(-0.25, 0.2) node{\normalsize $\bar{v}_1$};
\draw(0, -3.2) node{\normalsize $\bar{v}_2$};
\draw(1.2, 0.2) node{\normalsize $\bar{v}_3$};
\draw(1.2, 2.15) node{\normalsize $\bar{v}_4$};
\draw(0.2, 2.15) node{\normalsize $\bar{v}_5$};
\draw(-2.2, 0.2) node{\normalsize $\bar{v}_6$};
\draw(-2.2, -3.15) node{\normalsize $\bar{v}_7$};

\draw(-1.5, 2) node{\normalsize $\left(\bar B, \bar{\mathscr{P}}\right)$};

\end{scriptsize}
\end{scope}

\draw[->] (2.5,0) -- (4,0);

\end{tikzpicture}

\end{center}
\caption{Transformation of $\left(\bar B, \bar{\mathscr{P}}\right)$ corresponding to blowing up $\bar{D}_1$.} \label{blowingadivisoraff}
\end{figure}

\subsubsection{Further blowing up components} \label{furthercompblowups}

Blowing up $\bar{D}_i$ improves the singularities of $\bar{\mathfrak{X}} \to \mathcal{S}$ contained in the strict transform of $\bar{D}_i$ but (unless $l_{\rho} = 1$ for all $\rho \in \bar{\mathscr{P}}^{[n-1]}$ with $\bar{v}_i \subseteq \rho$) it does not completely resolve them. Therefore, we need to continue blowing up components of the proper transform of $\bar{D}_i$ after the blowup to resolve the singularities.

Let $\rho \in \bar{\mathscr{P}}^{[1]}$ be an edge with $l_{\rho} = 2$. Then to resolve the ODP singularities of $\bar{D}_{\rho, 1} \cap \bar{D}_i$ it is enough to blow up either $\bar{D}_{\rho, 1}$ or $\bar{D}_i$ once (and the blowup map is small on $\bar{D}_{\rho, 1} \cap \bar{D}_i$). If $l_{\rho} > 3$ instead, then the blowup is not small on $\bar{D}_{\rho, 1} \cap \bar{D}_i$ and introduces a \enquote{new} $\mathbb{P}^1$-bundle. Moreover, the index of the singularities (now contained in the intersection of $\bar{D}_i$ with the \enquote{new} $\mathbb{P}^1$-bundle) decreases by $1$. These observations follow from analyzing the local models of Sections \ref{toricmodelcodim2} and \ref{localmodelgen}, similarly to the analysis of the previous section. The transformation of the dual intersection complex $\left(B', \mathscr{P}' \right)$ corresponding to this blowup is similar to the description at the end of the previous section (using the blown up local models for $\bar{X}_{\sigma}, \ \sigma \in \bar{\mathscr{P}}^{\max}$ and $\bar{X}_{\rho}, \ \rho \in \bar{\mathscr{P}}^{[1]}$). Namely, $\left(B', \mathscr{P}' \right)$ gets (further) subdivided according to Figure \ref{admisspolytopesresolve} (note that all the cells $\sigma \in \mathscr{P}'^{\max}$ are of one of the types in Figure \ref{admisspolytopes} by Corollary \ref{subwellbehaved}) and the singularities are pulled into the vertex corresponding to the component that we blow up.
One may then continue to blow up the components of the new central fibre (in particular, one may blow up the \enquote{new} $\mathbb{P}^1$-bundle).

Suppose that we have blown up some sequence of irreducible components of proper transforms of $\bar{\mathfrak{X}}_0$. This gives rise to a generically log smooth partial resolution $\mathfrak{X}'' \to \mathcal{S}$ of $\bar{\mathfrak{X}} \to \mathcal{S}$. This resolution has a well-defined dual intersection complex $\left(B'', \mathscr{P}''\right)$ that is a natural subdivision of $\left(\bar{B}, \bar{\mathscr{P}}\right)$ and has an affine structure given by Construction \ref{affstrintgen} that extends across the vertices corresponding to the toric irreducible components. 

\begin{notations} \label{generalblowupnotations}
Let $\mathfrak{X}'' \to \mathcal{S}$ be a generically log smooth partial resolution of $\bar{\mathfrak{X}} \to \mathcal{S}$ with dual intersection complex $\left(B'', \mathscr{P}''\right)$.
We fix the notations naturally expanding Notations \ref{blowupnotationsgeneral} and Notation \ref{vertexnot}. 
\begin{enumerate}
\item We denote an irreducible component of $\mathfrak{X}''_0$ with a bar unless the degeneration is log smooth in a neighbourhood of that component (and we will use the bar if speaking of the general situation). We use similar notation for the corresponding vertices $\bar{v} \in \mathscr{P}''^{[0]}$. Note that every $\bar{v} \in \mathscr{P}''^{[0]}$ corresponds to an element of $\bar{B}(\mathbb{Z})$. 
\item We write $\bar{D}_i$ (or $D_i$) for the strict transform of a component $\bar{D}_i$ of $\bar{\mathfrak{X}}_0$ and we write $\bar{v}_i$ (or $v_i$) for the corresponding vertex in $\mathscr{P}''^{[0]}$.
\item For an edge $\rho = \langle \bar{v}_i, \bar{v}_j \rangle \in \bar{\mathscr{P}}$ of length $l_{\rho}$ we write $D_{p}^{ij}, \ 1 \leq p \leq l_{\rho} -1$ (possibly not all divisors and possibly with bars) for the divisors corresponding to integral points of $\rho$ at distance $p$ from $\bar{v}_{j}$. We write $F_{p}^{ij}, \ 1 \leq p \leq l_{\rho} -1$ for the exceptional curves of the $\mathbb{P}^1$-bundles $\bar{D}_{p}^{ij}$ (and their strict transforms after further blowups). We may also use alternative notations $D_{\rho, p}$ and $F_{\rho, p}$. We write $v_p^{ij} \in \mathscr{P}''^{[1]}$ (or $v_{\rho, p} \in \mathscr{P}''^{[1]}$) for the vertex corresponding to $D_{p}^{ij}$ (similarly with bars). 
\item We denote the exceptional curves of the last blowup of one of the $D_{p}^{ij}, \ 1 \leq p \leq l_{\rho} -1$ (it could also be a blowup of $\bar{D}_i$ or $\bar{D}_j$) resolving the ODP-s by $E_{k}^{ij}$ (or $E_{\rho, k}$) for $1 \leq k \leq r_{\rho}$.
\end{enumerate}
\end{notations}

\begin{figure}[t]
\begin{center}

\begin{tikzpicture}[line cap=round,line join=round,>=triangle 45,x=1cm,y=1cm]

\begin{scope}[scale = 0.8]
\clip(-4.5,-4.5) rectangle (4.5,4.5);
\draw  (-1.5,1.5)-- (1.5,1.5);
\draw  (-1.5,1.5)-- (-1.5,-1.5);
\draw  (-1.5,-1.5)-- (2,-1.5);
\draw  (-1.5,1.5)-- (-4,4);
\draw  (1.5,1.5)-- (1.5,4);
\draw  (2,1)-- (4.5,1);
\draw  (2,-1.5)-- (4.5,-4);
\draw  (-1.5,-1.5)-- (-1.5,-4);
\draw  (-1.5,-1.5)-- (-4,-1.5);

\draw  (-3.5, 3.5)-- (1.5, 3.5);
\draw  (-3.5, 3.5)-- (-3.5, -1.5);
\draw  (-1.5, -3.5)-- (4, -3.5);
\draw  (-1.5, -2.5)-- (3, -2.5);
\draw  (-1.5, -3.5)-- (-3.5, -1.5);

\draw  (1.5, 1.5)-- (2,1);
\draw  (2, 1)-- (2,-1.5);

\begin{scriptsize}
\draw(0,0) node{\normalsize $\bar{D}_1$};
\draw(0,-4) node{\normalsize $\bar{D}_2$};
\draw(4,0) node{\normalsize $\bar{D}_3$};
\draw(4,4) node{\normalsize $\bar{D}_4$};
\draw(0,4) node{\normalsize $\bar{D}_5$};
\draw(-4,0) node{\normalsize $\bar{D}_6$};
\draw(-4,-4) node{\normalsize $\bar{D}_7$};
\draw(-2, -2) node{\normalsize $\tilde{D}$};
\draw [color=red, line width=1pt] (2,0.5)-- ++(-12pt,0pt);
\draw(1.4, 0.55) node{\tiny $E_1^{13}$};
\draw [color=red, line width=1pt] (2,0)-- ++(-12pt,0pt);
\draw(1.4, 0.05) node{\tiny $E_2^{13}$};
\draw [color=red, line width=1pt] (2,-0.5)-- ++(-12pt,0pt);
\draw(1.4, -0.45) node{\tiny $E_3^{13}$};

\draw [color=red, line width=1pt] (0.75,1.5)-- (0.75,3.5);
\draw(1.1, 2.5) node{\tiny $F_1^{15}$};
\draw(0, 2.5) node{\tiny $D_1^{15}$};

\draw(0.4, 0.8) node{\tiny $E_2^{15}$};
\draw [color=red, line width=1pt] (0.33,1.5)-- ++(0pt,-12pt);
\draw(-0.26, 0.8) node{\tiny $E_1^{15}$};
\draw [color=red, line width=1pt] (-0.33,1.5)-- ++(0pt,-12pt);

\draw [color=red, line width=1pt] (-3.5,0.5)-- (-1.5,0.5);
\draw(-2.5, 0.8) node{\tiny $F_1^{16}$};
\draw(-2.5, 0) node{\tiny $D_1^{16}$};

\draw [color=blue, line width = 1pt] (0.5, -1.5)-- ++(-2.25pt,-2.25pt) -- ++(4.5pt,4.5pt) ++(-4.5pt,0) -- ++(4.5pt,-4.5pt);

\draw [color=red, line width=1pt] (-0.5,-1.5)-- (-0.5,-2.5);
\draw(-0.92, -2) node{\tiny $F_2^{12}$};
\draw(0.25, -2) node{\tiny $\bar{D}_2^{12}$};

\draw [color=red, line width=1pt] (0,-2.5)-- (0,-3.5);
\draw(-0.42, -3) node{\tiny $F_1^{12}$};
\draw(0.75, -3) node{\tiny $D_1^{12}$};

\draw [color=blue, line width = 1pt] (1.5, 3.9)-- ++(-2.25pt,-2.25pt) -- ++(4.5pt,4.5pt) ++(-4.5pt,0) -- ++(4.5pt,-4.5pt);

\draw [color=blue, line width = 1pt] (-3.9, -1.5)-- ++(-2.25pt,-2.25pt) -- ++(4.5pt,4.5pt) ++(-4.5pt,0) -- ++(4.5pt,-4.5pt);

\end{scriptsize}
\end{scope}

\begin{scope}[xshift=8cm, yshift=0.5cm, scale = 1.2]
\clip(-2.5,-3.5) rectangle (2,2.5);
\draw (0,2)-- (1,2);
\draw (1,2)-- (1,0);
\draw (1,0)-- (0,-3);
\draw (0,-3)-- (-2,-3);
\draw (-2,0)-- (-2,-3);
\draw (-2,0)-- (0,2);
\draw (0,2)-- (0,0);
\draw (0,0)-- (1,0);
\draw (0,0)-- (0,-3);
\draw (0,0)-- (-2,0);
\draw (1,2)-- (-1,0);
\draw (-1,0)-- (-1,-2);
\draw (-1,-2)-- (-2,-3);
\draw (0,-2)-- (-1,-2);
\draw (0,-2)-- (1,0);
\draw (-1,-2)-- (1,0);
\draw (0,0)-- (1,2);
\begin{scriptsize}
\draw [fill] (0,0) circle (1pt);
\draw [fill] (0,1) circle (1pt);
\draw [fill] (1,0) circle (1pt);
\draw [fill] (1,1) circle (1pt);
\draw [fill] (0,2) circle (1pt);
\draw [fill] (1,2) circle (1pt);
\draw [fill] (-1,0) circle (1pt);
\draw [fill] (-2,0) circle (1pt);
\draw [fill] (-1,1) circle (1pt);
\draw [fill] (0,-1) circle (1pt);
\draw [fill] (0,-2) circle (1pt);
\draw [fill] (0,-3) circle (1pt);
\draw [fill] (-1,-1) circle (1pt);

\draw [fill] (-2,-1) circle (1pt);
\draw [fill] (-2,-2) circle (1pt);
\draw [fill] (-2,-3) circle (1pt);
\draw [fill] (-1,-2) circle (1pt);
\draw [fill] (-1,-3) circle (1pt);

\draw [color=blue, line width = 1pt] (0,0)-- ++(-2.25pt,-2.25pt) -- ++(4.5pt,4.5pt) ++(-4.5pt,0) -- ++(4.5pt,-4.5pt);
\draw [color=blue, line width = 1pt] (0,-0.5)-- ++(-2.25pt,-2.25pt) -- ++(4.5pt,4.5pt) ++(-4.5pt,0) -- ++(4.5pt,-4.5pt);

\draw [color=blue, line width = 1pt] (0.5,2)-- ++(-2.25pt,-2.25pt) -- ++(4.5pt,4.5pt) ++(-4.5pt,0) -- ++(4.5pt,-4.5pt);
\draw [color=blue, line width = 1pt] (-2, -1.5)-- ++(-2.25pt,-2.25pt) -- ++(4.5pt,4.5pt) ++(-4.5pt,0) -- ++(4.5pt,-4.5pt);

\draw(-0.25, 0.2) node{\normalsize $\bar{v}_1$};
\draw(0, -3.2) node{\normalsize $\bar{v}_2$};
\draw(1.2, 0.2) node{\normalsize $\bar{v}_3$};
\draw(1.2, 2.15) node{\normalsize $\bar{v}_4$};
\draw(0.2, 2.15) node{\normalsize $\bar{v}_5$};
\draw(-2.2, 0.2) node{\normalsize $\bar{v}_6$};
\draw(-2.2, -3.15) node{\normalsize $\bar{v}_7$};

\draw(-1.2, 0.2) node{\tiny $v_1^{16}$};
\draw(0.2, 0.85) node{\tiny $v_1^{15}$};
\draw(-0.2, -1.85) node{\tiny $v_1^{12}$};
\draw(-0.2, -0.85) node{\tiny $\bar{v}_2^{12}$};

\end{scriptsize}
\end{scope}

\end{tikzpicture}

\end{center}
\caption{Blowing up $\bar{D}_1$ again and the corresponding transformation of the dual intersection complex.} \label{blowingadivisor1}
\end{figure}
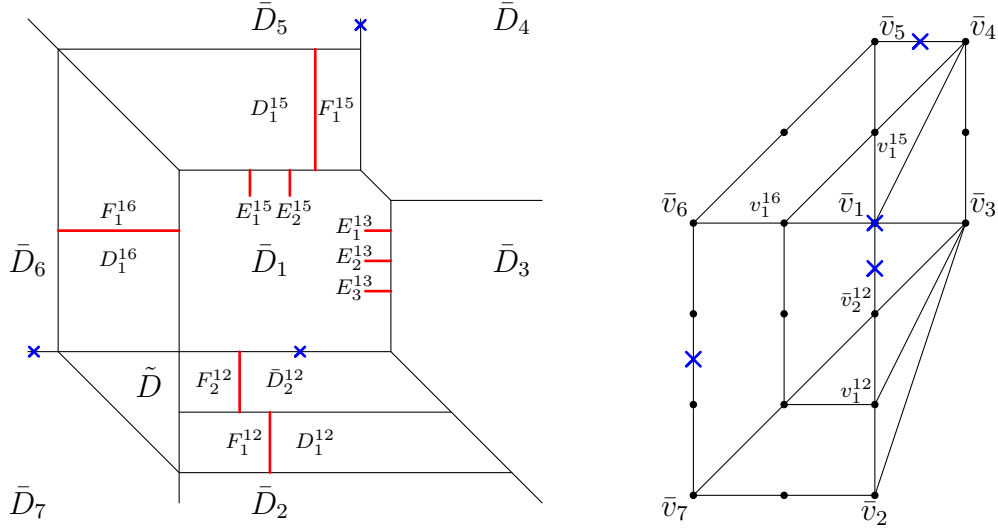

To give an example of these further blowups and illustrate the notation, we explain how to resolve the singularities on $\bar{D}_1^{15} \cap \bar{D}_1$ and $\bar{D}_1^{12} \cap \bar{D}_1$ in the resolution obtained after blowing up $\bar{D}_1$ in Figures \ref{blowingadivisor} and \ref{blowingadivisoraff}. Blow up $\bar{D}_1$ again, this resolves the ODP singularities of $\bar{D}_1^{15} \cap \bar{D}_1$ (since $l_{\rho_{15}} = 2$) producing exceptional curves $E_{1}^{15}, E_{2}^{15}$ and improves the singularity on $\bar{D}_1^{12} \cap \bar{D}_1$ to an ODP, introducing a $\mathbb{P}^1$-bundle $\bar{D}_2^{12}$ with exceptional curve $F_2^{12}$. It also separates $\bar{D}_1^{15}$ and $\bar{D}_3$, introducing another exceptional curve. The singularity of the dual intersection complex contained in the interior of  $\langle \bar{v}_1^{15}, \bar{v}_1 \rangle$ moves to $\bar{v}_1$ and the singularity contained in the interior of $\langle \bar{v}_1^{12}, \bar{v}_1 \rangle$ moves into the interior of $\langle \bar{v}_2^{12}, \bar{v}_1 \rangle$. We give the central fibre and the dual intersection complex of this partial resolution in Figure \ref{blowingadivisor1}. Note that we could have blown up $\bar{D}_{1}^{15}$ instead. Then the proper transforms of the singularities on $\bar{D}_1^{15} \cap \bar{D}_1$ would be contained in $D_1^{15}$ (and the singularity contained in the interior of  $\langle \bar{v}_1^{15}, \bar{v}_1 \rangle$ would move to $v_1^{15}$).

\begin{figure}[t]
\begin{center}

\begin{tikzpicture}[line cap=round,line join=round,>=triangle 45,x=1cm,y=1cm]

\begin{scope}[scale = 0.8]
\clip(-4.8,-4.5) rectangle (4.5,4.5);
\draw  (-2,1.5)-- (1.5,1.5);
\draw  (-1.5,-1.5)-- (2,-1.5);
\draw  (-2,1.5)-- (-4.5,4);
\draw  (1.5,1.5)-- (1.5,4);
\draw  (2,1)-- (4.5,1);
\draw  (2,-1.5)-- (4.5,-4);
\draw  (-1.5,-1.5)-- (-1.5,-4);
\draw  (-2,-1)-- (-4.5,-1);

\draw  (-4, 3.5)-- (1.5, 3.5);
\draw  (-4, 3.5)-- (-4, -1);
\draw  (-1.5, -3.5)-- (4, -3.5);
\draw  (-1.5, -2.5)-- (3, -2.5);
\draw  (-1.5, -3.5)-- (-4, -1);

\draw  (1.5, 1.5)-- (2,1);
\draw  (2, 1)-- (2,-1.5);

\draw  (-1.5, -1.5)-- (-2,-1);
\draw  (-2, -1)-- (-2,1.5);

\begin{scriptsize}
\draw(0,0) node{\normalsize $D_1$};
\draw(0,-4) node{\normalsize $\bar{D}_2$};
\draw(4,0) node{\normalsize $\bar{D}_3$};
\draw(4,4) node{\normalsize $\bar{D}_4$};
\draw(0,4) node{\normalsize $\bar{D}_5$};
\draw(-4.5,0) node{\normalsize $\bar{D}_6$};
\draw(-4,-4) node{\normalsize $\bar{D}_7$};
\draw(-2.25, -1.75) node{\normalsize $\tilde{D}$};
\draw [color=red, line width=1pt] (2,0.5)-- ++(-12pt,0pt);
\draw(1.4, 0.55) node{\tiny $E_1^{13}$};
\draw [color=red, line width=1pt] (2,0)-- ++(-12pt,0pt);
\draw(1.4, 0.05) node{\tiny $E_2^{13}$};
\draw [color=red, line width=1pt] (2,-0.5)-- ++(-12pt,0pt);
\draw(1.4, -0.45) node{\tiny $E_3^{13}$};

\draw [color=red, line width=1pt] (0.75,1.5)-- (0.75,3.5);
\draw(1.1, 2.5) node{\tiny $F_1^{15}$};
\draw(0, 2.5) node{\tiny $D_1^{15}$};

\draw(0.4, 0.8) node{\tiny $E_2^{15}$};
\draw [color=red, line width=1pt] (0.33,1.5)-- ++(0pt,-12pt);
\draw(-0.26, 0.8) node{\tiny $E_1^{15}$};
\draw [color=red, line width=1pt] (-0.33,1.5)-- ++(0pt,-12pt);

\draw [color=red, line width=1pt] (-4,0.5)-- (-2,0.5);
\draw(-3, 0.8) node{\tiny $F_1^{16}$};
\draw(-3, 0) node{\tiny $D_1^{16}$};

\draw(0.6, -0.9) node{\tiny $E_1^{12}$};
\draw [color=red, line width=1pt] (0.5,-1.5)-- ++(0pt,12pt);

\draw [color=red, line width=1pt] (-0.5,-1.5)-- (-0.5,-2.5);
\draw(-0.92, -2) node{\tiny $F_2^{12}$};
\draw(0.25, -2) node{\tiny $D_2^{12}$};

\draw [color=red, line width=1pt] (0,-2.5)-- (0,-3.5);
\draw(-0.42, -3) node{\tiny $F_1^{12}$};
\draw(0.75, -3) node{\tiny $D_1^{12}$};

\draw [color=blue, line width = 1pt] (1.5, 3.9)-- ++(-2.25pt,-2.25pt) -- ++(4.5pt,4.5pt) ++(-4.5pt,0) -- ++(4.5pt,-4.5pt);

\draw [color=blue, line width = 1pt] (-4.4, -1)-- ++(-2.25pt,-2.25pt) -- ++(4.5pt,4.5pt) ++(-4.5pt,0) -- ++(4.5pt,-4.5pt);

\end{scriptsize}
\end{scope}

\begin{scope}[xshift=8cm, yshift=0.5cm, scale = 1.2]
\clip(-2.5,-3.5) rectangle (2,2.5);
\draw (0,2)-- (1,2);
\draw (1,2)-- (1,0);
\draw (1,0)-- (0,-3);
\draw (0,-3)-- (-2,-3);
\draw (-2,0)-- (-2,-3);
\draw (-2,0)-- (0,2);
\draw (0,2)-- (0,0);
\draw (0,0)-- (1,0);
\draw (0,0)-- (0,-3);
\draw (0,0)-- (-2,0);
\draw (1,2)-- (-1,0);
\draw (-1,0)-- (-1,-2);
\draw (-1,-2)-- (-2,-3);
\draw (0,-2)-- (-1,-2);
\draw (0,-2)-- (1,0);
\draw (-1,-2)-- (1,0);
\draw (-1,-2)-- (1,2);
\begin{scriptsize}
\draw [fill] (0,0) circle (1pt);
\draw [fill] (0,1) circle (1pt);
\draw [fill] (1,0) circle (1pt);
\draw [fill] (1,1) circle (1pt);
\draw [fill] (0,2) circle (1pt);
\draw [fill] (1,2) circle (1pt);
\draw [fill] (-1,0) circle (1pt);
\draw [fill] (-2,0) circle (1pt);
\draw [fill] (-1,1) circle (1pt);
\draw [fill] (0,-1) circle (1pt);
\draw [fill] (0,-2) circle (1pt);
\draw [fill] (0,-3) circle (1pt);
\draw [fill] (-1,-1) circle (1pt);

\draw [fill] (-2,-1) circle (1pt);
\draw [fill] (-2,-2) circle (1pt);
\draw [fill] (-2,-3) circle (1pt);
\draw [fill] (-1,-2) circle (1pt);
\draw [fill] (-1,-3) circle (1pt);

\draw [color=blue, line width = 1pt] (0,0)-- ++(-2.25pt,-2.25pt) -- ++(4.5pt,4.5pt) ++(-4.5pt,0) -- ++(4.5pt,-4.5pt);

\draw [color=blue, line width = 1pt] (0.5,2)-- ++(-2.25pt,-2.25pt) -- ++(4.5pt,4.5pt) ++(-4.5pt,0) -- ++(4.5pt,-4.5pt);
\draw [color=blue, line width = 1pt] (-2, -1.5)-- ++(-2.25pt,-2.25pt) -- ++(4.5pt,4.5pt) ++(-4.5pt,0) -- ++(4.5pt,-4.5pt);

\draw(-0.25, 0.2) node{\normalsize $v_1$};
\draw(0, -3.2) node{\normalsize $\bar{v}_2$};
\draw(1.2, 0.2) node{\normalsize $\bar{v}_3$};
\draw(1.2, 2.15) node{\normalsize $\bar{v}_4$};
\draw(0.2, 2.15) node{\normalsize $\bar{v}_5$};
\draw(-2.2, 0.2) node{\normalsize $\bar{v}_6$};
\draw(-2.2, -3.15) node{\normalsize $\bar{v}_7$};

\draw(-1.2, 0.2) node{\tiny $v_1^{16}$};
\draw(0.2, 0.85) node{\tiny $v_1^{15}$};
\draw(-0.2, -1.85) node{\tiny $v_1^{12}$};
\draw(-0.2, -0.85) node{\tiny $v_2^{12}$};

\end{scriptsize}
\end{scope}

\end{tikzpicture}

\end{center}
\caption{Blowing up $\bar{D}_1$ (once again) and the corresponding transformation of the dual intersection complex.} \label{blowingadivisor2}
\end{figure}
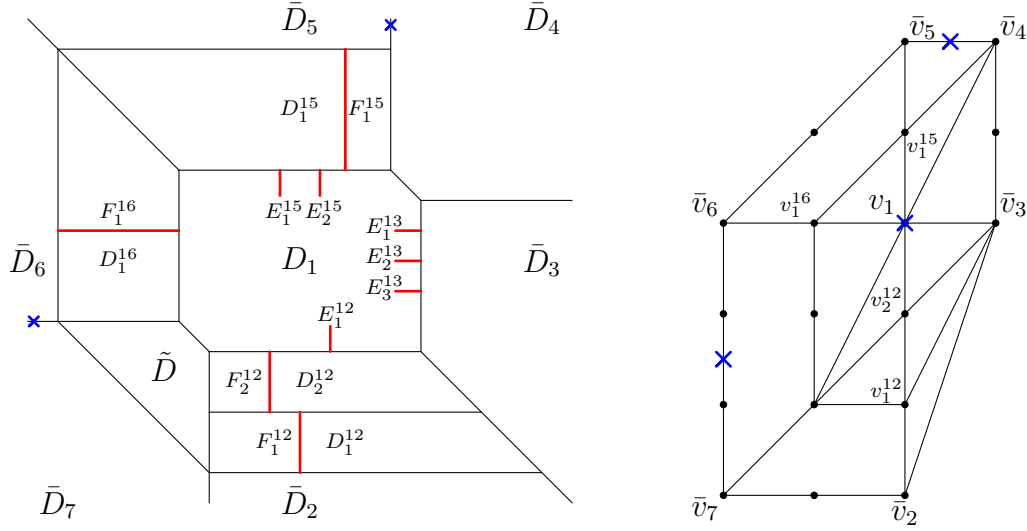

Now, to resolve the singularity on $\bar{D}_2^{12} \cap \bar{D}_1$, we can blow up either $\bar{D}_2^{12}$ or $\bar{D}_1$. Either way, this resolves the ODP singularity on $\bar{D}_2^{12} \cap \bar{D}_1$, producing an exceptional curve $E_1^{12}$. In the first case, the curve is contained in $D_2^{12}$ and intersects a curve of class $F_{2}^{12}-E_1^{12}$ at one point. In the second case, the curve is contained in $D_1$. The blowup also introduces an exceptional curve separating $\tilde{D}$ and $D_1$ in the first case and separating $D_{1}^{16}$ and $D_{2}^{12}$ in the second case. The singularity of the dual intersection complex contained in the interior of $\langle \bar{v}_2^{12}, \bar{v}_1 \rangle$ moves to $v_2^{12}$ in the first case and to $v_1$ in the second case. We give the central fibre and the dual intersection complex of the resolution obtained by blowing up $\bar{D}_1$ (once again) in Figure \ref{blowingadivisor2}.

After these blowups, the resolution becomes log smooth in a neighbourhood of the proper transform of $\bar{D}_1$.

\subsubsection{A global resolution} \label{global42}

By continuously performing the blowups of Sections \ref{divisorblowupgeneral} and \ref{furthercompblowups} in some order, we resolve all the singularities of $\bar{\mathfrak{X}} \to \mathcal{S}$ and obtain a resolution $\pi: \mathfrak{X} \to \bar{\mathfrak{X}}$ to a log smooth degeneration $\mathfrak{X} \to \mathcal{S}$. Here the log structure on $\mathfrak{X}$ is the divisorial log structure given by $D = D_1 + \dots + D_m$ where $D_i, \ 1 \leq i \leq m$ are the irreducible components of the central fibre $\mathfrak{X}_0$ and $D_i$ is the strict transform of $\bar{D}_i$ for $1 \leq i \leq \bar{m}$. Note that this construction agrees with the construction of Section \ref{resolvesmall} if all the $\sigma \in \bar{\mathscr{P}}^{\max}$ are standard triangles. 

The dual intersection complex $(B, \mathscr{P})$ of $\mathfrak{X} \to \mathcal{S}$ is a natural subdivision of $\left(\bar{B}, \bar{\mathscr{P}}\right)$ and has an affine structure (with singularities at the vertices) given by Construction \ref{affstrintgen}. One may view the resolution tropically as a composition of transformations as in Figure \ref{blowingadivisoraff}. Note that the affine structure extends across the toric components, so we only have singularities at vertices of the form $v_i$ for $ 1 \leq i \leq \bar{m}$ and $v_{\rho, p}$ for $\rho \in \bar{\mathscr{P}}^{[1]}, \ 1 \leq p \leq l_{\rho} -1$. Moreover, for any $\rho \in \bar{\mathscr{P}}^{[1]}$ there is at most one singularity at a vertex of the form $v_{\rho, p}, \ 1 \leq p \leq l_{\rho} -1$.

\begin{remarks} \label{restrictedclassofresolutions}
\begin{enumerate}
    \item One can always resolve the singularities of $\bar{\mathfrak{X}} \to \mathcal{S}$ by blowing up the original components $\bar{D}_i, \ 1 \leq i \leq \bar{m}$ and their strict transforms. Indeed, we can blow up $\bar{D}_1$ the number of times equal to the integral length of the longest edge adjacent to $\bar{v}_1 \in \bar{\mathscr{P}}^{[0]}$ and then do the same for all the other divisors. In this case, all the $E_{\rho, p}$ (for $\rho \in \bar{\mathscr{P}}^{[1]}$ and $1 \leq p \leq r_{\rho}$) are contained in $D_i, \ 1 \leq i \leq \bar{m}$ and the singularities are contained in $v_i, \ 1 \leq i \leq \bar{m}$.
    \item One can also always construct a resolution with $D$ simple normal crossings. To ensure this, after obtaining a log smooth resolution, one can keep blowing up components (corresponding to vertices contained in maximal cells that are not standard triangles) until all the cells $\sigma \in \mathscr{P}^{\max}$ are standard triangles.
    \item Resolutions satisfying (1) and (2) would suffice for our purposes, but we choose to allow some additional flexibility. The reader may prefer to think of the resolutions satisfying these additional requirements.
\end{enumerate}
\end{remarks}

The exceptional locus of $\pi: \mathfrak{X} \to \bar{\mathfrak{X}}$ is a union of exceptional curves and exceptional divisors. Note that the cone $NE(\pi) \subseteq NE(\mathfrak{X}_0)$ of curves contracted by $\pi$ is finitely generated and that we have
\begin{equation} 
\left\{ E_{\rho,k}  \ | \ \rho \in \bar{\mathscr{P}}^{[1]}, \ 1 \leq k \leq r_{\rho}\right\} \cup \left\{F_{\rho, p}, \ \rho \in \bar{\mathscr{P}}^{[1]}, \ 1 \leq p \leq l_{\rho} -1  \right\} \subseteq NE(\pi).
\end{equation}
(with possibly not all $F_{\rho,p}$ present if $r_{\rho} = 0$).
We will mostly be interested in these curve classes.

Unlike Section \ref{resolvesmall}, the divisor $D$ is not simple normal crossings unless all the $\sigma \in \mathscr{P}^{\max}$ are standard triangles. However, we can generalize the rest of Proposition \ref{sncmlcy}.

\begin{fact} \label{minlogcy42}
$\mathfrak{X} \to \mathcal{S}$ is minimal log CY. The dual intersection complex $(B, \mathscr{P})$ satisfies Assumption \ref{manifoldass}.
\end{fact}
	
\begin{proof}
The fact that $\mathfrak{X} \to \mathcal{S}$ is minimal log CY follows from the behaviour of the canonical class under blowup. Indeed, for any degeneration $\bar{\mathfrak{X}} \to \mathcal{S}$ with $K_{\bar{\mathfrak{X}}} + \bar{\mathfrak{X}}_0 \equiv 0$ and any blowup $\pi: \mathfrak{X} \to \bar{\mathfrak{X}}$ supported on $\bar{\mathfrak{X}}_0$, we have $K_{\mathfrak{X}} + (\mathfrak{X}_0)_{\reduced} \equiv 0$. 
Now, $D = D_1 + \dots + D_m$ where $D_i, \ 1 \leq i \leq m$ are the irreducible components of $\mathfrak{X}_0$, $\bar{D} = \bar{D}_1 + \dots + \bar{D}_{\bar{m}}$ where $\bar{D}_i, \ 1 \leq i \leq \bar{m}$ are the irreducible components of $\bar{\mathfrak{X}}_0$, and $K_{\bar{\mathfrak{X}}} + \bar{D} \equiv 0$ since $\bar{\mathfrak{X}} \to \mathcal{S}$ is a toric degeneration. So it suffices to check that $\mathfrak{X}_0$ is reduced. But this follows from the fact that $(B, \mathscr{P})$ is an integral subdivision of $\left( \bar{B}, \bar{\mathscr{P}} \right)$ by construction (see Remark \ref{nonreducedeness}).  So $\mathfrak{X} \to \mathcal{S}$ is minimal log CY.

The fact that $(B, \mathscr{P})$ satisfies Assumption \ref{manifoldass} follows from $\left(\bar B, \bar{\mathscr{P}}\right)$ satisfying Assumption \ref{manifoldass} and the description of $(B, \mathscr{P})$ as a subdivision of $\left(\bar B, \bar{\mathscr{P}}\right)$ (here it is crucial that the cells $\sigma \in \bar{\mathscr{P}}^{\max}$ are subdivided according to Figure \ref{admisspolytopesresolve}).
\end{proof}

As in Section \ref{globalresolvesmall}, let $P$ be a well-chosen monoid (see Definition \ref{wellchosendef}) with a face $K$ containing the classes of the contracted curves.

\begin{notation} \label{minnotation}
For a PA-function $\alpha$ on a polyhedral manifold $B$ (possibly with boundary), we denote by $\min(\alpha)$ the subset of $B$ where the minimum of $\alpha$ is achieved. We say that $\min(\alpha)$ is \emph{well-defined} if it consists of one point. 
\end{notation}

Proposition \ref{ampledivisorsmall} can be partially generalized as follows (again, we can't guarantee that the $\pi$-ample divisor $D'$ is simple normal crossings in general). 

\begin{fact} \label{ampledivisor}
The resolution $\pi: \mathfrak{X} \to \bar{\mathfrak{X}}$ satisfies the following assumptions:
\begin{enumerate}
    \item The log structure $\mathcal{M}_{\mathfrak{X}}$ on $\mathfrak{X}$ is fine, saturated, and Zariski.
    \item There exists a $\pi$-ample effective divisor $D' = \sum_{i=1}^{m} a_i D_i$ such that:
    \begin{enumerate}
        \item $D'$ is PA-generated i.e. corresponds to an integral PA-function on $B$ (see Definition \ref{pagenerated}).
        \item Let  $D_{\irrel} := \sum_{i=1}^m \varepsilon_i D_i$ where $\varepsilon_i = 0$ if $a_i > 0$ and $\varepsilon_i = 1$ if $a_i = 0$. Then $D_{\irrel}$ is $\mathbb{Q}$-Cartier and $\pi$-nef.
    \end{enumerate}
    \item $K \cap NE(\mathfrak{X}_0) = K \cap NE(\mathfrak{X}_0)_{\numer}$ (under the splittings of \eqref{standardsplittings}), i.e. $\pi$ only contracts numerical classes of curves $C \subseteq \mathfrak{X}_0$.
\end{enumerate}
In particular, $\pi: \mathfrak{X} \to \bar{\mathfrak{X}}$ satisfies the assumptions of Proposition \ref{P2general} if $D'$ is simple normal crossings.
\end{fact}

\begin{proof}
\textbf{Step 1.} First of all, $(\bar{\mathfrak{X}}, \bar{D})$ is a Zariski log scheme by Assumption \ref{toricdegenerationass}(3), so $(\mathfrak{X}, D)$ is also a Zariski log scheme. The log structure $\mathcal{M}_{\mathfrak{X}}$ on $\mathfrak{X}$ is fine and saturated away from the proper transform of the singular locus $Z \subseteq \bar{\mathfrak{X}}_0$ since the log structure $\mathcal{M}_{\bar{\mathfrak{X}}}$ on $\bar{\mathfrak{X}}$ is fine and saturated away from $Z$. It follows from the local models of Sections \ref{toricmodelcodim2} and \ref{localmodelgen} that it is also fine and saturated in the neighbourhood of any point of the proper transform of $Z$. So $\mathcal{M}_{\mathfrak{X}}$ is fine and saturated, and $\pi: \mathfrak{X} \to \bar{\mathfrak{X}}$ satisfies assumption (1). It is also clear from the analysis of the local models that the cone $NE(\pi) \subseteq NE(\mathfrak{X}_0)$ of curves contracted by $\pi$ is finitely generated and every curve $C \in NE(\mathfrak{X}_0)$ defines a non-zero class in $A_1(\mathfrak{X}_0, \mathbb{Z})_{\numer}$. So $\pi: \mathfrak{X} \to \bar{\mathfrak{X}}$ satisfies assumption (3). 

\textbf{Step 2.} It suffices to show that there exists a $\pi$-ample PA-generated divisor $D'$ on $\mathfrak{X}$. Indeed, if we find such a $D'$, then $D' + ND$ is PA-generated, $\pi$-ample, and effective for $N \gg 0$ (since $D \equiv 0$, $D$ is PA-generated by Corollary \ref{DisPLgen} and PA-generated divisors are supported on the central fibre). Moreover, we would have $(D' + ND)_{\irrel} = 0$ for $N \gg 0$, so $\pi: \mathfrak{X} \to \bar{\mathfrak{X}}$ would satisfy assumption (2). 

\textbf{Step 3.} Finding a $\pi$-ample PA-generated divisor $D'$ is equivalent to finding a PA-function $\alpha \in PA(B)$ such that:
\begin{enumerate}
    \item The restriction $\alpha|_{\sigma}$ of $\alpha$ to any $\sigma \in \bar{\mathscr{P}}^{\max}$ is strictly convex (on the induced subdivision of $\sigma$).
    \item The restriction $\alpha|_{\rho}$ of $\alpha$ to any $\rho \in \bar{\mathscr{P}}^{[1]}$ is strictly convex with $\min(\alpha|_{\rho}) = \{v_{\rho, p} \}$ (using Notation \ref{minnotation}) if the exceptional curves $E_{\rho,k}, \ 1 \leq k \leq r_{\rho}$ are contained in $D_{\rho, p}$ (as before, if $p = 0$ or $p = l_{\rho}$, we understand $D_{\rho, p}$ and $v_{\rho, p}$ as the strict transform of an irreducible component of $\bar{\mathfrak{X}}_0$ and the corresponding vertex respectively). There is no condition if $r_{\rho} = 0$.
\end{enumerate}
Indeed, let $D': = \divisor(\alpha)$ via \eqref{divmap}. Condition (1) is equivalent to the toric divisors corresponding to $\alpha|_{\sigma}$ in the resolved local models for $\bar{X}_{\sigma}, \sigma \in \bar{\mathscr{P}}^{\max}$ being relatively ample. Condition (2) is equivalent to requiring that the PA-functions $\alpha|^0_{\rho, x}$ that are the extensions of $\alpha|_{\rho}$ by $0$ to the subdivision of the triangle or trapezoid of Figure \ref{model2} (or Figure \ref{model1}) giving rise to the resolved local models at the points $x \in \bar{X}_{\rho}, \rho \in \bar{\mathscr{P}}^{[1]}$ are strictly convex. Here by an extension by $0$, we mean that $\alpha|_{\rho, x}^0$ is the PA-function that takes the same values as $\alpha$ at the points corresponding to the integer points of $\rho$ and takes value $0$ at $(0,1,1)$ and $(1,1,1)$ (or just $(0,1,1)$). Now condition (2) is equivalent to requiring that the toric divisors corresponding to $\alpha|^0_{\rho, x}$ in these local models are relatively ample. The claim follows from the fact that $NE(\pi) = \overline{NE(\pi)}$ (since $NE(\pi)$ is finitely generated), from the (relative) Kleiman's criterion for ampleness, and from the local nature of intersection numbers (see, e.g. \cite[Appendix A.1]{hartshorne}).
 
Note that condition (2) implies that $\alpha \in PA(B)$ determines the singular locus $\Delta$ of the affine structure on $(B, \mathscr{P})$. Indeed, since for every $\rho \in \bar{\mathscr{P}}^{[1]}$ the exceptional curves $E_{\rho,k}, \ 1 \leq k \leq r_{\rho}$ are contained in $D_{\rho, p}$ and all the other irreducible components of $\mathfrak{X}_0$ are toric, the singular locus of the affine structure on $(B, \mathscr{P})$ is $\Delta = \{ v_{\rho,p} \ | \ \rho \in \bar{\mathscr{P}}^{[1]}, \ r_{\rho} \neq 0\}$. 

\textbf{Step 4}. To find a PA-function $\alpha \in PA(B)$ satisfying the required properties, we proceed inductively.\footnote{The motivation for the argument below is that for a morphism $f: X \to Y$ of projective varieties with $D_Y$ an ample divisor on $Y$ and $D_{X/Y}$ a relatively ample divisor on $X$, $D_{X/Y} + Mf^*D_Y$ is an ample divisor on $X$ for $M \gg 0$ (see \cite[Part II, Proposition 7.10]{hartshorne}).} For a partial resolution $\mathfrak{X}' \to \mathcal{S}$ with dual intersection complex $\left(B', \mathscr{P}'\right)$, we require a PA-function $\alpha' \in PA(B')$ that satisfies the analogues of conditions (1) and (2) in Step 3 with the minimum requirement in (2) replaced by the requirement that
\begin{enumerate}
    \item[$(*)$] $\bar{v}_{\rho,p} \in \min(\alpha|_{\rho})$ for $\bar{D}_{\bar{v}_{\rho,p}}$ the strict transform of the last divisor (with $\bar{v}_{\rho,p} \subseteq \rho$) that has been blown up.
\end{enumerate}
in the case that the singularities of $\bar{X}_{\rho}$ have not been fully resolved.

The dual intersection complex $\left( \bar{B}, \bar{\mathscr{P}} \right)$ of $\bar{\mathfrak{X}} \to \mathcal{S}$ admits a PA-function $\bar{\alpha}$ satisfying these (modified) conditions. Namely, we can just take $\bar{\alpha} := 0 \in PA(\bar{B})$. Suppose that at some point we have a partial resolution $\mathfrak{X}'_1 \to \mathcal{S}$ with dual intersection complex $\left(B'_1, \mathscr{P}'_1 \right)$ and a PA-function $\alpha'_1 \in PA\left(B'_1\right)$ satisfying the conditions. Let $\mathfrak{X}'_2 \to \mathfrak{X}'_1$ be the blowup of a component $\bar{D}'_{1,i}$ of $\left(\mathfrak{X}'_1\right)_0$ (the next blowup in the resolution process). We need to construct a PA-function $\alpha'_2 \in PA\left(B'_2\right)$ on the the dual intersection complex $\left(B'_2, \mathscr{P}'_2 \right)$ of $\mathfrak{X}'_2 \to \mathcal{S}$ satisfying the conditions. Let $\alpha_{\new}  \in PA\left(B_2'\right)$ be a PA-function with value $-1$ at $\bar{v}'_{1,i}$ (the vertex corresponding to the strict transform of $\bar{D}'_{1,i}$ under the blowup) and at all the \enquote{new} vertices and value $0$ at all the other vertices. It is easy to confirm using Figure \ref{admisspolytopesresolve} that $\alpha_{\new} \in PA\left(B_2'\right)$ is well-defined, that for every $\sigma \in \mathscr{P}'^{\max}_1$ subdivided by the blowup, $(\alpha_{\new})|_{\sigma}$ is a strictly convex function on $\sigma$, and that we have $\bar{v}'_{1,i} \in \min ((\alpha_{\new})|_{\rho})$ for every $\rho \in \mathscr{P}'^{[1]}_1$ with $\bar{v}'_{1,i} \subseteq \rho$. Now for $M \gg 0$ the PA-function $\alpha_2':= \alpha_{\new} + M \alpha_1'$ on $\left(B'_2, \mathscr{P}'_2 \right)$ satisfies the conditions. Repeating the argument for all the blowups in the construction of $\pi: \mathfrak{X} \to \bar{\mathfrak{X}}$, we obtain an $\alpha \in PA(B)$ satisfying conditions (1) and (2) of Step 3. 
\end{proof}

We shall revisit the idea of Steps 3 and 4 of the proof of Proposition \ref{ampledivisor} in Section \ref{admissible}, where we construct admissible resolutions in general.

As a consequence of Proposition \ref{ampledivisor}, we obtain resolutions satisfying the assumptions of Proposition \ref{P2general}. 

\begin{cor} \label{furtherresolvetosatisfy}
For any resolution $\pi: \mathfrak{X} \to \bar{\mathfrak{X}}$ obtained by blowing up a sequence of irreducible components of $\bar{\mathfrak{X}}_0$, there exists a further sequence of blowups $\pi': \tilde{\mathfrak{X}} \to \mathfrak{X}$ of irreducible components of $\mathfrak{X}_0$ with the combined $\tilde{\mathfrak{X}} \to \mathfrak{X} \to \bar{\mathfrak{X}}$ satisfying the assumptions of Proposition \ref{P2general}.
\end{cor}

\begin{proof}
    We let $\tilde{D}:=(\pi')^{-1}D$ and denote by $\tilde{D}'$ any divisor on $\tilde{\mathfrak{X}}$ supported on $\tilde{D}$.
    By Proposition \ref{ampledivisor}, we just need to check that there exists a further sequence of blowups $\pi':\tilde{\mathfrak{X}} \to \mathfrak{X}$ of irreducible components of $\mathfrak{X}_0$ such that the PA-generated relatively ample divisor $\tilde{D}'$ on $\tilde{\mathfrak{X}}$ given by Proposition \ref{ampledivisor} is simple normal crossings. Indeed, by Remark \ref{restrictedclassofresolutions}(2), we can always resolve to a $\tilde{\mathfrak{X}}$ with $\tilde{D}$ simple normal crossings and $\tilde{D}'$ is supported on $\tilde{D}$.
\end{proof}

\begin{remark} \label{satisfybutnotsnc}
Note that one can often produce resolutions $\mathfrak{X} \to \mathcal{S}$ satisfying the assumptions of Proposition \ref{P2general} but with $D$ not simple normal crossings, especially in the case that $r_{\rho} = 0$ for some $\rho \in \bar{\mathscr{P}}^{[1]}$.
\end{remark}

We define the initial slab functions for  $\bar{\mathfrak{X}} \to \mathcal{S}$ as in \eqref{choiceslab} but require that the slabs $\underline{\rho}, \underline{\rho}' \in \tilde{\bar{\mathscr{P}}}^{[1]}$ are chosen so that the vertex $v_{\rho, p} \in \mathscr{P}^{[0]}$ corresponding to the component $D_{\rho, p}$ containing the exceptional curves $E_{\rho,k}, \ 1 \leq k \leq r_{\rho}$ is contained in $\underline{\rho}$. We have now defined all the necessary data in Basic Setup \ref{setupforproof}. 

Sometimes the extended intrinsic mirror is algebraic. 

\begin{observe} \label{extmirroralgebraic}
Let $A$ be the polarization on $\bar{\mathfrak{X}} \to \mathcal{S}$, $\pi: \mathfrak{X} \to \bar{\mathfrak{X}}$ be a resolution as before and suppose that $\pi^{*}A$ is a PA-generated divisor (supported on the central fibre). Then for $M \gg 0$, $D'+ M \pi^{*}A$ is a PA-generated divisor on $\mathfrak{X}$, relatively ample for $\mathfrak{X} \to \mathcal{S}$ (see the footnote in the proof of Proposition \ref{ampledivisor} above). The extended intrinsic mirror is defined as an algebraic family $\check{\mathfrak{X}} \to \Spec \kk [P]$ by arguing as in Proposition \ref{P1general} for the contraction $\mathfrak{X} \to \mathcal{S}$.
\end{observe}

\subsubsection{The PL-isomorphism \texorpdfstring{$\Phi: B \to \bar{B}$}{Phi}} 

We want to generalize Construction \ref{plisosmall} of the PL-isomorphism $\Phi: (B, \mathscr{P}) \to \left(\bar{B}, \bar{\mathscr{P}}\right)$ to the setup of this section. The difference with the situation of Section \ref{resolvesmall} is that we no longer have $(B, \mathscr{P}) \cong \left(\bar B, \bar{\mathscr{P}}\right)$ as polyhedral manifolds.

\begin{cons} \label{plmap}
We define a PL-isomorphism $\Phi: (B, \mathscr{P}) \to \left(\bar{B}, \bar{\mathscr{P}}\right)$, linear on the maximal cells of $\mathscr{P}$. Note that $(B, \mathscr{P})$ is a natural subdivision of $\left(\bar{B}, \bar{\mathscr{P}}\right)$. Moreover, the affine structure on $(B, \mathscr{P})$ extends across the toric components so we only have singularities at vertices of the form $v_i$ for $ 1 \leq i \leq \bar{m}$ and at most one vertex of the form $v_{\rho, p}, \ 1 \leq p \leq l_{\rho} -1$ for $\rho \in \bar{\mathscr{P}}^{[1]}$. Let $\mathscr{P}_{\coar}$ be the obvious coarsened polyhedral decomposition on $B$ such that $\left(B, \mathscr{P}_{\coar}\right) \cong \left(\bar B, \bar{\mathscr{P}}\right)$ as polyhedral manifolds. 

Now, $\left(B, \mathscr{P}_{\coar}\right)$ is an affine manifold with singularities that satisfies Construction \ref{affchartstool}. Indeed, it has the obvious structures of integral polyhedra on $\left\{\sigma \ | \ \sigma \in \mathscr{P}^{\max}_{\coar} \right\}$ and the structures of integral affine manifolds on $$\left\{W_v \setminus \Delta \ | \ v \in \mathscr{P}^{[0]}_{\coar} \right\}$$ coming from the extended affine structure on $(B, \mathscr{P})$. Note that $\left(B, \mathscr{P}_{\coar}\right)$ may have singularities at both the vertices $v \in \mathscr{P}_{\coar}^{[0]}$ and at the interiors of the edges $\rho \in \mathscr{P}_{\coar}^{[1]}$ (unless we consider a resolution $\pi: \mathfrak{X} \to \bar{\mathfrak{X}}$ satisfying Remark \ref{restrictedclassofresolutions}(2)). We have $(B, \mathscr{P}) \cong \left(B, \mathscr{P}_{\coar}\right)$ as affine manifolds with singularities, so it is enough to define a PL-isomorphism $\Phi: \left(B, \mathscr{P}_{\coar}\right) \to \left(\bar{B}, \bar{\mathscr{P}}\right)$, linear on the maximal cells of $\mathscr{P}_{\coar}$. 

We define the PL-isomorphism $\Phi: \left(B, \mathscr{P}_{\coar}\right) \to \left(\bar{B}, \bar{\mathscr{P}}\right)$ as in Construction \ref{plisosmall}, using the fact that for all the maximal cells $\sigma \in \mathscr{P}^{\max}_{\coar}$ there are canonical identity maps $\Id: \sigma \to \bar{\sigma}$ compatible with each other (here for $\tau \in \mathscr{P}_{\coar}$, we let $\bar{\tau} \in \bar{\mathscr{P}}$ be the corresponding cell of $\bar{\mathscr{P}}$) and the refined description of Construction \ref{affchartstool}$(2')$ for the affine structures.

We define $\Phi(\tau)$ for $\tau \in \mathscr{P}$ or $\tau \in \mathscr{P}_{\coar}$ (in which case $\Phi(\tau) \in \bar{\mathscr{P}}$), $\Phi(\mathfrak{p})$ for a codimension $0$ wall $\mathfrak{p}$ on $(B, \mathscr{P})$ and $\Phi(\mathfrak{b})$ for a slab $\mathfrak{b}$ on $(B, \mathscr{P})$ as in Construction \ref{plisosmall}. 
Note that the image of a codimension $0$ (resp. codimension $1$) rational polyhedral subset of $\sigma \in \mathscr{P}^{\max}$ is a codimension $0$ (resp. codimension $1$) rational polyhedral subset of the unique cell containing $\Phi(\sigma)$.
\end{cons}

\subsection{Tropical approach and admissible resolutions in general} \label{admissible}

We continue with the general setup of Section \ref{resolvesmall}. We remove the restrictions of Sections \ref{resolvesmall} and \ref{resolve} on the types of $\sigma \in \bar{\mathscr{P}}^{\max}$ and will treat the general case. As we have seen in Remark \ref{moregeneralsubdivisions}, a blowup of an irreducible component of the central fibre $\bar{\mathfrak{X}}_0$ of $\bar{\mathfrak{X}} \to \mathcal{S}$ (a toric degeneration of K3-s that is a divisorial log deformation with a smooth generic fibre) can give rise to a subdivision of $\left(\bar{B}, \bar{\mathscr{P}}\right)$ that does not satisfy Assumption \ref{manifoldass}. Therefore, to construct a resolution $\pi: \mathfrak{X} \to \bar{\mathfrak{X}}$ to a minimal log CY degeneration $\mathfrak{X} \to \mathcal{S}$ in the general case, we need to consider more general blowups to have more control over the induced subdivisions of $\left(\bar{B}, \bar{\mathscr{P}}\right)$.

\subsubsection{Tropical approach to resolutions} \label{tropicalresolutions}

As in Sections \ref{resolvesmall} and \ref{resolve}, we have étale local models for points in codimension $1$ and $2$ strata of $\bar{\mathfrak{X}}_0$ that are given by toric varieties. 

\begin{notations}
    \begin{enumerate}
        \item For $\sigma \in \bar{\mathscr{P}}^{\max}$, we denote the local toric model at the point $\bar{X}_{\sigma}$ (given by the cone over $\sigma$) by $\tilde{\bar{X}}_{\sigma}$. 
        \item For $\rho \in \bar{\mathscr{P}}^{[1]}$, we denote the local toric model at the point $x \in \bar{X}_{\rho}$ (given by the cone over the triangle as in Figure \ref{model1} if $x$ is a non-singular point or the cone over the trapezoid as in Figure \ref{model2} if $x$ is a singular point) by $\tilde{\bar{X}}_{\rho, x}$.
    \end{enumerate}    
\end{notations}

Suppose that we have obtained a resolution $\pi: \mathfrak{X} \to \bar{\mathfrak{X}}$ to a log smooth $\mathfrak{X} \to \mathcal{S}$. We want to understand $\pi: \mathfrak{X} \to \bar{\mathfrak{X}}$ in terms of log smooth resolutions of the étale local models $\tilde{\bar{X}}_{\sigma}$ and $\tilde{\bar{X}}_{\rho, x}$.  

For every $\sigma \in \bar{\mathscr{P}}^{\max}$, since $\tilde{\bar{X}}_{\sigma}$ is an étale local model for $\bar{\mathfrak{X}}$, we have a variety $\bar{\mathfrak{U}}_{\sigma}$ equipped with étale maps $\bar{\mathfrak{U}}_{\sigma} \to \bar{\mathfrak{X}}$ and $\bar{\mathfrak{U}}_{\sigma} \to \tilde{\bar{X}}_{\sigma}$. We let $\mathfrak{U}_{\sigma}:= \bar{\mathfrak{U}}_{\sigma} \times_{\bar{\mathfrak{X}}} \mathfrak{X}$ be the basechange of $\bar{\mathfrak{U}}_{\sigma} \to \bar{\mathfrak{X}}$ by $\pi: \mathfrak{X} \to \bar{\mathfrak{X}}$. For any morphism $\pi_{\sigma}: \tilde{X}_{\sigma} \to \tilde{\bar{X}}_{\sigma}$, we let $\tilde{\mathfrak{U}}_{\sigma, \pi_{\sigma}}:= \bar{\mathfrak{U}}_{\sigma} \times_{\tilde{\bar{X}}_{\sigma}} \tilde{X}_{\sigma}$ be the basechange of $\bar{\mathfrak{U}}_{\sigma} \to \tilde{\bar{X}}_{\sigma}$ by $\pi_{\sigma}$. 

Similarly, for every $\rho \in \bar{\mathscr{P}}^{[1]}$ and every $x \in \bar{X}_{\rho}$, since $\tilde{\bar{X}}_{\rho,x}$ is an étale local model for $\bar{\mathfrak{X}}$, we have a variety $\bar{\mathfrak{U}}_{\rho,x}$ equipped with étale maps $\bar{\mathfrak{U}}_{\rho,x} \to \bar{\mathfrak{X}}$ and $\bar{\mathfrak{U}}_{\rho,x} \to \tilde{\bar{X}}_{\rho,x}$. We let $\mathfrak{U}_{\rho,x}:= \bar{\mathfrak{U}}_{\rho,x} \times_{\bar{\mathfrak{X}}} \mathfrak{X}$ be the basechange of $\bar{\mathfrak{U}}_{\rho,x} \to \bar{\mathfrak{X}}$ by $\pi: \mathfrak{X} \to \bar{\mathfrak{X}}$. For any morphism $\pi_{\rho,x}: \tilde{X}_{\rho,x} \to \tilde{\bar{X}}_{\rho,x}$,we let $\tilde{\mathfrak{U}}_{\rho,x, \pi_{\rho,x}}:= \bar{\mathfrak{U}}_{\rho,x} \times_{\tilde{\bar{X}}_{\rho,x}} \tilde{X}_{\rho,x}$ be the basechange of $\bar{\mathfrak{U}}_{\rho,x} \to \tilde{\bar{X}}_{\rho,x}$ by $\pi_{\rho,x}$.

\begin{defn} \label{toricresolutions}
 We say that a resolution $\pi: \mathfrak{X} \to \bar{\mathfrak{X}}$ of a special toric degeneration $\bar{\mathfrak{X}} \to \mathcal{S}$ of K3-s to a log smooth degeneration $\mathfrak{X} \to \mathcal{S}$ is \emph{toric} if for every $\sigma \in \bar{\mathscr{P}}^{\max}$ there exists a toric blowup $\pi_{\sigma}: \tilde{X}_{\sigma} \to \tilde{\bar{X}}_{\sigma}$ such that $\mathfrak{U}_{\sigma} \cong \tilde{\mathfrak{U}}_{\sigma, \pi_{\sigma}}$, for every $\rho \in \bar{\mathscr{P}}^{[1]}$ and every $x \in \bar{X}_{\rho}$ there exists a toric blowup $\pi_{\rho,x}: \tilde{X}_{\rho,x} \to \tilde{\bar{X}}_{\rho,x}$  such that $\mathfrak{U}_{\rho,x} \cong \tilde{\mathfrak{U}}_{\rho,x, \pi_{\rho,x}}$, and $\pi$ is trivial at every point of $\bar{\mathfrak{X}}$ not contained in a codimension $1$ or $2$ stratum of $\bar{\mathfrak{X}}_0$. 

We say that a toric resolution $\pi: \mathfrak{X} \to \bar{\mathfrak{X}}$ is \emph{integral} if for every $\sigma \in \bar{\mathscr{P}}^{\max}$, the toric blowup $\pi_{\sigma}: \tilde{X}_{\sigma} \to \tilde{\bar{X}}_{\sigma}$ induces an integral subdivision of $\sigma$ and for every $\rho \in \bar{\mathscr{P}}^{[1]}$ and every $x \in \bar{X}_{\rho}$, the toric blowup $\pi_{\rho,x}: \tilde{X}_{\rho,x} \to \tilde{\bar{X}}_{\rho,x}$ induces an integral subdivision of:
 \begin{enumerate}
     \item If $x$ is a non-sigular point, the triangle defining $\tilde{\bar{X}}_{\rho, x}$ of Figure \ref{model1}. 
     \item If $x$ is a singular point, the trapezoid defining $\tilde{\bar{X}}_{\rho, x}$ of Figure \ref{model2}. 
 \end{enumerate} 
Moreover, we require that the subdivisions of (1) and (2) are of the form in the respective figures (possibly in an intermediate form if $\rho \in \bar{\mathscr{P}}^{[1]}$ has $r_{\rho} = 0$). We also require that for every $\rho \subseteq \sigma$ with $\rho \in \bar{\mathscr{P}}^{[1]}, \ \sigma \in \bar{\mathscr{P}}^{\max}$, the subdivision of $\rho$ induced by the subdivision of $\sigma$ agrees with the subdivision of $\rho$ induced by the subdivision of (1) or (2). Here and later, we assume that the cells of any subdivision don't self-intersect and that an intersection of any two cells is also a cell.
 
 We say that a toric and integral resolution is \emph{homogeneous} if for every $\rho \in \bar{\mathscr{P}}^{[1]}$ and $x,y \in \bar{X}_{\rho}$ singular points, the subdivisions of the trapezoid of (2) corresponding to the blowups $\pi_{\rho,x}: \tilde{X}_{\rho,x} \to \tilde{\bar{X}}_{\rho,x}$ and $\pi_{\rho,y}: \tilde{X}_{\rho,y} \to \tilde{\bar{X}}_{\rho,y}$ are the same. Note that the same is true for non-singular points, even when $r_{\rho} = 0$, by the definition of an integral resolution above.
\end{defn}

We are going to consider toric, integral, and homogeneous resolutions $\pi: \mathfrak{X} \to \bar{\mathfrak{X}}$ from now on. Sections \ref{resolvesmall} and \ref{resolve} give examples of such resolutions. We will use Notations \ref{generalblowupnotations} for exceptional curves and irreducible components of $\mathfrak{X}_0$, and for vertices $v \in \mathscr{P}^{[0]}$ of the dual intersection complex $(B, \mathscr{P})$ of $\mathfrak{X} \to \mathcal{S}$. 

From the local models, it is easy to see that the cone $NE(\pi) \subseteq NE(\mathfrak{X}_0)$ of curves contracted by $\pi$ is finitely generated and that we have
\begin{equation} \label{intnumbers}
\left\{ E_{\rho,k}  \ | \ \rho \in \bar{\mathscr{P}}^{[1]}, \ 1 \leq k \leq r_{\rho}\right\} \cup \left\{F_{\rho, p}, \ \rho \in \bar{\mathscr{P}}^{[1]}, \ 1 \leq p \leq l_{\rho} -1  \right\} \subseteq NE(\pi)
\end{equation}
(with possibly not all $F_{\rho,p}$ present if $r_{\rho} = 0$)
as in Section \ref{resolve}.

\begin{remark} \label{otherresolutions}
Considering non-homogeneous (but toric and integral) resolutions is also interesting. Such resolutions $\pi: \mathfrak{X} \to \bar{\mathfrak{X}}$ correspond to not requiring all the $E_{\rho,k},  \  1 \leq k \leq r_{\rho}$ (for a fixed $\rho \in \bar{\mathscr{P}}$) to be contained in the same irreducible component of $\mathfrak{X}_0$. These resolutions can't be constructed as in Proposition \ref{stronglogsmoothresolution} below, might not be projective (they are still proper), or might have $\mathfrak{X}$ an algebraic space. See Section \ref{relax} for a further discussion. At the level of dual intersection complexes, the singularity $x_{\rho}$ of the affine structure on $\left(\bar{B}, \bar{\mathscr{P}}\right)$ splits into the singularities of the affine structure on $(B, \mathscr{P})$ at the vertices corresponding to the irreducible components containing the $E_{\rho,k},  \  1 \leq k \leq r_{\rho}$. 

Assuming that $D$ is simple normal crossings, by considering intersection numbers of the curve classes of \eqref{intnumbers}  with the irreducible components, one can see that there can be no $\pi$-ample divisor $D'$ supported on $D$. If $D$ is not simple normal crossings, one can similarly argue that there is no $\pi$-ample PA-generated $D'$. Indeed, the extension by $0$ of the restriction $\alpha|_{\rho}$ of the corresponding $\alpha \in PA(B)$ (see Step 3 in the proof of Proposition \ref{ampledivisor}) would not be strictly convex in one of the subdivisions of the trapezoids corresponding to the resolutions of the local models at the singular points of $\bar{X}_{\rho}$. 
\end{remark}

In Sections \ref{resolvesmall} and \ref{resolve}, we defined the resolutions globally and looked at the local models to understand how the central fibre $\bar{\mathfrak{X}}_0$ and the dual intersection complex $\left(\bar{B}, \bar{\mathscr{P}}\right)$ transform after the blowup. We will now reverse this logic and define resolutions corresponding to certain subdivisions of $\left(\bar{B}, \bar{\mathscr{P}}\right)$ by gluing together the resolutions in the local models. We use the observations about PA-functions of Step 3 in the proof of Proposition \ref{ampledivisor}.

\begin{fact} \label{stronglogsmoothresolution}
Let $\bar{\mathfrak{X}} \to \mathcal{S}$ be a special toric degeneration of K3-s. Let $\left(\bar{B}, \bar{\mathscr{P}}\right)$ be its dual intersection complex and consider an integral subdivision $(B, \mathscr{P})$ (defined as a polyhedral manifold only) of $\left(\bar{B}, \bar{\mathscr{P}}\right)$ with every $\rho \in \bar{\mathscr{P}}^{[1]}$ such that $r_{\rho} \neq 0$ subdivided into edges of integral length $1$. We require that the cells of the subdivision don't self-intersect and that an intersection of any two cells is also a cell (i.e. $(B, \mathscr{P})$ satisfies Assumption \ref{manifoldass}). Suppose also that there exists an $\alpha \in PA(B)$ such that:
\begin{enumerate}[(a)]
    \item $\alpha(v) \geq 0$ for all $v \in \bar{\mathscr{P}}^{[0]}$.\footnote{This is a technical condition that is necessary to ensure that the ideal sheaves defining the blowups in the étale local models glue to an ideal sheaf on $\bar{\mathfrak{X}}$.}
    \item The restriction $\alpha|_{\sigma}$ of $\alpha$ to any $\sigma \in \bar{\mathscr{P}}^{\max}$ is strictly convex (on the induced subdivision of $\sigma$).
    \item The restriction $\alpha|_{\rho}$ of $\alpha$ to any $\rho \in \bar{\mathscr{P}}^{[1]}$ with $r_{\rho} \neq 0$ has a well-defined $\min(\alpha|_{\rho}) \subseteq \rho(\mathbb{Z})$ (see Notation \ref{minnotation}).
\end{enumerate}
Then:
\begin{enumerate}
    \item There exists a toric, integral, and homogeneous resolution $\pi: \mathfrak{X} \to \bar{\mathfrak{X}}$ to a log smooth and minimal log CY degeneration $\mathfrak{X} \to \mathcal{S}$ with dual intersection complex isomorphic to $(B, \mathscr{P})$ as a polyhedral manifold.
    \item The PA-function $\alpha$ defines a $\pi$-ample PA-generated divisor $\divisor(\alpha) = \sum_{v \in \mathscr{P}^{[0]}} \alpha(v) D_v$ on $\mathfrak{X}$ via \eqref{divmap} (with the coefficient at $D_v$ positive if $v \in \bar{\mathscr{P}}^{[0]}$). 
    \item Identifying the dual intersection complex of $\mathfrak{X} \to \mathcal{S}$ with $(B, \mathscr{P})$ under the isomorphism, the affine structure on $(B, \mathscr{P})$ extends to the complement of the vertices 
    \begin{equation} \label{extsinglocus}
    \Delta = \left\{ v_{\rho, p} \in \mathscr{P}^{[0]} \ | \ \min(\alpha|_{\rho}) = \{v_{\rho, p} \}, \ \rho \in \bar{\mathscr{P}}^{[1]}, \ r_{\rho} \neq 0  \right\}
    \end{equation}
    (as before, we understand $v_{\rho, p}$ as the vertex corresponding to the strict transform of an irreducible component of $\bar{\mathfrak{X}}_0$ if $p = 0, l_{\rho}$, so some of such vertices may appear more than once). 
\end{enumerate}
\end{fact}

\begin{proof}
By Proposition \ref{whenspecialk3}, $\bar{\mathfrak{X}} \to \mathcal{S}$ is special if and only if it is a divisorial log deformation and the generic fibre of $\bar{\mathfrak{X}} \to \mathcal{S}$ is smooth. The subdivision $(B, \mathscr{P})$ and the PA-function $\alpha \in PA(B)$ define the blowups in the local models and the corresponding ideal sheaves. We will construct these ideal sheaves and glue them to an ideal sheaf $\mathcal{I}$ on $\bar{\mathfrak{X}}$ supported on codimension $1$ and $2$ strata of $\bar{\mathfrak{X}}_0$. The blowup of $\mathcal{I}$ will give the desired resolution $\pi:\mathfrak{X} \to \bar{\mathfrak{X}}$.

\textbf{Step 1.} We consider the local models in codimension $2$ first. Fix a $\sigma \in \bar{\mathscr{P}}^{\max}$. Then the induced subdivision of the cone $\mathbf{C}\sigma$ over $\sigma$ gives rise to a toric blowup $\pi_{\sigma}: \tilde{X}_{\sigma} \to \tilde{\bar{X}}_{\sigma}$ of the local model $\tilde{\bar{X}}_{\sigma}$. Moreover, the canonical extension of $\alpha|_{\sigma}$ to a PL-function $\mathbf{C}\alpha|_{\mathbf{C}\sigma}$ on the cone $\mathbf{C}\sigma$ defines the support function for the induced subdivision of $\mathbf{C}\sigma$ since $\alpha|_{\sigma}$ is strictly convex by condition (b). By the discussion of \cite[Section 3]{T}, $\mathbf{C}\alpha|_{\mathbf{C}\sigma}$ gives rise to an ideal sheaf $\mathcal{I}_{\sigma}$ on $\tilde{\bar{X}}_{\sigma}$ such that the blowup of $\mathcal{I}_{\sigma}$ induces the given subdivision. 

Explicitly, the ideal sheaf $\mathcal{I}_{\sigma}$ can be described as follows. We let $I_{\sigma} \subseteq \widecheck{\mathbf{C}\sigma} \cap N$ be the monoid ideal $$I_{\sigma} := \left\langle -m_{\mathbf{C}\tau}, \ | \  \tau \in \mathscr{P}^{\max}, \ \tau \subseteq \sigma \right\rangle \cap \left(\widecheck{\mathbf{C}\sigma} \cap N \right)$$ where $m_{\mathbf{C}\tau}$ is the linear function defining $\mathbf{C}\alpha|_{\mathbf{C}\sigma}$ on $\mathbf{C}\tau$ (i.e. $-m_{\mathbf{C}\tau}, \ \tau \subseteq \sigma$ define the Cartier data of \cite[Section 3]{T} for the subdivision of $\sigma$). Then $\mathcal{I}_{\sigma}$ is the unique torus-invariant ideal sheaf on $\tilde{\bar{X}}_{\sigma}$ with $\Gamma(U_{\mathbf{C}\sigma}, \mathcal{I}_{\sigma}) \subseteq \kk[\widecheck{\mathbf{C}\sigma} \cap N]$ the monomial ideal corresponding to $I_{\sigma} \subseteq \widecheck{\mathbf{C}\sigma} \cap N$.

We note that the support of $\mathcal{I}_{\sigma}$ is contained in the union of codimension $1$ and $2$ toric strata of $\partial \tilde{\bar{X}}_{\sigma}$. Indeed, for every $v \in \bar{\mathscr{P}}^{[0]}$ with $v \subseteq \sigma$ let $u_v$ be the primitive generator of  $\mathbf{C}v$ and let $\tau \in \mathscr{P}^{\max}$ be such that $v \subseteq \tau \subseteq \sigma$. Then $\langle u_v, -m_{\mathbf{C}\tau} \rangle = - \alpha(v) \leq 0$ by condition (a). This implies that $I_{\sigma} \cap (\mathbf{C}v)^{\perp} \neq \varnothing$ so $\mathcal{I}_{\sigma}$ vanishes on the toric divisor corresponding to $v$. 

\textbf{Step 2.} We now consider the local models in codimension $1$. Fix a $\rho \in \bar{\mathscr{P}}^{[1]}$ and $x \in \bar{X}_{\rho}$. If $x$ is non-singular, we define the subdivision of the triangle $\rho'_x$ defining $\tilde{\bar{X}}_{\rho, x}$ of Figure \ref{model1} to be as in Figure \ref{model1} (possibly, we have an intermediate form if $r_{\rho} = 0$, in which case the subdivision is determined by the induced subdivision of $\rho$ and is the same for all non-singular $x \in \bar{X}_{\rho}$). If $x$ is singular, we define the subdivision of the trapezoid $\rho'_x$ defining $\tilde{\bar{X}}_{\rho, x}$ of Figure \ref{model2} to be as in Figure \ref{model2} with the vertex connected to both $(0,1,1)$ and $(1,1,1)$ corresponding to the $\min(\alpha|_{\rho}) \subseteq \rho(\mathbb{Z})$ of condition (c) (note that this does not depend on the choice of a singular point $x \in \bar{X}_{\rho}$). The subdivision of $\rho'_x$ gives rise to a toric blowup $\pi_{\rho,x}: \tilde{X}_{\rho,x} \to \tilde{\bar{X}}_{\rho,x}$ of the local model $\tilde{\bar{X}}_{\rho, x}$. 

Let $\alpha|_{\rho, x}^0$ be the extension by $0$\footnote{This is as in Step 3 in the proof of Proposition \ref{ampledivisor}.} of $\alpha|_{\rho}$ to $\rho'_x$. Here by an extension by $0$, we mean that $\alpha|_{\rho, x}^0$ is the PA-function that takes the same values as $\alpha$ at the points corresponding to the integer points of $\rho$ and takes value $0$ at $(0,1,1)$ and $(1,1,1)$ (or just $(0,1,1)$ if $x$ is non-singular). As in Step 3 in the proof of Proposition \ref{ampledivisor}, $\alpha|_{\rho, x}^0$ is strictly convex with respect to the subdivision of $\rho'_x$. So its canonical extension to a PL-function $\mathbf{C}\alpha|_{\mathbf{C}\rho, x}^0$ on the cone $\mathbf{C}\rho'_x$ over $\rho'_x$ defines the support function for the induced subdivision of $\mathbf{C}\rho'_x$. By the discussion of \cite[Section 3]{T}, $\mathbf{C}\alpha|_{\mathbf{C}\rho, x}^0$ gives rise to an ideal sheaf $\mathcal{I}_{\rho, x}$ on $\tilde{\bar{X}}_{\rho, x}$ such that the blowup of $\mathcal{I}_{\rho, x}$ induces the given subdivision.

An explicit description of $\mathcal{I}_{\rho, x}$ using a monoid ideal $I_{\rho,x} \subseteq \widecheck{\mathbf{C}\rho'_x} \cap N$ is similar to Step 1, and we can similarly show that the support of $\mathcal{I}_{\rho, x}$ is contained in the union of codimension $1$ and $2$ toric strata of $\partial \tilde{\bar{X}}_{\rho, x}$. However, we actually know more. Since $\alpha|_{\rho, x}^0$ is the extension by $0$, we have $\alpha|_{\rho, x}^0(\tilde{v}) = 0$ for any vertex $\tilde{v}$ of $\rho'_x$ that is not a vertex of $\rho$. Arguing as in Step 1, this implies that $I_{\rho,x} \cap (\mathbf{C}\tilde{\rho})^{\perp} \neq \varnothing$ for any edge $\tilde{\rho}$ of $\rho'_x$ that is not $\rho$ and that $I_{\rho,x}$ contains the vertex of $\widecheck{\mathbf{C}\rho'_x}$. So the support of $\mathcal{I}_{\rho,x}$ is contained in the toric stratum corresponding to $\rho$.

\textbf{Step 3.} We want to define an ideal sheaf $\mathcal{I}$ on $\bar{\mathfrak{X}}$ supported on codimension $1$ and $2$ strata of $\bar{\mathfrak{X}}_0$, that gives rise to the desired resolution. Let $$\bar{\mathfrak{U}}:= \bar{\mathfrak{X}} \ \Big\backslash \bigcup_{\sigma \in \bar{\mathscr{P}}^{\max}} \bar{X}_{\sigma} \cup \bigcup_{\rho \in \bar{\mathscr{P}}^{[1]}} \bar{X}_{\rho}$$ be the complement of codimension $1$ and $2$ strata of $\bar{\mathfrak{X}}_0$, $\bar{\mathfrak{U}} \xhookrightarrow{} \bar{\mathfrak{X}}$ be the natural inclusion, and $\mathcal{I}_{\bar{\mathfrak{U}}}$ be the trivial sheaf on $\bar{\mathfrak{U}}$. For every $\sigma \in \bar{\mathscr{P}}^{\max}$, let $\mathcal{I}_{\bar{\mathfrak{U}}_{\sigma}}$ be the ideal sheaf on $\bar{\mathfrak{U}}_{\sigma}$ that is the pullback of the ideal sheaf $\mathcal{I}_{\sigma}$ on $\tilde{\bar{X}}_{\sigma}$ constructed in Step 1 via the étale map $\bar{\mathfrak{U}}_{\sigma} \to \tilde{\bar{X}}_{\sigma}$ (using the notations of Definition \ref{toricresolutions}). Similarly, for every $\rho \in \bar{\mathscr{P}}^{[1]}$ and $x \in \bar{X}_{\rho}$, let $\mathcal{I}_{\bar{\mathfrak{U}}_{\rho,x}}$ be the ideal sheaf on $\bar{\mathfrak{U}}_{\rho,x}$ that is the pullback of the ideal sheaf $\mathcal{I}_{\rho,x}$ on $\tilde{\bar{X}}_{\rho,x}$ constructed in Step 2 via the étale map $\bar{\mathfrak{U}}_{\rho,x} \to \tilde{\bar{X}}_{\rho,x}$. 

The collection of maps $\bar{\mathfrak{U}} \xhookrightarrow{} \bar{\mathfrak{X}}$, $\bar{\mathfrak{U}}_{\sigma} \to \bar{\mathfrak{X}}$ for $\sigma \in \bar{\mathscr{P}}^{\max}$, and $\bar{\mathfrak{U}}_{\rho,x} \to \bar{\mathfrak{X}}$ for $\rho \in \bar{\mathscr{P}}^{[1]}$ and $x \in \bar{X}_{\rho}$ gives an étale covering of $\bar{\mathfrak{X}}$. Therefore, we can apply étale descent for quasi-coherent sheaves (see, e.g. \cite[Theorem 2.3]{decent}) to construct an ideal sheaf $\mathcal{I}$ on $\bar{\mathfrak{X}}$ such that its pullback via the maps in the covering gives the ideal sheaves $\mathcal{I}_{\bar{\mathfrak{U}}}, \ \mathcal{I}_{\bar{\mathfrak{U}}_{\sigma}}$, and $\mathcal{I}_{\bar{\mathfrak{U}}_{\rho,x}}$ respectively. To do that, we need to check that we can define isomorphisms of the pullbacks of the ideal sheaves $\mathcal{I}_{\bar{\mathfrak{U}}}, \ \mathcal{I}_{\bar{\mathfrak{U}}_{\sigma}}, \ \mathcal{I}_{\bar{\mathfrak{U}}_{\rho,x}}$ to the double fibre products that satisfy the cocycle condition on the triple fibre products (i.e. we can define a descent datum for quasi-coherent sheaves, see \cite[Definition 2.1]{decent}). This follows from the correspondence between the ideals in the local toric models and the support functions, the fact that $(\alpha|_{\sigma})|_{\rho} = (\alpha|_{\rho, x}^{0})|_{\rho}$ for every $\rho \subseteq \sigma$ with $\rho \in \bar{\mathscr{P}}^{[1]}, \ 
\sigma \in \bar{\mathscr{P}}^{\max}$, and every $x \in \bar{X}_{\rho}$\footnote{This is since $(\alpha|_{\rho, x}^{0})|_{\rho} = \alpha|_{\rho}$ and both $\alpha|_{\rho}$ and $\alpha|_{\sigma}$ arise as the restrictions of a globally defined $\alpha \in PA(B)$.} (for fibre products not involving $\bar{\mathfrak{U}}$), and the fact that $\mathcal{I}_{\rho, x}$ is supported on the union of codimension $1$ and $2$ toric strata of $\partial \tilde{\bar{X}}_{\rho, x}$ and $\mathcal{I}_{\rho,x}$ is supported on the toric stratum corresponding to $\rho$ (for fibre products involving $\bar{\mathfrak{U}}$).

\textbf{Step 4.} We let $\pi: \mathfrak{X} \to \bar{\mathfrak{X}}$ be the blowup of $\mathcal{I}$. Then the construction of Step 3 immediately implies that $\mathfrak{U}_{\sigma} \cong \tilde{\mathfrak{U}}_{\sigma, \pi_{\sigma}}$ for every $\sigma \in \bar{\mathscr{P}}^{\max}$ and $\pi_{\sigma}: \tilde{X}_{\sigma} \to \tilde{\bar{X}}_{\sigma}$ as in the construction and $\mathfrak{U}_{\rho,x} \cong \tilde{\mathfrak{U}}_{\rho,x, \pi_{\rho,x}}$ for every $\rho \in \bar{\mathscr{P}}^{[1]}$, every $x \in \bar{X}_{\rho}$, and $\pi_{\rho,x}: \tilde{X}_{\rho,x} \to \tilde{\bar{X}}_{\rho,x}$ as in the construction (using the notations of Definition \ref{toricresolutions}). Further, the fact that $\mathcal{I}_{\bar{\mathfrak{U}}}$ is trivial implies that $\pi$ is trivial at every point of $\bar{\mathfrak{X}}$ not contained in a codimension $1$ or $2$ stratum of $\bar{\mathfrak{X}}_0$. So $\pi:\mathfrak{X} \to \bar{\mathfrak{X}}$ is toric. Our definitions of the blowups of local models in Steps 1 and 2 imply that $\pi:\mathfrak{X} \to \bar{\mathfrak{X}}$ is integral and homogeneous. The fact that $\mathfrak{X} \to \mathcal{S}$ is minimal log CY follows as in Proposition \ref{minlogcy42}, so we have proved (1). 

Now (2) follows exactly as in Step 3 in the proof of Proposition \ref{ampledivisor} and (3) follows since the analysis of the local models of Sections \ref{toricmodelcodim2} and \ref{localmodelgen} implies that the irreducible components $D_{v_{\rho, p}}$ contain the curves $E_{\rho,k},  \  1 \leq k \leq r_{\rho}$ and all the other irreducible components of $\mathfrak{X}_0$ are toric. 
\end{proof}

Proposition \ref{stronglogsmoothresolution} gives us a tropical way to work with resolutions. Note that the correspondence in Proposition \ref{stronglogsmoothresolution} is actually one-to-one.

\begin{cor} \label{onetoonetropical}
There is a one-to-one correspondence $$\left\{(B, \mathscr{P}), \alpha \in PA(B)\right\} \Longleftrightarrow \left\{ \pi:\mathfrak{X} \to \bar{\mathfrak{X}}, D'\right\}$$
where $\bar{\mathfrak{X}} \to \mathcal{S}$ is a special toric degeneration of K3-s with dual intersection complex $\left(\bar{B}, \bar{\mathscr{P}}\right)$, $(B, \mathscr{P})$ is an integral subdivision of $\left(\bar{B}, \bar{\mathscr{P}}\right)$ with every $\rho \in \bar{\mathscr{P}}^{[1]}$ such that $r_{\rho} \neq 0$ subdivided into edges of integral length $1$, $\alpha \in PA(B)$ satisfies conditions (a), (b), (c) of Proposition \ref{stronglogsmoothresolution}, $\pi:\mathfrak{X} \to \bar{\mathfrak{X}}$ is a toric, integral, and homogeneous resolution, and $D' = \sum_{v \in \mathscr{P}^{[0]}} a_v D_v$ is a $\pi$-ample PA-generated divisor with $a_v \geq 0$ if $D_v$ is the strict transform of an irreducible component of $\bar{\mathfrak{X}}_0$.
\end{cor}

\begin{proof}
One direction is Proposition \ref{stronglogsmoothresolution}. For the other direction, note that for every toric, integral, and homogeneous resolution $\pi:\mathfrak{X} \to \bar{\mathfrak{X}}$, the dual intersection complex $(B, \mathscr{P})$ of $\mathfrak{X} \to \mathcal{S}$ is an integral subdivision of $\left(\bar{B}, \bar{\mathscr{P}}\right)$ by the analysis of the local models. By an argument as in Step 3 in the proof of Proposition \ref{ampledivisor}, every $\pi$-ample PA-generated divisor $D'$ corresponds to a PA-function $\alpha \in PA(B)$ such that the restriction $\alpha|_{\sigma}$ of $\alpha$ to any $\sigma \in \bar{\mathscr{P}}^{\max}$ is strictly convex (on the induced subdivision of $\sigma$) and the restriction $\alpha|_{\rho}$ of $\alpha$ to any $\rho \in \bar{\mathscr{P}}^{[1]}$ with $r_{\rho} \neq 0$ has a well-defined $\min(\alpha|_{\rho}) \subseteq \rho(\mathbb{Z})$. The fact that $a_v \geq 0$ if $D_v$ is the strict transform of an irreducible component of $\bar{\mathfrak{X}}_0$ implies that $\alpha(v) \geq 0$ for all $v \in \bar{\mathscr{P}}^{[0]}$.
\end{proof}

\subsubsection{Admissible and strongly admissible resolutions} \label{admissibleandstronglyso}

We are ready to define admissible resolutions (i.e. resolutions allowed in Conjecture \ref{conjecture}), the data of Basic Setup \ref{setupforproof}, and the PL-isomorphism $\Phi: B \to \bar{B}$ in general. In the next section, we will show that any special toric degeneration admits a (strongly) admissible resolution.

\begin{defn} \label{stronglyadmissible}
We say that a resolution $\pi: \mathfrak{X} \to \bar{\mathfrak{X}}$ of a special toric degeneration $\bar{\mathfrak{X}} \to \mathcal{S}$ of K3-s is \emph{strongly admissible} if it is a toric, integral, and homogeneous resolution admitting a $\pi$-ample PA-generated divisor $D'$ such that:
\begin{enumerate}
    \item The divisor $D'$ is simple normal crossings and effective.
    \item The divisor $D_{\irrel}$ (defined as in the proof of Proposition \ref{P2}) is $\mathbb{Q}$-Cartier and $\pi$-nef.
\end{enumerate}
\end{defn}

\begin{remarks}
\begin{enumerate}
    \item Note that conditions (1) and (2) in Definition \ref{stronglyadmissible} correspond to conditions (1) and (2) in Proposition \ref{P2general} that we use to define the extended intrinsic mirror. In particular, it follows from Proposition \ref{ampledivisorsmall} that the resolutions of Section \ref{resolvesmall} are strongly admissible. Similarly, the resolutions of Corollary \ref{furtherresolvetosatisfy} and Remark \ref{satisfybutnotsnc} in Section \ref{resolve} are strongly admissible.
    \item It will follow from the results of Section \ref{scatint} that we don't actually need the existence of a $\pi$-ample PA-generated divisor $D'$ in Definition \ref{stronglyadmissible} (satisfying conditions (1) and (2)) to construct the extended intrinsic mirror. As a result, Conjecture \ref{conjecture} will follow for all toric, integral, and homogeneous resolutions, see Remark \ref{stronglyadmissibleweakening}(1). 
\end{enumerate}
\end{remarks}

\begin{defn} \label{admissibledef}
We say that a resolution $\pi: \mathfrak{X}' \to \bar{\mathfrak{X}}$ of a special toric degeneration $\bar{\mathfrak{X}} \to \mathcal{S}$ of K3-s is \emph{admissible} if it factors as $\mathfrak{X}' \to \mathfrak{X} \to \bar{\mathfrak{X}}$ with $\mathfrak{X} \to \bar{\mathfrak{X}}$ strongly admissible and $\mathfrak{X}' \to \mathfrak{X}$ a logarithmic modification.
\end{defn}

By Proposition \ref{furtherresolution}, it is enough to prove Conjecture \ref{conjecture} for strongly admissible resolutions. Using the results of Section \ref{tropicalresolutions}, we can understand strongly admissible resolutions tropically. 

\begin{notation}
For $\alpha \in PA(B)$, denote by $(\supp \alpha, \mathscr{P}_{\alpha})$ the polyhedral subcomplex of $(B, \mathscr{P})$ consisting of the following cells of $\mathscr{P}$. A cell  $\tau \in \mathscr{P}$ is in $\mathscr{P}_{\alpha}$ if and only if for all vertices $v \subseteq \tau, \ v \in \mathscr{P}^{[0]}$, we have $\alpha(v) \neq 0$. Note that $(\supp \alpha, \mathscr{P}_{\alpha})$ might consist of multiple components of different dimensions.
\end{notation}

\begin{fact} \label{stronglyadmissibletosubdivision}
\begin{sloppypar}Under the correspondence of Corollary \ref{onetoonetropical} a pair $\left\{ \pi:\mathfrak{X} \to \bar{\mathfrak{X}}, D'\right\}$ with $\pi:\mathfrak{X} \to \bar{\mathfrak{X}}$ a strongly admissible resolution and $D'$ a $\pi$-ample PA-generated divisor corresponds to a pair $\left\{(B, \mathscr{P}), \alpha \in PA(B)\right\}$ with $(B, \mathscr{P})$ the dual intersection complex of $\mathfrak{X} \to \mathcal{S}$ and $\alpha$ such that $\divisor(\alpha) = D'$ via \eqref{divmap} satisfying the following conditions:\end{sloppypar}
\begin{enumerate}
    \item $\alpha(v) \geq 0$ for all $v \in \mathscr{P}^{[0]}$ and $\supp \alpha$ is such that every $\sigma \in \mathscr{P}_{\alpha}^{[2]}$ is a standard triangle and every $\rho \in \mathscr{P}_{\alpha}^{[1]}$ is of integral length $1$.
    \item Suppose in addition that $MD_{\irrel}$ is PA-generated for those $M \in \mathbb{N}$ for which it is Cartier. Then $\alpha_{\comp}:B \to \mathbb{R}$ defined by 
    $$\alpha_{\comp}(v) = \begin{cases}\begin{aligned} 1, &\text{ if } v \notin \supp{\alpha} \\ 0, &\text{ if } v \in \supp{\alpha} \end{aligned} \end{cases}$$ is a well-defined PA-function on $B$ with $N\alpha_{\comp} \in PA(B)$ (i.e. it is integral\footnote{Recall from Section \ref{generalextensionstuff} that $PA(B):=PA(B, \mathbb{Z})$ stands for the group of integral PA-functions on $B$.}) for some $N \in \mathbb{N}$. Moreover, it satisfies the following conditions:
    \begin{enumerate}
        \item The restriction $(N\alpha_{\comp})|_{\sigma}$ of $N\alpha_{\comp}$ to any $\sigma \in \bar{\mathscr{P}}^{\max}$ is convex (but not necessarily strictly convex).
        \item For every $\rho \in \bar{\mathscr{P}}^{[1]}$ with $r_{\rho} \neq 0$ we either have $\supp (N\alpha_{\comp}) \cap \rho = \varnothing$ or $\supp (N\alpha_{\comp}) \cap \rho = \rho$.
    \end{enumerate}
\end{enumerate}
\end{fact}

\begin{proof}
    The description of the correspondence follows from the proof of Corollary \ref{onetoonetropical}. Condition (1) on $\alpha$ corresponds to condition (1) on $D'$ in Definition \ref{stronglyadmissible} of being strongly admissible. For condition (2) note that $N\alpha_{\comp}$ corresponds to $ND_{\irrel}$ via \eqref{divmap}, provided that $ND_{\irrel}$ is Cartier. Further, similarly to Step 3 in the proof of Proposition \ref{ampledivisor}, one can show that a PA-generated divisor $\divisor(\beta)$ is $\pi$-nef if and only if $\beta \in PA(B)$ satisfies the analogues of conditions (1) and (2) in Step 3 on $\beta$ with strict convexity replaced by convexity and the minimum condition $\min(\beta|_{\rho}) = \{v_{\rho, p}\}$ replaced by $v_{\rho, p} \in \min(\beta|_{\rho})$. Now, $(N\alpha_{\comp})$ satisfying the analogue of condition (1) in Step 3 is equivalent to (2)(a). For any $\rho \in \bar{\mathscr{P}}^{[1]}$ with $r_{\rho} = 0$, the analogue of condition (2) in Step 3 for $(N\alpha_{\comp})|_{\sigma}$ is equivalent to requiring that $(N\alpha_{\comp})_{\rho}$ is convex on $\rho$. This is implied by (2)(a), so it suffices to check that for $\rho \in \bar{\mathscr{P}}^{[1]}$ with $r_{\rho} \neq 0$, the analogue of condition (2) in Step 3 is equivalent to (2)(b). Indeed, it is clear from the setup that (2)(b) implies the analogue of condition (2) in Step 3. Conversely, if neither $\supp (N\alpha_{\comp}) \cap \rho = \varnothing$ nor $\supp (N\alpha_{\comp}) \cap \rho = \rho$, then $v_{\rho, p} \notin \min(N\alpha_{\comp}|_{\rho})$ and the analogue of condition (2) in Step 3 is not satisfied.
\end{proof}

A generalization of Proposition \ref{minlogcy42} to strongly admissible resolutions in the general case follows from the setup.

\begin{fact}
    A strongly amissible resolution $\mathfrak{X} \to \mathcal{S}$ is minimal log CY. The dual intersection complex $(B, \mathscr{P})$ satisfies Assumption \ref{manifoldass}.
\end{fact}

\begin{proof}
The fact that $\mathfrak{X} \to \mathcal{S}$ is minimal log CY follows as in the proof of Proposition \ref{minlogcy42} (we have also already proved this for any toric, integral, and homogeneous resolution via Corollary \ref{onetoonetropical} and Proposition \ref{stronglogsmoothresolution}). 

$(B, \mathscr{P})$ satisfies Assumption \ref{manifoldass} since so does $\left(\bar B, \bar{\mathscr{P}}\right)$ and $(B, \mathscr{P})$ is a subdivision of $\left(\bar B, \bar{\mathscr{P}}\right)$ such that the cells of the subdivision don't self-intersect and that an intersection of any two cells is also a cell (see Definition \ref{toricresolutions} and Proposition \ref{stronglogsmoothresolution}).
\end{proof}

As in Section \ref{globalresolvesmall}, let $P$ be a well-chosen monoid (see Definition \ref{wellchosendef}) with a face $K$ containing the classes of the contracted curves. We also have an analogue of Proposition \ref{ampledivisor} and Corollary \ref{furtherresolvetosatisfy} by design.

\begin{fact} \label{canextendifstronglyadmissible}
     A strongly amissible resolution $\mathfrak{X} \to \bar{\mathfrak{X}}$ satisfies the assumptions of Proposition \ref{P2general}.
\end{fact}

\begin{proof}
    The resolution $\mathfrak{X} \to \bar{\mathfrak{X}}$ satisfies the assumptions of Proposition \ref{P1general} and assumption (3) of Proposition \ref{P2general} by Definition \ref{stronglyadmissible} and the same arguments as in Step 1 in the proof of Proposition \ref{ampledivisor}. Assumptions (1) and (2) of Proposition \ref{P2general} are the same as assumptions (1) and (2) of Definition \ref{stronglyadmissible}. 
\end{proof}

We define the initial slab functions for  $\bar{\mathfrak{X}} \to \mathcal{S}$ as in \eqref{choiceslab} with the same convention for the choice of slabs $\underline{\rho}, \underline{\rho}' \in \tilde{\bar{\mathscr{P}}}^{[1]}$ as in Section \ref{global42}. For any strongly admissible resolution $\pi: \mathfrak{X} \to \bar{\mathfrak{X}}$, we have now defined all the necessary data in Basic Setup \ref{setupforproof}. Observation \ref{extmirroralgebraic} on when the extended intrinsic mirror is algebraic holds in this context. 

Finally, for any strongly admissible resolution $\pi:\mathfrak{X} \to \bar{\mathfrak{X}}$, we need to define a PL-isomorphism $\Phi: (B, \mathscr{P}) \to \left(\bar{B}, \bar{\mathscr{P}}\right)$ generalizing Construction \ref{plmap}.

\begin{cons} \label{plmapadmissible}
By Proposition \ref{stronglogsmoothresolution}(3), the affine structure on $(B, \mathscr{P})$ extends to the complement of singularities at vertices of the form $v_i$ for $ 1 \leq i \leq \bar{m}$ and at most one vertex of the form $v_{\rho, p}, \ 1 \leq p \leq l_{\rho} -1$ for $\rho \in \bar{\mathscr{P}}^{[1]}$, so we are in a similar setup to Construction \ref{plmap}. As in Construction \ref{plmap}, let $\mathscr{P}_{\coar}$ be the obvious coarsened polyhedral decomposition on $B$ such that $\left(B, \mathscr{P}_{\coar}\right) \cong \left(\bar B, \bar{\mathscr{P}}\right)$ as polyhedral manifolds. Then $\left(B, \mathscr{P}_{\coar}\right)$ is an affine manifold with singularities that satisfies Construction \ref{affchartstool} by the same argument as in Construction \ref{plmap}, and we define $\Phi: \left(B, \mathscr{P}_{\coar}\right) \to \left(\bar{B}, \bar{\mathscr{P}}\right)$, $\Phi(\tau)$ for $\tau \in \mathscr{P}$ or $\tau \in \mathscr{P}_{\coar}$ (in which case $\Phi(\tau) \in \bar{\mathscr{P}}$), $\Phi(\mathfrak{p})$ for a codimension $0$ wall $\mathfrak{p}$ on $(B, \mathscr{P})$ and $\Phi(\mathfrak{b})$ for a slab $\mathfrak{b}$ on $(B, \mathscr{P})$ in exactly the same way. 
\end{cons}

\subsubsection{Existence of strongly admissible resolutions}

We now show that any special toric degeneration $\bar{\mathfrak{X}} \to \mathcal{S}$ of K3-s admits a strongly admissible resolution. First, we prove a weaker result to illustrate how Proposition \ref{stronglyadmissibletosubdivision} can be applied in practice. 

\begin{fact} \label{weakadmissibleresolutionresult}
Let $\bar{\mathfrak{X}} \to \mathcal{S}$ be a special toric degeneration of K3-s with dual intersection complex $\left(\bar{B}, \bar{\mathscr{P}}\right)$ and suppose that all the $\sigma \in \bar{\mathscr{P}}^{\max}$ are either a standard triangle or one of the 16 reflexive polytopes in dimension $2$ (up to an action of $AGL(2, \mathbb{Z})$). Then there exists a strongly admissible resolution $\pi: \mathfrak{X} \to \bar{\mathfrak{X}}$. Moreover, we may assume that the divisor $D$ of $\mathfrak{X} \to \mathcal{S}$ is simple normal crossings.
\end{fact}

\begin{proof}
\textbf{Step 1.} By Proposition \ref{stronglyadmissibletosubdivision}, it suffices to construct a subdivision $(B, \mathscr{P})$ of $\left(\bar{B}, \bar{\mathscr{P}}\right)$ and an $\alpha \in PA(B)$ such that:
\begin{enumerate}
    \item All the cells $\sigma \in \mathscr{P}^{\max}$ are standard triangles. 
    \item $\alpha(v) > 0$ for all $v \in \mathscr{P}^{[0]}$.
    \item The restriction $\alpha|_{\sigma}$ of $\alpha$ to any $\sigma \in \bar{\mathscr{P}}^{\max}$ is strictly convex (on the induced subdivision of $\sigma$).
    \item The restriction $\alpha|_{\rho}$ of $\alpha$ to any $\rho \in \bar{\mathscr{P}}^{[1]}$ with $r_{\rho} \neq 0$ has a well-defined $\min(\alpha|_{\rho}) \subseteq \rho(\mathbb{Z})$.
\end{enumerate}       
Indeed, assumptions (1) and (2) together imply condition (2) of Proposition \ref{stronglyadmissibletosubdivision}. Assumption (2) implies that $\alpha_{\comp} = 0 \in PA(B)$ which is a well-defined convex PA-function satisfying conditions (2)(a) and (2)(b) of Proposition \ref{stronglyadmissibletosubdivision} for any subdivision $(B, \mathscr{P})$ of $\left(\bar{B}, \bar{\mathscr{P}}\right)$. Assumption (1) implies that $D$ is simple normal crossings (so $D'$ is also simple normal crossings).

\textbf{Step 2.} Let $(B, \mathscr{P})$ be the subdivision of $\left(\bar{B}, \bar{\mathscr{P}}\right)$ obtained by introducing, for every $\sigma \in \bar{\mathscr{P}}^{\max}$ that is not a standard triangle, a vertex $v_{\sigma} \in \mathscr{P}^{[0]}$ at the unique interior integral point of $\sigma$ and edges $\rho = \langle v_{\sigma}, v \rangle \in \mathscr{P}^{[1]}$ for every $v \in \mathscr{P}^{[0]}$ an integral point of an edge $\rho \in \bar{\mathscr{P}}^{[1]}$ with $\rho \subseteq \sigma$. It follows from the classification of the 16 reflexive polytopes (see, e.g. \cite[Fig.1]{KOS}) that this procedure produces a subdivision satisfying assumption (1) above. As in Step 2 in the proof of Proposition \ref{ampledivisor}, it suffices to find an $\alpha \in PA(B)$ satisfying assumptions (3) and (4) above.

\textbf{Step 3.} To find an $\alpha \in PA(B)$ satisfying assumptions (3) and (4) above, it is enough to define a PA-function $\alpha_{\partial}$ on the 1-dimensional polyhedral subcomplex $(\partial, \mathscr{P}_{\partial})$ that is given by the induced subdivision of the codimension $1$ skeleton $(\partial, \bar{\mathscr{P}}_{\partial})$ (with $\bar{\mathscr{P}}_{\partial} = \bar{\mathscr{P}}^{[0]} \cup \bar{\mathscr{P}}^{[1]}$) of $\bar{\mathscr{P}}$, such that $(\alpha_{\partial})|_{\rho}$ is strictly convex for every $\rho \in \bar{\mathscr{P}}^{[1]}$ and has well-defined $\min((\alpha_{\partial})|_{\rho}) \subseteq \rho(\mathbb{Z})$. Indeed, given such a function $\alpha_{\partial}$, the PA-function $\alpha \in PA(B)$ given by $$\alpha(v) = \begin{cases}\begin{aligned} N\alpha_{\partial}(v), &\text{ if } v \neq v_{\sigma}, \ \sigma \in \bar{\mathscr{P}}^{\max} \\ 0, \ \ \ \ \ \ \ \  &\text{ if } v = v_{\sigma}, \ \sigma \in \bar{\mathscr{P}}^{\max} \end{aligned} \end{cases}$$  for $N \in \mathbb{N}, \ N \gg 0$ and chosen so that $\alpha$ is integral, satisfies assumptions (3) and (4) above.

\textbf{Step 4.} Let $\bar{G}$ be a graph supported on $(\partial, \bar{\mathscr{P}}_{\partial})$ and choose an acyclic orientation on $\bar{G}$.\footnote{It is classical that every graph admits an acyclic orientation. Indeed, to obtain an acyclic orientation of a graph $G$, suppose that $V(G) = \{v_1, \dots, v_n \}$ for some $n \in \mathbb{N}$ and orient every $\langle v_i, v_j \rangle \in E(G)$ with $1 \leq i < j \leq n$ from $v_i$ to $v_j$.} Then it gives rise to an acyclic orientation on a graph $G$ supported on $(\partial, \mathscr{P}_{\partial})$ by either orienting all the edges subdividing every $\rho \in \bar{\mathscr{P}}_{\partial}^{[1]}$ in the same direction as $\rho$\footnote{Then, similarly to Remark \ref{restrictedclassofresolutions}(1) for the resolutions of Section \ref{resolve}, the dual intersection complex of the corresponding resolution would have all the singularities contained at the vertices corresponding to the strict transforms of the original irreducible components.} or modifying the orientation in such a way that there is a unique vertex $v \subseteq \rho, \ v \in \mathscr{P}_{\partial}^{[0]}$ with two adjacent edges $\rho_1, \rho_2$ with orientation vectors pointing into $v$ and no vertices $v' \subseteq \rho, \ v' \in \mathscr{P}_{\partial}^{[0]}$ with two adjacent edges $\rho'_1, \rho'_2$ with orientation vectors pointing out of $v$. Then let $\alpha_{\partial}$ be given by $\alpha_{\partial}(v) = N_v \in \mathbb{N}$ for $v \in \mathscr{P}_{\partial}^{[0]}$ with the $N_v$ chosen so that for any $\rho = \langle v, v' \rangle \in \mathscr{P}_{\partial}^{[1]}$ with orientation vector on $\rho$ pointing from $v$ to $v'$, we have $N_v > 2 N_{v'}$. The fact that the orientation on $G$ is acyclic implies that this system of inequalities is consistent and gives rise to a well-defined $\alpha_{\partial}$. Moreover, $\left(\alpha_{\partial}\right)|_{\rho}$ is strictly convex for every $\rho \in \bar{\mathscr{P}}^{[1]}$ with well-defined $\min((\alpha_{\partial})|_{\rho}) \subseteq \rho(\mathbb{Z})$.
\end{proof}

\begin{obs}
    We make a few observations:
    \begin{enumerate}
        \item Proposition \ref{weakadmissibleresolutionresult} directly applies only in the case that all the $\sigma \in \bar{\mathscr{P}}^{\max}$ are either a standard triangle or one of the 16 reflexive polytopes in dimension $2$. Indeed, we need every $\sigma \in \bar{\mathscr{P}}^{\max}$ to have a unique (to get a subdivision with all maximal cells standard triangles) integral interior point and reflexive polytopes are precisely such polytopes in dimension $2$. 
        \item In the case that all the $\sigma \in \bar{\mathscr{P}}^{\max}$ are either a standard triangle or contain an integral interior point, we can argue as in the proof of  Proposition \ref{weakadmissibleresolutionresult} (choosing $v_{\sigma} \in \mathscr{P}^{[0]}$ to be any integral interior point) that there exists a toric, integral, and homogeneous resolution $\mathfrak{X} \to \bar{\mathfrak{X}}$ admitting a $\pi$-ample PA-generated divisor $D'$.
        \item In particular, (2) implies that any special toric degeneration $\bar{\mathfrak{X}} \to \mathcal{S}$ of K3-s admits a toric, integral, and homogeneous resolution $\mathfrak{X} \to \bar{\mathfrak{X}}$ with a $\pi$-ample PA-generated divisor $D'$ after a basechange by $R \to R, \ t \mapsto t^3$ (where $t$ is the uniformizer of $R$). Indeed, such a basechange rescales the integral affine structure on $\left(\bar{B}, \bar{\mathscr{P}}\right)$ by $3$ and every integral polygon $\sigma \in \bar{\mathscr{P}}^{\max}$ has an integral interior point after such a rescaling.\footnote{This follows from the classification of lattice polygons without interior lattice points, see, e.g. \cite[§5]{Kh}.}
    \end{enumerate}
\end{obs}

To obtain a strongly admissible resolution in general, we can argue inductively similarly to Step 4 in the proof of Proposition \ref{ampledivisor}. Unlike the setup of Proposition \ref{ampledivisor}, we have control over the subdivision by appealing to the correspondence of Proposition \ref{stronglyadmissibletosubdivision}.

\begin{defn}
Let $(\bar{B}, \bar{\mathscr{P}})$ be a polyhedral manifold and let $v \in B(\mathbb{Z})$ be an integral point. Consider a polyhedral manifold $(B, \mathscr{P})$ obtained by subdividing every $\sigma \in \bar{\mathscr{P}}^{\max}$ with $v \in \sigma$ by the cells 
$$\left\{\Conv(\tau, v) \ | \ \tau \in \bar{\mathscr{P}}, \ \tau \subseteq \sigma, \ v \notin \tau \right\}.$$ We say that $(B, \mathscr{P})$ is obtained from $(\bar{B}, \bar{\mathscr{P}})$ by \emph{pulling $v$}.
\end{defn}

\begin{lemma} \label{pullingresolution}
Let $(\bar{B}, \bar{\mathscr{P}})$ be a polyhedral manifold of dimension $2$. Then one can subdivide $(\bar{B}, \bar{\mathscr{P}})$ to a polyhedral manifold $(B, \mathscr{P})$ with all the cells $\sigma \in \mathscr{P}^{\max}$ standard triangles by pulling a sequence of integral points.
\end{lemma}

\begin{proof}
Start pulling integral points and producing more and more refined integral subdivisions of $(\bar{B}, \bar{\mathscr{P}})$. Since $\bar{B}$ is compact, it is enough to show that for any integral subdivision $\left(B', \mathscr{P}'\right)$ such that not all $\sigma \in \mathscr{P}'^{\max}$ are standard simplices, one has a non-trivial pull of an integral point of $B'(\mathbb{Z})$. Indeed, consider a $\sigma \in \mathscr{P}'^{\max}$ that is not a standard simplex or a standard square. Then $\sigma$ has an integral point $v \in \sigma$ that is not one of the vertices $v \in \sigma, \ v \in \mathscr{P}'^{[0]}$ (i.e. either $v \in \Int(\sigma)$ or $v \in \Int(\rho)$ for $\rho \subseteq \sigma, \ \rho \in \mathscr{P}'^{[1]}$). But then pulling $v$ induces a non-trivial subdivision of $\sigma$. If $\sigma$ is a standard square, then pulling any vertex of $\sigma$ produces a non-trivial subdivision.
\end{proof}

\begin{remark} \label{nononsingularsubdivision}
Note that the analogue of Lemma \ref{pullingresolution} does not hold in higher dimensions since in dimension $n \geq 3$ there exist integral simplices $\sigma$ with the only integral points at the vertices (called \emph{elementary} or \emph{empty} simplices, see, e.g. \cite{HZ2}) and not $AGL(n, \mathbb{Z})$-equivalent to a standard simplex of dimension $n$. 
\end{remark}

\begin{fact} \label{strongadmissibleresolutionresult}
Any special toric degeneration $\bar{\mathfrak{X}} \to \mathcal{S}$ of K3-s admits a strongly admissible resolution $\pi: \mathfrak{X} \to \bar{\mathfrak{X}}$. Moreover, we may assume that the divisor $D$ of $\mathfrak{X} \to \mathcal{S}$ is simple normal crossings.
\end{fact}

\begin{proof}
Let $\left( \bar{B}, \bar{\mathscr{P}}\right)$ be the dual intersection complex of $\bar{\mathfrak{X}} \to \mathcal{S}$ and let $(B, \mathscr{P})$ be an integral subdivision of $\left( \bar{B}, \bar{\mathscr{P}}\right)$ with all the cells $\sigma \in \mathscr{P}^{\max}$ standard triangles, obtained by pulling a sequence of integral points (such a subdivision exists by Lemma \ref{pullingresolution}). Similarly to Steps 1 and 2 in the proof of Proposition \ref{weakadmissibleresolutionresult} (using Proposition \ref{stronglyadmissibletosubdivision} and an argument as in Step 2 in the proof of Proposition \ref{ampledivisor}), it suffices to find an $\alpha \in PA(B)$ that satisfies:
\begin{enumerate}
    \item[$(*)$] The restriction $\alpha|_{\sigma}$ of $\alpha$ to any $\sigma \in \bar{\mathscr{P}}^{\max}$ is strictly convex (on the induced subdivision of $\sigma$).
    \item[$(\dagger)$] The restriction $\alpha|_{\rho}$ of $\alpha$ to any $\rho \in \bar{\mathscr{P}}^{[1]}$ with $r_{\rho} \neq 0$ has a well-defined $\min(\alpha|_{\rho}) \subseteq \rho(\mathbb{Z})$.
\end{enumerate}

To find an $\alpha \in PA(B)$ satisfying $(*)$, we argue inductively similarly to Step 4 in the proof of Proposition \ref{ampledivisor}. 
$\left( \bar{B}, \bar{\mathscr{P}}\right)$ admits a PA-function satisfying $(*)$. Namely, we can just take $\bar{\alpha} := 0 \in PA(\bar{B})$. Suppose that an intermediate subdivision $\left(B'_1, \mathscr{P}'_1 \right)$ (obtained after pulling a sequence of integral points) of $\left( \bar{B}, \bar{\mathscr{P}}\right)$ admits a PA-function $\alpha'_1 \in PA\left(B'_1\right)$ satisfying $(*)$ and that we next pull an integral point $v \in B'(\mathbb{Z})$, obtaining a further subdivision $\left(B'_2, \mathscr{P}'_2 \right)$. We need to construct a PA-function $\alpha'_2 \in PA\left(B'_2\right)$ satisfying $(*)$. Let $\alpha_{\new}  \in PA\left(B_2'\right)$ be a PA-function with value $-1$ at $v$ and value $0$ at all the other vertices $v' \in \mathscr{P}_2'^{[0]}$. Then $k\alpha_{\new} \in PA\left(B_2'\right)$ is a well-defined integral PA-function for some $k \in \mathbb{N}$, and for every $\sigma \in \mathscr{P}'^{\max}_1$ that is subdivided by pulling $v$, $(k\alpha_{\new})|_{\sigma}$ is a strictly convex function on $\sigma$. Now for $M \gg k$, the PA-function $\alpha_2':= k\alpha_{\new} + M \alpha_1'$ on $\left(B'_2, \mathscr{P}'_2 \right)$ satisfies $(*)$. Repeating the argument over the whole sequence of pulls of integral points giving rise to $(B, \mathscr{P})$, we obtain an $\alpha \in PA(B)$ satisfying $(*)$. There are multiple ways to ensure that $\alpha \in PA(B)$ obtained by this procedure satisfies $(\dagger)$. The easiest is to just pull all the vertices $v \in \bar{\mathscr{P}}^{[0]}$ (in some order) first.
\end{proof}

\begin{observe} \label{singularitiesatthevertices}
    Similarly to Remark \ref{restrictedclassofresolutions}(1), one can always guarantee that the singularities of $\left(B, \mathscr{P}_{\coar}\right)$ of Construction \ref{plmapadmissible} lie at the vertices. To do that, it is enough to pull all the vertices $v \in \bar{\mathscr{P}}^{[0]}$ (in some order) first in the proof of Proposition \ref{strongadmissibleresolutionresult}. 
\end{observe}

We also have an analogue of Corollary \ref{furtherresolvetosatisfy}.

\begin{cor} \label{furtherresolutionofadmissible}
Let $\bar{\mathfrak{X}} \to \mathcal{S}$ be a special toric degeneration of K3-s. For any strongly admissible resolution $\pi: \mathfrak{X} \to \bar{\mathfrak{X}}$, there exists a further logarithmic modification $\pi': \tilde{\mathfrak{X}} \to \mathfrak{X}$ with the combined $\tilde{\mathfrak{X}} \to \mathfrak{X} \to \bar{\mathfrak{X}}$ a strongly admissible resolution with $\tilde{D}:=(\pi')^{-1}D$ simple normal crossings.
\end{cor}

\begin{proof}
Let $(B, \mathscr{P})$ be the dual intersection complex of $\mathfrak{X} \to \mathcal{S}$ and let $(\tilde{B}, \tilde{\mathscr{P}})$ be an integral subdivision of $(B, \mathscr{P})$ with all the cells $\sigma \in \tilde{\mathscr{P}}^{\max}$ standard triangles, obtained by pulling a sequence of integral points (such a subdivision exists by Lemma \ref{pullingresolution}). Recall that by \cite[Theorem 2.4.1 and Corollary 2.6.7]{AW} there is a one-to-one correspondence between logarithmic modifications of $\mathfrak{X} \to \mathcal{S}$ and subdivisions of $(B, \mathscr{P})$, and let $\pi': \tilde{\mathfrak{X}} \to \mathfrak{X}$ be the logarithmic modification corresponding to the subdivision $(\tilde{B}, \tilde{\mathscr{P}})$. Alternatively, consider $(\tilde{B}, \tilde{\mathscr{P}})$ as a subdivision of $\left(\bar{B}, \bar{\mathscr{P}}\right)$ and use the correspondence of Corollary \ref{onetoonetropical} to obtain the resolution $\tilde{\mathfrak{X}} \to \bar{\mathfrak{X}}$. Then $\tilde{\mathfrak{X}} \to \bar{\mathfrak{X}}$ factors through $\mathfrak{X}$ since $(\tilde{B}, \tilde{\mathscr{P}})$ is an integral subdivision of $(B, \mathscr{P})$. Note that $\tilde{D}$ is simple normal crossings since all the cells $\sigma \in \tilde{\mathscr{P}}^{\max}$ are standard triangles.

As in the proof of Proposition \ref{strongadmissibleresolutionresult}, it is enough to find an $\tilde{\alpha} \in PA(\tilde{B})$ that satisfies assumptions $(*)$ and $(\dagger)$ in that proof. By Proposition \ref{stronglyadmissibletosubdivision} and Corollary \ref{onetoonetropical}, we have an $\alpha \in PA(B)$ satisfying $(*)$ and $(\dagger)$. But then we can construct an $\tilde{\alpha} \in PA(\tilde{B})$ that satisfies assumptions $(*)$ and $(\dagger)$ by the same argument as in the proof of Proposition \ref{strongadmissibleresolutionresult}.
\end{proof}

\begin{remark} \label{satisfybutnotsncgeneral}
As in Remark \ref{satisfybutnotsnc}, note that one can often produce strongly admissible resolutions $\mathfrak{X} \to \mathcal{S}$ with $D$ not simple normal crossings using Proposition \ref{stronglyadmissibletosubdivision}, especially in the case that $r_{\rho} = 0$ for some $\rho \in \bar{\mathscr{P}}^{[1]}$.
\end{remark}

\begin{remark}
It is classical that any lattice polygon $P$ admits a subdivision that is \emph{unimodular} (i.e. all the maximal cells of the subdivision are standard triangles) and \emph{regular} (i.e. the subdivision supports a strictly convex PA-function), and such a subdivision can be obtained by pulling integral points, see, e.g. \cite[Lemma 2.1]{HPPS}.  
Lemma \ref{pullingresolution} (showing that there exists a $(B, \mathscr{P})$ with all the cells $\sigma \in \mathscr{P}^{\max}$ standard triangles) along with the argument in the proof of Proposition \ref{strongadmissibleresolutionresult} (giving an $\alpha \in PA(B)$ such that the restriction $\alpha|_{\sigma}$ of $\alpha$ to any $\sigma \in \bar{\mathscr{P}}^{\max}$ is strictly convex) provide a natural generalization of this fact to polyhedral manifolds of dimension $2$.
\end{remark}

\subsection[Scattering diagram interpretation]{Scattering diagram interpretation of the extended intrinsic mirror} \label{scatint}

For any special toric degeneration $\bar{\mathfrak{X}} \to \mathcal{S}$ of K3-s and a strongly admissible resolution $\pi:\mathfrak{X} \to \bar{\mathfrak{X}}$, we have defined the necessary data of Basic Setup \ref{setupforproof} in Section \ref{admissibleandstronglyso}. In particular, this gives us the toric degeneration mirror $\check{\bar{\mathfrak{X}}} \to \Spec \kk \llbracket t \rrbracket$, the extended intrinsic mirror $\check{\mathfrak{X}} \to \Spec \widehat{\kk [P]}_{J}$ (here $J := P \setminus K$ for $P$ a well-chosen monoid and $K \subseteq P$ the face containing the classes of curves contracted by $\pi$), and the basechange of the extended intrinsic mirror by $P \to \mathbb{N}, \ \beta \mapsto \pi^{*} A \cdot \beta$ makes sense.
The toric degeneration mirror is defined via a collection of algorithmic scattering diagrams $\bar{\mathfrak{D}}=\left\{\bar{\mathfrak{D}}_k, \ k \geq 0\right\}$ on $\left(\bar B, \bar{\mathscr{P}}\right)$. As promised in Remark \ref{scatforext}, we are going to give a collection of scattering diagrams $\mathfrak{D}^J: = \left\{ \mathfrak{D}_{J^{k+1}}, \ k \geq 0  \right\}$ (that we still call \emph{canonical scattering diagrams}) on $(B, \mathscr{P})$ giving rise to the extended intrinsic mirror. We actually won't use the existence of a $\pi$-ample PA-generated divisor $D'$ in Definition \ref{stronglyadmissible} (satisfying conditions (1) and (2)) of a strongly admissible resolution, see Remark \ref{stronglyadmissibleweakening}(1). Our analysis will generalize \cite[Section 5.3]{GHKS}, and we shall often refer to loc. cit. for the details. 

The setup of \cite{GHKS} is that one has a normal crossings degeneration $\mathcal{Y} \to \mathcal{S}$, a small contraction $\mathcal{Y} \to \bar{\mathcal{Y}}$, and a log étale blowup $\tilde{\mathcal{Y}} \to \mathcal{Y}$ such that $\tilde{\mathcal{Y}} \to \mathcal{S}$ is simple normal crossings. One then constructs the mirror to $\tilde{\mathcal{Y}} \to \mathcal{S}$ by defining, for every ideal $\tilde{I}$ with radical $\tilde{J} = \tilde{P} \setminus \tilde{K}$ (for, in our terminology, a well-chosen ideal $\tilde{P}$ with a face $\tilde{K}$) a scattering diagram $\mathfrak{D}_{\tilde{I}}$, consistent if $\mathcal{Y} \to \mathcal{S}$ is projective. 

The construction of \cite[Section 5.3]{GHKS} is quite general since it allows an arbitrary small contraction $\mathcal{Y} \to \bar{\mathcal{Y}}$. In particular, it covers our case of interest if the resolution $\pi: \mathfrak{X} \to \bar{\mathfrak{X}}$ factors as a composition of a small contraction and a log étale blowup. For a general $\pi: \mathfrak{X} \to \bar{\mathfrak{X}}$, the main difference between our setup and the setup of \cite{GHKS} is that $\left(B, \mathscr{P}_{\coar} \right)$ might have a singularity at the interior of an edge $\rho \in \mathscr{P}_{\coar}^{[1]}$, corresponding to the case that the exceptional curves $E_{\rho,k}, \ 1 \leq k \leq r_{\rho}$ of \eqref{intnumbers} are contained in some $D_{\rho, p}$ for $1 \leq p \leq l_{\rho}-1$. We shall need to treat this case separately. We will write $\Delta$ for the singular locus of the extended affine structure of \eqref{extsinglocus} on $(B, \mathscr{P})$.

\subsubsection{Canonical scattering modulo \texorpdfstring{$J$}{J}} \label{canonicalmodJ}

We first define a scattering diagram $\mathfrak{D}_{J}$ that is consistent modulo $J$ and such that $\check{\mathfrak{X}}_{\mathfrak{D}_{J}} \to \Spec \kk [P]/J$ is isomorphic to the family $\check{\mathfrak{X}}_{J} \to \Spec \kk [P]/J$ obtained by reducing the extended intrinsic mirror $\check{\mathfrak{X}} \to \Spec \widehat{\kk [P]}_{J}$  modulo $J$.

\begin{notation} \label{tropabusenotation}
Let $\btau$ be a decorated wall type. Following \cite[Section 5.3]{GHKS}, we will abuse the notation by writing $h_{\btau}:G \to B$ (or just $h: G \to B$) for the corresponding tropical map to $B$ (see Remark \ref{dicinstead}).
\end{notation}

\begin{remark}
Following \cite[Remark 5.20]{GHKS}, note that if $\btau$ is a decorated wall type (see Definition \ref{walltype}) and $\mathscr{M}(\mathfrak{X}, \btau)$ is non-empty, then the curve classes $\mathbf{A}(v)$ are determined by the underlying type $\tau$ provided $v$ is a vertex of $G$ with $h_{\btau}(v) \notin \Delta$. We say that a decorated wall type $\btau$ is \emph{well-decorated} if the decorations $\mathbf{A}(v)$ are as required by $\tau$ whenever $h_{\btau}(v) \notin \Delta$. See \cite[Remark 5.20]{GHKS} for details.
\end{remark}

First, we modify \cite[Definition 5.23]{GHKS} of a \emph{slab twig}. If there is a singularity $v_{\rho,p} \in \mathscr{P}^{[1]}$ of the affine structure in the interior of an edge $\tilde{\rho} \in \mathscr{P}_{\coar}^{[1]}$ (with $\mathscr{P}_{\coar}$ defined as in Construction \ref{plmapadmissible}) and $\rho$ is one of the connected components of $\tilde{\rho} \setminus v_{\rho,p}$, we say that $\rho$ is a \emph{half-edge of $\mathscr{P}_{\coar}$}. 

\begin{notation}
For $\rho$ a half-edge of $\mathscr{P}_{\coar}$ and $\tilde{\rho} \in \mathscr{P}_{\coar}^{[1]}$ the unique edge containing $\rho$, we will write $r_{\rho}$, $E_{\rho,k}$, and $F_{\rho,p}$ for $r_{\tilde{\rho}}$, $E_{\tilde{\rho},k}$, and $F_{\tilde{\rho},p}$ respectively. 
\end{notation}

Note that if there is a singularity $v_{\rho,p} \in \mathscr{P}^{[1]}$ of the affine structure in the interior of an edge $\rho \in \mathscr{P}_{\coar}^{[1]}$, then $D_{v_{\rho,p}}$ contains $(-1)$-curves with curve classes $E_{\rho,k}$ and $F_{\rho,p} - E_{\rho,k}$ for $1 \leq k \leq r_{\rho}$. If there is no such singularity, then we only have $(-1)$-curves with curve classes $E_{\rho,k}, \ 1 \leq k \leq r_{\rho}$ contained in $D_v$ for $v \in \mathscr{P}^{[0]}_{\coar}$ one of the endpoints of $\rho$.

\begin{defn} \label{slabtwig}
We say that a well-decorated wall type $\btau = (G, \bsigma, \mathbf{u}, \mathbf{A})$ is a \emph{slab twig with curve class $E_{\rho,k}$ (resp. either $E_{\rho,k}$ or $F_{\rho,p} - E_{\rho,k}$), weight $w$, and length $n \geq 0$} for 
\begin{enumerate}[(a)]
    \item $\rho$ an edge (resp. half-edge) of $\mathscr{P}_{\coar}$ such that $E_{\rho,k}$ (resp. either $E_{\rho,k}$ or $F_{\rho,p} - E_{\rho,k}$)\footnote{Here we confuse the curve classes with the underlying curves.} intersects $X_{\rho'}$ for a $\rho' \in \mathscr{P}^{[1]}$ with $\rho' \subseteq \rho$.
    \item $v \in \Delta \cap \mathscr{P}_{\coar}^{[0]}$ (resp. $v = v_{\rho,p} \in \Delta \setminus \mathscr{P}_{\coar}^{[0]}$) such that $E_{\rho,k} \subseteq D_{v}$ (resp. either $E_{\rho,k} \subseteq D_v$ or $F_{\rho,p} - E_{\rho,k} \subseteq D_{v}$).
\end{enumerate}
if:
\begin{enumerate}
    \item $G$ has vertices $v_0, \dots, v_n = v_{\out}$ and edges $E_i$ connecting $v_i$ to $v_{i+1}$ for $i \geq 0$, so that $v_0$ is univalent and all the other vertices are bivalent.
    \item $\bsigma(v_0) = \mathbf{C}v \in \mathbf{C}\Delta$ and $\textbf{A}(v_0) = wE_{\rho,k}$ (resp. either $\textbf{A}(v_0) = wE_{\rho,k}$ or $\textbf{A}(v_0) = w(F_{\rho,p} -E_{\rho,k})$).
    \item $h_{\btau}(v_i) \in \Int(\rho)$ for $1 \leq i \leq n$.
    \item For each $i$, $\mathbf{u}(E_i)$ is $w\nu_{\rho}$ where $\nu_{\rho}$ is a primitive tangent vector to $\rho$ pointing away from $v$, when $E_i$ is oriented from $v_i$ to $v_{i+1}$. The same is true for $\mathbf{u}(L_{\out})$.
\end{enumerate}
\end{defn}

\begin{remarks} \label{slabtwigremarks}
\begin{enumerate}
    \item As in \cite[Remark 5.24]{GHKS}, a slab twig is completely determined by its curve class, length, and weight. If $n > 0$, the vertices $v_1, \dots, v_n$ map to successive vertices of $\mathscr{P}^{[0]}$ along $\rho$ by rigidity. The curve classes $\mathbf{A}(v_i)$ for $i \geq 1$ are determined by the well-decorated condition. Note that the length of a slab twig for an edge or half-edge $\rho$ is bounded by $l_{\rho} - 1$ where $l_{\rho}$ is the length of $\rho$.
    \item Definition \ref{slabtwig} agrees with \cite[Definition 5.23]{GHKS} (up to different notations) if $\Delta \subseteq \mathscr{P}_{\coar}^{[0]}$. By Observation \ref{singularitiesatthevertices}, we can always construct a strongly admissible resolution $\pi: \mathfrak{X} \to \bar{\mathfrak{X}}$ satisfying this property. The results of \cite[Section 5.3]{GHKS} suffice for such resolutions. 
\end{enumerate}
\end{remarks}

We also need the following.

\begin{defn}[{\cite[Definition 5.25]{GHKS}}] \label{subwalltype}
Let $\btau = (G, \bsigma, \mathbf{u}, \mathbf{A})$ be a decorated wall type. We define a \emph{sub-wall type} of $\btau$ as a type $\btau'$ obtained as follows. There is an edge $E \in E(G)$ such that if we split $G$ at $E$, we obtain two connected components $G', G''$. Here, splitting at $E$ turns $E$ into a leg for both $G'$ and $G''$. Choose $G''$ so that $v_{\out} \in V(G'')$. Then $\btau'$ is obtained by restricting $\bsigma, \mathbf{u}$ and $\mathbf{A}$ to $G'$.\footnote{It follows from \cite[Lemma 8.12]{G3} that a sub-wall type is itself a decorated wall type.} We call $E$ the \emph{output edge} of the sub-wall type.
\end{defn}

Construction \ref{canonical scattering} describes the canonical scattering diagram modulo any monoid ideal $I \subseteq P$ with $P \setminus I$ finite. We compute the walls this construction gives for $\mathfrak{X} \to \mathcal{S}$ modulo $J$. Since $P \setminus J = K$ is not finite, it is reasonable to expect infinitely many non-trivial walls. 

\begin{fact} \label{diag}
The walls given by Construction \ref{canonical scattering} for $\mathfrak{X} \to \mathcal{S}$ modulo $J$ can be described as follows. Let $\rho \in \mathscr{P}_{\coar}^{[1]}$ be such that the corresponding edge $\rho \in \bar{\mathscr{P}}^{[1]}$ has $r_{\rho} \neq 0$. Suppose that $\rho = \langle v, v' \rangle$ for some $v,v' \in \mathscr{P}_{\coar}^{[0]}$ and let 
$$\rho_1 := \langle v, v_{\rho, 1} \rangle, \  \rho_2 := \langle v_{\rho, 1}, v_{\rho, 2} \rangle, \ \dots, \ \rho_{l_{\rho}}:= \langle v_{\rho, l_{\rho}-1}, v' \rangle$$ be the edges of $\mathscr{P}^{[1]}$ subdividing $\rho$. Suppose also that the exceptional curves $E_{\rho,k}, \ 1 \leq k \leq r_{\rho}$ are contained in $D_{\rho, p_0}$ for some $0 \leq p_0 \leq l_{\rho} -1$ (here we set $D_{\rho, 0}:= D_v, D_{\rho, l_{\rho}}:=D_{v'}$ and use similar notations for the corresponding vertices, as usual) and intersect $D_{\rho, p_0+1}$ at one point (this can be achieved by permuting $v$ and $v'$). Let $w_{\rho} = z^{m_{\rho}}$ where $m_{\rho} \in \Lambda_{\rho}$ is the integral generator pointing from $v'$ to $v$. Then the walls supported on every $\rho \in \mathscr{P}_{\coar}^{[1]}$ with $r_{\rho} \neq 0$ are
\begin{equation} \label{infmanywallfunctions}
\begin{aligned}
&\left\{\mathfrak{b}_{\rho_p,k,b} := \left( \rho_p, \exp\left(\frac{(-1)^{b-1}}{b} t^{b\left(F_{\rho,p} + \dots + F_{\rho,p_0} - E_{\rho, k}\right)} w_{\rho}^{-b}\right)   \right) \Big|\,
\substack{1 \leq p \leq p_0\\1 \leq k \leq r_{\rho}\\b \in \mathbb{Z}_{>0}}\right\}, \\
&\left\{\mathfrak{b}_{\rho_{p_0+1},k,b} := \left( \rho_{p_0+1}, \exp\left(\frac{(-1)^{b-1}}{b} t^{bE_{\rho, k}} w_{\rho}^{b}\right)   \right) \Big|\,
\substack{1 \leq k \leq r_{\rho}\\b \in \mathbb{Z}_{>0}}\right\}, \\
&\left\{\mathfrak{b}_{\rho_{p},k,b} := \left( \rho_p, \exp\left(\frac{(-1)^{b-1}}{b} t^{b\left(E_{\rho, k} + F_{\rho, p_0+1} + \dots + F_{\rho, p-1}\right)} w_{\rho}^{b}\right)   \right) \Big|\,
\substack{p_0 + 1 < p \leq l_{\rho}\\1 \leq k \leq r_{\rho}\\b \in \mathbb{Z}_{>0}}\right\},
\end{aligned}
\end{equation}
and all the other walls are trivial. 
\end{fact}

\begin{proof}
By birational invariance of punctured log Gromov-Witten invariants (see \cite[Theorem 1.4]{J}) and Corollary \ref{furtherresolutionofadmissible}, we may assume that $D$ is simple normal crossings (note that the map $\tilde{\mathfrak{X}} \to \mathfrak{X}$ of Corollary \ref{furtherresolutionofadmissible} does not resolve any curve classes appearing in the description of walls \eqref{infmanywallfunctions}). In particular, all the irreducible components of $\mathfrak{X}_0$ are Cartier, and we have well-defined intersection numbers $\beta \cdot D_i$ for $\beta \in P$ and $D_i, \ 1 \leq i \leq m$ an irreducible component of $\mathfrak{X}_0$. We use the notations of Construction \ref{canonical scattering}.

Let $\btau = (\tau, \mathbf{A})$ be a decorated wall type with total curve class $A \in P \setminus J = K$ and $W_{\btau} \neq 0$. We will show that $\btau$ is necessarily a slab twig. Let $h: G \to B$ be the corresponding tropical map (using Notation \ref{tropabusenotation}). We proceed inductively on the number of vertices of $G$. This is similar to the argument in Step I in the proof of \cite[Theorem 5.7]{GHKS}. 

The base case is that $G$ consists of a single vertex $v_{\out}$, necessarily with adjacent leg $L_{\out}$. By the balancing condition, we must have $v:=h(v_{\out}) \in \Delta$. There are two cases: either $v \in \Delta \cap \mathscr{P}_{\coar}^{[0]}$ or $v = v_{\rho,p_0} \in \Delta \setminus \mathscr{P}_{\coar}^{[0]}$. In the first case, $\btau$ is a slab twig by an argument as in Step I in the proof of \cite[Theorem 5.7]{GHKS} (note that case (1) in that argument does not occur, and in case (3), $C'$ is necessarily the unique curve of class $E_{\rho,k}$ for some $1 \leq k \leq r_{\rho}$). So we may assume $v = v_{\rho,p_0} \in \Delta \setminus \mathscr{P}_{\coar}^{[0]}$. Let $f: C^{o} \to \mathfrak{X}$ be a punctured curve of class $A \in K$ contributing to $W_{\btau}$. Then its image $C'$ is a genus $0$ connected curve contracted by $\pi: \mathfrak{X} \to \bar{\mathfrak{X}}$ and satisfying $C' \subseteq D_v$. Let $\rho_1, \rho_2$ be the two half-edges of $\mathscr{P}_{\coar}$ containing $v = v_{\rho,p_0}$ and let $S_{\rho_1}, S_{\rho_2}$ be the corresponding irreducible components of $\partial D_v$. Let $S_{\rho_1}, S_{1}, \dots, S_r, S_{\rho_2}, S'_1, \dots, S'_s$ be a cyclic ordering of the irreducible components of $\partial D_v$.\footnote{Our notations here are different from Step I in the proof of \cite[Theorem 5.7]{GHKS}.} Then $C'$ is one of the following:
\begin{enumerate}
    \item $C' \subseteq \bigcup_{i=1}^r S_i$ or $C' \subseteq \bigcup_{j=1}^s S'_j$ is a chain of curves contained in $\partial D_v$. 
    \item $C'$ is a fibre of $\pi|_{D_v}: D_v \to \mathbb{P}^1$ and $A$ is a multiple of the fibre class $F_{\rho, p}$.
    \item $C'$ is of the form $C' = C'_1 \cup C'_2$ where $C'_1$ is the unique curve of class $E_{\rho,k}$ and $C'_2$ is the unique curve of class $F_{\rho, p_0} - E_{\rho,k}$ for some $1 \leq k \leq r_{\rho}$.
    \item $C'$ is the unique curve of class $E_{\rho,k}$ for some $1 \leq k \leq r_{\rho}$.
    \item $C'$ is the unique curve of class $F_{\rho,p_0} - E_{\rho,k}$ for some $1 \leq k \leq r_{\rho}$.
\end{enumerate}
In cases (1-3), there exist (distinct) irreducible components $D_1$ and $D_2$ of $\mathfrak{X}_0$ such that $D_1 \cap C'$ and $D_2 \cap C'$ are one-point sets. As in case (2) of the argument in Step I in the proof of \cite[Theorem 5.7]{GHKS}, \cite[Corollary 1.14]{GS3} (or its analogue \cite[Proposition 2.1]{G3} for irreducible components of $C'$) implies that $f: C^{o} \to \mathfrak{X}$ must have at least two punctured points, mapping to $D_1 \cap C'$ and $D_2 \cap C'$. So these cases do not occur. Therefore, $C'$ is as in case (4), or $C'$ is as in case (5). In particular, $\btau$ is a slab twig.

Suppose our induction hypothesis is true when $G$ has $<n$ vertices, and assume $\btau$ is now given with $G$ having $n$ vertices. If $\btau_1, \dots, \btau_r$ are the sub-wall types adjacent to $v_{\out}$, these are all slab twigs by the induction hypothesis. We have two cases:
\begin{enumerate}
    \item $h(v_{\out})$ lies in the interior of a half-edge $\rho$.
    \item $h(v_{\out}) = v \in \Delta \cap \mathscr{P}_{\coar}^{[0]}$.
\end{enumerate}
Here it is crucial that our setup does not allow a slab twig $\btau'$ of Definition \ref{slabtwig} to \enquote{end} at $v_{\rho,p_0} \in \Delta \setminus \mathscr{P}_{\coar}^{[0]}$ (i.e. to have $h(\tau'_{\out})$ an edge of $\mathscr{P}^{[1]}$ containing $v_{\rho,p} \in \Delta \setminus \mathscr{P}_{\coar}^{[0]}$). But then $\btau$ is a slab twig by arguing in (1) and (2) exactly as in cases (1) and (2) of the induction argument in Step I in the proof of \cite[Theorem 5.7]{GHKS} respectively (one needs to adjust \cite[Construction 5.34]{GHKS} and \cite[Theorem 5.37]{GHKS} to our setup, see Steps 2 and 3 in the proof of Proposition \ref{finitenessresult} below).

It remains to compute the contributions of the slab twigs. It is clear that \cite[Lemma 5.36]{GHKS} holds in our setup (with all $\mu_j = 1$ since $\mathfrak{X}_0$ is reduced). We have the following cases:
\begin{enumerate}[(I)]
        \item \textbf{$\btau$ is a slab twig with curve class $E_{\rho,k}$ for some $1 \leq k \leq r_{\rho}$, weight $b \in \mathbb{Z}_{>0}$ and length $n=0$.} We have $k_{\tau} W_{\btau} = \frac{(-1)^{b-1}}{b}$ and $A = b E_{\rho,k}$ by \cite[Lemma 5.36]{GHKS}. Further, it is immediate from Definition \ref{slabtwig} of a slab twig that $h(\tau_{\out}) = \rho_{p_0+1}$ and $u_{\tau}:=\mathbf{u}(L_{\out}) = -b m_{\rho}$. We have $\ind(\btau) = 1$ by Observation \ref{nnmn}. By Construction \ref{canonical scattering}, the wall type $\btau$ defines a wall $$\mathfrak{b}_{\rho_{p_0+1},k, b} := \left( \rho_{p_0+1}, \exp\left(\frac{(-1)^{b-1}}{b} t^{bE_{\rho, k}} w_{\rho}^{b}\right)   \right).$$
        \item \textbf{$\btau$ is a slab twig with curve class $E_{\rho,k}$ for some $1 \leq k \leq r_{\rho}$, weight $b \in \mathbb{Z}_{>0}$ and length $1 \leq n \leq l_{\rho} - (p_0 + 1)$.} We have $k_{\tau} W_{\btau} = \frac{(-1)^{b-1}}{b}$ and $A = b\left(E_{\rho, k} + F_{\rho, p_0+1} + \dots + F_{\rho, p_0 + n} \right)$ by \cite[Lemma 5.36]{GHKS}. Further, it is immediate from Definition \ref{slabtwig} of a slab twig that $h(\tau_{\out}) = \rho_{p_0+n+ 1}$ and $u_{\tau}:=\mathbf{u}(L_{\out}) = -b m_{\rho}$. We have $\ind(\btau) = 1$ by Observation \ref{nnmn}. Set $p:=p_0+n+1$.
        By Construction \ref{canonical scattering}, the wall type $\btau$ defines a wall $$\mathfrak{b}_{\rho_{p},k,b} := \left( \rho_p, \exp\left(\frac{(-1)^{b-1}}{b} t^{b\left(E_{\rho, k} + F_{\rho, p_0+1} + \dots + F_{\rho, p-1}\right)} w_{\rho}^{b}\right)   \right).$$
        \item \textbf{$\btau$ is a slab twig with curve class $F_{\rho,p_0} - E_{\rho,k}$ for some $1 \leq k \leq r_{\rho}$, weight $b \in \mathbb{Z}_{>0}$ and length $n=0$.} Arguing as in (I), this case contributes a wall $$\mathfrak{b}_{\rho_{p_0},k,b} := \left( \rho_{p_0}, \exp\left(\frac{(-1)^{b-1}}{b} t^{b\left(F_{\rho,p_0} - E_{\rho, k}\right)} w_{\rho}^{-b}\right)\right).$$ 
        \item \textbf{$\btau$ is a slab twig with curve class $F_{\rho,p_0} - E_{\rho,k}$ for some $1 \leq k \leq r_{\rho}$, weight $b \in \mathbb{Z}_{>0}$ and length $1 \leq n \leq p_0 - 1$.} Arguing as in (II), this case contributes a wall $$\mathfrak{b}_{\rho_p,k,b} := \left( \rho_p, \exp\left(\frac{(-1)^{b-1}}{b} t^{b\left(F_{\rho,p} + \dots + F_{\rho,p_0} - E_{\rho, k}\right)} w_{\rho}^{-b}\right)   \right)$$ for $p:=p_0-n$.
    \end{enumerate}
The collection of walls of (I-IV), varied in all the relevant parameters, defines exactly the walls of \eqref{infmanywallfunctions}.
\end{proof}

We define $\mathfrak{D}_J$ by grouping the walls of Proposition \ref{diag} with the same support together. This procedure is similar to Remark \ref{dec}(2), but now we have countably many walls with the same support.
	
\begin{cons} \label{newwalls}
In the notations of Proposition \ref{diag}, we define $\mathfrak{D}_J$ as follows. Let $\rho \in \mathscr{P}_{\coar}^{[1]}$ with $r_{\rho} \neq 0$ and fix a $p$ with $1 \leq p \leq p_0$. We compute the infinite product of the wall functions over the walls $\left\{\mathfrak{b}_{\rho_p,k,b} \ | \ 1 \leq k \leq r_{\rho}, \ b \in \mathbb{Z}_{>0}\right\}$ supported on $\rho_p$:
\begin{multline*}
\prod_{1 \leq k \leq r_{\rho}, \ b \in \mathbb{Z}_{>0}} f_{\mathfrak{b}_{\rho_p,k,b}} = \prod_{1 \leq k \leq r_{\rho}, \ b \in \mathbb{Z}_{>0}} \exp\left( \frac{(-1)^{b-1}}{b} t^{b\left(F_{\rho,p} + \dots + F_{\rho, p_0} - E_{\rho, k}\right)} w_{\rho}^{-b}\right) = \\ = \prod_{1 \leq k \leq r_{\rho}} \exp\left(\sum_{b \in \mathbb{Z}_{>0}}  \frac{(-1)^{b-1}}{b} t^{b\left(F_{\rho,p} + \dots + F_{\rho,p_0} - E_{\rho, k}\right)} w_{\rho}^{-b}\right) = \\ = \prod_{1 \leq k \leq r_{\rho}} \exp \left(\log \left(1+t^{F_{\rho,p} + \dots + F_{\rho,p_0} - E_{\rho, k}} w_{\rho}^{-1}\right)\right) =  \prod_{1 \leq k \leq r_{\rho}} \left(1+t^{F_{\rho,p} + \dots + F_{\rho,p_0} - E_{\rho, k}} w_{\rho}^{-1}\right).
\end{multline*}
Note that the product is a polynomial function, so we can use it as a wall function. Instead of the infinitely many walls $\left\{\mathfrak{b}_{\rho_p,k,b} \ | \ 1 \leq k \leq r_{\rho}, \ b \in \mathbb{Z}_{>0}\right\}$, we define a single wall supported on $\rho_p$:
$$\mathfrak{b}_{\rho_p}:=\left(\rho_p, \prod_{1 \leq k \leq r_{\rho}} \left(1+t^{F_{\rho,p} + \dots + F_{\rho,p_0} - E_{\rho, k}} w_{\rho}^{-1}\right) \right).$$
Similarly, for $p = p_0+1$, the infinite product of the wall functions over the walls $\left\{\mathfrak{b}_{\rho_{p_0},k,b} \ | \ 1 \leq k \leq r_{\rho}, \ b \in \mathbb{Z}_{>0}\right\}$ supported on $\rho_{p_0}$ is:
$$\prod_{1 \leq k \leq r_{\rho}, \ b \in \mathbb{Z}_{>0}} f_{\mathfrak{b}_{\rho_{p_0},k,b}} = \prod_{1 \leq k \leq r_{\rho}} \left(1+t^{E_{\rho,k}} w_{\rho}\right)$$ and we define a single wall supported on $\rho_{p_0}$:
$$\mathfrak{b}_{\rho_{p_0}}:=\left(\rho_{p_0}, \prod_{1 \leq k \leq r_{\rho}} \left(1+t^{E_{\rho, k}} w_{\rho}\right) \right).$$
Finally, for $p$ with $p_0+1 \leq p \leq l_{\rho}$, the infinite product of the wall functions over the walls $\left\{\mathfrak{b}_{\rho_{p},k,b} \ | \ 1 \leq k \leq r_{\rho}, \ b \in \mathbb{Z}_{>0}\right\}$ supported on $\rho_{p}$ is:
$$\prod_{1 \leq k \leq r_{\rho}, \ b \in \mathbb{Z}_{>0}} f_{\mathfrak{b}_{\rho_{p},k,b}} = \prod_{1 \leq k \leq r_{\rho}} \left(1+t^{E_{\rho, k} + F_{\rho, p_0+1} + \dots + F_{\rho, p-1}} w_{\rho}\right)$$ and we define a single wall supported on $\rho_{p}$:
$$\mathfrak{b}_{\rho_{p}}:=\left(\rho_{p_0}, \prod_{1 \leq k \leq r_{\rho}} \left(1+t^{E_{\rho, k} + F_{\rho, p_0+1} + \dots + F_{\rho, p-1}} w_{\rho}\right) \right).$$
We set: 
$$\mathfrak{D}_J:= \left\{ \mathfrak{b}_{\rho_p} \ | \ \rho \in \mathscr{P}_{\coar}^{[1]} \text{ with } r_{\rho} \neq 0, \ 1 \leq p \leq l_{\rho} \right\}.$$ See Figure \ref{scattering} for an example of $\mathfrak{D}_J$.  
\end{cons}

Note that $\mathfrak{D}_J$ has finitely many walls, so it is a well-defined scattering diagram. We need to check that $\mathfrak{D}_J$ is consistent and recovers the extended intrinsic mirror $\check{\mathfrak{X}} \to \Spec \widehat{\kk [P]}_{J}$  modulo $J$. 

\begin{sloppypar}
\begin{fact} \label{consistencyandiso}
$\mathfrak{D}_J$ is a consistent scattering diagram and $\check{\mathfrak{X}}_{\mathfrak{D}_{J}} \to \Spec \kk [P]/J$ is isomorphic to the family $\check{\mathfrak{X}}_{J} \to \Spec \kk [P]/J$ obtained by reducing the extended intrinsic mirror $\check{\mathfrak{X}} \to \Spec \widehat{\kk [P]}_{J}$  modulo $J$.  
\end{fact}
\end{sloppypar}

\begin{proof}
Consistency of $\mathfrak{D}_J$ follows from the fact that its only walls are slabs. Alternatively, note that for any $l \geq 1$, the ideal $J + \mathfrak{m}^l$ has a finite complement in $P$ and Construction \ref{newwalls} implies that $\mathfrak{D}_J$ agrees with the scattering diagram $\mathfrak{D}_{J + \mathfrak{m}^l}$ of Construction \ref{canonical scattering} modulo $J + \mathfrak{m}^l$. Moreover, the ideals $J + \mathfrak{m}^l$ for $l \geq 1$ form an inverse system with limit $J$ so consistency of $\mathfrak{D}_J$ follows from consistency of all the $\mathfrak{D}_{J + \mathfrak{m}^l}$ for $l \geq 1$ (see Theorem \ref{canonicalconsistency} or \cite[Theorem 5.2]{GS4}). 

The families $\check{\mathfrak{X}}_{\mathfrak{D}_{J+ \mathfrak{m}^l}} \to \Spec \kk [P]/(J + \mathfrak{m}^l)$ for $l \geq 1$ form an inverse system with limit $\check{\mathfrak{X}}_{\mathfrak{D}_{J}} \to \Spec \kk [P]/J$, so to prove the second claim it is enough to check that $\check{\mathfrak{X}}_{\mathfrak{D}_{J+ \mathfrak{m}^l}} \to \Spec \kk [P]/(J + \mathfrak{m}^l)$ is isomorphic to the family $\check{\mathfrak{X}}_{J + \mathfrak{m}^l} \to \Spec \kk [P]/{\left(J + \mathfrak{m}^l\right)}$ obtained by reducing the extended intrinsic mirror $\check{\mathfrak{X}} \to \Spec \widehat{\kk [P]}_{J}$  modulo $J + \mathfrak{m}^l$. But $J + \mathfrak{m}^l$ is an ideal with a finite complement in $P$, so this follows by \cite[Theorem 6.1]{GS4}.
\end{proof}

\subsubsection{Canonical scattering modulo \texorpdfstring{$J^{k+1}, \ k \geq 0$}{Jk}} \label{modjkscatter}

By the same reasoning as in Section \ref{canonicalmodJ}, to define $\mathfrak{D}_{J^{k+1}}, \ k \geq 0$, we need to combine certain walls with the same support to get a finite scattering diagram. The idea is to group the wall types that only differ by slab twigs. We copy the precise definitions from \cite[Section 5.3]{GHKS} (with our modified definition of a slab twig and different notations).

\begin{defn}
Let $\btau$ be a well-decorated wall type. A \emph{slab twig of $\btau$} is a sub-wall type of $\btau$ (see Definition \ref{subwalltype}) which is:
\begin{enumerate}
    \item A slab twig in the sense of Definition \ref{slabtwig}.
    \item Maximal in the sense that it is not a sub-wall type of another sub-wall type of $\btau$, which is also a slab twig.
\end{enumerate}
\end{defn}

\begin{notation}
For a (decorated) slab type $\tau$ and $v \in V(G)$, we denote $\tilde{v}:=h_{\tau}(v)$. 
\end{notation}

\begin{defn} \label{slabequiv}
We say that two decorated wall types $\btau_1, \ \btau_2$ are \emph{slab equivalent} if the following holds. There exist subgraphs $G'_i \subseteq G_i$ obtained by deleting all slab twigs of $\btau_i$ and their output edges, and an isomorphism $\phi:G'_1 \to G'_2$ preserving $\bsigma, \mathbf{u}$, and $\mathbf{A}$. Further, if $v \in V(G'_1)$, $\tilde{v} \in \Delta$, $\rho$ is either a half-edge or an edge (if no half-edge exists) of $\mathscr{P}_{\coar}$ containing $\tilde{v}$, and $k$ is a fixed integer with $1 \leq k \leq r_{\rho}$, let $w_j^1, \ 1 \leq j \leq n^1$ be the weights of the slab twigs for $v$ and $\rho$ of curve class $E_{\rho, k}$ (or $F_{\rho,p} - E_{\rho,k}$ if $\tilde{v} = v_{\rho,p}$ and $E_{\rho,k}$ intersects $X_{\rho'}$ for $\rho' \in \mathscr{P}^{[1]}$ such that $\rho' \not \subseteq \rho$). Similarly, let $w_j^2, \ 1 \leq j \leq n^2$ be the weights of the slab twigs adjacent to $\phi(v)$ of the same curve class. Then we require $$\sum_{j=1}^{n^1} w_j^1 = \sum_{j=1}^{n^2} w_j^2.$$ Given a decorated wall type $\btau$, we write $[\btau]$ for its slab equivalence class.
\end{defn}

\begin{defn} \label{decslabtype}
\begin{sloppypar}
A \emph{decorated slab type} is a decorated wall type $\tilde{\btau} = (\tilde{G}, \tilde{\bsigma}, \tilde{\mathbf{u}}, \tilde{\mathbf{A}})$ such that, for each vertex $v \in V(\tilde{G})$ not contained in a slab twig of $\tilde{\btau}$, there is at most one slab twig adjacent to $v$ with any fixed curve class (either $E_{\rho,k}$ or $F_{\rho,p} - E_{\rho,k}$). 
\end{sloppypar}

If a decorated slab type is itself a slab twig, we call the decorated slab type \emph{trivial}.
\end{defn}

\begin{remark}
It is easy to see that for every decorated wall type $\btau$, there exists a unique slab type $\tilde{\btau}$ representing the slab equivalence class $[\btau]$. Conversely, given a decorated slab type $\tilde{\btau}$ it is easy to describe all the decorated wall types $\btau$ giving rise to $\tilde{\btau}$, see \cite[Construction 5.30]{GHKS}.
\end{remark}

We define $\mathfrak{D}_{J^{k+1}}, \ k \geq 0$ following \cite[Definition 5.31]{GHKS} (with the same remarks as in Remarks \ref{diagramandcone}, i.e. we define a scattering diagram on $B$ instead of $\mathbf{C}B$).

\begin{cons} \label{extendcanonicalscat}
For $\tilde{\btau}$ a decorated slab type, we write $[\tilde{\btau}]$ for the set of isomorphism classes of decorated wall types which give rise to $\tilde{\btau}$. For a non-trivial decorated slab type $\tilde{\btau}$, set 
\begin{equation} \label{virtslabtypecount}
W_{\tilde{\btau}} :=\sum_{\btau \in [\tilde{\btau}]} k_{\tau}W_{\btau}
\end{equation}
and 
$$\mathfrak{d}_{\btau}:=
\left(h(\tau_{\out}) \cap g_{\trop}^{-1}(1), \ \exp\left(W_{\tilde{\btau}}t^A z^{-\tilde{\mathbf{u}}(L_{\out})}\right)^{\ind(\btau)}\right)$$
where $A$ is the total class of any $\btau \in [\tilde{\btau}]$ and $\ind(\btau)$ is as in Definition \ref{index}. 

Let $k \geq 0$. We define
\begin{equation} \label{canonicalscatteringmodJk}
\mathfrak{D}_{J^{k+1}}:=\mathfrak{D}_J \cup \left\{\mathfrak{d}_{\tilde{\btau}}\,\Big|\,
\substack{\hbox{$\tilde{\btau}$ an isomorphism class of a non-trivial decorated slab}\\ \hbox{type with total curve class lying in 
$P\setminus J^{k+1}$, $W_{\tilde{\btau}} \neq 0$}}\right\}.
\end{equation}
Note that this is consistent with Construction \ref{newwalls} for $\mathfrak{D}_J$. Indeed, as shown in the proof of Proposition \ref{diag}, for a decorated wall type $\btau$ with total curve class $A \in P \setminus J = K$ one has $W_{\btau} = 0$ unless $\btau$ is a slab twig (i.e. a trivial slab type).
\end{cons}

\begin{fact} \label{finitenessresult}
For every $k \geq 0$, $\mathfrak{D}_{J^{k+1}}$ is a well-defined consistent scattering diagram. 
\end{fact}

\begin{proof}\textbf{Step 1.}
Since we allow scattering diagrams that don't satisfy conditions (2) and (3) of Definition \ref{scatter} (see Remark \ref{dec}), it is enough to show that $\mathfrak{D}_{J^{k+1}}$ is finite and consistent. Note that Construction \ref{extendcanonicalscat} implies that $\mathfrak{D}_{J^{k+1}}$ agrees with the scattering diagram $\mathfrak{D}_{J^{k+1} + \mathfrak{m}^l}$ of Construction \ref{canonical scattering} modulo $J^{k+1} + \mathfrak{m}^l$ for every $l \geq 1$. So, assuming that $\mathfrak{D}_{J^{k+1}}$ is finite, consistency follows as in the proof of Proposition \ref{consistencyandiso}, replacing $J$ with $J^{k+1}$. It is enough to prove that $\mathfrak{D}_{J^{k+1}}$ is finite. Moreover, by birational invariance of punctured log Gromov-Witten invariants (see \cite[Theorem 1.4]{J}) and Corollary \ref{furtherresolutionofadmissible}, we may assume that $D$ is simple normal crossings. 

\textbf{Step 2.} The crucial ingredients of the proof of the corresponding \cite[Theorem 5.7]{GHKS} are \cite[Construction 5.34]{GHKS} and \cite[Theorem 5.37]{GHKS}, which give an inductive way to compute the invariants $W_{\tilde{\btau}}$ of \eqref{virtslabtypecount}. It is easy to generalize them to our setting. We first explain how to modify \cite[Construction 5.34]{GHKS}.

Fix a non-trivial decorated slab type $\tilde{\btau} = (\tilde{G}, \tilde{\bsigma}, \tilde{\mathbf{u}}, \tilde{\mathbf{A}})$ with $W_{\tilde{\btau}} \neq 0$ and let $\tilde{h}: \tilde{G} \to B$ be the unique tropical map of type $\tilde{\btau}$. Let $v_{\out}$ be the vertex in $\tilde{G}$ attached to $L_{\out}$. 
Since $\tilde{\btau}$ is non-trivial, it is not a slab twig. Write the sub-wall types adjacent to $v_{\out}$ as $\tilde{\btau}'_1, \dots, \tilde{\btau}'_{p'}, \tilde{\btau}_1, \dots, \tilde{\btau}_{p}$ where the $\tilde{\btau}'_i$ are all slab twigs and the $\tilde{\btau}_j$ are not slab twigs. We assume that $\tilde{\btau}'_i$ has associated curve class $E_{\rho_i, k_i}$ (or $F_{\rho_i, p_i} - E_{\rho_i, k_i}$) and weight $w_i$. Let $E'_1, \dots, E'_{p'}, E_1, \dots, E_p$ be the edges connecting $v_{\out}$ to $\tilde{\btau}'_1, \dots, \tilde{\btau}'_{p'}, \tilde{\btau}_1, \dots, \tilde{\btau}_p$ respectively. Let $x:=\tilde{h}(v_{\out})$ and denote by $\sigma_x$ the minimal cell of $\mathscr{P}_{\coar}$ containing $x$. We use similar notations for a decorated wall type $\btau \in [\tilde{\btau}]$ with $\btau = (G, \bsigma, \mathbf{u}, \mathbf{A})$, denoting the adjacent sub-wall types to $v_{\out}$ by $\btau_1, \dots, \btau_{r}$ and the edges attaching them to $v_{\out}$ by $E_1, \dots, E_{r}$.  

\begin{sloppypar}
\cite[Construction 5.34]{GHKS} associates to $\tilde{\btau}$ and $\btau$ Looijenga pairs $(\hat{Y}_x, \hat{D}_x)$ and $(\hat{Y}'_x, \hat{D}'_x)$ respectively. Moreover, it defines decorated types $\tilde{\btau}_{\out} = (\tilde{G}_{\out}, \tilde{\bsigma}_{\out}, \tilde{\mathbf{u}}_{\out},  \hat{A}_x)$ (with $V(\tilde{G}_{\out}) = \{v_{\out}\}$, $E(\tilde{G}_{\out}) = \varnothing$,  $L(\tilde{G}_{\out}) = \{E_1, \dots, E_p, L_{\out} \}$) and $\btau_{\out} = (G_{\out}, \bsigma_{\out}, \mathbf{u}_{\out}, \hat{A}'_x)$ (with $V(G_{\out}) = \{v_{\out}\}$, $E(G_{\out}) = \varnothing$, $L(G_{\out}) = \{E_1, \dots, E_{r}, L_{\out} \}$) of log maps to $(\hat{Y}_x, \hat{D}_x)$ and $(\hat{Y}'_x, \hat{D}'_x)$ respectively. For $\btau$, the construction follows \cite[Construction 8.13]{G3}, and for $\tilde{\btau}$, one gets a modification that takes the slab twigs into account.
\end{sloppypar}

We define $\tilde{\tau}_{\out}$ and $\tau_{\out}$ by restricting the corresponding data for $\tilde{\tau}$ and $\tau$ respectively, as in \cite[Construction 5.34]{GHKS}. There are four cases to consider for the construction of $(\hat{Y}_x, \hat{D}_x)$ and $(\hat{Y}'_x, \hat{D}'_x)$, and the curve classes $\hat{A}_x \in A_1(\hat{Y}_x), \ \hat{A}'_x \in A_1(\hat{Y}'_x)$.

\textbf{Case I: $x \in B_0$ and there are no slab twigs adjacent to $v_{\out}$.} This is exactly as Case I of \cite[Construction 5.34]{GHKS}.

\textbf{Case II: $x \in B_0$ but there are slab twigs adjacent to $v_{\out}$.} Note that we have $\dim \sigma_x = 1$ and there exists a half-edge (or edge) $\rho_x$ with $x \subseteq \rho_x \subseteq \sigma_x$. Let $\rho, \rho' \in \mathscr{P}^{[1]}$ be the two edges contained in $\rho_x$ and adjacent to the endpoints of $\rho_x$. Note that either every slab twig $\tilde{\btau}_i, \ 1 \leq i \leq p'$ has an edge (or leg if it is length $0$) mapping to $\rho$ or every slab twig $\tilde{\btau}_i, \ 1 \leq i \leq p'$ has an edge (or leg if it is length $0$) mapping to $\rho'$. We proceed as in Case II of \cite[Construction 5.34]{GHKS} with the modified definitions of $\rho$ and $\rho'$ (Case II of \cite[Construction 5.34]{GHKS} corresponds to the situation that $\rho_x = \sigma_x$ is an edge). Note that Case II of \cite[Construction 5.34]{GHKS} is more general in the sense that it allows some slab twigs to have an edge mapping to $\rho$ and some to have an edge mapping to $\rho'$. Our construction of $\pi: \mathfrak{X} \to \bar{\mathfrak{X}}$ does not allow that to happen. 

\textbf{Case III: $x \in \Delta \cap \mathscr{P}_{\coar}^{[0]}$.} For every half-edge or edge $\rho$ of $\mathscr{P}_{\coar}$ with $x$ as an endpoint (here we only take edges $\rho$ with $\Int(\rho) \cap \Delta = \varnothing$), let $v_{\rho}$ be the other endpoint of $\rho$. Recall that the curve class associated to the twig $\tilde{\btau}'_i$ is $E_{\rho_i, k_i}$ (or $F_{\rho_i, p_i} - E_{\rho_i, k_i}$). Let $\tilde{v}_i$ be the image of the univalent vertex of $\tilde{\btau}'_i$ and note that $\tilde{v}_i = v_{\rho_i}$ as otherwise the edge $E'_i$ is contracted. Let $\tilde{\rho} \in \mathscr{P}^{[1]}$ be the edge contained in $\rho \in \mathscr{P}_{\coar}$ and containing $x$. We set $V_{\tilde{\rho}} := \left\{ \tilde{\btau}'_{i} \ | \ \rho = \rho_i\right\}$ and proceed as in Case III of \cite[Construction 5.34]{GHKS} with the modified definition of $V_{\tilde{\rho}}$ (again, Case III of \cite[Construction 5.34]{GHKS} corresponds to the situation that all the $\rho$ containing $x$ are edges).

\textbf{Case IV: $x = v_{\rho,p} \in \Delta \setminus \mathscr{P}_{\coar}^{[0]}$.} This case does not appear in \cite{GHKS}. Note that there are no slab twigs adjacent to $v_{\out}$ as otherwise the edges $E'_1, \dots, E'_{p'}$ are contracted. We can again proceed as in Case III of \cite[Construction 5.34]{GHKS}. 

Explicitly, following the notations of \cite[Construction 5.34]{GHKS}, let $Y_x$\footnote{We usually denote this by $D_x$.} be the irreducible component of $\mathfrak{X}_0$ corresponding to $x \in \mathscr{P}^{[0]}$ and let $(B_x, \Sigma_x)$ be the tropicalization of $Y_x$ (with the log structure coming from the divisor $D_x:=\partial Y_x$). Note that by construction and the assumption that $\tilde{\btau}$ is a wall type, for any adjacent edge or leg $E$ to $v_{\out} \in V(\tilde{G})$, $\tilde{h}(E)$ is one-dimensional and $\rho_E:=T_x(\tilde{h}(E))$ is a ray in $B_x$. Then let $\hat{\Sigma}_x$ be any refinement of $\Sigma_x$ such that:
\begin{enumerate}
    \item Every top-dimensional cone in $\hat{\Sigma}_x$ is a strictly convex rational polyhedral cone integral affine isomorphic to the standard cone $\mathbb{R}_{\geq 0}^2$.
    \item The rays $\rho_{E_i}, \ 1 \leq i \leq q$ and $\rho_{L_{\out}}$ are one-dimensional cones of $\hat{\Sigma}_x$.
\end{enumerate}

The refinement $\hat{\Sigma}_x$ of $\Sigma_x$ leads to a toric blowup $\pi': \bar{Y}_x \to Y_x$, see \cite[Lemma 1.6]{GHK}. This comes with the boundary divisor $\bar{D}_x$, which is the strict transform of $D_x$. For $L \in L(\tilde{G}_{\out})$, write $w_L$ for the index of $\tilde{\mathbf{u}}_{\out}(L)$, i.e. $\tilde{\mathbf{u}}_{\out}(L)$ is $w_L$ times a primitive tangent vector. This represents the order of tangency imposed by the contact order $\tilde{\mathbf{u}}_{\out}(L)$ with the corresponding divisor. Define $\bar{A}_x \in A_1(\bar{Y}_x)$ to be a class such that $\pi'_{*}(\bar{A}_x) = \tilde{\textbf{A}}(v_{\out}) \in A_1(Y_x)$ and, for $\tau$ a ray in $\hat{\Sigma}_x$, 
\begin{equation} \label{controlledbyintnumbers}
\bar{A}_x \cdot \bar{D}_{\tau} = \sum_{L \in L(\tilde{G}_{\out}), \ \tilde{\bsigma}_{\out}(L) = \tau} w_L.
\end{equation}

By \cite[Lemma 8.14]{G3}, $\bar{A}_x$ is unique if it exists. If $\bar{A}_x$ exists, we define $(\hat{Y}_x, \hat{D}_x) = (\hat{Y}'_x, \hat{D}'_x) := (\bar{Y}_x, \bar{D}_x)$ and $\hat{A}_x = \hat{A}'_x := \bar{A}_x$. This corresponds to the construction of Case III of \cite[Construction 5.34]{GHKS} since there are no slab twigs adjacent to $v_{\out}$.

\textbf{Step 3.} As in \cite[Construction 5.34]{GHKS}, the Construction of Step 2 gives rise to invariants $N_{\tilde{\btau}_{\out}}$ and $N_{\btau_{\out}}$ via \cite[Definition 5.32]{GHKS} or \cite[(8.3)]{G3}. In cases III and IV, if a curve class $\bar{A}_x$ satisfying \cite[(5.21)]{GHKS} (which reduces to \eqref{controlledbyintnumbers} in Case IV) does not exist, we define $N_{\tilde{\btau}_{\out}} = N_{\btau_{\out}} := 0$. There are some strong vanishing results on $N_{\tilde{\btau}_{\out}}$ and $N_{\btau_{\out}}$ (see \cite[Lemma 5.33]{GHKS}) that are crucial to the proof of finiteness.

\cite[Theorem 5.37]{GHKS} gives an inductive way to compute the invariants $W_{\tilde{\btau}}$ of \eqref{virtslabtypecount}. It is the analogue of \cite[Theorem 8.15]{G3}, which does the same for $\btau$. We claim that the theorem still works in our setup. Let $\tilde{\btau}$ be a non-trivial decorated slab type. Let $\tilde{w}_{\out}$ be the index of $\tilde{\mathbf{u}}_{\out}(L_{\out})$ in the tangent lattice to $\tilde{\bsigma}_{\out}(L_{\out})$ in the notations of Step 2. Then 
\begin{equation} \label{indvanishing}
\ind(\tilde{\btau})W_{\tilde{\btau}} = \frac{\tilde{w}_{\out}N_{\tilde{\btau}_{\out}}\prod_{i=1}^p (\ind(\tilde{\btau}_i)W_{\tilde{\btau}_i})}{|\Aut(\tilde{\btau}_1, \dots, \tilde{\btau}_p)|}
\end{equation}
where the group $\Aut(\tilde{\btau}_1, \dots, \tilde{\btau}_p)$ is the set of permutations $\sigma$ of $\{1, \dots, p\}$ such that $\tilde{\btau}_i$ and $\tilde{\btau}_{\sigma(i)}$ are isomorphic decorated slab types.
To check that \eqref{indvanishing} applies in our setup, we proceed as in the proof of \cite[Theorem 5.37]{GHKS} until the case-by-case analysis. Then our Cases I, II, and III follow exactly as Cases I, II, and III in the proof of \cite[Theorem 5.37]{GHKS} respectively. In Case IV, note that $x$ is integral point of $(B, \mathscr{P})$, so $\ind(\tilde{\btau}) = 1$ and $\ind(\tilde{\btau}_i) = 1$ for $1 \leq i \leq p$. Further, there are no slab twigs, so it is enough to check that $w_{\out} N_{\btau_{\out}} =\tilde{w}_{\out}N_{\tilde{\btau}_{\out}}$ where $w_{\out}$ is the index of $\mathbf{u}_{\out}(L_{\out})$ in the tangent lattice to $\bsigma_{\out}(L_{\out})$. Since $x$ is integral point of $(B, \mathscr{P})$, we have $w_{\out} = \tilde{w}_{\out}$. Further, $\hat{Y}_x = \hat{Y}'_x$, so $N_{\btau_{\out}} = N_{\tilde{\btau}_{\out}}$. Therefore, \eqref{indvanishing} applies in our setup.

\textbf{Step 4.} To prove that $\mathfrak{D}_{J^{k+1}}$ is finite, we will argue as in the proof of \cite[Theorem 5.7]{GHKS} but will need to consider Case IV separately. Denote the polarization on $\bar{\mathfrak{X}}$ by $H$ to keep with the notations of \cite[Theorem 5.7]{GHKS} (we usually denote it by $A$). There exists a $k \in \mathbb{N}$ such that if $A \in P$ and $A \cdot \pi^{*}(H) > k$, then $A \in J^{k+1}$ (here $\pi$ is the resolution $\pi: \mathfrak{X} \to \bar{\mathfrak{X}}$ as usual). Call $A \cdot \pi^{*}(H)$ the \emph{degree} of the curve class $A$. We then go by induction on degree, showing that there are at most a finite number of non-trivial decorated slab types $\tilde{\btau}$ with $W_{\tilde{\btau}} \neq 0$ and the total curve class $\tilde{A}$ of $\tilde{\btau}$ of degree at most $k$. 

Note that if the degree of $A \in P$ is $0$, then $A \in K$. As shown in the proof of Proposition \ref{diag}, for a decorated wall type $\btau$ with total curve class $A \in  K$ one has $W_{\btau} = 0$ unless $\btau$ is a slab twig (i.e. a trivial slab type). Since the $\mathfrak{D}_J$ of Construction \ref{newwalls} is finite, it suffices to prove the induction step.

Assume a finite number of non-trivial decorated slab types $\tilde{\btau}$ with $W_{\tilde{\btau}} \neq 0$ and total curve class of degree $<k$. We need to show that there are only a finite number of non-trivial decorated slab types $\tilde{\btau}$ of degree $k$ with $W_{\tilde{\btau}} \neq 0$. 

First, we say that $\tilde{\btau}$ has a \emph{tail} if there exists a sequence of vertices $v_0, \dots, v_n = v_{\out}$ of $\tilde{G}$ with $n \geq 1$ and edges $E_1, \dots, E_n$ with $E_i$ having endpoints $v_{i-1}$ and $v_i$, satisfying the following properties:
\begin{enumerate}
    \item $v_i$ is bivalent for $i \geq 1$.
    \item $\tilde{h}(v_i) \neq \Delta$ for $i \geq 1$ (unlike \cite{GHKS}, $\Delta$ includes singularities contained in the interiors of the edges $\rho \in \mathscr{P}_{\coar}^{[1]}$).
    \item There exists a cell $\sigma \in \mathscr{P}_{\coar}$ such that $\tilde{\bsigma}(E_i) \subseteq \mathbf{C}\sigma$ for each $i$ and $\tilde{\bsigma}(L_{\out}) \subseteq \mathbf{C}\sigma$.
\end{enumerate}
As in Step II of the proof of \cite[Theorem 5.7]{GHKS}, it follows that it is enough to consider slab types without a tail. Let $\tilde{\btau}$ be such a decorated slab type. We set $x := \tilde{h}(v_{\out})$ and let  $\sigma_x$ be the minimal cell of $\mathscr{P}_{\coar}$ containing $x$ as before. We have a case-by-case analysis.

\textbf{Case I: $\dim \sigma_x = 2$.} This corresponds to Case I in Step 2 above and follows exactly as in Case I of the proof of \cite[Theorem 5.7]{GHKS}.

\textbf{Case II: $\dim \sigma_x = 1$ and $x \in B_0$.} This corresponds to Case II in Step 2 above and follows exactly as in Case II of the proof of \cite[Theorem 5.7]{GHKS}.

\textbf{Case III: $\dim \sigma_x = 2$.} This corresponds to Case III in Step 2 above and follows exactly as in Case III of the proof of \cite[Theorem 5.7]{GHKS}.

\textbf{Case IV: $x = v_{\rho,p} \in \Delta \setminus \mathscr{P}_{\coar}^{[0]}$.} This corresponds to Case IV in Step 2 above. The argument will use ideas from both Case II and Case III of the proof of \cite[Theorem 5.7]{GHKS}. Let $\tilde{\btau}_1, \dots, \tilde{\btau}_p$ be the adjacent sub-wall types to $v_{\out}$ as before (recall that there are no slab twigs adjacent to $v_{\out}$ in this case). We need to show that there are a finite number of choices for $\tilde{\btau}_1, \dots, \tilde{\btau}_p$ and $\tilde{\mathbf{A}}(v_{\out})$. Note that we have $p \geq 1$ since otherwise $\tilde{\btau}$ is a trivial slab type. Note also that the total degree of each $\tilde{\btau}_i$ is positive (indeed, otherwise $\tilde{\btau}_i$ is a trivial slab type, a contradiction). 

If $p \geq 2$, then the degree of each $\tilde{\btau}_i$ is smaller than $k$ so there are a finite number of choices for $\tilde{\btau}_1, \dots, \tilde{\btau}_p$. If $p = 1$ and the degree of $\tilde{\btau}_1$ is $<k$, again there are a finite number of choices for $\tilde{\btau}_1$. So, in these cases, it suffices to show that there are a finite number of choices for $\tilde{\mathbf{A}}(v_{\out})$ (we will consider the case that $p = 1$ and the degree of $\tilde{\btau}_1$ is $k$ separately). To show that there are a finite number of choices for $\tilde{\mathbf{A}}(v_{\out})$ it is enough to prove that there are a finite number of choices for $\bar{A}_x$ of degree $\leq k$ (with degree measured with respect to the pullback of $H$ to $\bar{Y}_x$). 

Consider the contraction $g: \bar{Y}_x \to \mathbb{P}^1$ induced by $\pi': \bar{Y}_x \to Y_x$ and $\pi: \mathfrak{X} \to \bar{\mathfrak{X}}$, and note that $g^{*}H$ is a multiple of the fibre class $F_{\rho,p}$. We claim that there exist $l \in \mathbb{N}$ and $a_{\tau} \in \mathbb{N}$ for $\tau$ a ray in $\hat{\Sigma}_x$ such that $W := l g^{*}H + \sum_{\tau} a_{\tau} \bar{D}_{\tau}$ is ample. Let $\rho_1, \rho_2$ be the two half-edges of $\mathscr{P}_{\coar}$ containing $x = v_{\rho,p}$ and let $\tau_{\rho_1}, \tau_{\rho_2}$ be the corresponding rays of $\hat{\Sigma}_x$. Note that either $\bar{D}_{\tau_{\rho_1}} \cdot E_{\rho,k} > 0$ or $\bar{D}_{\tau_{\rho_2}} \cdot E_{\rho,k} > 0$ for all $1 \leq k \leq r_{\rho}$, and that $g^{*}H \cdot \bar{D}_{\tau_{\rho_1}} > 0, \ g^{*}H \cdot \bar{D}_{\tau_{\rho_2}} > 0$. Let $\tau_{\rho_1}, \tau_{1}, \dots, \tau_r, \tau_{\rho_2}, \tau'_1, \dots, \tau'_s$ be a cyclic ordering of the rays of $\hat{\Sigma}_x$. Choose $a_{\tau} \in \mathbb{N}$ so that $a_{\tau_{1}} \gg \dots \gg  a_{\tau_{r}} > 0$, $a_{\tau'_{1}} \gg \dots \gg  a_{\tau'_{s}} > 0$, $a_{\tau_{\rho_1}} = a_{\tau_{\rho_2}} \gg \max\{a_{\tau_{1}}, a_{\tau'_{1}}\}$, and choose $l \gg a_{\tau_{\rho_1}}$. It follows that $W = l g^{*}H + \sum_{\tau} a_{\tau} \bar{D}_{\tau}$ is ample. Moreover, $$W \cdot \bar{A}_x = l g^{*}H \cdot \bar{A}_x +   \left(\sum_{\tau} a_{\tau} \bar{D}_{\tau}\right) \cdot \bar{A}_x \leq lk + \sum_{\tau} a_{\tau} \left( \bar{D}_{\tau} \cdot \bar{A}_x \right)$$ is bounded for a fixed decorated slab type $\tilde{\btau}$ since $\tilde{\btau}$ controls every $\bar{D}_{\tau} \cdot \bar{A}_x$ via \eqref{controlledbyintnumbers}. By the fact that the Hilbert scheme is of finite type, there are at most a finite number of possibilities for $\bar{A}_x$, as desired.

Suppose now that $p = 1$ and the degree of $\tilde{\btau}_1$ is $k$. Necessarily, the degree of $\tilde{\mathbf{A}}(v_{\out})$ is $0$. A priori, we have $$\tilde{\mathbf{A}}(v_{\out}) = \sum_{1 \leq k \leq r_{\rho}} a_{\rho,k} E_{\rho,k} + \sum_{\tau} b_{\tau} \bar{D}_{\tau} +  b_{\rho,p} F_{\rho,p}$$ for some $a_{\rho,k} \in \mathbb{Z}$ and $ b_{\tau}, b_{\rho,p} \in \mathbb{N}$. However, if $b_{\tau} > 0$ for some ray $\tau$ in $\hat{\Sigma}_x$, the image of every stable log map in $\mathscr{M}(\bar{Y}_x, \tilde{\btau}_{\out})$ contains $\bar{D}_{\tau}$ and $N_{\tilde{\btau}_{\out}} = 0$ by \cite[Lemma 5.33(1)]{GHKS}. Hence $W_{\tilde{\btau}} = 0$ by \eqref{indvanishing}. If $a_{\rho,k} > 0$, the image of every stable log map in $\mathscr{M}(\bar{Y}_x, \tilde{\btau}_{\out})$ contains the unique curve $C'$ of class $E_{\rho, k}$ in $\bar{Y}_x$ (the preimage under $\pi': \bar{Y}_x \to Y_x$ of the unique curve of class $E_{\rho, k}$ in $Y_x$). Moreover, $C'$ meets $\bar{D}_x$ at one point. But then $N_{\tilde{\btau}_{\out}} = 0$ by \cite[Lemma 5.33(2)]{GHKS} and $W_{\tilde{\btau}} = 0$ by \eqref{indvanishing}. Similarly, if $a_{\rho,k} < 0$, the image of every stable log map in $\mathscr{M}(\bar{Y}_x, \tilde{\btau}_{\out})$ contains the unique curve $C'$ of class $F_{\rho,p} - E_{\rho, k}$ in $\bar{Y}_x$ and $W_{\tilde{\btau}} = 0$. 

\begin{sloppypar}
So $\tilde{\mathbf{A}}(v_{\out}) = b_{\rho,p} F_{\rho,p}$ is a multiple of the fibre class. This implies that $\tilde{\bsigma}_{\out} (L_{\out})$ and $\tilde{\bsigma}_{\out} (E_1)$ are the rays corresponding to $\langle v_{\rho, p},  v_{\rho, p - 1}\rangle$ and $\langle v_{\rho, p+1},  v_{\rho, p}\rangle$ respectively (in the notations of Proposition \ref{diag}). Further, $\tilde{\mathbf{u}}_{\out}(L_{\out}) = b_{\rho,p} m_{\rho}$ and $\tilde{\mathbf{u}}_{\out}(E_{1}) = -b_{\rho,p} m_{\rho}$. But then we have $\tilde{\bsigma}(L_{\out}) = \mathbf{C}\langle v_{\rho, p},  v_{\rho, p - 1}\rangle, \  \tilde{\bsigma}(E_1) = \mathbf{C}\langle v_{\rho, p+1},  v_{\rho, p }\rangle$ and $\tilde{\mathbf{u}}(L_{\out}) = b_{\rho,p} m_{\rho}, \tilde{\mathbf{u}}(E_{1}) = -b_{\rho,p} m_{\rho}$ by the construction of $\tilde{\btau}_{\out}$. This immediately implies that for every $\tilde{\btau}_1$ there exists a unique $\tilde{\btau}$ inducing $\tilde{\btau}_1$. Note that this is similar to the argument that it is enough to consider slab types without a tail in Step II of the proof of \cite[Theorem 5.7]{GHKS}. The point is that even though there is no notion of $\tilde{\btau}$ being well-decorated at $x \in \Delta$, $\tilde{\mathbf{A}}(v_{\out}) = b_{\rho,p} F_{\rho,p}$ is still controlled by the underlying wall type $\tilde{\tau}$ (since $b_{\rho,p} F_{\rho,p}$ is the pullback of a curve class in a toric variety). 
\end{sloppypar}

To conclude the proof, it is enough to show that there are a finite number of choices for $\tilde{\btau}_1$. Indeed, we have $x_1: = \tilde{h}_1(v_{1,\out}) = v_{\rho,p+1}$ (with obvious notations). So either $\dim \sigma_{x_1} = 1$ and $x_1 \in B_0$ (if $p+1 \neq l_{\rho}$) or $\dim \sigma_{x_1} = 2$  (if $p+1 = l_{\rho}$). But then there are a finite number of choices for $\tilde{\btau}_1$ by Case II or Case III respectively.  
\end{proof}

It follows from Construction \ref{extendcanonicalscat} that the scattering diagram $\mathfrak{D}_{J^{k+1}}$ is compatible with $\mathfrak{D}_{J^{k}}$ (in the sense of Definition \ref{compatibility}) for $k \geq 1$, so the families $\check{\mathfrak{X}}_{\mathfrak{D}_{J^{k+1}}} \to \Spec \kk[P]/J^{k+1}$ for $k \geq 0$ form an inverse system and taking the limit over this system gives 
\begin{equation} \label{extendedmirrorscatint}
\check{\mathfrak{X}}_{\mathfrak{D}^J} \to \Spec \widehat{\kk [P]}_J.
\end{equation}

Finally, we need to check that \eqref{extendedmirrorscatint} agrees with the extended intrinsic mirror.

\begin{fact} \label{samemirrors}
The family $\check{\mathfrak{X}}_{\mathfrak{D}^J} \to \Spec \widehat{\kk [P]}_J$ is isomorphic to the extended intrinsic mirror family $\check{\mathfrak{X}} \to \Spec \widehat{\kk [P]}_J$. 
\end{fact}

\begin{proof}
It is enough to show that for every $k \geq 0$, $\check{\mathfrak{X}}_{\mathfrak{D}_{J^{k+1}}} \to \Spec \kk [P]/J^{k+1}$ is isomorphic to the family $\check{\mathfrak{X}}_{J^{k+1}} \to \Spec \kk [P]/J^{k+1}$ obtained by reducing the extended intrinsic mirror $\check{\mathfrak{X}} \to \Spec \widehat{\kk [P]}_{J}$  modulo $J^{k+1}$. This follows as in the proof of Proposition \ref{consistencyandiso}, replacing $J$ with $J^{k+1}$ (and using the fact that for every $l \geq 1$, $\mathfrak{D}_{J^{k+1}}$ agrees with the scattering diagram $\mathfrak{D}_{J^{k+1} + \mathfrak{m}^l}$ of Construction \ref{canonical scattering} modulo $J^{k+1} + \mathfrak{m}^l$). 
\end{proof}

\begin{remarks} \label{stronglyadmissibleweakening}
The results of this section allow us to extend the setup of Conjecture \ref{conjecture}.
\begin{enumerate}
\item Note that we have not used the existence of a $\pi$-ample PA-generated divisor $D'$ in Definition \ref{stronglyadmissible} (satisfying conditions (1) and (2)) of a strongly admissible resolution in this section. So the family $\check{\mathfrak{X}}_{\mathfrak{D}^J} \to \Spec \widehat{\kk [P]}_J$ is well-defined for any toric, integral, and homogeneous $\pi: \mathfrak{X} \to \bar{\mathfrak{X}}$. Moreover, it follows from Proposition \ref{samemirrors} that the extended intrinsic mirror $\check{\mathfrak{X}} \to \Spec \widehat{\kk [P]}_J$ of \eqref{extendedmirror} is also well-defined and isomorphic to $\check{\mathfrak{X}}_{\mathfrak{D}^J} \to \Spec \widehat{\kk [P]}_J$. So the results of this section are stronger than the extension results of Section \ref{extresults} (at least for toric, integral, and homogeneous $\pi$).  As a consequence, we will actually prove Conjecture \ref{conjecture} for any toric, integral, and homogeneous resolution $\pi: \mathfrak{X} \to \bar{\mathfrak{X}}$.
\item One can also consider working with non-homogeneous resolutions of Remark \ref{otherresolutions}, i.e. allowing $E_{\rho,k},  \  1 \leq k \leq r_{\rho}$ (for a fixed $\rho \in \bar{\mathscr{P}}$) to be contained in different irreducible components of $\mathfrak{X}_0$. This requires a more general definition of a slab twig, and the argument in Case IV in the proof of Proposition \ref{finitenessresult} (as well as the induction step in the proof of Proposition \ref{diag}) is more involved. Since non-homogeneous resolutions are often non-projective, we can't appeal to Construction \ref{canonical scattering} to claim that $\mathfrak{D}_{J^{k+1}}$ is consistent in this case.
\end{enumerate}
\end{remarks}

\subsection[Relating the scattering diagrams]{Relating the canonical and algorithmic scattering diagrams} \label{main}

We are ready to prove Conjecture \ref{conjecture} for a special toric degeneration $\bar{\mathfrak{X}} \to \mathcal{S}$ of K3-s. Let $\pi:\mathfrak{X} \to \bar{\mathfrak{X}}$ be a strongly admissible resolution (or just a toric, integral, and homogeneous resolution, see Remark \ref{stronglyadmissibleweakening}(1)). We are going to relate the collection of scattering diagrams $\mathfrak{D}^J: = \left\{ \mathfrak{D}_{J^{k+1}}, \ k \geq 0  \right\}$ giving rise to the extended intrinsic mirror
$\check{\mathfrak{X}} \to \Spec \widehat{\kk [P]}_J$ and the collection of scattering diagrams $\bar{\mathfrak{D}} = \left\{ \bar{\mathfrak{D}}_k, \ k \geq 0 \right\}$ giving rise to the toric degeneration mirror $\check{\bar{\mathfrak{X}}} \to \Spec \kk \llbracket t \rrbracket$.

\subsubsection{Basechange of the canonical scattering diagram}
First, we interpret the basechange $h: P \to \mathbb{N}, \ \beta \mapsto \pi^{*}A \cdot \beta$ at the level of scattering diagrams.

\begin{cons} \label{ccc}
We define a collection of scattering diagrams $\mathfrak{D} = \left\{ \mathfrak{D}_k, \ k \geq 0 \right\}$ on $(B, \mathscr{P})$ with monoid $\mathbb{N}$ and $I_0 = \mathfrak{m} \subseteq \mathbb{N}$ as follows.

We let $\mathfrak{D}_0 := \mathfrak{D}_{(t)}$ have walls $$ \left(\mathfrak{b}_{\rho_p}, h\left(f_{\mathfrak{b}_{\rho_p}}\right) \right)$$ for $\mathfrak{b}_{\rho_p}$ a wall of $\mathfrak{D}_J$ of Construction \ref{newwalls}. We still write $\mathfrak{b}_{\rho_p}$ for these walls.
Note that we have 
\begin{equation} \label{functionsareok}
\begin{aligned}
h(\mathfrak{b}_{\rho_p}) &= \left( 1+ w_{\rho}^{-1}\right)^{r_{\rho}}, \quad  &1 \leq p \leq p_0 \\
h(\mathfrak{b}_{\rho_p}) &= \left( 1+ w_{\rho}\right)^{r_{\rho}}, \quad &p_0+1 \leq p \leq l_{\rho} 
\end{aligned}
\end{equation}
since $\pi^{*}A \cdot \beta = 0$ for $\beta \in K$.

More generally, we let $\mathfrak{D}_k := \mathfrak{D}_{(t^{k+1})}$ have walls $$(\mathfrak{p}, h(f_{\mathfrak{p}}) \mod (t^{k+1}))$$ for $\mathfrak{p}$ a wall of $\mathfrak{D}_{J^{k+1}}$ of Construction \ref{extendcanonicalscat}. Again, we still denote such a wall by $\mathfrak{p}$. Here the reduction modulo $(t^{k+1})$ is necessary to guarantee that $\mathfrak{D}_k$ is compatible with $\mathfrak{D}_{k-1}$ for $k \geq 1$. 

We define the MPA function $\varphi_{A}$ on $(B, \mathscr{P})$ via its kinks by setting $\kappa_{\rho} := X_{\rho} \cdot \pi^{*} A$. Note that this is compatible with the basechange.

Now $\mathfrak{D}_k$ is a consistent scattering diagram since consistency in codimensions $0$, $1$ and $2$ follows trivially from the corresponding consistency statement for $\mathfrak{D}_{J^{k+1}}$ by interpreting all the relevant monomials as elements of $\kk[t]$ via the basechange and reducing the equations modulo $(t^{k+1})$. By taking the inverse limit over $\check{\mathfrak{X}}_{\mathfrak{D}_k} \to \Spec \kk [t]/(t^{k+1})$ for $k \geq 0$ we define a family $\check{\mathfrak{X}}_{\mathfrak{D}} \to \Spec \kk \llbracket t \rrbracket$.
\end{cons}

\begin{fact} \label{pr}
The basechange of the extended intrinsic mirror family $\check{\mathfrak{X}} \to \Spec \widehat{\kk [P]}_J$ by $h: P \to \mathbb{N}, \ \beta \mapsto \pi^{*}A \cdot \beta$ is isomorphic to the family $\check{\mathfrak{X}}_{\mathfrak{D}} \to \Spec \kk \llbracket t \rrbracket$ of Construction \ref{ccc}.    
\end{fact}

\begin{proof}
By the construction of $\check{\mathfrak{X}}^o$ (see Section \ref{gluedfamily} and \cite[Proposition 2.4.1]{GHS}), the basechange of the family $\check{\mathfrak{X}}^o_{\mathfrak{D}_{J^{k+1}}}$ via $h$ is isomorphic to the family $\check{\mathfrak{X}}^o_{h\left(\mathfrak{D}_{J^{k+1}}\right)}$ with $h\left(\mathfrak{D}_{J^{k+1}}\right)$ a scattering diagram on $(B, \mathscr{P})$ with monoid $\mathbb{N}$, the MPA function $\varphi_A$ of Construction \ref{ccc}, and walls $(\mathfrak{p}, h(f_{\mathfrak{p}}))$ for $\mathfrak{p}$ a wall of $\mathfrak{D}_{J^{k+1}}$ (with no reduction modulo $(t^{k+1})$). Indeed, changing the scattering diagram from $\mathfrak{D}_{J^{k+1}}$ to $h\left(\mathfrak{D}_{J^{k+1}}\right)$ corresponds to interpreting the monomials in $R_{\mathfrak{u}}$ and $R_{\mathfrak{b}}$ as elements of $\kk [t]$ in the gluing construction of $\check{\mathfrak{X}}^{o}$.

Note that the scattering diagram $h(\mathfrak{D}_{J^{k+1}})$ is defined over the ideal $(t^m)$ for $$m := \max \left\{ \pi^{*}A \cdot \beta \ | \ \beta \in P \setminus J^{k+1} \right\}$$ but it agrees with $\mathfrak{D}_k$ modulo $(t^{k+1})$. So the family $\check{\mathfrak{X}}^o_{h\left(\mathfrak{D}_{J^{k+1}}\right)}$ agrees with $\check{\mathfrak{X}}^o_{\mathfrak{D}_k}$ modulo $(t^{k+1})$. But then as in the proof of Proposition \ref{help}, we see that also $\check{\mathfrak{X}}_{h\left(\mathfrak{D}_{J^{k+1}}\right)}$ agrees with $\check{\mathfrak{X}}_{\mathfrak{D}_k}$ modulo $(t^{k+1})$. This implies the result by taking the limit of both families over $(t^{k+1})$ for $k \geq 0$.
\end{proof}

\subsubsection{The image of \texorpdfstring{$\mathfrak{D}$}{D} under \texorpdfstring{$\Phi$}{Phi}}

To prove Conjecture \ref{conjecture}, we show that $(B, \mathfrak{D}_k)$ is equivalent to the algorithmically constructed $(\bar{B}, \bar{\mathfrak{D}}_{k})$ of Theorem \ref{reconstruction}. To do that, we will define a scattering diagram $\Phi (\mathfrak{D}_k)$ on $\left(\bar B, \bar{\mathscr{P}}\right)$ that is consistent and equivalent to $(B, \mathfrak{D}_k)$, and appeal to the uniqueness statement of Theorem \ref{reconstruction}. It will be more convenient to work with a scattering diagram on a refinement of $\left(\bar B, \bar{\mathscr{P}}\right)$.\footnote{If all the singularities of $\left(B, \mathscr{P}_{\coar}\right)$ are contained at the vertices, then one can argue in terms of $\Phi(\mathfrak{D}_k)$ directly. However, singularities at the edges make it more convenient to work with a refinement to check consistency at the joints corresponding to these singularities. Note from Observation \ref{singularitiesatthevertices}, that one can always construct a strongly admissible resolution $\pi:\mathfrak{X} \to \bar{\mathfrak{X}}$ of $\bar{\mathfrak{X}} \to \mathcal{S}$ such that the singularities of $\left(B, \mathscr{P}_{\coar}\right)$ are at the vertices.} 

\begin{notation} \label{nonsubdividing}
Let $(B, \mathscr{P})$ be a polyhedral manifold of dimension $2$ and let $\mathscr{P}'$ be a refinement of the polyhedral decomposition $\mathscr{P}$ on $B$. Then we denote by $\mathscr{P}'^{[0]} \setminus \mathscr{P}^{[0]}$ the set of vertices of $\mathscr{P}'^{[0]}$ that are not vertices of $\mathscr{P}^{[0]}$ and denote by $\mathscr{P}'^{[1]} \setminus \mathscr{P}^{[1]}$ the set of edges of $\mathscr{P}'^{[1]}$ that do not subdivide an edge of $\mathscr{P}^{[1]}$.
\end{notation}

\begin{cons} \label{refinement}
We let $\bar{\mathscr{P}}'$ be the obvious refined polyhedral decomposition on $\bar{B}$ such that $\left(B, \mathscr{P}\right) \cong \left(\bar B, \bar{\mathscr{P}}'\right)$ as polyhedral complexes. Note that $\left(\bar B, \bar{\mathscr{P}}'\right)$ is an affine manifold with singularities, and the singularities only lie at the interiors of the edges $\rho' \in \bar{\mathscr{P}}'^{[1]}$ (since we assumed that the singularities of $\left(\bar B, \bar{\mathscr{P}}\right)$ lie at irrational points of the edges). Moreover, by Remark \ref{indofsinglocus}, for every singularity $x_{\rho} \in \Int(\rho), \ \rho \in \bar{\mathscr{P}}^{[1]}$, we are free to choose a $\rho' \in \bar{\mathscr{P}}'^{[1]}, \ \rho' \subseteq \rho$ such that $x_{\rho} \in \Int(\rho')$. We will fix particular choices of such $\rho' \subseteq \rho$ in Construction \ref{chs}. The PL-isomorphism $\Phi: (B, \mathscr{P}) \to \left(\bar{B}, \bar{\mathscr{P}}\right)$ can be viewed as a map to $\left(\bar{B}, \bar{\mathscr{P}}'\right)$. This procedure is inverse to the coarsening $\mathscr{P}_{\coar}$ of $\mathscr{P}$ of Constructions \ref{plmap} and \ref{plmapadmissible}.

Define the MPA function $\varphi'_A$ on $\left(\bar{B}, \bar{\mathscr{P}}'\right)$ via its kinks by setting $\kappa_{\rho_p}(\varphi'_A) := \kappa_{\rho}(\varphi_A)$ for $\rho_p \in \bar{\mathscr{P}}'^{[1]}, \ 1 \leq p \leq l_{\rho}$ subdividing a $\rho \in \bar{\mathscr{P}}^{[1]}$ (using the notations of Construction \ref{newwalls}) and $\kappa_{\rho'}(\varphi'_A) := 0$ for $\rho' \in \bar{\mathscr{P}}'^{[1]} \setminus \bar{\mathscr{P}}^{[1]}$ (using Notation \ref{nonsubdividing}). We will typically denote the edges of $\bar{\mathscr{P}}'^{[1]}$ and the corresponding edges of $\mathscr{P}^{[1]}$ by $\rho'$ or $\rho_p$, and reserve the notation $\rho$ for edges of $\bar{\mathscr{P}}^{[1]}$ and $\mathscr{P}^{[1]}_{\coar}$.

Now, let $\bar{\mathfrak{D}}$ be a scattering diagram on $\left(\bar{B}, \bar{\mathscr{P}}\right)$. Then $\bar{\mathfrak{D}}$ defines a scattering diagram $\bar{\mathfrak{D}}'$ on $\left(\bar{B}, \bar{\mathscr{P}}'\right)$ by subdividing the walls of $\bar{\mathfrak{D}}$ (similarly to Remark \ref{dec}(1)). Note that since $\kappa_{\rho}(\varphi'_A) = 0$ for $\rho' \in \bar{\mathscr{P}}'^{[1]} \setminus \bar{\mathscr{P}}^{[1]}$, crossing a $\rho' \in \bar{\mathscr{P}}'^{[1]} \setminus \bar{\mathscr{P}}^{[1]}$ is equivalent (for gluing purposes) to crossing a trivial wall with support $\rho'$ and crossing a slab $\mathfrak{b} \subseteq \rho' \in \bar{\mathscr{P}}'^{[1]} \setminus \bar{\mathscr{P}}^{[1]}$ is equivalent to crossing a codimension $0$ wall with the same support and wall function. Moreover, since there are no singularities at the vertices $v \in \bar{\mathscr{P}}'^{[0]} \setminus \bar{\mathscr{P}}^{[0]}$, consistency for the newly-formed codimension $2$ joints can be reinterpreted similarly to consistency in codimension $0$, see the proof of Proposition \ref{consistencydiff}. So consistency of $\bar{\mathfrak{D}}$ is equivalent to consistency of $\bar{\mathfrak{D}}'$. Moreover, $\bar{\mathfrak{D}}'$ is equivalent to a scattering diagram on $\left(\bar{B}, \bar{\mathscr{P}}\right)$ obtained by adding trivial walls with support $\rho' \in \bar{\mathscr{P}}'^{[1]} \setminus \bar{\mathscr{P}}^{[1]}$ to $\bar{\mathfrak{D}}$ (assuming that $\bar{\mathfrak{D}}'$ is consistent). But then $\bar{\mathfrak{D}}'$ is equivalent to $\bar{\mathfrak{D}}$ by Remark \ref{combimpliesequiv}(3), and the two give rise to the same family by an argument as in Proposition \ref{help}. 
\end{cons}

We define a scattering diagram $\Phi (\mathfrak{D}_k)$ on $\left(\bar B, \bar{\mathscr{P}}\right)$ by defining $\Phi (\mathfrak{D}_k)'$ on $\left(\bar B, \bar{\mathscr{P}}'\right)$ first.

\begin{cons} \label{chs}
We define a scattering diagram $\Phi (\mathfrak{D}_k)'$ on $\left(\bar B, \bar{\mathscr{P}}'\right)$ as follows. For every codimension $0$ wall $\mathfrak{p} \in \mathfrak{D}_k$, we introduce a wall $$\Phi(\mathfrak{p}):= (\Phi(\mathfrak{p}), \Phi(f_{\mathfrak{p}})) \in \Phi(\mathfrak{D}_k)'$$ defined as in Construction \ref{plmapadmissible} via Construction \ref{plisosmall}. Similarly, for every slab $\mathfrak{b} \in \mathfrak{D}_k$ with $\mathfrak{b} \subseteq \rho' \in \mathscr{P}^{[1]} \setminus \mathscr{P}^{[1]}_{\coar}$ (using Notation \ref{nonsubdividing}), we introduce a slab 
\begin{equation} \label{irrslabimage}
\Phi(\mathfrak{b}):=(\Phi(\mathfrak{b}), \Phi(f_{\mathfrak{b}})) \in \Phi(\mathfrak{D}_k)'.
\end{equation}

By Remark \ref{indofsinglocus}, we may assume that for every $\rho \in \mathscr{P}_{\coar}^{[1]}$, the singularity $x_{\Phi(\rho)}$ of $\Phi(\rho) \in \bar{\mathscr{P}}^{[1]}$ is contained in $\Phi(\rho_{p_0+1}) \in \bar{\mathscr{P}}'^{[1]}$ (using the notations of Proposition \ref{diag}).\footnote{This is not a crucial requirement for the construction but allows for a less cumbersome definition for images of the slabs $\mathfrak{b} \in \mathfrak{D}_k, \ \mathfrak{b} \subseteq \rho_p$.} We have $\Phi(\rho) = \underline{\Phi(\rho)} \cup \underline{\Phi(\rho)}'$ for $\underline{\Phi(\rho)}, \underline{\Phi(\rho)}'\in \tilde{\bar{\mathscr{P}}}^{[1]}$ chosen as in Sections \ref{admissibleandstronglyso} and \ref{global42}. This choice, along with \eqref{choiceslab}, ensures that $\Phi(m_{\rho}) = m_{\Phi(\rho)}$ (for $\Phi(m_{\rho})$ defined as in Construction \ref{plmapadmissible} via Construction \ref{plisosmall}) and we have 
$$f_{\underline{\Phi(\rho)}} = \left(1+w_{\Phi(\rho)}\right)^{r_{\Phi(\rho)}}, \quad f_{\underline{\Phi(\rho)}'} = z^{m_{\underline{\Phi(\rho)}'\underline{\Phi(\rho)} }}f_{\underline{\Phi(\rho)}} = \left(1+w_{\Phi(\rho)}^{-1}\right)^{r_{\Phi(\rho)}}.$$ Moreover, we have $\Phi(\rho_p) \subseteq \underline{\Phi(\rho)}'$ for $1 \leq p \leq p_0$ and $\Phi(\rho_p) \subseteq \underline{\Phi(\rho)}$ for $p_0+2 \leq p \leq l_{\rho}$ (since for every $\rho \in \mathscr{P}_{\coar}^{[1]}$, we required that the singularity $x_{\Phi(\rho)}$ of $\Phi(\rho) \in \bar{\mathscr{P}}^{[1]}$ is contained in $\Phi(\rho_{p_0+1}) \in \bar{\mathscr{P}}'^{[1]}$).

Now, suppose that $\mathfrak{b} \in \mathfrak{D}_k$ is a slab with $\mathfrak{b} \subseteq \rho_{p} \in \mathscr{P}^{[1]}$. 
There are three cases: 
\begin{enumerate}
    \item Either $\Phi(\mathfrak{b}) \subseteq \underline{\Phi(\rho)}$, or $\Phi(\mathfrak{b}) \subseteq \underline{\Phi(\rho)}'$ and $1 \leq p \leq p_0$. We introduce a slab $$\Phi(\mathfrak{b}):=(\Phi(\mathfrak{b}), \Phi(f_{\mathfrak{b}})) \in \Phi(\mathfrak{D}_k)'$$ defined as in Construction \ref{plmapadmissible} via Construction \ref{plisosmall}. 
    \item $\Phi(\mathfrak{b}) \subseteq \underline{\Phi(\rho)}'$ and $p=p_0+1$. We introduce a slab $$\Phi(\mathfrak{b})' :=(\Phi(\mathfrak{b}), z^{m_{\underline{\Phi(\rho)}'\underline{\Phi(\rho)} }}f_{\Phi(\mathfrak{b})}) \in \Phi(\mathfrak{D}_k)'.$$
    \item $x_{\Phi(\rho)} \in \Int(\Phi(\mathfrak{b}))$. Then $\Phi(\mathfrak{b}) \setminus x_{\Phi(\rho)} = \bar{\Phi}(\mathfrak{b}) \cup \bar{\Phi}(\mathfrak{b})'$ for some $$\bar{\Phi}(\mathfrak{b}) \subseteq \underline{\Phi(\rho)} \cap \Phi(\rho_{p_0+1}), \quad \bar{\Phi}(\mathfrak{b})'  \subseteq \underline{\Phi(\rho)}' \cap \Phi(\rho_{p_0+1}).$$ We introduce two slabs as follows\footnote{Here we confuse $\bar{\Phi}(\mathfrak{b})$ and $\bar{\Phi}(\mathfrak{b})'$ with their closures.}:
    $$
    \begin{aligned}
    \bar{\Phi}(\mathfrak{b}) &:=(\bar{\Phi}(\mathfrak{b}), f_{\Phi(\mathfrak{b})}) \in \Phi(\mathfrak{D}_k)' \\
    \bar{\Phi}(\mathfrak{b})' &:=(\bar{\Phi}(\mathfrak{b})', z^{m_{\underline{\Phi(\rho)}'\underline{\Phi(\rho)} }}f_{\Phi(\mathfrak{b})}) \in \Phi(\mathfrak{D}_k)'
    \end{aligned}
    $$
\end{enumerate}
This defines $\Phi (\mathfrak{D}_k)'$. We define the MPA function $\bar{\varphi}'_A$ on $\left(\bar{B}, \bar{\mathscr{P}}'\right)$ as in Construction \ref{refinement}. Note that for any $\rho' \in \bar{\mathscr{P}}'^{[1]} \setminus \bar{\mathscr{P}^{[1]}}$ we have $\kappa_{\rho'} = \pi^{*} A \cdot X_{\rho'} = 0$ since $X_{\rho'}$ is contracted by $\pi:\mathfrak{X} \to \bar{\mathfrak{X}}$. It follows that $\kappa_{\rho'} (\varphi_A) = \kappa_{\Phi(\rho')} (\bar{\varphi}_A')$ for any $\rho' \in \bar{\mathscr{P}}'^{[1]}$, so the construction is compatible with the MPA functions. 

We call slabs of the form $\Phi(\mathfrak{b})$ or $\bar{\Phi}(\mathfrak{b})$ \emph{slabs of the first type}. We call slabs of the form $\Phi(\mathfrak{b})'$ or $\bar{\Phi}(\mathfrak{b})'$ \emph{slabs of the second type}. We will often drop the bar in the notation for $\bar{\Phi}(\mathfrak{b})$ and $\bar{\Phi}(\mathfrak{b})'$.

Let $\Phi(\mathfrak{D}_k)$ be the scattering diagram on $\left(\bar B, \bar{\mathscr{P}}\right)$ that is obtained from $\Phi(\mathfrak{D}_k)'$ by replacing the edges $\rho' \in \bar{\mathscr{P}}'^{[1]} \setminus \bar{\mathscr{P}}^{[1]}$ with trivial walls with support $\rho'$ and replacing the slabs $\mathfrak{b} \subseteq \rho' \in \bar{\mathscr{P}}'^{[1]} \setminus \bar{\mathscr{P}}^{[1]}$ with codimension $0$ walls with the same support and wall functions.
\end{cons}

\begin{obs} \label{obsplimage}
It is immediate from Construction \ref{chs} that:
\begin{enumerate}
    \item $\Phi(\mathfrak{D}_k)'$ arises from $\Phi(\mathfrak{D}_k)$ by the procedure of Construction \ref{refinement} (here it is crucial that the walls of $\mathfrak{D}_0$ satisfy \eqref{functionsareok}). So if $\Phi(\mathfrak{D}_k)'$ is consistent, then so is $\Phi(\mathfrak{D}_k)$. Moreover, the two are equivalent and give rise to the same family. 
    \begin{sloppypar}\item $\Phi(\mathfrak{D}_k)'$ is compatible with $\Phi(\mathfrak{D}_{k-1})'$ and $\Phi(\mathfrak{D}_k)$ is compatible with $\Phi(\mathfrak{D}_{k-1})$ for $k \geq 1$ (since the same is true for $\mathfrak{D}_k$).\end{sloppypar}
    \item $\Phi(\mathfrak{D}_0)$ is combinatorially equivalent to $\bar{\mathfrak{D}}_0$. This is the main motivation for Construction \ref{chs} and will be crucial when we appeal to the uniqueness statement of Theorem \ref{reconstruction}.
\end{enumerate} 
\end{obs}

In Figure \ref{scattering} on the next page, we give the scattering diagrams $\mathfrak{D}_J$, $\mathfrak{D}_{0}$, $\Phi(\mathfrak{D}_{0})'$ (to obtain $\Phi(\mathfrak{D}_{0})$, remove the vertex $\bar{v}_1^{15}$ and replace $\langle v_1^{15}, v_4 \rangle$ and $\langle v_1^{15}, v_1^{16} \rangle$ with trivial walls), and $\bar{\mathfrak{D}}_0$ for $\pi: \mathfrak{X} \to \bar{\mathfrak{X}}$ of Figure \ref{blowingadivisor2} in the neighbourhood of the edge $\rho = \langle v_1, v_5 \rangle$ (or $\langle \bar{v}_1, \bar{v}_5 \rangle$ for the diagrams on $\left(\bar{B}, \bar{\mathscr{P}}\right)$ and $\left(\bar{B}, \bar{\mathscr{P}}'\right)$). Note that we have $\rho_{p_0+1} = \langle v_1, v_1^{15} \rangle$. 

\newpage

\begin{figure}[H]
\begin{center}

\begin{tikzpicture}[line cap=round,line join=round,>=triangle 45,x=1cm,y=1cm]

\begin{scope}[yshift = 15cm, scale = 2]
\clip(-3,-1.8) rectangle (2,1.8);
\draw  (-1,0)-- (1,0);
\draw  (1,0)-- (2,1);
\draw  (1,0)-- (2,-1);
\draw  (-1,0)-- (-2,1);
\draw  (-1,0)-- (-2,-1);
\draw  (0,0)-- (0.1,0.8);
\draw  (0,0)-- (0.1,-0.8);

\begin{scriptsize}

\draw [fill] (0,0) circle (1pt);
\draw [fill] (1,0) circle (1pt);
\draw [fill] (-1,0) circle (1pt);

\draw[latex-,dashed] 
        (-0.5,0) -- (-0.5, -1)
         node[fill=white]{\tiny $(1+t^{E_1^{15}}w_{\rho})  (1+t^{E_2^{15}}w_{\rho})$};
         
\draw[latex-,dashed] 
        (0.5,0) -- (0.5, 1)
         node[fill=white]{\tiny $(1+t^{E_1^{15} + F_1^{15}}w_{\rho})  (1+t^{E_2^{15} + F_1^{15}}w_{\rho})$};
         
\draw(-1.2,0) node{\normalsize $v_1$};
\draw(1.2,0) node{\normalsize $v_5$};
\draw(0.2,0.15) node{\normalsize $v_1^{15}$};
\draw(-0.3,0.5) node{\tiny $\langle v_1^{15}, v_4 \rangle$};
\draw(0.45,-0.5) node{\tiny $\langle v_1^{15}, v_1^{16} \rangle$};

\draw [color=blue, line width = 1pt] (-1, 0)-- ++(-2.25pt,-2.25pt) -- ++(4.5pt,4.5pt) ++(-4.5pt,0) -- ++(4.5pt,-4.5pt);
\draw [color=blue, line width = 1pt] (1,0)-- ++(-2.25pt,-2.25pt) -- ++(4.5pt,4.5pt) ++(-4.5pt,0) -- ++(4.5pt,-4.5pt);

\draw(-2.5,0) node{\normalsize $\mathfrak{D}_J$};

\end{scriptsize}
\end{scope}

\begin{scope}[yshift = 10cm, scale = 2]
\clip(-3,-1.8) rectangle (2,1.8);
\draw  (-1,0)-- (1,0);
\draw  (1,0)-- (2,1);
\draw  (1,0)-- (2,-1);
\draw  (-1,0)-- (-2,1);
\draw  (-1,0)-- (-2,-1);
\draw  (0,0)-- (0.1,0.8);
\draw  (0,0)-- (0.1,-0.8);

\draw [fill] (0,0) circle (1pt);

\begin{scriptsize}

\draw [fill] (1,0) circle (1pt);
\draw [fill] (-1,0) circle (1pt);

\draw(-1.2,0) node{\normalsize $v_1$};
\draw(1.2,0) node{\normalsize $v_5$};  

\draw(-0.3,0.5) node{\tiny $\langle v_1^{15}, v_4 \rangle$};
\draw(0.45,-0.5) node{\tiny $\langle v_1^{15}, v_1^{16} \rangle$};
\draw(0.2,0.15) node{\normalsize $v_1^{15}$};

\draw(0.5,-0.15) node{\tiny $(1+w_{\rho})^{2}$};
\draw(-0.5,-0.15) node{\tiny $(1+w_{\rho})^{2}$};

\draw [color=blue, line width = 1pt] (-1, 0)-- ++(-2.25pt,-2.25pt) -- ++(4.5pt,4.5pt) ++(-4.5pt,0) -- ++(4.5pt,-4.5pt);
\draw [color=blue, line width = 1pt] (1,0)-- ++(-2.25pt,-2.25pt) -- ++(4.5pt,4.5pt) ++(-4.5pt,0) -- ++(4.5pt,-4.5pt);

\draw(-2.5,0) node{\normalsize $\mathfrak{D}_0$};

\end{scriptsize}
\end{scope}

\begin{scope}[yshift = 5cm, scale = 2]
\clip(-3,-1.8) rectangle (2,1.8);
\draw  (-1,0)-- (1,0);
\draw  (1,0)-- (2,1);
\draw  (1,0)-- (2,-1);
\draw  (-1,0)-- (-2,1);
\draw  (-1,0)-- (-2,-1);
\draw  (0,0)-- (0.1,0.8);
\draw  (0,0)-- (0.1,-0.8);

\begin{scriptsize}

\draw [fill] (0,0) circle (1pt);
\draw [fill] (1,0) circle (1pt);
\draw [fill] (-1,0) circle (1pt);

\draw[latex-,dashed] 
        (-0.75,0) -- (-0.75, -1)
         node[fill=white]{\tiny $(1+w_{\rho}^{-1})^2$};
         
\draw[latex-,dashed] 
        (-0.25,0) -- (-0.25, 1)
         node[fill=white]{\tiny $(1+w_{\rho})^2$};
         
\draw(0.5,-0.15) node{\tiny $(1+w_{\rho})^{2}$};

\draw(-1.2,0) node{\normalsize $\bar{v}_1$};
\draw(1.2,0) node{\normalsize $\bar{v}_5$};
\draw(0.45,0.5) node{\tiny $\langle v_1^{15}, v_4 \rangle$};
\draw(0.45,-0.5) node{\tiny $\langle v_1^{15}, v_1^{16} \rangle$};
\draw(0.2,0.15) node{\normalsize $\bar{v}_1^{15}$};

\draw [color=blue, line width = 1pt] (-0.5, 0)-- ++(-2.25pt,-2.25pt) -- ++(4.5pt,4.5pt) ++(-4.5pt,0) -- ++(4.5pt,-4.5pt);

\draw(-2.5,0) node{\normalsize $\Phi(\mathfrak{D}_0)'$};

\end{scriptsize}
\end{scope}

\begin{scope}[scale = 2]
\clip(-3,-1.8) rectangle (2,1.8);
\draw  (-1,0)-- (1,0);
\draw  (1,0)-- (2,1);
\draw  (1,0)-- (2,-1);
\draw  (-1,0)-- (-2,1);
\draw  (-1,0)-- (-2,-1);

\begin{scriptsize}

\draw [fill] (1,0) circle (1pt);
\draw [fill] (-1,0) circle (1pt);

\draw[latex-,dashed] 
        (-0.75,0) -- (-0.75, -1)
         node[fill=white]{\tiny $(1+w_{\rho}^{-1})^2$};
         
\draw(0.25,-0.15) node{\tiny $(1+w_{\rho})^{2}$};

\draw(-1.2,0) node{\normalsize $\bar{v}_1$};
\draw(1.2,0) node{\normalsize $\bar{v}_5$};

\draw [color=blue, line width = 1pt] (-0.5, 0)-- ++(-2.25pt,-2.25pt) -- ++(4.5pt,4.5pt) ++(-4.5pt,0) -- ++(4.5pt,-4.5pt);

\draw(-2.5,0) node{\normalsize $\bar{\mathfrak{D}}_0$};

\end{scriptsize}
\end{scope}

\end{tikzpicture}

\end{center}
\caption{$\mathfrak{D}_J, \mathfrak{D}_{0}, \Phi(\mathfrak{D}_{0})'$, and $\bar{\mathfrak{D}}_0$ in the neighbourhood of $\langle v_1, v_5 \rangle$ of Figure \ref{blowingadivisor2}.} \label{scattering}
\end{figure}
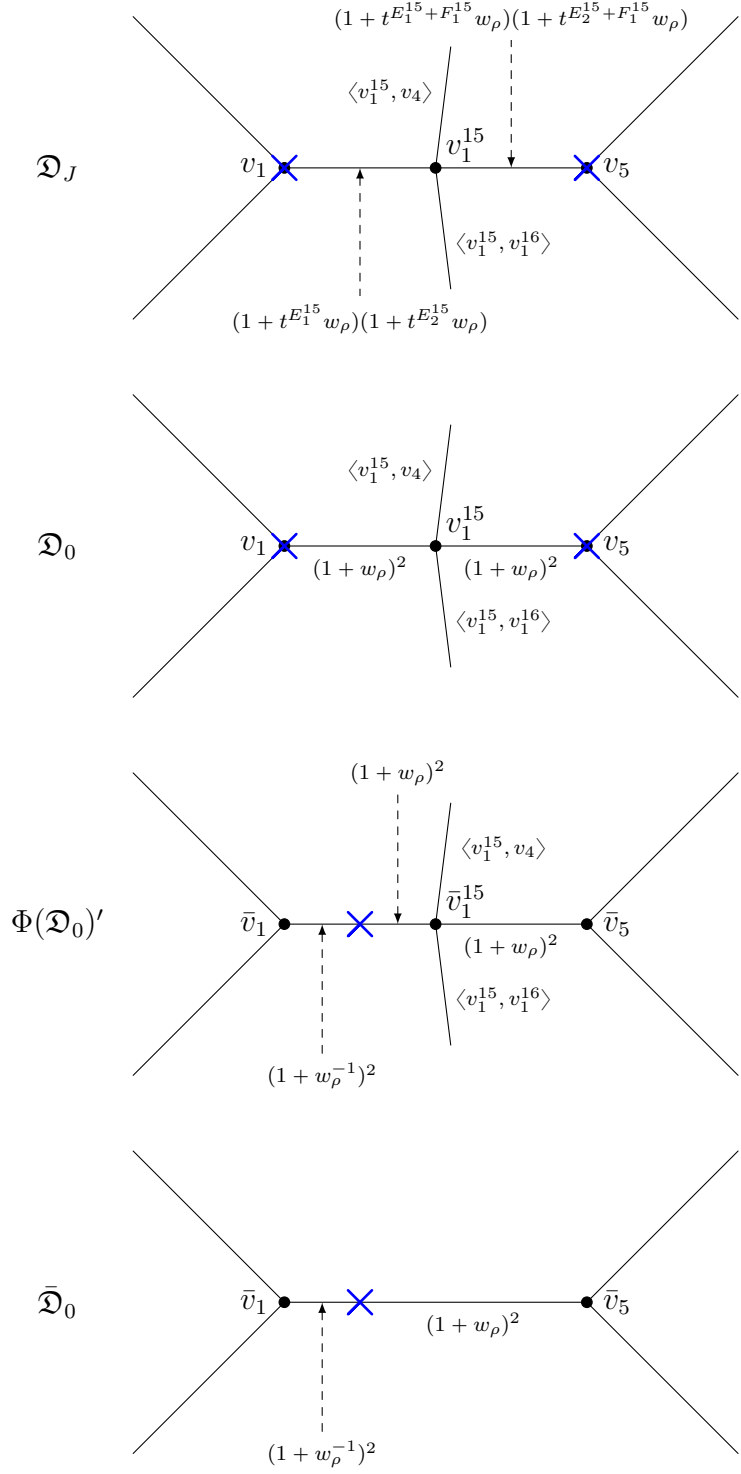

We check consistency of $\Phi(\mathfrak{D}_k)$.

\begin{fact} \label{consis}
For every $k \geq 0$, the scattering diagram $\Phi(\mathfrak{D}_k)$ is consistent.
\end{fact}

\begin{proof}
By Remark \ref{obsplimage}(1), it is enough to show that $\Phi(\mathfrak{D}_k)'$ is consistent.
By Definition \ref{consistency} of consistency, we need to check that $\Phi(\mathfrak{D}_k)$ is consistent in codimensions $0$, $1$, and $2$. This will follow from consistency of $\mathfrak{D}_k$ and Construction \ref{chs}. 

\textbf{Consistency in codimension $\mathbf{0}$.} Every codimension $0$ joint of $\bar{\mathscr{P}}_{\Phi(\mathfrak{D}_k)'}'^{[n-2]} = \bar{\mathscr{P}}_{\Phi(\mathfrak{D}_k)'}'^{[0]}$ is of the form $\Phi(\mathfrak{j})$ for $\mathfrak{j}$ a codimension $0$ joint of $\mathscr{P}_{\mathfrak{D}_k}^{[n-2]} = \mathscr{P}_{\mathfrak{D}_k}^{[0]}$. Moreover, the walls containing $\Phi(\mathfrak{j})$ are of the form $\Phi(\mathfrak{p}_1), \dots, \Phi(\mathfrak{p}_r)$ for $\mathfrak{p}_1, \dots, \mathfrak{p}_r$ the walls containing $\mathfrak{j}$ and if $\sigma \in \mathscr{P}^{\max}$ is the maximal cell containing $\mathfrak{j}$, then $\Phi(\sigma) \in \bar{\mathscr{P}}'^{\max}$ is the maximal cell containing $\Phi(\mathfrak{j})$.  Recall from Construction \ref{plmap} (via Construction \ref{plisosmall}) that $\Phi$ induces an isomorphism $\mathcal{P}_x^{+} \cong \mathcal{\bar{P}}_{\Phi(x)}^{+}$ for every $x \in \Int(\sigma), \ \sigma \in \mathscr{P}^{\max}$. This implies that we have an isomorphism $R_{\sigma} \cong R_{\Phi(\sigma)}$ and for any $1 \leq i \leq r$, the automorphism $\theta_{\Phi(\mathfrak{p}_i)}$ of $R_{\Phi(\sigma)}$ is the same as the automorphism $\theta_{\mathfrak{p}_i}$ of $R_{\sigma}$ under this identification. But then $\theta_{\mathfrak{p}_1} \circ \dots  \circ \theta_{\mathfrak{p}_r} = \Id$ implies that also $\theta_{\Phi(\mathfrak{p}_1)} \circ \dots  \circ \theta_{\Phi(\mathfrak{p}_r)} = \Id$ and $\Phi(\mathfrak{D}_k)'$ is consistent in codimension $0$. 

\textbf{Consistency in codimension $\mathbf{1}$.} Every joint $\mathfrak{j} \in \mathscr{P}_{\mathfrak{D}_k}^{[0]}$ of codimension $1$ gives a joint $\Phi(\mathfrak{j}) \in \bar{\mathscr{P}}_{\Phi(\mathfrak{D}_k)'}'^{[0]}$ of codimension $1$ and we have $\Phi(\mathfrak{j}) \subseteq \Phi(\rho') \in \bar{\mathscr{P}}'^{[1]}$ if $\mathfrak{j} \subseteq \rho' \in \mathscr{P}^{[1]}$. We also have joints $\mathfrak{j}^{\speci}_{\rho} := x_{\Phi(\rho)}$ for every edge $\rho'$ of the form $\rho_{p_0+1} \in \mathscr{P}^{[1]}$ and $x_{\Phi(\rho)}$ the singular point of $\Phi(\rho_{p_0+1})$. Moreover, if $\sigma, \sigma' \in \mathscr{P}^{\max}$ are the maximal cells containing $\rho'$, then $\Phi(\sigma), \Phi(\sigma') \in \bar{\mathscr{P}}'^{\max}$ are the maximal cells containing $\Phi(\rho')$. Using the conventions of Construction \ref{chs}, all the joints $\Phi(\mathfrak{j}) \in \bar{\mathscr{P}}_{\Phi(\mathfrak{D}_k)'}'^{[0]}$ of codimension $1$ fall into one of the following cases:
\begin{enumerate}
    \item[Case 1.] One of the following holds:
    \begin{enumerate}
        \item $\mathfrak{j} \subseteq \Int(\rho'), \ \rho' \in \mathscr{P}^{[1]} \setminus \mathscr{P}^{[1]}_{\coar}$ (using Notation \ref{nonsubdividing}).
        \item $\mathfrak{j} \subseteq \Int(\rho_p), \ \rho_p \in \mathscr{P}^{[1]}$ and $\Phi(\mathfrak{j}) \subseteq \Int(\underline{\Phi(\rho)})$ for $\rho \in \mathscr{P}^{[1]}_{\coar}$ the edge containing $\rho_p$.
        \item $\mathfrak{j} \subseteq \Int(\rho_p), \ \rho_p \in \mathscr{P}^{[1]}$ for $1 \leq p \leq p_0$ and $\Phi(\mathfrak{j}) \subseteq \Int(\Phi(\rho_p)) \subseteq \Int(\underline{\Phi(\rho)}')$.
    \end{enumerate}
    \item[Case 2.] $\mathfrak{j} \subseteq \Int(\rho_{p_0+1}), \ \rho_{p_0+1} \in \mathscr{P}^{[1]}$ and $\Phi(\mathfrak{j}) \subseteq \Int(\Phi(\rho_{p_0+1}) \cap \underline{\Phi(\rho)}')$.
    \item[Case 3.] $\mathfrak{j}^{\speci}_{\rho} \subseteq \Int(\Phi(\rho_{p_0+1}))$ of the form $\mathfrak{j}^{\speci}_{\rho} = x_{\Phi(\rho)}$. 
\end{enumerate}
Here in Cases 1, 2, and 3, the wall functions for the slabs of $\Phi(\mathfrak{D}_k)'$ containing $\mathfrak{j}$ arise as in Construction \ref{chs}(1) (or \eqref{irrslabimage}), Construction \ref{chs}(2), and Construction \ref{chs}(3) respectively. 
We consider these cases one by one.

\textbf{Case 1.} In this case, the proof of consistency is similar to the proof of consistency in codimension $0$. We have $R_{\sigma} \cong R_{\Phi(\sigma)}$ and $R_{\sigma'} \cong R_{\Phi(\sigma')}$. Moreover, the two slabs $\Phi(\mathfrak{b}_1), \Phi(\mathfrak{b}_2)$ containing $\Phi(\mathfrak{j})$ are of the first type. This implies that we have natural isomorphisms $R_{\mathfrak{b}_1} \cong R_{\Phi(\mathfrak{b}_1)}, \ R_{\mathfrak{b}_2} \cong R_{\Phi(\mathfrak{b}_2)}$. Then for $\mathfrak{b}_1, \mathfrak{p}_1, \dots, \mathfrak{p}_r, \mathfrak{b}_2, \mathfrak{p}'_1, \dots, \mathfrak{p}'_s$ a cyclic ordering of walls around $\mathfrak{j}$ as in Definition \ref{codim1consistency}, the automorphism $\theta_{\Phi(\mathfrak{p}_i)}$ of $R_{\Phi(\sigma)}$ is the same as the automorphism $\theta_{\mathfrak{p}_i}$ of $R_{\sigma}$ and the automorphism $\theta_{\Phi(\mathfrak{p}'_j)}$ of $R_{\Phi(\sigma')}$ is the same as the automorphism $\theta_{\mathfrak{p}'_j}$ of $R_{\sigma'}$ under the identifications. Further, the localization homomorphisms $\chi_{\Phi(\mathfrak{b}_i), \Phi(\sigma)}, \chi_{\Phi(\mathfrak{b}_i), \Phi(\sigma')}$ for $i=1,2$ are the same as $\chi_{\mathfrak{b}_i, \sigma}, \chi_{\mathfrak{b}_i, \sigma'}$  under the identifications of rings. Let
$$
\begin{aligned}
\Phi(\theta) &:= \theta_{\Phi(\mathfrak{p}_r)} \circ \theta_{\Phi(\mathfrak{p}_{r-1})} \circ \dots \circ \theta_{\Phi(\mathfrak{p}_{1})}: \ &R_{\Phi(\sigma)} \to R_{\Phi(\sigma)} \\
\Phi(\theta') &:= \theta_{\Phi(\mathfrak{p}'_1)} \circ \theta_{\Phi(\mathfrak{p}'_{2})} \circ \dots \circ \theta_{\Phi(\mathfrak{p}'_{s})}: \ &R_{\Phi(\sigma')} \to R_{\Phi(\sigma')}
\end{aligned}
$$
be defined similarly to $\theta$ and $\theta'$ in Definition \ref{codim1consistency}. Consistency of $\mathfrak{D}_k$ around $\mathfrak{j}$ means that 
\begin{equation} \label{proof1}
(\theta \times \theta')\left((\chi_{\mathfrak{b}_1, \sigma}, \chi_{\mathfrak{b}_1, \sigma'})(R_{\mathfrak{b}_1})\right) = (\chi_{\mathfrak{b}_2, \sigma}, \chi_{\mathfrak{b}_2, \sigma'})(R_{\mathfrak{b}_2})
\end{equation}
which implies that $$(\Phi(\theta) \times \Phi(\theta'))\left((\chi_{\Phi(\mathfrak{b}_1), \Phi(\sigma)}, \chi_{\Phi(\mathfrak{b}_1), \Phi(\sigma')})(R_{\Phi(\mathfrak{b}_1)})\right) = (\chi_{\Phi(\mathfrak{b}_2), \sigma}, \chi_{\Phi(\mathfrak{b}_2), \Phi(\sigma')})(R_{\Phi(\mathfrak{b}_2)})$$ by the discussion above. So $\Phi(\mathfrak{D}_k)'$ is consistent around $\Phi(\mathfrak{j})$.

\begin{sloppypar}
\textbf{Case 2.} To prove consistency in this case, note that the two slabs $\Phi(\mathfrak{b}_1)'$, $\Phi(\mathfrak{b}_2)'$ containing $\Phi(\mathfrak{j})$ are of the second type. The only difference with the first case is that we no longer have natural isomorphisms $R_{\mathfrak{b}_1} \cong R_{\Phi(\mathfrak{b}_1)'}, \ R_{\mathfrak{b}_2} \cong R_{\Phi(\mathfrak{b}_2)'}$ and need to choose the isomorphisms. We have \end{sloppypar}
$$
\begin{aligned}
R_{\mathfrak{b}_i} &= \kk[t]/(t^{k+1})[\Lambda_{\rho}][Z_+, Z_-]/\left(Z_+Z_- - f_{\mathfrak{b}_i} \cdot z^{\kappa_{\rho_{p_0+1}}}\right)  \\
R_{\Phi(\mathfrak{b}_i)'} &= \kk[t]/(t^{k+1})[\Lambda_{\Phi(\rho)}][Z_+, Z_-]/\left(Z_+Z_- - z^{m_{\underline{\Phi(\rho)}' \underline{\Phi(\rho)}}}\Phi(f_{\mathfrak{b}_i}) \cdot z^{\kappa_{\Phi(\rho_{p_0+1})}}\right)
\end{aligned}
$$ for $i=1,2$ and $\kappa_{\rho_{p_0+1}} = \kappa_{\Phi(\rho_{p_0+1})}$ by construction. We define an isomorphism\footnote{This is similar to the isomorphism $R_{\underline{\rho}} \to R_{\underline{\rho}'}$ for $\underline{\rho}, \underline{\rho}' \subseteq \rho \in \mathscr{P}^{[1]}$ two slabs constructed in \cite[Lemma 2.2.3]{GHS}.} 
\begin{equation} \label{proof2}
\Phi': R_{\mathfrak{b}_i} \to R_{\Phi(\mathfrak{b}_i)'}, \ Z_{+} \mapsto Z_{+}, \ Z_{-} \mapsto  z^{m_{\underline{\Phi(\rho)} \, \underline{\Phi(\rho)}'}} Z_{-}
\end{equation} 
with $\Lambda_{\rho}$ and $\Lambda_{\Phi(\rho)}$ identified by the PL-isomorphism. We now have 
\begin{equation} \label{compatibilityofchi}
(\chi_{\Phi(\mathfrak{b}_i)', \Phi(\sigma)}, \chi_{\Phi(\mathfrak{b}_i)', \Phi(\sigma')})(R_{\Phi(\mathfrak{b}_i)'}) = (\chi_{\mathfrak{b}_i, \sigma}, \chi_{\mathfrak{b}_i, \sigma'})(R_{\mathfrak{b}_i}) 
\end{equation}
for $i=1,2$
under the identifications $R_{\sigma} \cong R_{\Phi(\sigma)}, \ R_{\sigma'} \cong R_{\Phi(\sigma')}$. But then \eqref{proof1} implies that 
$$(\Phi(\theta) \times \Phi(\theta'))\left((\chi_{\Phi(\mathfrak{b}_1)', \Phi(\sigma)}, \chi_{\Phi(\mathfrak{b}_1)', \Phi(\sigma')})(R_{\Phi(\mathfrak{b}_1)'})\right) = (\chi_{\Phi(\mathfrak{b}_2)', \sigma}, \chi_{\Phi(\mathfrak{b}_2)', \Phi(\sigma')})(R_{\Phi(\mathfrak{b}_2)'}).$$ So $\Phi(\mathfrak{D}_k)'$ is consistent around $\Phi(\mathfrak{j})$.

\textbf{Case 3.} Recall that $\mathfrak{j}_{\rho}^{\speci} = x_{\Phi(\rho)}$ is an irrational point of $\Phi(\rho_{p_0+1})$ and since both the walls on $\mathfrak{D}_k$ and the PL-isomorphism $\Phi$ are rationally defined, the only walls containing $\mathfrak{j}_{\rho}^{\speci}$ are the two slabs $\bar{\Phi}(\mathfrak{b})=(\bar{\Phi}(\mathfrak{b}), f_{\Phi(\mathfrak{b})})$ and $\bar{\Phi}(\mathfrak{b})' =(\bar{\Phi}(\mathfrak{b})', z^{m_{\underline{\Phi(\rho)}'\underline{\Phi(\rho)} }}f_{\Phi(\mathfrak{b})})$ defined as in Construction \ref{chs}(3) for some slab $\mathfrak{b}$ of $\mathfrak{D}_k$. Now, we have a natural isomorphism $R_{\mathfrak{b}} \cong R_{\bar{\Phi}(\mathfrak{b})}$ and an isomorphism $\Phi': R_{\mathfrak{b}} \to R_{\bar{\Phi}(\mathfrak{b})'}$ defined as in \eqref{proof2}.
But then 
$$(\chi_{\bar{\Phi}(\mathfrak{b}), \Phi(\sigma)}, \chi_{\bar{\Phi}(\mathfrak{b}), \Phi(\sigma')})(R_{\bar{\Phi}(\mathfrak{b})}) = 
(\chi_{\mathfrak{b}, \sigma}, \chi_{\mathfrak{b}, \sigma'})(R_{\mathfrak{b}})=
(\chi_{\bar{\Phi}(\mathfrak{b})', \Phi(\sigma)}, \chi_{\bar{\Phi}(\mathfrak{b})', \Phi(\sigma')})(R_{\bar{\Phi}(\mathfrak{b})'})
$$
under the identifications $R_{\sigma} \cong R_{\Phi(\sigma)}, \ R_{\sigma'} \cong R_{\Phi(\sigma')}$ where the second equality follows as in \eqref{compatibilityofchi}. Since there are no other walls containing $\mathfrak{j}_{\rho}^{\speci}$, this implies that $\Phi(\mathfrak{D}_k)'$ is consistent around $\mathfrak{j}_{\rho}^{\speci}$.\footnote{Note that for $k=0$, we recover the fact that every scattering diagram is consistent modulo $I_0$.}

\textbf{Consistency in codimension $\mathbf{2}$.} Every codimension $2$ joint of $\bar{\mathscr{P}}_{\Phi(\mathfrak{D}_k)'}'^{[0]}$ is of the form $\Phi(\mathfrak{j})$ for $\mathfrak{j}$ a codimension $2$ joint of $\mathscr{P}_{\mathfrak{D}_k}^{[0]}$. First, we note that the PL-isomorphism of Construction \ref{plmapadmissible} induces a PL-isomorphism $\Phi: \left(B_{\mathfrak{j}}, \mathscr{P}_{\mathfrak{j}, \mathfrak{D}_{k, \mathfrak{j}}} \right) \to \left(\bar{B}_{\Phi(\mathfrak{j})}, \bar{\mathscr{P}}'_{\Phi(\mathfrak{j}), \Phi(\mathfrak{D}_{k})'_\mathfrak{j}} \right)$. Consistency of $\mathfrak{D}_k$ along $\mathfrak{j}$ means that for a general point $p \in B_{\mathfrak{j}}$ contained in a chamber $\mathfrak{u}$ of $\mathscr{P}_{\mathfrak{j}, \mathfrak{D}_{k, \mathfrak{j}}}$, the theta functions $\vartheta_m^{\mathfrak{j}}(p)=\sum_{\beta} a_{\beta} z^{m_{\beta}} \in R_{\mathfrak{u}}$ of \eqref{thetafunc} satisfy conditions (1) and (2) of Definition \ref{consistency}. We need to show that a similar property holds for $B_{\Phi(\mathfrak{j})}$.

As a set, we define the image of a broken line $\beta: (-\infty, 0] \to B_{\mathfrak{j}}^0$ with endpoint $p$ to be the image $\Phi(\beta): (-\infty, 0] \to B_{\Phi(\mathfrak{j})}^0$ under the PL-isomorphism (with endpoint $\Phi(p)$). We also define the attached monomial on $\Phi(\beta)([t_{i-1}, t_i])$ to be $a_i z^{\Phi(m_i)}$ if the attached monomial on $\beta([t_{i-1}, t_i])$ is $a_i z^{m_i}$. This gives a one-to-one correspondence since $\Phi$ is a PL-isomorphism. We need to check that $\Phi(\beta)$ satisfies the properties of a broken line of Definition \ref{brokenlines}. 

The only non-trivial check is that $a_{i+1} z^{\Phi(m_{i+1})}$ is the result of transport of $a_i z^{\Phi(m_i)}$ from $\Phi(\mathfrak{u})$ to $\Phi(\mathfrak{u}')$ (here $\Phi(\mathfrak{u})$ and $\Phi(\mathfrak{u}')$ are the chambers of  $\bar{\mathscr{P}}'_{\Phi(\mathfrak{j}), \Phi(\mathfrak{D}_{k})'_{\mathfrak{j}}}$ corresponding to the chambers $\mathfrak{u}$ and $\mathfrak{u}'$ of $\mathscr{P}_{\mathfrak{j}, \mathfrak{D}_{k, \mathfrak{j}}}$ respectively). If $\Phi(\mathfrak{u}) \cap \Phi(\mathfrak{u}') = \Phi(\mathfrak{p})$ is a codimension $0$ wall, this follows from the canonical identifications $R_{\sigma} \cong R_{\Phi(\sigma)}, \ R_{\sigma'} \cong R_{\Phi(\sigma')}$ for $\sigma \in \mathscr{P}_{\mathfrak{j}}^{\max}$ and $\sigma' \in \mathscr{P}_{\mathfrak{j}}^{\max}$ the unique cells containing $\mathfrak{u}$ and $ \mathfrak{u}'$  respectively. Similarly, if $\Phi(\mathfrak{u}) \cap \Phi(\mathfrak{u}') = \Phi(\mathfrak{b})$ (or $\Phi(\mathfrak{b})=\bar{\Phi}(\mathfrak{b})$ as in Construction \ref{chs}(3)) is a slab of the first type, this follows from the identifications $R_{\sigma} \cong R_{\Phi(\sigma)}, \ R_{\sigma'} \cong R_{\Phi(\sigma')}$ along with the natural isomorphism $R_{\mathfrak{b}} \cong R_{\Phi(\mathfrak{b})}$. If $\Phi(\mathfrak{u}) \cap \Phi(\mathfrak{u}') = \Phi(\mathfrak{b})'$ (or $\bar{\Phi}(\mathfrak{b})'$ as in Construction \ref{chs}(3)) is a slab of the second type, the condition follows from the fact that the diagram 
\begin{equation} \label{commutativityofchangeofchamber}
        \begin{tikzcd}
        R_{\mathfrak{b}} \arrow[hookrightarrow]{r} \arrow{d}{\Phi'} & R_{\mathfrak{u}}^{\mathfrak{b}} \arrow{r}{\theta_{\mathfrak{b}}} & R_{\mathfrak{u}'} \arrow{d}{\cong}   \\
         R_{\Phi(\mathfrak{b})'} \arrow[hookrightarrow]{r} & R_{\Phi(\mathfrak{u})}^{\Phi(\mathfrak{b})'} \arrow{r}{{\theta_{\Phi(\mathfrak{b})'}}} & R_{\Phi(\mathfrak{u}')}
        \end{tikzcd}
        \end{equation}
with the rows defined by \eqref{slabcrossing} and $\Phi'$ the isomorphism of \eqref{proof2} is commutative.

Every chamber of $\bar{\mathscr{P}}'_{\Phi(\mathfrak{j}), \Phi(\mathfrak{D}_{k})'_{\mathfrak{j}}}$ is of the form $\Phi(\mathfrak{u})$ for $\mathfrak{u}$ a chamber of $\mathscr{P}_{\mathfrak{j}, \mathfrak{D}_{k, \mathfrak{j}}}$ and every asymptotic monomial on $(\bar B_{\Phi(\mathfrak{j})}, \bar{\mathscr{P}}'_{\Phi(\mathfrak{j})})$ is of the form $\Phi(m)$ for $m$ an asymptotic monomial on $\left( B_{\mathfrak{j}}, \mathscr{P}_{\mathfrak{j}} \right)$. So the theta functions $\vartheta_{\Phi(m)}^{\Phi(\mathfrak{j})}(\Phi(p))$ on $(\bar B_{\Phi(\mathfrak{j})}, \bar{\mathscr{P}}'_{\Phi(\mathfrak{j})})$ are of the form  $\vartheta_{\Phi(m)}^{\Phi(\mathfrak{j})}(\Phi(p))=\sum_{\beta} a_{\beta} z^{\Phi(m_{\beta})} \in R_{\Phi(\mathfrak{u})}$ for $\vartheta_m^{\mathfrak{j}}(p)=\sum_{\beta} a_{\beta} z^{m_{\beta}} \in R_{\mathfrak{u}}$ the corresponding theta function on $\left( B_{\mathfrak{j}}, \mathscr{P}_{\mathfrak{j}} \right)$. 

This immediately implies that $\vartheta_{\Phi(m)}^{\Phi(\mathfrak{j})}(\Phi(p))$ satisfy condition (1) of Definition \ref{consistency} and satisfy condition (2) of Definition \ref{consistency} for all the change of chamber homomorphisms $\theta_{\Phi(\mathfrak{u}') \Phi(\mathfrak{u})}$ with $\Phi(\mathfrak{u}') \cap \Phi(\mathfrak{u})$ either a codimension $0$ wall or a slab of the first type. Compatibility with the change of chamber homomorphism $\theta_{\Phi(\mathfrak{u}') \Phi(\mathfrak{u})}$ for $\Phi(\mathfrak{u}') \cap \Phi(\mathfrak{u}) = \Phi(\mathfrak{b})'$ a slab of the second type again follows from commutativity of \eqref{commutativityofchangeofchamber}.   
\end{proof}

Now the construction of $\check{\mathfrak{X}}^o$ implies that $\Phi(\mathfrak{D}_k)$ is equivalent to $\mathfrak{D}_k$.

\begin{fact} \label{equivalence}
For every $k \geq 0$, the scattering diagram $\Phi(\mathfrak{D}_k)$ is equivalent to $\mathfrak{D}_k$. 
\end{fact}

\begin{proof}
By Remark \ref{obsplimage}(1), it is enough to show that $\Phi(\mathfrak{D}_k)'$ is equivalent to $\mathfrak{D}_k$. We need to check that the families $\check{\mathfrak{X}}^o_{(B, \mathfrak{D}_k)}$ and $\check{\mathfrak{X}}^o_{(\bar{B}, \Phi(\mathfrak{D}_k)')}$ are isomorphic.

As in the proof of \cite[Proposition 2.4.1]{GHS}, consistency of $\mathfrak{D}_k$  implies that all the rings $R_{\mathfrak{u}}$ for chambers $\mathfrak{u} \in \mathscr{P}_{\mathfrak{D}_k}$ contained in the same $  \sigma \in \mathscr{P}^{\max}$ are canonically isomorphic and all the rings $R_{\mathfrak{b}}$ for slabs $\mathfrak{b} \in \mathfrak{D}_k$ contained in the same $ \rho' \in \mathscr{P}^{[1]}$ are canonically isomorphic. Fixing one chamber $\mathfrak{u}$ for every $\sigma \in \mathscr{P}^{\max}$ and one slab $\mathfrak{b}$ for every $\rho' \in \mathscr{P}^{[1]}$, we can glue $\check{\mathfrak{X}}^o_{(B, \mathfrak{D}_k)}$ from $\left\{\Spec R_{\mathfrak{u}}, \sigma \in \mathscr{P}^{\max} \right\}$ and $\left\{\Spec R_{\mathfrak{b}}, \rho' \in \mathscr{P}^{[1]} \right\}$ via the maps $\chi_{\mathfrak{b}, \mathfrak{u}}, \chi_{\mathfrak{b}, \mathfrak{u}'}$ that are compositions of the localization homomorphisms with isomorphisms of $R_{\mathfrak{u}}$ for $\mathfrak{u} \subseteq \sigma \in \mathscr{P}^{\max}$. The same is true for $\check{\mathfrak{X}}^o_{(\bar{B}, \Phi(\mathfrak{D}_k)')}$ by consistency of $\Phi(\mathfrak{D}_k)'$, so we are free to choose the chambers and slabs in both constructions. 

We make an arbitrary choice of chambers $\mathfrak{u}$ for every $\sigma \in \mathscr{P}^{\max}$ and use the choice $\Phi(\mathfrak{u})$ for $\Phi(\sigma) \in \bar{\mathscr{P}}'^{\max}$. We choose $\mathfrak{b}$ for every $\rho' \in \mathscr{P}^{[1]}$ so that its image under the PL-isomorphism $\Phi$ is either $\Phi(\mathfrak{b})$ or contains $\bar{\Phi}(\mathfrak{b})$ using the setup of Construction \ref{chs}. We use the choice $\Phi(\mathfrak{b})$ or $\bar{\Phi}(\mathfrak{b})$ respectively for $\Phi(\rho') \in \bar{\mathscr{P}}'^{[1]}$ and drop the bar in the second case.

As in the proof of Proposition \ref{consis}, we have natural isomorphisms $R_{\mathfrak{u}} \cong R_{\Phi(\mathfrak{u})}$ and $R_{\mathfrak{b}} \cong R_{\Phi(\mathfrak{b})}$. Moreover, $\chi_{\Phi(\mathfrak{b}), \Phi(\mathfrak{u})}, \chi_{\Phi(\mathfrak{b}), \Phi(\mathfrak{u}')}$  are the same as $\chi_{\mathfrak{b}, \mathfrak{u}}, \chi_{\mathfrak{b}, \mathfrak{u}'}$  under these identifications. This implies that $\check{\mathfrak{X}}^o_{(B, \mathfrak{D}_k)} \cong \check{\mathfrak{X}}^o_{(\bar{B}, \Phi(\mathfrak{D}_k))}$ by the description of the gluing above.
\end{proof}

\subsubsection{Proof of the conjecture}

We combine the results of this chapter to prove the conjecture.

\begin{theorem} \label{proofconj}
Conjecture \ref{conjecture} holds for special toric degenerations of K3-s. 
\end{theorem}

\begin{proof}
By Proposition \ref{whenspecialk3}, a toric degeneration $\bar{\mathfrak{X}} \to \mathcal{S}$ of K3-s is special if and only if it is a divisorial log deformation and the generic fibre of $\bar{\mathfrak{X}} \to \mathcal{S}$ is smooth. Sections \ref{resolvesmall}, \ref{resolve}, and \ref{admissible} construct resolutions in increasing generality. Proposition \ref{furtherresolution} implies that it is enough to prove the conjecture for strongly admissible resolutions, and Proposition \ref{strongadmissibleresolutionresult} constructs a strongly admissible resolution $\pi:\mathfrak{X} \to \bar{\mathfrak{X}}$ of $\bar{\mathfrak{X}} \to \mathcal{S}$ to a minimal log CY degeneration $\mathfrak{X} \to \mathcal{S}$ (with $D$ simple normal crossings) in general. It is enough to prove the conjecture for well-chosen monoids $P$ by Proposition \ref{enoughforwellchosen} (well-chosen monoids exist by Proposition \ref{wellchosen}). For every strongly admissible resolution $\pi:\mathfrak{X} \to \bar{\mathfrak{X}}$ and a well-chosen monoid $P$, Section \ref{admissibleandstronglyso} defines the data of Basic Setup \ref{setupforproof}. 

In Section \ref{scatint}, we constructed a collection of canonical scattering diagrams $\mathfrak{D}^J: = \left\{ \mathfrak{D}_{J^{k+1}}, \ k \geq 0  \right\}$ on $(B, \mathscr{P})$ giving rise to the extended intrinsic mirror $\check{\mathfrak{X}} \to \Spec \widehat{\kk [P]}_{J}$. Construction \ref{ccc} defines a collection of scattering diagrams $\mathfrak{D} = \left\{ \mathfrak{D}_k, \ k \geq 0 \right\}$ on $(B, \mathscr{P})$ with monoid $\mathbb{N}$ and $\mathfrak{D}_k:=\mathfrak{D}_{(t^{k+1})}$. By Proposition \ref{pr}, the basechange of the extended intrinsic mirror $\check{\mathfrak{X}} \to \Spec \widehat{\kk [P]}_J$ by $h: P \to \mathbb{N}, \ \beta \mapsto \pi^{*}A \cdot \beta$ is isomorphic to the family $\check{\mathfrak{X}}_{\mathfrak{D}} \to \Spec \kk \llbracket t \rrbracket$ of Construction \ref{ccc}. 

It is enough to show that the families $\check{\mathfrak{X}}_{\mathfrak{D}_k} \to \Spec \kk [t]/(t^{k+1})$ and $\check{\bar{\mathfrak{X}}}_{\bar{\mathfrak{D}}_k} \to \Spec \kk [t]/(t^{k+1})$ are isomorphic for every $k \geq 0$. Here the second family is defined in \eqref{reconstructedfamily} using Theorem \ref{reconstruction} that can be applied since local rigidity is an empty condition in dimension $2$. By Proposition \ref{help} and Observation \ref{obsplimage}(1), it is enough to show that $(B, \mathfrak{D}_k)$ is equivalent to $(\bar{B}, \bar{\mathfrak{D}}_k)$ for every $k \geq 0$. Since $\Phi(\mathfrak{D}_k)$ is consistent by Proposition \ref{consis} and $(B, \mathfrak{D}_k)$ is equivalent to $(\bar{B}, \Phi(\mathfrak{D}_k))$ by Proposition \ref{equivalence}, it is enough to prove that $\Phi(\mathfrak{D}_k)$ is equivalent to $\bar{\mathfrak{D}}_k$. By Observation \ref{combimpliesequiv}(2), it is enough to prove that they are combinatorially equivalent.  

We want to use the uniqueness statement of Theorem \ref{reconstruction}. Consider two sequences of scattering diagrams: $$\bar{\mathfrak{D}} = \left\{\bar{\mathfrak{D}}_k, \ k \geq 0\right\}, \quad \Phi(\mathfrak{D}):=\left\{\Phi(\mathfrak{D}_k), \ k \geq 0 \right\}.$$ By the existence part of Theorem \ref{reconstruction}, $\bar{\mathfrak{D}}_k$ is compatible with $\bar{\mathfrak{D}}_{k-1}$ for $k \geq 1$. By Observation \ref{obsplimage}(2), $\Phi(\mathfrak{D}_k)$ is compatible with $\Phi(\mathfrak{D}_{k-1})$ for $k \geq 1$. Moreover, $\Phi(\mathfrak{D}_0)$ is combinatorially equivalent to $\bar{\mathfrak{D}}_0$ by Observation \ref{obsplimage}(3). Every $\bar{\mathfrak{D}}_k, \ k \geq 0$ is consistent in the sense of \cite[Definition 2.28]{GS1} by the existence part of Theorem \ref{reconstruction}. It remains to prove that every $\Phi(\mathfrak{D}_k), \ k \geq 0$ is consistent in the sense of \cite[Definition 2.28]{GS1}.

By Proposition \ref{consis}, $\Phi(\mathfrak{D}_k)$ is consistent for every $k \geq 0$. By Proposition \ref{consistencydiff}, it suffices to check that $\Phi(\mathfrak{D}_k)$ is consistent in the sense of \cite[Definition 2.28]{GS1} around every codimension one joint $\mathfrak{j} \subseteq \bar \Delta$. Every such joint is of the form $\mathfrak{j}_{\rho}^{\speci} = x_{\Phi(\rho)}$ for $x_{\Phi(\rho)}$ the singularity of $\Phi(\rho)$. As in Case 3 in the proof of Proposition \ref{consis}, the fact that $x_{\Phi(\rho)}$ is an irrational point implies that $\mathfrak{j}_{\rho}^{\speci}$ is contained in exactly two walls that are the slabs $\bar{\Phi}(\mathfrak{b})=(\bar{\Phi}(\mathfrak{b}), f_{\Phi(\mathfrak{b})})$ and $\bar{\Phi}(\mathfrak{b})' =(\bar{\Phi}(\mathfrak{b})', z^{m_{\underline{\rho}'\underline{\rho} }}f_{\Phi(\mathfrak{b})})$. Since there are no other walls containing $\mathfrak{j}_{\rho}^{\speci}$ and the slab functions satisfy $f_{\bar{\Phi}(\mathfrak{b})'} = z^{m_{\Phi(\underline{\rho})'\Phi(\underline{\rho}) }}f_{\bar{\Phi}(\mathfrak{b})}$, consistency around $ \mathfrak{j}_{\rho}^{\speci}$ in the sense of \cite[Definition 2.28]{GS1} follows exactly as in the proof of consistency of $\bar{\mathfrak{D}}_0$ in the sense of \cite[Definition 2.28]{GS1} (in the case of trivial gluing data), given in \cite[Proposition 3.2]{GS1}.
\end{proof}

\section[Intrinsic mirror over the Gross-Siebert locus]{Intrinsic mirror over the minimal relative Gross-Siebert locus} \label{gslocus}

In this chapter, we generalize Conjecture \ref{conjecture} following the plan of Section \ref{problem}.

\subsection{Setup for the generalizations} \label{gensetup}

Let $\bar{\mathfrak{X}} \to \mathcal{S}$ be a special toric degeneration of K3-s as before. In Basic Setup \ref{setupforproof}, we assumed a choice of polarization $A$ on $\bar{\mathfrak{X}} \to \mathcal{S}$. In this chapter, we just require that there exists some polarization $A$ on $\bar{\mathfrak{X}} \to \mathcal{S}$. We also require that $\bar{\mathfrak{X}} \to \mathcal{S}$ satisfies Assumption \ref{techass}. Let $\pi: \mathfrak{X} \to \bar{\mathfrak{X}}$ be a strongly admissible resolution (see Definition \ref{stronglyadmissible}) to a log smooth minimal log CY degeneration $\mathfrak{X} \to \mathcal{S}$. We further require that $A_1(\mathfrak{X}_0, \mathbb{Z}) = A_1(\mathfrak{X}_0, \mathbb{Z})_{\numer}$. 

\begin{remark} \label{restrictingtonumer}
Note that for any well-chosen monoid $NE(\mathfrak{X}_0) \subseteq P \subseteq A_1 (\mathfrak{X}_0, \mathbb{Z})$ with a face $K$ containing the classes of the contracted curves (see Definition \ref{wellchosendef}), we can construct a well-chosen monoid $NE(\mathfrak{X}_0)_{\numer} \subseteq P_{\numer} \subseteq A_1 (\mathfrak{X}_0, \mathbb{Z})_{\numer}$ by setting $P_{\numer} := P \cap A_1(\mathfrak{X}_0, \mathbb{Z})_{\numer}$ and $ K_{\numer} := K \cap A_1(\mathfrak{X}_0, \mathbb{Z})_{\numer}$ using the splitting $A_1(\mathfrak{X}_0, \mathbb{Z})_{\numer} \oplus G$ of Construction \ref{splittings}. Then for $J_{\numer}:= P_{\numer} \setminus K_{\numer}$, we have a map $\kk[P]_{J} \to \kk[P_{\numer}]_{J_{\numer}}$ by sending  $t^{\beta} \mapsto t^{\beta}$ for $\beta \in A_1(\mathfrak{X}_0, \mathbb{Z})_{\numer}$ and sending $t^{\beta} \mapsto 0$ for $\beta \in G$. The extended intrinsic mirror $\check{\mathfrak{X}} \to \Spec \kk[P_{\numer}]_{J_{\numer}}$ is the basechange of $\check{\mathfrak{X}} \to \Spec \kk[P]_{J}$ via this map (this is easy to see both from the construction of the extended intrinsic mirror of Section \ref{extresults} and from the scattering diagram description of Section \ref{scatint}) that we call the \emph{numerical locus} of $\check{\mathfrak{X}} \to \Spec \kk[P]_{J}$. 

Conversely, given a well-chosen monoid $NE(\mathfrak{X}_0)_{\numer} \subseteq P_{\numer} \subseteq A_1 (\mathfrak{X}_0, \mathbb{Z})_{\numer}$ with a face $K_{\numer}$, we can construct a well-chosen monoid $NE(\mathfrak{X}_0) \subseteq P \subseteq A_1 (\mathfrak{X}_0, \mathbb{Z})$ and a face $K$ such that $\check{\mathfrak{X}} \to \Spec \kk[P_{\numer}]_{J_{\numer}}$ is the basechange of $\check{\mathfrak{X}} \to \Spec \kk[P]_{J}$ by fixing a well-chosen monoid $NE(\mathfrak{X}_0) \subseteq P' \subseteq A_1 (\mathfrak{X}_0, \mathbb{Z})$ with a face $K'$ and setting $P := P' \cap (P \oplus G)$ and $K := K' \cap (P \oplus G)$. So working over $A_1(\mathfrak{X}_0, \mathbb{Z}) = A_1(\mathfrak{X}_0, \mathbb{Z})_{\numer}$ is equivalent to studying the numerical locus of the extended intrinsic mirror $\check{\mathfrak{X}} \to \Spec \kk[P]_{J}$ for arbitrary $A_1(\mathfrak{X}_0, \mathbb{Z})$. 
\end{remark}

\subsubsection{A compatible choice of base monoid}  \label{basemonoidcompatible}

In Conjecture \ref{conjecture}, we use an arbitrary choice of a finitely generated saturated monoid $P$ satisfying $NE(\mathfrak{X}_0)_{\numer} \subseteq P \subseteq A_1(\mathfrak{X}_0, \mathbb{Z})_{\numer}$. For generalizations, we will use a fixed choice, compatible with $\pi: \mathfrak{X} \to \bar{\mathfrak{X}}$. In \cite[Appendix A.2]{GHS}, the authors construct a universal toric degeneration mirror to $\bar{\mathfrak{X}} \to \mathcal{S}$ using a (finitely generated, saturated, and sharp) \emph{universal monoid} $Q$. 

\begin{fact} \label{unitoric} 
The universal monoid $Q$ of \cite[Appendix A.2]{GHS} naturally satisfies $NE\left(\bar{\mathfrak{X}}_0\right)_{\numer} \subseteq Q \subseteq A_1\left(\bar{\mathfrak{X}}_0, \mathbb{Z}\right)_{\numer}$.
\end{fact}

\begin{proof}
We use the notation $\MPA(\bar{B}, Q)$ of \cite[Definition 1.2.5]{GHS} for the monoid of convex integral MPA functions on $\left(\bar{B}, \bar{\mathscr{P}} \right)$with values in $Q$ (without the requirement of Section \ref{mpafunc} that $\kappa_{\underline{\rho}} = \kappa_{\underline{\rho}'}=: \kappa_{\rho}$ for any slabs $\underline{\rho}, \underline{\rho}' \subseteq \rho \in \bar{\mathscr{P}}^{[n-1]}$). We also use the notation $\bMPA(\bar{B}, Q)$ for the monoid of the restricted convex integral MPA functions that are used in \cite{GS2, GS3} (see \cite[Example 1.2.8(1)]{GHS}). These functions satisfy $\kappa_{\underline{\rho}} = \kappa_{\underline{\rho}'}=: \kappa_{\rho}$ for any slabs $\underline{\rho}, \underline{\rho}' \subseteq \rho \in \bar{\mathscr{P}}^{[n-1]}$ as well as an additional balancing condition at each vertex. 

The monoid $Q$ is defined as $Q:=r(Q_0)^{\sat}$ where $Q_0 := \Hom (\MPA(\bar{B}, \mathbb{N}), \mathbb{N})$, $r$ is the surjective restriction map $$r: \Hom (\MPA(\bar{B}, \mathbb{N}), \mathbb{Z}) \to \Hom (\bMPA(\bar{B}, \mathbb{N}), \mathbb{Z}),$$ and we use the notation $M^{\sat}$ for the saturation of a monoid $M$. See \cite[Appendix A.2]{GHS} for details.

\textbf{Step 1.}
We first show that $\Nef\left(\bar{\mathfrak{X}}_0\right) \cong \bMPA(\bar{B}, \mathbb{N})$. Indeed, a nef line bundle $L$ on $\bar{\mathfrak{X}}_0$ restricts to a nef line bundle $L_{\bar{v}}$ on each irreducible component $\bar{D}_{\bar{v}}, \ \bar{v} \in \bar{\mathscr{P}}^{[0]}$ of $\bar{\mathfrak{X}}_0$. By standard toric geometry, every $L_{\bar{v}}, \ \bar{v} \in \bar{\mathscr{P}}^{[0]}$ defines a convex PL-function $\varphi_{\bar{v}}$ with integral slopes on the fan $\Sigma_{\bar{v}}$ of $\bar{D}_{\bar{v}}$. The $\varphi_{\bar{v}}$ patch to a global multi-valued convex integral $\MPA$ function $\varphi_L$ on $(\bar{B}, \bar{\mathscr{P}})$ with kinks $\kappa_{\underline{\rho}}(\varphi_L) = \kappa_{\underline{\rho}'}(\varphi_L) = \deg L|_{\bar{X}_{\rho}}$ for any slabs $\underline{\rho}, \underline{\rho}' \subseteq \rho \in \bar{\mathscr{P}}^{[1]}$. The fact that the $\bar{D}_{\bar{v}}, \ \bar{v} \in \bar{\mathscr{P}}^{[0]}$ are toric varieties implies that the balancing condition is satisfied. This defines a homomorphism $$\Nef\left(\bar{\mathfrak{X}}_0\right) \to \bMPA(\bar{B}, \mathbb{N}).$$ 

\begin{sloppypar}
Conversely, given a convex multi-valued integral PL-function $\psi$ with $\kappa_{\underline{\rho}}(\psi) = \kappa_{\underline{\rho}'}(\psi)$ for any slabs $\underline{\rho}, \underline{\rho}' \subseteq \rho \in \bar{\mathscr{P}}^{[1]}$ and satisfying a balancing condition, a piecewise-linear representative $\psi_{\bar{v}}$ at a vertex $\bar{v} \in \bar{\mathscr{P}}^{[0]}$ defines a line bundle $L_{\bar{v}}$ on $\bar{D}_{\bar{v}}$. These are isomorphic on the double curves $\bar{X}_{\rho}, \ \rho \in \bar{\mathscr{P}}^{[1]}$ since for $\rho$ of the form $\rho = \left \langle \bar{v}, \bar{v}' \right\rangle$ for $\bar{v}, \bar{v}' \in \bar{\mathscr{P}}^{[0]}$, we have $$\deg L_{\bar{v}} |_{\bar{X}_{\rho}} = \kappa_{\underline{\rho}}(\psi) = \kappa_{\underline{\rho}'}(\psi) = \deg L_{\bar{v}'} |_{\bar{X}_{\rho}}.$$ Choosing an isomorphism on the overlaps defines an obstruction class in $H^2(\bar{B}, \kk^{\times})$ that is precisely the $o(s)$ of \cite[Theorem 2.34]{GS2}. So it vanishes by Assumption \ref{techass}, and we get a line bundle $L$ on $\bar{\mathfrak{X}}_0$. This proves $\Nef\left(\bar{\mathfrak{X}}_0\right) \cong \bMPA(\bar{B}, \mathbb{N})$.
\end{sloppypar}

\textbf{Step 2.}
Note that we have a natural pairing 
\begin{equation} \label{nondegeneratepairing}
\deg: A_1\left(\bar{\mathfrak{X}}_0, \mathbb{Z}\right)_{\numer} \times \Pic\left(\bar{\mathfrak{X}}_0\right) \to \mathbb{Z}
\end{equation}
that is non-degenerate.\footnote{Indeed, it is enough to show that $\Pic\left(\bar{\mathfrak{X}}_0\right)$ is isomorphic to $ \Num\left(\bar{\mathfrak{X}}_0\right)$, the quotient of $\Pic\left(\bar{\mathfrak{X}}_0\right)$ by the numerically trivial line bundles. For a smooth K3-surface $X$, it is classical that $\Pic(X) \cong \Num(X)$ since the fact that $H^1(X, \mathcal{O}_X) = 0$ and the Riemann-Roch theorem for line bundles on surfaces imply that every numerically trivial line bundle is trivial, see, e.g. \cite[Chapter 1, Proposition 2.4]{Huy}. However, a similar proof applies in the case of $\bar{\mathfrak{X}}_0$ since $H^1(\bar{\mathfrak{X}}_0, \mathcal{O}_{\bar{\mathfrak{X}}_0}) \cong H^1(\bar{B}, \kk) = 0$, see the proof of Proposition \ref{spheredic}.} We define a map 
\begin{equation}
\begin{aligned}
i: NE\left(\bar{\mathfrak{X}}_0\right) &\to \Hom (\bMPA(\bar{B}, \mathbb{N}), \mathbb{Z}) \cong \Hom \left(\Nef\left(\bar{\mathfrak{X}}_0\right), \mathbb{Z}\right) \\
\beta &\mapsto (L \mapsto \deg_L(\beta)) 
\end{aligned}
\end{equation}
Non-degeneracy of the pairing \eqref{nondegeneratepairing} and the fact that for every $L \in \Nef\left(\bar{\mathfrak{X}}_0\right)$ and every $\beta \in NE\left(\bar{\mathfrak{X}}_0\right)$ we have $\deg_L(\beta) \geq 0$ imply that $i$ is injective with image contained in $r(Q_0)$. But $r(Q_0) \subseteq r(Q_0)^{\sat} = Q$ which shows the first inclusion.

\textbf{Step 3.} To define the inclusion $Q \subseteq A_1\left(\bar{\mathfrak{X}}_0, \mathbb{Z}\right)_{\numer}$, note that we have $Q \subseteq \Hom (\bMPA(\bar{B}, \mathbb{N}), \mathbb{Z})$ and 
\begin{multline*}
\Hom (\bMPA(\bar{B}, \mathbb{N}), \mathbb{Z})  \cong \Hom \left(\Nef\left(\bar{\mathfrak{X}}_0\right), \mathbb{Z}\right) \cong \\ \cong  \Hom \left(\Nef\left(\bar{\mathfrak{X}}_0\right)^{\gp}, \mathbb{Z}\right) = \Hom \left(\Pic\left(\bar{\mathfrak{X}}_0\right), \mathbb{Z}\right) \cong A_1\left(\bar{\mathfrak{X}}_0, \mathbb{Z}\right)_{\numer}
\end{multline*}
where the last isomorphism is again due to non-degeneracy of \eqref{nondegeneratepairing}.
\end{proof}

Consider the surjective pushforward map 
\begin{equation} \label{pushforward}
A_1(\mathfrak{X}_0, \mathbb{Z}) \overset{\pi_{*}}{\longrightarrow}  A_1\left(\bar{\mathfrak{X}}_0, \mathbb{Z}\right)
\end{equation}
and let $K^{\gp}:= \ker{\pi_{*}}$, so that there is a short exact sequence
\begin{equation} \label{exactseqcurveclasses}
0 \longrightarrow K^{\gp} \longrightarrow A_1(\mathfrak{X}_0, \mathbb{Z}) \overset{\pi_{*}}{\longrightarrow}  A_1\left(\bar{\mathfrak{X}}_0, \mathbb{Z}\right) \longrightarrow 0.
\end{equation}
Note that $K^{\gp}$ is generated by the classes of curves contracted by $\pi$. In particular, we have $E_{\rho,k} \in K^{\gp}$ and $F_{\rho, p} \in K^{\gp}$ for $E_{\rho,k}, \ \rho \in \bar{\mathscr{P}}^{[1]}, \ 1 \leq k \leq r_{\rho}$ and $F_{\rho, p}, \ \rho \in \bar{\mathscr{P}}^{[1]}, \ 1 \leq p \leq l_{\rho} -1$ the curve classes of \eqref{intnumbers}.

\begin{cor}
We have natural inclusions $NE(\mathfrak{X}_0)_{\numer} \subseteq K^{\gp} \oplus Q  \subseteq A_1(\mathfrak{X}_0, \mathbb{Z})_{\numer}$
for $Q$ the universal monoid of \cite[Appendix A.2]{GHS}.
\end{cor}

\begin{sloppypar}
\begin{proof}
Since $A_1\left(\bar{\mathfrak{X}}_0, \mathbb{Z}\right)$ is a finitely generated free abelian group, we have a splitting of $\pi_{*}$ in \eqref{exactseqcurveclasses} inducing an isomorphism $A_1(\mathfrak{X}_0, \mathbb{Z})_{\numer} \cong K^{\gp} \oplus A_1\left(\bar{\mathfrak{X}}_0, \mathbb{Z}\right)_{\numer}$. We have $NE(\mathfrak{X}_0)_{\numer} \subseteq K^{\gp} \oplus NE\left(\bar{\mathfrak{X}}_0\right)_{\numer}$ under this isomorphism and  $NE\left(\bar{\mathfrak{X}}_0\right)_{\numer} \subseteq Q \subseteq A_1\left(\bar{\mathfrak{X}}_0, \mathbb{Z}\right)_{\numer}$ by Proposition \ref{unitoric}. The claim follows from
\begin{multline}
NE(\mathfrak{X}_0)_{\numer} \subseteq K^{\gp} \oplus NE\left(\bar{\mathfrak{X}}_0\right)_{\numer} \subseteq  K^{\gp} \oplus Q  \subseteq \\ \subseteq K^{\gp} \oplus A_1\left(\bar{\mathfrak{X}}_0, \mathbb{Z}\right)_{\numer}  \cong A_1(\mathfrak{X}_0, \mathbb{Z})_{\numer}.
\end{multline}
\end{proof}
\end{sloppypar}

The intrinsic mirror is well-defined over $K^{\gp} \oplus Q$. Indeed, Proposition \ref{canextendifstronglyadmissible} and the proof of Proposition \ref{P2general} (via the proof of Proposition \ref{P1general}) imply that one can find a well-chosen monoid $P'$ and a face $K' \subseteq P'$ such that $K'^{\gp} = K^{\gp}$. But then the claim follows by applying Remark \ref{extmonoidchange} to the monoid $P' \cap (K^{\gp} \oplus Q)$ and the face $K' \cap (K^{\gp} \oplus Q)$. 

In this chapter, we shall work with the base monoid $K^{\gp} \oplus Q$. Note that $K^{\gp} \oplus Q$ does not satisfy condition (3) of Definition \ref{monoidass}, but Remark \ref{restrictingtonumer} still makes sense in this setup. 

\subsubsection{The minimal relative Gross-Siebert locus} \label{minrelgslocus}

The intrinsic mirror over $P:=K^{\gp} \oplus Q$ is of the form $$\check{\mathfrak{X}} \to \Spec \widehat{\kk [P]}_J = \Spec \kk [K^{\gp}] \llbracket Q \rrbracket$$ (where the notation means that the completion is only with respect to the second factor) for  $J := P \setminus K^{\gp} = Q \setminus \left\{0 \right\}$ (thinking of $Q \setminus \left\{0\right\}$ as an ideal of $P$). We still call it the \emph{extended intrinsic mirror}.

Note from \eqref{intnumbers} that we have $E^{\gp} \subseteq K^{\gp}$ for
\begin{equation} \label{egp}
E^{\gp} := \langle E_{\rho, k} \ , \ -E_{\rho,k} \ | \  \rho \in \bar{\mathscr{P}}^{[1]}, 1 \leq k \leq r_{\rho}  \rangle.
\end{equation}
It is clear from the construction that the curve classes $E_{\rho,k}$ for $\rho \in \bar{\mathscr{P}}^{[1]}, \ 1 \leq k \leq r_{\rho}$ don't have any relations with the other curve classes of $K^{\gp}$. On the other hand, in general, there will be relations between the $E_{\rho,k}$ so we have $\dim E^{\gp} \leq \sum_{\rho \in \bar{\mathscr{P}}^{[1]}} r_{\rho}$. One can check that both $\dim E^{\gp} = \sum_{\rho \in \bar{\mathscr{P}}^{[1]}} r_{\rho}$ and $\dim E^{\gp} < \sum_{\rho \in \bar{\mathscr{P}}^{[1]}} r_{\rho}$ occur.\footnote{For the equality, consider a special toric degeneration $\bar{\mathfrak{X}} \to \mathcal{S}$ with a simple dual intersection complex $\left( \bar{B}, \bar{\mathscr{P}}\right)$ and a strongly admissible small resolution $\pi: \mathfrak{X} \to \bar{\mathfrak{X}}$ to a log smooth $\mathfrak{X} \to \mathcal{S}$ with the divisor $D$ simple normal crossings. Simplicity of $\left( \bar{B}, \bar{\mathscr{P}}\right)$ implies that $r_{\rho} \leq 1$ for every $\rho \in \bar{\mathscr{P}}^{[1]}$. Suppose further that if $\rho_1, \rho_1 \in \bar{\mathscr{P}}^{[1]}$ are such that $r_{\rho_1} = r_{\rho_2} = 1$, then $\rho_1 \cap \rho_2 = \varnothing$. Then the equality $\dim E^{\gp} = \sum_{\rho \in \bar{\mathscr{P}}^{[1]}} r_{\rho}$ can be deduced by considering the intersection numbers of the curves $E_{\rho,1}$ (for $\rho \in \bar{\mathscr{P}}^{[1]}$ with $ r_{\rho} = 1$) with the irreducible components of $\mathfrak{X}_0$.

For the inequality, consider the toric degeneration $\bar{\mathfrak{X}} \to \mathcal{S}$ of Example \ref{runningexample} and its small log smooth resolution $\mathfrak{X} \to \mathcal{S}$. Then there are 24 curves $E_{\rho,k}$ of \eqref{24curves} corresponding to $6$ singularities $x_{\rho}, \ \rho \in \bar{\mathscr{P}}^{[1]}$ with $r_{\rho} = 4$ of the affine structure on $\left( \bar{B}, \bar{\mathscr{P}}\right)$. On the other hand, we have $H_2(\mathfrak{X}_0, \mathbb{Z}) \cong \mathbb{Z}^{22}$. This follows from an easy computation using a Mayer-Vietoris type spectral sequence for homology similar to the spectral sequence of \cite[Chapter VI, §1]{TTAG} (by David R. Morrison) for cohomology. So $\dim E^{\gp} \leq \dim H_2(\mathfrak{X}_0, \mathbb{Z}) = 22 < 24 = \sum_{\rho \in \bar{\mathscr{P}}^{[1]}} r_{\rho}$.}

Write $K^{\gp}$ as $K^{\gp} = E^{\gp} \oplus G^{\gp}$ for a finitely generated free abelian group $G^{\gp}$. In particular, $G^{\gp}$ contains the curve classes $F_{\rho, p}, \ \rho \in \bar{\mathscr{P}}^{[1]}, \ 1 \leq p \leq l_{\rho} -1$ of \eqref{intnumbers}. 

\begin{defn} \label{gslocusdef}
Define $h^{\GS}: K^{\gp} \oplus Q = E^{\gp} \oplus G^{\gp} \oplus Q \to E^{\gp} \oplus Q$ by sending $G^{\gp} \to 0$ and by the identity on $E^{\gp}$ and $Q$. We call the basechange $\check{\mathfrak{X}} \to \Spec \kk [E^{\gp}]\llbracket Q \rrbracket$ of the extended intrinsic mirror $\check{\mathfrak{X}} \to \Spec \kk [K^{\gp}] \llbracket Q \rrbracket$ by $h^{\GS}$ the \emph{(numerical) minimal relative Gross-Siebert locus}.\footnote{Note that this is defined via Remark \ref{restrictingtonumer} for any extended intrinsic mirror, even when using $A_1(\mathfrak{X}_0, \mathbb{Z}) \neq A_1(\mathfrak{X}_0, \mathbb{Z})_{\numer}$.} 
\end{defn}

\begin{notation}
We shall often write the ring $\kk [E^{\gp}]\llbracket Q \rrbracket$ as $\kk [t^{\pm E_{\rho, k} }]\llbracket Q \rrbracket$ using \eqref{egp}. 
\end{notation}

\begin{remark}
The \emph{Gross-Siebert locus} of the extended intrinsic mirror $\check{\mathfrak{X}} \to \Spec \widehat{\kk[P]}_{J}$ to a log CY surface $(\mathfrak{X}, D)$ was defined in \cite[Section 3.2]{GHK}. Our definition is similar in the case that $\pi: \mathfrak{X} \to \bar{\mathfrak{X}}$ is small. By \enquote{minimality}, we mean that the minimal relative Gross-Siebert locus only captures the essential curve classes $E^{\gp}$ of $K^{\gp}$ (i.e. the curve classes that would be present if we could construct a minimal resolution). See Section \ref{gslocusstuff} for the discussion on the strata of $\check{\mathfrak{X}} \to \Spec \kk [K^{\gp}] \llbracket Q \rrbracket$ that are not in the minimal relative Gross-Siebert locus. 
\end{remark}

We will show, by extending the correspondence of Conjecture \ref{conjecture} to larger and larger strata, that the minimal relative Gross-Siebert locus $\check{\mathfrak{X}} \to \Spec \kk [E^{\gp}]\llbracket Q \rrbracket$ of the extended intrinsic mirror can be recovered (up to a basechange) from  $\bar{\mathfrak{X}} \to \mathcal{S}$. First, we need to interpret the basechange by $h^{\GS}$ at the level of scattering diagrams. The following is a universal version of Construction \ref{ccc}.

\begin{cons} \label{gsscattering}
We define a collection of scattering diagrams $$\mathfrak{D}^{\GS} = \left\{ \mathfrak{D}^{\GS}_{Q^{k+1}}, \ k \geq 0 \right\}$$ on $(B, \mathscr{P})$ with monoid $E^{\gp} \oplus Q$ and $I_0 = (E^{\gp} \oplus Q) \setminus E^{\gp} = Q \setminus \left\{0\right\}$ as follows. We let $\mathfrak{D}^{\GS}_{Q}:= \mathfrak{D}^{\GS}_{I_0} $ have walls $$ \left(\mathfrak{b}_{\rho_p}, h^{\GS}(f_{\mathfrak{b}_{\rho_p}}) \right)$$ for $\mathfrak{b}_{\rho_p}$ a wall of $\mathfrak{D}_J$ defined in Construction \ref{newwalls}. We still write $\mathfrak{b}_{\rho_p}$ for these walls. 
Note that we have 
\begin{equation} 
\begin{aligned}
h^{\GS}(f_{\mathfrak{b}_{\rho_p}}) &= \prod_{1 \leq k \leq r_{\rho}} \left(1+t^{-E_{\rho,k}} w_{\rho}^{-1}\right) , \quad  &1 \leq p \leq p_0 \\
h^{\GS}(f_{\mathfrak{b}_{\rho_p}}) &= \prod_{1 \leq k \leq r_{\rho}}\left(1+t^{E_{\rho,k}} w_{\rho}\right), \quad &p_0+1 \leq p \leq l_{\rho} 
\end{aligned}
\end{equation}
since $h^{\GS}(F_{\rho,p}) = 0$ for all $\rho \in \bar{\mathscr{P}}^{[1]}, \ 1 \leq p \leq l_{\rho} -1$. 

More generally, we let $\mathfrak{D}^{\GS}_{Q^{k+1}}:= \mathfrak{D}^{\GS}_{I_0^{k+1}}$ have walls $$(\mathfrak{p}, h^{\GS}(f_{\mathfrak{p}}))$$ for $\mathfrak{p}$ a wall of $\mathfrak{D}_{J^{k+1}}$ of Construction \ref{extendcanonicalscat}. Again, we still denote such a wall by $\mathfrak{p}$. 

\begin{sloppypar}
We define the MPA function $\varphi^{\GS}$ on $(B, \mathscr{P}_{\coar})$ via its kinks by setting $\kappa_{\rho}(\varphi^{\GS}) := h^{\GS}(\kappa_{\rho}(\varphi))$. Note that this is clearly compatible with the basechange and that $\kappa_{\rho_p}(\varphi^{\GS}) = \kappa_{\rho_p}(\varphi)$ for an edge $\rho_p \in \mathscr{P}^{[1]}$ subdividing $\rho \in \mathscr{P}^{[1]}_{\coar}$ and $\kappa_{\rho'}(\varphi^{\GS})=0$ for $\rho' \in \mathscr{P}^{[1]} \setminus \mathscr{P}^{[1]}_{\coar}$.
\end{sloppypar}

Now $\mathfrak{D}^{\GS}_{Q^{k+1}}$ is a consistent scattering diagram since consistency in codimensions $0$, $1$ and $2$ follows trivially from the corresponding consistency of $\mathfrak{D}_{J^{k+1}}$ by interpreting all the relevant monomials as elements of $\kk[E^{\gp}][Q]$ via the basechange. By taking the inverse limit over $\check{\mathfrak{X}}_{\mathfrak{D}^{\GS}_{Q^{k+1}}} \to \Spec \kk [E^{\gp}][Q]/I_0^{k+1}$ for $k \geq 0$ (using the MPA function $\varphi^{\GS}$ in the construction) we define a family $\check{\mathfrak{X}}_{\mathfrak{D}^{\GS}} \to  \Spec \kk [E^{\gp}]\llbracket Q \rrbracket$. It is evident from the construction that this is the same family as the minimal relative Gross-Siebert locus $\check{\mathfrak{X}} \to \Spec \kk [E^{\gp}]\llbracket Q \rrbracket$.
\end{cons}

\begin{remark} \label{remring1}
One can equivalently view $\mathfrak{D}^{GS}$ as a collection of scattering diagrams with $A=\kk [E^{\gp}]$, monoid $Q$, and $I_0=\mathfrak{m}$. 
\end{remark}

We are ready to start generalizing Conjecture \ref{conjecture}.

\subsection{Universal version} \label{universal}

First, we would like to remove the dependency on the choice of polarization $A$ on $\bar{\mathfrak{X}} \to \mathcal{S}$. We use a version of the toric degeneration mirror defined over the universal monoid $Q$. 

As explained in \cite[Appendix A.2]{GHS}, \cite[Theorem 5.2]{GS1} implies that one can do a universal version of the reconstruction algorithm of Theorem \ref{reconstruction}. Namely, given a choice of the initial slab functions $$\left\{f_{\underline{\rho}} \in \kk[\Lambda_{\rho}] \ \Big| \ \underline{\rho} \in \tilde{\bar{\mathscr{P}}}^{[n-1]} \right\}$$ as in Proposition \ref{struct}, one can produce a collection of scattering diagrams $\bar{\mathfrak{D}} := \left\{ \bar{\mathfrak{D}}_{k}, \ k \geq 0 \right\}$ with $\bar{\mathfrak{D}}_{k}:= \bar{\mathfrak{D}}_{\mathfrak{m}^{k+1}}$  using monoid $Q$, the maximal ideal $ \mathfrak{m}= Q \setminus \left\{0\right\}$, and the universal MPA function $\bar{\varphi}$\footnote{This function is denoted by $\breve{\varphi}$ in \cite[Appendix A.2]{GHS} but $\bar{\varphi}$ fits better with our notation. One needs to have some polarization $A$ on $\bar{\mathfrak{X}} \to \mathcal{S}$ to define $\bar{\varphi}$ (which is why in this chapter we require that some such polarization exists).} with values in $Q$ (constructed in \cite[Appendix A.2]{GHS}). Moreover, $\bar{\mathfrak{D}}$ satisfies the compatibility and uniqueness properties of Theorem \ref{reconstruction}. We let the initial slab functions be given by \eqref{choiceslab} as before (with the convention for slabs $\underline{\rho}, \underline{\rho'} \subseteq \rho \in \bar{\mathscr{P}}^{[1]}$ as in Sections \ref{admissibleandstronglyso} and \ref{global42}).

The restriction of the universal toric degeneration mirror of
\cite[Theorem A.2.4]{GHS} to the trivial gluing data is the family $\check{\bar{\mathfrak{X}}}_{\bar{\mathfrak{D}}} \to \Spec \kk \llbracket Q \rrbracket$ defined as the inverse limit over $\check{\bar{\mathfrak{X}}}_{\bar{\mathfrak{D}}} \to \Spec \kk [Q]/ \mathfrak{m}^{k+1}$ for $k \geq 0$.

\begin{remarks} \label{univrem1}
\begin{enumerate}
    \item Note a slight abuse of notation in the definition of $\bar{\mathfrak{D}}$. In this chapter, we will denote by $\bar{\mathfrak{D}}$ any collection of scattering diagrams defining the toric degeneration mirror of interest, regardless of the exact situation (i.e. of the ring $A$, monoid $Q$, and ideal $I_0$).
    \item Similarly, regardless of the exact situation, we will denote by $\mathfrak{D}$ the collection of scattering diagrams defining the basechange $\check{\mathfrak{X}}_{\mathfrak{D}}$ of the intrinsic mirror family $\check{\mathfrak{X}}$, for which we compare $\check{\mathfrak{X}}_{\mathfrak{D}}$ and $\check{\bar{\mathfrak{X}}}_{\bar{\mathfrak{D}}}$.
    \item The family $\check{\bar{\mathfrak{X}}}_{\bar{\mathfrak{D}}} \to \Spec \kk \llbracket Q \rrbracket$ is universal in the sense that every toric degeneration mirror family $\check{\bar{\mathfrak{X}}} \to \Spec \kk \llbracket t \rrbracket$ constructed using trivial gluing data, the initial slab functions of \eqref{choiceslab}, and any polarization $A$ on $\bar{\mathfrak{X}} \to \mathcal{S}$ can be obtained from $\check{\bar{\mathfrak{X}}}_{\bar{\mathfrak{D}}} \to \Spec \kk \llbracket Q \rrbracket$ via basechange by $Q \to \mathbb{N}, \ \beta \mapsto A \cdot \beta$. Note that the composition of this map with $h^{\GS}$ is exactly the map $h$ of Conjecture \ref{conjecture}.  
\end{enumerate}
\end{remarks}

Instead of taking the basechange of the extended intrinsic mirror by $h: P \to \mathbb{N}, \ \beta \mapsto \pi^{*}A \cdot \beta$ as in Conjecture \ref{conjecture}, we consider the basechange of the minimal relative Gross-Siebert locus by $h: E^{\gp} \oplus Q \to Q$ that sends $E^{\gp} \to 0$ and is the identity on $Q$ (note that the composition $h \circ h^{GS}$ sends $K^{\gp} \to 0$ and is the identity on $Q$). At the level of rings, this corresponds to restricting to $t^{E_{\rho, k}} =1$ for all $1 \leq k \leq r_{\rho}, \ \rho \in \bar{\mathscr{P}}^{[1]}$. We need to define the appropriate version of $\mathfrak{D}$ giving rise to the basechanged family. The following (combined with Construction \ref{gsscattering}) is the analogue of Construction \ref{ccc}.

\begin{cons} \label{univ1}
We define a collection of scattering diagrams $\mathfrak{D} = \left\{ \mathfrak{D}_k, \ k \geq 0 \right\}$ on $(B, \mathscr{P})$ with monoid $Q$ and $I_0 = \mathfrak{m}$ as follows. We let $\mathfrak{D}_k := \mathfrak{D}_{\mathfrak{m}^{k+1}}$ have walls $$(\mathfrak{p}, h(f_{\mathfrak{p}}))$$ for $\mathfrak{p}$ a wall of $\mathfrak{D}^{\GS}_{Q^{k+1}}$. We still denote such a wall by $\mathfrak{p}$. We use the MPA function $\varphi$ with the same kinks as $\varphi^{\GS}$ (seen as elements of $Q$), this is clearly compatible with the basechange. 

As in Constructions \ref{ccc} and \ref{gsscattering}, $\mathfrak{D}_k$ is a consistent scattering diagram and taking the inverse limit over $\check{\mathfrak{X}}_{\mathfrak{D}_k} \to \Spec \kk [Q]/\mathfrak{m}^{k+1}$ for $k \geq 0$ we get a family $\check{\mathfrak{X}}_{\mathfrak{D}} \to \Spec \kk \llbracket Q \rrbracket$, which is the same family as the basechange of the minimal relative Gross-Siebert locus $\check{\mathfrak{X}} \to \Spec \kk [E^{\gp}]\llbracket Q \rrbracket$ by $h: E^{\gp} \oplus Q \to Q$.
\end{cons}

\begin{cons} \label{univ2}
We use the analogue of Construction \ref{refinement} for a scattering diagram $\bar{\mathfrak{D}}'$ on $\left(\bar{B}, \bar{\mathscr{P}}'\right)$ equivalent to $\bar{\mathfrak{D}}$, and define the MPA function $\bar{\varphi}'$ on $\left(\bar{B}, \bar{\mathscr{P}}' \right)$ from $\bar{\varphi}$ in the same way as the MPA function $\varphi'_A$ is defined from $\varphi_A$ in Construction \ref{refinement}. 

We define $\Phi(\mathfrak{D}_k)'$ and $\Phi(\mathfrak{D}_k)$ exactly as in Construction \ref{chs}. It is easy to see that $\kappa_{\rho'} (\varphi) = \kappa_{\Phi(\rho')} (\bar{\varphi}')$ for any $\rho' \in \bar{\mathscr{P}}'^{[1]}$. So the constructions are compatible with the MPA functions. Moreover, $\Phi(\mathfrak{D}_0)$ is still combinatorially equivalent to $\bar{\mathfrak{D}}_0$.
\end{cons}

\begin{fact} \label{proofuniversal}
The basechange of the minimal relative Gross-Siebert locus $\check{\mathfrak{X}} \to \Spec \kk [E^{\gp}]\llbracket Q \rrbracket$ by $h: E^{\gp} \oplus Q \to Q$ is isomorphic to the restriction $\check{\bar{\mathfrak{X}}}_{\bar{\mathfrak{D}}} \to \Spec \kk \llbracket Q \rrbracket$ of the universal toric degeneration mirror of
\cite[Theorem A.2.4]{GHS} to the trivial gluing data.
\end{fact}
	
\begin{proof}
Consistency of $\Phi(\mathfrak{D}_k)$ and equivalence of $\Phi(\mathfrak{D}_k)$ and $\mathfrak{D}_k$ follow exactly as in the proofs of Propositions \ref{consis} and \ref{equivalence} respectively since the proofs do not depend on the base monoid. Similarly, changing $\mathbb{N}$ to $Q$, $(t^{k+1})$ to $\mathfrak{m}^{k+1}$, and using $h: E^{\gp} \oplus Q \to Q$, the same argument as in the proof of Theorem \ref{proofconj} implies the result. 
\end{proof}

Viewing $\Spec \kk [t^{\pm E_{\rho, k} }]$ as the parameter space, we have interpreted the fibre over the point $$\left\{t^{E_{\rho, k}} =1, \ \rho \in \bar{\mathscr{P}}^{[1]}, \ 1 \leq k \leq r_{\rho} \right\}$$ as a (universal in polarization) toric degeneration. We will now extend this correspondence to the restriction of the extended intrinsic mirror to a family over a certain subvariety of $\Spec \kk [t^{\pm E_{\rho, k} }]$.
	
\subsection{Extension in the free parameters} \label{freeparam}

Recall from Proposition \ref{struct} that if $\left(\bar B, \bar{\mathscr{P}}\right)$ is not simple, there are coefficients $a_{\rho, i} \in \kk$ for $ 1 \leq i \leq r_{\rho} - 1, \ \rho \in \bar{\mathscr{P}}^{[1]}$ in the initial slab functions $f_{\underline{\rho}}, f_{\underline{\rho}'}$ that we are free to choose when constructing the toric degeneration mirror. So far, we have always made the choice specified by \eqref{choiceslab}. We can make any other choice of the initial slab functions $f_{\underline{\rho}}, f_{\underline{\rho}'}$ satisfying \eqref{slabfunctions}. Proving the analogues of Theorem \ref{proofconj} and Proposition \ref{proofuniversal} is straightforward. However, we can in fact do better. As shown in \cite[Appendix A.4]{GHS}, in the case of toric degenerations of K3-s, one can make the toric degeneration mirror universal in the choice of slab functions.

Indeed, let $K := \sum_{\rho \in \bar{\mathscr{P}}^{[1]}, r_{\rho} > 0} (r_{\rho} -1)$ and consider $A = \kk[\mathbb{N}^{K}] =: \kk [a_{\rho, i}]$ where we now think of $a_{\rho, i}$ as variables. Then the initial slab functions 
\begin{equation} \label{extslabfunc}
\begin{aligned}f_{\underline{\rho}} = &1 + a_{\rho, 1}w_{\rho} + \dots + a_{\rho, r_{\rho}-1}w_{\rho}^{r_{\rho}-1} + w_{\rho}^{r_{\rho}}& \\ f_{\underline{\rho}'} = &1 + a_{\rho, r_{\rho}-1}w^{-1} + \dots + a_{\rho,1}w_{\rho}^{-r_{\rho}+1} + w_{\rho}^{-r_{\rho}}&
\end{aligned}
\end{equation} of \eqref{slabfunctions} 
make sense as elements of $A[\Lambda_{\rho}]$. Again, \cite[Theorem 5.2]{GS1} implies that the reconstruction algorithm of Theorem \ref{reconstruction} generalizes to this setting  
and we can produce a collection of scattering diagrams $\bar{\mathfrak{D}} := \left\{ \bar{\mathfrak{D}}_{k}, \ k \geq 0 \right\}$ with $\bar{\mathfrak{D}}_{k}:= \bar{\mathfrak{D}}_{\mathfrak{m}^{k+1}}$  using the ring $A$ and the same monoid, ideal and MPA function as in Section \ref{universal}.  Moreover, $\bar{\mathfrak{D}}$ satisfies the compatibility and uniqueness properties of Theorem \ref{reconstruction}.

The restriction of the universal toric degeneration mirror of
\cite[Theorem A.4.2]{GHS} to the trivial gluing data is the family $\check{\bar{\mathfrak{X}}}_{\bar{\mathfrak{D}}} \to \Spec \kk [a_{\rho, i}]\llbracket Q \rrbracket$ defined as the inverse limit over $\check{\bar{\mathfrak{X}}}_{\bar{\mathfrak{D}}} \to \Spec \kk [a_{\rho, i}][Q]/ \mathfrak{m}^{k+1}$ for $k \geq 0$.

\begin{remark} \label{freeremark1}
Similarly to Remark \ref{univrem1}(3), the family $\check{\bar{\mathfrak{X}}}_{\bar{\mathfrak{D}}} \to \Spec \kk [a_{\rho, i}]\llbracket Q \rrbracket$ is universal in the sense that every one-parameter toric degeneration mirror family $\check{\bar{\mathfrak{X}}} \to \Spec \kk \llbracket t \rrbracket$ constructed using the trivial gluing data and polarization $A$ can be obtained from $\check{\bar{\mathfrak{X}}}_{\bar{\mathfrak{D}}} \to \Spec \kk [a_{\rho, i}]\llbracket Q \rrbracket$ by a basechange that is induced by $Q \to \mathbb{N}, \ \beta \mapsto A \cdot \beta$ and by sending the variables $a_{\rho, i}$ to the values chosen for the initial slab functions $f_{\underline\rho}, f_{\underline{\rho}'}$ satisfying \eqref{slabfunctions}.
\end{remark}

The slab functions $f_{\mathfrak{b}_{\rho_p}}$ for the walls $\mathfrak{b}_{\rho_p}$ of $\mathfrak{D}_Q^{GS}$ are of the form 
\begin{equation} 
\begin{aligned}
f_{\mathfrak{b}_{\rho_p}} &= \prod_{1 \leq k \leq r_{\rho}} \left(1+t^{-E_{\rho,k}} w_{\rho}^{-1}\right) = 1+ \sigma_1(t^{-E_{\rho}})w_{\rho} + \dots + \sigma_{r_{\rho}-1} (t^{-E_{\rho}})w_{\rho}^{r-1} + \sigma_{r_{\rho}} (t^{-E_{\rho}})w_{\rho}^r  \\
f_{\mathfrak{b}_{\rho_p}} &= \prod_{1 \leq k \leq r_{\rho}}\left(1+t^{E_{\rho,k}} w_{\rho}\right) = 1+ \sigma_1(t^{E_{\rho}})w_{\rho} + \dots + \sigma_{r_{\rho}-1} (t^{E_{\rho}})w_{\rho}^{r-1} + \sigma_{r_{\rho}} (t^{E_{\rho}})w_{\rho}^r 
\end{aligned}
\end{equation}
where $\sigma_i(t^{-E_{\rho}})$ is the $i$-th symmetric polynomial in $t^{-E_{\rho, k}}, \ 1 \leq k \leq r_{\rho}$ and $\sigma_i(t^{E_{\rho}})$ is the $i$-th symmetric polynomial in $t^{E_{\rho, k}}, \ 1 \leq k \leq r_{\rho}$. Since the initial slab functions $f_{\underline{\rho}}$ of \eqref{slabfunctions} have coefficient $1$ at $w_{\rho}^{r_{\rho}}$ and we want to construct a collection of scattering diagrams $\mathfrak{D}:= \left\{ \mathfrak{D}_{k}, \ k \geq 0 \right\}$ with $\Phi(\mathfrak{D}_0)$ combinatorially equivalent to $\bar{\mathfrak{D}}_0$, we need to restrict to the family over the locus where $\sigma_{r_{\rho}} (t^{E_{\rho}}) = \prod_{k=1}^{r_{\rho}} t^{E_{\rho, k}} = 1$ for all $\rho \in \bar{\mathscr{P}}^{[1]}$ (clearly, this forces $\sigma_{r_{\rho}} (t^{-E_{\rho}})=1$ for all $\rho \in \bar{\mathscr{P}}^{[1]}$ as well).

Let $E^{\gp}_1:= E^{\gp} / \left \langle \sum_{k=1}^{r_{\rho}} E_{\rho, k}, \ \rho \in \bar{\mathscr{P}}^{[1]} \right \rangle$ and consider the basechange of the minimal relative Gross-Siebert locus by $$h: E^{\gp} \oplus Q \to E^{\gp}_1 \oplus Q$$ that is the natural projection on $E^{\gp} \to E_1^{\gp}$ and the identity on $Q$. At the level of rings, this corresponds to the natural projection
\begin{equation} \label{projmap1}
h: \kk [t^{\pm E_{\rho, k} }]\llbracket Q \rrbracket \to \kk [t^{\pm E_{\rho, k} }]\llbracket Q \rrbracket \Bigg/ \left \langle\prod_{k=1}^{r_{\rho}} t^{E_{\rho, k}} = 1, \ \rho \in \bar{\mathscr{P}}^{[1]} \right \rangle.
\end{equation}

The scattering diagram interpretation of this basechange is similar to Construction \ref{univ1}.

\begin{cons} \label{freeinitdiag}
We define a collection of scattering diagrams $\mathfrak{D} = \left\{ \mathfrak{D}_k, \ k \geq 0 \right\}$ on $(B, \mathscr{P})$ with monoid $E_1^{\gp} \oplus Q$ and $I_0 = (E_1^{\gp} \oplus Q) \setminus E_1^{\gp} = Q \setminus \left\{0\right\}$ as follows. We let $\mathfrak{D}_k := \mathfrak{D}_{I_0^{k+1}}$ have walls $$(\mathfrak{p}, h(f_{\mathfrak{p}}))$$ for $\mathfrak{p}$ a wall of $\mathfrak{D}^{\GS}_{Q^{k+1}}$. We still denote such a wall by $\mathfrak{p}$. We use the MPA function $\varphi$ with the same kinks as $\varphi^{\GS}$, this is clearly compatible with the basechange. 

Taking the inverse limit over $\check{\mathfrak{X}}_{\mathfrak{D}_k} \to \Spec \kk [E_1^{\gp}][Q]/I_0^{k+1}$ for $k \geq 0$ we get a family $\check{\mathfrak{X}}_{\mathfrak{D}} \to \Spec \kk [E_1^{\gp}]\llbracket Q \rrbracket$, which is the same family as the basechange of the minimal relative Gross-Siebert locus $\check{\mathfrak{X}} \to \Spec \kk [E^{\gp}]\llbracket Q \rrbracket$ by $h: E^{\gp} \oplus Q \to E^{\gp}_1 \oplus Q$.
\end{cons}

\begin{remark} \label{remring2}
As in Remark \ref{remring1}, we can view $\mathfrak{D}$ as a collection of scattering diagrams with $A=\kk [E^{\gp}_1]$, monoid $Q$, and $I_0=\mathfrak{m}$.
\end{remark}

The slab functions $f_{\mathfrak{b}_{\rho_p}}$ for the walls $\mathfrak{b}_{\rho_p}$ of $\mathfrak{D}_0$ are of the form 
\begin{equation} 
\begin{aligned}
f_{\mathfrak{b}_{\rho_p}} &= \prod_{1 \leq k \leq r_{\rho}} \left(1+t^{-E_{\rho,k}} w_{\rho}^{-1}\right) = 1+ \sigma_1(t^{-E_{\rho}})w_{\rho} + \dots + \sigma_{r_{\rho}-1} (t^{-E_{\rho}})w_{\rho}^{r-1} + w_{\rho}^r  \\
f_{\mathfrak{b}_{\rho_p}} &= \prod_{1 \leq k \leq r_{\rho}}\left(1+t^{E_{\rho,k}} w_{\rho}\right)= 1+ \sigma_1(t^{E_{\rho}})w_{\rho} + \dots + \sigma_{r_{\rho}-1} (t^{E_{\rho}})w_{\rho}^{r-1} + w_{\rho}^r 
\end{aligned}
\end{equation}
where we still write $t^{-E_{\rho, k}}$ and $t^{E_{\rho, k}}$ for the images of $t^{-E_{\rho, k}}$ and $t^{E_{\rho, k}}$ under $h$. 

We need to interpret the variables $a_{\rho, i}$ as elements of $\kk[E_1^{\gp}]$ to set up the correspondence. This requires a basechange of $\check{\bar{\mathfrak{X}}}_{\bar{\mathfrak{D}}} \to \Spec \kk [a_{\rho, i}]\llbracket Q \rrbracket$.  

\begin{cons} \label{freecons2}
Consider the ring map
\begin{equation} \label{sigmabasechange}
\sigma: \kk [a_{\rho, i}] \to \kk [E_1^{\gp}], \quad  a_{\rho, i} \mapsto \sigma_i (t^{E_{\rho}}), \ 1 \leq i \leq r_{\rho} -1, \ \rho \in \bar{\mathscr{P}}^{[1]}.
\end{equation}

We can define a new collection of scattering diagrams $\bar{\mathfrak{D}}^{\sigma} := \left\{ \bar{\mathfrak{D}}^{\sigma}_{k}, \ k \geq 0 \right\}$ with ring $A=\kk[E_1^{\gp}]$, monoid $Q$, and ideal $I_0 = \mathfrak{m}$ as follows. We let $\bar{\mathfrak{D}}^{\sigma}_{k}$ have walls $$(\mathfrak{p}, \sigma(f_{\mathfrak{p}}))$$ for $\mathfrak{p}$ a wall of $\bar{\mathfrak{D}}_k$. We also use the same MPA function $\bar{\varphi}$. 

Taking the inverse limit over $\check{\bar{\mathfrak{X}}}_{\bar{\mathfrak{D}}^{\sigma}_k} \to \Spec \kk [E_1^{\gp}][Q]/\mathfrak{m}^{k+1}$ for $k \geq 0$ we get a family $\check{\bar{\mathfrak{X}}}_{\bar{\mathfrak{D}}^{\sigma}} \to \Spec \kk [E_1^{\gp}]\llbracket Q \rrbracket$, which is the same family as the basechange of $\check{\bar{\mathfrak{X}}}_{\bar{\mathfrak{D}}} \to \Spec \kk [a_{\rho, i}]\llbracket Q \rrbracket$ by $\sigma$. If $\dim E^{\gp} = \sum_{\rho \in \bar{\mathscr{P}}^{[1]}} r_{\rho}$, then the induced map $\check{\bar{\mathfrak{X}}}_{\bar{\mathfrak{D}}^{\sigma}} \to \check{\bar{\mathfrak{X}}}_{\bar{\mathfrak{D}}}$ is generically a covering of index $\prod_{\rho \in \bar{\mathscr{P}}^{[1]}} r_{\rho}!$ ramified at the points where $t^{E_{\rho, i}} = t^{E_{\rho, j}}$ for some $\rho \in \bar{\mathscr{P}}^{[1]}$ and $1 \leq i < j \leq r_{\rho}$. If $\dim E^{\gp} < \sum_{\rho \in \bar{\mathscr{P}}^{[1]}} r_{\rho}$, then the morphism $\check{\bar{\mathfrak{X}}}_{\bar{\mathfrak{D}}^{\sigma}} \to \check{\bar{\mathfrak{X}}}_{\bar{\mathfrak{D}}}$ is not surjective.
\end{cons}

Now that the collections of scattering diagrams $\mathfrak{D}$ and $\bar{\mathfrak{D}}^{\sigma}$ are defined over the same ring (via Remark \ref{remring2} for $\mathfrak{D}$). We define $\Phi(\mathfrak{D}_k)'$ and $\Phi(\mathfrak{D}_k)$ as before.

\begin{cons} \label{freecons1}
As in Construction \ref{univ2}, we use the analogue of Construction \ref{refinement} for a scattering diagram $\bar{\mathfrak{D}}'$ on $\left(\bar{B}, \bar{\mathscr{P}}'\right)$.
We define $\Phi(\mathfrak{D}_k)'$ and $\Phi(\mathfrak{D}_k)$ exactly as in Construction \ref{chs}. This is compatible with the MPA functions and $\Phi(\mathfrak{D}_0)$ is combinatorially equivalent to $\bar{\mathfrak{D}}^{\sigma}_0$.
\end{cons}

\begin{fact} \label{prooffree}
The basechange of the minimal relative Gross-Siebert locus $\check{\mathfrak{X}} \to \Spec \kk [E^{\gp}]\llbracket Q \rrbracket$ by $h: E^{\gp} \oplus Q \to E^{\gp} \oplus Q$ is isomorphic to the basechange by $\sigma: \kk [a_{\rho, i}] \to \kk [E_1^{\gp}], \  a_{\rho, i} \mapsto \sigma_i (t^{E_{\rho}})$ of the restriction $\check{\bar{\mathfrak{X}}}_{\bar{\mathfrak{D}}} \to \Spec \kk [a_{\rho, i}]\llbracket Q \rrbracket$ of the universal toric degeneration mirror of
\cite[Theorem A.4.2]{GHS} to the trivial gluing data.
\end{fact}

\begin{proof}
We need to show that the families $\check{\mathfrak{X}}_{\mathfrak{D}} \to \Spec \kk [E_1^{\gp}]\llbracket Q \rrbracket$ and $\check{\bar{\mathfrak{X}}}_{\bar{\mathfrak{D}}^{\sigma}} \to \Spec \kk [E_1^{\gp}]\llbracket Q \rrbracket$ are isomorphic.
As in the proof of Proposition \ref{proofuniversal}, consistency of $\Phi(\mathfrak{D}_k)$ and equivalence of $\Phi(\mathfrak{D}_k)$ and $\mathfrak{D}_k$ follow exactly as in the proofs of Propositions \ref{consis} and \ref{equivalence} respectively. Similarly, changing $\bar{\mathfrak{D}}$ to $\bar{\mathfrak{D}}^{\sigma}$, $\kk$ to $A=\kk[E_1^{\gp}]$, $\mathbb{N}$ to $Q$, $(t^{k+1})$ to $\mathfrak{m}^{k+1}$, and using $h: E^{\gp} \oplus Q \to E_1^{\gp} \oplus Q$, the same argument as in the proof of Theorem \ref{proofconj} implies the result.
\end{proof}

\begin{observe}
If $\left(\bar B, \bar{\mathscr{P}}\right)$ is simple, Proposition \ref{prooffree} reduces to Proposition \ref{proofuniversal}.
\end{observe}

Viewing $\Spec \kk [t^{\pm E_{\rho, k} }]$ as the parameter space, we have interpreted the restriction of the extended intrinsic mirror to the family over the subvariety $$\left\{\prod_{k=1}^{r_{\rho}} t^{E_{\rho, k}} = 1, \ \rho \in \bar{\mathscr{P}}^{[1]} \right\}$$ as a natural basechange of a (universal in polarization and slab functions) toric degeneration. It remains to understand the fibres over the points with $\prod_{k=1}^{r_{\rho}} t^{E_{\rho, k}} \neq 1$ for some $\rho \in \bar{\mathscr{P}}^{[1]}$, which requires introducing gluing data for toric degenerations.

\subsection{Extension over non-normalized fibres via gluing data} \label{glue}

Let $c_{\rho} \in \kk^{\times}$ for $\rho \in \bar{\mathscr{P}}^{[1]}$ be fixed constants and suppose that we tried to do the same construction as in Section \ref{freeparam} but instead of restricting to $\prod_{k=1}^{r_{\rho}} t^{E_{\rho, k}} = 1, \ \rho \in \bar{\mathscr{P}}^{[1]}$ we restricted to $\prod_{k=1}^{r_{\rho}} t^{E_{\rho, k}} = c_{\rho}, \ \rho \in \bar{\mathscr{P}}^{[1]}$ (we assume that $c_{\rho}, \ \rho \in \bar{\mathscr{P}}^{[1]}$ are chosen so that this system of equations is consistent\footnote{Recall that one can have relations between the $E_{\rho, k}$, see the discussion after \eqref{egp}.}). Then for a wall $\mathfrak{b}_{\rho_p}$ on $(B, \mathscr{P})$, the constant coefficient of $z^{m_{\underline{\Phi(\rho)}'\underline{\Phi(\rho)} }}f_{\Phi(\mathfrak{b}_{\rho_p})}$ is $c_\rho$. So unless $c_{\rho} = 1$ for all $\rho \in \bar{\mathscr{P}}^{[1]}$, Construction \ref{freecons1} would give some slabs $\mathfrak{b}$ of $\Phi(\mathfrak{D}_0)$ with the constant coefficient of $f_{\mathfrak{b}}$ not equal to $1$. However, recall from Section \ref{tdmpafunc} that the initial slab functions for a toric degeneration are assumed to be normalized. As explained in \cite[Definition 4.23]{GS2} (in a rather different language), one can always normalize the slab functions by introducing gluing data.

We will give an overview of gluing data for toric degenerations and prove the analogue of Proposition \ref{prooffree} for the restriction of the extended intrinsic mirror to the family over the subvariety $$\left\{\prod_{k=1}^{r_{\rho}} t^{E_{\rho, k}} = c_{\rho}, \ \rho \in \bar{\mathscr{P}}^{[1]} \ | \ c_{\rho} \in \kk^{\times}\right\}.$$ Finally, in Theorem \ref{compgen}, we do the construction universally in a subfamily of gluing data and extend the correspondence to the whole minimal relative Gross-Siebert locus. 

It is important to emphasize that we keep the gluing data on the intrinsic mirror side trivial in this chapter. See Section \ref{gdint} for a discussion on introducing gluing data for intrinsic mirrors.

\subsubsection{Overview of gluing data} \label{gluingdataintro}

We review gluing data in our case of interest following \cite[Section 5.2]{GHS}, see loc. cit. for more details. The treatment of gluing data in \cite{GS2, GS1} is rather different, but any gluing data in the sense of \cite[Section 5.2]{GHS} induces gluing data in the sense of \cite{GS2, GS1} (see \cite[Appendix A.1]{GHS}).

Let $\left(\bar B, \bar{\mathscr{P}}\right)$ be the dual intersection complex of a toric degeneration of K3-s. We work over an arbitrary ring $A$, use monoid $Q$, ideal $I_0$, and fix an ideal $I$ with $\sqrt{I} = I_0$, as in Sections \ref{mpafunc} and \ref{monomialsaffstr}. 

\begin{defn} \label{gddef}
For every $\rho \in \bar{\mathscr{P}}^{[1]}$ with two (or one, if $r_{\rho}=0$) slabs $\underline{\rho}, \underline{\rho}' \subseteq \rho \in \bar{\mathscr{P}}^{[1]}$ and a maximal cell $\sigma \in \bar{\mathscr{P}}^{\max}$ adjacent to $\rho$, choose homomorphisms of abelian groups $$\begin{aligned}
&s_{\sigma \underline{\rho}}: \Lambda_{\sigma} \to A^{\times}  \\ &s_{\sigma \underline{\rho}'}: \Lambda_{\sigma} \to A^{\times}  \end{aligned}$$ subject to the constraint 
\begin{equation} \label{opengd}
s_{\sigma \underline{\rho}}|_{\Lambda_{\rho}} \cdot (s_{\sigma \underline{\rho}'}|_{\Lambda_{\rho}})^{-1} = s_{\sigma' \underline{\rho}}|_{\Lambda_{\rho}} \cdot (s_{\sigma' \underline{\rho}'}|_{\Lambda_{\rho}})^{-1}    
\end{equation}
if $\rho = \sigma \cap \sigma'$ for $\sigma, \sigma' \in \bar{\mathscr{P}}^{\max}$ (if there is just one slab $\rho$, the condition simplifies to $s_{\sigma \underline{\rho}}|_{\Lambda_{\rho}} = s_{\sigma' \underline{\rho}}|_{\Lambda_{\rho}}$). We call this collection $s$ of homomorphisms \emph{(open) gluing data}. Setting all homomorphisms to be trivial defines \emph{trivial gluing data}.
\end{defn}

For a slab $\mathfrak{b} \subseteq \underline{\rho} \in \tilde{\bar{\mathscr{P}}}^{[1]}$ and a chamber $\mathfrak{u} \subseteq \sigma \in \bar{\mathscr{P}}^{\max}$, we define the \emph{localization homomorphism $\chi^s_{\mathfrak{b}, \mathfrak{u}}$ twisted by $s$} as the composition of the canonical $\chi_{\mathfrak{b}, \mathfrak{u}}:R_{\mathfrak{b}} \to R_{\mathfrak{u}}$ of \eqref{lochom} with the map
\begin{equation} \label{gdhomomorphism}
s_{\sigma \underline{\rho}}: R_{\mathfrak{u}} \to R_{\mathfrak{u}}, \quad z^m \mapsto s_{\sigma \underline{\rho}} (m) z^m.
\end{equation}

We also need to modify the equation \eqref{comp} relating the initial slab functions $f_{\underline{\rho}}, f_{\underline{\rho}'}  \in (A[Q]/I_0)[\Lambda_{\rho}]$ for slabs $\underline{\rho}, \underline{\rho}' \subseteq \rho \in \bar{\mathscr{P}}^{[n-1]}$ to ensure that $R_{\underline{\rho}}$ and $R_{\underline{\rho}'}$ are compatible with the localization homomorphisms twisted by $s$. We replace \eqref{comp} with  
\begin{equation} \label{compgd}
f_{\underline{\rho}'} = \frac{s_{\sigma \underline{\rho}'}(\xi) s_{\sigma' \underline{\rho}'}(\xi')}{s_{\sigma \underline{\rho}}(\xi) s_{\sigma' \underline{\rho}}(\xi')}s^{-1}_{\sigma \underline{\rho}'}(s_{\sigma \underline{\rho}}(f_{\underline{\rho}})) z^{m_{\underline{\rho}' \underline{\rho}}}
\end{equation}
where $\xi$ is the normal generator pointing into $\sigma$ as in \eqref{posstalk1} and $\xi'$ is the parallel transport of $\xi$ through $\underline{\rho}$ into $\sigma'$. As in Section \ref{tdmpafunc}, we require that the initial slab functions have no poles and are normalized.

We have an analogue of Proposition \ref{struct}. 

\begin{fact}
Let $\underline{\rho}, \underline{\rho}' \in \tilde{\bar{\mathscr{P}}}^{[1]}$ be two slabs with $ \underline{\rho}, \underline{\rho}' \subseteq \rho \in \bar{\mathscr{P}}^{[1]}$ and let $f_{\underline{\rho}}, f_{\underline{\rho}'}$ be the corresponding slab functions. Let $w_{\rho} := z^{m_{\rho}}$ where $m_{\rho}$ is the integral generator of $\Lambda_{\rho}$ that points towards the vertex endpoint of $\underline{\rho}'$. We have
\begin{equation} \label{slabfuncgd}
\begin{aligned}f_{\underline{\rho}} = &1 + a_{\rho, 1}w_{\rho} + \dots + a_{\rho, r_{\rho}-1}w_{\rho}^{r_{\rho}-1} + c_{s, \underline{\rho}} w_{\rho}^{r_{\rho}}& \\ f_{\underline{\rho}'} = &1 + a_{\rho, 1}'w_{\rho} + \dots + a_{\rho, r_{\rho}-1}'w_{\rho}^{-r_{\rho}+1} + c_{s, \underline{\rho}'} w_{\rho}^{-r_{\rho}}&
\end{aligned}
\end{equation}
where the $a'_{\rho, i} \in \kk[a_{\rho, i}]$ are polynomials in $a_{\rho, i}, \ 1 \leq i \leq r_{\rho} -1$ such that $f_{\underline{\rho}}, f_{\underline{\rho}'}$ satisfy \eqref{compgd} and $c_{s, \underline{\rho}}, c_{s, \underline{\rho}'} \in \kk^{\times}$ are certain fixed constants determined by the gluing data $s$. A similar statement holds over $\kk$ in which case we view $a_{\rho, i}, a_{\rho, i}' \in \kk, \ \ 1 \leq i \leq r_{\rho} -1$ as constants satisfying \eqref{compgd}.
\end{fact}

\begin{proof}
We refer to the proof of \cite[Proposition A.4.1]{GHS}.
\end{proof}

The notion of consistency of a scattering diagram $\bar{\mathfrak{D}}_I$ on $\left(\bar B, \bar{\mathscr{P}}\right)$ in codimension $0$ is as in Definition \ref{consistency}. Consistency in codimension $1$ is similar, but we replace $\chi_{\mathfrak{b}, \mathfrak{u}}$ with  $\chi^s_{\mathfrak{b}, \mathfrak{u}}$. 

We also replace $\chi_{\mathfrak{b}, \mathfrak{u}}$ with $\chi^s_{\mathfrak{b}, \mathfrak{u}}$ in the construction of $\check{\bar{\mathfrak{X}}}^{o}_{\mathfrak{D}_J}$.  In general, there is an obstruction to doing that, but it is empty in our case of interest.  

\begin{fact} \label{obsglue}
The obstruction to gluing $\check{\bar{\mathfrak{X}}}^{o}_{\mathfrak{D}_I}$ using the localization homomorphisms $\chi^s_{\mathfrak{b}, \mathfrak{u}}$ vanishes if $\bar{B}$ is a manifold and $\dim \bar{B} = 2$.
\end{fact}

\begin{proof}
We need to show that $s$ is consistent for $\left(\bar B, \bar{\mathscr{P}}\right)$ in the sense of \cite[Definition 5.2.10]{GHS}. Since $\bar{B}$ is a topological manifold, by \cite[Proposition 5.2.9]{GHS} it is enough to show that every vertex $v \in \bar{\mathscr{P}}^{[0]}$ satisfies the equation \cite[(5.4)]{GHS}. But that equation is trivially satisfied since $\dim \bar{B} = 2$. 
\end{proof}
Given a scattering diagram $\bar{\mathfrak{D}}_I$ on $(\bar B, \bar{\mathscr{P}})$ consistent in codimensions $0$ and $1$ for gluing data $s$, we have a well-defined $\check{\bar{\mathfrak{X}}}^{o,s}_{\bar{\mathfrak{D}}_I}$. Here we introduce $s$ in the notation to indicate the use of gluing data in the construction. We will extend this notation to other objects obtained using $s$.\footnote{Our notation here differs from \cite{GHS} where $s$ is implicit in the notation.} 

To define consistency in codimension $2$, one needs to modify the requirement $a_1=1$ on the starting coefficient $a_1$ in the Definition \ref{brokenlines} of a broken line to ensure that theta functions are compatible with slab-crossings. We require that $a_1 = 1$ in a chosen cell $\sigma \in \bar{\mathscr{P}}_{\mathfrak{j}}^{\max}$ such that the asymptotic monomial $m$ of the broken line satisfies $m \in \Lambda_{\sigma}$. Then for a cell $\sigma' \in \bar{\mathscr{P}}_{\mathfrak{j}}^{\max}$ with $\sigma \cap \sigma' = \rho \in \bar{\mathscr{P}}_{\mathfrak{j}}^{[1]}$ and $m \in \Lambda_{\sigma'}$, we need to modify the initial coefficient of a broken line with asymptotic monomial $m$ by $s^{-1}_{\sigma \underline{\rho}}(\bar m) s_{\sigma \underline{\rho}}(\bar m)$. We can define the initial coefficient in all cells $\sigma'' \in \bar{\mathscr{P}}_{\mathfrak{j}}^{\max}$ with $m \in \Lambda_{\sigma''}$ in this fashion. The vanishing of the obstruction of Proposition \ref{obsglue} and contractability of the subcomplex of cells $\sigma'' \in \bar{\mathscr{P}}_{\mathfrak{j}}^{\max}$ with $m \in \Lambda_{\sigma''}$ along with their faces imply that this is well-defined. With this modification, consistency in codimension $2$ is defined as in Definition \ref{consistency} using $\chi^s_{\mathfrak{b}, \mathfrak{u}}$ instead of $\chi_{\mathfrak{b}, \mathfrak{u}}$ to define the slab-crossing homomorphisms of \eqref{slabcrossing}.

Suppose that $\bar{\mathfrak{D}}_I$ is consistent for gluing data $s$. There is an additional obstruction to extending the mirror family over the codimension $2$ strata. Due to Proposition \ref{obsglue}, every open gluing data $s$ induces a cohomology class $\bar{s} \in H^1\left(\bar{B}, \mathcal{Q} \otimes \underline{A}^{\times}\right)$\footnote{Here $\mathcal{Q}$ is a certain constructible sheaf with stalks $\Hom(\Lambda_{\tau}, \mathbb{Z})$ along $\Int \check{\tau}$ for every cell $\check{\tau} \in \check{\bar{\mathscr{P}}}$ of the Legendre dual decomposition $\check{\bar{\mathscr{P}}}$ of $\bar{B}$. See \cite[Section 5.2]{GHS} for details. } that we call \emph{closed gluing data}.

\begin{defn}
We say that $s$ is \emph{projective} if the induced closed gluing data $\bar{s}$ lifts to a \emph{lifted closed gluing data} $\tilde{\bar{s}}$ for $\mathbf{C}B$ in the sense of \cite[Definition 5.2.12]{GHS}.
\end{defn}

\begin{sloppypar}
In general, there might be different lifts (the set of lifts is a torsor for $H^1\left(\bar{B}, \underline{A}^{\times}\right)$) producing non-isomorphic families $$\check{\bar{\mathfrak{X}}}^s_{\bar{\mathfrak{D}}_I} := \Proj \Gamma(\check{\bar{\mathfrak{N}}}_{\bar{\mathbf{C}\mathfrak{D}}_I}^{o,\tilde{s}}, \mathcal{O}_{\check{\bar{\mathfrak{N}}}_{\bar{\mathbf{C}\mathfrak{D}}_I}^{o,\tilde{s}}})$$ (here $\tilde{s}$ is an open gluing data inducing $\tilde{\bar{s}}$, see the proof of \cite[Proposition 5.2.13]{GHS}). However, as explained in the proof of \cite[Proposition A.4.1]{GHS}, reexamining the proof of \cite[Theorem 5.4]{GS2} shows that there is always a unique lift $\tilde{\bar{s}}$ of $\bar{s}$ in the case that $\left(\bar B, \bar{\mathscr{P}}\right)$ is the dual intersection complex of a toric degeneration. So $\check{\bar{\mathfrak{X}}}^s_{\bar{\mathfrak{D}}_I}$ is well-defined.
\end{sloppypar}

Theorem \ref{reconstruction} admits a generalization to non-trivial gluing data $s$ (see \cite[Proposition 3.9]{GS1}). Using this generalization, one can obtain a collection of scattering diagrams $\bar{\mathfrak{D}} := \left\{ \bar{\mathfrak{D}}_{k}, \ k \geq 0 \right\}$ consistent for gluing data $s$. 
Fixing the initial slab functions $$\left\{f_{\underline{\rho}} \in \kk[a_{\rho, i}][\Lambda_{\rho}] \ | \ \underline{\rho} \in \tilde{\mathscr{P}}^{[n-1]} \right\}$$ satisfying \eqref{slabfuncgd}, we define the \emph{initial scattering diagram} $\bar{\mathfrak{D}}_0:=\bar{\mathfrak{D}}_{\mathfrak{m}}$ to have slabs with support $\underline{\rho} \in \tilde{\bar{\mathscr{P}}}^{[1]}$ and the attached slab functions $f_{\underline{\rho}}$ as before. Consistency of $\bar{\mathfrak{D}}_0$ follows from equation \eqref{compgd} relating $f_{\underline{\rho}}, f_{\underline{\rho}'}$ for two slabs $\underline{\rho}, \underline{\rho}' \subseteq \rho \in \bar{\mathscr{P}}^{[1]}$. Then one can produce a collection $\bar{\mathfrak{D}} := \left\{ \bar{\mathfrak{D}}_{k}, \ k \geq 0 \right\}$ of scattering diagrams with $\bar{\mathfrak{D}}_{k}:= \bar{\mathfrak{D}}_{\mathfrak{m}^{k+1}}$ such that $\bar{\mathfrak{D}}_{k}, \ k \geq 0$ are consistent for $s$ in the sense of \cite[Definition 2.28]{GS1} and satisfy the compatibility and uniqueness properties of Theorem \ref{reconstruction}.

\subsubsection{Non-normalized fibres via gluing data} \label{nonnormalized}

We want to generalize the result of Section \ref{freeparam} to the restriction of the extended intrinsic mirror to the family over the subvariety $$\left\{\prod_{k=1}^{r_{\rho}} t^{E_{\rho, k}} = c_{\rho}, \ \rho \in \bar{\mathscr{P}}^{[1]} \ | \ c_{\rho} \in \kk^{\times}\right\}.$$ The crucial extra step is finding the gluing data $s$ normalizing the slab functions of all the slabs in $\Phi(\mathfrak{D}_0)$. Using the ring $\kk[a_{\rho, i}]$, for a maximal cell $\sigma \in \bar{\mathscr{P}}^{\max}$ adjacent to $\rho \in \bar{\mathscr{P}}^{[1]}$ we define the gluing data $s$ via:
\begin{equation} \label{gluingdatak}
\begin{aligned}
&s_{\sigma \underline{\rho}}: \Lambda_{\sigma} \to (\kk[a_{\rho, i}])^{\times} = \kk^{\times}, \quad \begin{cases} \begin{aligned} m_{\rho} &\mapsto c_{\rho}^{-\frac{1}{r_{\rho}}} \\ \xi &\mapsto 1 \end{aligned} \end{cases} \\ &s_{\sigma \underline{\rho}'}: \Lambda_{\sigma} \to (\kk[a_{\rho, i}])^{\times} = \kk^{\times}, \quad \ \ \ m \mapsto 1 \end{aligned}
\end{equation}
where the slabs $\underline{\rho}, \underline{\rho}' \in \tilde{\bar{\mathscr{P}}}^{[1]}$ are chosen as in Sections \ref{admissibleandstronglyso} and \ref{global42}, $\xi$ is the normal generator pointing into $\sigma$ as in \eqref{posstalk1}, and we fix some choice of primitive $r_{\rho}$-th roots for $r_{\rho} \geq 1$. If $r_{\rho} = 0$, we interpret $s_{\sigma \underline{\rho}}$ as the trivial homomorphism. It is easy to check using the long exact cohomology sequence for \cite[(5.9)]{GHS} (as in the proof of \cite[Proposition 5.2.13]{GHS}) that this gluing data lifts. So $s$ is projective, and there is a unique lift of $\bar{s}$ to a lifted closed gluing data $\tilde{\bar{s}}$ on $\mathbf{C}B$.

The corresponding initial slab function $f_{\underline{\rho}}$ of \eqref{slabfuncgd} is
\begin{equation} \label{glueslabfunc1}
f_{\underline{\rho}} = 1 + a_{\rho, 1}w_{\rho} + \dots + a_{\rho, r_{\rho}-1}w_{\rho}^{r_{\rho}-1} + c_{\rho} w_{\rho}^{r_{\rho}}  
\end{equation}
and from equation \eqref{compgd} relating $f_{\underline{\rho}}$ and $f_{\underline{\rho}'}$ we have
\begin{equation} \label{glueslabfunc2}
f_{\underline{\rho}'} = \frac{s_{\sigma \underline{\rho}'}(\xi) s_{\sigma' \underline{\rho}'}(\xi')}{s_{\sigma \underline{\rho}}(\xi) s_{\sigma' \underline{\rho}}(\xi')}s^{-1}_{\sigma \underline{\rho}'}(s_{\sigma \underline{\rho}}(f_{\underline{\rho}})) z^{m_{\underline{\rho}' \underline{\rho}}} = s_{\sigma \underline{\rho}}(f_{\underline{\rho}}) z^{m_{\underline{\rho}' \underline{\rho}}}
\end{equation}
which has constant coefficient $1$, so it is normalized. The initial slab functions define $\bar{\mathfrak{D}}_0$ that is consistent for $s$ and by the discussion at the end of Section \ref{gluingdataintro} we get a collection of scattering diagrams $\bar{\mathfrak{D}} := \left\{ \bar{\mathfrak{D}}_{k}, \ k \geq 0 \right\}$ and a family $\check{\bar{\mathfrak{X}}}^s_{\bar{\mathfrak{D}}} \to \Spec \kk [a_{\rho, i}]\llbracket Q \rrbracket$ defined as the inverse limit over $\check{\bar{\mathfrak{X}}}^s_{\bar{\mathfrak{D}}} \to \Spec \kk [a_{\rho, i}][Q]/ \mathfrak{m}^{k+1}$ as usual. This family is isomorphic to the restriction of the universal toric degeneration mirror of
\cite[Theorem A.4.2]{GHS} to closed gluing data $\bar{s}$.

\begin{remark} \label{glueremarkcorr}
Similarly to Remark \ref{freeremark1}, the family $\check{\bar{\mathfrak{X}}}_{\bar{\mathfrak{D}}}^s \to \Spec \kk [a_{\rho, i}]\llbracket Q \rrbracket$ is universal in the sense that every one-parameter toric degeneration mirror family $\check{\bar{\mathfrak{X}}}^s \to \Spec \kk \llbracket t \rrbracket$ constructed using gluing data $s$ and polarization $A$ can be obtained from $\check{\bar{\mathfrak{X}}}_{\bar{\mathfrak{D}}}^s \to \Spec \kk [a_{\rho, i}]\llbracket Q \rrbracket$ via a basechange that is induced by $Q \to \mathbb{N}, \ \beta \mapsto A \cdot \beta$ and by sending the variables $a_{\rho, i}$ to the values chosen for the initial slab function $f_{\underline\rho}$ of \eqref{slabfuncgd}.
\end{remark}

We proceed as in Section \ref{freeparam}. First, we restrict the minimal relative Gross-Siebert locus to a family over $$\left\{\sigma_{r_{\rho}} (t^{E_{\rho}}) = \prod_{k=1}^{r_{\rho}} t^{E_{\rho, k}} = c_{\rho}, \ \rho \in \bar{\mathscr{P}}^{[1]}\right\}.$$ We generalize \eqref{projmap1} by setting $$A^{GS} := \kk [t^{\pm E_{\rho, k} }]\llbracket Q \rrbracket, \quad A^{\GS}_{c_{\rho}}:= \kk [t^{\pm E_{\rho, k} }]\llbracket Q \rrbracket \Bigg/ \left \langle\prod_{k=1}^{r_{\rho}} t^{E_{\rho, k}} = c_{\rho}, \ \rho \in \bar{\mathscr{P}}^{[1]} \right \rangle$$ and considering the natural projection
$h_{c_{\rho}}: A^{\GS} \to A^{\GS}_{c_{\rho}}$. We still write $t^{E_{\rho, k}}$ for the image of $t^{E_{\rho, k}}$ under $h_{c_{\rho}}$.

\begin{cons}
 We define a collection of scattering diagrams $\mathfrak{D} = \left\{ \mathfrak{D}_k, \ k \geq 0 \right\}$ in the same way as in Construction \ref{freeinitdiag} (via Remark \ref{remring2} and using $h_{c_{\rho}}$ instead of $h$). The resulting family $\check{\mathfrak{X}}_{\mathfrak{D}} \to \Spec A^{\GS}_{c_{\rho}}\llbracket Q \rrbracket$ is the same family as the basechange of the minimal relative Gross-Siebert locus $\check{\mathfrak{X}} \to \Spec A^{\GS}\llbracket Q \rrbracket$ by $h_{c_{\rho}}: A^{\GS} \to A^{\GS}_{c_{\rho}}$.
\end{cons}

Now we need to interpret the variables $a_{\rho, i}$ as elements of $A^{\GS}_{c_{\rho}}$ to set up the correspondence. Similarly to Construction \ref{freecons2}, this is achieved by a basechange. 

\begin{sloppypar}
\begin{cons} \label{basechangeglue}
We consider the map $$\sigma: \kk [a_{\rho, i}] \to  A^{\GS}_{c_{\rho}}, \quad  a_{\rho, i} \mapsto \sigma_i (t^{E_{\rho}}), \ 1 \leq i \leq r_{\rho} -1, \ \rho \in \bar{\mathscr{P}}^{[1]}$$ as in \eqref{sigmabasechange} and define a new collection of scattering diagrams $\bar{\mathfrak{D}}^{\sigma} := \left\{ \bar{\mathfrak{D}}_{k}^{\sigma}, \ k \geq 0 \right\}$ with ring $A^{\GS}_{c_{\rho}}$, monoid $Q$, and ideal $I_0 = \mathfrak{m}$ in the same way as in Construction \ref{freecons2}. The family $\check{\bar{\mathfrak{X}}}^s_{\bar{\mathfrak{D}}^{\sigma}} \to \Spec A^{\GS}_{c_{\rho}}\llbracket Q \rrbracket$ constructed using the gluing data (that we still call $s$) defined via the homomorphisms $$s_{\sigma \underline{\rho}}: \Lambda_{\sigma} \to \kk^{\times} \subseteq \left(A_{c_{\rho}}^{\GS}\right)^{\times}$$ induced by the $s$ of \eqref{gluingdatak} is the same family as the basechange of $\check{\bar{\mathfrak{X}}}^s_{\bar{\mathfrak{D}}} \to \Spec \kk [a_{\rho, i}]\llbracket Q \rrbracket$ by $\sigma$. If $\dim E^{\gp} = \sum_{\rho \in \bar{\mathscr{P}}^{[1]}} r_{\rho}$, then the induced map $\check{\bar{\mathfrak{X}}}^s_{\bar{\mathfrak{D}}^{\sigma}} \to \check{\bar{\mathfrak{X}}}^s_{\bar{\mathfrak{D}}}$ is generically a covering of index $\prod_{\rho \in \bar{\mathscr{P}}^{[1]}} r_{\rho}!$ ramified at the points where $t^{E_{\rho, i}} = t^{E_{\rho, j}}$ for some $\rho \in \bar{\mathscr{P}}^{[1]}$ and $1 \leq i < j \leq r_{\rho}$. If $\dim E^{\gp} < \sum_{\rho \in \bar{\mathscr{P}}^{[1]}} r_{\rho}$, then the morphism $\check{\bar{\mathfrak{X}}}^s_{\bar{\mathfrak{D}}^{\sigma}} \to \check{\bar{\mathfrak{X}}}^s_{\bar{\mathfrak{D}}}$ is not surjective.
\end{cons}
\end{sloppypar}

Now that the scattering diagrams $\mathfrak{D}$ and $\bar{\mathfrak{D}}^{\sigma}$ are defined over the same ring, we define $\Phi(\mathfrak{D}_k)'$ and  $\Phi(\mathfrak{D}_k)$ in a way that ensures that $\Phi(\mathfrak{D}_0)$ is combinatorially equivalent to $\bar{\mathfrak{D}}^{\sigma}_0$.

\begin{cons} \label{chsgd}
As in Construction \ref{univ2}, we use the analogue of Construction \ref{refinement} to get a scattering diagram $\bar{\mathfrak{D}}^{\sigma'}$ on $\left(\bar{B}, \bar{\mathscr{P}}'\right)$. Note that the gluing data $s$ of \eqref{gluingdatak} defines a gluing data $s'$ on $\left(\bar{B}, \bar{\mathscr{P}}'\right)$ by defining the gluing data for every slab $\underline{\rho_p} \subseteq \rho_p \in \bar{\mathscr{P}}'^{[1]}$ to be the same as the gluing data on the unique slab $\underline{\rho} \subseteq \rho \in \bar{\mathscr{P}}^{[1]}$ containing $\rho_p$ and defining the gluing data for $\rho' \subseteq \bar{\mathscr{P}}'^{[1]} \setminus \bar{\mathscr{P}}^{[1]}$ (using Notation \ref{nonsubdividing}) to be trivial.

We modify Construction \ref{chs} of $\Phi(\mathfrak{D}_k)'$ and $\Phi(\mathfrak{D}_k)$ by replacing  $z^{m_{\underline{\Phi(\rho)}'\underline{\Phi(\rho)} }}f_{\Phi(\mathfrak{b})}$ in (2) and (3) with $z^{m_{\underline{\Phi(\rho)}'\underline{\Phi(\rho)}}} s'_{\Phi(\sigma) \underline{\Phi(\rho_p)}}(f_{\Phi(\mathfrak{b})})$.
Thus defined $\Phi(\mathfrak{D}_k)'$ and $\Phi(\mathfrak{D}_k)$ are compatible with the MPA functions and it follows immediately from the construction of $s$ in \eqref{gluingdatak} and $f_{\underline{\rho}}, f_{\underline{\rho}'}$ in \eqref{glueslabfunc1}, \eqref{glueslabfunc2} that $\Phi(\mathfrak{D}_0)$ is combinatorially equivalent to $\bar{\mathfrak{D}}^{\sigma}_0$ (one needs to verify that $f_{\Phi(\mathfrak{b}_1)} = z^{m_{\underline{\Phi(\rho)}'\underline{\Phi(\rho)}}} s'_{\Phi(\sigma) \underline{\Phi(\rho_p)}}(f_{\Phi(\mathfrak{b}_2)})$ for $\mathfrak{b}_1$ as in Construction \ref{chs}(1) with $\mathfrak{b}_1 \subseteq \rho_{p_1} \in \mathscr{P}^{[1]}$ for $1 \leq p_1 \leq p_0$ and $\mathfrak{b}_2$ as in Construction \ref{chs}(2,3) with $\mathfrak{b}_2 \subseteq \rho_{p_2} \in \mathscr{P}^{[1]}$ for $p_0+1 \leq p_2 \leq l_{\rho}$). Here we extend Definition \ref{combequiv} of combinatorial equivalence to non-trivial gluing data by requiring that the gluing data used to construct the family is the same for the two scattering diagrams.
\end{cons}

We can now prove an analogue of Proposition \ref{prooffree}.

\begin{fact} \label{proofglue}
For any gluing data $s$ as in \eqref{gluingdatak}, the basechange of the minimal relative Gross-Siebert locus $\check{\mathfrak{X}} \to \Spec  A\llbracket Q \rrbracket$ by $h_{c_{\rho}}: A^{\GS} \to A^{\GS}_{c_{\rho}}$ is isomorphic to the basechange by $\sigma: \kk [a_{\rho, i}] \to A^{\GS}_{c_{\rho}}, \  a_{\rho, i} \mapsto \sigma_i (t^{E_{\rho}})$ of the restriction $\check{\bar{\mathfrak{X}}}^s_{\bar{\mathfrak{D}}} \to \Spec \kk [a_{\rho, i}]\llbracket Q \rrbracket$ of the universal toric degeneration mirror of
\cite[Theorem A.4.2]{GHS} to closed gluing data $\bar{s}$ (induced from $s$).
\end{fact}

\begin{proof}
We need to show that the families $\check{\mathfrak{X}}_{\mathfrak{D}} \to \Spec A^{\GS}_{c_{\rho}}\llbracket Q \rrbracket$ and $\check{\bar{\mathfrak{X}}}^s_{\bar{\mathfrak{D}}^{\sigma}} \to \Spec A^{\GS}_{c_{\rho}}\llbracket Q \rrbracket$ are the same. First, we need to check consistency of $\Phi(\mathfrak{D}_k)$. Remark \ref{obsplimage}(1) still applies in this setting, so it is enough to show consistency of $\Phi(\mathfrak{D}_k)'$. This follows as in the proof of Proposition \ref{consis} with the following modifications:
\begin{enumerate}
    \item Work over the ring $A^{\GS}_{c_{\rho}}$, use monoid $Q$, and replace Constructions \ref{refinement} and \ref{chs} with Construction \ref{chsgd}. The proof of consistency in codimension $0$ is the same.
    \item Replace the canonical localization homomorphisms $\chi_{\Phi(\mathfrak{b}), \Phi(\mathfrak{u})}$ with 
    $\chi^{s'}_{\Phi(\mathfrak{b}), \Phi(\mathfrak{u})}$ and the canonical isomorphisms $R_{\mathfrak{b}} \cong R_{\Phi(\mathfrak{b})}$ with 
    \begin{equation} \label{Phirho}
    \Phi_{\underline{\rho}'}:  R_{\mathfrak{b}} \to  R_{\Phi(\mathfrak{b})}^{s'}, \quad Z_+ \mapsto Z_+, \ Z_- \mapsto Z_-, \ w_{\rho} \mapsto  {s'}^{-1}_{\Phi(\sigma) \underline{\Phi(\rho')}}(w_{\rho}).
    \end{equation}
    for $\mathfrak{b} \subseteq \rho' \in \mathscr{P}^{[1]}$ a slab of the first type. 
    Modify the argument in Step 1 of the proof of consistency in codimension $1$ to an argument similar to the one in Step 2.
    \item Replace the isomorphism $\Phi': R_{\mathfrak{b}_i} \to R_{\Phi(\mathfrak{b}_i)'}$ of \eqref{proof2} with, for 
    $$
    \begin{aligned}
    R_{\mathfrak{b}_i} &= \left( A^{\GS}_{c_{\rho}}[Q]/\mathfrak{m}^{k+1}\right)[\Lambda_{\rho}][Z_+, Z_-]/\left(Z_+Z_- - f_{\mathfrak{b}_i} \cdot z^{\kappa_{\rho_{p_0+1}}}\right)  \\
    R^{s'}_{\Phi(\mathfrak{b}_i)'} &= \left(A^{\GS}_{c_{\rho}}[Q]/\mathfrak{m}^{k+1}\right)[\Lambda_{\Phi(\rho)}][Z_+, Z_-] \Big/ \left(Z_+Z_- - z^{m_{\underline{\Phi(\rho)}'\underline{\Phi(\rho)}}} s'_{\Phi(\sigma) \underline{\Phi(\rho)}}(f_{\Phi(\mathfrak{b}_i)}) \cdot z^{\kappa_{\Phi(\rho_{p_0+1})}}\right)
    \end{aligned}
    $$ (for $i=1,2$ and with $\kappa_{\rho_{p_0+1}} = \kappa_{\Phi(\rho_{p_0+1})}$ by construction), the map
    \begin{equation} 
    \Phi': R_{\mathfrak{b}_i} \to R^{s'}_{\Phi(\mathfrak{b}_i)'}, \quad Z_{+} \mapsto Z_{+}, \ Z_{-} \mapsto  z^{m_{\underline{\Phi(\rho)} \, \underline{\Phi(\rho)}'}} Z_{-}, \ w_{\rho} \mapsto {s'}^{-1}_{\Phi(\sigma) \underline{\Phi(\rho_p)}}(w_{\rho}).
    \end{equation} The arguments of Steps 2 and 3 in the proof of consistency in codimension $1$ go through with this modification. 
    \item In the proof of consistency in codimension $2$, use broken lines in the modified sense of Section \ref{gluingdataintro}. The proof that $\Phi(\beta)$ is a well-defined broken line goes through (using the modified slab-crossing homomorphisms). Compatibility with the change of chamber homomorphism $\theta_{\Phi(\mathfrak{u}') \Phi(\mathfrak{u})}$ for $\Phi(\mathfrak{u}') \cap \Phi(\mathfrak{u}) = \Phi(\mathfrak{b})'$ a slab of the second type follows in the same way (using the modified slab-crossing homomorphisms and $\Phi'$). 
\end{enumerate}

Equivalence of $\left(\bar B, \Phi(\mathfrak{D}_k)\right)$ and $(B, \mathfrak{D}_k)$ should be understood in the sense of using trivial gluing data on $(B, \mathscr{P})$ and gluing data $s$ on $(\bar{B}, \bar{\mathscr{P}})$. Remark \ref{obsplimage}(1) still applies in this setting, so it is enough to show the equivalence of $\left(\bar B, \Phi(\mathfrak{D}_k)'\right)$ and $(B, \mathfrak{D}_k)$ (using gluing data $s'$ on $\left(\bar{B}, \bar{\mathscr{P}}'\right)$).
This follows as in the proof of Proposition \ref{equivalence} with the same modifications as in (1) and (2) above. We need to slightly modify the last paragraph of the proof since the isomorphism $\Phi_{\underline{\rho}'}$ of \eqref{Phirho} is not canonical. However, it follows from the constructions of $\check{\mathfrak{X}}^{o, s'}_{(\bar{B}, \Phi(\mathfrak{D}_k)')}$ and $\check{\mathfrak{X}}^o_{(B, \mathfrak{D}_k)}$ that it is enough to check that the diagrams 
$$\begin{tikzcd}
R_{\rho'} \arrow{r}{\chi_{\rho', \sigma}} \arrow{d}{\Phi_{\underline{\rho}}} & {R_{\sigma}} \arrow{d}{\Id} \\
R^{s'}_{\underline{\Phi(\rho')}} \arrow{r}{\chi^{s'}_{\underline{\Phi(\rho')}, \Phi(\sigma)}} & R^{s'}_{\Phi(\sigma)}
\end{tikzcd}
$$
commute for all $\rho' \in \mathscr{P}^{[1]}$ and $\sigma \in \mathscr{P}^{\max}$ with $\rho' \subseteq \sigma$. This is immediate from the definitions.

Now, changing $\bar{\mathfrak{D}}$ to $\bar{\mathfrak{D}}^{\sigma}$, $\kk$ to $A^{\GS}_{c_{\rho}}$, $\mathbb{N}$ to $Q$, $(t^{k+1})$ to $\mathfrak{m}^{k+1}$, the trivial gluing data for $\left(\bar B, \bar{\mathscr{P}}\right)$ to the gluing data $s$ of \eqref{gluingdatak}, and using $h_{c_{\rho}}: A^{\GS} \to A^{\GS}_{c_{\rho}}$, the same argument as in the proof of Theorem \ref{proofconj} implies the result. 

There is one place where an additional argument is required. We can't directly apply Proposition \ref{help} that says that equivalent scattering diagrams induce isomorphic families to claim that it is enough to prove equivalence of $(B, \mathfrak{D}_k)$ and $\left(\bar B, \Phi(\mathfrak{D}_k)'\right)$ in the sense discussed above to ensure that $\check{\mathfrak{X}}_{(B, \mathfrak{D}_k)} \cong \check{\mathfrak{X}}^{s'}_{(\bar{B}, \Phi(\mathfrak{D}_k)')}$. Indeed, the definition of $\check{\bar{X}}^{s'}_0$ is different from the definition of $\check{X}_0$. Instead of the ring $\left( A_{c_{\rho}}^{\GS} [Q]/\mathfrak{m} \right)[B]$ of \eqref{ringcentralfibre}, one uses the ring $\left( A_{c_{\rho}}^{\GS} [Q]/\mathfrak{m} \right)[\bar{B}]\left(\tilde{\bar{s}}'\right)$ where the multiplication rule depends on the lifted closed gluing data $\tilde{\bar{s}}'$ (see \cite[Section 5.2]{GHS}). However, $\check{\bar{X}}^{s'}_0:= \Proj \left( A_{c_{\rho}}^{\GS} [Q]/\mathfrak{m} \right)[\bar{B}]\left(\tilde{\bar{s}}'\right)$ is isomorphic to $\check{X}_0:= \Proj \left( A_{c_{\rho}}^{\GS} [Q]/\mathfrak{m} \right)[B]$ via
$$
\begin{aligned}
\alpha: \left( A_{c_{\rho}}^{\GS} [Q]/\mathfrak{m} \right)[B] &\to \left( A_{c_{\rho}}^{\GS} [Q]/\mathfrak{m} \right)[\bar{B}]\left(\tilde{\bar{s}}'\right) \\
z^{m} &\mapsto {\tilde{\bar{s}}'}^{-1}_{\Phi(\sigma) \underline{\Phi(\rho')}}(z^{m}), \quad   m \in \mathbf{C}\rho', \ \rho' \in \mathscr{P}^{[1]} \\
z^m &\mapsto z^m, \quad \quad \quad \quad \quad \  m \notin \mathbf{C}\rho', \ \rho' \in \mathscr{P}^{[1]}
\end{aligned}
$$
which is clearly compatible with the isomorphisms $\Phi_{\underline{\rho'}}$ of \eqref{Phirho}. Using the isomorphism $\alpha$, an argument similar to the one in the proof of Proposition \ref{help} implies that it is enough to prove equivalence of $(B, \mathfrak{D}_k)$ and $\left(\bar B, \Phi(\mathfrak{D}_k)'\right)$ in the sense discussed above to ensure that $\check{\mathfrak{X}}_{(B, \mathfrak{D}_k)} \cong \check{\mathfrak{X}}^{s'}_{(\bar{B}, \Phi(\mathfrak{D}_k)')}$. 
\end{proof}

\begin{observe}
If $s$ is the \emph{trivial gluing data}, then $c_{\rho} = 1$ for all $\rho \in \bar{\mathscr{P}}^{[1]}$ and  Proposition \ref{proofglue} reduces to Proposition \ref{prooffree}.
\end{observe}

\subsubsection{The minimal relative Gross-Siebert locus via gluing data} \label{compgencons}

The subvarieties $$\left\{\prod_{k=1}^{r_{\rho}} t^{E_{\rho, k}} = c_{\rho}, \ \rho \in \bar{\mathscr{P}}^{[1]} \ | \ c_{\rho} \in \kk^{\times}\right\}$$ of $\Spec \kk [t^{\pm E_{\rho, k} }]$ cover $\Spec \kk [t^{\pm E_{\rho, k} }]$. In particular, Proposition \ref{proofglue} implies that we can realize any fibre of the minimal relative Gross-Siebert locus $\check{\mathfrak{X}} \to \Spec \kk [t^{\pm E_{\rho, k} }]\llbracket Q \rrbracket$ as a (universal in polarization) toric degeneration mirror family $\check{\bar{\mathfrak{X}}}^s \to \Spec \kk \llbracket Q \rrbracket$ constructed using gluing data $s$ of the form \eqref{gluingdatak} and a choice of slab functions $f_{\underline{\rho}}, f_{\underline{\rho}'}$ as in \eqref{glueslabfunc1}, \eqref{glueslabfunc2} (with $a_{
\rho, i}, a_{\rho, i}' \in \kk$). Note that not all (universal in polarization) toric degeneration mirror families of this form arise as fibres of $\check{\mathfrak{X}} \to \Spec \kk [t^{\pm E_{\rho, k} }]\llbracket Q \rrbracket$ unless $\dim E^{\gp} = \sum_{\rho \in \bar{\mathscr{P}}^{[1]}} r_{\rho}$.  

We would like to do the toric degeneration mirror construction universally in gluing data of the form \eqref{gluingdatak} and exhibit a connection between the minimal relative Gross-Siebert locus of the extended intrinsic mirror and a certain toric degeneration mirror family. To do that, we need to think of $c_{
\rho} \in \kk^{\times}, \ \rho \in \bar{\mathscr{P}}^{[1]}$ as variables and work over $\kk [a_{\rho, i}, c^{\pm 1}_{\rho}]$ instead of $\kk [a_{\rho, i}]$. The slab functions $f_{\underline{\rho}}, f_{\underline{\rho}'}$ of \eqref{glueslabfunc1}, \eqref{glueslabfunc2} make sense as elements of $\kk [a_{\rho, i},c^{\pm 1}_{\rho}][\Lambda_{\rho}]$. However, unless $\left(\bar B, \bar{\mathscr{P}}\right)$ is simple (that is unless $r_{\rho} \leq 1$ for all $\rho \in \bar{\mathscr{P}}^{[1]}$), the maps defined in \eqref{gluingdatak} can't be viewed as homomorphisms $\Lambda_{\sigma} \to (\kk [a_{\rho, i},c^{\pm 1}_{\rho}])^{\times} = \kk^{\times} [c^{\pm 1}_{\rho}]$ since $c_{\rho}^{-\frac{1}{r_{\rho}}}$ does not make sense as an element of $\kk [a_{\rho, i},c^{\pm 1}_{\rho}]$ for $r_{\rho} \geq 2$. To circumvent this issue, we do an additional basechange.

\begin{cons} \label{compgen1}
We work over $\kk [a_{\rho, i},c^{\pm 1}_{\rho}]$, use monoid $Q$, and set $I_0 = \mathfrak{m}$.
Consider the map 
\begin{equation} \label{mapr}
r: \kk [a_{\rho, i},c^{\pm 1}_{\rho}] \to \kk [a_{\rho, i},c^{\pm 1}_{\rho}], \quad a_{\rho, i} \mapsto a_{\rho, i}, \ c_{\rho} \mapsto c_{\rho}^{r_{\rho}}
\end{equation}
and define the initial slab functions $f^r_{\underline{\rho}}:= r(f_{
\underline{\rho}}), f^r_{\underline{\rho}'}:= r(f_{
\underline{\rho}'})$ where $f_{\underline{\rho}}, f_{\underline{\rho}'}$ are as in \eqref{glueslabfunc1}, \eqref{glueslabfunc2}. Let $s^r$ be the gluing data on $\left(\bar B, \bar{\mathscr{P}}\right)$ defined via 
\begin{equation} \label{gluingdatakuniv}
\begin{aligned}
&s^r_{\sigma \underline{\rho}}: \Lambda_{\sigma} \to \kk^{\times} [c^{\pm 1}_{\rho}], \quad \begin{cases} \begin{aligned} m_{\rho} &\mapsto c_{\rho}^{-1} \\ \xi &\mapsto 1 \end{aligned} \end{cases} \\ &s^r_{\sigma \underline{\rho}'}: \Lambda_{\sigma} \to \kk^{\times} [c^{\pm 1}_{\rho}], \quad \ \ \ m \mapsto 1 \end{aligned}
\end{equation}
for every $\rho \in \bar{\mathscr{P}}^{[1]}$ with $r_{\rho} = 1$ where the slabs $\underline{\rho}, \underline{\rho}' \in \tilde{\bar{\mathscr{P}}}^{[1]}$ are chosen as in Sections \ref{admissibleandstronglyso} and \ref{global42}, and $\xi$ is the normal generator pointing into $\sigma$ as in \eqref{posstalk1}. For $\rho \in \bar{\mathscr{P}}^{[1]}$ with $r_{\rho} = 0$, we set $s^r_{\sigma \underline{\rho}}$ to be the trivial homomorphism. Similarly to \eqref{gluingdatak}, one can check that $s^r$ is projective. Note that $f_{\underline{\rho}}^r, f_{\underline{\rho}'}^r$ satisfy \eqref{compgd} with gluing data $s^r$ since $f_{\underline{\rho}}$ and $f_{\underline{\rho}'}$ satisfy \eqref{compgd} with gluing data $s$ (treating $c_{\rho}$ as constants). 

We define $\bar{\mathfrak{D}}^r_0:=\bar{\mathfrak{D}}^r_{\mathfrak{m}}$ to have the slabs with support $\underline{\rho} \in \tilde{\bar{\mathscr{P}}}^{[1]}$ and the attached slab functions $f^r_{\underline{\rho}}, f^r_{\underline{\rho}'}$. Consistency of $\bar{\mathfrak{D}}_0^r$ for $s^r$ follows from the fact that \eqref{compgd} is satisfied. As before, \cite[Theorem 5.2]{GS1} implies that the reconstruction algorithm of Theorem \ref{reconstruction} generalizes to this setting  
and we can produce a collection $\bar{\mathfrak{D}}^r := \left\{ \bar{\mathfrak{D}}^r_{k}, \ k \geq 0 \right\}$ of scattering diagrams such that the scattering diagrams $\bar{\mathfrak{D}}^r_{k}:= \bar{\mathfrak{D}}^r_{\mathfrak{m}^{k+1}}$ are consistent for $s^r$ in the sense of \cite[Definition 2.28]{GS1} and satisfy the compatibility and uniqueness properties of Theorem \ref{reconstruction}.

We can now construct a family $\check{\bar{\mathfrak{X}}}^{s^r}_{\bar{\mathfrak{D}}^r} \to \Spec \kk [a_{\rho, i}, c_{\rho}^{\pm 1}]\llbracket Q \rrbracket$ via the recipe of Section \ref{gluingdataintro} (using the universal MPA function $\bar{\varphi}$) by taking the inverse limit as usual. We need to understand the connection with the universal toric degeneration mirror $\check{\bar{\mathfrak{X}}}_{\mathbb{P}} \to \Spec \kk_{\mathbb{P}}[a_{\rho, i}]\llbracket Q \rrbracket$ of \cite[Theorem A.4.2]{GHS}, which is defined over $\kk_{\mathbb{P}}[a_{\rho, i}]\llbracket Q \rrbracket$ for a subring $\kk_{\mathbb{P}} \subseteq \kk [H^1(\bar{B}, i_{*} \check{\Lambda})_f^{*}]$. Here $\kk [H^1(\bar{B}, i_{*} \check{\Lambda})_f^{*}]$ parameterizes closed gluing data and $\kk_{\mathbb{P}}$ is a universal choice of a subring parametrizing projective gluing data, see \cite[Appendix A.2]{GHS} for details. For every choice of $c_{\rho} \in \kk^{\times}, \ \rho \in \bar{\mathscr{P}}^{[1]}$ we obtain a gluing data (with values in $\kk^{\times}$) via \eqref{gluingdatak}. It is easy to see from the universality of $\kk_{\mathbb{P}}$ that the corresponding closed gluing data is contained in the subfamily parameterized by $\kk_{\mathbb{P}}$. Let $\kk'_{\mathbb{P}} \subseteq \kk_{\mathbb{P}}$ be the subring generated by $\bar{s}$ for $s$ as in \eqref{gluingdatak} and let $\check{\bar{\mathfrak{X}}}'_{\mathbb{P}} \to \Spec \kk'_{\mathbb{P}}[a_{\rho, i}]\llbracket Q \rrbracket$ be the restriction of $\check{\bar{\mathfrak{X}}}_{\mathbb{P}} \to \Spec \kk_{\mathbb{P}}[a_{\rho, i}]\llbracket Q \rrbracket$ to $\Spec \kk'_{\mathbb{P}}[a_{\rho, i}]\llbracket Q \rrbracket \subseteq \Spec \kk_{\mathbb{P}}[a_{\rho, i}]\llbracket Q \rrbracket$. 

We define a universal analogue of the map $r$ of \eqref{mapr}:
\begin{equation} \label{mapruniv}
r_{\univ}: \kk'_{\mathbb{P}}[a_{\rho, i}] \to \kk'_{\mathbb{P}}[a_{\rho, i}], \quad a_{\rho, i} \mapsto a_{\rho, i}, \ \bar{s} \mapsto \bar{s}^r.
\end{equation}
Here $\bar{s}^r$ is defined as follows. Let $s$ be any open gluing data giving rise to $\bar{s}$ and define an open gluing data $s^r$ by setting $s_{\sigma \underline{\rho}}^r(m) := (s_{\sigma \underline{\rho}}(m))^{r_{\rho}}$ for every $\sigma \in \bar{\mathscr{P}}^{\max}$, every slab $\underline{\rho} \subseteq \rho \in \bar{\mathscr{P}}^{[1]}$ with $\rho \subseteq \sigma$, and any $m \in \Lambda_{\sigma}$. Then $\bar{s}^r$ is the closed gluing data corresponding to $s^r$ (it is easy to check that $\bar{s}^r$ is independent of the choice of $s$ giving rise to $\bar{s}$). 
Let $\check{\bar{\mathfrak{X}}}'^{r_{\univ}}_{\mathbb{P}} \to \Spec \kk'_{\mathbb{P}}[a_{\rho, i}]\llbracket Q \rrbracket$ be the basechange of $\check{\bar{\mathfrak{X}}}'_{\mathbb{P}} \to \Spec \kk'_{\mathbb{P}}[a_{\rho, i}]\llbracket Q \rrbracket$ by $r_{\univ}$. 

\begin{sloppypar}
Now, for every $\rho \in \bar{\mathscr{P}}^{[1]}$, define a gluing data $s^{r,\rho}$ by setting $s^{r,\rho}_{\sigma, \underline{\rho}}(m) := s^{r}_{\sigma, \underline{\rho}}(m)$ (for $s^{r}$ and $\underline{\rho}$ as in \eqref{gluingdatakuniv} and any $m \in \Lambda_{\sigma}$) and setting $s^{r,\rho}_{\sigma, \underline{\rho}'}$ to be the trivial homomorphism for any other slab $\underline{\rho'} \subseteq \rho' \in \bar{\mathscr{P}}^{[1]}$ and $\sigma \in \bar{\mathscr{P}}^{\max}$ containing $\rho'$. Let $\bar{s}^{r,\rho}$ be the corresponding closed gluing data. Then the map 
$$c_{\univ}: \kk [a_{\rho, i}, c_{\rho}^{\pm 1}] \to \kk'_{\mathbb{P}}[a_{\rho, i}], \quad a_{\rho,i} \mapsto a_{\rho,i}, \ c_{\rho}  \mapsto \bar{s}^{r, \rho}$$ is well-defined and surjective. 
Our construction implies that $\check{\bar{\mathfrak{X}}}'^{r_{\univ}}_{\mathbb{P}} \to \Spec \kk'_{\mathbb{P}}[a_{\rho, i}]\llbracket Q \rrbracket$ is isomorphic to the basechange of $\check{\bar{\mathfrak{X}}}^{s^r}_{\bar{\mathfrak{D}}^r} \to \Spec \kk [a_{\rho, i}, c_{\rho}^{\pm 1}]\llbracket Q \rrbracket$ by $c_{\univ}$, i.e. we have a commutative diagram as follows:
\begin{equation} \label{mirrorfamiliescd}
\hspace{-1cm}
    \begin{tikzcd}
     \check{\bar{\mathfrak{X}}}^{s^r}_{\bar{\mathfrak{D}}^r} \arrow[r, leftarrow] \arrow{d} & \check{\bar{\mathfrak{X}}}'^{r_{\univ}}_{\mathbb{P}} \arrow{r} \arrow{d} & \check{\bar{\mathfrak{X}}}'_{\mathbb{P}} \arrow[hookrightarrow]{r} \arrow{d} & \check{\bar{\mathfrak{X}}}_{\mathbb{P}}  \arrow{d}  \\
     \Spec \kk [a_{\rho, i}, c_{\rho}^{\pm 1}]\llbracket Q \rrbracket \arrow[leftarrow]{r}{c_{\univ}} & \Spec \kk'_{\mathbb{P}}[a_{\rho, i}]\llbracket Q \rrbracket \arrow{r}{r_{\univ}} & \Spec \kk'_{\mathbb{P}}[a_{\rho, i}]\llbracket Q \rrbracket \arrow[hookrightarrow]{r} & \Spec \kk_{\mathbb{P}}[a_{\rho, i}]\llbracket Q \rrbracket   
    \end{tikzcd}
\end{equation}
\end{sloppypar}

\end{cons}

\begin{remark} \label{completeuniversality}
Similarly to Remark \ref{glueremarkcorr}, the family $\check{\bar{\mathfrak{X}}}^{s^r}_{\bar{\mathfrak{D}}^r} \to \Spec \kk [a_{\rho, i}, c_{\rho}^{\pm 1}]\llbracket Q \rrbracket$ is universal in the sense that every one-parameter toric degeneration mirror family $\check{\bar{\mathfrak{X}}}^s \to \Spec \kk \llbracket t \rrbracket$ constructed using gluing data $s$ as in \eqref{gluingdatak} and polarization $A$ can be obtained from $\check{\bar{\mathfrak{X}}}^{s^r}_{\bar{\mathfrak{D}}^r} \to \Spec \kk [a_{\rho, i}, c_{\rho}^{\pm 1}]\llbracket Q \rrbracket$ via a basechange that is induced by $Q \to \mathbb{N}, \ \beta \mapsto A \cdot \beta$ and by sending the variables $a_{\rho, i}$ to the values chosen for the initial slab function $f_{\underline\rho}$ of \eqref{slabfuncgd}, and sending the variables $c_{\rho}$ to the $\frac{1}{r_{\rho}}$-th powers of the constants used in defining $s$ (for some fixed choice of primitive $r_{\rho}$-th roots for $r_{\rho} \geq 1$ as in \eqref{gluingdatak}).
\end{remark}

\begin{sloppypar}
We will now relate the minimal relative Gross-Siebert locus $\check{\mathfrak{X}} \to \Spec \kk [t^{\pm E_{\rho, k} }]\llbracket Q \rrbracket$ to $\check{\bar{\mathfrak{X}}}^{s^r}_{\bar{\mathfrak{D}}^r} \to \Spec \kk [a_{\rho, i}, c_{\rho}^{\pm 1}]\llbracket Q \rrbracket$. This requires a basechange by a map similar to the map $r$ of \eqref{mapr}.
\end{sloppypar}

\begin{cons} \label{compgen2}
Consider the map 
\begin{equation} \label{maprintrinsic}
r: \kk [t^{\pm E_{\rho, k} }] \to \kk [t^{\pm E_{\rho, k} }], \quad t^{E_{\rho, k}} \mapsto (t^{E_{\rho, k}})^{r_{\rho}}
\end{equation}
and define a collection of scattering diagrams $\mathfrak{D}^r = \left\{ \mathfrak{D}^r_{\mathfrak{m}^{k+1}}, \ k \geq 0 \right\}$ on $(B, \mathscr{P})$ with $A=\kk [E^{\gp}]$, monoid $Q$, and $I_0=\mathfrak{m}$ by letting $\mathfrak{D}^r_{\mathfrak{m}^{k+1}}$ have walls $$(\mathfrak{p}, r(f_{\mathfrak{p}}))$$ for $\mathfrak{p}$ a wall of $\mathfrak{D}^{\GS}_{\mathfrak{m}^{k+1}}$ (using Remark \ref{remring1}). We use the MPA function $\varphi^{\GS}$. Taking the inverse limit as usual, we define a family $\check{\mathfrak{X}}_{\mathfrak{D}^{r}} \to  \Spec \kk [t^{\pm E_{\rho, k} }]\llbracket Q \rrbracket$ (using the trivial gluing data) which is just the basechange of the minimal relative Gross-Siebert locus $\check{\mathfrak{X}}_{\mathfrak{D}^{\GS}} \to  \Spec \kk [t^{\pm E_{\rho, k} }]\llbracket Q \rrbracket$ by $r$.
\end{cons}

After these basechanges, we can proceed as in Section \ref{nonnormalized}. Similarly to Construction \ref{basechangeglue}, we interpret $a_{\rho, i}$ and $c_{\rho}$ as elements of $\kk [t^{\pm E_{\rho, k} }]$ and the gluing data $s^r$ as gluing data with homomorphisms taking values in $\left(\kk [t^{\pm E_{\rho, k} }]\right)^{\times} = \kk^{\times} [t^{\pm E_{\rho, k} }]$.

\begin{cons} \label{compgen3}
We consider the map 
\begin{equation} \label{mapsigma}
\begin{aligned}
\sigma: \kk [a_{\rho, i},c^{\pm 1}_{\rho}] &\to  \kk [t^{\pm E_{\rho, k} }] \\ a_{\rho, i} &\mapsto \sigma_i (t^{r_{\rho} E_{\rho}}), \quad 1 \leq i \leq r_{\rho} -1, \ \rho \in \bar{\mathscr{P}}^{[1]} \\ c_{\rho} &\mapsto \sigma_{r_{\rho}} (t^{E_{\rho}}) = \prod_{k=1}^{r_{\rho}} t^{E_{\rho, k}}, \quad \quad \quad \ \ \rho \in \bar{\mathscr{P}}^{[1]}
\end{aligned}
\end{equation} 
where $\sigma_i(t^{r_{\rho} E_{\rho}})$ is the $i$-th symmetric polynomial in $t^{r_{\rho} E_{\rho, k}}, \ 1 \leq k \leq r_{\rho}$. Define a new collection of scattering diagrams $\bar{\mathfrak{D}}^{r\sigma} := \left\{ \bar{\mathfrak{D}}_{k}^{r\sigma}, \ k \geq 0 \right\}$ with ring $\kk[E^{\gp}]$, monoid $Q$, and ideal $I_0 = \mathfrak{m}$ by setting the walls of $\bar{\mathfrak{D}}_{k}^{r\sigma}$ to be $$(\mathfrak{p}, \sigma(f_{\mathfrak{p}}))$$ for $\mathfrak{p}$ a wall of $\bar{\mathfrak{D}}_k^r$.

Similarly, we define gluing data $s^{r\sigma}$ on $\left(\bar B, \bar{\mathscr{P}}\right)$ by setting $s^{r\sigma}_{\sigma \underline{\rho}}(m) := \sigma(s^{r}_{\sigma \underline{\rho}}(m))$ (where $s^r$ is the gluing data of \eqref{gluingdatakuniv}) for any slab $\underline{\rho} \in \tilde{\bar{\mathscr{P}}}^{[1]}$ with adjacent maximal cell $\sigma \in \bar{\mathscr{P}}^{\max}$ and any $m \in \Lambda_{\sigma}$. The gluing data $s^{r\sigma}$ is well-defined since the $t^{E_{
\rho, i}}$ are invertible in $\kk [t^{\pm E_{\rho, k} }]$. Explicitly, $s^{r\sigma}$ given by homomorphisms
\begin{equation}
\begin{aligned}
&s^{r\sigma}_{\sigma \underline{\rho}}: \Lambda_{\sigma} \to \kk^{\times} [t^{\pm E_{\rho, k} }], \quad \begin{cases} \begin{aligned} m_{\rho} &\mapsto \left(\prod_{k=1}^{r_{\rho}} t^{E_{\rho, k}}\right)^{-1} \\ \xi &\mapsto 1 \end{aligned} \end{cases} \\ &s^{r\sigma}_{\sigma \underline{\rho}'}: \Lambda_{\sigma} \to \kk^{\times} [t^{\pm E_{\rho, k} }], \quad \ \ \ m \mapsto 1 \end{aligned}
\end{equation}
using the same conventions as in \eqref{gluingdatakuniv}.

We use the MPA function $\bar{\varphi}$. Taking the inverse limit as usual, we get a family $\check{\bar{\mathfrak{X}}}^{s^{r\sigma}}_{\bar{\mathfrak{D}}^{r\sigma}} \to \Spec \kk [t^{\pm E_{\rho, k} }]\llbracket Q \rrbracket$ which is the same family as the basechange of $\check{\bar{\mathfrak{X}}}^{s^{r}}_{\bar{\mathfrak{D}}^{r}} \to \Spec \kk [a_{\rho, i}, c_{\rho}^{\pm 1}]\llbracket Q \rrbracket$ by $\sigma$. 
\end{cons}

Now that the scattering diagrams $\mathfrak{D}^r$ and $\bar{\mathfrak{D}}^{r\sigma}$ are defined over the same ring, we define $\Phi(\mathfrak{D}^r_k)'$ and $\Phi(\mathfrak{D}^r_k)$ in a way that ensures that $\Phi(\mathfrak{D}^r_0)$ is combinatorially equivalent to $\bar{\mathfrak{D}}^{r\sigma}_0$. 

\begin{cons} 
As in Construction \ref{chsgd}, we use the analogue of Construction \ref{refinement} for a scattering diagram $\bar{\mathfrak{D}}^{r\sigma'}$ on $\left(\bar{B}, \bar{\mathscr{P}}'\right)$ with gluing data $(s^{r\sigma})'$ on $\left(\bar{B}, \bar{\mathscr{P}}'\right)$ defined similarly to the $s'$ of Construction \ref{chsgd}.

We define $\Phi(\mathfrak{D}^r_k)'$ and $\Phi(\mathfrak{D}^r_k)$ in the same way as in Construction \ref{chsgd}, but using $(s^{r\sigma})'$ instead of $s'$. As in Construction \ref{chsgd}, it follows that $\Phi(\mathfrak{D}^r_k)'$ and $\Phi(\mathfrak{D}^r_k)$ are compatible with the MPA functions and $\Phi(\mathfrak{D}^r_0)$ is combinatorially equivalent to $\bar{\mathfrak{D}}^{r\sigma}_0$.
\end{cons}

Finally, we obtain a correspondence between the minimal relative Gross-Siebert locus $\check{\mathfrak{X}} \to \Spec \kk [t^{\pm E_{\rho, k} }]\llbracket Q \rrbracket$ and the subfamily $\check{\bar{\mathfrak{X}}}'_{\mathbb{P}} \to \Spec \kk'_{\mathbb{P}}[a_{\rho, i}]\llbracket Q \rrbracket$ of the universal toric degeneration mirror $\check{\bar{\mathfrak{X}}}_{\mathbb{P}} \to \Spec \kk_{\mathbb{P}}[a_{\rho, i}]\llbracket Q \rrbracket$ of \cite[Theorem A.4.2]{GHS}.

\begin{theorem} \label{compgen}
The basechange of the minimal relative Gross-Siebert locus  $$\check{\mathfrak{X}} \to \Spec \kk [t^{\pm E_{\rho, k} }]\llbracket Q \rrbracket$$ by $$r: \kk [t^{\pm E_{\rho, k} }] \to \kk [t^{\pm E_{\rho, k} }], \quad t^{E_{\rho, k}} \mapsto (t^{E_{\rho, k}})^{r_{\rho}}$$ is isomorphic to the basechange of the toric degeneration mirror family $$\check{\bar{\mathfrak{X}}}^{s^r}_{\bar{\mathfrak{D}}^r} \to \Spec \kk [a_{\rho, i}, c_{\rho}^{\pm 1}]\llbracket Q \rrbracket$$ of Construction \ref{compgen1} by 
$$
\begin{aligned}
\sigma: \kk [a_{\rho, i},c^{\pm 1}_{\rho}] &\to  \kk [t^{\pm E_{\rho, k} }] \\ a_{\rho, i} &\mapsto \sigma_i (t^{r_{\rho} E_{\rho}}), \quad 1 \leq i \leq r_{\rho} -1, \ \rho \in \bar{\mathscr{P}}^{[1]} \\ c_{\rho} &\mapsto \sigma_{r_{\rho}} (t^{E_{\rho}}) = \prod_{k=1}^{r_{\rho}} t^{E_{\rho, k}}, \quad \quad \quad \ \ \rho \in \bar{\mathscr{P}}^{[1]}
\end{aligned}
$$
Therefore, we obtain a correspondence between $\check{\mathfrak{X}} \to \Spec \kk [t^{\pm E_{\rho, k} }]\llbracket Q \rrbracket$ and $\check{\bar{\mathfrak{X}}}'_{\mathbb{P}} \to \Spec \kk'_{\mathbb{P}}[a_{\rho, i}]\llbracket Q \rrbracket$ via \eqref{mirrorfamiliescd}. 
\end{theorem}

\begin{proof}
We need to show that the families $\check{\mathfrak{X}}_{\mathfrak{D}^{r}} \to  \Spec \kk [t^{\pm E_{\rho, k} }]\llbracket Q \rrbracket$ and $\check{\bar{\mathfrak{X}}}^{s^{r\sigma}}_{\bar{\mathfrak{D}}^{r\sigma}} \to \Spec \kk [t^{\pm E_{\rho, k} }]\llbracket Q \rrbracket$ are the same. The proof is the same as the proof of Proposition \ref{proofglue} after replacing $A^{\GS}_{c_{\rho}}$ with $\kk [t^{\pm E_{\rho, k} }]$ and $s$ with $s^{r\sigma}$.
\end{proof}

If the dual intersection complex $\left(\bar B, \bar{\mathscr{P}}\right)$ of $\bar{\mathfrak{X}} \to \mathcal{S}$ is simple, the correspondence is more direct. Note that in this case we have $\kk[a_{\rho,i}] \cong \kk$ since $K := \sum_{\rho \in \bar{\mathscr{P}}^{[1]}, r_{\rho} > 0} (r_{\rho} -1) = 0$. We also have $\kk [t^{\pm E_{\rho, k} }] \cong \kk [t^{\pm E_{\rho, 1} }]$.

\begin{cor} \label{compgensimple}
Suppose that $\left(\bar B, \bar{\mathscr{P}}\right)$ is simple. Let $\kk'_{\mathbb{P}} \subseteq \kk_{\mathbb{P}}$ be the subring generated by $\bar{s}$ for $s$ as in \eqref{gluingdatak} and let $\check{\bar{\mathfrak{X}}}'_{\mathbb{P}} \to \Spec \kk'_{\mathbb{P}}\llbracket Q \rrbracket$ be the corresponding subfamily of the universal toric degeneration mirror $\check{\bar{\mathfrak{X}}}_{\mathbb{P}} \to \Spec \kk_{\mathbb{P}}\llbracket Q \rrbracket$ of \cite[Theorem A.2.4]{GHS}. 
The map 
$$c_{\univ}: \kk [c_{\rho}^{\pm 1}] \to \kk'_{\mathbb{P}},  \quad c_{\rho} \mapsto \bar{s}^{\rho}$$ (where the gluing data $s^{\rho}$ is defined similarly to $s^{r,\rho}$ in Construction \ref{compgen1} and $\bar{s}^{\rho}$ is the corresponding closed gluing data) is well-defined and surjective.  

Let 
$\check{\bar{\mathfrak{X}}}^{s}_{\bar{\mathfrak{D}}} \to \Spec \kk [c_{\rho}^{\pm 1}]\llbracket Q \rrbracket$
be the toric degeneration mirror family constructed using the gluing data $s$ of \eqref{gluingdatak} (viewing $c_{\rho}$ as a variable) and let $\sigma$ be the map
\begin{equation} \label{samsigmasimple}
\sigma: \kk [c^{\pm 1}_{\rho}] \to  \kk [t^{\pm E_{\rho, 1} }], \quad  c_{\rho} \mapsto t^{E_{\rho, 1}}.
\end{equation} 
Note that $\sigma$ is not an isomorphism in general, see the discussion after \eqref{egp}. Let $\check{\mathfrak{X}} \to \Spec \kk [t^{\pm E_{\rho, 1} }]\llbracket Q \rrbracket$ be the minimal relative Gross-Siebert locus. We have a commutative diagram as follows:
\begin{equation}
    \begin{tikzcd}
     \check{\mathfrak{X}} \arrow[r] \arrow{d} & \check{\bar{\mathfrak{X}}}^{s}_{\bar{\mathfrak{D}}} \arrow[leftarrow]{r} \arrow{d} & \check{\bar{\mathfrak{X}}}'_{\mathbb{P}} \arrow[hookrightarrow]{r} \arrow{d} & \check{\bar{\mathfrak{X}}}_{\mathbb{P}}  \arrow{d}  \\
     \Spec \kk [t^{\pm E_{\rho, 1} }]\llbracket Q \rrbracket \arrow{r}{\sigma} & \Spec \kk[c_{\rho}^{\pm1}]\llbracket Q \rrbracket \arrow[leftarrow]{r}{c_{\univ}} & \Spec \kk'_{\mathbb{P}}\llbracket Q \rrbracket \arrow[hookrightarrow]{r} & \Spec \kk_{\mathbb{P}}\llbracket Q \rrbracket   
    \end{tikzcd}
\end{equation}
\end{cor}

\begin{proof}
Since $\left(\bar B, \bar{\mathscr{P}}\right)$ is simple, the universal toric degeneration mirror of \cite[Theorem A.4.2]{GHS} is just the universal toric degeneration mirror of \cite[Theorem A.2.4]{GHS}. The maps $r$ of \eqref{mapr} and \eqref{maprintrinsic} and the map $r_{\univ}$ of \eqref{mapruniv} are trivial since $r_{\rho} \leq 1$ for all $\rho \in \bar{\mathscr{P}}^{[1]}$. The map $\sigma$ of \eqref{mapsigma} reduces to the map $\sigma$ of \eqref{samsigmasimple}. So Theorem \ref{compgen} reduces to the claimed result.   
\end{proof}

\begin{remark}
Even in the case of simple $\left(\bar B, \bar{\mathscr{P}}\right)$, there is no isomorphism between $\check{\mathfrak{X}} \to \Spec \kk [t^{\pm E_{\rho, 1} }]\llbracket Q \rrbracket$ and $\check{\bar{\mathfrak{X}}}'_{\mathbb{P}} \to \Spec \kk'_{\mathbb{P}}\llbracket Q \rrbracket$ in general. Indeed, we have $\kk'_{\mathbb{P}} \subseteq \kk [H^1(\bar{B}, i_{*} \check{\Lambda})_f^{*}]$ and $\dim H^1(\bar{B}, i_{*} \check{\Lambda})_f^{*} = 20$ (see \cite[Example 5.22(2)]{GS2}). On the other hand, there are examples when $\dim E^{\gp} = \sum_{\rho \in \bar{\mathscr{P}}^{[1]}} r_{\rho} =24$ (see the discussion after \eqref{egp} and \cite[Example 5.22(2)]{GS2}).
\end{remark}

\subsection{Discussion of the results} \label{discussion}

We highlight certain features of our construction along with some philosophy. 

Note that the fact that there is an algorithmic construction for the toric degeneration mirror means that we might be able to compute certain intrinsic mirrors to minimal log CY resolutions of special toric degenerations of K3-s explicitly (and in particular, compute certain punctured log Gromov-Witten invariants) similarly to the calculations of \cite{A} (also see Observation \ref{imsmireq} for the equation of the intrinsic mirror to a small resolution of the toric degeneration of Example \ref{runningexample}). For this reason, it would be interesting to study minimal log CY degenerations of K3-s that admit a blowdown to a special toric degeneration. 

We now discuss relaxing the underlying assumptions, the minimal relative Gross-Siebert locus modulo $I_0$ and mirrors to generically log smooth resolutions, and gluing data for intrinsic mirrors.

\subsubsection{Relaxing the underlying assumptions} \label{relax}

We have made certain restricting assumptions throughout the paper. First, we have Assumption \ref{toricdegenerationass} on the toric degeneration $\bar{\mathfrak{X}} \to \mathcal{S}$. The assumption that $\bar{\mathfrak{X}}$ is a variety is only needed to ensure that so is $\mathfrak{X}$ (for a resolution $\pi:\mathfrak{X} \to \bar{\mathfrak{X}}$ that we consider), which is necessary to fit with the assumptions of the literature \cite{ACGS, ACGS2} on punctured log Gromov-Witten invariants. As explained in \cite[Remark 5.16(2)]{GHKS}, one can work with an algebraic space $\mathfrak{X}$ instead, and this only requires a slight modification of Definition \ref{walltype} of a decorated wall type $\btau$. Assumptions \ref{toricdegenerationass}(3) and \ref{toricdegenerationass}(4) ensure that both $\left(\bar{B}, \bar{\mathscr{P}}\right)$ and $(B, \mathscr{P})$ are polyhedral manifolds in the sense of Definition \ref{polymani}. These assumptions can likely be removed by working with the toric polyhedral decompositions of \cite[Definitions 1.21 and 1.22]{GS2} and revising \cite{GS3, GS4, GHS} for this more general setup, see Remark \ref{polymaniremarks}(2). The assumption that $\bar{\mathfrak{X}} \to \mathcal{S}$ is projective (and not just proper) is only needed to ensure that so is $\mathfrak{X} \to \mathcal{S}$, a condition we now discuss. 

It may sometimes seem more natural to work with non-projective but proper resolutions $\mathfrak{X} \to \mathcal{S}$. For instance, consider the case of Example \ref{runningexample}. We can blow up two singular points of $D_i \cap D_j$ in $D_i$ and the other two points in $D_j$ for $1 \leq i < j \leq 4$ to obtain a more symmetric resolution $\mathfrak{X} \to \Spec \kk\llbracket t \rrbracket$. There are certain obstacles to extending the results of this thesis to such a resolution. 

It is clear that $\mathfrak{X} \to \Spec \kk\llbracket t \rrbracket$ should be the classical degeneration of a Calabi-Yau hypersurface to four rational surfaces that are $\mathbb{P}^2$-s with $6$ boundary points blown up. However, such $\mathfrak{X} \to \Spec \kk\llbracket t \rrbracket$ can't be obtained as a sequence of blowups of irreducible components of the central fibre. Moreover, by considering the intersection numbers of $E_{\rho,k}, \ \rho \in \bar{\mathscr{P}}^{[1]}, 1 \leq k \leq r_{\rho}$ with the irreducible components of $\mathfrak{X}_0$ it is easy to see that there is no relatively ample divisor $D'$ supported on $D$ for this resolution. Note that we actually prove Conjecture \ref{conjecture} for any toric, integral, and homogeneous resolution (see Remark \ref{stronglyadmissibleweakening}(1)), so this is not crucial. More seriously, the projectivity assumption is currently required in \cite{GS4} to ensure that the canonical scattering diagram is consistent. 

This discussion generalizes to the setup of the tropical approach to resolutions of Section \ref{tropicalresolutions} and corresponds to considering non-homogeneous resolutions, see Remark \ref{otherresolutions}. These resolutions are often not projective. Generalizing the results of this thesis to toric, integral, and non-homogeneous projective resolutions should be straightforward. This only requires a generalization of Section \ref{scatint}, as explained in Remark \ref{stronglyadmissibleweakening}(2), and a corresponding generalization of Construction \ref{chs}.  

We believe that the results of this thesis can also be extended to certain proper, non-projective resolutions. One reason this is desirable is the results of \cite{H} that construct log smooth families directly from polytopes (with certain specified decompositions), similarly to how we constructed Batyrev degenerations in Section \ref{batyrev}. These families are often non-projective but proper. It would be interesting to see how the mirror Batyrev degenerations of Section \ref{batyrev} are related to the intrinsic mirrors to the degenerations of \cite{H}.

The assumption that a strongly admissible resolution $\pi:\mathfrak{X} \to \bar{\mathfrak{X}}$ (see Definition \ref{stronglyadmissible}) is integral is necessary to guarantee that $(B, \mathscr{P})$ is an integral subdivision of $\left( \bar{B}, \bar{\mathscr{P}}\right)$, which implies that $(B, \mathscr{P})$ is a polyhedral manifold. In particular, the central fibre of $\mathfrak{X} \to \mathcal{S}$ is reduced, see Remark \ref{nonreducedeness}. We chose to require integrality of the polytopes in the polyhedral decomposition $\mathscr{P}$ on $B$ in Definition \ref{polymani} of a polyhedral manifold to be in line with the conventions of \cite{GHS}. However, as we mention in Remark \ref{nonreducedeness}, one can allow rational polyhedra in the construction, and all the results still make sense (this is the approach of \cite{GHKS}). Therefore, one may weaken the requirement that a strongly admissible resolution is integral and allow rational subdivisions of $\sigma \in \bar{\mathscr{P}}^{\max}$ (they should still induce integral subdivisions of $\rho \in \bar{\mathscr{P}}^{[1]}$ since in Definition \ref{toricresolutions} we require compatibility with the subdivisions of the cones defining the local models for the singular points $x \in \bar{X}_{\rho}$). This extends the results of this thesis, particularly those of Chapter \ref{gslocus} where we only work with strongly admissible resolutions.

We discuss generalizing Conjecture \ref{conjecture} to non-special toric degenerations $\bar{\mathfrak{X}} \to \mathcal{S}$ of K3-s in Remarks \ref{basechangeconj}. In particular, Remark \ref{basechangeconj}(1) gives a positive result in this direction.

\subsubsection{The minimal relative Gross-Siebert locus modulo \texorpdfstring{$I_0$}{I0} and mirrors to generically log smooth degenerations} \label{gslocusstuff}

In Theorem \ref{compgen}, we proved that there is a correspondence between the minimal relative Gross-Siebert locus and a subfamily of the universal toric degeneration mirror family of \cite[Theorem A.4.2]{GHS}. Proposition \ref{proofglue} and Remark \ref{glueremarkcorr} imply that we can obtain certain one-parameter toric degeneration mirrors from the minimal relative Gross-Siebert locus $\check{\mathfrak{X}} \to \Spec \kk [t^{\pm E_{\rho, k} }]\llbracket Q \rrbracket$ by basechange. In particular, $\Spec \kk[t^{\pm E_{\rho, k}}]$ can be viewed as the parameter space of toric log CY structures on $\check{\bar{X}}_0$ of a certain form.

\begin{fact} \label{fibres}
The minimal relative Gross-Siebert locus modulo $I_0 = Q \setminus \left\{0\right\}$ is the family $\check{\mathfrak{X}}_{\mathfrak{D}^{\GS}_{Q}} \to \Spec \kk [t^{\pm E_{\rho, k} }]$ (here $\mathfrak{D}^{\GS}_{Q}$ is defined in Construction \ref{gsscattering}). It is a trivial deformation of schemes (i.e. $\check{\mathfrak{X}}_{\mathfrak{D}^{\GS}_{Q}} \cong \check{\bar{X}}_0 \times \Spec \kk [t^{\pm E_{\rho, k} }]$). Suppose further that $\dim E^{\gp} = \sum_{\rho \in \bar{\mathscr{P}}^{[1]}} r_{\rho}$. Then all the toric log CY structures on $\check{\bar{X}}_0$ with gluing data in the subring $\kk'_{\mathbb{P}} \subseteq \kk [H^1(\bar{B}, i_{*} \check{\Lambda})_f^{*}]$ parametrizing closed gluing data generated by $\bar{s}$ for $s$ as in \eqref{gluingdatak} appear as fibres of $\check{\mathfrak{X}}_{\mathfrak{D}^{\GS}_{Q}} \to \Spec \kk [t^{\pm E_{\rho, k} }]$. Here the log structures on the fibres are induced from the natural log structure on $\check{\mathfrak{X}}_{\mathfrak{D}^{\GS}_{Q}} \to \Spec \kk [t^{\pm E_{\rho, k} }]$ (see Appendix \ref{log}) by inclusions.
\end{fact}

\begin{proof}
To show that 
$\check{\mathfrak{X}}_{\mathfrak{D}^{\GS}_{Q}} \to \Spec \kk [t^{\pm E_{\rho, k} }]$ is a trivial deformation of the central fibre, it is enough to show that the family $\check{\mathfrak{X}}^o_{\mathfrak{D}^{\GS}_{Q}}$ outside codimension 2 is a trivial deformation (by Proposition \ref{help} and since the empty scattering diagram provides the trivial deformation). $\check{\mathfrak{X}}^o_{\mathfrak{D}^{\GS}_{Q}}$ is obtained by gluing together the affine schemes $\Spec R_{\mathfrak{b}}$ with 
\begin{equation} \label{localeq}
R_{\mathfrak{b}} = (\kk [E^{\gp} \oplus Q]/I_0)[\Lambda_{\rho}][Z_+, Z_-]/(Z_+Z_- - f_{\mathfrak{b}} \cdot z^{\kappa_{\rho}})
\end{equation}
for choices of slabs $\mathfrak{b} \subseteq \rho \in \mathscr{P}^{[1]}$ for every $\rho \in \mathscr{P}^{[1]}$ along the $\Spec R_{\mathfrak{u}}$ for choices of chambers $\mathfrak{u} \subseteq \sigma \in \mathscr{P}^{\max}$ for every $\sigma \in \mathscr{P}^{\max}$ (see the proof of Proposition \ref{equivalence}). Note that modulo $I_0$, all the non-trivial walls are slabs. 

Construction \ref{gsscattering} implies that for every non-trivial slab $\mathfrak{b} \subseteq \rho \in \mathscr{P}^{[1]}$, we have $\kappa_{\rho} := [X_{\rho}]$. Then either $X_{\rho}$ is not contracted by $\pi:\mathfrak{X} \to \bar{\mathfrak{X}}$ and $z^{\kappa_{\rho}} = 0 \in \kk [E^{\gp} \oplus Q]/I_0$, or $\rho \in \mathscr{P}^{[1]} \setminus \mathscr{P}^{[1]}_{\coar}$ (using Notation \ref{nonsubdividing}) and $z^{\kappa_{\rho}} = f_{\mathfrak{b}}=1$. Thus  $\Spec R_{\mathfrak{b}}$ are trivial deformations of their counterparts over the maximal ideal. From the construction of $\check{\mathfrak{X}}^o_{\mathfrak{D}^{\GS}_{Q}}$ in \cite[Proposition 2.4.1]{GHS}, we immediately see that the gluing respects this triviality and produces a trivial deformation.

A toric log CY structure with gluing data $s$ as in \eqref{gluingdatak} specified by $c_{\rho} \in \kk^{\times}$ is determined by a choice of the initial slab functions as in \eqref{slabfuncgd} with $a_{\rho, i} \in \kk$ and $c_{s, \underline{\rho}} = c_{\rho}, \ c_{s, \underline{\rho}'} = 1$ for $\rho \in \bar{\mathscr{P}}^{[1]}$. Proposition \ref{proofglue} and Remark \ref{glueremarkcorr} imply that $\check{\bar{X}}_0$ with such a toric log CY structure is isomorphic as a log scheme to the fibre of $\check{\mathfrak{X}}_{\mathfrak{D}^{\GS}_{Q}} \to \Spec \kk [t^{\pm E_{\rho, k} }]$ (with the log structure induced from the natural log structure on $\check{\mathfrak{X}}_{\mathfrak{D}^{\GS}_{Q}} \to \Spec \kk [t^{\pm E_{\rho, k} }]$ by inclusion) over the point 
\begin{equation} \label{solvefort}
\left\{t^{E_{\rho, k}} = \zeta_{\rho, k} \ | \ \rho \in \bar{\mathscr{P}}^{[1]}, \  1 \leq k \leq r_{\rho} \right\}
\end{equation}
where $\zeta_{\rho, k}$ is the $k$-th root of 
\begin{equation} \label{polyeq}
1 + a_{\rho, 1}w_{\rho} + \dots + a_{\rho, r_{\rho}-1}w_{\rho}^{r_{\rho}-1} + c_{\rho} w_{\rho}^{r_{\rho}} = 0
\end{equation}
(for any ordering of the roots, note that $c_{\rho} \in \kk^{\times}$ implies that $\zeta_{\rho, k} \neq 0$ for all $1 \leq k \leq r_{\rho}$). Here the fact that $\dim E^{\gp} = \sum_{\rho \in \bar{\mathscr{P}}^{[1]}} r_{\rho}$ implies that the system \eqref{solvefort} has a solution in $E^{\gp}$. 
\end{proof}

\begin{remarks} \label{fibresremarks}
\begin{enumerate}
    \item Note that the correspondence between the fibres of $\check{\mathfrak{X}}_{\mathfrak{D}^{\GS}_{Q}} \to \Spec \kk [t^{\pm E_{\rho, k} }]$ and toric log CY structures on $\check{\bar{X}}_0$ is not one-to-one. If the equations \eqref{polyeq} have no multiple roots for all $\rho \in \bar{\mathscr{P}}^{[1]}$, then there are at least $\prod_{\rho \in \bar{\mathscr{P}}^{[1]}} r_{\rho}!$ fibres of $\check{\mathfrak{X}}_{\mathfrak{D}^{\GS}_{Q}} \to \Spec \kk [t^{\pm E_{\rho, k} }]$ corresponding to the toric log CY structure on $\check{\bar{X}}_0$ specified by $c_{\rho} \in \kk^{\times}$ and $a_{\rho, i} \in \kk$ (note that some of these log CY structures are equivalent). 
    This follows from the fact that the basechange by $\sigma$ in Construction \ref{basechangeglue} gives a covering of index $\prod_{\rho \in \bar{\mathscr{P}}^{[1]}} r_{\rho}!$. 
    \item If $\dim E^{\gp} < \sum_{\rho \in \bar{\mathscr{P}}^{[1]}} r_{\rho}$, then the system \eqref{solvefort} does not always have a solution in $E^{\gp}$. Therefore, $\check{\mathfrak{X}}_{\mathfrak{D}^{\GS}_{Q}} \to \Spec \kk [t^{\pm E_{\rho, k} }]$ only parametrizes the log CY structures on $\check{\bar{X}}_0$ corresponding to choosing $c_{\rho} \in \kk^{\times}$ and $a_{\rho, i} \in \kk$ so that \eqref{solvefort} has a solution.
\end{enumerate}
\end{remarks}

Suppose that we are in the case of a small resolution $\pi: \mathfrak{X} \to \bar{\mathfrak{X}}$ of a special toric degeneration $\bar{\mathfrak{X}} \to \mathcal{S}$ of Section \ref{resolvesmall} (e.g. the case of Example \ref{runningexample}). Then \eqref{canonicalsplitting} implies that the intrinsic mirror $\check{\mathfrak{X}} \to \Spec \kk [t^{ E_{\rho, k} }]\llbracket Q \rrbracket$ is defined over $\kk [E] = \kk[t^{E_{\rho, k}}]$. 
By Theorem \ref{compgen}, we have a correspondence between the restriction of the intrinsic mirror to the minimal relative Gross-Siebert locus $\check{\mathfrak{X}} \to \Spec \kk [t^{\pm E_{\rho, k} }]\llbracket Q \rrbracket$ (the family over the large torus of $\Spec \kk [t^{ E_{\rho, k} }]$) and a certain universal toric degeneration mirror. It is natural to ask what the restrictions of the intrinsic mirror to families over deeper toric strata of $\Spec \kk [t^{ E_{\rho, k} }]$ correspond to. Since $t^{E_{\rho, k}}=0$ implies that $\sigma_{r_{\rho}} (t^{E_{\rho}}) = \prod_{k=1}^{r_{\rho}} t^{E_{\rho, k}} = 0$, the analogue of the constructions of Chapters \ref{proof} and \ref{gslocus} can't work for such strata since $\prod_{i=1}^{r_{\rho}} (1+ t^{E_{\rho, k}}w_{\rho})$ has coefficient $0$ at $w_{\rho}^{r_{\rho}}$, so it is not of the form $f_{\underline{\rho}}$ in \eqref{slabfuncgd}. In particular, restricting modulo $I_0$, a fibre of $\check{\mathfrak{X}} \to \Spec \kk [t^{ E_{\rho, k} }]$ with some $t^{E_{\rho, k}}=0$ and with the log structure induced from the natural log structure on $\check{\mathfrak{X}} \to \Spec \kk [t^{ E_{\rho, k} }]$ (see Appendix \ref{log}) by inclusion is not toric log CY (it follows from the local equations \eqref{localeq} that the singularities of the log structure fall into codimension $\geq 2$). So the restrictions of the intrinsic mirror to deeper strata are not toric degenerations.

The original non-extended intrinsic mirror is a family $\check{\mathfrak{X}} \to \Spec \kk \llbracket t^{ E_{\rho, k} } \rrbracket\llbracket Q \rrbracket$ defined over the completion of the minimal stratum $\left\{0\right\} \subseteq \Spec \kk [t^{ E_{\rho, k} }]$. Recall the generically log smooth families of Definition \ref{genlogsmooth}. The mirror to a generically log smooth partial resolution $\mathfrak{X}' \to \mathcal{S}$ of $\bar{\mathfrak{X}} \to \mathcal{S}$ should correspond to the restriction of $\check{\mathfrak{X}} \to \Spec \kk [t^{ E_{\rho, k} }]\llbracket Q \rrbracket$ to the completion of the stratum corresponding to the resolved $E_{\rho, k}$. This suggests that there should be a general construction of mirrors to generically log smooth degenerations.

From this point of view, the extended intrinsic mirror $\check{\mathfrak{X}} \to \Spec \widehat{\kk [P]}_{I_{\min}}$ of \eqref{extendedmirror} obtained using the ideal $I_{\min}$ of Proposition \ref{imin} can be viewed as a universal mirror such that the restrictions to the completions of the toric strata of $\Spec \widehat{\kk [P]}_{I_{\min}}$ correspond to mirrors to the families obtained by contracting the corresponding curves in the central fibre of $\mathfrak{X} \to \mathcal{S}$.
It is not clear how to construct such contractions or if they can always be constructed. Some contractions may be non-projective or only defined as algebraic spaces. 

This philosophy extends to the case when the resolution is not small. In that case, we have some exceptional divisors. In Section \ref{minrelgslocus}, we restricted to the minimal relative Gross-Siebert locus, which is the family over the stratum corresponding to setting all the curve classes except for $E_{\rho, k}, \ \rho \in \mathscr{P}^{[1]}, \ 1 \leq k \leq r_{\rho}$ to $0$. The reason that we restricted to this locus is that we were only interested in the mirror to $\bar{\mathfrak{X}} \to \mathcal{S}$ and the curves $E_{\rho, k}, \ \rho \in \mathscr{P}^{[1]}, \ 1 \leq k \leq r_{\rho}$ correspond to the singularities of $\bar{\mathfrak{X}} \to \mathcal{S}$. This restriction smooths out the components of the intrinsic mirror corresponding to the other curve classes.
If we include all the strata, $\check{\mathfrak{X}} \to \Spec \widehat{\kk [P]}_{I_{\min}}$ should parameterize the mirrors to generically log smooth contractions in this case as well (these contractions don't need to be small). In particular, there may be multiple non-isomorphic toric degenerations among the contractions where we know how to construct mirrors, and the described property is just Theorem \ref{compgen}. To summarize, contracting on the side of the minimal log CY degeneration $\mathfrak{X} \to \mathcal{S}$ should correspond to smoothing out the components on the mirror side $\check{\mathfrak{X}} \to \Spec \widehat{\kk [P]}_{I_{\min}}$.

Conversely, Theorem \ref{compgen} implies that for any two strongly admissible resolutions $\pi_1: \mathfrak{X}_1 \to \bar{\mathfrak{X}}$ and $\pi_2: \mathfrak{X}_2 \to \bar{\mathfrak{X}}$ of a toric degeneration $\bar{\mathfrak{X}} \to \mathcal{S}$ of K3-s, the minimal relative Gross-Siebert locus $\check{\mathfrak{X}}_1 \to \Spec \kk [t^{\pm E_{\rho, k} }]\llbracket Q \rrbracket$ agrees with the minimal relative Gross-Siebert locus $\check{\mathfrak{X}}_2 \to \Spec \kk [t^{\pm E_{\rho, k} }]\llbracket Q \rrbracket$. So the difference between the two mirrors comes from the curve classes that are not common between them. In the case of non-small resolutions, the mirrors might have different components that smooth out to the same components of the minimal relative Gross-Siebert locus. In the case of small resolutions, $\mathfrak{X}_1$ and $\mathfrak{X}_2$ are related by a sequence of flops which correspond to swapping the wall functions attached to the corresponding slabs of $\mathfrak{D}_J$.

It would be interesting to find a mirror construction for generically log smooth degenerations generalizing both the toric degeneration mirror construction \cite{GS1} and the intrinsic mirror construction  \cite{GS3, GS4} and make the philosophy of this section precise.

\subsubsection{Gluing data for intrinsic mirrors} \label{gdint}

In Section \ref{glue}, we used gluing data on $\left(\bar B, \bar{\mathscr{P}}\right)$ to obtain the correspondence between the minimal relative Gross-Siebert locus and the restriction of the universal toric degeneration mirror of \cite[Theorem A.4.2]{GHS} to the subfamily $\kk'_{\mathbb{P}} \subseteq \kk_{\mathbb{P}}$ of gluing data generated by $\bar{s}$ for $s$ as in \eqref{gluingdatak}. It is natural to wonder if we can extend this correspondence to a larger subfamily of the universal toric degeneration mirror of \cite[Theorem A.4.2]{GHS} by introducing gluing data on $(B, \mathscr{P})$ and modifying the wall functions of Construction \ref{canonical scattering} of the canonical scattering diagram to account for gluing data (changing Construction \ref{extendcanonicalscat} of $\mathfrak{D}^{J}$ accordingly). Introducing gluing data into the construction of the intrinsic mirror does not seem very natural and requires extending the results of \cite{GS4} to the case of gluing data. Therefore, we refrain from giving a formal exposition and just outline the possible extension.

First, we introduce gluing data on $(B, \mathscr{P})$. We work over the ring $\kk [t^{\pm E_{\rho, k} }]$ as in Section \ref{compgencons}. Since there are no singularities at the interiors of the edges, we can simplify Definition \ref{gddef} and define gluing data $s$ on $(B, \mathscr{P})$ as a collection of homomorphisms of abelian groups
$s_{\sigma \rho'}: \Lambda_{\sigma} \to \left(\kk [t^{\pm E_{\rho, k} }]\right)^{\times} = \kk^{\times} [t^{\pm E_{\rho, k} }]$ for every $\rho' \in \mathscr{P}^{[1]}$ and an adjacent maximal cell $\sigma \in \mathscr{P}^{\max}$ subject to the constraint 
\begin{equation} \label{slabrestint}
s_{\sigma \rho'}|_{\Lambda_{\rho'}} = s_{\sigma' \rho'}|_{\Lambda_{\rho'}}
\end{equation}
if $\rho' = \sigma \cap \sigma'$ for $\sigma, \sigma' \in \bar{\mathscr{P}}^{\max}$.

Now, consider the map $\sigma: \kk [a_{\rho, i},c^{\pm 1}_{\rho}] \to  \kk [t^{\pm E_{\rho, k} }] $ of Construction \ref{compgen3} and let $\kk[\sigma_i]: = \im \sigma \cong \kk [a_{\rho, i},c^{\pm 1}_{\rho}]$. Note that we have $$\left(\kk[\sigma_i]\right)^{\times} = \kk^{\times} \left[\left(\prod_{k=1}^{r_{\rho}} t^{E_{\rho, k}}\right)^{\pm 1} \ \Bigg| \ \rho \in \bar{\mathscr{P}}^{[1]}\right].$$ We shall only consider gluing data $s$ on $(B, \mathscr{P})$ with all the $
s_{\sigma \rho'}: \Lambda_{\sigma} \to \kk^{\times} [t^{\pm E_{\rho, k} }]$ taking values in $(\kk[\sigma_i])^{\times} \subseteq \kk^{\times} [t^{\pm E_{\rho, k} }]$. Further, we require that $s_{\sigma_1 \rho_{p_1}} = s_{\sigma_2 \rho_{p_2}}:= s_{\sigma, \rho}$ for any $\rho_{p_1}, \rho_{p_2} \in \mathscr{P}^{[1]}$ supported on the same $\rho \in \mathscr{P}^{[1]}_{\coar}$ and contained in the maximal cells $\sigma_1, \sigma_2 \in \mathscr{P}^{\max}$ with $\sigma_1, \sigma_2 \subseteq \sigma \in \mathscr{P}^{\max}_{\coar}$ and require the gluing data to be trivial on any $\rho' \in \mathscr{P}^{[1]} \setminus \mathscr{P}^{[1]}_{\coar}$ (using Notation \ref{nonsubdividing}). Indeed, only such gluing data corresponds to gluing data on $\left(\bar B, \bar{\mathscr{P}}\right)$ with values in $(\kk [a_{\rho, i},c^{\pm 1}_{\rho}])^{\times} = \kk^{\times} [c^{\pm 1}_{\rho}]$. Given a gluing data $s$ on $(B, \mathscr{P})$ with values in $(\kk[\sigma_i])^{\times}$, we define the corresponding gluing data $s'$ on $\left(\bar B, \bar{\mathscr{P}}\right)$ with values in $ \kk^{\times} [c^{\pm 1}_{\rho}]$ by setting, for every two slabs $\underline{\rho}, \underline{\rho}' \subseteq \rho \in \bar{\mathscr{P}}^{[1]}$ and an adjacent maximal cell $\sigma \in \bar{\mathscr{P}}^{\max}$:
\begin{equation} \label{gdcorr}
s'_{\sigma \underline{\rho}} = s'_{\sigma \underline{\rho}'} := s_{\sigma \rho}
\end{equation}
(using the isomorphism $\kk[\sigma_i] \cong \kk [a_{\rho, i},c^{\pm 1}_{\rho}]$ and the same notations $\rho, \sigma$ for the corresponding cells of $\mathscr{P}_{\coar}$). 

Unless $s$ is a trivial gluing data, the canonical scattering diagram $\mathfrak{D}_{I}$ of Construction \ref{canonical scattering} is not consistent for $s$ and the wall functions of \eqref{refinedwalleq} need to be modified to account for gluing data. We set the new wall functions to be $$f^s_{\mathfrak{p}_{\btau}} := \left(\prod_{\rho \in \mathscr{P}^{[1]}, \Int \rho \cap \mathfrak{p}_{\btau} \neq \varnothing } s_{\sigma \rho}  \right) (f_{\mathfrak{p}_{\btau}})$$ where $\sigma$ is the maximal cell containing $\mathfrak{p}_{\btau}$ (if $\mathfrak{p}_{\btau}$ is a slab, then we choose any maximal cell containing $\mathfrak{p}_{\btau}$ and the choice does not matter by \eqref{slabrestint}), $f_{\mathfrak{p}_{\btau}}$ is the wall function of \eqref{refinedwalleq}, and we set $s_{\sigma \rho}(z^m):= s_{\sigma \rho}(m) z^m$ as in \eqref{gdhomomorphism}. One needs to check that this modification defines a scattering diagram $\mathfrak{D}^s_I$ consistent for $s$ by modifying the argument of \cite{GS4} to account for gluing data. The modification also requires changing Construction \ref{extendcanonicalscat} of $\mathfrak{D}^J$ accordingly. Now, one should be able to construct a universal in gluing data extended intrinsic mirror $\check{\mathfrak{X}} \to \Spec \widehat{\kk^{\IMS}_{\mathbb{P}} [P]}_J$, where $\kk^{\IMS}_{\mathbb{P}}$ is the universal choice of a subfamily of projective closed gluing data with values in $(\kk[\sigma_i])^{\times}$ and satisfying the requirements above, by using the techniques similar to \cite[Appendix A.2]{GHS}. One can then define the minimal relative Gross-Siebert locus varied in gluing data $\check{\mathfrak{X}} \to \Spec \kk^{\IMS}_{\mathbb{P}} [t^{\pm E_{\rho, k} }]\llbracket Q \rrbracket$ in the same way as in Definition \ref{gslocusdef}.

Consider the subfamily $$ \kk^{\TD}_{\mathbb{P}} := \left\langle \bar{s}' \cdot \bar{s}'' \ | \ s' \text{ as in } \eqref{gdcorr}, \ s'' \text{ as in } \eqref{gluingdatak} \right\rangle \subseteq \kk [H^1(\bar{B}, i_{*} \check{\Lambda})_f^{*}]$$ of projective gluing data on $\left(\bar B, \bar{\mathscr{P}}\right)$. It is easy to see that $\kk^{\TD}_{\mathbb{P}}$ is a subring of the universal ring $\kk_{\mathbb{P}}$ parametrizing projective gluing data used in the construction of the universal toric degeneration mirror of \cite[Theorem A.4.2]{GHS}. Arguing similarly to Section \ref{compgencons}, one should be able to prove the following.

\begin{sloppypar}
\begin{conj} \label{conjclueims}
The correspondence of Theorem \ref{compgen} naturally extends to the minimal relative Gross-Siebert locus varied in gluing data $\check{\mathfrak{X}} \to \Spec \kk^{\IMS}_{\mathbb{P}} [t^{\pm E_{\rho, k} }]\llbracket Q \rrbracket$ and the restriction of the universal toric degeneration mirror of \cite[Theorem A.4.2]{GHS} to $\Spec \kk^{\TD}_{\mathbb{P}}\llbracket Q \rrbracket \subseteq \Spec \kk_{\mathbb{P}}\llbracket Q \rrbracket$.
\end{conj}
\end{sloppypar}

\begin{remark}
One may wonder if (possibly under some additional assumptions) one has $\kk^{\TD}_{\mathbb{P}} = \kk_{\mathbb{P}}$ (or even $\kk'_{\mathbb{P}} = \kk_{\mathbb{P}}$) so that Conjecture \ref{conjclueims} (resp. Theorem \ref{compgen}) gives a correspondence with the whole universal toric degeneration mirror of \cite[Theorem A.4.2]{GHS}. Answering this question requires a better understanding of the structure of the ring $\kk [H^1(\bar{B}, i_{*} \check{\Lambda})_f^{*}]$ parametrizing closed gluing data.
\end{remark}

\section{Generalizing to higher dimensions} \label{higherdimgen}

We outline the challenges of generalizing the results of this thesis to higher dimensions which is a subject of ongoing work. To prove Conjecture \ref{conjecture} for a special toric degeneration $\bar{\mathfrak{X}} \to \mathcal{S}$ of relative dimension $n \geq 3$ we need to:
\begin{enumerate}
    \item Define strongly admissible and admissible resolutions $\pi: \mathfrak{X} \to \bar{\mathfrak{X}}$ of $\bar{\mathfrak{X}} \to \mathcal{S}$ and show that $\bar{\mathfrak{X}} \to \mathcal{S}$ admits a strongly admissible resolution, generalizing the results of Section \ref{admissible}. For any well-chosen monoid $P$ (see Definition \ref{wellchosendef}) and $J := P \setminus K$ the complement of the face containing the curve classes contracted by $\pi$ as before, we obtain the extended intrinsic mirror $\check{\mathfrak{X}} \to \Spec \widehat{\kk [P]}_J$ via \eqref{extendedmirror} (alternatively, we may define the extended intrinsic mirror via the collection $\mathfrak{D}^J$ of scattering diagrams as in (2) to allow a more general definition of a strongly admissible resolution).
    \item Provide a collection of scattering diagrams $\mathfrak{D}^J: = \left\{ \mathfrak{D}_{J^{k+1}}, \ k \geq 0  \right\}$ giving rise to $\check{\mathfrak{X}} \to \Spec \widehat{\kk [P]}_J$, generalizing the results of Section \ref{scatint}.
    \item Relate the collection $\mathfrak{D}^J$ of the canonical scattering diagrams giving rise to $\check{\mathfrak{X}} \to \Spec \widehat{\kk [P]}_J$ to the collection $\bar{\mathfrak{D}} = \left\{ \bar{\mathfrak{D}}_k, \ k \geq 0 \right\}$ of the algorithmic scattering diagrams giving rise to the toric degeneration mirror $\check{\bar{\mathfrak{X}}} \to \Spec \kk\llbracket t \rrbracket$ via a PL-isomorphism $\Phi: \left(B, \mathscr{P}\right) \to \left(\bar{B}, \bar{\mathscr{P}}\right)$. Here $\check{\bar{\mathfrak{X}}} \to \Spec \kk\llbracket t \rrbracket$ is the toric degeneration mirror constructed using the initial slab functions obtained from the walls of $\mathfrak{D}_{J}$ by basechange and trivial gluing data. In particular, check that this choice of slab functions defines a structure of a toric log CY on $\check{\bar{X}}_0$. This should be a straightforward generalization of the argument in Section \ref{main}. 
\end{enumerate}

To do (1), we need good local models of $\bar{\mathfrak{X}} \to \mathcal{S}$ near the points of the singular locus (for special toric degenerations of K3-s these are given by Observation \ref{modelslocal}(2)). Although condition (3) in Definition \ref{special} of a special toric degeneration guarantees that we have such local models (see Remark \ref{specialass}(3)), in general, these local models can be quite complicated (especially when the toric degeneration satisfies part (b) of condition (3)) and hard to control. Moreover, in general, it is challenging to globalize the resolution (using the arguments as in the proof of Proposition \ref{stronglogsmoothresolution}). 

By Proposition \ref{distinguishedspecial}, we can restrict to the natural case of a distinguished toric degeneration $\bar{\mathfrak{X}} \to \mathcal{S}$ with a simple dual intersection complex $\left(\bar B, \bar{\mathscr{P}}\right)$ and a smooth generic fibre. In this case, the local models for the singular points of $\bar{\mathfrak{X}} \to \mathcal{S}$ are controlled by $\left(\bar B, \bar{\mathscr{P}}\right)$. Another reason to restrict to this case is that simplicity of $\left(\bar B, \bar{\mathscr{P}}\right)$ guarantees a unique choice of the initial slab functions independent of the gluing data (see Proposition \ref{prop65}). Therefore, unlike the general case, the space of toric log CY structures on $\check{\bar{X}}_0$ is well-behaved (see \cite[Theorem 5.4]{GS2}). This should make it possible to generalize the results of Chapter \ref{gslocus} to this setting.

\begin{ass} \label{naturalcase}
Let $\bar{\mathfrak{X}} \to \mathcal{S}$ be a special toric degeneration. We assume that $\bar{\mathfrak{X}} \to \mathcal{S}$ is distinguished in the sense of Definition \ref{distinguished} and the dual intersection complex $\left(\bar B, \bar{\mathscr{P}}\right)$ of $\bar{\mathfrak{X}} \to \mathcal{S}$ is simple in the sense of \cite[Definition 1.60]{GS2}.\footnote{Note that by Proposition \ref{distinguishedspecial}, under these assumptions $\bar{\mathfrak{X}} \to \mathcal{S}$ is special if and only if the generic fibre of $\bar{\mathfrak{X}} \to \mathcal{S}$ is smooth.} 
\end{ass}

In this chapter, we construct (strongly) admissible resolutions of special toric degenerations of CY threefolds satisfying Assumption \ref{naturalcase} and sketch a generalization to relative dimension $n \geq 4$ (under an additional Assumption \ref{facelocalmodelass} and modulo a combinatorial Conjecture \ref{combconjecture}). We also discuss generalizing the scattering results of Section \ref{scatint} (see Conjectures \ref{canonicalscathd} and \ref{scatdiaghigherpowershd}) and reduce Conjecture \ref{conjecture} for a special toric degeneration $\bar{\mathfrak{X}} \to \mathcal{S}$ satisfying Assumption \ref{naturalcase} and Assumption \ref{facelocalmodelass} (if $n \geq 4$) to Conjectures \ref{combconjecture}, \ref{canonicalscathd}, and \ref{scatdiaghigherpowershd}. Finally, we conjecture (under the same assumptions) a generalization of Theorem \ref{compgen} (the main result of Chapter \ref{gslocus}) to higher dimensions (see Conjecture \ref{conjecturehd}).

\subsection{Resolutions in higher dimensions} \label{resolvehd}

As in the case of special toric degenerations of K3-s, to construct resolutions, we need to first understand how to resolve the local models. 
Let $\bar{\mathfrak{X}} \to \mathcal{S}$ be a special toric degeneration of relative dimension $n$ satisfying Assumption \ref{naturalcase}. Let $\left( \bar{B}, \bar{\mathscr{P}} \right)$ be the dual intersection complex of $\bar{\mathfrak{X}} \to \mathcal{S}$ and let $x \in Z$ be a singular point contained in the minimal stratum $\bar{X}_{\tau}, \ \tau \in \bar{\mathscr{P}}$ (note that since $x$ is singular, $1 \leq \dim \tau \leq n-1$). Then the étale local model $\tilde{\bar{X}}_{\tau, x}$ in the neighbourhood of $x$ is an affine toric variety defined by the cone over the convex hull $$\Delta_{\tau} := \Conv\left( \bigcup_{i=0}^q \left(\Delta_{\tau,i} \times \left\{e_i\right\}\right) \right)$$ as in \eqref{convexhull} by \cite[Theorem 2.6 and Corollary 2.18]{models}. Here $\Delta_{\tau,i}$ are integral polytopes in a lattice $N$ (not necessarily of maximal dimension), $e_i$ for $ \ 0 \leq i \leq q$ are the standard generators of the second factor of $N' := N \oplus \mathbb{Z}^{q+1}$ and the convex hull is taken in $N'_{\mathbb{R}}$. Moreover, \cite[Theorem 2.6 and Corollary 2.18]{models} imply that $\Delta_{\tau,0} := \tau$ is a cell of dimension $1 \leq \dim \tau \leq n-1$ (with its integral structure), $q := n - \dim \tau$, and $\Delta_{\tau,i}, \ 1 \leq i \leq q$ are the monodromy simplices associated to $\tau$ (see \cite[Definition 1.60]{GS3}) that are certain  Minkowski summands of $\tau$. The map $\pi: \tilde{\bar{X}}_{\tau, x} \to \mathbb{A}^1 \cong \Spec \kk[\mathbb{N}]$ is induced by the dual of $e_0$ and we equip $\tilde{\bar{X}}_{\tau, x}$ with the divisorial log structure with divisor $\pi^{-1}(0)$ as usual.

The condition that the generic fibre of $\bar{\mathfrak{X}} \to \mathcal{S}$ is smooth implies (see \cite[Proposition 2.2]{models}) that $$\Delta_{\tau,g}: = \Conv\left( \bigcup_{i=1}^q \left(\Delta_{\tau,i} \times \left\{e_i\right\}\right) \right)$$ is a standard simplex. Note that we have 
\begin{equation} \label{conv2}
\Delta_{\tau} \cong \Conv\left((\tau \times \{e_0\}) \cup  \Delta_{\tau,g})\right).
\end{equation}

In what follows, we will mostly consider local models with $\tau \in \bar{\mathscr{P}}^{[n-1]}$. We assume that the local models with $\tau \in \bar{\mathscr{P}}^{[k]}, \ 1 \leq k \leq n-2$ are determined by the local models with $\tau \in \bar{\mathscr{P}}^{[n-1]}$.

\begin{ass} \label{facelocalmodelass}
For cells $\tau', \tau \in \bar{\mathscr{P}}$ with $\tau' \subseteq \tau, \ 1 \leq \dim \tau' < \dim \tau \leq n-1$, we assume that the local model at a singular point $x'$ with minimal stratum $\bar{X}_{\tau'}$ is given by the cone over the unique face $F_{\tau'}$ of $\Delta_{\tau}$ such that $F_{\tau'} \cap (\tau \times \{e_0\}) = \tau'$ and $F_{\tau'} \cap \Delta_{\tau,g}  \neq \varnothing$ (using the description of \eqref{conv2}).
\end{ass}

\begin{remarks} \label{facelocalmodel}
\begin{enumerate}
    \item If $n = 3$, Assumption \ref{facelocalmodelass} is satisfied. Indeed, in this case, we have $\dim \tau = 2, \ \dim \tau' = 1$, and Assumption \ref{facelocalmodelass} follows from the construction in the proof of \cite[Theorem 2.6]{models} and the analysis of the local models with $\dim \tau' = 1$ of Sections \ref{locmod} and \ref{localmodelgen}.
    \item If $n \geq 4$, Assumption \ref{facelocalmodelass} is not always satisfied, and one can have more complicated local models with $\tau \in \bar{\mathscr{P}}^{[k]}, \ 1 \leq k \leq n-2$. This follows from the behaviour of the singular locus, see \cite[Example 1.62]{GS2}.
\end{enumerate}
\end{remarks}

We will now construct resolutions of a special toric degeneration $\bar{\mathfrak{X}} \to \mathcal{S}$ of CY threefolds satisfying Assumption \ref{naturalcase} and sketch the approach in higher dimensions. This work is related to \cite[Section 4.1.1]{kasa}, which constructed resolutions of the local models $\tilde{\bar{X}}_{\tau, x}$ in general (these resolutions do not patch to a global resolution of $\bar{\mathfrak{X}} \to \mathcal{S}$), and to \cite[Section 5.5]{ruddat} (in German), which constructed a partial resolution of a toric log CY space (i.e. of the central fibre $\bar{\mathfrak{X}}_0$ of a toric degeneration). 

\subsubsection{Resolutions in relative dimension \texorpdfstring{$n=3$}{n3} assuming all \texorpdfstring{$\sigma \in \bar{\mathscr{P}}^{\max}$}{sigmamax} are standard simplices} \label{genericallysncgeneral}

Suppose that $\bar{\mathfrak{X}} \to \mathcal{S}$ is a special toric degeneration of relative dimension $3$ satisfying Assumption \ref{naturalcase}. Let us first consider the case that all the maximal cells $\sigma \in \bar{\mathscr{P}}^{\max}$ are standard simplices (this corresponds to the case of Section \ref{resolvesmall} for special toric degenerations of K3-s). 

First, we need to understand how to resolve the local models. Let $x$ be a singular point contained in the minimal stratum $\bar{X}_{\tau}$ for $\tau \in \mathscr{P}^{[2]}$ a codimension $1$ cell. Then our assumption on the maximal cells implies that $\tau$ is a standard triangle. Moreover, we have $q = 1$ so $\Delta_{\tau,g} = \Delta_{\tau,1}$ is a Minkowski summand of $\tau$. However, the only possible Minkowski summands of a standard triangle are a point or the triangle itself. If $\Delta_{\tau,g}$ is a point, the corresponding local model (given by the cone over the $\Delta_{\tau}$ of \eqref{conv2}) is log smooth. So we need to consider the case that $\Delta_{\tau,g} \cong \tau$. Using Lemma \ref{blowuplemma}, one can check that this local model can be resolved by blowing up any two of the three divisors corresponding to the vertices of $\tau$. See Figure \ref{locmodeltriangle} for the corresponding subdivisions of the local model (the subdivisions are given in blue, and the vertices corresponding to the blown up divisors are highlighted in red). Note that the local models for the edges of $\tau$ are as in Figure \ref{coneoversquare} (see Remark \ref{facelocalmodel}(1)) and are subdivided accordingly by the blowups.

\begin{figure}[b]
\begin{center}

\begin{minipage}{.3\textwidth}
\begin{tikzpicture}[line cap=round,line join=round,>=triangle 45,x=1cm,y=1cm, scale = 3]
\clip(-0.2,-0.2) rectangle (1.5,1.5);
\draw (0,1)-- (1,1);
\draw (1,1)-- (1,0);
\draw (0,0)-- (1,0);
\draw (0,0)-- (0,1);
\draw (0,1)-- (0.66,1.33);
\draw (1,1)-- (0.66,1.33);
\draw[dashed] (0,0)-- (0.66,0.33);
\draw[dashed] (1,0)-- (0.66,0.33);
\draw[dashed] (0.66,0.33)-- (0.66,1.33);
\draw [fill] (0,0) circle (0.5pt);
\draw [fill] (0.66,0.33) circle (0.5pt);
\draw [fill] (0.66,1.33) circle (0.5pt);
\draw [fill] (0,1) circle (0.5pt);
\draw [fill] (1,0) circle (0.5pt);
\draw [fill] (1,1) circle (0.5pt);

\draw(0.5, 0.1) node{\scriptsize $\tau$};
\draw(0.5, 1.1) node{\scriptsize $\Delta_{\tau,g}$};

\draw(-0.06, 0.1) node{\scriptsize $v_2$};
\draw(1.07, 0.1) node{\scriptsize $v_3$};
\draw(0.75, 0.33) node{\scriptsize $v_1$};

\end{tikzpicture}
\captionsetup{labelformat=empty}
\caption{Local model}
\end{minipage}
\begin{minipage}{.3\textwidth}
\begin{tikzpicture}[line cap=round,line join=round,>=triangle 45,x=1cm,y=1cm, scale = 3]
\clip(-0.2,-0.2) rectangle (1.5,1.5);
\draw (0,1)-- (1,1);
\draw (1,1)-- (1,0);
\draw (0,0)-- (1,0);
\draw (0,0)-- (0,1);
\draw (0,1)-- (0.66,1.33);
\draw (1,1)-- (0.66,1.33);
\draw[dashed] (0,0)-- (0.66,0.33);
\draw[dashed] (1,0)-- (0.66,0.33);
\draw[dashed] (0.66,0.33)-- (0.66,1.33);
\draw [fill] (0,0) circle (0.5pt);
\draw [fill] (0.66,1.33) circle (0.5pt);
\draw [fill, color = red] (0.66,0.33) circle (1pt);
\draw [fill] (0,1) circle (0.5pt);
\draw [fill] (1,0) circle (0.5pt);
\draw [fill] (1,1) circle (0.5pt);

\filldraw[draw=black, fill=blue, opacity=.2] (0.66,0.33) -- (0,1) -- (1,1) -- cycle;

\draw(0.5, 0.1) node{\scriptsize $\tau$};
\draw(0.5, 1.1) node{\scriptsize $\Delta_{\tau,g}$};

\draw(-0.06, 0.1) node{\scriptsize $v_2$};
\draw(1.07, 0.1) node{\scriptsize $v_3$};
\draw(0.75, 0.33) node{\scriptsize $v_1$};

\end{tikzpicture}
\captionsetup{labelformat=empty}
\caption{After one blowup}
\end{minipage}
\begin{minipage}{.3\textwidth}
\begin{tikzpicture}[line cap=round,line join=round,>=triangle 45,x=1cm,y=1cm, scale = 3]
\clip(-0.2,-0.2) rectangle (1.5,1.5);
\draw (0,1)-- (1,1);
\draw (1,1)-- (1,0);
\draw (0,0)-- (1,0);
\draw (0,0)-- (0,1);
\draw (0,1)-- (0.66,1.33);
\draw (1,1)-- (0.66,1.33);
\draw[dashed] (0,0)-- (0.66,0.33);
\draw[dashed] (1,0)-- (0.66,0.33);
\draw[dashed] (0.66,0.33)-- (0.66,1.33);
\draw [fill, color=red] (0,0) circle (1pt);
\draw [fill] (0.66,1.33) circle (0.5pt);
\draw [fill] (0.66,0.33) circle (0.5pt);
\draw [fill] (0,1) circle (0.5pt);
\draw [fill] (1,0) circle (0.5pt);
\draw [fill] (1,1) circle (0.5pt);

\filldraw[draw=black, fill=blue, opacity=.2] (0.66,0.33) -- (0,1) -- (1,1) -- cycle;

\filldraw[draw=black, fill=blue, opacity=.2] (0,0) -- (0.66,0.33) -- (1,1) -- cycle;

\draw(0.5, 0.1) node{\scriptsize $\tau$};
\draw(0.5, 1.1) node{\scriptsize $\Delta_{\tau,g}$};

\draw(-0.06, 0.1) node{\scriptsize $v_2$};
\draw(1.07, 0.1) node{\scriptsize $v_3$};
\draw(0.75, 0.33) node{\scriptsize $v_1$};

\end{tikzpicture}
\captionsetup{labelformat=empty}
\caption{After two blowups}
\end{minipage}
\end{center}
\setcounter{figure}{0}
\caption{Resolution of the local model $\tilde{\bar{X}}_{\tau,x}$ for $\tau$ a standard triangle.}
\label{locmodeltriangle}
\end{figure}

The singular locus of the local model $\tilde{\bar{X}}_{\tau,x}$ of Figure \ref{locmodeltriangle} is a union of three lines $\ell_{1,2} \subseteq D_{v_1} \cup D_{v_2}, \ell_{13} \subseteq D_{v_1} \cup D_{v_2}, \ell_{23} \subseteq D_{v_2} \cup D_{v_3}$ intersecting at a point. The exceptional locus of the blowup is the union of a $\mathbb{P}^1$-bundle over the strict transform of $\ell_{13}$ with exceptional curve $E_{v_1}$, a $\mathbb{P}^1$-bundle over the strict transform of $\ell_{23}$ with exceptional curve $E_{v_2}$, and a $\mathbb{P}^1$-bundle over the strict transform of $\ell_{12}$ with exceptional curve $E_{v_1}-E_{v_2}$. See Figure \ref{locmodeltriangleresolve} for a sketch of the central fibre of $\tilde{\bar{X}}_{\tau,x}$ and the central fibres of the blowups of Figure \ref{locmodeltriangle}.

\begin{figure}[H]
\begin{center}

\begin{minipage}{.3\textwidth}
\begin{tikzpicture}[line cap=round,line join=round,>=triangle 45,x=1cm,y=1cm, scale = 3]
\clip(-0.2,-0.2) rectangle (1.25,1.25);
\draw  (0.5,0.1)-- (0.5,0.9);

\draw  (0.5,0.1)-- (0.9,0.1);
\draw  (0.5,0.1)-- (0.1,0.3);
\draw  (0.5,0.1)-- (0.3,-0.1);

\draw  (0.5,0.9)-- (0.9,0.9);
\draw  (0.5,0.9)-- (0.1,1.1);
\draw  (0.5,0.9)-- (0.3,0.7);

\draw [color=blue, line width=1pt] (0.5,0.5)-- (0.9,0.5);
\draw [color=blue, line width=1pt] (0.5,0.5)-- (0.1,0.7);
\draw [color=blue, line width=1pt] (0.5,0.5)-- (0.3,0.3);

\draw(0.4, 0.3) node{\scriptsize $\ell_{12}$};
\draw(0.2, 0.75) node{\scriptsize $\ell_{13}$};
\draw(0.85, 0.57) node{\scriptsize $\ell_{23}$};

\draw(0.2, 0) node{\scriptsize $D_{v_1}$};
\draw(0.7, 0) node{\scriptsize $D_{v_2}$};
\draw(0.7, 0.2) node{\scriptsize $D_{v_3}$};

\end{tikzpicture}
\captionsetup{labelformat=empty}
\caption{Local model}
\end{minipage}
\begin{minipage}{.3\textwidth}
\begin{tikzpicture}[line cap=round,line join=round,>=triangle 45,x=1cm,y=1cm, scale = 3]
\clip(-0.2,-0.2) rectangle (1.25,1.25);
\draw  (0.5,0.1)-- (0.5,0.9);

\draw  (0.5,0.1)-- (0.9,0.1);
\draw  (0.5,0.1)-- (0.1,0.3);
\draw  (0.5,0.1)-- (0.3,-0.1);

\draw  (0.5,0.9)-- (0.9,0.9);
\draw  (0.5,0.9)-- (0.1,1.1);
\draw  (0.5,0.9)-- (0.3,0.7);

\draw [color=blue, line width=1pt] (0.5,0.5)-- (0.9,0.5);
\draw  (0.5,0.5)-- (0.1,0.7);
\draw  (0.5,0.5)-- (0.3,0.3);

\draw(0.4, 0.3) node{\scriptsize $\ell_{12}$};
\draw(0.2, 0.75) node{\scriptsize $\ell_{13}$};
\draw(0.85, 0.57) node{\scriptsize $\ell_{23}$};

\draw(0.2, 0) node{\scriptsize $D_{v_1}$};
\draw(0.7, 0) node{\scriptsize $D_{v_2}$};
\draw(0.7, 0.2) node{\scriptsize $D_{v_3}$};

\draw [color=red, line width=1pt] (0.5,0.5)-- ++(-4pt,0pt);
\draw [color=red, line width=1pt] (0.4,0.55)-- ++(-4pt,0pt);
\draw [color=red, line width=1pt] (0.3,0.6)-- ++(-4pt,0pt);
\draw [color=red, line width=1pt] (0.2,0.65)-- ++(-4pt,0pt);
\draw [color=red, line width=1pt] (0.1,0.7)-- ++(-4pt,0pt);
\draw [color=red, line width=1pt] (0.45,0.45)-- ++(-4pt,0pt);
\draw [color=red, line width=1pt] (0.4,0.4)-- ++(-4pt,0pt);
\draw [color=red, line width=1pt] (0.35,0.35)-- ++(-4pt,0pt);
\draw [color=red, line width=1pt] (0.3,0.3)-- ++(-4pt,0pt);

\draw(0, 0.61) node{\tiny $E_{v_1}$};

\end{tikzpicture}
\captionsetup{labelformat=empty}
\caption{After one blowup}
\end{minipage}
\begin{minipage}{.3\textwidth}
\begin{tikzpicture}[line cap=round,line join=round,>=triangle 45,x=1cm,y=1cm, scale = 3]
\clip(-0.3,-0.2) rectangle (1.25,1.25);
\draw  (0.5,0.1)-- (0.5,0.9);

\draw  (0.5,0.1)-- (0.9,0.1);
\draw  (0.5,0.1)-- (0.1,0.3);
\draw  (0.5,0.1)-- (0.3,-0.1);

\draw  (0.5,0.9)-- (0.9,0.9);
\draw  (0.5,0.9)-- (0.1,1.1);
\draw  (0.5,0.9)-- (0.3,0.7);

\draw  (0.5,0.5)-- (0.9,0.5);
\draw  (0.5,0.5)-- (0.1,0.7);
\draw  (0.5,0.5)-- (0.3,0.3);

\draw(0.4, 0.3) node{\scriptsize $\ell_{12}$};
\draw(0.2, 0.75) node{\scriptsize $\ell_{13}$};
\draw(0.85, 0.57) node{\scriptsize $\ell_{23}$};

\draw(0.2, 0) node{\scriptsize $D_{v_1}$};
\draw(0.7, 0) node{\scriptsize $D_{v_2}$};
\draw(0.7, 0.2) node{\scriptsize $D_{v_3}$};

\draw [color=red, line width=1pt] (0.5,0.5)-- ++(-4pt,0pt);
\draw [color=red, line width=1pt] (0.4,0.55)-- ++(-4pt,0pt);
\draw [color=red, line width=1pt] (0.3,0.6)-- ++(-4pt,0pt);
\draw [color=red, line width=1pt] (0.2,0.65)-- ++(-4pt,0pt);
\draw [color=red, line width=1pt] (0.1,0.7)-- ++(-4pt,0pt);
\draw [color=red, line width=1pt] (0.45,0.45)-- ++(-3pt,1pt);
\draw [color=red, line width=1pt] (0.4,0.4)-- ++(-3pt,1pt);
\draw [color=red, line width=1pt] (0.35,0.35)-- ++(-3pt,1pt);
\draw [color=red, line width=1pt] (0.3,0.3)-- ++(-3pt,1pt);

\draw(0, 0.61) node{\tiny $E_{v_1}$};

\draw(-0.05, 0.35) node{\tiny $E_{v_1} - E_{v_2}$};

\draw [color=red, line width=1pt] (0.5,0.5)-- ++(-1.5pt,-1.5pt);
\draw [color=red, line width=1pt] (0.6,0.5)-- ++(-1.5pt,-1.5pt);
\draw [color=red, line width=1pt] (0.7,0.5)-- ++(-1.5pt,-1.5pt);
\draw [color=red, line width=1pt] (0.8,0.5)-- ++(-1.5pt,-1.5pt);
\draw [color=red, line width=1pt] (0.9,0.5)-- ++(-1.5pt,-1.5pt);

\draw(0.95, 0.4) node{\tiny $E_{v_2}$};

\end{tikzpicture}
\captionsetup{labelformat=empty}
\caption{After two blowups}
\end{minipage}
\end{center}
\setcounter{figure}{1}
\caption{Blowing up the local model $\tilde{\bar{X}}_{\tau,x}$.} \label{locmodeltriangleresolve}
\end{figure}

Similarly to the analysis of Section \ref{resolvesmall}, blowing up a component $\bar{D}_{v}, \ v \in \bar{\mathscr{P}}^{[0]}$ of the central fibre $\bar{\mathfrak{X}}_0$ of $\bar{\mathfrak{X}}$ gives a small partial resolution $\mathfrak{X}' \to \bar{\mathfrak{X}}$ that is trivial in the local models for the codimension $3$ strata $\bar{X}_{\sigma}, \ \sigma \in \bar{\mathscr{P}}^{\max}$ (given by the cones over $\sigma \in \bar{\mathscr{P}}^{\max}$). In codimension $2$, the blowup is described via the local models for $\tau \in \bar{\mathscr{P}}^{[2]}$ with $v \subseteq \tau$ of Figure \ref{locmodeltriangleresolve} (after one blowup) and the behaviour in codimension $1$ follows from Remark \ref{facelocalmodel}(1). Note that the singular locus of $\bar{D}_v$ is a connected union of curves $C_{\rho} \subseteq \bar{X}_{\rho}$ for $\rho \in \bar{\mathscr{P}}^{[1]}$ with $v \subseteq \rho$. Blowing up $\bar{D}_v$ introduces a $\mathbb{P}^1$-bundle over the strict transform of $\bigcup C_{\rho}$ with exceptional curve $E_{v}$. This globalizes the first blowup of Figure \ref{locmodeltriangleresolve}.

As usual, we define the dual intersection complex $\left(B', \mathscr{P}'\right)$ of $\mathfrak{X}' \overset{g'}{\to} \mathcal{S}$ as $(g')^{-1}_{\trop}(1)$ for $g'_{\trop}: \Sigma(\mathfrak{X}') \to \mathbb{R}_{\geq 0}$ the tropicalization of $g'$ of Construction \ref{tropicalize}. The fact that all the maximal cells $\sigma \in \bar{\mathscr{P}}^{\max}$ are standard simplices implies that $\left(B', \mathscr{P}'\right) \cong \left(\bar{B}, \bar{\mathscr{P}}\right)$ as polyhedral manifolds. Further, $\left(B', \mathscr{P}'\right)$ possesses the structure of an affine manifold with singularities via Construction \ref{affstrint}. 

\begin{sloppypar}
A similar analysis applies for further blowups of irreducible components $\bar{D}_{v'}, \ v' \in \bar{\mathscr{P}}^{[0]}$ of $\bar{\mathfrak{X}}_0$. In the case that a component corresponding to a vertex of $\tau \in \bar{\mathscr{P}}^{[2]}$ (with $v' \subseteq \tau$) has already been blown up, we use the local model of Figure \ref{locmodeltriangleresolve} after two blowups to describe the behaviour in codimension $2$. In this case, the singular locus of $\bar{D}_{v'}$ is smaller than in the analysis above. The blowup still introduces a $\mathbb{P}^1$-bundle over the singular locus.
\end{sloppypar}

In Figure \ref{singularlocusltriangle}, we give the transformation of the singular locus (in blue) of a cell $\tau = \langle v_1,v_2,v_3 \rangle \in \bar{\mathscr{P}}^{[2]}$ corresponding to first blowing up $\bar{D}_{v_1}$ and then blowing up $\bar{D}_{v_2}$. Note that, similarly to the analysis of Section \ref{resolvesmall}, at the level of the affine manifolds, the blowup of a component $\bar{D}_v, \ v \in \bar{\mathscr{P}}^{[0]}$ can be visualized as pulling the singular loci of the cells of $\mathscr{P}$ adjacent to $v$ into $v$.

\begin{figure}[H]
\begin{center}

\begin{minipage}{.3\textwidth}
\begin{tikzpicture}[line cap=round,line join=round,>=triangle 45,x=1cm,y=1cm, scale = 3]
\clip(-0.2,-0.2) rectangle (1.2,1.2);
\draw (0,1)-- (1,0);
\draw (0,0)-- (0,1);
\draw (0,0)-- (1,0);

\draw [fill] (0,0) circle (0.5pt);
\draw [fill] (0,1) circle (0.5pt);
\draw [fill] (1,0) circle (0.5pt);

\draw(-0.07, 0.07) node{$v_2$};
\draw(1.06, 0.07) node{$v_3$};
\draw(-0.06, 1.07) node{$v_1$};

\draw [color=blue, line width=1pt] (0.33,0.33)-- (0.33,0);
\draw [color=blue, line width=1pt] (0.33,0.33)-- (0,0.33);
\draw [color=blue, line width=1pt] (0.33,0.33)-- (0.5,0.5);

\end{tikzpicture}
\captionsetup{labelformat=empty}
\caption{Codim one cell $\tau$}
\end{minipage}
\begin{minipage}{.3\textwidth}
\begin{tikzpicture}[line cap=round,line join=round,>=triangle 45,x=1cm,y=1cm, scale = 3]
\clip(-0.2,-0.2) rectangle (1.2,1.2);
\draw (0,1)-- (1,0);
\draw (0,0)-- (0,1);
\draw (0,0)-- (1,0);

\draw [fill] (0,0) circle (0.5pt);
\draw [fill] (0,1) circle (0.5pt);
\draw [fill] (1,0) circle (0.5pt);

\draw(-0.07, 0.07) node{$v_2$};
\draw(1.06, 0.07) node{$v_3$};
\draw(-0.06, 1.07) node{$v_1$};

\draw [color=blue, line width=1pt] (0,1)-- (0.33,0);

\end{tikzpicture}
\captionsetup{labelformat=empty}
\caption{Blowing up $\bar{D}_{v_1}$}
\end{minipage}
\begin{minipage}{.3\textwidth}
\begin{tikzpicture}[line cap=round,line join=round,>=triangle 45,x=1cm,y=1cm, scale = 3]
\clip(-0.2,-0.2) rectangle (1.2,1.2);
\draw (0,1)-- (1,0);
\draw (0,0)-- (0,1);
\draw (0,0)-- (1,0);

\draw [fill] (0,0) circle (0.5pt);
\draw [fill] (0,1) circle (0.5pt);
\draw [fill] (1,0) circle (0.5pt);

\draw(-0.07, 0.07) node{$v_2$};
\draw(1.06, 0.07) node{$v_3$};
\draw(-0.06, 1.07) node{$v_1$};

\draw [color=blue, line width=1pt] (0,1)-- (0,0);

\end{tikzpicture}
\captionsetup{labelformat=empty}
\caption{Blowing up $\bar{D}_{v_2}$}
\end{minipage}
\end{center}
\setcounter{figure}{2}
\caption{Transformation of a cell $\tau \in \bar{\mathscr{P}}^{[2]}$ under blowing up components.}
\label{singularlocusltriangle}
\end{figure}

Blow up all the irreducible components of the central fibre $\bar{\mathfrak{X}}_0$ of $\bar{\mathfrak{X}} \to \mathcal{S}$ in any order. By the analysis above, this gives a resolution $\pi: \mathfrak{X} \to \bar{\mathfrak{X}}$ to a log smooth degeneration $\mathfrak{X} \to \mathcal{S}$. Here the log structure on $\mathfrak{X}$ is the divisorial log structure given by $D = D_1 + \dots + D_m$ where $D_i, \ 1 \leq i \leq m$ are the irreducible components of the central fibre $\mathfrak{X}_0$. The exceptional locus of $\pi: \mathfrak{X} \to \bar{\mathfrak{X}}$ is the union of $\mathbb{P}^1$-bundles with exceptional curves  
\begin{equation} \label{contcurvesgeneral}
\left\{ E_{v}  \ | \ v \in \bar{\mathscr{P}}^{[0]}\right\} \cup \left\{ E_{v} - E_{v'} \ | \ \rho = \langle v, v' \rangle \in \bar{\mathscr{P}}^{[1]}\right\}
\end{equation}
(it is possible that all the singularities of some $\bar{D}_v$ are resolved before we blow up $\bar{D}_v$ in which case there is no $E_v$). Similarly to Section \ref{resolvesmall}, $\pi$ is a small morphism.
The singularities of the affine structure on the dual intersection complex $(B, \mathscr{P})$ of $\mathfrak{X} \to \mathcal{S}$ are contained in cells of codimension $2$, and one may view the resolution tropically as a composition of transformations as in Figure \ref{singularlocusltriangle}. Note that we have $(B, \mathscr{P}) \cong \left(\bar B, \bar{\mathscr{P}}\right)$ as polyhedral manifolds. 

Propositions \ref{sncmlcy} and \ref{ampledivisorsmall} of Section \ref{globalresolvesmall} admit a direct generalization.

\begin{fact}
\begin{enumerate}
    \item $\mathfrak{X} \to \mathcal{S}$ is minimal log CY and $D$ is simple normal crossings. The dual intersection complex $(B, \mathscr{P})$ of $\mathfrak{X} \to \mathcal{S}$ satisfies Assumption \ref{manifoldass}.
    \item The resolution $\pi: \mathfrak{X} \to \bar{\mathfrak{X}}$ satisfies the assumptions of Proposition \ref{P2}.    
\end{enumerate}
\end{fact}

\begin{proof}
(a) is similar to Proposition \ref{sncmlcy} using the local models of Figure \ref{locmodeltriangle} to check that $D$ is simple normal crossings in codimension $\geq 2$. (b) is similar to Proposition \ref{ampledivisorsmall} using the fact that for any $v \in \bar{\mathscr{P}}^{[0]}$ we have $E_v \cdot D_v = -1$, $E_v \cdot D_{v'} = 1$ for any  $v' \in \bar{\mathscr{P}}^{[0]}$ that is connected to $v$ by an edge $\rho \in \bar{\mathscr{P}}^{[1]}$, and $E_v \cdot D_{v'} = 0$ for all the other $v' \in \bar{\mathscr{P}}^{[0]}$.
\end{proof}

\subsubsection{Resolutions of the local models in the easily decomposable case} \label{resolveed}

One can generalize the analysis of Section \ref{resolve} to special toric degenerations of relative dimension $3$ satisfying Assumption \ref{naturalcase}. To do that one needs to restrict the types of $\sigma \in \bar{\mathscr{P}}^{\max}$ to natural generalizations of the types (1-4) in Figure \ref{admisspolytopes} (one might be able to take a larger class) and analyze the corresponding local models for $\tau \in \bar{\mathscr{P}}^{[2]}$. We do not provide further details and instead consider the general case, arguing similarly to Section \ref{admissible}. 

First, we need to understand how to resolve the local models $\tilde{\bar{X}}_{\tau,x}$ in general. Let $x$ be a singular point contained in the minimal stratum $\bar{X}_{\tau}$ for $\tau \in \mathscr{P}^{[2]}$ a codimension $1$ cell. Since $\left( \bar{B}, \bar{\mathscr{P}} \right)$ is simple by Assumption \ref{naturalcase}, we can use the classification of the singular loci of \cite[Example 1.62]{GS2} (here one should understand the singular locus of $\left( \bar{B}, \bar{\mathscr{P}} \right)$ to mean the smallest subset $\Delta$ such that the affine structure on $\left( \bar{B}, \bar{\mathscr{P}} \right)$ extends to $B \setminus \Delta$). This classification, along with the fact that $\Delta_{\tau,g} = \Delta_{\tau,1}$ is a standard triangle and a Minkowski summand of $\tau$, implies that there are two cases. The first case is that $\tau$ has two parallel edges $\rho_1, \rho_2 \in \bar{\mathscr{P}}^{[1]}$ and the singular locus of $\tau$ is as in Figure \ref{singlocitypes}(1).  The second case is that $\tau$ has three edges $\rho_1, \rho_2, \rho_3 \in \bar{\mathscr{P}}^{[1]}$ that are parallel to the corresponding edges of some standard triangle and the singular locus of $\tau$ is as in Figure \ref{singlocitypes}(2). Recall that the precise location of the singular locus does not matter as long as it respects the cell structure.

In the case of Figure \ref{singlocitypes}(1), $\Delta_{\tau,g} = \Delta_{\tau,1}$ is the line segment connecting the points $(0,0)$ and $(0,1)$ (i.e. a standard simplex of dimension $1$). In the case of Figure \ref{singlocitypes}(2), $\Delta_{\tau,g} = \Delta_{\tau,1}$ is a standard triangle. This describes the local models $\tilde{\bar{X}}_{\tau,x}$ via \eqref{conv2}. 

\newpage

\begin{figure}[H]
\begin{center}
\begin{minipage}{.4\textwidth}
\begin{tikzpicture}[scale=0.7]

\draw[fill] (0,0) circle (1.5pt) coordinate (m-0-0);

\draw[fill] (0,1) circle (1.5pt) coordinate (m-0-1);
\draw[fill] (1,1) circle (1.5pt) coordinate (m-1-1);
\draw[fill] (2,1) circle (1.5pt) coordinate (m-2-1);
\draw[fill] (3,1) circle (1.5pt) coordinate (m-3-1);
\draw[fill] (4,1) circle (1.5pt) coordinate (m-4-1);
\draw[fill] (5,1) circle (1.5pt) coordinate (m-5-1);
\draw[fill] (6,1) circle (1.5pt) coordinate (m-6-1);

\draw[fill] (0,2) circle (1.5pt) coordinate (m-0-2);
\draw[fill] (1,2) circle (1.5pt) coordinate (m-1-2);
\draw[fill] (2,2) circle (1.5pt) coordinate (m-2-2);
\draw[fill] (3,2) circle (1.5pt) coordinate (m-3-2);
\draw[fill] (4,2) circle (1.5pt) coordinate (m-4-2);
\draw[fill] (5,2) circle (1.5pt) coordinate (m-5-2);
\draw[fill] (6,2) circle (1.5pt) coordinate (m-6-2);

\draw[fill] (0,3) circle (1.5pt) coordinate (m-0-3);
\draw[fill] (1,3) circle (1.5pt) coordinate (m-1-3);
\draw[fill] (2,3) circle (1.5pt) coordinate (m-2-3);
\draw[fill] (3,3) circle (1.5pt) coordinate (m-3-3);
\draw[fill] (4,3) circle (1.5pt) coordinate (m-4-3);
\draw[fill] (5,3) circle (1.5pt) coordinate (m-5-3);
\draw[fill] (6,3) circle (1.5pt) coordinate (m-6-3);

\draw[fill] (0,4) circle (1.5pt) coordinate (m-0-4);
\draw[fill] (1,4) circle (1.5pt) coordinate (m-1-4);
\draw[fill] (2,4) circle (1.5pt) coordinate (m-2-4);
\draw[fill] (3,4) circle (1.5pt) coordinate (m-3-4);
\draw[fill] (4,4) circle (1.5pt) coordinate (m-4-4);
\draw[fill] (5,4) circle (1.5pt) coordinate (m-5-4);
\draw[fill] (6,4) circle (1.5pt) coordinate (m-6-4);

\draw[fill] (1,5) circle (1.5pt) coordinate (m-1-5);
\draw[fill] (2,5) circle (1.5pt) coordinate (m-2-5);
\draw[fill] (3,5) circle (1.5pt) coordinate (m-3-5);
\draw[fill] (4,5) circle (1.5pt) coordinate (m-4-5);
\draw[fill] (5,5) circle (1.5pt) coordinate (m-5-5);
\draw[fill] (6,5) circle (1.5pt) coordinate (m-6-5);

\draw[fill] (5,6) circle (1.5pt) coordinate (m-5-6);

\draw (m-0-0) -- (m-6-1) -- (m-6-5) -- (m-5-6) -- (m-1-5) -- (m-0-4) -- cycle;

\draw(3, 6) node{\large $\tau$};

\draw(-0.5, 2) node{\large $\rho_1$};
\draw(6.5, 3) node{\large $\rho_2$};

\draw [color=blue, line width=1pt] (0,2.5)-- (6,2.5);

\end{tikzpicture}
\captionsetup{labelformat=empty}
\caption{(1)}
\end{minipage}
\hspace{1cm}
\begin{minipage}{.4\textwidth}
\begin{tikzpicture}[scale=0.7]

\draw[fill] (1,0) circle (1.5pt) coordinate (m-1-0);
\draw[fill] (2,0) circle (1.5pt) coordinate (m-2-0);
\draw[fill] (3,0) circle (1.5pt) coordinate (m-3-0);
\draw[fill] (4,0) circle (1.5pt) coordinate (m-4-0);

\draw[fill] (1,1) circle (1.5pt) coordinate (m-1-1);
\draw[fill] (2,1) circle (1.5pt) coordinate (m-2-1);
\draw[fill] (3,1) circle (1.5pt) coordinate (m-3-1);
\draw[fill] (4,1) circle (1.5pt) coordinate (m-4-1);
\draw[fill] (5,1) circle (1.5pt) coordinate (m-5-1);
\draw[fill] (6,1) circle (1.5pt) coordinate (m-6-1);

\draw[fill] (0,2) circle (1.5pt) coordinate (m-0-2);
\draw[fill] (1,2) circle (1.5pt) coordinate (m-1-2);
\draw[fill] (2,2) circle (1.5pt) coordinate (m-2-2);
\draw[fill] (3,2) circle (1.5pt) coordinate (m-3-2);
\draw[fill] (4,2) circle (1.5pt) coordinate (m-4-2);
\draw[fill] (5,2) circle (1.5pt) coordinate (m-5-2);
\draw[fill] (6,2) circle (1.5pt) coordinate (m-6-2);

\draw[fill] (0,3) circle (1.5pt) coordinate (m-0-3);
\draw[fill] (1,3) circle (1.5pt) coordinate (m-1-3);
\draw[fill] (2,3) circle (1.5pt) coordinate (m-2-3);
\draw[fill] (3,3) circle (1.5pt) coordinate (m-3-3);
\draw[fill] (4,3) circle (1.5pt) coordinate (m-4-3);
\draw[fill] (5,3) circle (1.5pt) coordinate (m-5-3);

\draw[fill] (0,4) circle (1.5pt) coordinate (m-0-4);
\draw[fill] (1,4) circle (1.5pt) coordinate (m-1-4);
\draw[fill] (2,4) circle (1.5pt) coordinate (m-2-4);
\draw[fill] (3,4) circle (1.5pt) coordinate (m-3-4);
\draw[fill] (4,4) circle (1.5pt) coordinate (m-4-4);

\draw[fill] (0,5) circle (1.5pt) coordinate (m-0-5);
\draw[fill] (1,5) circle (1.5pt) coordinate (m-1-5);
\draw[fill] (2,5) circle (1.5pt) coordinate (m-2-5);
\draw[fill] (3,5) circle (1.5pt) coordinate (m-3-5);

\draw[fill] (2,6) circle (1.5pt) coordinate (m-2-6);

\draw (m-1-0) -- (m-4-0) -- (m-6-1) -- (m-6-2) -- (m-2-6) -- (m-0-5) -- (m-0-2) -- cycle;

\draw(3, 6) node{\large $\tau$};

\draw(-0.5, 3) node{\large $\rho_1$};
\draw(2, -0.5) node{\large $\rho_3$};
\draw(4.5, 4.5) node{\large $\rho_2$};

\draw [color=blue, line width=1pt] (0,3.5)-- (2.5,3.5);
\draw [color=blue, line width=1pt] (2.5,0)-- (2.5,3.5);
\draw [color=blue, line width=1pt] (4.5,3.5)-- (2.5,3.5);

\end{tikzpicture}
\captionsetup{labelformat=empty}
\caption{(2)}
\end{minipage}
\end{center}
\setcounter{figure}{3}
\caption{Types of singular loci of $\tau \in \bar{\mathscr{P}}^{[2]}$.} 
\label{singlocitypes}
\end{figure}

We will now construct a subdivision of $\tau$ (or some rescaling of $\tau$) and an induced subdivision of $\Delta_{\tau}$ that resolves the local model $\tilde{\bar{X}}_{\tau,x}$. We consider a special case first.

\begin{defn} 
Let $\tau \in \bar{\mathscr{P}}^{[2]}$ be a cell of codimension $1$ (contained in a lattice $N \cong \mathbb{Z}^2$). For any edge $\rho \in \bar{\mathscr{P}}^{[1]}$ with $\rho \subseteq \tau$ and a point $y \in N_{\mathbb{R}}$, denote by $n_{\rho,y}$ any line such that:
\begin{enumerate}
    \item $n_{\rho,y}$ passes through $y$ and intersects the line $l_{\rho}$ generated by $\rho$.
    \item The primitive generator of $n_{\rho,y}$ and the primitive generator of $\rho$ generate the lattice $N$. 
\end{enumerate}
We call such a $n_{\rho,y}$ an \emph{integral normal}.
\end{defn}

\begin{defn} \label{easydecompose}
We say that $\tau$ is \emph{easily decomposable} if there exists a point $y \in \Int(\tau)$ with irrational coordinates and a choice of integral normals such that:
\begin{enumerate}
    \item If $\tau$ is as in Figure \ref{singlocitypes}(1), $n_{\rho_1, y} \cap l_{\rho_1} \in \Int(\rho_1)$ and $n_{\rho_2, y} \cap l_{\rho_2} \in \Int(\rho_2)$.
    \item If $\tau$ is as in Figure \ref{singlocitypes}(2), $n_{\rho_1, y} \cap l_{\rho_1} \in \Int(\rho_1), \ n_{\rho_2, y} \cap l_{\rho_2} \in \Int(\rho_2)$, and $n_{\rho_3, y} \cap l_{\rho_3} \in \Int(\rho_3)$.
\end{enumerate}
\end{defn}

Suppose that $\tau \in \bar{\mathscr{P}}^{[2]}$ is easily decomposable. We first do a partial subdivision as in Figure \ref{subdivisiontypes} for the types of $\tau$ in Figure \ref{singlocitypes}. Such subdivisions always exist for suitable rescalings $k\tau$ of $\tau$ with $k \in \mathbb{Z}_{>0}$ large enough (note that we have $\Delta_{k\tau,1} \cong \Delta_{\tau,1}$ under the rescalings). We denote by a blue point the point $y$ of Definition \ref{easydecompose} and draw the integral normals as dashed lines. We require that every edge of $\tau$ subdivided by our subdivision is a union of three (not two) edges of the subdivision. Again, this can clearly be achieved by a suitable choice of $y$ in Definition \ref{easydecompose} and a large enough rescaling factor $k$. The idea behind this subdivision is to create a tubular neighbourhood of the singular locus of $\tau$ composed of standard squares (and a standard triangle in the case of Figure \ref{subdivisiontypes}(2)). 

\begin{figure}[H]
\begin{center}
\begin{minipage}{.4\textwidth}
\begin{tikzpicture}[scale=0.7]

\draw[fill] (0,0) circle (1.5pt) coordinate (m-0-0);

\draw[fill] (0,1) circle (1.5pt) coordinate (m-0-1);
\draw[fill] (1,1) circle (1.5pt) coordinate (m-1-1);
\draw[fill] (2,1) circle (1.5pt) coordinate (m-2-1);
\draw[fill] (3,1) circle (1.5pt) coordinate (m-3-1);
\draw[fill] (4,1) circle (1.5pt) coordinate (m-4-1);
\draw[fill] (5,1) circle (1.5pt) coordinate (m-5-1);
\draw[fill] (6,1) circle (1.5pt) coordinate (m-6-1);

\draw[fill] (0,2) circle (1.5pt) coordinate (m-0-2);
\draw[fill] (1,2) circle (1.5pt) coordinate (m-1-2);
\draw[fill] (2,2) circle (1.5pt) coordinate (m-2-2);
\draw[fill] (3,2) circle (1.5pt) coordinate (m-3-2);
\draw[fill] (4,2) circle (1.5pt) coordinate (m-4-2);
\draw[fill] (5,2) circle (1.5pt) coordinate (m-5-2);
\draw[fill] (6,2) circle (1.5pt) coordinate (m-6-2);

\draw[fill] (0,3) circle (1.5pt) coordinate (m-0-3);
\draw[fill] (1,3) circle (1.5pt) coordinate (m-1-3);
\draw[fill] (2,3) circle (1.5pt) coordinate (m-2-3);
\draw[fill] (3,3) circle (1.5pt) coordinate (m-3-3);
\draw[fill] (4,3) circle (1.5pt) coordinate (m-4-3);
\draw[fill] (5,3) circle (1.5pt) coordinate (m-5-3);
\draw[fill] (6,3) circle (1.5pt) coordinate (m-6-3);

\draw[fill] (0,4) circle (1.5pt) coordinate (m-0-4);
\draw[fill] (1,4) circle (1.5pt) coordinate (m-1-4);
\draw[fill] (2,4) circle (1.5pt) coordinate (m-2-4);
\draw[fill] (3,4) circle (1.5pt) coordinate (m-3-4);
\draw[fill] (4,4) circle (1.5pt) coordinate (m-4-4);
\draw[fill] (5,4) circle (1.5pt) coordinate (m-5-4);
\draw[fill] (6,4) circle (1.5pt) coordinate (m-6-4);

\draw[fill] (1,5) circle (1.5pt) coordinate (m-1-5);
\draw[fill] (2,5) circle (1.5pt) coordinate (m-2-5);
\draw[fill] (3,5) circle (1.5pt) coordinate (m-3-5);
\draw[fill] (4,5) circle (1.5pt) coordinate (m-4-5);
\draw[fill] (5,5) circle (1.5pt) coordinate (m-5-5);
\draw[fill] (6,5) circle (1.5pt) coordinate (m-6-5);

\draw[fill] (5,6) circle (1.5pt) coordinate (m-5-6);

\draw (m-0-0) -- (m-6-1) -- (m-6-5) -- (m-5-6) -- (m-1-5) -- (m-0-4) -- cycle;

\draw(3, 6) node{\large $\tau$};

\draw(-0.5, 2) node{\large $\rho_1$};
\draw(6.5, 3) node{\large $\rho_2$};

\draw(3, 4.5) node{$\tau_2$};
\draw(3, 1.5) node{$\tau_1$};

\draw(0.5, 2.5) node{$\sigma_1$};
\draw(1.5, 2.5) node{$\sigma_2$};
\draw(2.5, 2.5) node{$\sigma_3$};
\draw(3.5, 2.5) node{$\sigma_4$};
\draw(4.5, 2.5) node{$\sigma_5$};
\draw(5.5, 2.5) node{$\sigma_6$};

\draw[fill, color=blue] (2.3,2.2) circle (1.5pt);

\draw [color=red, line width=1pt] (0,3)-- (6,3);
\draw [color=red, line width=1pt] (0,2)-- (6,2);

\draw [color=red, line width=1pt] (1,2)-- (1,3);
\draw [color=red, line width=1pt] (2,2)-- (2,3);
\draw [color=red, line width=1pt] (3,2)-- (3,3);
\draw [color=red, line width=1pt] (4,2)-- (4,3);
\draw [color=red, line width=1pt] (5,2)-- (5,3);

\draw [dashed] (0,2.2)-- (6,2.2);

\draw[fill] (0,2.2) circle (1pt);
\draw[fill] (6,2.2) circle (1pt);

\end{tikzpicture}
\captionsetup{labelformat=empty}
\caption{(1)}
\end{minipage}
\hspace{1cm}
\begin{minipage}{.4\textwidth}
\begin{tikzpicture}[scale=0.7]

\draw[fill] (1,0) circle (1.5pt) coordinate (m-1-0);
\draw[fill] (2,0) circle (1.5pt) coordinate (m-2-0);
\draw[fill] (3,0) circle (1.5pt) coordinate (m-3-0);
\draw[fill] (4,0) circle (1.5pt) coordinate (m-4-0);

\draw[fill] (1,1) circle (1.5pt) coordinate (m-1-1);
\draw[fill] (2,1) circle (1.5pt) coordinate (m-2-1);
\draw[fill] (3,1) circle (1.5pt) coordinate (m-3-1);
\draw[fill] (4,1) circle (1.5pt) coordinate (m-4-1);
\draw[fill] (5,1) circle (1.5pt) coordinate (m-5-1);
\draw[fill] (6,1) circle (1.5pt) coordinate (m-6-1);

\draw[fill] (0,2) circle (1.5pt) coordinate (m-0-2);
\draw[fill] (1,2) circle (1.5pt) coordinate (m-1-2);
\draw[fill] (2,2) circle (1.5pt) coordinate (m-2-2);
\draw[fill] (3,2) circle (1.5pt) coordinate (m-3-2);
\draw[fill] (4,2) circle (1.5pt) coordinate (m-4-2);
\draw[fill] (5,2) circle (1.5pt) coordinate (m-5-2);
\draw[fill] (6,2) circle (1.5pt) coordinate (m-6-2);

\draw[fill] (0,3) circle (1.5pt) coordinate (m-0-3);
\draw[fill] (1,3) circle (1.5pt) coordinate (m-1-3);
\draw[fill] (2,3) circle (1.5pt) coordinate (m-2-3);
\draw[fill] (3,3) circle (1.5pt) coordinate (m-3-3);
\draw[fill] (4,3) circle (1.5pt) coordinate (m-4-3);
\draw[fill] (5,3) circle (1.5pt) coordinate (m-5-3);

\draw[fill] (0,4) circle (1.5pt) coordinate (m-0-4);
\draw[fill] (1,4) circle (1.5pt) coordinate (m-1-4);
\draw[fill] (2,4) circle (1.5pt) coordinate (m-2-4);
\draw[fill] (3,4) circle (1.5pt) coordinate (m-3-4);
\draw[fill] (4,4) circle (1.5pt) coordinate (m-4-4);

\draw[fill] (0,5) circle (1.5pt) coordinate (m-0-5);
\draw[fill] (1,5) circle (1.5pt) coordinate (m-1-5);
\draw[fill] (2,5) circle (1.5pt) coordinate (m-2-5);
\draw[fill] (3,5) circle (1.5pt) coordinate (m-3-5);

\draw[fill] (2,6) circle (1.5pt) coordinate (m-2-6);

\draw (m-1-0) -- (m-4-0) -- (m-6-1) -- (m-6-2) -- (m-2-6) -- (m-0-5) -- (m-0-2) -- cycle;

\draw(3, 6) node{\large $\tau$};

\draw(-0.5, 3) node{\large $\rho_1$};
\draw(2, -0.5) node{\large $\rho_3$};
\draw(4.5, 4.5) node{\large $\rho_2$};

\draw(2, 4.5) node{$\tau_2$};
\draw(0.5, 2) node{$\tau_1$};
\draw(4.5, 2) node{$\tau_3$};

\draw(0.5, 3.5) node{$\sigma_{1}$};
\draw(1.5, 3.5) node{$\sigma_{2}$};
\draw(2.2, 3.5) node{\tiny$\delta$};
\draw(3, 3.5) node{$\sigma'_{2}$};
\draw(4, 3.5) node{$\sigma'_{1}$};
\draw(2.7, 0.5) node{$\sigma''_{1}$};
\draw(2.7, 1.5) node{$\sigma''_{2}$};
\draw(2.7, 2.5) node{$\sigma''_{3}$};

\draw[fill, color = blue] (2.3,3.2) circle (1.5pt);

\draw [dashed] (4.8,3.2)-- (0,3.2);
\draw [dashed] (2.3,3.2)-- (2.3,0);

\draw[fill] (0,3.2) circle (1pt);
\draw[fill] (4.8,3.2) circle (1pt);
\draw[fill] (2.3,0) circle (1pt);

\draw [color=red, line width=1pt] (0,3)-- (5,3);
\draw [color=red, line width=1pt] (0,4)-- (4,4);
\draw [color=red, line width=1pt] (1,3)-- (1,4);

\draw [color=red, line width=1pt] (2,3)-- (2,4);
\draw [color=red, line width=1pt] (2,3)-- (3,3);
\draw [color=red, line width=1pt] (2,4)-- (3,3);

\draw [color=red, line width=1pt] (2,0)-- (2,3);
\draw [color=red, line width=1pt] (3,0)-- (3,3);
\draw [color=red, line width=1pt] (2,1)-- (3,1);
\draw [color=red, line width=1pt] (2,2)-- (3,2);

\draw [color=red, line width=1pt] (4,3)-- (3,4);

\end{tikzpicture}
\captionsetup{labelformat=empty}
\caption{(2)}
\end{minipage}
\end{center}
\setcounter{figure}{4}
\caption{Subdivisions of easily decomposable $\tau \in \bar{\mathscr{P}}^{[2]}$.} 
\label{subdivisiontypes}
\end{figure}

\begin{notation}
We confuse $\Delta_{\tau,1}$ and the polyhedra of Figure \ref{subdivisiontypes} with the corresponding subsets of $\Delta_{\tau}$. We will follow similar conventions for cells of $\Delta_{\tau,1}$ and further subdivisions of $\tau$.
\end{notation}

The subdivisions of $\tau$ of Figure \ref{subdivisiontypes} lift to subdivisions of $\Delta_{\tau}$ as follows. In the case of Figure \ref{subdivisiontypes}(1), let $v_1$ and $v_2$ be the two endpoints $(0,0)$ and $(0,1)$ of $\Delta_{\tau,1}$. We subdivide $\Delta_{\tau}$ (following the notations of Figure \ref{subdivisiontypes}(1)) into $\Conv(v_1, \tau_1)$, $\Conv(v_2,\tau_2)$, and $\Conv(\Delta_{\tau,1}, \sigma_i)$ for $1 \leq i \leq 6$.
In the case of Figure \ref{subdivisiontypes}(2), let $v_1$, $v_2$ and $v_3$ be the vertices $(0,0)$, $(0,1)$, and $(1,0)$ of $\Delta_{\tau,1}$ respectively. Then we subdivide $\Delta_{\tau}$ (following the notations of Figure \ref{subdivisiontypes}(2)) into $\Conv(v_1, \tau_1)$, $\Conv(v_2, \tau_2)$, $\Conv(v_2, \tau_3)$, $\Conv(\langle v_1, v_2 \rangle, \sigma_i)$ for $i = 1,2$, $\Conv(\langle v_2, v_3 \rangle, \sigma'_j)$ for $j = 1,2$, $\Conv(\langle v_3, v_1 \rangle, \sigma''_k)$ for $k = 1,2,3$, and $\Conv(\Delta_{\tau,1}, \delta)$. The same construction can be done in general.

Note that the cells $\tau$ with the subdivisions of Figure \ref{subdivisiontypes} are not polyhedral manifolds (with boundary) in the sense of Definition \ref{polymani}. For instance, in Figure \ref{subdivisiontypes}(1), the intersection of $\tau_1$ and $\sigma_1$ is not a cell of $\tau_1$. To fix this, we choose integral points in the interiors of $\tau_i$ (this can always be done after a rescaling) and do star subdivisions as in Figure \ref{refinedsubdivisiontypes}. If $\tau_i$ is subdivided into $\tau_{ij}, \ 1 \leq j \leq l$, we lift this to a subdivision of $\Delta_{\tau}$ by subdividing $\Conv(v_{i}, \tau_i)$ into $\Conv(v_{i}, \tau_{ij}), \ 1 \leq j \leq l$.

\begin{figure}[t]
\begin{center}
\begin{minipage}{.4\textwidth}
\begin{tikzpicture}[scale=0.7]

\draw[fill] (0,0) circle (1.5pt) coordinate (m-0-0);

\draw[fill] (0,1) circle (1.5pt) coordinate (m-0-1);
\draw[fill] (1,1) circle (1.5pt) coordinate (m-1-1);
\draw[fill] (2,1) circle (1.5pt) coordinate (m-2-1);
\draw[fill] (3,1) circle (1.5pt) coordinate (m-3-1);
\draw[fill] (4,1) circle (1.5pt) coordinate (m-4-1);
\draw[fill] (5,1) circle (1.5pt) coordinate (m-5-1);
\draw[fill] (6,1) circle (1.5pt) coordinate (m-6-1);

\draw[fill] (0,2) circle (1.5pt) coordinate (m-0-2);
\draw[fill] (1,2) circle (1.5pt) coordinate (m-1-2);
\draw[fill] (2,2) circle (1.5pt) coordinate (m-2-2);
\draw[fill] (3,2) circle (1.5pt) coordinate (m-3-2);
\draw[fill] (4,2) circle (1.5pt) coordinate (m-4-2);
\draw[fill] (5,2) circle (1.5pt) coordinate (m-5-2);
\draw[fill] (6,2) circle (1.5pt) coordinate (m-6-2);

\draw[fill] (0,3) circle (1.5pt) coordinate (m-0-3);
\draw[fill] (1,3) circle (1.5pt) coordinate (m-1-3);
\draw[fill] (2,3) circle (1.5pt) coordinate (m-2-3);
\draw[fill] (3,3) circle (1.5pt) coordinate (m-3-3);
\draw[fill] (4,3) circle (1.5pt) coordinate (m-4-3);
\draw[fill] (5,3) circle (1.5pt) coordinate (m-5-3);
\draw[fill] (6,3) circle (1.5pt) coordinate (m-6-3);

\draw[fill] (0,4) circle (1.5pt) coordinate (m-0-4);
\draw[fill] (1,4) circle (1.5pt) coordinate (m-1-4);
\draw[fill] (2,4) circle (1.5pt) coordinate (m-2-4);
\draw[fill] (3,4) circle (1.5pt) coordinate (m-3-4);
\draw[fill] (4,4) circle (1.5pt) coordinate (m-4-4);
\draw[fill] (5,4) circle (1.5pt) coordinate (m-5-4);
\draw[fill] (6,4) circle (1.5pt) coordinate (m-6-4);

\draw[fill] (1,5) circle (1.5pt) coordinate (m-1-5);
\draw[fill] (2,5) circle (1.5pt) coordinate (m-2-5);
\draw[fill] (3,5) circle (1.5pt) coordinate (m-3-5);
\draw[fill] (4,5) circle (1.5pt) coordinate (m-4-5);
\draw[fill] (5,5) circle (1.5pt) coordinate (m-5-5);
\draw[fill] (6,5) circle (1.5pt) coordinate (m-6-5);

\draw[fill] (5,6) circle (1.5pt) coordinate (m-5-6);

\draw (m-0-0) -- (m-6-1) -- (m-6-5) -- (m-5-6) -- (m-1-5) -- (m-0-4) -- cycle;

\draw(3, 6) node{\large $\tau$};

\draw(-0.5, 2) node{\large $\rho_1$};
\draw(6.5, 3) node{\large $\rho_2$};

\draw(0.5, 2.5) node{$\sigma_1$};
\draw(1.5, 2.5) node{$\sigma_2$};
\draw(2.5, 2.5) node{$\sigma_3$};
\draw(3.5, 2.5) node{$\sigma_4$};
\draw(4.5, 2.5) node{$\sigma_5$};
\draw(5.5, 2.5) node{$\sigma_6$};

\draw [color=red, line width=1pt] (0,3)-- (6,3);
\draw [color=red, line width=1pt] (0,2)-- (6,2);

\draw [color=red, line width=1pt] (1,2)-- (1,3);
\draw [color=red, line width=1pt] (2,2)-- (2,3);
\draw [color=red, line width=1pt] (3,2)-- (3,3);
\draw [color=red, line width=1pt] (4,2)-- (4,3);
\draw [color=red, line width=1pt] (5,2)-- (5,3);

\draw [color=red] (2,1)-- (0,0);
\draw [color=red] (2,1)-- (0,2);
\draw [color=red] (2,1)-- (6,1);
\draw [color=red] (2,1)-- (1,2);
\draw [color=red] (2,1)-- (2,2);
\draw [color=red] (2,1)-- (3,2);
\draw [color=red] (2,1)-- (4,2);
\draw [color=red] (2,1)-- (5,2);
\draw [color=red] (2,1)-- (6,2);

\draw [color=red] (3,4)-- (0,3);
\draw [color=red] (3,4)-- (0,4);
\draw [color=red] (3,4)-- (1,5);
\draw [color=red] (3,4)-- (5,6);
\draw [color=red] (3,4)-- (6,5);
\draw [color=red] (3,4)-- (6,3);
\draw [color=red] (3,4)-- (1,3);
\draw [color=red] (3,4)-- (2,3);
\draw [color=red] (3,4)-- (3,3);
\draw [color=red] (3,4)-- (4,3);
\draw [color=red] (3,4)-- (5,3);

\end{tikzpicture}
\captionsetup{labelformat=empty}
\caption{(1)}
\end{minipage}
\hspace{1cm}
\begin{minipage}{.4\textwidth}
\begin{tikzpicture}[scale=0.7]

\draw[fill] (1,0) circle (1.5pt) coordinate (m-1-0);
\draw[fill] (2,0) circle (1.5pt) coordinate (m-2-0);
\draw[fill] (3,0) circle (1.5pt) coordinate (m-3-0);
\draw[fill] (4,0) circle (1.5pt) coordinate (m-4-0);

\draw[fill] (1,1) circle (1.5pt) coordinate (m-1-1);
\draw[fill] (2,1) circle (1.5pt) coordinate (m-2-1);
\draw[fill] (3,1) circle (1.5pt) coordinate (m-3-1);
\draw[fill] (4,1) circle (1.5pt) coordinate (m-4-1);
\draw[fill] (5,1) circle (1.5pt) coordinate (m-5-1);
\draw[fill] (6,1) circle (1.5pt) coordinate (m-6-1);

\draw[fill] (0,2) circle (1.5pt) coordinate (m-0-2);
\draw[fill] (1,2) circle (1.5pt) coordinate (m-1-2);
\draw[fill] (2,2) circle (1.5pt) coordinate (m-2-2);
\draw[fill] (3,2) circle (1.5pt) coordinate (m-3-2);
\draw[fill] (4,2) circle (1.5pt) coordinate (m-4-2);
\draw[fill] (5,2) circle (1.5pt) coordinate (m-5-2);
\draw[fill] (6,2) circle (1.5pt) coordinate (m-6-2);

\draw[fill] (0,3) circle (1.5pt) coordinate (m-0-3);
\draw[fill] (1,3) circle (1.5pt) coordinate (m-1-3);
\draw[fill] (2,3) circle (1.5pt) coordinate (m-2-3);
\draw[fill] (3,3) circle (1.5pt) coordinate (m-3-3);
\draw[fill] (4,3) circle (1.5pt) coordinate (m-4-3);
\draw[fill] (5,3) circle (1.5pt) coordinate (m-5-3);

\draw[fill] (0,4) circle (1.5pt) coordinate (m-0-4);
\draw[fill] (1,4) circle (1.5pt) coordinate (m-1-4);
\draw[fill] (2,4) circle (1.5pt) coordinate (m-2-4);
\draw[fill] (3,4) circle (1.5pt) coordinate (m-3-4);
\draw[fill] (4,4) circle (1.5pt) coordinate (m-4-4);

\draw[fill] (0,5) circle (1.5pt) coordinate (m-0-5);
\draw[fill] (1,5) circle (1.5pt) coordinate (m-1-5);
\draw[fill] (2,5) circle (1.5pt) coordinate (m-2-5);
\draw[fill] (3,5) circle (1.5pt) coordinate (m-3-5);

\draw[fill] (2,6) circle (1.5pt) coordinate (m-2-6);

\draw (m-1-0) -- (m-4-0) -- (m-6-1) -- (m-6-2) -- (m-2-6) -- (m-0-5) -- (m-0-2) -- cycle;

\draw(3, 6) node{\large $\tau$};

\draw(-0.5, 3) node{\large $\rho_1$};
\draw(2, -0.5) node{\large $\rho_3$};
\draw(4.5, 4.5) node{\large $\rho_2$};

\draw(0.5, 3.5) node{$\sigma_{1}$};
\draw(1.5, 3.5) node{$\sigma_{2}$};
\draw(2.2, 3.5) node{\tiny$\delta$};
\draw(3, 3.5) node{$\sigma'_{2}$};
\draw(4, 3.5) node{$\sigma'_{1}$};
\draw(2.7, 0.5) node{$\sigma''_{1}$};
\draw(2.7, 1.5) node{$\sigma''_{2}$};
\draw(2.7, 2.5) node{$\sigma''_{3}$};

\draw [color=red, line width=1pt] (0,3)-- (5,3);
\draw [color=red, line width=1pt] (0,4)-- (4,4);
\draw [color=red, line width=1pt] (1,3)-- (1,4);

\draw [color=red, line width=1pt] (2,3)-- (2,4);
\draw [color=red, line width=1pt] (2,3)-- (3,3);
\draw [color=red, line width=1pt] (2,4)-- (3,3);

\draw [color=red, line width=1pt] (2,0)-- (2,3);
\draw [color=red, line width=1pt] (3,0)-- (3,3);
\draw [color=red, line width=1pt] (2,1)-- (3,1);
\draw [color=red, line width=1pt] (2,2)-- (3,2);

\draw [color=red, line width=1pt] (4,3)-- (3,4);

\draw [color=red] (1,2)-- (0,2);
\draw [color=red] (1,2)-- (0,3);
\draw [color=red] (1,2)-- (1,0);
\draw [color=red] (1,2)-- (2,0);
\draw [color=red] (1,2)-- (2,1);
\draw [color=red] (1,2)-- (2,2);
\draw [color=red] (1,2)-- (2,3);
\draw [color=red] (1,2)-- (1,3);

\draw [color=red] (4,1)-- (3,0);
\draw [color=red] (4,1)-- (4,0);
\draw [color=red] (4,1)-- (6,1);
\draw [color=red] (4,1)-- (6,2);
\draw [color=red] (4,1)-- (5,3);
\draw [color=red] (4,1)-- (3,1);
\draw [color=red] (4,1)-- (3,2);
\draw [color=red] (4,1)-- (3,3);
\draw [color=red] (4,1)-- (4,3);

\draw [color=red] (2,5)-- (0,5);
\draw [color=red] (2,5)-- (2,6);
\draw [color=red] (2,5)-- (4,4);
\draw [color=red] (2,5)-- (3,4);
\draw [color=red] (2,5)-- (2,4);
\draw [color=red] (2,5)-- (1,4);
\draw [color=red] (2,5)-- (0,4);

\end{tikzpicture}
\captionsetup{labelformat=empty}
\caption{(2)}
\end{minipage}
\end{center}
\setcounter{figure}{5}
\caption{Refined subdivisions of easily decomposable $\tau \in \bar{\mathscr{P}}^{[2]}$.} 
\label{refinedsubdivisiontypes}
\end{figure}

The subdivisions of Figure \ref{refinedsubdivisiontypes} give rise to partial resolutions that are log smooth in the neighbourhoods of strata corresponding to $\Conv(v_i, \tau_{ij})$.\footnote{Indeed, the toric variety $X$ defined by the cone over $\Conv(v_i, \tau_{ij})$ splits as a product $X\cong  Y \times \mathbb{A}^1$ where $Y$ is the toric variety defined by the cone over the face $\tau_{ij}$ of $\Conv(v_i, \tau_{ij})$. The map $X \to \mathbb{A}^1$ is just the trivial fibration $Y \times \mathbb{A}^1 \to \mathbb{A}^1$, which is clearly log smooth.}  We note that the subdivisions do not change the singular locus in the sense that the singular loci of the polytopes in the subdivision of $\tau$ can be understood as the dashed lines in Figure \ref{subdivisiontypes}.

We still need to resolve the singularities of $\Conv(\Delta_{\tau,1}, \sigma_i)$ in the case of Figure \ref{refinedsubdivisiontypes}(1) and of $\Conv(\langle v_1, v_2 \rangle, \sigma_i)$, $\Conv(\langle v_2, v_3 \rangle, \sigma'_j)$, $\Conv(\langle v_3, v_1 \rangle, \sigma''_k)$, and $\Conv(\Delta_{\tau,1}, \delta)$ in the case of Figure \ref{refinedsubdivisiontypes}(2). Note that $\Conv(\Delta_{\tau,1}, \delta)$ in the case of Figure \ref{refinedsubdivisiontypes}(2) is as in Figure \ref{locmodeltriangle} and all the other convex hulls are $AGL(3,\mathbb{Z})$-equivalent to $$\Delta_{\sigma} \cong \Conv\left((\sigma \times \{e_0\}) \cup  \Delta_{\sigma,1} \right)$$ for $\sigma$ a standard square (i.e. the convex hull of $(0,0), (1,0), (0,1)$, and $(1,1)$) and $\Delta_{\sigma,1}$ the line segment connecting $(0,0)$ and $(0,1)$. One can resolve the corresponding local model by subdividing $\Delta_{\sigma}$ with $2$ hyperplane sections, similarly to the subdivisions of Figure \ref{locmodeltriangle}. We give a possible subdivision in Figure \ref{locmodelsquare}.

\begin{figure}[t]
\begin{center}

\begin{minipage}{.4\textwidth}
\begin{tikzpicture}[line cap=round,line join=round,>=triangle 45,x=1cm,y=1cm, scale = 3]
\clip(-0.2,-0.2) rectangle (1.5,1.5);
\draw (0,0)-- (1,0);
\draw (0,0)-- (0,1);
\draw (0,1)-- (1,0);
\draw (0,1)-- (0.33,1.33);
\draw[dashed] (0,0)-- (0.33,0.33);
\draw[dashed] (1.33,0.33)-- (0.33,0.33);
\draw[dashed] (0.33,0.33)-- (0.33,1.33);

\draw (1.33,0.33)-- (0.33,1.33);
\draw (1.33,0.33)-- (1,0);

\draw [fill] (0,0) circle (0.5pt);
\draw [fill] (0.33,0.33) circle (0.5pt);
\draw [fill] (0.33,1.33) circle (0.5pt);
\draw [fill] (0,1) circle (0.5pt);
\draw [fill] (1,0) circle (0.5pt);

\draw [fill] (1.33,0.33) circle (0.5pt);

\draw(0.5, 0.15) node{\scriptsize $\tau$};
\draw(0.065, 1.2) node{\scriptsize $\Delta_{\tau,1}$};

\draw(-0.06, 0.1) node{\scriptsize $v_1$};
\draw(1, 0.1) node{\scriptsize $v_2$};
\draw(0.25, 0.43) node{\scriptsize $v_4$};
\draw(1.33, 0.43) node{\scriptsize $v_3$};

\end{tikzpicture}
\captionsetup{labelformat=empty}
\caption{Local model}
\end{minipage}
\begin{minipage}{.4\textwidth}
\begin{tikzpicture}[line cap=round,line join=round,>=triangle 45,x=1cm,y=1cm, scale = 3]
\clip(-0.2,-0.2) rectangle (1.5,1.5);
\draw (0,0)-- (1,0);
\draw (0,0)-- (0,1);
\draw (0,1)-- (1,0);
\draw (0,1)-- (0.33,1.33);
\draw[dashed] (0,0)-- (0.33,0.33);
\draw[dashed] (1.33,0.33)-- (0.33,0.33);
\draw[dashed] (0.33,0.33)-- (0.33,1.33);

\draw (1.33,0.33)-- (0.33,1.33);
\draw (1.33,0.33)-- (1,0);

\draw [fill] (0,0) circle (0.5pt);
\draw [fill] (0.33,0.33) circle (0.5pt);
\draw [fill] (0.33,1.33) circle (0.5pt);
\draw [fill] (0,1) circle (0.5pt);
\draw [fill] (1,0) circle (0.5pt);

\draw [fill] (1.33,0.33) circle (0.5pt);

\draw(0.5, 0.15) node{\scriptsize $\tau$};
\draw(0.065, 1.2) node{\scriptsize $\Delta_{\tau,1}$};

\draw(-0.06, 0.1) node{\scriptsize $v_1$};
\draw(1, 0.1) node{\scriptsize $v_2$};
\draw(0.25, 0.43) node{\scriptsize $v_4$};
\draw(1.33, 0.43) node{\scriptsize $v_3$};

\filldraw[draw=black, fill=blue, opacity=.2] (0.33,0.33) -- (0,1) -- (1.33,0.33) -- cycle;

\filldraw[draw=black, fill=blue, opacity=.2] (0.33,0.33) -- (0,1) -- (1,0) -- cycle;

\end{tikzpicture}
\captionsetup{labelformat=empty}
\caption{Resolution}
\end{minipage}
\end{center}
\setcounter{figure}{6}
\caption{Resolution of the local model $\tilde{\bar{X}}_{\sigma,x}$ for $\sigma$ a standard square.}
\label{locmodelsquare}
\end{figure}

As in the case of Figure \ref{locmodeltriangle}, the subdivision of Figure \ref{locmodelsquare} corresponds to a sequence of blowups of the irreducible components corresponding to the vertices of $\sigma$. Namely, the chosen subdivision corresponds to blowing up the divisor corresponding to $v_4$ and then blowing up the divisor corresponding to $v_2$.

In Figure \ref{singularlocuslsquare}, we give the transformation of the singular locus (in blue) of $\sigma$ corresponding to the subdivision of Figure \ref{locmodelsquare}. Again, at the level of the affine manifolds, the blowup of an irreducible component corresponding to $v_i$ can be visualized as pulling the singular locus into $v_i$.

\begin{figure}[H]
\begin{center}

\begin{minipage}{.4\textwidth}
\begin{tikzpicture}[line cap=round,line join=round,>=triangle 45,x=1cm,y=1cm, scale = 3]
\clip(-0.2,-0.2) rectangle (1.2,1.2);
\draw (0,1)-- (1,1);
\draw (0,0)-- (0,1);
\draw (0,0)-- (1,0);
\draw (1,0)-- (1,1);

\draw [fill] (0,0) circle (0.5pt);
\draw [fill] (0,1) circle (0.5pt);
\draw [fill] (1,0) circle (0.5pt);
\draw [fill] (1,1) circle (0.5pt);

\draw(-0.07, 0.07) node{$v_1$};
\draw(1.07, 0.07) node{$v_2$};
\draw(1.07, 1.07) node{$v_3$};
\draw(-0.06, 1.07) node{$v_4$};

\draw [color=blue, line width=1pt] (0,0.5)-- (1,0.5);

\end{tikzpicture}
\captionsetup{labelformat=empty}
\caption{$\sigma$ before \\ resolution}
\end{minipage}
\begin{minipage}{.3\textwidth}
\begin{tikzpicture}[line cap=round,line join=round,>=triangle 45,x=1cm,y=1cm, scale = 3]
\clip(-0.2,-0.2) rectangle (1.2,1.2);
\draw (0,1)-- (1,1);
\draw (0,0)-- (0,1);
\draw (0,0)-- (1,0);
\draw (1,0)-- (1,1);

\draw [fill] (0,0) circle (0.5pt);
\draw [fill] (0,1) circle (0.5pt);
\draw [fill] (1,0) circle (0.5pt);
\draw [fill] (1,1) circle (0.5pt);

\draw(-0.07, 0.07) node{$v_1$};
\draw(1.07, 0.07) node{$v_2$};
\draw(1.07, 1.07) node{$v_3$};
\draw(-0.06, 1.07) node{$v_4$};

\draw [color=blue, line width=1pt] (0,1)-- (1,0);

\end{tikzpicture}
\captionsetup{labelformat=empty}
\caption{$\sigma$ after \\ resolution}
\end{minipage}
\end{center}
\setcounter{figure}{7}
\caption{Transformation of $\sigma$ corresponding to Figure \ref{locmodelsquare}.}
\label{singularlocuslsquare}
\end{figure}

There are $6$ possible resolutions of the local model corresponding to Figure \ref{locmodelsquare}. They correspond to choosing a different order of blowups of the irreducible components or, tropically, to taking the singular locus in Figure \ref{singularlocuslsquare} to be one of the $\langle v_1, v_2 \rangle, \ \langle v_1, v_3 \rangle, \ \langle v_4, v_2 \rangle, \ \langle v_4, v_3 \rangle$ and subdividing the square by one of the diagonals (if the singular locus is itself a diagonal, the subdivision has to coincide with it). In particular, one can always arrange any choice of diagonals subdividing the two square faces of $\Delta_{\sigma}$ that intersect $\Delta_{\sigma,1}$.

Let us return to the partial resolutions of the (possibly rescaled) local models with an easily decomposable $\tau$ we constructed above. To resolve the local model $\tilde{\bar{X}}_{\tau,x}$ to a log smooth degeneration in the case of Figure \ref{refinedsubdivisiontypes}(1), we just subdivide the $\Conv(\Delta_{\tau,1}, \sigma_i)$ for $1 \leq i \leq 6$ as in Figure \ref{locmodelsquare}. The subdivisions of $\Conv(\Delta_{\tau,1}, \sigma_i)$ and $\Conv(\Delta_{\tau,1}, \sigma_{i+1})$ for $1 \leq i \leq 5$ have to be compatible, but it is clear from the above that this can be arranged. In the case of Figure \ref{refinedsubdivisiontypes}(2), we choose a subdivision of $\Conv(\Delta_{\tau,1}, \delta)$ as in Figure \ref{locmodeltriangle} and choose the subdivisions of the convex hulls $AGL(3,\mathbb{Z})$-equivalent to the convex hull of Figure \ref{locmodelsquare} as in Figure \ref{locmodelsquare}. Again, this should be compatible on the intersections of cells of the subdivision, but this can be arranged. The fact that the corresponding resolution is log smooth follows from the fact that all the cells in the subdivisions of Figures \ref{locmodeltriangle} and \ref{locmodelsquare} are standard simplices. The general construction is similar. 

The singular locus of the subdivision of (a rescaling of) $\tau$ can be recovered from the singular loci of Figures \ref{locmodeltriangle} and \ref{locmodelsquare} and depends on the particular subdivisions chosen. For example, continuing with the setup of Figure \ref{refinedsubdivisiontypes}, the resulting subdivisions and singular loci of $\tau$ can look as in Figure \ref{finalsubdivisiontypes}. Note that the singular loci are connected.

\begin{figure}[H]
\begin{center}
\begin{minipage}{.4\textwidth}
\begin{tikzpicture}[scale=0.7]

\draw[fill] (0,0) circle (1.5pt) coordinate (m-0-0);

\draw[fill] (0,1) circle (1.5pt) coordinate (m-0-1);
\draw[fill] (1,1) circle (1.5pt) coordinate (m-1-1);
\draw[fill] (2,1) circle (1.5pt) coordinate (m-2-1);
\draw[fill] (3,1) circle (1.5pt) coordinate (m-3-1);
\draw[fill] (4,1) circle (1.5pt) coordinate (m-4-1);
\draw[fill] (5,1) circle (1.5pt) coordinate (m-5-1);
\draw[fill] (6,1) circle (1.5pt) coordinate (m-6-1);

\draw[fill] (0,2) circle (1.5pt) coordinate (m-0-2);
\draw[fill] (1,2) circle (1.5pt) coordinate (m-1-2);
\draw[fill] (2,2) circle (1.5pt) coordinate (m-2-2);
\draw[fill] (3,2) circle (1.5pt) coordinate (m-3-2);
\draw[fill] (4,2) circle (1.5pt) coordinate (m-4-2);
\draw[fill] (5,2) circle (1.5pt) coordinate (m-5-2);
\draw[fill] (6,2) circle (1.5pt) coordinate (m-6-2);

\draw[fill] (0,3) circle (1.5pt) coordinate (m-0-3);
\draw[fill] (1,3) circle (1.5pt) coordinate (m-1-3);
\draw[fill] (2,3) circle (1.5pt) coordinate (m-2-3);
\draw[fill] (3,3) circle (1.5pt) coordinate (m-3-3);
\draw[fill] (4,3) circle (1.5pt) coordinate (m-4-3);
\draw[fill] (5,3) circle (1.5pt) coordinate (m-5-3);
\draw[fill] (6,3) circle (1.5pt) coordinate (m-6-3);

\draw[fill] (0,4) circle (1.5pt) coordinate (m-0-4);
\draw[fill] (1,4) circle (1.5pt) coordinate (m-1-4);
\draw[fill] (2,4) circle (1.5pt) coordinate (m-2-4);
\draw[fill] (3,4) circle (1.5pt) coordinate (m-3-4);
\draw[fill] (4,4) circle (1.5pt) coordinate (m-4-4);
\draw[fill] (5,4) circle (1.5pt) coordinate (m-5-4);
\draw[fill] (6,4) circle (1.5pt) coordinate (m-6-4);

\draw[fill] (1,5) circle (1.5pt) coordinate (m-1-5);
\draw[fill] (2,5) circle (1.5pt) coordinate (m-2-5);
\draw[fill] (3,5) circle (1.5pt) coordinate (m-3-5);
\draw[fill] (4,5) circle (1.5pt) coordinate (m-4-5);
\draw[fill] (5,5) circle (1.5pt) coordinate (m-5-5);
\draw[fill] (6,5) circle (1.5pt) coordinate (m-6-5);

\draw[fill] (5,6) circle (1.5pt) coordinate (m-5-6);

\draw (m-0-0) -- (m-6-1) -- (m-6-5) -- (m-5-6) -- (m-1-5) -- (m-0-4) -- cycle;

\draw(3, 6) node{\large $\tau$};

\draw(-0.5, 2) node{\large $\rho_1$};
\draw(6.5, 3) node{\large $\rho_2$};

\draw [color=red, line width=1pt] (0,3)-- (6,3);
\draw [color=red, line width=1pt] (0,2)-- (6,2);

\draw [color=red, line width=1pt] (1,2)-- (1,3);
\draw [color=red, line width=1pt] (2,2)-- (2,3);
\draw [color=red, line width=1pt] (3,2)-- (3,3);
\draw [color=red, line width=1pt] (4,2)-- (4,3);
\draw [color=red, line width=1pt] (5,2)-- (5,3);

\draw [color=red] (2,1)-- (0,0);
\draw [color=red] (2,1)-- (0,2);
\draw [color=red] (2,1)-- (6,1);
\draw [color=red] (2,1)-- (1,2);
\draw [color=red] (2,1)-- (2,2);
\draw [color=red] (2,1)-- (3,2);
\draw [color=red] (2,1)-- (4,2);
\draw [color=red] (2,1)-- (5,2);
\draw [color=red] (2,1)-- (6,2);

\draw [color=red] (3,4)-- (0,3);
\draw [color=red] (3,4)-- (0,4);
\draw [color=red] (3,4)-- (1,5);
\draw [color=red] (3,4)-- (5,6);
\draw [color=red] (3,4)-- (6,5);
\draw [color=red] (3,4)-- (6,3);
\draw [color=red] (3,4)-- (1,3);
\draw [color=red] (3,4)-- (2,3);
\draw [color=red] (3,4)-- (3,3);
\draw [color=red] (3,4)-- (4,3);
\draw [color=red] (3,4)-- (5,3);

\draw [color=blue, line width=1pt] (0,3)-- (1,2) -- (3,2) -- (4,3) -- (5,2) -- (6,2);

\draw [color=red] (1,2)-- (2,3);
\draw [color=red] (2,2)-- (3,3);
\draw [color=red] (5,2)-- (6,3);

\end{tikzpicture}
\captionsetup{labelformat=empty}
\caption{(1)}
\end{minipage}
\hspace{1cm}
\begin{minipage}{.4\textwidth}
\begin{tikzpicture}[scale=0.7]

\draw[fill] (1,0) circle (1.5pt) coordinate (m-1-0);
\draw[fill] (2,0) circle (1.5pt) coordinate (m-2-0);
\draw[fill] (3,0) circle (1.5pt) coordinate (m-3-0);
\draw[fill] (4,0) circle (1.5pt) coordinate (m-4-0);

\draw[fill] (1,1) circle (1.5pt) coordinate (m-1-1);
\draw[fill] (2,1) circle (1.5pt) coordinate (m-2-1);
\draw[fill] (3,1) circle (1.5pt) coordinate (m-3-1);
\draw[fill] (4,1) circle (1.5pt) coordinate (m-4-1);
\draw[fill] (5,1) circle (1.5pt) coordinate (m-5-1);
\draw[fill] (6,1) circle (1.5pt) coordinate (m-6-1);

\draw[fill] (0,2) circle (1.5pt) coordinate (m-0-2);
\draw[fill] (1,2) circle (1.5pt) coordinate (m-1-2);
\draw[fill] (2,2) circle (1.5pt) coordinate (m-2-2);
\draw[fill] (3,2) circle (1.5pt) coordinate (m-3-2);
\draw[fill] (4,2) circle (1.5pt) coordinate (m-4-2);
\draw[fill] (5,2) circle (1.5pt) coordinate (m-5-2);
\draw[fill] (6,2) circle (1.5pt) coordinate (m-6-2);

\draw[fill] (0,3) circle (1.5pt) coordinate (m-0-3);
\draw[fill] (1,3) circle (1.5pt) coordinate (m-1-3);
\draw[fill] (2,3) circle (1.5pt) coordinate (m-2-3);
\draw[fill] (3,3) circle (1.5pt) coordinate (m-3-3);
\draw[fill] (4,3) circle (1.5pt) coordinate (m-4-3);
\draw[fill] (5,3) circle (1.5pt) coordinate (m-5-3);

\draw[fill] (0,4) circle (1.5pt) coordinate (m-0-4);
\draw[fill] (1,4) circle (1.5pt) coordinate (m-1-4);
\draw[fill] (2,4) circle (1.5pt) coordinate (m-2-4);
\draw[fill] (3,4) circle (1.5pt) coordinate (m-3-4);
\draw[fill] (4,4) circle (1.5pt) coordinate (m-4-4);

\draw[fill] (0,5) circle (1.5pt) coordinate (m-0-5);
\draw[fill] (1,5) circle (1.5pt) coordinate (m-1-5);
\draw[fill] (2,5) circle (1.5pt) coordinate (m-2-5);
\draw[fill] (3,5) circle (1.5pt) coordinate (m-3-5);

\draw[fill] (2,6) circle (1.5pt) coordinate (m-2-6);

\draw (m-1-0) -- (m-4-0) -- (m-6-1) -- (m-6-2) -- (m-2-6) -- (m-0-5) -- (m-0-2) -- cycle;

\draw(3, 6) node{\large $\tau$};

\draw(-0.5, 3) node{\large $\rho_1$};
\draw(2, -0.5) node{\large $\rho_3$};
\draw(4.5, 4.5) node{\large $\rho_2$};

\draw [color=red, line width=1pt] (0,3)-- (5,3);
\draw [color=red, line width=1pt] (0,4)-- (4,4);
\draw [color=red, line width=1pt] (1,3)-- (1,4);

\draw [color=red, line width=1pt] (2,3)-- (2,4);
\draw [color=red, line width=1pt] (2,3)-- (3,3);
\draw [color=red, line width=1pt] (2,4)-- (3,3);

\draw [color=red, line width=1pt] (2,0)-- (2,3);
\draw [color=red, line width=1pt] (3,0)-- (3,3);
\draw [color=red, line width=1pt] (2,1)-- (3,1);
\draw [color=red, line width=1pt] (2,2)-- (3,2);

\draw [color=red, line width=1pt] (4,3)-- (3,4);

\draw [color=red] (1,2)-- (0,2);
\draw [color=red] (1,2)-- (0,3);
\draw [color=red] (1,2)-- (1,0);
\draw [color=red] (1,2)-- (2,0);
\draw [color=red] (1,2)-- (2,1);
\draw [color=red] (1,2)-- (2,2);
\draw [color=red] (1,2)-- (2,3);
\draw [color=red] (1,2)-- (1,3);

\draw [color=red] (4,1)-- (3,0);
\draw [color=red] (4,1)-- (4,0);
\draw [color=red] (4,1)-- (6,1);
\draw [color=red] (4,1)-- (6,2);
\draw [color=red] (4,1)-- (5,3);
\draw [color=red] (4,1)-- (3,1);
\draw [color=red] (4,1)-- (3,2);
\draw [color=red] (4,1)-- (3,3);
\draw [color=red] (4,1)-- (4,3);

\draw [color=red] (2,5)-- (0,5);
\draw [color=red] (2,5)-- (2,6);
\draw [color=red] (2,5)-- (4,4);
\draw [color=red] (2,5)-- (3,4);
\draw [color=red] (2,5)-- (2,4);
\draw [color=red] (2,5)-- (1,4);
\draw [color=red] (2,5)-- (0,4);

\draw [color=blue, line width=1pt] (0,3)-- (1,4) -- (2,4) -- (3,3) -- (2,2) -- (2,1) -- (2,0);

\draw [color=blue, line width=1pt] (2,4) -- (4,3) -- (5,3);

\draw [color=red] (1,3) -- (2,4);

\draw [color=red] (2,0) -- (3,1) -- (2,2);

\draw [color=red] (4,3) -- (4,4);

\end{tikzpicture}
\captionsetup{labelformat=empty}
\caption{(2)}
\end{minipage}
\end{center}
\setcounter{figure}{8}
\caption{Final subdivisions and singular loci of easily decomposable $\tau \in \bar{\mathscr{P}}^{[2]}$.} 
\label{finalsubdivisiontypes}
\end{figure}

\subsubsection{Resolutions of the local models in general}

Let $x$ be a singular point contained in the minimal stratum $\bar{X}_{\tau}$ for $\tau \in \mathscr{P}^{[2]}$ a codimension $1$ cell (we no longer assume that $\tau$ is easily decomposable). By the same analysis as in Section \ref{resolveed}, the singular locus of $\tau$ is of one of the two types in Figure \ref{singlocitypes}. Still, in the case of Figure \ref{singlocitypes}(1), $\Delta_{\tau,g} = \Delta_{\tau,1}$ is the line segment connecting the points $(0,0)$ and $(0,1)$. In the case of Figure \ref{singlocitypes}(2), $\Delta_{\tau,g} = \Delta_{\tau,1}$ is a standard triangle. We want to produce a subdivision of $\tau$ and lift it to a subdivision of $\Delta_{\tau}$. Since we no longer assume that $\tau$ is easily decomposable, the subdivision will be more involved. 

First, we do partial subdivisions similar to Figure \ref{refinedsubdivisiontypes}, see Figure \ref{singlocitypesgeneral}. Here we choose new $\tau$-s to illustrate the concept (it is easy to check that they are not easily decomposable). We still require that every edge of $\tau$ subdivided by our subdivision is a union of three (not two) edges of the subdivision. The only difference with Figure \ref{refinedsubdivisiontypes} is that we no longer require that $\sigma_i, \sigma'_j, \sigma''_k$ are standard squares and allow them to be general parallelograms. 

It is clear that for all $\tau$, such subdivisions exist for suitable rescalings $k\tau$ of $\tau$ with $k \in \mathbb{Z}_{>0}$ large enough. Indeed, in the case of Figure \ref{singlocitypesgeneral}(1), we can choose length-one segments on $\rho_1$ and $\rho_2$ (after a rescaling by $3$ if one of the edges is length $1$ or $2$) and connect them as in Figure \ref{singlocitypesgeneral}(1) creating a parallelogram $\sigma$. Possibly after an additional rescaling (and a choice of the \enquote{new} $\sigma$ inside the rescaling of the \enquote{old} $\sigma$), we can ensure that the other two polygons $\tau_1$ and $\tau_2$ of the subdivision have interior points. Then we can do a star subdivision of $\tau_1$ and $\tau_2$ into triangles $\tau_{1j}$ and $\tau_{2j}$, and subdivide the parallelogram $\sigma$ as in Figure \ref{singlocitypesgeneral}(1). In the case of Figure \ref{singlocitypesgeneral}(2), we can choose a standard triangle $\delta$ contained in $\tau$ with edges parallel to $\rho_1, \rho_2, \rho_3$ and connect the edges of the triangle to length-one segments in the corresponding edges (again, possibly after some rescalings). This subdivides $\tau$ into $\delta$, parallelograms $\sigma, \sigma', \sigma''$ (with obvious notations), and $3$ polygons $\tau_1, \tau_2, \tau_3$ one of which can be non-convex. Despite this non-convexity, it is easy to see that after an additional rescaling (and a choice of the \enquote{new} $\delta,\sigma, \sigma', \sigma''$ inside the rescaling of the \enquote{old} $\delta,\sigma, \sigma', \sigma''$) we can ensure that $\tau_i, \ 1 \leq i \leq 3$ admit star subdivisions into triangles $\tau_{ij}$ as in Figure \ref{singlocitypesgeneral}(2). We perform these star subdivisions and subdivide the parallelograms $\sigma, \sigma', \sigma''$ as in Figure \ref{singlocitypesgeneral}(2). This gives the required subdivisions.

\begin{figure}[t]
\begin{center}
\begin{minipage}{.4\textwidth}
\begin{tikzpicture}[scale=0.7]

\draw[fill] (0,0) circle (1.5pt) coordinate (m-0-0);
\draw[fill] (1,0) circle (1.5pt) coordinate (m-1-0);
\draw[fill] (2,0) circle (1.5pt) coordinate (m-2-0);
\draw[fill] (3,0) circle (1.5pt) coordinate (m-3-0);

\draw[fill] (0,1) circle (1.5pt) coordinate (m-0-1);
\draw[fill] (1,1) circle (1.5pt) coordinate (m-1-1);
\draw[fill] (2,1) circle (1.5pt) coordinate (m-2-1);
\draw[fill] (3,1) circle (1.5pt) coordinate (m-3-1);

\draw[fill] (0,2) circle (1.5pt) coordinate (m-0-2);
\draw[fill] (1,2) circle (1.5pt) coordinate (m-1-2);
\draw[fill] (2,2) circle (1.5pt) coordinate (m-2-2);
\draw[fill] (3,2) circle (1.5pt) coordinate (m-3-2);
\draw[fill] (4,2) circle (1.5pt) coordinate (m-4-2);

\draw[fill] (0,3) circle (1.5pt) coordinate (m-0-3);
\draw[fill] (1,3) circle (1.5pt) coordinate (m-1-3);
\draw[fill] (2,3) circle (1.5pt) coordinate (m-2-3);
\draw[fill] (3,3) circle (1.5pt) coordinate (m-3-3);
\draw[fill] (4,3) circle (1.5pt) coordinate (m-4-3);
\draw[fill] (5,3) circle (1.5pt) coordinate (m-5-3);

\draw[fill] (1,4) circle (1.5pt) coordinate (m-1-4);
\draw[fill] (2,4) circle (1.5pt) coordinate (m-2-4);
\draw[fill] (3,4) circle (1.5pt) coordinate (m-3-4);
\draw[fill] (4,4) circle (1.5pt) coordinate (m-4-4);
\draw[fill] (5,4) circle (1.5pt) coordinate (m-5-4);
\draw[fill] (6,4) circle (1.5pt) coordinate (m-6-4);

\draw[fill] (1,5) circle (1.5pt) coordinate (m-1-5);
\draw[fill] (2,5) circle (1.5pt) coordinate (m-2-5);
\draw[fill] (3,5) circle (1.5pt) coordinate (m-3-5);
\draw[fill] (4,5) circle (1.5pt) coordinate (m-4-5);
\draw[fill] (5,5) circle (1.5pt) coordinate (m-5-5);
\draw[fill] (6,5) circle (1.5pt) coordinate (m-6-5);

\draw[fill] (4,6) circle (1.5pt) coordinate (m-4-6);
\draw[fill] (5,6) circle (1.5pt) coordinate (m-5-6);
\draw[fill] (6,6) circle (1.5pt) coordinate (m-6-6);

\draw[fill] (6,7) circle (1.5pt) coordinate (m-6-7);

\draw (m-0-0) -- (m-3-0) -- (m-6-4) -- (m-6-7) -- (m-1-5) -- (m-0-3) -- cycle;

\draw(3, 7) node{\large $\tau$};

\draw(-0.5, 2) node{\large $\rho_1$};
\draw(6.5, 5) node{\large $\rho_2$};

\draw(1.5, 2.5) node{$\sigma_1$};

\draw(4.5, 4.5) node{$\sigma_2$};

\draw [color=red, line width=1pt] (0,2)-- (6,6);
\draw [color=red, line width=1pt] (0,1)-- (6,5);

\draw [color=red, line width=1pt] (3,3)-- (3,4);

\draw [color=red] (2,1)-- (0,0);
\draw [color=red] (2,1)-- (0,1);
\draw [color=red] (2,1)-- (3,3);
\draw [color=red] (2,1)-- (6,5);
\draw [color=red] (2,1)-- (6,4);
\draw [color=red] (2,1)-- (3,0);

\draw [color=red] (2,4)-- (0,2);
\draw [color=red] (2,4)-- (0,3);
\draw [color=red] (2,4)-- (1,5);
\draw [color=red] (2,4)-- (3,4);
\draw [color=red] (2,4)-- (6,6);
\draw [color=red] (2,4)-- (6,7);

\end{tikzpicture}
\captionsetup{labelformat=empty}
\caption{(1)}
\end{minipage}
\hspace{1cm}
\begin{minipage}{.4\textwidth}
\begin{tikzpicture}[scale=0.7]

\draw[fill] (0,0) circle (1.5pt) coordinate (m-0-0);
\draw[fill] (1,0) circle (1.5pt) coordinate (m-1-0);
\draw[fill] (2,0) circle (1.5pt) coordinate (m-2-0);
\draw[fill] (3,0) circle (1.5pt) coordinate (m-3-0);

\draw[fill] (0,1) circle (1.5pt) coordinate (m-0-1);
\draw[fill] (1,1) circle (1.5pt) coordinate (m-1-1);
\draw[fill] (2,1) circle (1.5pt) coordinate (m-2-1);
\draw[fill] (3,1) circle (1.5pt) coordinate (m-3-1);

\draw[fill] (0,2) circle (1.5pt) coordinate (m-0-2);
\draw[fill] (1,2) circle (1.5pt) coordinate (m-1-2);
\draw[fill] (2,2) circle (1.5pt) coordinate (m-2-2);
\draw[fill] (3,2) circle (1.5pt) coordinate (m-3-2);
\draw[fill] (4,2) circle (1.5pt) coordinate (m-4-2);

\draw[fill] (0,3) circle (1.5pt) coordinate (m-0-3);
\draw[fill] (1,3) circle (1.5pt) coordinate (m-1-3);
\draw[fill] (2,3) circle (1.5pt) coordinate (m-2-3);
\draw[fill] (3,3) circle (1.5pt) coordinate (m-3-3);
\draw[fill] (4,3) circle (1.5pt) coordinate (m-4-3);
\draw[fill] (5,3) circle (1.5pt) coordinate (m-5-3);

\draw[fill] (1,4) circle (1.5pt) coordinate (m-1-4);
\draw[fill] (2,4) circle (1.5pt) coordinate (m-2-4);
\draw[fill] (3,4) circle (1.5pt) coordinate (m-3-4);
\draw[fill] (4,4) circle (1.5pt) coordinate (m-4-4);
\draw[fill] (5,4) circle (1.5pt) coordinate (m-5-4);
\draw[fill] (6,4) circle (1.5pt) coordinate (m-6-4);

\draw[fill] (1,5) circle (1.5pt) coordinate (m-1-5);
\draw[fill] (2,5) circle (1.5pt) coordinate (m-2-5);
\draw[fill] (3,5) circle (1.5pt) coordinate (m-3-5);
\draw[fill] (4,5) circle (1.5pt) coordinate (m-4-5);
\draw[fill] (5,5) circle (1.5pt) coordinate (m-5-5);

\draw[fill] (1,6) circle (1.5pt) coordinate (m-1-6);
\draw[fill] (2,6) circle (1.5pt) coordinate (m-2-6);
\draw[fill] (3,6) circle (1.5pt) coordinate (m-3-6);
\draw[fill] (4,6) circle (1.5pt) coordinate (m-4-6);

\draw[fill] (1,7) circle (1.5pt) coordinate (m-1-7);
\draw[fill] (2,7) circle (1.5pt) coordinate (m-2-7);
\draw[fill] (3,7) circle (1.5pt) coordinate (m-3-7);
\draw[fill] (4,7) circle (1.5pt) coordinate (m-4-7);

\draw[fill] (1,8) circle (1.5pt) coordinate (m-1-8);
\draw[fill] (2,8) circle (1.5pt) coordinate (m-2-8);
\draw[fill] (3,8) circle (1.5pt) coordinate (m-3-8);

\draw[fill] (1,9) circle (1.5pt) coordinate (m-1-9);
\draw[fill] (2,9) circle (1.5pt) coordinate (m-2-9);

\draw[fill] (1,10) circle (1.5pt) coordinate (m-1-10);

\draw (m-0-0) -- (m-3-0) -- (m-6-4) -- (m-4-7) -- (m-1-10) -- (m-0-3) -- cycle;

\draw(5, 10) node{\large $\tau$};

\draw(-0.5, 2) node{\large $\rho_1$};
\draw(2, -0.5) node{\large $\rho_3$};
\draw(3.1, 8.5) node{\large $\rho_2$};

\draw [color=red, line width=1pt] (3,4)-- (3,5);
\draw [color=red, line width=1pt] (3,4)-- (4,4);
\draw [color=red, line width=1pt] (3,5)-- (4,4);

\draw(3.3, 4.3) node{\tiny$\delta$};

\draw [color=red, line width=1pt] (0,1)-- (3,4);
\draw [color=red, line width=1pt] (0,2)-- (3,5);
\draw [color=red, line width=1pt] (1,2)-- (1,3);
\draw [color=red, line width=1pt] (2,3)-- (2,4);

\draw(0.5, 2) node{$\sigma_{1}$};
\draw(1.5, 3) node{$\sigma_{2}$};
\draw(2.5, 4) node{$\sigma_{3}$};

\draw [color=red, line width=1pt] (1,0)-- (3,4);
\draw [color=red, line width=1pt] (2,0)-- (4,4);
\draw [color=red, line width=1pt] (2,2)-- (3,2);

\draw(1.7, 0.5) node{$\sigma''_{1}$};
\draw(2.7, 2.5) node{$\sigma''_{2}$};

\draw [color=red, line width=1pt] (3,5)-- (2,9);
\draw [color=red, line width=1pt] (4,4)-- (3,8);

\draw(3, 6.5) node{$\sigma'_{1}$};

\draw [color=red] (1,1)-- (0,0);
\draw [color=red] (1,1)-- (1,0);
\draw [color=red] (1,1)-- (0,1);
\draw [color=red] (1,1)-- (1,2);
\draw [color=red] (1,1)-- (2,2);
\draw [color=red] (1,1)-- (2,3);
\draw [color=red] (1,1)-- (3,4);

\draw [color=red] (5,4)-- (4,4);
\draw [color=red] (5,4)-- (3,2);
\draw [color=red] (5,4)-- (2,0);
\draw [color=red] (5,4)-- (3,0);
\draw [color=red] (5,4)-- (6,4);
\draw [color=red] (5,4)-- (3,8);
\draw [color=red] (5,4)-- (4,7);

\draw [color=red] (2,6)-- (0,2);
\draw [color=red] (2,6)-- (0,3);
\draw [color=red] (2,6)-- (2,4);
\draw [color=red] (2,6)-- (3,5);
\draw [color=red] (2,6)-- (2,9);
\draw [color=red] (2,6)-- (1,10);

\end{tikzpicture}
\captionsetup{labelformat=empty}
\caption{(2)}
\end{minipage}
\end{center}
\setcounter{figure}{9}
\caption{Subdivisions of general $\tau \in \bar{\mathscr{P}}^{[2]}$.} 
\label{singlocitypesgeneral}
\end{figure}

We lift the subdivisions of $\tau$ of Figure \ref{singlocitypesgeneral} to subdivisions of $\Delta_{\tau}$ in the same way that we lifted the subdivisions of $\tau$ of Figure \ref{refinedsubdivisiontypes} in Section \ref{resolveed} (with the only difference that $\sigma_i, \sigma'_j, \sigma''_k$ are no longer standard squares). Using the same notations as in Section \ref{resolveed}, the subdivisions of Figure \ref{singlocitypesgeneral} give rise to partial resolutions that are log smooth in the neighbourhoods of strata corresponding to $\Conv(v_i, \tau_{ij})$. We still need to resolve the singularities of $\Conv(\Delta_{\tau,1}, \sigma_i)$ in the case of Figure \ref{singlocitypesgeneral}(1) and of $\Conv(\langle v_1, v_2 \rangle, \sigma_i)$, $\Conv(\langle v_2, v_3 \rangle, \sigma'_j)$, $\Conv(\langle v_3, v_1 \rangle, \sigma''_k)$, and $\Conv(\Delta_{\tau,1}, \delta)$ in the case of Figure \ref{singlocitypesgeneral}(2). Again, $\Conv(\Delta_{\tau,1}, \delta)$ in the case of Figure \ref{singlocitypesgeneral}(2) is as in Figure \ref{locmodeltriangle}. The other convex hulls are $AGL(3,\mathbb{Z})$-equivalent to $$\Delta_{\sigma} \cong \Conv\left((\sigma \times \{e_0\}) \cup  \Delta_{\sigma,1} \right)$$ for $\sigma$ a parallelogram and $\Delta_{\sigma,1}$ the line segment connecting $(0,0)$ and $(1,0)$.\footnote{Note a different choice of representation for $\Delta_{\sigma,1}$ compared to Figure \ref{locmodelsquare}. Our choice here will make it more convenient to draw pictures.} Therefore, we need to find resolutions of such local models. The general idea is to rescale $\sigma$ to some $k\sigma$, subdivide $k\sigma$ into standard squares and triangles, and resolve the local models for the standard squares as we did in Section \ref{resolveed}. Of course, to resolve the local model $\tilde{\bar{X}}_{\tau,x}$, the rescaling will have to be done on the whole $\tau$, as we explain later.

Suppose that a parallelogram $\sigma$ is $AGL(2,\mathbb{Z})$-equivalent to the convex hull of $(0,0), (1,0), (p,q)$ and $(p+1, q)$ for some $p,q \in \mathbb{Z}_{>0}$ with $(p,q) = 1$ (the parallelograms of our subdivisions of $\tau$ satisfy this requirement). We will show how to resolve the local model defined by the cone over $\Delta_{k\sigma}$ for a suitable rescaling $k\sigma$ of $\sigma$. Namely, we require that $k \in \mathbb{Z}_{>0}$ is such that $k \geq p+3$. First, we subdivide $k\sigma$ as in Figure \ref{rescalingby3} (we ignore some interior lattice points of $\sigma$ and $k \sigma$ in the pictures).

\begin{figure}[t]
\begin{center}
\hspace{-2cm}
\begin{minipage}{.2\textwidth}
\begin{tikzpicture}[scale=0.7]

\draw[fill] (0,0) circle (1.5pt) coordinate (m-0-0);

\draw[fill] (2,0) circle (1.5pt) coordinate (m-2-0);

\draw[fill] (2,3) circle (1.5pt) coordinate (m-2-3);

\draw[fill] (4,3) circle (1.5pt) coordinate (m-4-3);

\draw (m-0-0) -- (m-2-0) -- (m-4-3) -- (m-2-3) -- cycle;

\draw(-0.4, 0.3) node{\tiny $(0,0)$};
\draw(-0.4+2, 0.3+3) node{\tiny $(p,q)$};
\draw(2.8, 0.3) node{\tiny $(1,0)$};
\draw(2.8 + 2, 0.3 + 3) node{\tiny $(p+1,q)$};

\end{tikzpicture}
\captionsetup{labelformat=empty}
\caption{$\sigma$}
\end{minipage}
\hspace{1cm}
\begin{minipage}{.6\textwidth}
\begin{tikzpicture}[scale=0.7]

\draw[fill] (0,0) circle (1.5pt) coordinate (m-0-0);
\draw[fill] (2,0) circle (1.5pt) coordinate (m-2-0);
\draw[fill] (4,0) circle (1.5pt) coordinate (m-4-0);
\draw[fill] (6,0) circle (1.5pt) coordinate (m-6-0);

\draw[fill] (2,3) circle (1.5pt) coordinate (m-2-3);
\draw[fill] (4,3) circle (1.5pt) coordinate (m-4-3);
\draw[fill] (6,3) circle (1.5pt) coordinate (m-6-3);
\draw[fill] (8,3) circle (1.5pt) coordinate (m-8-3);

\draw[fill] (4,6) circle (1.5pt) coordinate (m-4-6);
\draw[fill] (6,6) circle (1.5pt) coordinate (m-6-6);
\draw[fill] (8,6) circle (1.5pt) coordinate (m-8-6);
\draw[fill] (10,6) circle (1.5pt) coordinate (m-10-6);

\draw[fill] (6,9) circle (1.5pt) coordinate (m-6-9);
\draw[fill] (8,9) circle (1.5pt) coordinate (m-8-9);
\draw[fill] (10,9) circle (1.5pt) coordinate (m-10-9);
\draw[fill] (12,9) circle (1.5pt) coordinate (m-12-9);

\draw (m-0-0) -- (m-6-0) -- (m-12-9) -- (m-6-9) -- cycle;

\draw(-0.4, 0.3) node{\tiny $(0,0)$};
\draw(-0.4+2, 0.3+3) node{\tiny $(p,q)$};
\draw(2.6, 0.3) node{\tiny $(1,0)$};
\draw(2.8 + 2, 0.3 + 3) node{\tiny $(p+1,q)$};

\draw(6.8, 0.3) node{\tiny $(k,0)$};
\draw(7.2 + 2, 0.3 + 3) node{\tiny $(p+k,q)$};

\draw(-0.4+6, 0.3+9) node{\tiny $(kp,kq)$};
\draw(7.2 + 6, 0.3 + 9) node{\tiny $(k(p+1),kq)$};

\draw [color=red, line width=1pt] (2,3)-- (8,3);
\draw [color=red, line width=1pt] (4,6)-- (10,6);

\draw(4, 1.5) node{$k\sigma_1$};
\draw(4 + 2, 1.5 + 3) node{$\vdots$};
\draw(4 + 4, 1.5 + 6) node{$k\sigma_k$};

\end{tikzpicture}
\captionsetup{labelformat=empty}
\caption{$k\sigma$}
\end{minipage}
\end{center}
\setcounter{figure}{10}
\caption{Subdivision of $k\sigma$.} 
\label{rescalingby3}
\end{figure}

The subdivision of Figure \ref{rescalingby3} lifts to a subdivision of $\Delta_{k\sigma}$ into $\Conv(\Delta_{\sigma,1}, k\sigma_i)$, $ 1 \leq i \leq k$. We now need to subdivide each of the $k\sigma_i$ and lift the subdivisions to $\Delta_{k\sigma_i} \cong \Conv(\Delta_{\sigma,1}, k\sigma_i)$ for $ 1 \leq i \leq k$. We show how to subdivide $k\sigma_1$ in Figure \ref{subdivideparallelogram} (clearly, every $k\sigma_i$ is $AGL(2,\mathbb{Z})$-equivalent to $k\sigma_1$ and we will use the same subdivision for them). We may choose $\eta_1$ and $\eta_q$ as we wish (again, we require that every edge of $k\sigma_1$ subdivided by the subdivision is a union of three edges of the subdivision). The assumption that $k \geq p+3$ is necessary to obtain this subdivision.

\begin{figure}[H]
\begin{center}
\begin{tikzpicture}[scale=0.7]

\draw[fill] (0,0) circle (1.5pt) coordinate (m-0-0);
\draw[fill] (1,0) circle (1.5pt) coordinate (m-1-0);
\draw[fill] (2,0) circle (1.5pt) coordinate (m-2-0);
\draw[fill] (3,0) circle (1.5pt) coordinate (m-3-0);
\draw[fill] (4,0) circle (1.5pt) coordinate (m-4-0);
\draw[fill] (5,0) circle (1.5pt) coordinate (m-5-0);
\draw[fill] (6,0) circle (1.5pt) coordinate (m-6-0);

\draw[fill] (1,1) circle (1.5pt) coordinate (m-1-1);
\draw[fill] (2,1) circle (1.5pt) coordinate (m-2-1);
\draw[fill] (3,1) circle (1.5pt) coordinate (m-3-1);
\draw[fill] (4,1) circle (1.5pt) coordinate (m-4-1);
\draw[fill] (5,1) circle (1.5pt) coordinate (m-5-1);
\draw[fill] (6,1) circle (1.5pt) coordinate (m-6-1);

\draw[fill] (2,2) circle (1.5pt) coordinate (m-2-2);
\draw[fill] (3,2) circle (1.5pt) coordinate (m-3-2);
\draw[fill] (4,2) circle (1.5pt) coordinate (m-4-2);
\draw[fill] (5,2) circle (1.5pt) coordinate (m-5-2);
\draw[fill] (6,2) circle (1.5pt) coordinate (m-6-2);
\draw[fill] (7,2) circle (1.5pt) coordinate (m-7-2);

\draw[fill] (3,3) circle (1.5pt) coordinate (m-3-3);
\draw[fill] (4,3) circle (1.5pt) coordinate (m-4-3);
\draw[fill] (5,3) circle (1.5pt) coordinate (m-5-3);
\draw[fill] (6,3) circle (1.5pt) coordinate (m-6-3);
\draw[fill] (7,3) circle (1.5pt) coordinate (m-7-3);
\draw[fill] (8,3) circle (1.5pt) coordinate (m-8-3);

\draw[fill] (3,4) circle (1.5pt) coordinate (m-3-4);
\draw[fill] (4,4) circle (1.5pt) coordinate (m-4-4);
\draw[fill] (5,4) circle (1.5pt) coordinate (m-5-4);
\draw[fill] (6,4) circle (1.5pt) coordinate (m-6-4);
\draw[fill] (7,4) circle (1.5pt) coordinate (m-7-4);
\draw[fill] (8,4) circle (1.5pt) coordinate (m-8-4);
\draw[fill] (9,4) circle (1.5pt) coordinate (m-9-4);

\draw (m-0-0) -- (m-6-0) -- (m-9-4) -- (m-3-4) -- cycle;

\draw(-0.4, 0.3) node{\tiny $(0,0)$};
\draw(-0.4+3, 0.3+4) node{\tiny $(p,q)$};

\draw(6.9, 0.3) node{\tiny $(k,0)$};
\draw(7 + 3, 0.2 + 4) node{\tiny $(p+k,q)$};

\draw(4.5,1.6) node{\tiny $\vdots$};
\draw(4.5,2.6) node{\tiny $\vdots$};
\draw(5,3.45) node{\scriptsize $\eta_{q}$};

\draw [color=red, line width=1pt] (1,0)-- (4,1);
\draw [color=red, line width=1pt] (2,0)-- (5,1);
\draw [color=red, line width=1pt] (4,1)-- (4,3);
\draw [color=red, line width=1pt] (5,1)-- (5,3);

\draw [color=red, line width=1pt] (5,4)-- (4,3);
\draw [color=red, line width=1pt] (6,4)-- (5,3);

\draw [color=red, line width=1pt] (5,1)-- (4,1);
\draw [color=red, line width=1pt] (5,2)-- (4,2);
\draw [color=red, line width=1pt] (5,3)-- (4,3);

\draw [color=red] (0,0)-- (4,1);
\draw [color=red] (0,0)-- (4,2);
\draw [color=red] (0,0)-- (4,3);
\draw [color=red] (3,4)-- (4,3);

\draw [color=red] (6,0)-- (5,1);
\draw [color=red] (6,0)-- (5,2);
\draw [color=red] (6,0)-- (5,3);
\draw [color=red] (9,4)-- (5,3);

\draw[latex-,dashed] 
        (4,3) -- (-0.3, 3)
         node[fill=white]{\tiny $(p+1, q-1)$};

\draw[latex-,dashed] 
        (4,1) -- (-0.7, 1)
         node[fill=white]{\tiny $(p+1, 1)$};

\draw[latex-,dashed] 
        (5,3) -- (10, 3)
         node[fill=white]{\tiny $(p+2, q-1)$};

\draw[latex-,dashed] 
        (5,1) -- (9.7, 1)
         node[fill=white]{\tiny $(p+2, 1)$};

\draw[latex-,dashed] 
        (3.5,0.7) -- (3.5, -0.7)
         node[fill=white]{\scriptsize $\eta_{1}$};

\end{tikzpicture}
\end{center}
\setcounter{figure}{11}
\caption{Subdivision of $k\sigma_1$.} 
\label{subdivideparallelogram}
\end{figure}
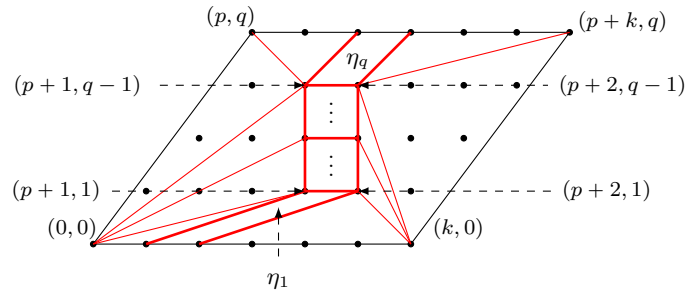

Note that all the $\eta_j, \ 1 \leq j \leq q$ in Figure \ref{subdivideparallelogram} are $AGL(2,\mathbb{Z})$-equivalent to standard squares. The subdivision of $k\sigma_1$ of Figure \ref{subdivideparallelogram} lifts to a subdivision of $\Delta_{k\sigma_1} \cong \Conv(\Delta_{\sigma,1}, k\sigma_i)$ as in the case of Figure \ref{refinedsubdivisiontypes}(1). Namely, if $v_1$ and $v_2$ are the vertices of $\Delta_{\sigma,1}$ corresponding to the lattice points $(0,0)$ and $(1,0)$, we subdivide $\Delta_{k\sigma_1}$ by $\Conv(v_1, \tau_{1i})$ for $\tau_{1i}$ the triangles \enquote{to the left of $\eta_j, \ 1 \leq j \leq q$} in Figure \ref{subdivideparallelogram}, $\Conv(v_2, \tau_{2i})$ for $\tau_{2i}$ the triangles \enquote{to the right of $\eta_j, \ 1 \leq j \leq q$} in Figure \ref{subdivideparallelogram}, and $\Conv(\Delta_{\sigma,1}, \eta_j), \ 1 \leq j \leq q$. 

Now each $\Conv(\Delta_{\sigma,1}, \eta_j), \ 1 \leq j \leq q$ can be subdivided as in Figure \ref{locmodelsquare}. Similarly to Section \ref{resolveed}, we can ensure this is compatible on the intersections of cells. This gives a resolution of the local model corresponding to $\Delta_{k\sigma_1}$. By performing similar subdivisions of $\Delta_{k\sigma_i}$ for $1 \leq i \leq k$ (again, ensuring that they are compatible on the intersections $\Delta_{k\sigma_i} \cap \Delta_{k\sigma_{i+1}}$ for $1 \leq i \leq k-1$), we obtain a resolution of the local model corresponding to $\Delta_{k\sigma}$. 

\begin{cons} \label{generalsubdivisionrecipe}
To obtain a subdivision of a rescaling $L \tau$ (for all $L \in \mathbb{Z}_{>0}$ sufficiently large) of a general $\tau \in \bar{\mathscr{P}}^{[2]}$ and lift it to a subdivision of $\Delta_{L\tau}$ one can proceed as follows:

\begin{enumerate}
    \item Obtain a subdivision of $\tau$ as in Figure \ref{singlocitypesgeneral} (possibly after doing a rescaling) and a lifting to a subdivision of $\Delta_{\tau}$ as explained after Figure \ref{singlocitypesgeneral}.
    \item Rescale the whole subdivision of $\tau$ to a subdivision of $L\tau$ so that all the $L\sigma_i$ (in the case of Figure \ref{singlocitypesgeneral}(1)) or $L\sigma_i,L\sigma'_j, L\sigma''_k$ (in the case of Figure \ref{singlocitypesgeneral}(2)) can be subdivided as in Figures \ref{rescalingby3} and \ref{subdivideparallelogram} (i.e. so that the condition $k \geq p+3$ is satisfied for all of them).
    \item Subdivide all the $L\sigma_i$ (or the $L\sigma_i,L\sigma'_j, L\sigma''_k$) as in Figure \ref{rescalingby3} into $L\sigma_{il}$ (or into $L\sigma_{il},L\sigma'_{jl}, L\sigma''_{kl}$) and do a refinement to a subdivision of $L\tau$ as follows. Suppose that we have a triangle $L\tau_{ij}$ with an edge $\rho$ and the subdivision of $L\sigma_i$ (or one of the $L\sigma_i,L\sigma'_j, L\sigma''_k$) subdivides $\rho$ into $\rho_1, \dots \rho_L$. Let $v_{\rho}$ be the vertex of $L\tau_{ij}$ that is not contained in $\rho$. Then we subdivide $\tau_{ij}$ into $\Conv(v_{\rho}, \rho_l), \ 1 \leq l \leq L$.  
    \begin{sloppypar}\item The subdivision of (3) lifts to $\Delta_{L\tau}$ by subdividing $\Conv(\Delta_{\tau,1}, L\sigma_i)$ (or $\Conv(\langle v_1, v_2 \rangle, L\sigma_i)$, $\Conv(\langle v_2, v_3 \rangle, L\sigma'_j)$, and $\Conv(\langle v_3, v_1 \rangle, L\sigma''_k)$) into  $\Conv(\Delta_{\tau,1}, L\sigma_{il})$ (or into  $\Conv(\langle v_1, v_2 \rangle, L\sigma_{il})$, $\Conv(\langle v_2, v_3 \rangle, L\sigma'_{jl})$, and  $\Conv(\langle v_3, v_1 \rangle, L\sigma''_{kl})$) and subdividing $\Conv (v_i, \tau_{ij})$ into $$\Conv (v_i, \Conv(v_{\rho}, \rho_l)), \quad 1 \leq l \leq L$$ using the same notations as in (3).\end{sloppypar}
    \item In the case of Figure \ref{singlocitypesgeneral}(2), note that $L\delta$ is easily decomposable and subdivide it as in Figure \ref{refinedsubdivisiontypes}(2). We denote by $\delta'$ the standard triangle of Figure \ref{refinedsubdivisiontypes}(2) contained in $L\delta$. This subdivision lifts to a subdivision of $\Conv(\Delta_{\tau,1}, L\delta)$ in the same way as we lift in Figure \ref{refinedsubdivisiontypes}(2) (or Figure \ref{singlocitypesgeneral}(2)). 
    \begin{sloppypar}\item Subdivide the $L\sigma_{il}$ (or $L\sigma_{il},L\sigma'_{jl}, L\sigma''_{kl}$) as in Figure \ref{subdivideparallelogram} so that the subdivisions are compatible with each other and with the subdivision of $L\delta$ (in the case of Figure \ref{singlocitypesgeneral}(2)). We can ensure compatibility by choosing the appropriate $\eta_1$ and $\eta_q$ in each case. These subdivisions lift to subdivisions of $\Conv(\Delta_{\tau,1}, L\sigma_{il})$ (or $\Conv(\langle v_1, v_2 \rangle, L\sigma_{il})$, $\Conv(\langle v_2, v_3 \rangle, L\sigma'_{jl})$, and $\Conv(\langle v_3, v_1 \rangle, L\sigma''_{kl})$) as explained after Figure \ref{subdivideparallelogram}. Compatibility ensures that these subdivisions (along with the subdivision of $\Conv(\Delta_{\tau,1}, L\delta)$ in the case of Figure \ref{singlocitypesgeneral}(2)) give a subdivision of $\Delta_{L\tau}$.\end{sloppypar}
    \item Subdivide the polyhedra of the subdivision of $\Delta_{L\tau}$ containing the standard squares (and containing $\delta'$ in the case of Figure \ref{singlocitypesgeneral}(2)) as in Figure \ref{locmodelsquare} (and Figure \ref{locmodeltriangle}). Again, this should be done compatibly, but similarly to Section \ref{resolveed}, this can be arranged. This gives the subdivision of $L \tau$ and a lifting to a subdivision of $\Delta_{L\tau}$.
\end{enumerate}
The subdivision of $\Delta_{L\tau}$ gives rise to a log smooth resolution of the corresponding local model $\tilde{\bar{X}}_{L\tau,x}$ given by the cone over $\Delta_{L\tau}$. Log smoothness follows from the fact that all the cells in the subdivisions of Figures \ref{locmodeltriangle} and \ref{locmodelsquare} are standard simplices, and all the other cells in the subdivision of $\Delta_{L\tau}$ are of the form $\Conv(v_i, \tau')$ for $\tau' \subseteq L\tau$ a cell in the subdivision of $L\tau$ (also see the footnote after Figure \ref{refinedsubdivisiontypes}).  
\end{cons}

\begin{remark} \label{generalsubdivisionreciperemark}
In the notations of Construction \ref{generalsubdivisionrecipe}, let $L\tau_{\rel}$ be the union of all the cells in the subdivision of $L\tau$ that subdivide the standard squares obtained after step (6) of Construction \ref{generalsubdivisionrecipe} (and $\delta'$ in the case of Figure \ref{singlocitypesgeneral}(2)). We note some flexibility of Construction \ref{generalsubdivisionrecipe}:
\begin{enumerate}
    \item In step (6) of Construction \ref{generalsubdivisionrecipe} one can allow subdivisions more general than those of Figure \ref{subdivideparallelogram}, i.e. one can vary not just $\eta_1$ and $\eta_q$ in Figure \ref{subdivideparallelogram} but also of all the $\eta_i, \ 2 \leq i \leq q-1$ (to the degree that the form of the subdivision in Figure \ref{subdivideparallelogram} allows). Similarly, in step (5) of Construction \ref{generalsubdivisionrecipe}, one can allow more general subdivisions of $L\delta$. This modifies the location of $L\tau_{\rel}$.
    \item Construction \ref{generalsubdivisionrecipe} gives an explicit subdivision of the connected components of $L\tau \setminus L\tau_{\rel}$ into triangular cells. However, one can allow any general subdivision of these connected components (into convex cells and inducing the given subdivision of $\partial(L\tau_{\rel})$), in particular, allow non-triangular cells.
\end{enumerate}
One can lift these more general subdivisions to subdivisions of $\Delta_{L\tau}$ in the same way as in Construction \ref{generalsubdivisionrecipe}. This still gives rise to log smooth resolutions of the corresponding local models $\tilde{\bar{X}}_{L\tau,x}$ given by the cones over $\Delta_{L\tau}$.  
\end{remark}

\subsubsection{Admissible resolutions in relative dimension \texorpdfstring{$n=3$}{n3}}

We are ready to define admissible resolutions of a special toric degeneration $\bar{\mathfrak{X}} \to \mathcal{S}$ of CY threefolds satisfying Assumption \ref{naturalcase}.

Suppose that we have obtained a resolution $\pi: \mathfrak{X} \to \bar{\mathfrak{X}}$ to a log smooth $\mathfrak{X} \to \mathcal{S}$.
First, we generalize the notation of Section \ref{tropicalresolutions} for the local models and their resolutions. For every $\tau \in \bar{\mathscr{P}}^{[k]}, \ 1 \leq k \leq 3$ and every $x \in \bar{X}_{\tau}$ we have an étale local model at $x$ given by a toric variety $\tilde{\bar{X}}_{\tau,x}$. This local model is defined by the cone over the $\Delta_{\tau}$ of \eqref{conv2} in the case that $x \in Z$ is a singular point and by the cone over $\tau$ in the case that $x \notin Z$ is non-singular. Since $\tilde{\bar{X}}_{\tau,x}$ is an étale local model for $\bar{\mathfrak{X}}$, we have a variety $\bar{\mathfrak{U}}_{\tau,x}$ equipped with étale maps $\bar{\mathfrak{U}}_{\tau,x} \to \bar{\mathfrak{X}}$ and $\bar{\mathfrak{U}}_{\tau,x} \to \tilde{\bar{X}}_{\sigma}$. We let $\mathfrak{U}_{\tau,x}:= \bar{\mathfrak{U}}_{\tau,x} \times_{\bar{\mathfrak{X}}} \mathfrak{X}$ be the basechange of $\bar{\mathfrak{U}}_{\tau,x} \to \bar{\mathfrak{X}}$ by $\pi: \mathfrak{X} \to \bar{\mathfrak{X}}$. For any morphism $\pi_{\tau,x}: \tilde{X}_{\tau,x} \to \tilde{\bar{X}}_{\tau,x}$, we let $\tilde{\mathfrak{U}}_{\tau,x, \pi_{\tau,x}}:= \bar{\mathfrak{U}}_{\tau,x} \times_{\tilde{\bar{X}}_{\tau,x}} \tilde{X}_{\tau,x}$ be the basechange of $\bar{\mathfrak{U}}_{\tau,x} \to \tilde{\bar{X}}_{\tau,x}$ by $\pi_{\tau,x}$. The following is a direct generalization of Definition \ref{toricresolutions}.

\begin{defn} 
We say that a resolution $\pi: \mathfrak{X} \to \bar{\mathfrak{X}}$ of a special toric degeneration $\bar{\mathfrak{X}} \to \mathcal{S}$ of CY threefolds satisfying Assumption \ref{naturalcase} to a log smooth degeneration $\mathfrak{X} \to \mathcal{S}$ is \emph{toric} if for every $\tau \in \bar{\mathscr{P}}^{[k]}, \ 1 \leq k \leq 3$ and every $x \in \bar{X}_{\tau}$ there exists a toric blowup $\pi_{\tau,x}: \tilde{X}_{\tau,x} \to \tilde{\bar{X}}_{\tau,x}$ such that $\mathfrak{U}_{\tau,x} \cong \tilde{\mathfrak{U}}_{\tau,x, \pi_{\tau,x}}$ and $\pi$ is trivial at every point of $\bar{\mathfrak{X}}$ not contained in a codimension $1$, $2$, or $3$ stratum of $\bar{\mathfrak{X}}_0$. 

We say that a toric resolution $\pi: \mathfrak{X} \to \bar{\mathfrak{X}}$ is \emph{integral} if for every $\tau \in \bar{\mathscr{P}}^{[k]}, \ 1 \leq k \leq 3$ and every $x \in \bar{X}_{\tau}$ the toric blowup $\pi_{\tau,x}: \tilde{X}_{\tau,x} \to \tilde{\bar{X}}_{\tau,x}$ induces an integral subdivision of:
 \begin{enumerate}
     \item If $x$ is a non-sigular point, the cone over $\tau$.
     \item If $x$ is a singular point, the cone over $\Delta_{\tau}$.
 \end{enumerate} 
Moreover, we require that the subdivisions of (1) are obtained by taking the cone over an integral subdivision of $\tau$, the subdivisions of (2) for $\Delta_{\tau}, \ \tau \in \bar{\mathscr{P}}^{[2]}$ are as in Construction \ref{generalsubdivisionrecipe} or Remark \ref{generalsubdivisionreciperemark} (with $\tau$ corresponding to the $L\tau$ of Construction \ref{generalsubdivisionrecipe}). 
We also require that for cells $\tau', \tau \in \bar{\mathscr{P}}$ with $\tau' \subseteq \tau$, the subdivision of $\tau'$ agrees with the one induced by the subdivision of $\tau$. Here and later, we assume that the cells of any subdivision don't self-intersect and that an intersection of any two cells is also a cell.
 
 We say that a toric and integral resolution is \emph{homogeneous} if for every $\tau \in \bar{\mathscr{P}}^{[k]}, \ 1 \leq k \leq 3$ and $x,y \in \bar{X}_{\tau}$ singular points, the subdivisions of (2) corresponding to the blowups $\pi_{\tau,x}: \tilde{X}_{\tau,x} \to \tilde{\bar{X}}_{\tau,x}$ and $\pi_{\tau,y}: \tilde{X}_{\tau,y} \to \tilde{\bar{X}}_{\tau,y}$ are the same. Note that the same is true for non-singular points by the definition of an integral resolution above. Further, if $\tau', \tau \in \bar{\mathscr{P}}$ are cells of dimension $1$ and $2$ respectively, $\tau' \subseteq \tau$, and $\bar{\mathfrak{X}} \to \mathcal{S}$ is not log smooth on $\bar{X}_{\tau'}$, we require that the subdivision of $\Delta_{\tau'}$ is the one induced by the subdivision of $\Delta_{\tau}$ by restricting to the face $F_{\tau'}$ of $\Delta_{\tau}$ (see Remark \ref{facelocalmodel}(1)). Note that the same is true for the subdivisions of the cones over $\tau'$ and $\tau$.
\end{defn}

We can now define strongly admissible resolutions as follows.

\begin{defn} \label{stronglyadmissiblegeneral}
We say that a resolution $\pi: \mathfrak{X} \to \bar{\mathfrak{X}}$ of a special toric degeneration $\bar{\mathfrak{X}} \to \mathcal{S}$ of CY threefolds satisfying Assumption \ref{naturalcase} is \emph{strongly admissible} if it is a toric, integral, and homogeneous resolution admitting a $\pi$-ample PA-generated divisor $D'= \sum_{v \in \mathscr{P}^{[0]}} a_v D_v$ such that $a_v \geq 0$ if $D_v$ is the strict transform of an irreducible component of $\bar{\mathfrak{X}}_0$ (the last condition can always be achieved by adding a positive multiple of $D$ to $D'$, see Corollary \ref{DisPLgen}).
\end{defn}

\begin{remarks}
\begin{enumerate}
    \item Definition \ref{stronglyadmissiblegeneral} is more general than the corresponding Definition \ref{stronglyadmissible} that we gave in the case of toric degenerations of K3-s. Indeed, in the case of special toric degenerations of K3-s, the results of Section \ref{scatint} imply that the extended intrinsic mirror $\check{\mathfrak{X}} \to \Spec \widehat{\kk [P]}_J$ can be defined via the collection $\mathfrak{D}^{J}$ of scattering diagrams of Construction \ref{extendcanonicalscat}  (so Conjecture \ref{conjecture} holds all toric, integral, and homogeneous resolutions of a special toric degeneration of K3-s, see Remark \ref{stronglyadmissibleweakening}(1)). It is reasonable to expect the same to be true in higher dimensions. If one wishes to define $\check{\mathfrak{X}} \to \Spec \widehat{\kk [P]}_J$ via \eqref{extendedmirror} (using Proposition \ref{P2general}), one can require a strongly admissible resolution to satisfy conditions (1) and (2) of the corresponding Definition \ref{stronglyadmissible}. Alternatively, one can restrict to strongly admissible resolutions with the divisor $\sum_{v \in \mathscr{P}^{[0]}} a_v D_v$ simple normal crossings and such that $a_v > 0$ for all $v \in \mathscr{P}^{[0]}$ (this clearly ensures that the assumptions of Proposition \ref{P2general} are satisfied). We show that such strongly admissible resolutions always exist (after a basechange) in Proposition \ref{basechangstronglyadmissiblesnc} below.  
    \item We still require the technical condition that $a_v \geq 0$ if $D_v$ is the strict transform of an irreducible component of $\bar{\mathfrak{X}}_0$ since we want to have a tropical way to construct strongly admissible resolutions, see Proposition \ref{onetoonetropicalgeneral} below (this is the natural analogue of Proposition \ref{stronglogsmoothresolution} and Corollary \ref{onetoonetropical} in the case of special toric degenerations of K3-s).
\end{enumerate}

\end{remarks}

Note that the resolutions of Section \ref{genericallysncgeneral} are strongly admissible.
We define admissible resolutions as in the case of special toric degenerations of K3-s (see Definition \ref{admissibledef}).

\begin{defn} \label{admissibledefgeneral}
We say that a resolution $\pi: \mathfrak{X}' \to \bar{\mathfrak{X}}$ of a special toric degeneration $\bar{\mathfrak{X}} \to \mathcal{S}$ of CY threefolds satisfying Assumption \ref{naturalcase} is \emph{admissible} if it factors as $\mathfrak{X}' \to \mathfrak{X} \to \bar{\mathfrak{X}}$ with $\mathfrak{X} \to \bar{\mathfrak{X}}$ strongly admissible and $\mathfrak{X}' \to \mathfrak{X}$ a logarithmic modification.
\end{defn}

\begin{remarks} \label{furtheradmissiblegenremarks}
\begin{enumerate}
    \item For special toric degenerations of relative dimension $n \geq 4$ satisfying Assumptions \ref{naturalcase} and \ref{facelocalmodelass}, the definition of a strongly admissible resolution should be similar to Definition \ref{stronglyadmissiblegeneral} provided that we have a suitable analogue of Construction \ref{generalsubdivisionrecipe} and Remark \ref{generalsubdivisionreciperemark}. The analogue should give rise to, for every $\tau \in \bar{\mathscr{P}}^{[n-1]}$ and every singular point $x \in \bar{X}_{\tau}$, a toric blowup $\pi_{\tau,x}: \tilde{X}_{\tau,x} \to \tilde{\bar{X}}_{\tau,x}$ resolving the local model (that does not depend on the choice of a singular point $x \in \bar{X}_{\tau}$). The definition of an admissible resolution is then as in Definition \ref{admissibledefgeneral}. 
    \item If a special toric degeneration does not satisfy Assumption \ref{facelocalmodelass}, one only needs to modify the compatibility condition in the definition of a homogeneous resolution. For special toric degenerations not satisfying Assumption \ref{naturalcase}, we still expect the definition of a strongly admissible resolution to follow the same lines. The condition that the resolution is toric should be replaced by a condition that there exist certain blowups of the local models of Remark \ref{specialass}(3) (provided by condition (3) of Definition \ref{special} of a special toric degeneration) that satisfy a similar compatibility condition with the resolution. The local models at the singular points $x, x' \in Z$ with $x'$ the specialization of $x$ should also be compatible (corresponding to the compatibility requirements in the integral and homogeneous conditions of Definition \ref{stronglyadmissiblegeneral}). These requirements should allow one to glue the blowups of the local models to a log smooth resolution, generalizing Proposition \ref{stronglogsmoothresolution}. The definition of an admissible resolution is then as in Definition \ref{admissibledefgeneral}.
\end{enumerate}
\end{remarks}

We have an analogue of Proposition \ref{stronglogsmoothresolution} and Corollary \ref{onetoonetropical} that allows us to construct strongly admissible resolutions of $\bar{\mathfrak{X}} \to \mathcal{S}$ tropically.

\begin{fact} \label{onetoonetropicalgeneral}
There is a one-to-one correspondence $$\left\{(B, \mathscr{P}), \alpha \in PA(B)\right\} \Longleftrightarrow \left\{ \pi:\mathfrak{X} \to \bar{\mathfrak{X}}, D'\right\}$$
where:
\begin{enumerate}
    \item $\bar{\mathfrak{X}} \to \mathcal{S}$ is a special toric degeneration of CY threefolds satisfying Assumption \ref{naturalcase} with dual intersection complex $\left(\bar{B}, \bar{\mathscr{P}}\right)$.
    \item $(B, \mathscr{P})$ is an integral subdivision of $\left(\bar{B}, \bar{\mathscr{P}}\right)$ such that the induced subdivision of every $\tau \in \bar{\mathscr{P}}^{[2]}$ is as in Construction \ref{generalsubdivisionrecipe} or Remark \ref{generalsubdivisionreciperemark} (or arbitrary if $\bar{\mathfrak{X}} \to \mathcal{S}$ is log smooth on $\bar{X}_{\tau}$) and for every $\rho \in \bar{\mathscr{P}}^{[1]}$ and $\tau, \tau' \in \bar{\mathscr{P}}^{[2]}$ with $\rho \subseteq \tau, \tau'$, the standard squares of the subdivisions of $\tau$ and $\tau'$ adjacent to $\rho$ intersect in an edge $\rho' \subseteq \rho, \ \rho' \in \mathscr{P}^{[1]}$. We require that the cells of the subdivision don't self-intersect and that an intersection of any two cells is also a cell (i.e. $(B, \mathscr{P})$ satisfies Assumption \ref{manifoldass}).
    \item $\pi:\mathfrak{X} \to \bar{\mathfrak{X}}$ is a strongly admissible resolution to a log smooth and minimal log CY degeneration $\mathfrak{X} \to \mathcal{S}$ with a $\pi$-ample PA-generated divisor $D' = \sum_{v \in \mathscr{P}^{[0]}} a_v D_v$ such that $a_v \geq 0$ if $D_v$ is the strict transform of an irreducible component of $\bar{\mathfrak{X}}_0$.
    \item $\alpha \in PA(B)$ is such that:
    \begin{enumerate}
        \item $\alpha(v) \geq 0$ for all $v \in \bar{\mathscr{P}}^{[0]}$.\footnote{This is a technical condition that is necessary to ensure that the ideal sheaves defining the blowups in the étale local models glue to an ideal sheaf on $\bar{\mathfrak{X}}$.}
        \item The restriction $\alpha|_{\sigma}$ of $\alpha$ to any $\sigma \in \bar{\mathscr{P}}^{\max}$ is strictly convex (on the induced subdivision of $\sigma$).
        \item The restriction $\alpha|_{\rho}$ of $\alpha$ to any $\rho \in \bar{\mathscr{P}}^{[1]}$ such that $\bar{\mathfrak{X}} \to \mathcal{S}$ is not log smooth on $\bar{X}_{\rho}$ has a well-defined $\min(\alpha|_{\rho}) \subseteq \rho(\mathbb{Z})$ (see Notation \ref{minnotation}). Further, $\min(\alpha|_{\rho})$ is one of the vertices of the edge $\rho' \in \mathscr{P}^{[1]}$ of a standard square of Construction \ref{generalsubdivisionrecipe} or Remark \ref{generalsubdivisionreciperemark} subdividing $\rho$.
        \item Consider an edge $\rho' \in \mathscr{P}^{[1]}$ of a standard square in the subdivision of a cell $\tau \in \bar{\mathscr{P}}^{[2]}$ of Construction \ref{generalsubdivisionrecipe} or Remark \ref{generalsubdivisionreciperemark}. Suppose further that $\rho'$ is parallel to the edge (or edges) of $\tau$ connected to $\rho$ by standard squares. Then $\alpha$ has different values at the two  vertices of $\rho'$ (so the restriction $\alpha|_{\rho'}$ of $\alpha$ to $\rho'$ has a well-defined $\min(\alpha|_{\rho'}) \subseteq \rho'(\mathbb{Z})$)
        \item Suppose that $\tau \in \bar{\mathscr{P}}^{[2]}$ is as in Figure \ref{singlocitypesgeneral}(2) and let $\delta'$ be the standard triangle of Construction \ref{generalsubdivisionrecipe} or Remark \ref{generalsubdivisionreciperemark}. Let $v_1, v_2, v_3$ be the vertices of $\delta'$. Then $\alpha(v_1) < \alpha(v_2) <\alpha(v_3)$ for some rearrangement of $v_1, v_2, v_3$. 
    \end{enumerate} 
\end{enumerate}
Under the correspondence, we have $\divisor(\alpha) = D'$ via \eqref{divmap}. Moreover, the dual intersection complex of $\mathfrak{X} \to \mathcal{S}$ is isomorphic to $(B, \mathscr{P})$ as a polyhedral manifold. Identifying the dual intersection complex of $\mathfrak{X} \to \mathcal{S}$ with $(B, \mathscr{P})$ under the isomorphism, the affine structure on $(B, \mathscr{P})$ extends to the complement of the singular locus  
\begin{multline} \label{extsinglocusgeneral}
\Delta = 
\left\{\langle \min(\alpha|_{\rho'}), \min(\alpha|_{\rho''})  \rangle\,\Bigg|\,
\substack{\hbox{$\rho', \rho'' \in \mathscr{P}^{[1]}$ two parallel edges as in (4)(c)}\\\hbox{ or (4)(d) of the same standard square in the}\\\hbox{subdivision of a cell $\tau \in \bar{\mathscr{P}}^{[2]}$}}\right\} \bigcup \\  
\bigcup \left\{\langle v_1,v_2  \rangle\,\Bigg|\,
\substack{\hbox{$\tau \in \bar{\mathscr{P}}^{[2]}$ as in Figure \ref{singlocitypesgeneral}(2),}\\\hbox{$\delta'$ the standard triangle of $\tau$,}\\\hbox{and $v_1, v_2, v_3$ as in (4)(e)}}  \right\}.  \ \ \ \ \ \ \ \ \ \ \, 
\end{multline}
\end{fact}

\begin{proof}
We argue similarly to the proofs of Proposition \ref{stronglogsmoothresolution} and Corollary \ref{onetoonetropical}. 

The subdivision $(B, \mathscr{P})$ and the PA-function $\alpha \in PA(B)$ define the blowups in the local models and the corresponding ideal sheaves. The analysis of the local models in codimension $3$ is similar to Step 1 in the proof of Proposition \ref{stronglogsmoothresolution} using conditions (4)(a) and (4)(b). For every $\sigma \in \bar{\mathscr{P}}^{\max}$ we obtain a toric blowup $\pi_{\sigma, \bar{X}_{\sigma}}: \tilde{X}_{\sigma, \bar{X}_{\sigma}} \to \tilde{\bar{X}}_{\sigma, \bar{X}_{\sigma}}$ of the local model $\tilde{\bar{X}}_{\sigma, \bar{X}_{\sigma}}$ defined by an ideal sheaf $\mathcal{I}_{\sigma, \bar{X}_{\sigma}}$ supported on codimension $1$, $2$, and $3$ toric strata of $\partial \tilde{\bar{X}}_{\sigma, \bar{X}_{\sigma}}$. The analysis of the local models in codimension $1$ is similar to Step 2 in the proof of Proposition \ref{stronglogsmoothresolution} using conditions (4)(a) and (4)(c). For every $\rho \in \bar{\mathscr{P}}^{[1]}$ and $x \in \bar{X}_{\rho}$ we obtain a toric blowup $\pi_{\rho, x}: \tilde{X}_{\rho,x} \to \tilde{\bar{X}}_{\rho,x}$ of the local model $\tilde{\bar{X}}_{\rho,x}$ defined by an ideal sheaf $\mathcal{I}_{\rho,x}$ supported on the toric stratum corresponding to $\rho$. 

The analysis of the local models in codimension $2$ is similar to Step 2 in the proof of Proposition \ref{stronglogsmoothresolution} but is more involved. Let $\tau \in \bar{\mathscr{P}}^{[2]}$ and $x \in \bar{X}_{\tau}$. 

Assume that $x$ is a singular point. Then conditions (4)(c), (4)(d), and (4)(e) determine the subdivision of $\Delta_{\tau}$. Indeed, the subdivision of $\Delta_{\tau}$ away from the polyhedra containing the cells of $\tau_{\rel}$ (using the notation of Remark \ref{generalsubdivisionreciperemark}) is already determined by Construction \ref{generalsubdivisionrecipe} or Remark \ref{generalsubdivisionreciperemark}. So we need to explain how to subdivide $\Conv(\rho, \eta)$ for $\eta$ a standard square and $\rho \subseteq \Delta_{\tau,1}$ the corresponding edge of $\Delta_{\tau,1}$ and how to subdivide $\Conv(\Delta_{\tau,1}, \delta')$ for the standard triangles $\Delta_{\tau,1}, \delta'$ in the case of Figure \ref{singlocitypesgeneral}(2). We use the unique subdivision of $\Conv(\rho, \eta)$ as in Figure \ref{locmodelsquare} such that its restriction to $\eta$ is the given subdivision of $\eta$ and the corresponding singular locus in Figure \ref{singularlocuslsquare} is given by $\langle \min(\alpha|_{\rho'}), \min(\alpha|_{\rho''}) \rangle$ for $\rho', \rho''$ as in \eqref{extsinglocusgeneral}. We use the unique subdivision of $\Conv(\Delta_{\tau,1}, \delta')$ as in Figure \ref{locmodeltriangle} such that the corresponding singular locus in Figure \ref{singularlocusltriangle} is given by $\langle v_1, v_2 \rangle$ as in \eqref{extsinglocusgeneral}. The fact that the involved PA-functions are restrictions of the globally defined $\alpha \in PA(B)$ implies that the defined subdivisions agree on the common faces and define a subdivision of $\Delta_{\tau}$ as in Construction \ref{generalsubdivisionrecipe} or Remark \ref{generalsubdivisionreciperemark}. The subdivision gives rise to a toric blowup $\pi_{\tau, x}: \tilde{X}_{\tau,x} \to \tilde{\bar{X}}_{\tau,x}$ of the local model $\tilde{\bar{X}}_{\tau,x}$. 

Let $\alpha|_{\tau,x}^{0}$ be the extension by $0$ of $\alpha|_{\tau}$ to $\Delta_{\tau}$. Here by an extension by $0$, we mean that $\alpha|_{\tau,x}^{0}$ is the PA-function that takes the same values as $\alpha|_{\tau}$ at the points corresponding to the integer points of $\tau$ and takes value $0$ at the points corresponding to the integer points of $\Delta_{\tau,1}$. Now conditions (4)(c), (4)(d), (4)(e) are equivalent to strict convexity of $\alpha|_{\tau,x}^{0}$ with respect to the described subdivision of $\Delta_{\tau}$. This follows by considering the subdivisions of Figure \ref{locmodelsquare} and Figure \ref{locmodeltriangle}. As in Step 2 in the proof of Proposition \ref{stronglogsmoothresolution}, $\alpha|_{\tau,x}^{0}$ gives rise to an ideal sheaf $\mathcal{I}_{\tau,x}$ on $\tilde{\bar{X}}_{\tau,x}$ inducing the subdivision of $\Delta_{\tau}$. Moreover, by an argument similar to the one in Step 2 in the proof of Proposition \ref{stronglogsmoothresolution}, the fact that $\alpha|_{\tau,x}^{0}$ is the extension by $0$ implies that $\mathcal{I}_{\tau,x}$ is supported on the union of the toric stratum corresponding to $\tau$ and the codimension $1$ toric strata containing the stratum corresponding to $\tau$. 

If $x$ is a non-singular point, the subdivision of the cone over $\tau$ is determined by the subdivision of $\tau$. This gives rise to a toric blowup $\pi_{\tau, x}: \tilde{X}_{\tau,x} \to \tilde{\bar{X}}_{\tau,x}$ of the local model $\tilde{\bar{X}}_{\tau,x}$. The corresponding ideal sheaf $\mathcal{I}_{\tau,x}$ is constructed as in the non-singular case, and the statements on the support still hold.

Now we can glue the ideal sheaves $\mathcal{I}_{\tau, x}$ for every $\tau \in \bar{\mathscr{P}}^{[k]}, \ 1 \leq k \leq 3$ and every $x \in \bar{X}_{\tau}$ by an argument as in Step 3 in the proof of Proposition \ref{stronglogsmoothresolution}, setting $\bar{\mathfrak{U}}$ to be the complement of codimension $1$, $2$, and $3$ strata of $\bar{\mathfrak{X}}_0$. The fact that we can define a descent datum for the ideal sheaves $\mathcal{I}_{\bar{\mathfrak{U}}_{\tau, x}}$ and $\mathcal{I}_{\bar{\mathfrak{U}}}$ (defined as in Step 3 in the proof of Proposition \ref{stronglogsmoothresolution}) follows from:
\begin{enumerate}
    \item The correspondence between the ideals in the local toric models and the support functions.
    \item The fact that $(\alpha|_{\sigma})|_{\tau} = (\alpha|_{\tau, x}^{0})|_{\tau}$ for every $\tau \subseteq \sigma$ with $\tau \in \bar{\mathscr{P}}^{[1]} \cup \bar{\mathscr{P}}^{[2]}, \ 
    \sigma \in \bar{\mathscr{P}}^{\max}$ and every $x \in \bar{X}_{\tau}$.
    \item The fact that $(\alpha|^0_{\tau, x})|_{F_{\tau'}} = \alpha|_{\tau', x'}^{0}$ for every $\tau', \tau \in \bar{\mathscr{P}}$ cells of dimension $1$ and $2$ respectively with $\tau \subseteq \tau'$ (see Remark \ref{facelocalmodel}(1)), and $x$ and $x'$ singular points with minimal strata $\bar{X}_{\tau'}$ and $\bar{X}_{\tau}$ respectively.
    \item The vanishing statements on $\mathcal{I}_{\sigma, \bar{X}_{\sigma}}$ for $\sigma \in \bar{\mathscr{P}}^{\max}$,  $\mathcal{I}_{\tau,x}$ for $\tau \in \bar{\mathscr{P}}^{[2]}$ and $\mathcal{I}_{\rho,x}$ for $\rho \in \bar{\mathscr{P}}^{[1]}$.
\end{enumerate}
Here claims (2) and (3) hold since both $\alpha|_{\rho}$ and $\alpha|_{\sigma}$ arise as the restrictions of a globally defined $\alpha \in PA(B)$. We obtain an ideal sheaf $\mathcal{I}$ on $\bar{\mathfrak{X}}$ and let $\pi: \mathfrak{X} \to \bar{\mathfrak{X}}$ be the blowup of $\mathcal{I}$.

The fact that $\pi:\mathfrak{X} \to \bar{\mathfrak{X}}$ is toric and integral and $\mathfrak{X} \to \mathcal{S}$ is minimal log CY follows as in Step 4 of the proof of Proposition \ref{stronglogsmoothresolution}. $\pi:\mathfrak{X} \to \bar{\mathfrak{X}}$ is homogeneous by construction and condition (3) in the previous paragraph. As we mentioned above, conditions (4)(c), (4)(d), (4)(e) are equivalent to strict convexity of $\alpha|_{\tau,x}^{0}$ with respect to the described subdivision of $\Delta_{\tau}$. As in Step 3 of the proof of Proposition \ref{ampledivisor}, condition (4)(c) is equivalent to strict convexity of the corresponding $\alpha|_{\rho}^{0}$. Then the fact that $\divisor(\alpha)$ is ample if and only if $\alpha$ satisfies conditions (4)(b), (4)(c), (4)(d), (4)(e) follows by an argument as in Step 3 of the proof of Proposition \ref{ampledivisor}. The structure of the singular locus follows from the described subdivisions of $\Delta_{\tau}, \ \tau \in \bar{\mathscr{P}}^{[2]}$ and the analysis of the local models of Figures \ref{singularlocusltriangle} and \ref{singularlocuslsquare}. The rest follows as in the proof of Corollary \ref{onetoonetropical}.
\end{proof}

We can now use Proposition \ref{onetoonetropicalgeneral} to deduce an existence result for strongly admissible resolutions. 

\begin{theorem} \label{basechangstronglyadmissible}
Let $\bar{\mathfrak{X}} \to \mathcal{S}$ be a special toric degeneration of CY threefolds satisfying Assumption \ref{naturalcase}. Then there exists an $L_{\bar{\mathfrak{X}}} \in \mathbb{Z}_{>0}$ such that for every $L \in \mathbb{Z}_{>0}$ with $L \geq L_{\bar{\mathfrak{X}}}$ the basechange $\bar{\mathfrak{X}}' \to \mathcal{S}$ of $\bar{\mathfrak{X}} \to \mathcal{S}$ by $R \to R, \ t \mapsto t^L$ (where $t$ is the uniformizer of $R$) admits a strongly admissible resolution $\pi: \mathfrak{X}' \to \bar{\mathfrak{X}}'$.
\end{theorem}

\begin{proof}
Let $\left(\bar{B}, \bar{\mathscr{P}}\right)$ be the dual intersection complex of $\bar{\mathfrak{X}} \to \mathcal{S}$. For every $\tau \in \bar{\mathscr{P}}^{[2]}$ such that $\bar{\mathfrak{X}} \to \mathcal{S}$ is not log smooth on $\bar{X}_{\tau}$, we let $L_{\tau} \in \mathbb{Z}_{>0}$ be such that for every $L \in \mathbb{Z}_{>0}$ with $L \geq L_{\tau}$, the rescaling $L\tau$ of $\tau$ admits a subdivision as in Construction \ref{generalsubdivisionrecipe}. For every $\sigma \in \bar{\mathscr{P}}^{\max}$ we let $L_{\sigma} \in \mathbb{Z}_{>0}$ be such that for every $L \in \mathbb{Z}_{>0}$ with $L \geq L_{\sigma}$, the rescaling $L\sigma$ of $\sigma$ has an interior point.\footnote{The fact that $\sigma$ is an integral polytope of dimension $3$ and the classification of lattice polytops of dimension $3$ without interior lattice points (see \cite[Theorem 2.2]{AWW} and \cite[Theorem 1]{AKW}) implies that we can choose $L_{\sigma} \leq 4$.} We let $$L_{\bar{\mathfrak{X}}}:=\max\left\{L_{\tau}, L_{\sigma} \ | \ \tau \in \bar{\mathscr{P}}^{[2]}, \ \sigma \in \bar{\mathscr{P}}^{\max} \right\}.$$
Fix an $L \in \mathbb{Z}_{>0}$ with $L \geq L_{\bar{\mathfrak{X}}}$ and let $\bar{\mathfrak{X}}' \to \mathcal{S}$ be the basechange of $\bar{\mathfrak{X}} \to \mathcal{S}$ by $R \to R, \ t \mapsto t^L$.

Let $\left(\bar{B}', \bar{\mathscr{P}}'\right)$ be the dual intersection complex of $\bar{\mathfrak{X}}' \to \mathcal{S}$. By the correspondence of Proposition \ref{onetoonetropicalgeneral}, it is enough to find an integral subdivision $(B', \mathscr{P}')$ of $\left(\bar{B}', \bar{\mathscr{P}}'\right)$ satisfying condition (2) of Proposition \ref{onetoonetropicalgeneral} and an $\alpha' \in PA(B')$ satisfying condition (4) of Proposition \ref{onetoonetropicalgeneral}. We will obtain $(B', \mathscr{P}')$ and $\alpha' \in PA(B')$ by arguing similarly to the proof of Proposition \ref{weakadmissibleresolutionresult}.

We define the integral subdivision $(B', \mathscr{P}')$ of $\left(\bar{B}', \bar{\mathscr{P}}'\right)$ as follows. First, we subdivide every $\tau \in \bar{\mathscr{P}}'^{[2]}$ as in Construction \ref{generalsubdivisionrecipe} (or trivially if $\bar{\mathfrak{X}}' \to \mathcal{S}$ is log smooth on $\bar{X}_{\tau}$) ensuring that for every $\rho \in \bar{\mathscr{P}}'^{[1]}$ and $\tau, \tau' \in \bar{\mathscr{P}}'^{[2]}$ with $\rho \subseteq \tau, \tau'$, the standard squares of the subdivisions of $\tau$ and $\tau'$ adjacent to $\rho$ intersect in an edge $\rho' \subseteq \rho$ of the subdivision. This is possible since in Construction \ref{generalsubdivisionrecipe}(1) we can choose any length $1$ segments on the edges for the subdivision of Figure \ref{singlocitypesgeneral} and in Construction \ref{generalsubdivisionrecipe}(6) we can choose appropriate $\eta_1$ and $\eta_q$ in the subdivisions of Figure \ref{subdivideparallelogram}. This gives an integral subdivision $(\partial, \mathscr{P}'_{\partial})$ of the codimension $1$ skeleton $(\partial, \bar{\mathscr{P}}'_{\partial})$ of $\bar{\mathscr{P}}'$ (here $\bar{\mathscr{P}}'_{\partial} = \bar{\mathscr{P}}'^{[0]} \cup \bar{\mathscr{P}}'^{[1]} \cup \bar{\mathscr{P}}'^{[2]}$). To lift it to an integral subdivision $(B', \mathscr{P}')$ of $\left(\bar{B}', \bar{\mathscr{P}}'\right)$, choose an interior integral point $x_{\sigma}$ for every $\sigma \in \bar{\mathscr{P}}'^{\max}$ and subdivide $\sigma$ into $\Conv(x_{\sigma}, \tilde{\tau})$ for every $\tilde{\tau} \subseteq \tau \subseteq \sigma$ with $\tilde{\tau} \in \mathscr{P}_{\partial}'^{[2]}$ and $\tau \in \bar{\mathscr{P}}'^{[2]}$ (this is similar to Step 2 in the proof of Proposition \ref{weakadmissibleresolutionresult}). This defines an integral subdivision $(B', \mathscr{P}')$ of $\left(\bar{B}', \bar{\mathscr{P}}'\right)$ which satisfies condition (2) of Proposition \ref{onetoonetropicalgeneral} by construction.

We now want to find an $\alpha' \in PA(B')$ satisfying condition (4) of Proposition \ref{onetoonetropicalgeneral}.
By an argument as in Steps 2 and 3 in the proof of Proposition \ref{weakadmissibleresolutionresult}, it is enough to define a PA-function $\alpha'_{\partial}$ on $(\partial, \mathscr{P}'_{\partial})$ that satisfies the natural analogues of conditions (4)(c), (4)(d), (4)(e) of Proposition \ref{onetoonetropicalgeneral}. Let $(\partial, \mathscr{P}''_{\partial})$ be a partial subdivision of $(\partial, \bar{\mathscr{P}}'_{\partial})$ such that every cell $\tau \in \bar{\mathscr{P}}'^{[2]}_{\partial}$ is subdivided as after Construction \ref{generalsubdivisionrecipe}(3). For every $\tau \in \bar{\mathscr{P}}'^{[2]}_{\partial}$, let $\tau'_{\rel}$ be the union of parallelograms (and a rescaling $L\delta$ of the standard triangle $\delta$ in the case of Figure \ref{singlocitypesgeneral}(2)) in the subdivision of $\tau$ after Construction \ref{generalsubdivisionrecipe}(3). Note that $\tau_{\rel} \subseteq \tau'_{\rel}$ using the notations of Remark \ref{generalsubdivisionreciperemark}. By an argument as in Step 4 in the proof of Proposition \ref{ampledivisor}, it follows that it is enough to:
\begin{enumerate}
    \item Define a PA-function $\alpha''_{\partial}$ on $(\partial, \mathscr{P}''_{\partial})$ such that $\alpha''_{\partial}(v) = 0$ for $v \in \bar{\mathscr{P}}_{\partial}'^{[0]}$,  and the restriction $(\alpha''_{\partial})|_{\tau}$ of $\alpha''_{\partial}$ to any $\tau \in \bar{\mathscr{P}}_{\partial}'^{[2]}$ is strictly convex (on the induced subdivision of $\tau$). 
    \item For every cell $\tau \in \bar{\mathscr{P}}_{\partial}'^{[2]}$, define a strictly convex PA-function $\alpha_{\tau'_{\rel}}$ on the subdivision of $\tau'_{\rel}$ induced by Steps (4-7) of Construction \ref{generalsubdivisionrecipe} such that:
    \begin{enumerate}
        \item $\alpha_{\tau'_{\rel}}$ satisfies the natural analogues of conditions (4)(d) and (4)(e) of Proposition \ref{onetoonetropicalgeneral}.
        \item $\alpha_{\tau'_{\rel}}(v) = 0$ for any $v \in \mathscr{P}'^{[0]}$ such that $v \in \partial \tau'_{\rel}$.
        \item $\alpha_{\tau'_{\rel}}$ has $\alpha_{\tau'_{\rel}}(v_{\rho}) = -2, \ \alpha_{\tau'_{\rel}}(v'_{\rho}) = -3$ for any $\rho \in \bar{\mathscr{P}}_{\partial}'^{[1]}$ with $\rho \cap \tau_{\rel} \neq \varnothing$ and any chosen ordering of vertices $v_{\rho}, v'_{\rho}$ of the intersection of $\rho$ with the standard square of $\tau_{\rel}$ (this condition ensures that the PA-functions $\alpha_{\tau'_{\rel}}$ can be patched to a PA-function on $\cup_{\tau \in \bar{\mathscr{P}}_{\partial}'^{[2]}} \tau'_{\rel}$).
    \end{enumerate}
\end{enumerate}
It is not difficult to check explicitly that one can construct PA-functions $\alpha''_{\partial}$ and $\alpha_{\tau'_{\rel}}$ as in (1) and (2) respectively. By an argument as in Step 4 in the proof of Proposition \ref{weakadmissibleresolutionresult}, this amounts to checking that certain systems of linear inequalities have solutions. In the case of $\alpha''_{\partial}$, the system is given by inequalities between the values of $\alpha''_{\partial}$ at the vertices of the parallelograms in the subdivision of every $\tau \in \bar{\mathscr{P}}_{\partial}'^{[2]}$ after Construction \ref{generalsubdivisionrecipe}(3) (note that this includes the vertices of the rescaling $L\delta$ of the standard triangle $\delta$ in the case of Figure \ref{singlocitypesgeneral}(2)). In the case of $\alpha_{\tau'_{\rel}}$, the system is given by inequalities between the values of $\alpha_{\tau'_{\rel}}$ at the vertices of the standard squares in the induced subdivision of $\tau'_{\rel}$ after Steps (4-7) of Construction \ref{generalsubdivisionrecipe} (note that this includes the vertices of the standard triangle $\delta'$ of Construction \ref{generalsubdivisionrecipe}(5) in the case of Figure \ref{singlocitypesgeneral}(2)). 
\end{proof}

One can also obtain admissible resolutions with $D'$ simple normal crossings. We can partially generalize Corollary \ref{furtherresolutionofadmissible} as follows.

\begin{cor} \label{sncafterfurtherblowupgeneral}
Let $\bar{\mathfrak{X}} \to \mathcal{S}$ be a special toric degeneration of CY threefolds satisfying Assumption \ref{naturalcase} and admitting a strongly admissible resolution. For any strongly admissible resolution $\pi: \mathfrak{X} \to \bar{\mathfrak{X}}$, there exists a further logarithmic modification $\pi': \tilde{\mathfrak{X}} \to \mathfrak{X}$ with the combined $\tilde{\mathfrak{X}} \to \mathfrak{X} \to \bar{\mathfrak{X}}$ an admissible resolution with $\tilde{D}:=(\pi')^{-1}D$ simple normal crossings.
\end{cor}

\begin{proof}
Let $(B, \mathscr{P})$ be the dual intersection complex of $\mathfrak{X} \to \mathcal{S}$. A classical result of Knudsen, Mumford, and Waterman \cite[Theorem 4.1]{KKMS} states that every polyhedral complex $\mathscr{C}$ admits a rescaling $L\mathscr{C}$ for some $L \in \mathbb{Z}_{>0}$ such that $L\mathscr{C}$ admits a unimodular triangulation, i.e. a subdivision with all the maximal cells standard simplices. Therefore, for some $L \in \mathbb{Z}_{>0}$, $(B, \mathscr{P})$ admits a subdivision $(\tilde{B}, \tilde{\mathscr{P}})$ such that all the cells $\tilde{\sigma} \in \tilde{\mathscr{P}}^{\max}$ are rescalings of the standard simplex by $\frac{1}{L}$. Recall that there is a one-to-one correspondence between logarithmic modifications of $\mathfrak{X} \to \mathcal{S}$ and subdivisions of $(B, \mathscr{P})$ and let $\pi': \tilde{\mathfrak{X}} \to \mathfrak{X}$ be the logarithmic modification corresponding to $(\tilde{B}, \tilde{\mathscr{P}})$. Then the dual intersection complex of the combined $\tilde{\mathfrak{X}} \to \mathcal{S}$ is isomorphic to $(\tilde{B}, \tilde{\mathscr{P}})$ and the condition on the maximal cells  $\tilde{\sigma} \in \tilde{\mathscr{P}}^{\max}$ implies that $\tilde{D}$ is simple normal crossings (note that $\tilde{\mathfrak{X}}_0$ is not reduced).
\end{proof}

Unlike Corollary \ref{furtherresolutionofadmissible}, we can't ensure that the resolution $\tilde{\mathfrak{X}} \to \bar{\mathfrak{X}}$ in Corollary \ref{sncafterfurtherblowupgeneral} is strongly admissible (see Remark \ref{nononsingularsubdivision}). However, one can guarantee the existence of a strongly admissible resolution with a simple normal crossings divisor by performing an additional basechange of $\bar{\mathfrak{X}} \to \mathcal{S}$.

\begin{fact} \label{basechangstronglyadmissiblesnc}
Let $\bar{\mathfrak{X}} \to \mathcal{S}$ be a special toric degeneration of CY threefolds satisfying Assumption \ref{naturalcase} and admitting a strongly admissible resolution. There exists an $L_{\snc} \in \mathbb{Z}_{>0}$ such that the basechange $\bar{g}': \bar{\mathfrak{X}}' \to \mathcal{S}$ of $\bar{\mathfrak{X}} \to \mathcal{S}$ by $R \to R, \ t \mapsto t^{L_{\snc}}$ (where $t$ is the uniformizer of $R$) admits a strongly admissible resolution $\pi': \mathfrak{X}' \to \bar{\mathfrak{X}}'$ with the divisor $D'' := (\pi' \circ \bar{g}')^{-1}(0)$ simple normal crossings. 
\end{fact}

\begin{proof}
Let $\pi: \mathfrak{X} \to \bar{\mathfrak{X}}$ be a strongly admissible resolution and let $(B, \mathscr{P})$ be the dual intersection complex of $\mathfrak{X} \to \mathcal{S}$ with a function $\alpha \in PA(B)$ satisfying the assumptions of condition (4) of Proposition \ref{onetoonetropicalgeneral}. Again, by \cite[Theorem 4.1]{KKMS}, there exists a rescaling $(B, L_{\snc}\mathscr{P})$ of $(B, \mathscr{P})$ that admits a subdivision $(\tilde{B}, \widetilde{L_{\snc}\mathscr{P}})$ with all $\tilde{\sigma} \in \widetilde{L_{\snc}\mathscr{P}}^{\max}$ standard simplices. Note that we may still view $(B, L_{\snc}\mathscr{P})$ as a subdivision of the corresponding rescaling $\left(\bar{B}, \overline{L_{\snc}{\mathscr{P}}}\right)$ of $\left( \bar{B}, \bar{\mathscr{P}}\right)$.  Further, as explained in \cite{KKMS}, the subdivision admits an $\alpha' \in PA(\tilde{B})$ that is strictly convex on the induced subdivision of every $\sigma \in L_{\snc} \mathscr{P}^{\max}$ and it is easy to check that it can be chosen to satisfy $\alpha'(v) = 0$ for any $v \in \bar{\mathscr{P}}^{[0]}$. Then by an argument similar to Step 4 of Proposition \ref{ampledivisor}, $\alpha' \in PA(\tilde{B})$ and $\alpha \in PA(B)$ give rise to an $\tilde{\alpha} \in PA(\tilde{B})$ that is strictly convex on the induced subdivision of every $\bar{\sigma} \in \overline{L_{\snc}{\mathscr{P}}}^{\max}$. 

It suffices to check that the subdivision $(\tilde{B}, \widetilde{L_{\snc}\mathscr{P}})$ of $\left(\bar{B}, \overline{L_{\snc}{\mathscr{P}}}\right)$ satisfies condition (2) of Proposition \ref{onetoonetropicalgeneral} and $\tilde{\alpha} \in PA(\tilde{B})$ satisfies condition (4) of Proposition \ref{onetoonetropicalgeneral}. Indeed, since $\left(\bar{B}, \overline{L_{\snc}{\mathscr{P}}}\right)$ can be seen as the dual intersection complex of the basechange $\bar{\mathfrak{X}}' \to \mathcal{S}$ of $\bar{\mathfrak{X}} \to \mathcal{S}$ by $R \to R, \ t \mapsto t^{L_{\snc}}$, the correspondence of Proposition \ref{onetoonetropicalgeneral} would give rise to a strongly admissible resolution $\pi': \mathfrak{X}' \to \bar{\mathfrak{X}}'$. The condition that $D''$ is simple normal crossings would follow from the fact that all $\tilde{\sigma} \in \widetilde{L_{\snc}\mathscr{P}}^{\max}$ are standard simplices.

To check that the subdivision $(\tilde{B}, \widetilde{L_{\snc}\mathscr{P}})$ of $\left(\bar{B}, \overline{L_{\snc}{\mathscr{P}}}\right)$ satisfies condition (2) of Proposition \ref{onetoonetropicalgeneral} note that since for every $\tau \in \bar{\mathscr{P}}^{[2]}$ the cells $\nu_{\tau_{\rel}} \in \mathscr{P}^{[2]}$ forming the subdivision of $\tau_{\rel}$ (using the notation of Remark \ref{generalsubdivisionreciperemark}) are standard triangles, the induced subdivision of $L_{\snc}\tau_{\rel}$ is given on every $L_{\snc}\nu_{\tau_{\rel}}$ by the standard subdivision of a rescaling of a standard triangle (i.e. the subdivision is by the lines parallel to the edges of the triangle). Now it is easy to check that one can choose a new $(L_{\snc}\tau)_{\rel}$ of the form required by Remark \ref{generalsubdivisionreciperemark} to lie inside the rescaling $L_{\snc}\tau_{\rel}$. Further, one can choose $\alpha' \in PA(\tilde{B})$ to be such that for every rescaling $L \sigma_{\tau_{\rel}}$ of a standard square $\sigma_{\tau_{\rel}}$ of $\tau_{\rel}$, $\alpha'|_{L \sigma_{\tau_{\rel}}}$ satisfies the natural analogues of conditions (4)(d) and (4)(e) of Proposition \ref{onetoonetropicalgeneral}. This implies that the corresponding $\tilde{\alpha} \in PA(\tilde{B})$ satisfies condition (4) of Proposition \ref{onetoonetropicalgeneral}.
\end{proof}

\begin{remark}
Note the difference between the statements of Theorem \ref{basechangstronglyadmissible} and Proposition \ref{basechangstronglyadmissiblesnc}. In the first case, we allow a basechange using any $L \in \mathbb{Z}_{>0}$ sufficiently large, whereas in the second case, the basechange uses a fixed $L_{\snc} \in \mathbb{Z}_{>0}$. The reason is purely combinatorial. It is known that for every polyhedral complex $\mathscr{C}$ there exists an $L_{\mathscr{C}} \in \mathbb{Z}_{>0}$ such that for every $L \in \mathbb{Z}_{>0}$ with $L \geq L_{\mathscr{C}}$, the rescaling $L\mathscr{C}$ admits a unimodular triangulation, see \cite{L2}. However, it is only conjectured (see the discussion of \cite[Section 1.4]{HHY}) that the same statement is true if one additionally requires that for every $L \in \mathbb{Z}_{>0}$ with $L \geq L_{\mathscr{C}}$ there exists a PA-function $\alpha'_L$ that is strictly convex on the induced subdivision of every $\sigma \in L\mathscr{C}^{\max}$. 
\end{remark}

\subsubsection{Resolutions in relative dimension \texorpdfstring{$n \geq 4$}{n4}} \label{resolveindimatleast4}

Suppose that $\bar{\mathfrak{X}} \to \mathcal{S}$ is a special toric degeneration of relative dimension $n \geq 4$ satisfying Assumptions \ref{naturalcase} and \ref{facelocalmodelass}. If all the maximal cells $\sigma \in \bar{\mathscr{P}}^{\max}$ are standard simplices, then the discussion of Section \ref{genericallysncgeneral} generalizes directly. One can also generalize the analysis of Section \ref{resolve} to special toric degenerations of relative dimension $n$ satisfying Assumptions \ref{naturalcase} and \ref{facelocalmodelass} by restricting the types of $\sigma \in \bar{\mathscr{P}}^{\max}$ to natural generalizations of the types (1-4) in Figure \ref{admisspolytopes}.

In general, one needs to understand how to resolve, for every $\tau \in \bar{\mathscr{P}}^{[n-1]}$ and every singular point $x \in \bar{X}_{\tau}$, the local toric model $\tilde{\bar{X}}_{\tau,x}$ (if we wanted to remove Assumption \ref{facelocalmodelass}, we would also have to consider local models with $\tau \in \bar{\mathscr{P}}^{[k]}, \ 1 \leq k \leq n-2$, see Remark \ref{facelocalmodel}(2)).

\begin{defn}
Let $n \geq 2$. We say that an integral polytope $\tau \subseteq \mathbb{Z}^{n-1}$ is a standard \emph{$(m,k)$-quplex} for some $m \in \mathbb{Z}_{\geq 0}, \ k \in \mathbb{Z}_{\geq 1}$ with $m + k = n-1$ if it is $AGL(n-1, \mathbb{Z})$-equivalent to $\sigma_m \times \Delta_k \subseteq \mathbb{Z}^{m} \times \mathbb{Z}^k$ for $\sigma_m$ the standard $m$-cube (i.e. the convex hull of $\{(\varepsilon_{1}, \dots, \varepsilon_{m}) \ | \ \varepsilon_i = 0,1 \text{ for } 1 \leq i \leq m\}$ in $\mathbb{Z}^m$) and $\Delta_k$ the standard $k$-simplex. Note that a standard $(n-2,1)$-quplex is a standard $(n-1)$-cube.
\end{defn}

If $n = 2$, the standard $(0,1)$-quplex is a line segment of length $1$. If $n=3$, the standard $(0,2)$-quplex is a standard triangle and the standard $(1,1)$-quplex is a standard square. If $n=4$, the standard $(0,3)$-quplex is a standard $3$-simplex, the standard $(1,2)$-quplex is a triangular prism, and the standard $(2,1)$-quplex is a standard $3$-cube.

It is straightforward to show that if $\tau \in \bar{\mathscr{P}}^{[n-1]}$ is a standard $(m,k)$-quplex, then $\Delta_{\tau} \cong \Conv\left((\tau \times \{e_0\}) \cup  \Delta_{\tau,g} \right)$ for $\Delta_{\tau,g} = \Delta_{\tau,1}$ a standard $k$-simplex using the description of \eqref{conv2}. Further, one can also generalize Figures \ref{locmodeltriangle} and \ref{locmodelsquare} to obtain subdivisions of $\Delta_{\tau}$ into standard $n$-simplices. This gives rise to log smooth resolutions of the local toric models $\tilde{\bar{X}}_{\tau,x}$. One can also generalize Figures \ref{singularlocusltriangle} and \ref{singularlocuslsquare} to describe the transformations of the singular loci of $\tau$ corresponding to these resolutions.

\begin{conj} \label{combconjecture}
Consider a cell $\tau \in \bar{\mathscr{P}}^{[n-1]}$ and suppose that $$\Delta_{\tau} \cong \Conv\left((\tau \times \{e_0\}) \cup  \Delta_{\tau,g} \right)$$ for $\Delta_{\tau,g}$ a standard $k$-simplex. Then for all $L \in \mathbb{Z}_{>0}$ sufficiently large, there exists a subdivision of $L\tau$ into $(n-1)$-simplices and standard $(m,k')$-quplexes for $1 \leq k' \leq k$ (if $k = n-1$ then there is exactly one $(0,k)$-quplex in the subdivision) such that no $(m,k')$-quplex intersects cells of $L\tau$ (considered as a polyhedral manifold with boundary) of dimension less than $k'$. 

Further, let $L\tau_{\rel}$ be the union of all the $(m,k')$-quplexes (for $1 \leq k' \leq k$) in the subdivision. Then there exists a choice of a singular locus $Z_{L\tau}$ on $L\tau$ respecting the cell structure such that $L\tau_{\rel}$ is a tubular neighbourhood of $Z_{L\tau}$ and for every $(m,k')$-quplex $\tau'$ in the subdivision, the singular locus of $\tau'$ is given by $Z_{L\tau} \cap \tau'$.
\end{conj}

Note that Conjecture \ref{combconjecture} holds trivially if $\tau \in \bar{\mathscr{P}}^{[n-1]}$ is a standard $(m,k)$-quplex. It is also easy to check that it holds for a rescaling of a standard $(m,k)$-quplex. Conjecture \ref{combconjecture} is trivial for $n=2$ and Steps (1-6) of Construction \ref{generalsubdivisionrecipe} prove that Conjecture \ref{combconjecture} holds for $n = 3$. In general, we expect Conjecture \ref{combconjecture} to follow by an inductive argument. Indeed, for $n=4$, Steps (1-6) of Construction \ref{generalsubdivisionrecipe} imply, by an argument as in the construction of the integral subdivision $(B', \mathscr{P}')$ in the proof of Theorem \ref{basechangstronglyadmissible}, that for every $\tau \in \bar{\mathscr{P}}^{[3]}$ one can construct a subdivision of $\partial(L \tau)$ as required by Conjecture \ref{combconjecture}. The fact that one can construct a subdivision of $L\tau$ restricting to this subdivision of $\partial(L \tau)$ (possibly after rescaling by some $k \in \mathbb{Z}_{>0}$ and choosing a new subdivision of $\partial(kL \tau)$ of this form) should follow from the fact that $\Delta_{\tau,g} = \Delta_{\tau,1}$ is the monodromy simplex associated to $\tau$ (in particular, it is a Minkowski summand of $\tau$) and from \cite[Definition 1.60]{GS2} of the simplicity of $\left(\bar{B}, \bar{\mathscr{P}}\right)$ (satisfied by Assumption \ref{naturalcase}).

Suppose that $\tau \in \bar{\mathscr{P}}^{[n-1]}$ satisfies Conjecture \ref{combconjecture}. Then the subdivision of $L\tau$ of Conjecture \ref{combconjecture} lifts to a subdivision of $\Delta_{L\tau}$ by directly generalizing the lifting of Construction \ref{generalsubdivisionrecipe}. To obtain a log smooth resolution, it remains, similarly to Construction \ref{generalsubdivisionrecipe}(7), to subdivide the polygons of the subdivision of $\Delta_{L\tau}$ containing the standard $(m,k')$-quplex (for $1 \leq k' \leq k$). We subdivide them using the natural generalization of Figures \ref{locmodeltriangle} and \ref{locmodelsquare} (as usual, this should be done compatibly). This gives both the subdivision of $L\tau$ and a lifting to a subdivision of $\Delta_{L\tau}$ defining a toric blowup $\pi_{\tau,x}: \tilde{X}_{\tau,x} \to \tilde{\bar{X}}_{\tau,x}$ resolving the local model.

Suppose that $\bar{\mathfrak{X}} \to \mathcal{S}$ is a special toric degeneration of relative dimension $n \geq 4$ satisfying Assumptions \ref{naturalcase} and \ref{facelocalmodelass}, such that every $\tau \in \bar{\mathscr{P}}^{[n-1]}$ with $\bar{\mathfrak{X}} \to \mathcal{S}$ not log smooth on $\bar{X}_{\tau}$ satisfies Conjecture \ref{combconjecture}. Then Definition \ref{stronglyadmissiblegeneral} of a strongly admissible resolution and Definition \ref{admissibledefgeneral} of an admissible resolution directly generalize to this case, see Remark \ref{furtheradmissiblegenremarks}(1). Moreover, it is straightforward to generalize Proposition \ref{onetoonetropicalgeneral} giving a tropical way to construct strongly admissible resolutions, Theorem \ref{basechangstronglyadmissible} giving an existence result for strongly admissible resolutions and Propositions \ref{sncafterfurtherblowupgeneral} and \ref{basechangstronglyadmissiblesnc} giving admissible and strongly admissible (after a further basechange) resolutions with a simple normal crossings central fibre.

\subsection[Scattering and the Gross-Siebert locus]{Scattering and the minimal relative Gross-Siebert locus}

Let $\bar{\mathfrak{X}} \to \mathcal{S}$ be a special toric degeneration of relative dimension $n \geq 3$ satisfying Assumption \ref{naturalcase} and Assumption \ref{facelocalmodelass} (if $n \geq 4$). By the discussion of Section \ref{resolvehd}, we can construct strongly admissible resolutions of $\bar{\mathfrak{X}} \to \mathcal{S}$ (assuming Conjecture \ref{combconjecture} if $n \geq 4$). 

Suppose that we want to prove Conjecture \ref{conjecture} for $\bar{\mathfrak{X}} \to \mathcal{S}$. As in the proof of Theorem \ref{proofconj}, Proposition \ref{furtherresolution} implies that it is enough to prove it for strongly admissible resolutions, and Proposition \ref{enoughforwellchosen} implies that it is enough to prove it for well-chosen monoids. So we may fix a strongly admissible resolution $\pi: \mathfrak{X} \to \bar{\mathfrak{X}}$ and a well-chosen monoid $P$. Let $J := P \setminus K$ be the complement of the face containing the curve classes contracted by $\pi$ as before. By the discussion at the start of this chapter, we need to provide a collection of scattering diagrams $\mathfrak{D}^J: = \left\{ \mathfrak{D}_{J^{k+1}}, \ k \geq 0  \right\}$ giving rise to the extended intrinsic mirror $\check{\mathfrak{X}} \to \Spec \widehat{\kk [P]}_J$, generalizing the results of Section \ref{scatint}.

\subsubsection{Canonical scattering modulo \texorpdfstring{$J$}{J}}

Let $\left( \bar{B}, \bar{\mathscr{P}} \right)$ be the dual intersection complex of $\bar{\mathfrak{X}} \to \mathcal{S}$ and $(B, \mathscr{P})$ be the dual intersection complex of $\mathfrak{X} \to \mathcal{S}$. First, consider the collection $\bar{\mathfrak{D}} = \left\{ \bar{\mathfrak{D}}_k, \ k \geq 0 \right\}$ of the algorithmic scattering diagrams on $\left( \bar{B}, \bar{\mathscr{P}} \right)$. Recall from Theorem \ref{reconstruction} that $\bar{\mathfrak{D}}$ is determined (up to combinatorial equivalence) by $\bar{\mathfrak{D}}_0 = \bar{\mathfrak{D}}_{(t)}$ whose only walls are the slabs $(\underline{\tau}, f_{\underline{\tau}})$ with support $\underline{\tau} \in \tilde{\bar{\mathscr{P}}}^{[n-1]}$ and the attached initial slab function $f_{\underline{\tau}}$.\footnote{Here we switch the notation for elements of $\bar{\mathscr{P}}^{[n-1]}$ (resp. $\tilde{\bar{\mathscr{P}}}^{[n-1]}$) from $\rho$ (resp. $\underline{\rho}$) to $\tau$ (resp. $\underline{\tau}$) to align with the notations of Section \ref{resolvehd}.} Moreover, since $\left(\bar B, \bar{\mathscr{P}}\right)$ is simple (in the sense of \cite[Definition 1.60]{GS2}) by Assumption \ref{naturalcase}, the initial slab functions are fixed and independent of the gluing data.

\begin{fact}[{\cite[Lemma A.1.1]{GHS}}] \label{prop65}
Let $\underline{\tau} \in \tilde{\bar{\mathscr{P}}}^{[n-1]}$ be a slab with $\underline{\tau} \subseteq \tau \in \bar{\mathscr{P}}^{[n-1]}$ and let $v$ be the vertex of $\tau$ contained in $\underline{\tau}$. Let $$\bar{\Delta} (\tau, v) := \left\{ m_{\underline{\tau} \, \underline{\tau}'} \ | \ \underline{\tau}' \in \tilde{\bar{\mathscr{P}}}'^{[n-1]}, \ \underline{\tau}' \subseteq \tau \right\}$$ be the set of monodromy vectors of closed paths in $W_{\rho} \setminus \bar{\Delta}$ as in \cite[Appendix A.1]{GHS} (since $\left(\bar B, \bar{\mathscr{P}}\right)$ is simple, $\bar{\Delta} (\tau, v)$ is the set of vertices of an elementary simplex). Then $f_{\underline{\tau}} = \sum_{m \in \bar{\Delta} (\tau, v)} z^m$. 
\end{fact}

Now, note that $(B, \mathscr{P})$ is an integral subdivision of $\left( \bar{B}, \bar{\mathscr{P}} \right)$ (see Proposition \ref{onetoonetropicalgeneral}(2)). The construction of a PL-isomorphism $\Phi: \left(B, \mathscr{P}\right) \to \left(\bar{B}, \bar{\mathscr{P}}\right)$, linear on the maximal cells of $\mathscr{P}$, is a direct generalization of Construction \ref{plmapadmissible} (using Constructions \ref{plisosmall} and \ref{plmap}). In particular, for every $x \in \Int(\sigma), \ \sigma \in \mathscr{P}^{\max}$ we can identify the monomials via the canonical isomorphism $\mathcal{P}_x^{+} \cong \mathcal{\bar{P}}_{\Phi(x)}^{+}$. The following conjecture generalizes Proposition \ref{diag} and Construction \ref{newwalls} of the scattering diagram $\mathfrak{D}_J$ in the case that $n=2$.

\begin{conj} \label{canonicalscathd}
The decorated wall types given by Construction \ref{canonical scattering} of the canonical scattering diagram modulo $J$ for $\mathfrak{X} \to \mathcal{S}$ satisfy the following conditions:
\begin{enumerate}
    \item Let $\btau$ be a decorated wall type with total curve class $A \in K$ and $W_{\btau} \neq 0$. Then the support of the corresponding wall $\mathfrak{b}_{\btau}$ given by Construction \ref{canonical scattering} is a cell $\tau' \in \mathscr{P}^{[n-1]}$. Moreover, if $\tau' \in \mathscr{P}^{[n-1]}$ subdivides a cell $\tau \in \bar{\mathscr{P}}^{[n-1]}$, then $\dim(\tau \cap \bar{\Delta}) = n-2$. 
    \item Let $\tau' \in \mathscr{P}^{[n-1]}$ be a cell of $(B, \mathscr{P})$ subdividing a cell $\tau \in \bar{\mathscr{P}}^{[n-1]}$ with $\dim(\tau \cap \bar{\Delta}) = n-2$. Then there are countably many decorated wall types $\btau$ with total curve class $A \in K$ and $W_{\btau} \neq 0$ producing walls $\{\mathfrak{b}_{\tau',b} \ | \ b \in \mathbb{Z}_{>0} \}$ supported on $\tau'$ via Construction \ref{canonical scattering}. The infinite product 
    \begin{equation} \label{infproductwallshd}
    f_{\mathfrak{b}_{\tau'}} := \prod_{b \in \mathbb{Z}_{>0}} f_{\mathfrak{b}_{\tau',b}} \in (\kk[P]/J)[[\Lambda_{\tau'}]]
    \end{equation}
    of the wall functions is polynomial, i.e. $f_{\mathfrak{b}_{\tau'}} \in (\kk[P]/J)[\Lambda_{\tau'}]$.
    \item Let $A$ be a choice of polarization on $\bar{\mathfrak{X}} \to \mathcal{S}$ and let $h: P \to \mathbb{N}, \ \beta \mapsto \pi^{*}A \cdot \beta$ be the map of Conjecture \ref{conjecture}. Then we have $h(f_{\mathfrak{b}_{\tau'}}) = f_{\underline{\tau}}$ for some choice of a slab $\underline{\tau} \in \tilde{\bar{\mathscr{P}}}^{[n-1]}$ with $\underline{\tau} \subseteq \tau \in \bar{\mathscr{P}}^{[n-1]}$ and $f_{\underline{\tau}}$ as in Proposition \ref{prop65} (using the identification of monomials). 
\end{enumerate}
Assuming (1),(2),(3) are satisfied, we let $\mathfrak{b}_{\tau'}:=(\tau', f_{\mathfrak{b}_{\tau'}})$ for any $\tau' \in \mathscr{P}^{[n-1]}$ and set $$\mathfrak{D}_J:=\{\mathfrak{b}_{\tau'} \ | \ \tau' \subseteq \tau, \ \tau' \in \mathscr{P}^{[n-1]}, \ \tau \in \bar{\mathscr{P}}^{[n-1]}, \ \dim(\tau \cap \bar{\Delta}) = n-2\}.$$
By an argument as in the proof of Proposition \ref{consistencyandiso}, $\mathfrak{D}_J$ is a well-defined consistent scattering diagram.
\end{conj}

We will now discuss how we expect Conjecture \ref{canonicalscathd} to follow.
Consider the case of Section \ref{genericallysncgeneral}, i.e. assume that $n=3$ and all the $\sigma \in \bar{\mathscr{P}}^{\max}$ are standard simplices. Let $\tau \in \bar{\mathscr{P}}^{[2]}$ be a standard triangle with $\dim(\tau \cap \bar{\Delta}) = 1$ and note that it may be viewed as a cell of $\mathscr{P}^{[2]}$ since $(B, \mathscr{P}) \cong \left(\bar B, \bar{\mathscr{P}}\right)$ as polyhedral manifolds. Let $\underline{\tau}$ be as in Conjecture \ref{canonicalscathd}(3) and write $f_{\underline{\tau}}$ as $f_{\underline{\tau}} = 1 + z^{m_1} + z^{m_2}$ where $m_1, m_2$ are the non-zero vertices of $\bar{\Delta} (\tau, v)$. Then we expect that $f_{\mathfrak{b}_{\tau'}} = 1 + t^{E_{v_1} - E_{v_2}}z^{m_1} + t^{E_{v_2}} z^{m_2}$ where $E_{v_1}$ and $E_{v_2}$ (for $v_1, v_2, \in \bar{\mathscr{P}}^{[0]}$ with $v_1, v_2 \subseteq \tau$) are the exceptional curves of the first and second blowup resolving the local model $\tilde{\bar{X}}_{\tau,x}$ (permuting $m_1$ and $m_2$ if necessary). Now, note that we have 
\begin{multline}
f_{\mathfrak{b}_{\tau'}} = 1 + t^{E_{v_1}- E_{v_2}}z^{m_1} + t^{E_{v_2}} z^{m_2} = \exp(\log(1 + t^{E_{v_1}- E_{v_2}}z^{m_1} + t^{E_{v_2}} z^{m_2})) = \\ = \exp \left(\sum_{\substack{a_1,a_2 \in \mathbb{N}, \\ a_1+a_2 > 0}} \frac{(-1)^{a_1+a_2-1}}{a_1+a_2} \binom{a_1+a_2}{a_1} t^{a_1(E_{v_1}- E_{v_2})+a_2E_{v_2}} z^{a_1m_1+a_2m_2} \right) = \\ = \prod_{\substack{a_1,a_2 \in \mathbb{N}, \\ a_1+a_2 > 0}} \exp \left( \frac{(-1)^{a_1+a_2-1}}{a_1+a_2} \binom{a_1+a_2}{a_1} t^{a_1(E_{v_1}- E_{v_2})+a_2E_{v_2}} z^{a_1m_1+a_2m_2} \right).
\end{multline}

We expect that this is precisely the infinite product of \eqref{infproductwallshd} giving rise to $f_{\mathfrak{b}_{\tau'}}$. Specifically, we conjecture the following:

\begin{conj} \label{walltypesdim3hd}
 Suppose that $n=3$ and all the $\sigma \in \bar{\mathscr{P}}^{\max}$ are standard simplices. Let $\tau \in \bar{\mathscr{P}}^{[2]}$ be a standard triangle with $\dim(\tau \cap \bar{\Delta}) = 1$. Let $E_{v_1}$ and $E_{v_2}$ (for $v_1, v_2, \in \bar{\mathscr{P}}^{[0]}$ with $v_1, v_2 \subseteq \tau$) be the exceptional curves of the first and second blowup resolving the local model $\tilde{\bar{X}}_{\tau,x}$. Let $v \in \bar{\mathscr{P}}^{[0]}$ be the third vertex of $\tau$ (distinct from $v_1,v_2$) and let $m_1 \in \Lambda_{\tau}$ and $m_2 \in \Lambda_{\tau}$ be the primitive generators of $\langle v, v_1 \rangle$ and $\langle v, v_2 \rangle$ respectively that point out of $v$. Then the decorated wall types $\btau$ with total curve class $A \in K$ and $W_{\btau} \neq 0$ producing walls supported on $\tau$ are as follows. 
 
 Let $a_1, a_2 \in \mathbb{N}$ with $a_1+a_2 > 0$. Then there is a decorated wall type $\btau_{a_1,a_2}$ producing a wall supported on $\tau$ with 
 \begin{equation} \label{contributionofshrub}
 k_{\tau_{a_1,a_2}}W_{\btau_{a_1,a_2}} = \frac{(-1)^{a_1+a_2-1}}{a_1+a_2} \binom{a_1+a_2}{a_1}
 \end{equation}
 Further, the underlying graph $G_{a_1,a_2}$ of $\btau_{a_1,a_2}$ has two vertices $\tilde{v}_0, \tilde{v}_1$, one edge $\langle \tilde{v}_0, \tilde{v}_1 \rangle$, and a leg $L_{\out}$ attached to $\tilde{v}_0$. We have $\bsigma(\tilde{v}_0) = \mathbf{C}\langle v_1, v_2 \rangle$ and $\bsigma(\tilde{v}_1) = \mathbf{C}v_1$. Let $m_{12} \in \Lambda_{\tau}$ be the generator of $\langle v_1, v_2 \rangle$ pointing from $v_1$ to $v_2$.  We have $\mathbf{u}(L_{\out}) = -a_1m_1-a_2m_2$ and $\mathbf{u}(\langle \tilde{v}_0, \tilde{v}_1 \rangle) = -a_1m_{12}$. Finally, we have $\mathbf{A}(\tilde{v}_0) = a_2 E_{v_2}$ and $\mathbf{A}(\tilde{v}_1) = a_1(E_{v_1} - E_{v_2})$. We give a sketch (in blue) of the corresponding family of tropical curves in $(B, \mathscr{P})$ in Figure \ref{typepicture} (the image of $\tilde{v}_0$ is free to move inside the edge). 
\end{conj}

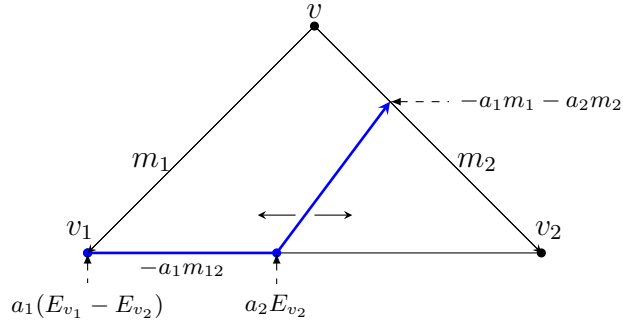
\begin{figure}[t]
\begin{center}
\begin{tikzpicture}

\draw[fill] (0,0) circle (1.5pt) coordinate (m-0-0);
\draw[fill] (6,0) circle (1.5pt) coordinate (m-6-0);
\draw[fill] (3,3) circle (1.5pt) coordinate (m-3-3);

\draw (m-0-0) -- (m-6-0) -- (m-3-3) -- cycle;

\draw(-0.1, 0.3) node{$v_1$};
\draw(6.1, 0.3) node{$v_2$};
\draw(3, 3.2) node{$v$};

\draw(0.85, 1.2) node{$m_1$};
\draw(5.15, 1.2) node{$m_2$};

\draw(1.25, -0.2) node{\scriptsize $-a_1m_{12}$};

\draw [-stealth, color=blue, line width=1pt] (2.5,0)-- (4,2);
\draw [color=blue, line width=1pt] (2.5,0)-- (0,0);

\draw[fill, color=blue] (2.5,0) circle (1.5pt);
\draw[fill, color=blue] (0,0) circle (1.5pt);

\draw [-stealth] (3,3)-- (0,0);
\draw [-stealth] (3,3)-- (6,0);

\draw [-stealth] (3,0.5)-- (3.5,0.5);
\draw [-stealth] (2.75,0.5)-- (2.25,0.5);

\draw[latex-,dashed] 
        (4,2) -- (6, 2)
         node[fill=white]{\scriptsize $-a_1m_1-a_2m_2$};

\draw[latex-,dashed] 
        (0,0) -- (0, -0.7)
         node[fill=white]{\scriptsize $a_1(E_{v_1} - E_{v_2})$};

\draw[latex-,dashed] 
        (2.5,0) -- (2.5, -0.7)
         node[fill=white]{\scriptsize $a_2E_{v_2}$};

\end{tikzpicture}
\end{center}
\setcounter{figure}{12}
\caption{The family of tropical curves in $(B, \mathscr{P})$ corresponding to $\btau_{a_1,a_2}$.} 
\label{typepicture}
\end{figure}

In the notations of Conjecture \ref{walltypesdim3hd}, let $C_{v_2}$ be the unique curve of class $E_{v_2}$ that intersects $X_{\tau}$ and let $C_{v_1}$ be the unique curve of class $E_{v_1} - E_{v_2}$ that intersects $C_{v_2}$ (see Figure \ref{locmodeltriangleresolve}). The decorated wall type $\btau_{a_1,a_2}$ corresponds to a stable map $C \to \mathfrak{X}_0$. Here $C$ consists of two rational curves $C_0, C_1$ (corresponding to $\tilde{v}_0, \tilde{v}_1$) and the map is an $a_1$-fold cover of $C_{v_1}$ on $C_1$ and an $a_2$-fold cover of $C_{v_2}$ on $C_2$.

By the discussion above and Construction \ref{canonical scattering} of the canonical scattering diagram, Conjecture \ref{walltypesdim3hd} implies Conjecture \ref{canonicalscathd} in the case that $n=3$ and all the $\sigma \in \bar{\mathscr{P}}^{\max}$ are standard simplices. The fact that the only wall types $\btau$ with total curve class $A \in K$ and $W_{\btau} \neq 0$ producing walls supported on $\tau \in \mathscr{P}^{[2]}$ are the $\btau_{a_1,a_2}$ should follow from a tropical analysis (generalizing the argument in the proof of Proposition \ref{diag}). The fact that \eqref{contributionofshrub} holds should follow by a suitable gluing formula, generalizing \cite[Theorem 8.15]{G3}. Note that we are in the case of toric gluing strata, so one should be able to use the convolution-style formula of \cite[Corollary 1.6]{W}. In any case, the gluing formula will reduce the computation of $k_{\tau_{a_1,a_2}}W_{\btau_{a_1,a_2}}$ to the computation of $k_{\tau^1_{a_1,a_2}}W_{\btau^1_{a_1,a_2}}$, where $\btau^1_{a_1,a_2}$ is the sub-wall type of $\btau_{a_1,a_2}$ obtained by splitting at $\tilde{v}_0$. Note that we have $k_{\tau^1_{a_1,a_2}} = a_1$. 

Now $\btau^1_{a_1,a_2}$ can be viewed as a decorated wall type $\tilde{\btau}^1_{a_1,a_2}$ of a punctured map to a log CY surface $\tilde{D}_{v_1}:=(D_{v_{1}}, \partial D_{v_{1}})$ (using the notation of Observation \ref{altaffinecons}). Further, $\tilde{\btau}^1_{a_1,a_2}$ has positive contact order $a_1$ with the component of $\partial D_{v_{1}}$ that intersects the exceptional $\mathbb{P}^1$-bundle with exceptional curve $E_{v_1}$, so it can be viewed as a type of an ordinary logarithmic map. It also has an induced curve class $\tilde{A}_{\tilde{\btau}^1_{a_1,a_2}}$ that we can write as $\tilde{A}_{\tilde{\btau}^1_{a_1,a_2}} =a_1(E_{v_1} - E_{v_2})$. We have $\mathscr{M}(\mathfrak{X},\btau^1_{a_1,a_2}) \simeq \mathscr{M}(\tilde{D}_{v_1},\tilde{\btau}^1_{a_1,a_2})$ by \cite[Theorem 6.1]{G3}, so $\deg [\mathscr{M}(\mathfrak{X},\btau^1_{a_1,a_2})]^{\virt} = \deg [\mathscr{M}(\tilde{D}_{v_{1}},\tilde{\btau}^1_{a_1,a_2})]^{\virt}$. Therefore, we have $W_{\btau^1_{a_1,a_2}} = N_{\tilde{\btau}^1_{a_1,a_2}}$ where $N_{\tilde{\btau}^1_{a_1,a_2}}$ is the virtual count corresponding to $\tilde{\btau}^1_{a_1,a_2}$. We conjecture that $N_{\tilde{\btau}^1_{a_1,a_2}} = \frac{(-1)^{a_1-1}}{a_1^2}$. More generally, we conjecture the following.

\begin{conj} \label{toricresolvecount}
Let $X_{\Sigma}$ be a toric variety and let $(X, D)$ be a log CY variety obtained by blowing up a sequence of hypersurfaces $H_i \subseteq D_{\rho_i}$ with $\cup H_i$ connected for some toric divisors $D_{\rho_i}$ corresponding to rays $\rho_i \in \Sigma$. Let $E$ be the exceptional curve of the last blowup. For every $a \in \mathbb{N}$, there is a unique decorated wall type $\tilde{\btau}_a$ with total curve class $aE$. Moreover, the corresponding virtual count is $N_{\tilde{\btau}_a} = \frac{(-1)^{a-1}}{a^2}$. 
\end{conj}

If $X_{\Sigma}$ is a toric surface, then by the comparison of the logarithmic and relative invariants of \cite{AMW} and the fact that $(X, D)$ is a Looijenga pair, Conjecture \ref{toricresolvecount} follows from \cite[Proposition 5.2]{GPS}. In higher dimensions, if one blows up only one hypersurface and $D$ is simple normal crossings, the setup falls into that of \cite{AG}. Conjecture \ref{toricresolvecount} should follow from \cite[Theorem 6.1]{AG} in this case (see the discussion at the end of \cite[Section 7]{AG}). In general, Conjecture \ref{toricresolvecount} should follow by generalizing the results of \cite{AG} to the case when $D$ is not simple normal crossings and the hypersurfaces $H_i$ are not disjoint, and performing a similar analysis. 

The whole discussion directly generalizes to the case when $n \geq 4$ and $\sigma \in \bar{\mathscr{P}}^{\max}$ are standard simplices. For general $\sigma \in \bar{\mathscr{P}}^{\max}$, let $\tau \in \bar{\mathscr{P}}^{[n-1]}$ and suppose that $k=n-1$ in Conjecture \ref{combconjecture}. Then there is exactly one $(0, k)$-quplex $\delta' \in \mathscr{P}^{[n-1]}$ in the subdivision of $\tau$ which is just the $(n-1)$-simplex. So the same reasoning applies to showing Conjecture \ref{canonicalscathd} for $\delta'$. For the other cells in the subdivision of $\tau$, Conjecture \ref{canonicalscathd} should follow by analyzing the local models for those cells and by inductively reducing the computation of the relevant virtual invariants to those contributing to the $f_{\mathfrak{b}_{\delta'}}$ of \eqref{infproductwallshd} by using a suitable gluing formula. Similarly, if $k < n-1$ in Conjecture \ref{combconjecture}, one needs to analyse the $(n-1-k, k)$-quplexes subdividing $\tau \in \bar{\mathscr{P}}^{[n-1]}$ first and then proceed inductively.

\subsubsection{Canonical scattering modulo \texorpdfstring{$J^{k+1}, \ k \geq 0$}{Jk} and Conjecture \ref{conjecture} in higher dimensions}

We have defined $\mathfrak{D}_J$ in Conjecture \ref{canonicalscathd}. Now, we also need to define scattering diagrams $\mathfrak{D}_{J^{k+1}}, \ k \geq 1$.

\begin{conj} \label{scatdiaghigherpowershd}
\begin{sloppypar}
There exist finite consistent scattering diagrams $\mathfrak{D}_{J^{k+1}}$, $\ k \geq 1$ such that:
\end{sloppypar}
\begin{enumerate}
    \item $\mathfrak{D}_{J^{k+1}}$ is compatible with $\mathfrak{D}_{J^{k}}$ (in the sense of Definition \ref{compatibility}) for $k \geq 1$.
    \item For every $k \geq 1$ and $l \geq 1$, the scattering diagram $\mathfrak{D}_{J^{k+1}}$ agrees with the scattering diagram $\mathfrak{D}_{J^{k+1} + \mathfrak{m}^l}$ of Construction \ref{canonical scattering} modulo $J^{k+1} + \mathfrak{m}^l$.
\end{enumerate}
\end{conj}

Conjecture \ref{scatdiaghigherpowershd} should follow by generalizing the results of \cite[Section 5.3]{GHKS} and Section \ref{modjkscatter} to higher dimensions. We expect the generalization to follow the same lines, i.e. the notion of slab twig of Definition \ref{slabtwig} should generalize to mean the wall types contributing to $\mathfrak{D}_{J}$. Then one can define slab types by grouping wall types that differ by the (generalized) slab twigs similarly to Definitions \ref{slabequiv} and \ref{decslabtype}. After this, the scattering diagrams $\mathfrak{D}_{J^{k+1}}, \ k \geq 1$ can be defined in the same way as in Definition \ref{extendcanonicalscat}. Consistency of $\mathfrak{D}_{J^{k+1}}, \ k \geq 1$ follows as in Step 1 in the proof of Proposition \ref{finitenessresult} and the hard part is to check that $\mathfrak{D}_{J^{k+1}}, \ k \geq 1$ are finite, generalizing Steps 2-4 in the proof of Proposition \ref{finitenessresult}.

Assuming Conjectures \ref{canonicalscathd} and \ref{scatdiaghigherpowershd}, we can define the extended intrinsic mirror $\check{\mathfrak{X}} \to \Spec \widehat{\kk [P]}_J$ as the inverse limit of the families $\check{\mathfrak{X}}_{\mathfrak{D}_{J^{k+1}}} \to \Spec \kk[P]/J^{k+1}$ for $k \geq 0$. Note that, as in Remark \ref{stronglyadmissibleweakening}(1), Conjecture \ref{scatdiaghigherpowershd}(2) implies that the extended intrinsic mirror of \eqref{extendedmirror} is well-defined and isomorphic to $\check{\mathfrak{X}} \to \Spec \widehat{\kk [P]}_J$. We can reduce Conjecture \ref{conjecture} to Conjectures \ref{combconjecture}, \ref{canonicalscathd}, and \ref{scatdiaghigherpowershd}.

\begin{fact}
Let $\bar{\mathfrak{X}} \to \mathcal{S}$ be a special toric degeneration of relative dimension $n \geq 3$ satisfying Assumption \ref{naturalcase} and Assumption \ref{facelocalmodelass} (if $n \geq 4$). Suppose that Conjecture \ref{combconjecture} holds for $\bar{\mathfrak{X}} \to \mathcal{S}$ (we have proved this for $n=3$) and Conjectures \ref{canonicalscathd} and \ref{scatdiaghigherpowershd} hold for any strongly admissible resolution $\mathfrak{X} \to \bar{\mathfrak{X}}$. Then Conjecture \ref{conjecture} holds for $\bar{\mathfrak{X}} \to \mathcal{S}$.
\end{fact}

\begin{proof}
We need to generalize the results of Section \ref{main}. Construction \ref{ccc} generalizes directly to give a collection of scattering diagrams $\mathfrak{D} = \left\{ \mathfrak{D}_k, \ k \geq 0 \right\}$ on $(B, \mathscr{P})$ such that the inverse limit over $\check{\mathfrak{X}}_{\mathfrak{D}_k} \to \Spec \kk [t]/(t^{k+1})$ for $k \geq 0$ is isomorphic to the basechange of the extended intrinsic mirror family $\check{\mathfrak{X}} \to \Spec \widehat{\kk [P]}_J$ (defined via  Conjectures \ref{canonicalscathd} and \ref{scatdiaghigherpowershd}) by $h: P \to \mathbb{N}, \ \beta \mapsto \pi^{*}A \cdot \beta$. 

Construction \ref{refinement} also admits a direct generalization which gives, for every scattering diagram $\bar{\mathfrak{D}}$ on $\left(\bar{B}, \bar{\mathscr{P}}\right)$, a scattering diagram $\bar{\mathfrak{D}}'$ on $\left(\bar{B}, \bar{\mathscr{P}}'\right)$ (where $\bar{\mathscr{P}}'$ is the obvious refined polyhedral decomposition on $\bar{B}$ such that $\left(B, \mathscr{P}\right) \cong \left(\bar B, \bar{\mathscr{P}}'\right)$ as polyhedral complexes) that is equivalent to $\bar{\mathfrak{D}}$. 

To generalize Construction \ref{chs}, we need to define a scattering diagram $\Phi (\mathfrak{D}_k)'$ on $\left(\bar B, \bar{\mathscr{P}}'\right)$. First, we define the images of codimension $0$ walls as in Construction \ref{chs}. Suppose that $\mathfrak{b} \in \mathfrak{D}_k$ is a slab with $\mathfrak{b} \subseteq \tau' \in \mathscr{P}^{[n-1]}$ for $\tau' \subseteq \tau \in \bar{\mathscr{P}}^{[n-1]}$ and let $\bar{\Phi}(\mathfrak{b})^{i}, \ 1 \leq i \leq m$ be the connected components of $\Phi(\mathfrak{b}) \setminus \Delta$. Let $\underline{\tau} \in \tilde{\bar{\mathscr{P}}}^{[n-1]}$ be a slab of $\tau$ such that $h(f_{\mathfrak{b}_{\tau'}}) = f_{\underline{\tau}}$ (see Conjecture \ref{canonicalscathd}(3)) and let $\underline{\tau}^i \in \tilde{\bar{\mathscr{P}}}^{[n-1]}$ be the slabs of $\tau$ containing $\bar{\Phi}(\mathfrak{b})^{i}$ for $1 \leq i \leq m$ respectively. Then we introduce $m$ slabs as follows:
$$\bar{\Phi}(\mathfrak{b})^{i} := \left(\bar{\Phi}(\mathfrak{b})^{i}, z^{m_{\underline{\tau}^i\underline{\tau} }}f_{\Phi(\mathfrak{b})}\right) \in \Phi (\mathfrak{D}_k)', \quad 1 \leq i \leq m$$ (here we set $m_{\underline{\tau}\,\underline{\tau}}:=0$). This defines a scattering diagram $\Phi (\mathfrak{D}_k)'$ on $\left(\bar B, \bar{\mathscr{P}}'\right)$. Then we also have an equivalent scattering diagram $\Phi (\mathfrak{D}_k)$ on $\left(\bar B, \bar{\mathscr{P}}\right)$ similarly to Construction \ref{chs}.

Propositions \ref{consis} and \ref{equivalence} generalize directly to show that for every $k \geq 0$, the scattering diagram $\Phi(\mathfrak{D}_k)$ is consistent and equivalent to $\mathfrak{D}_k$. Now Conjecture \ref{conjecture} follows by an argument as in the proof of Theorem \ref{proofconj}. 
\end{proof}

\subsubsection{A generalization of Conjecture \ref{conjecture} to the minimal relative Gross-Siebert locus}

We conclude by explaining how to generalize Theorem \ref{compgen} (the main result of Chapter \ref{gslocus}) to higher dimensions. The desire to have such a generalization is another reason that we restrict to the case of simple $\left(\bar B, \bar{\mathscr{P}}\right)$. Indeed, the local rigidity condition is not empty for $n \geq 3$, so for general $\left(\bar B, \bar{\mathscr{P}}\right)$ the space of toric log CY spaces with intersection complex $\left(\bar B, \bar{\mathscr{P}}\right)$ is not well-behaved. In particular, it may be singular, see \cite[Example 4.28]{GS2}. So, in general, there is no good analogue of \cite[Theorem A.2.4]{GHS} (or \cite[Theorem A.4.2]{GHS}) constructing the universal toric degeneration mirror varied in projective gluing data (or projective gluing data and free parameters of the initial slab functions). This means that in the general case, we can't hope to generalize Conjecture \ref{conjecture} further than removing the dependence on polarization $A$ by using the universal monoid $Q$ (which would give an analogue of Proposition \ref{proofuniversal}). 

If $\bar{\mathfrak{X}} \to \mathcal{S}$ is a special toric degeneration satisfying Assumption \ref{naturalcase},  $\left(\bar B, \bar{\mathscr{P}}\right)$ is simple. In this case, Theorem \ref{compgen} reduces to Corollary \ref{compgensimple}, so we need to generalize Corollary \ref{compgensimple}.

\begin{conj} \label{conjecturehd}
Let $\bar{\mathfrak{X}} \to \mathcal{S}$ be a special toric degeneration of relative dimension $n \geq 3$ satisfying Assumption \ref{naturalcase} and Assumption \ref{facelocalmodelass} (if $n \geq 4$), and supporting some polarization $A$. We further require that $A_1(\mathfrak{X}_0, \mathbb{Z}) = A_1(\mathfrak{X}_0, \mathbb{Z})_{\numer}$. Then (possibly after a finite basechange) there exists a strongly admissible resolution $\pi: \mathfrak{X} \to \bar{\mathfrak{X}}$ to a minimal log CY degeneration $\mathfrak{X} \to \mathcal{S}$. Suppose that Conjecture \ref{combconjecture} holds for $\bar{\mathfrak{X}} \to \mathcal{S}$ (we have proved this for $n=3$) and Conjectures \ref{canonicalscathd} and \ref{scatdiaghigherpowershd} hold for $\pi: \mathfrak{X} \to \bar{\mathfrak{X}}$. Let $P$ be a well-chosen monoid and let $J:= P \setminus K$ be the complement of the face containing the curve classes contracted by $\pi$ as before. Arguing as in Section \ref{basemonoidcompatible}, we may set $P:= K^{\gp} \oplus Q$ for $Q$ the universal monoid of \cite[Appendix A.2]{GHS} and the intrinsic mirror $$\check{\mathfrak{X}} \to \Spec \widehat{\kk [P]}_J = \Spec \kk [K^{\gp}] \llbracket Q \rrbracket$$ is well-defined.
\begin{enumerate}
    \item Let $\tau \in \bar{\mathscr{P}}^{[n-1]}$ and let $k \leq n-1$ be as in Conjecture \ref{combconjecture}. Choose a fixed $(n-1-k, k)$-quplex subdividing $\tau \in \bar{\mathscr{P}}^{[n-1]}$ (unique if $k=n-1$) and let $\delta'_{\tau} \in \mathscr{P}^{[n-1]}$ be a fixed cell subdividing the $(n-1-k, k)$-quplex. Then the $f_{\mathfrak{b}_{\delta'_{\tau}}}$ of Conjecture \ref{canonicalscathd}(2) is of the form $$f_{\mathfrak{b}_{\delta'_{\tau}}} = 1 + \sum_{m \in \bar{\Delta} (\tau, v) \setminus \{0 \}} t^{E_{\tau,m}}z^m$$ for some curve classes $E_{\tau,m} \in K^{\gp}$. Moreover, $K^{\gp}$ splits as $K^{\gp} = E^{\gp} \oplus G^{\gp}$ for $$E^{\gp} := \left\langle E_{\tau, m} \ , \ -E_{\tau,m} \ | \  \tau \in \bar{\mathscr{P}}^{[n-1]}, \ m \in \bar{\Delta} (\tau, v) \setminus \{0 \}  \right\rangle$$ and some finitely generated free abelian group $G^{\gp}$. 

    \item Define $h^{\GS}: K^{\gp} \oplus Q = E^{\gp} \oplus G^{\gp} \oplus Q \to E^{\gp} \oplus Q$ by sending $G^{\gp} \to 0$ and by the identity on $E^{\gp}$ and $Q$. We call the basechange $\check{\mathfrak{X}} \to \Spec \kk [E^{\gp}]\llbracket Q \rrbracket$ of the extended intrinsic mirror $\check{\mathfrak{X}} \to \Spec \kk [K^{\gp}] \llbracket Q \rrbracket$ by $h^{\GS}$ the (numerical) minimal relative Gross-Siebert locus. Then $\check{\mathfrak{X}} \to \Spec \kk [E^{\gp}]\llbracket Q \rrbracket$ is independent of the choice of $\delta'_{\tau} \in \mathscr{P}^{[n-1]}$ in (1).

    \item The correspondence of Corollary \ref{compgensimple} directly generalizes to a correspondence between $\check{\mathfrak{X}} \to \Spec \kk [t^{\pm E_{m} }]\llbracket Q \rrbracket$ and the subfamily $\check{\bar{\mathfrak{X}}}'_{\mathbb{P}} \to \Spec \kk'_{\mathbb{P}}\llbracket Q \rrbracket$ of the universal toric degeneration mirror $\check{\bar{\mathfrak{X}}}_{\mathbb{P}} \to \Spec \kk_{\mathbb{P}}\llbracket Q \rrbracket$ of \cite[Theorem A.2.4]{GHS}. Here $\kk_{\mathbb{P}}$ is the universal ring parametrizing projective gluing data and $\kk'_{\mathbb{P}} \subseteq \kk_{\mathbb{P}}$ is the subring generated by $\bar{s}$ for $s$ as follows. For every $\tau \in \bar{\mathscr{P}}^{[n-1]}$, let $\underline{\tau} \in \tilde{\bar{\mathscr{P}}}^{[n-1]}$ be the slab of $\tau$ as in Conjecture \ref{canonicalscathd}(3) for $\tau' = \delta'_{\tau}$ (for any choice of $\delta'_{\tau} \in \mathscr{P}^{[n-1]}$ in (1)). Choose a set $m_{\tau, i}, \ 1 \leq i \leq n-1$ of primitive generators of $\Lambda_{\tau}$ such that $m_{\tau, i}, \ 1 \leq i \leq \dim \bar{\Delta} (\tau, v)$ are the generators of $\bar{\Delta} (\tau, v)$. Then $s$ is given by, for $\sigma \in \bar{\mathscr{P}}^{\max}$ adjacent to $\tau \in \bar{\mathscr{P}}^{[n-1]}$, the homomorphisms:
\begin{equation} 
\begin{aligned}
&s_{\sigma \underline{\tau}}: \Lambda_{\sigma} \to \kk^{\times}, \quad \begin{cases} \begin{aligned} m_{\tau, i} &\mapsto c_{\rho, i}^{-1}, \ &1 \leq i \leq \dim \bar{\Delta} (\tau, v)  \\ m_{\tau, i} &\mapsto 1, \ &\dim \Delta (\tau, v) < i \leq n-1 \\ \xi &\mapsto 1 \end{aligned} \end{cases} \\ &s_{\sigma \underline{\tau}'}: \Lambda_{\sigma} \to \kk^{\times}, \quad \ \ \ \ m \mapsto 1 \end{aligned}
\end{equation} 
where $\xi$ is the normal generator pointing into $\sigma$ as in \eqref{posstalk1}, $\underline{\tau}' \subseteq \tau \in \bar{\mathscr{P}}^{[n-1]}$ is any slab of $\tau \in \bar{\mathscr{P}}^{[n-1]}$ distinct from $\underline{\tau}$, and $c_{\rho, i} \in \kk^{\times}$ for $1 \leq i \leq \dim \bar{\Delta} (\tau, v)$ are fixed constants. 
\end{enumerate}
\end{conj}

We expect that the proof of Conjecture \ref{canonicalscathd} will require an analysis of the curve classes contracted by $\pi$ and yield explicit expressions for the $f_{\mathfrak{b}_{\tau'}}$ of \eqref{infproductwallshd}. In particular, the analysis will imply conditions (1) and (2) of Conjecture \ref{conjecturehd}. The extension of Conjecture \ref{conjecture} to Conjecture \ref{conjecturehd}(3)  should then be a straightforward generalization of the arguments of Chapter \ref{gslocus}.

\appendix

\section[Toric degeneration mirror to a Batyrev degeneration in $\mathbb{P}^n$]{Toric degeneration mirror to a Batyrev degeneration in \texorpdfstring{$\mathbb{P}^n$}{Pn}} \label{projcase}

We discuss how to construct the toric degeneration mirror $\check{\bar{\mathfrak{X}}}_{\Delta, \TD}$ to the Batyrev degeneration $$\bar{\mathfrak{X}}_{\Delta}:=\left\{tf_{n+1}+\prod_{i=0}^n x_i\right\} \subseteq \mathbb{P}^n \times \Spec \kk\llbracket t \rrbracket$$ (in particular, it is a toric degeneration) with the natural projection to $\Spec \kk\llbracket t \rrbracket$. Here $f_{n+1}$ is a general homogeneous polynomial of degree $n+1$ and $\Delta$ is the polytope defining $\mathbb{P}^n, \ n \geq 3$. Following Chapter \ref{ell}, we will use the polarization by $-K_{\mathbb{P}_{\Delta}}$ to construct the toric degeneration mirror. 

We will give a direct generalization of Proposition \ref{ldr}(1) (the case of a Batyrev degeneration of elliptic curves) to $n = 3$ (the case of a Batyrev degeneration of K3-s), providing a direct connection between the toric degeneration mirror $\check{\bar{\mathfrak{X}}}_{\Delta, \TD}$ and the dual Batyrev degeneration $\bar{\mathfrak{X}}^{\eta(t)}_{\Delta^*}$ of \eqref{gentoricbat}. For $n \geq 4$, one can't construct the toric degeneration mirror $\check{\bar{\mathfrak{X}}}_{\Delta, \TD}$ by referring to the reconstruction algorithm of \cite[Proposition 3.9]{GS1} (see Theorem \ref{reconstruction}) since the local rigidity condition is not satisfied. However, we conjecture that one can still perform the algorithm and that Proposition \ref{ldr}(1) generalizes to this setting. We work in the notations of Chapter \ref{ell} (except for having a bar in the notation for toric degenerations in line with the rest of the thesis). 

It is easy to see that the Cox coordinate description of $\mathfrak{X}_{\Delta_1^*}$ of Figure \ref{gsequations} generalizes to arbitrary dimension. That is the degeneration $$  \bar{\mathfrak{X}}^{\eta(t)}_{\Delta^{*}} = \left\{ \eta(t)s+s_0 =0 \right\} \subseteq \mathbb{P}_{\Delta^{*}} \times \Spec \kk\llbracket t \rrbracket$$ of \eqref{gentoricbat} has the Cox coordinate description of Figure \ref{coxPn} (with $n-1$ copies of $\mathbb{Z}_{n+1}$). 

\begin{figure}[t]
$$
\begin{aligned}
&\begin{bmatrix}
& x_0 & x_1 & x_2 & x_3 & x_4 & \dots & x_n  \\ \hline
\mathbb{Z} & 1 & 1 & 1 &  1 &  \dots & 1  & 1 \\ \mathbb{Z}_{n+1} & 0 & 1 & 2 & 1 & \dots & 1   & 1 \\
\mathbb{Z}_{n+1} & 0 & 1 & 1 & 2 & \dots & 1   & 1 \\
\vdots & \vdots & \vdots & \vdots & \vdots & \ddots & \vdots   & \vdots \\
\mathbb{Z}_{n+1} & 0 & 1 & 1 & 1 & \dots & 2   & 1 \\
\mathbb{Z}_{n+1} & 0 & 1 & 1 & 1 & \dots & 1   & 2 \\
\end{bmatrix} \\  &s_0=\prod_{i=0}^n x_i \\ &s=\sum_{i=0}^n x_i^{n+1}
\end{aligned}
$$
\caption{Cox coordinate description of $\bar{\mathfrak{X}}^{\eta(t)}_{\Delta^{*}}$.} \label{coxPn}
\end{figure}

Constructing the toric degeneration mirror $\check{\bar{\mathfrak{X}}}_{\Delta, \TD}$ is more challenging. The dual intersection complex $\left(\bar{B}, \bar{\mathscr{P}}\right)$ of $\bar{\mathfrak{X}}_{\Delta} \to \Spec \kk \llbracket t \rrbracket$ can be identified with $\partial \Delta^*$ with the natural polyhedral decomposition and integral structure. Let the singular locus $\Gamma$\footnote{In Section \ref{tdaffstr}, we denote the singular locus by $\bar{\Delta}$.} be the union of all simplices in the barycentric\footnote{Unlike Section \ref{tdaffstr}, we don't require that $\Gamma$ contains no rational point. By the discussion of \cite[Appendix A.1]{GHS}, the barycentric choice gives rise to an isomorphic family.} subdivision $\tilde{\bar{\mathscr{P}}}$ of $\bar{\mathscr{P}}$ not containing a vertex or intersecting the interior of a maximal cell. Then $\bar{B} \setminus \Gamma$ carries an affine structure as usual (see \cite[Definition 2.10]{G2}).

We first compute the monodromy around the components of $\Gamma$. Let $v$ and $v'$ be two vertices of $\bar{B} \cong \partial \Delta^{*}$ and let $\sigma$ and $\sigma'$ be two maximal cells containing $v$ and $v'$ respectively. Let $\gamma$ be a simple loop based at $v$, passing successively into $\Int \sigma$, through $v'$, into $\Int \sigma'$, and back to $v$. We have $\Delta^* \subseteq N \cong \mathbb{Z}^{n}$. By choosing an integral basis $(e_1, \dots, e_n)$ of $N$, we may assume that $v=e_1, \ v' = e_2$, and 
$$
\sigma = \left\langle e_1, \dots, e_n \right\rangle, \quad
\sigma' = \left\langle e_1, \dots, e_{n-1}, -e_1 - \dots - e_n  \right\rangle.
$$
Now the MPA function $\varphi_{-K_{\mathbb{P}_{\Delta}}}$ has value $1$ at all the vertices of $\bar{B} \cong \partial \Delta^{*}$ and a calculation shows that it is defined by $k = e_1 + \dots + e_n \in N^*$ on $\sigma$ and $k'=e_1 + \dots + e_{n-1} -n e_n \in N^*$ on $\sigma'$. Identifying $\Lambda_v$ with $N / \langle v  \rangle \cong (\bar{e}_2, \dots, \bar{e}_n)$, by  \cite[Proposition 2.13]{G2} the monodromy $T_{\gamma}: \Lambda_v \to \Lambda_v$ takes form (compare with \eqref{monodromy}):
\begin{equation} \label{monodromyt}
T_{\gamma}(m) = m + \left\langle k'-k, m \right\rangle (v'-v) = m - \left\langle (n+1)\bar{e}_n, m \right\rangle \bar{e}_2.
\end{equation} So the monodromy matrix is the  $(n-1) \times (n-1)$ matrix $\Id - (n+1)I_{1, n-1}$ where $I_{1,n-1}$ is the matrix with $1$ in the $(1, n-1)$ position and zeroes elsewhere.

We would like to apply the recipe of \cite{GS1} to construct the toric degeneration mirror $\check{\bar{\mathfrak{X}}}_{\Delta, \TD} \to \Spec \kk \llbracket t \rrbracket$ from $\left(\bar{B}, \bar{\mathscr{P}}\right)$. Since $\left(\bar{B}, \bar{\mathscr{P}}\right)$ is not simple (see \cite[Definition 1.60]{GS2}), the existence of a toric log CY structure on $\check{\bar{X}}_0$ (the central fibre of $\check{\bar{\mathfrak{X}}}_{\Delta, \TD}$) is not automatic, and we need to specify\footnote{See the discussion of Section \ref{tdmpafunc}.} normalized initial slab functions $f_{\underline{\rho}}$ with no poles for every slab $\underline{\rho} \in \tilde{\bar{\mathscr{P}}}^{[n-1]}$ that satisfy for any $\underline{\rho}, \underline{\rho}' \subseteq \rho \in \bar{\mathscr{P}}^{[n-1]}$ the equation
\begin{equation} \label{compt}
f_{\underline{\rho}'} = z^{m_{\underline{\rho}'\underline{\rho} }}f_{\underline{\rho}}   
\end{equation}
as in \eqref{comp} where $m_{\underline{\rho}'\underline{\rho}}$ is the monodromy variable corresponding to $v'-v$ in \eqref{monodromyt}. Using the same choices as above, suppose that $v \subseteq \underline{\rho} \subseteq \sigma$. Then we may assume that $\Lambda_{\rho} \subseteq \Lambda_v$ is given by $(\bar{e}_3, \dots, \bar{e}_n)$ and let 
\begin{equation} \label{degeninprojslabfunc}
f_{\underline{\rho}} = 1+z^{(n+1)\bar{e}_3}+ \dots + z^{(n+1)\bar{e}_n}.
\end{equation}
Now \eqref{monodromyt} implies that this choice is consistent with \eqref{compt}. So $\check{\bar{X}}_0$ carries a toric log CY structure.

Suppose that $n = 3$ (i.e. $\bar{\mathfrak{X}}_{\Delta} \to \Spec \kk \llbracket t \rrbracket$ is the toric degeneration of Example \ref{runningexample}). Since $\left(\bar{B}, \bar{\mathscr{P}}\right)$ is of dimension $2$, $\check{\bar{X}}_0$ is locally rigid (see \cite[Definition 1.26]{GS1}). So by the reconstruction algorithm of \cite{GS1} (see Theorem \ref{reconstruction}) $\left(\bar{B}, \bar{\mathscr{P}}\right)$ gives rise to the toric degeneration mirror $\check{\bar{\mathfrak{X}}}_{\Delta, \TD} \to \Spec \kk\llbracket t \rrbracket$ along with a very ample line bundle $\mathcal{L}$ on $\check{\bar{\mathfrak{X}}}_{\Delta, \TD}$. Note that $\bar{B}(\mathbb{Z})$ consists of the vertices $v \in \bar{\mathscr{P}}^{[0]}$ that give rise to theta functions $\vartheta_0, \vartheta_1, \vartheta_2, \vartheta_3$, the canonically defined sections forming a basis of $\mathcal{L}$. We can generalize Proposition \ref{ldr}(1).

\begin{fact} \label{toricdegmirrorhddim3}
Let $n = 3$. There is a rational map $\check{\bar{\mathfrak{X}}}_{\Delta, \TD} \to \bar{\mathfrak{X}}^{\eta(t)}_{\Delta^*}$ given by $x_i \mapsto \vartheta_i, \ 0 \leq i \leq 3$ in the Cox coordinates (for a certain choice of $\eta(t) \in \kk\llbracket t \rrbracket$).
\end{fact}

\begin{proof}
Similarly to the proof of Proposition \ref{ldr}(1), it is enough to show that the toric degeneration mirror is given by
\begin{equation} \label{mirroreqform}
\check{\bar{\mathfrak{X}}}_{\Delta, \TD} = \left\{ \eta(t) (\vartheta_0^{4}+ \vartheta_1^{4} + \vartheta_2^{4} + \vartheta_3^{4}) +  \vartheta_0 \vartheta_1 \vartheta_2 \vartheta_3 = 0 \right\} \subseteq \mathbb{P}^3 \times \Spec \kk\llbracket t \rrbracket.
\end{equation}
We argue similarly to \cite[Example 6.0.3]{GHS} (on the discrete Legendre dual side), see there for more details. Namely, we use the symmetries of the dual intersection complex $\bar{B} \cong \partial \Delta^{*}$ of $\bar{\mathfrak{X}}_{\Delta}$ to restrict the form of the equation for the toric degeneration mirror. 

We can view the central fibre $\check{\bar{X}}_0$ of the mirror family (see Section \ref{gluedfamily} for the construction) as given by $\{\vartheta_0 \vartheta_1 \vartheta_2 \vartheta_3 = 0 \} \subseteq \mathbb{P}^3$. Now the fact that $\mathcal{L}$ is very ample implies that $\check{\bar{\mathfrak{X}}}_{\Delta, \TD}$ embeds into $\mathbb{P}_{\kk\llbracket t \rrbracket}^3$ as a hypersurface of degree $4$. So there is a single equation of degree $4$ satisfied by all the $\vartheta_i$ that restricts to $\vartheta_0 \vartheta_1 \vartheta_2 \vartheta_3 = 0$ at $t=0$.

We observe that $\bar{B} \cong \partial \Delta^{*}$ (and thus $\check{\bar{\mathfrak{X}}}_{\Delta, \TD}$) has a large symmetry group. First, $\bar{B} \cong \partial \Delta^{*}$ has an action given by permuting the vertices that preserves the initial structure. This action lifts to an action on $\mathcal{L}$ permuting the $\vartheta_i$. So the equation for $\check{\bar{\mathfrak{X}}}_{\Delta, \TD} \subseteq \mathbb{P}_{\kk\llbracket t \rrbracket}^3$ has to be a symmetric polynomial in $\vartheta_i$.

There is also an action of multiplication by $4$-th roots of unity. Indeed, note that the local system $\Lambda \otimes_{\mathbb{Z}} \mathbb{Z}_{4}$ has no monodromy, so it is trivial. Fix an isomorphism between $\Lambda \otimes_{\mathbb{Z}} \mathbb{Z}_4$ and a constant sheaf with stalk $(\mathbb{Z}_{4})^2$ by choosing an isomorphism $\Lambda_x \cong \mathbb{Z}^{2}$ at any $x \in \bar{B} \setminus \Gamma$. Now any character $\chi: (\mathbb{Z}_4)^{2} \to \kk^{\times}$ gives rise to a well-defined action on monomials via $$\chi: \mathcal{P} \to \Lambda \to \Lambda \otimes_{\mathbb{Z}} \mathbb{Z}_{4} \cong (\mathbb{Z}_{4})^{2} \to \kk^{\times}$$ (here $\mathcal{P}$ is the local system of \eqref{monomialexactseq} defining the monomials). Because of the form of the initial walls of the algorithmic scattering diagram $\bar{\mathfrak{D}}$, the scattering diagram is left invariant under this action, so $\chi$ acts on $\check{\bar{\mathfrak{X}}}_{\Delta, \TD}$. One can lift this action to an action on $\mathcal{L}$ in such a way that $\chi(\vartheta_0) = \vartheta_0$. 
Now the action on the other $\vartheta_i$ is determined. Let $v_i$ be the vertex corresponding to $\vartheta_i$. Then using the local fan structure at $v_0$ we can choose the coordinates $(\bar{e}_2, \bar{e}_3)$ as before so that $$v_1 = \bar{e}_{2}, \quad v_2 = \bar{e}_3, \quad v_3 = -\bar{e}_2-\bar{e}_3.$$
The lifted action of $\chi$ is given by (viewing $(\bar{e}_2, \bar{e}_3)$ as the basis for $(\mathbb{Z}_{4})^{2}$)
$$\chi(\vartheta_0) = \vartheta_0, \quad \chi(\vartheta_1) = \chi(\bar{e}_{2}) \cdot  \vartheta_1, \quad \chi(\vartheta_2) = \chi(\bar{e}_{3}) \cdot \vartheta_2,  \quad \chi(\vartheta_3) = \chi(-\bar{e}_2- \bar{e}_3) \cdot \vartheta_3.$$

A degree $4$ symmetric polynomial equation in $\vartheta_i$ that restricts to $\vartheta_0 \vartheta_1 \vartheta_2 \vartheta_3= 0$ at $t=0$ and is invariant under the action of $\Hom\left((\mathbb{Z}_{4})^{2}, \kk^{\times}\right)$ described above is necessarily of the form $$\eta(t) (\vartheta_0^{4}+ \vartheta_1^{4} + \vartheta_2^{4} + \vartheta_3^{4}) +  \vartheta_0 \vartheta_1 \vartheta_2 \vartheta_3 = 0$$ for some $\eta(t) \in \kk\llbracket t \rrbracket$.
\end{proof}

The correspondence of Theorem \ref{compgen} allows restricting the equation of the (extended) intrinsic mirror to a small resolution of $\mathfrak{\bar{X}}_{\Delta}$ in the case $n=3$.

\begin{observe} \label{imsmireq}
Note from \eqref{degeninprojslabfunc} that $\check{\bar{\mathfrak{X}}}_{\Delta, \TD}$ in Proposition \ref{toricdegmirrorhddim3} is the toric degeneration mirror to the toric degeneration of Example \ref{runningexample} constructed using the initial slab functions \begin{equation} 
f_{\underline{\rho}} = 1+w_{\rho}^4, \quad f_{\underline{\rho}'} = 1+w_{\rho}^{-4},
\end{equation}
where $w_{\rho} = z^{m_{\rho}}$ for $m_{\rho}$ an integral generator of $\Lambda_{\rho}$. It is easy to see that the argument in the proof of Proposition \ref{toricdegmirrorhddim3} directly generalizes to the universal setting of Section \ref{compgencons}. In particular, for the universal choice of slab functions as in \eqref{slabfuncgd}, the action given by permuting the vertices preserves the initial structure. Therefore, the universal toric degeneration mirror $\check{\bar{\mathfrak{X}}}^{s^r}_{\bar{\mathfrak{D}}^r} \to \Spec \kk [a_{\rho, i}, c_{\rho}^{\pm 1}]\llbracket Q \rrbracket$ of Construction \ref{compgen1} is given by an equation similar to \eqref{mirroreqform} of the form $$\left\{ \eta (\vartheta_0^{4}+ \vartheta_1^{4} + \vartheta_2^{4} + \vartheta_3^{4}) +  \vartheta_0 \vartheta_1 \vartheta_2 \vartheta_3 = 0 \right\} \subseteq \mathbb{P}^3 \times \Spec \kk [a_{\rho, i}, c_{\rho}^{\pm 1}]\llbracket Q \rrbracket$$ for some $\eta \in \kk [a_{\rho, i}, c_{\rho}^{\pm 1}]\llbracket Q \rrbracket$. Theorem \ref{compgen} implies that the minimal relative Gross-Siebert locus of the extended intrinsic mirror to a small resolution of the toric degeneration of Example \ref{runningexample} (see Sections \ref{problem} and \ref{resolvesmall}) is also given by an equation similar to \eqref{mirroreqform} of the form $$\left\{ \eta (\vartheta_0^{4}+ \vartheta_1^{4} + \vartheta_2^{4} + \vartheta_3^{4}) +  \vartheta_0 \vartheta_1 \vartheta_2 \vartheta_3 = 0 \right\} \subseteq \mathbb{P}^3 \times \Spec \kk [t^{\pm E_{\rho, k} }]\llbracket Q \rrbracket$$ for some $\eta \in \kk [t^{\pm E_{\rho, k} }]\llbracket Q \rrbracket$. 

\begin{sloppypar}
It follows from \eqref{canonicalsplitting} that the intrinsic mirror is actually defined over $\Spec \kk [t^{E_{\rho, k} }]\llbracket Q \rrbracket$ in this case since the resolution is small. The intrinsic mirror over $\Spec \kk [t^{E_{\rho, k} }]\llbracket Q \rrbracket$ is given by the same equation as it is given by it generically on the dense torus $\Spec \kk [t^{\pm E_{\rho, k} }]\llbracket Q \rrbracket \subseteq \Spec \kk [t^{E_{\rho, k} }]\llbracket Q \rrbracket$.
\end{sloppypar}
\end{observe}

Suppose now that $n \geq 4$. In this case, it is easy to check that $\check{\bar{X}}_0$ is not locally rigid. Indeed, the monodromy computation of \eqref{monodromyt} implies that $\check{\bar{X}}_0$ does not satisfy condition (i) of \cite[Definition 1.26]{GS1} for $n \geq 4$. It also does not satisfy condition (ii) of \cite[Definition 1.26]{GS1} for $n=4$ (but satisfies it for $n \geq 5$). We conjecture that Proposition \ref{ldr}(1) can still be generalized to this setting. 

\begin{conj} \label{conjprojhd}
Let $n \geq 4$. One can run the reconstruction algorithm of \cite{GS1} (see Theorem \ref{reconstruction}) using the initial slab functions of \eqref{degeninprojslabfunc} to obtain a toric degeneration mirror $\check{\bar{\mathfrak{X}}}_{\Delta, \TD} \to \Spec \kk\llbracket t \rrbracket$. Consequently, $\check{\bar{\mathfrak{X}}}_{\Delta, \TD}$ is given by $$\check{\bar{\mathfrak{X}}}_{\Delta, \TD} = \left\{ \eta(t) \sum_{i=0}^n \vartheta_i^{n+1} + \prod_{i=0}^n \vartheta_i = 0 \right\} \subseteq \mathbb{P}^n \times \Spec \kk\llbracket t \rrbracket$$ for some $\eta(t) \in \kk\llbracket t \rrbracket$ and there is a rational map $\check{\bar{\mathfrak{X}}}_{\Delta, \TD} \to \bar{\mathfrak{X}}^{\eta(t)}_{\Delta^*}$ given by $x_i \mapsto \vartheta_i$ in the Cox coordinates. 
\end{conj}

Provided that one can construct $\check{\bar{\mathfrak{X}}}_{\Delta, \TD} \to \Spec \kk\llbracket t \rrbracket$, the second part of Conjecture \ref{conjprojhd} follows as in the proof of Proposition \ref{toricdegmirrorhddim3} by using symmetries of the dual intersection complex $\bar{B} \cong \partial \Delta^{*}$. Local rigidity guarantees that one always has unique choices in the reconstruction algorithm of \cite{GS1}. We do not expect this to be the case here. Rather, we expect that one can always make natural choices yielding $\check{\bar{\mathfrak{X}}}_{\Delta, \TD} \to \Spec \kk\llbracket t \rrbracket$. Note that the map $\check{\bar{\mathfrak{X}}}_{\Delta, \TD} \to \bar{\mathfrak{X}}^{\eta(t)}_{\Delta^*}$ recovers the classical construction of Greene and Plesser \cite{GP} for $n = 4$.

\begin{remark}
Even assuming that Conjecture \ref{conjprojhd} holds, we can't directly generalize the argument of Observation \ref{imsmireq} to $n \geq 4$ using Conjecture \ref{conjecturehd} since the dual intersection complex $\left(\bar B, \bar{\mathscr{P}}\right)$ of $\bar{\mathfrak{X}}_{\Delta}$ is not simple. However, in this particular case, one should be able to check that the space of toric log CY spaces with intersection complex $\left(\bar B, \bar{\mathscr{P}}\right)$ is well-behaved, state the analogue of Theorem \ref{compgen} in higher dimensions, and do the same argument.
\end{remark}

\section{Log structures for canonical families} \label{log}

Suppose that we are in the setup of \cite{GHS} (see Section \ref{scatteringsetup}) working over a ring $A$, using monoid $Q$ and a monoid ideal $I_0 \subseteq Q$. Fix an ideal $I$ with $\sqrt{I} = I_0$ and suppose that $\mathfrak{D}_I$ is a consistent scattering diagram on $(B, \mathscr{P})$. Then $\mathfrak{D}_I$ gives rise to a family $\check{\mathfrak{X}}_{\mathfrak{D}_I} \to \Spec A[Q]/I$ by \cite[Theorem 4.3.2]{GHS}. We are going to construct log structures on $\check{\mathfrak{X}}_{\mathfrak{D}_I}$ and $\Spec A[Q]/I$, and an enhancement of $\check{\mathfrak{X}}_{\mathfrak{D}_I} \to \Spec A[Q]/I$ to a log morphism log smooth away from codimension $2$. Since both toric degeneration mirrors of Section \ref{tdsetup} and intrinsic mirrors of Section \ref{intsetup} fall into this framework, this will make the toric degeneration and intrinsic mirrors log smooth morphisms away from codimension $2$.

We endow $\Spec A[Q]/I$ with the canonical log structure defined by the global chart $Q \to A[Q]/I, \ q \mapsto z^q $. It is enough to give a log structure $\mathscr{M}_{\check{\mathfrak{X}}^o_{\mathfrak{D}_I}}$ on $\check{\mathfrak{X}}^o_{\mathfrak{D}_I}$ and a log morphism $(\check{\mathfrak{X}}^o_{\mathfrak{D}_I}, \mathscr{M}_{\check{\mathfrak{X}}^o_{\mathfrak{D}_I}}) \to \Spec A[Q]/I$ (with the underlying morphism of schemes the usual $\check{\mathfrak{X}}^o_{\mathfrak{D}_I} \to \Spec A[Q]/I$) that is log smooth away from codimension $2$. Indeed, $\check{\mathfrak{X}}_{\mathfrak{D}_I} \setminus \check{\mathfrak{X}}^o_{\mathfrak{D}_I}$ is codimension $2$ in $\check{\mathfrak{X}}_{\mathfrak{D}_I}$ and given a log structure $\mathscr{M}_{\check{\mathfrak{X}}^o_{\mathfrak{D}_I}}$ on $\check{\mathfrak{X}}^o_{\mathfrak{D}_I}$ we can define the log structure $\mathscr{M}_{\check{\mathfrak{X}}_{\mathfrak{D}_I}}$ on $\check{\mathfrak{X}}_{\mathfrak{D}_I}$ as $\mathscr{M}_{\check{\mathfrak{X}}_{\mathfrak{D}_I}}:=j_{*} \mathscr{M}_{\check{\mathfrak{X}}^o_{\mathfrak{D}_I}}$ for $j: \check{\mathfrak{X}}^o_{\mathfrak{D}_I} \hookrightarrow \check{\mathfrak{X}}_{\mathfrak{D}_I}$ the canonical inclusion of \cite[Theorem 4.3.2(c)]{GHS} (we have $j_{*} \mathcal{O}_{\check{\mathfrak{X}}^o_{\mathfrak{D}_I}} \cong \mathcal{O}_{\check{\mathfrak{X}}_{\mathfrak{D}_I}}$ since $\check{\mathfrak{X}}_{\mathfrak{D}_I}$ satisfies Serre's $S_2$ condition by \cite[Proposition 2.1.6]{GHS}). Then the log morphism $(\check{\mathfrak{X}}^o_{\mathfrak{D}_I}, \mathscr{M}_{\check{\mathfrak{X}}^o_{\mathfrak{D}_I}})  \to \Spec A[Q]/I$ extends to a log morphism $(\check{\mathfrak{X}}_{\mathfrak{D}_I}, \mathscr{M}_{\check{\mathfrak{X}}_{\mathfrak{D}_I}})  \to \Spec A[Q]/I$ that is log smooth away from the union of a codimension $2$ subset of $\check{\mathfrak{X}}^o_{\mathfrak{D}_I}$ and $\check{\mathfrak{X}}_{\mathfrak{D}_I} \setminus \check{\mathfrak{X}}^o_{\mathfrak{D}_I}$, which is a codimension $2$ subset of $\check{\mathfrak{X}}_{\mathfrak{D}_I}$. By the same argument, it is enough to give a log structure $\mathscr{M}_{\check{\mathfrak{X}}^{o'}_{\mathfrak{D}_I}}$ on a subset $\check{\mathfrak{X}}^{o'}_{\mathfrak{D}_I}$ of $\check{\mathfrak{X}}^{o}_{\mathfrak{D}_I}$ with $\check{\mathfrak{X}}^{o}_{\mathfrak{D}_I} \setminus \check{\mathfrak{X}}^{o'}_{\mathfrak{D}_I}$ of codimension $2$ and a log smooth morphism $(\check{\mathfrak{X}}^{o'}_{\mathfrak{D}_I}, \mathscr{M}_{\check{\mathfrak{X}}^{o'}_{\mathfrak{D}_I}})  \to \Spec A[Q]/I$ (with the underlying morphism of schemes the restriction of $\check{\mathfrak{X}}^o_{\mathfrak{D}_I} \to \Spec A[Q]/I$ to $\check{\mathfrak{X}}^{o'}_{\mathfrak{D}_I}$).

From the proof of \cite[Proposition 2.4.1]{GHS}, $\mathfrak{X}^o_{\mathfrak{D}_I}$ is constructed by gluing together $\Spec R_{\mathfrak{b}}$ along $\Spec R_{\mathfrak{u}}$ for choices of slabs $\mathfrak{b} \subseteq \underline{\rho} \in \tilde{\mathscr{P}}^{[n-1]}$ for every $\rho \in \mathscr{P}^{[n-1]}$ and choices of chambers $\mathfrak{u} \subseteq \sigma \in \mathscr{P}^{\max}$ for every $\sigma \in \mathscr{P}^{\max}$. Fix a slab $\mathfrak{b}$ and consider the subset $$ \left\{ Z_+ = Z_-=f_{\mathfrak{b}} = 0\right\} \subseteq \Spec (A[Q]/I)[\Lambda_{\rho}][Z_+, Z_-]/(Z_+Z_- - f_{\mathfrak{b}} \cdot z^{\kappa_{\rho}}) = \Spec R_{\mathfrak{b}}.$$ It defines a subset of $\mathfrak{X}^o_{\mathfrak{D}_I}$ of codimension $2$ via the canonical inclusion $\Spec R_{\mathfrak{b}} \hookrightarrow \mathfrak{X}^o_{\mathfrak{D}_I}$ (see \cite[Proposition 2.4.1]{GHS}) that we still denote by $\left\{ Z_+ = Z_-=f_{\mathfrak{b}} = 0\right\}$. Let $$\check{\mathfrak{X}}^{o'}_{\mathfrak{D}_I} := \check{\mathfrak{X}}^{o}_{\mathfrak{D}_I} \Big\backslash \bigcup_{\mathfrak{b} \subseteq \underline{\rho} \in \tilde{\mathscr{P}}^{[n-1]}} \left\{ Z_+ = Z_-=f_{\mathfrak{b}} = 0\right\}$$ and note that $\check{\mathfrak{X}}^{o'}_{\mathfrak{D}_I}$ can be constructed similarly to the construction of $\mathfrak{X}^o_{\mathfrak{D}_I}$ in \cite[Proposition 2.4.1]{GHS} by using appropriate localizations of $R_{\mathfrak{b}}$ and $R_{\mathfrak{u}}$. We will construct a fine saturated log structure on $\check{\mathfrak{X}}^{o'}_{\mathfrak{D}_I}$ and a log smooth morphism to $\Spec A[Q]/I$.

Away from the subsets $\left\{ Z_+ = Z_- = 0\right\} \cap \check{\mathfrak{X}}^{o'}_{\mathfrak{D}_I} \subseteq \Spec R_{\mathfrak{b}}$, the fine saturated log structure is induced by the log structures on $\Spec R_{\mathfrak{u}}$ defined via the global charts $h_{\mathfrak{u}}: Q \to R_{\mathfrak{u}}, \ q \mapsto z^q$. Moreover, for every $\Spec R_{\mathfrak{u}}$, we can define a log morphism $\Spec R_{\mathfrak{u}} \to \Spec A[Q]/I$ (with the underlying morphism of schemes defined by the natural inclusion $\frac{A [P]}{J} \to R_{\mathfrak{u}}$) via the global chart $\Id: Q \to Q, \ q \mapsto q$ that fits into the commutative diagram:
\begin{equation} \label{commdiagtriviality}
\begin{tikzcd}
{Q} \arrow{r}{h_{\mathfrak{u}}} \arrow[leftarrow]{d}{\Id} & {R_{\mathfrak{u}}}  \arrow[leftarrow]{d}{} \\
Q \arrow{r}{} & \frac{A [Q]}{I}
\end{tikzcd}
\end{equation}
The log morphism $\Spec R_{\mathfrak{u}} \to \Spec A[Q]/I$ is log smooth by Kato's criterion (see \cite[Theorem (3.5)]{Kato}) since $\Id: Q \to Q$ is injective.

In the neighbourhood of a subset $\left\{ Z_+ = Z_- = 0\right\} \cap \check{\mathfrak{X}}^{o'}_{\mathfrak{D}_I} \subseteq \Spec R_{\mathfrak{b}}$, the fine saturated log structure is induced by the log structure on $\Spec (R_{\mathfrak{b}})_{f_{\mathfrak{b}}}$ defined as follows. Let $P := Q \oplus_{\mathbb{N}} \mathbb{N}^2$ be the pushout via the maps $\mathbb{N} \to Q, \ 1 \mapsto \kappa_{\rho}$ and $\mathbb{N} \to \mathbb{N}^2, \ 1 \mapsto (1,1)$.  We will write the elements of $P$ as triplets $[q, (a, b)]$ for $q \in Q, \ (a, b) \in \mathbb{N}^2$ (some of them correspond to the same equivalence classes in $P$). The global chart $h_{\mathfrak{b}}: P \to (R_{\mathfrak{b}})_{f_{\mathfrak{b}}}, \ [q, (a,b)] \mapsto z^q Z_{+}^a Z_{-}^b f_{\mathfrak{b}}^{-b}$ defines a fine saturated (since $P$ is a finitely generated and saturated monoid) log structure on $\Spec (R_{\mathfrak{b}})_{f_{\mathfrak{b}}}$ (one needs to check that $h_{\mathfrak{b}}$ is well-defined on the equivalence classes of triplets $[q, (a, b)]$). Moreover, for every $\Spec (R_{\mathfrak{b}})_{f_{\mathfrak{b}}}$, we can define a log morphism $\Spec (R_{\mathfrak{b}})_{f_{\mathfrak{b}}} \to \Spec A[Q]/I$ (with the underlying morphism of schemes defined by the map $\frac{A [P]}{J} \to (R_{\mathfrak{b}})_{f_{\mathfrak{b}}}$ that is the natural inclusion in $R_{\mathfrak{b}}$ followed by localization) via the global chart $\pi: Q \to P, \ q \mapsto [q, (0, 0)]$ that fits into the commutative diagram:
\begin{equation} \label{commdiag}
\begin{tikzcd}
{P} \arrow{r}{h_{\mathfrak{b}}} \arrow[leftarrow]{d}{\pi} & {(R_{\mathfrak{b}})_{f_{\mathfrak{b}}}}  \arrow[leftarrow]{d}{} \\
Q \arrow{r}{} & \frac{A [Q]}{I}
\end{tikzcd}
\end{equation}
Again, the log morphism $\Spec (R_{\mathfrak{b}})_{f_{\mathfrak{b}}} \to \Spec A[Q]/I$ is log smooth by Kato's criterion since $\pi: Q \to P$ is injective.

To define a fine saturated log structure on $\check{\mathfrak{X}}^{o'}_{\mathfrak{D}_I}$ and a log smooth morphism to $\Spec A[Q]/I$ it is enough to check that the gluing in the proof of \cite[Proposition 2.4.1]{GHS} respects the log structures of \eqref{commdiagtriviality} and \eqref{commdiag}. Fix a chamber $\mathfrak{u} \subseteq \sigma \in \mathscr{P}^{\max}$ and a slab $\mathfrak{b} \subseteq \underline{\rho} \in \tilde{\mathscr{P}}^{[n-1]}$ with $\rho \subseteq \sigma$. We have a commutative diagram as follows:
\begin{equation} \label{commdiag2}
\begin{tikzcd}
{Q \oplus \mathbb{Z}} \arrow{r}{h'_{\mathfrak{u}}} \arrow[leftarrow]{d}{g_{\rho}} & {(R_{\mathfrak{u}})_{\chi_{\mathfrak{b}, \mathfrak{u}}(f_{\mathfrak{b}})}}  \arrow[leftarrow]{d}{\chi_{\mathfrak{b}, \mathfrak{u}}} \\
P=Q \oplus_{\mathbb{N}} \mathbb{N}^2 \arrow{r}{h_{\mathfrak{b}}} & {(R_{\mathfrak{b}})_{f_{\mathfrak{b}}}}
\end{tikzcd}
\end{equation}
Here $h_{\mathfrak{b}}$ is the chart of \eqref{commdiag}, $\chi_{\mathfrak{b}, \mathfrak{u}}$ is induced by the localization homomorphism of \eqref{lochom}, $g_{\rho}: Q \oplus_{\mathbb{N}} \mathbb{N}^2 \to Q \oplus \mathbb{Z}$ is given by $[q, (a,b)] \mapsto (q + b \kappa_{\rho}, a-b)$ (one needs to check that $g_{\rho}$ is well-defined on the equivalence classes of triplets $[q, (a, b)]$), and $h'_{\mathfrak{u}}: Q \oplus \mathbb{Z} \to (R_{\mathfrak{u}})_{\chi_{\mathfrak{b}, \mathfrak{u}}(f_{\mathfrak{b}})}$ is given by $(q, a) \mapsto z^q z^{a\xi}$. Moreover, the log structure on $\Spec (R_{\mathfrak{u}})_{\chi_{\mathfrak{b}, \mathfrak{u}}(f_{\mathfrak{b}})}$ defined by $h'_{\mathfrak{u}}$ agrees with the one induced by $h_{\mathfrak{u}}$ since $z^{a\xi}$ is invertible in $(R_{\mathfrak{u}})_{\chi_{\mathfrak{b}, \mathfrak{u}}(f_{\mathfrak{b}})}$. So the log structures of \eqref{commdiagtriviality} and \eqref{commdiag} are compatible and the gluing in the proof of \cite[Proposition 2.4.1]{GHS} defines a fine saturated log structure $\mathscr{M}_{\check{\mathfrak{X}}^{o'}_{\mathfrak{D}_I}}$ on $\check{\mathfrak{X}}^{o'}_{\mathfrak{D}_I}$ and a log morphism $(\check{\mathfrak{X}}^{o'}_{\mathfrak{D}_I}, \mathscr{M}_{\check{\mathfrak{X}}^{o'}_{\mathfrak{D}_I}})  \to \Spec A[Q]/I$. This morphism is log smooth by Kato's criterion since in the neighbourhood of every point of $\check{\mathfrak{X}}^{o'}_{\mathfrak{D}_I}$ we have a chart for the morphism of the form \eqref{commdiagtriviality} or of the form \eqref{commdiag}. By the discussion above, we have a canonical extension to a log morphism $(\check{\mathfrak{X}}_{\mathfrak{D}_I}, \mathscr{M}_{\check{\mathfrak{X}}_{\mathfrak{D}_I}})  \to \Spec A[Q]/I$ log smooth away from the subset $\check{\mathfrak{X}}_{\mathfrak{D}_I} \setminus \check{\mathfrak{X}}^{o'}_{\mathfrak{D}_I}$ of codimension $2$. Note that $\mathscr{M}_{\check{\mathfrak{X}}_{\mathfrak{D}_I}}$ is only fine saturated on  $\check{\mathfrak{X}}^{o'}_{\mathfrak{D}_I} \subseteq \check{\mathfrak{X}}_{\mathfrak{D}_I}$.

\begin{sloppypar}
Suppose that we have a collection of scattering diagrams $\mathfrak{D}^{I_0}:= \left\{ \mathfrak{D}_{I_0^{k+1}}, \ k \geq 0 \right\}$ such that $\mathfrak{D}_{I_0^{k}}$ is compatible with $\mathfrak{D}_{I_0^{k-1}}$ for $k \geq 1$. As usual, it gives rise to an inverse system of families $\check{\mathfrak{X}}_{\mathfrak{D}_{I_0^{k+1}}} \to \Spec A[Q]/I_0^{k+1}$ for $k \geq 0$. We define the log structure on the inverse limit 
$\check{\mathfrak{X}}_{\mathfrak{D}^{I_0}} \to \Spec \widehat{A[Q]}_{I_0}$ by taking the inverse limit in the category of log schemes. In particular, this defines log structures on the toric degeneration mirrors of \eqref{reconstructedfamily} and the (extended) intrinsic mirrors of \eqref{extendedmirrorscatint}.
\end{sloppypar}

\begin{remark}
In the presence of non-trivial gluing data $s$ on $(B, \mathscr{P})$ (see Section \ref{gluingdataintro} and \cite[Section 5.2]{GHS}) the construction of an enhancement of $\check{\mathfrak{X}}^s_{\mathfrak{D}_I} \to \Spec A[Q]/I$ to a log morphism log smooth away from codimension $2$ is similar, replacing the localization homomorphisms $\chi_{\mathfrak{b}, \mathfrak{u}}$ with the localization homomorphisms $\chi^s_{\mathfrak{b}, \mathfrak{u}}$ twisted by $s$ (see Section \ref{gluingdataintro}). 
\end{remark}

\begin{observe}
It is easy to check explicitly that with the above definition of log structures on canonical families, all the basechanges considered in Chapters \ref{proof}, \ref{gslocus}, and \ref{higherdimgen} are also basechanges in the category of log schemes. 
\end{observe}

\cleardoublepage

\end{document}